\documentclass[Afour,sageh,times]{sagej}

\usepackage{amsmath,amsfonts,amsthm}
\usepackage{amssymb}
\usepackage{prettyref}
\usepackage{refcount}
\usepackage{mathrsfs}
\usepackage{graphicx}
\usepackage{wrapfig}
\usepackage{subfig}
\usepackage{framed}
\usepackage{multirow}
\usepackage{booktabs}
\usepackage{tabularx}
\usepackage{algorithm}
\usepackage[table]{xcolor}
\usepackage[export]{adjustbox}
\usepackage[noend]{algpseudocode}
\usepackage{enumitem}
\usepackage{comment}
\usepackage{environ}
\usepackage[colorlinks,allcolors=gray,hypertexnames=true]{hyperref}
\newcommand{\citewithauthor}[1]{\citeauthor{#1} \cite{#1}}
\usepackage[nomessages]{fp}
\usepackage[overload]{empheq}
\usepackage{tikz}
\usepackage{pgfplots}
\pgfplotsset{compat=1.18}

\newtheorem{theorem}{Theorem}
\newtheorem{lemma}[theorem]{\TE{Lemma}}
\newtheorem{corollary}[theorem]{\TE{Corollary}}
\newtheorem{remark}{\TE{Remark}}
\newtheorem{assume}{\TE{Assumption}}
\newtheorem{definition}{\TE{Definition}}

\makeatletter
\newtheorem*{rep@theorem}{\rep@title}
\newcommand{\newreptheorem}[2]{%
  \newenvironment{rep#1}[1]{%
    \def\rep@title{#2~\ref{##1}}%
    \def\@currentlabel{\getrefnumber{##1}}%
    \begin{rep@theorem}}%
    {\end{rep@theorem}}}
\makeatother
\newreptheorem{theorem}{Theorem}
\newreptheorem{lemma}{Lemma}

\newenvironment{repalgorithm}[1]{%
  \renewcommand{\thealgorithm}{\getrefnumber{#1}}%
  \begin{algorithm}}%
 {\end{algorithm}%
  \addtocounter{algorithm}{-1}}

\algnewcommand{\LineComment}[1]{\State \(\triangleright\) #1}
\algdef{SE}[DOWHILE]{Do}{doWhile}{\algorithmicdo}[1]
{\algorithmicwhile\ #1}
\newcommand*{\colorboxed}{}
\def\colorboxed#1#{%
  \colorboxedAux{#1}%
}
\newcommand*{\colorboxedAux}[3]{%
  \begingroup
    \colorlet{cb@saved}{.}%
    \color#1{#2}%
    \boxed{%
      \color{cb@saved}%
      #3%
    }%
  \endgroup
}

\newrefformat{fig}{Figure~\ref{#1}}
\newrefformat{appen}{Appendix~\ref{#1}}
\newrefformat{par}{Section~\nameref{#1}}
\newrefformat{sec}{Section~\nameref{#1}}
\newrefformat{sub}{Section~\nameref{#1}}
\newrefformat{table}{Table~\ref{#1}}
\newrefformat{ass}{Assumption~\ref{#1}}
\newrefformat{alg}{Algorithm~\ref{#1}}
\newrefformat{def}{Definition~\ref{#1}}
\newrefformat{rem}{Remark~\ref{#1}}
\newrefformat{thm}{Theorem~\ref{#1}}
\newrefformat{lem}{Lemma~\ref{#1}}
\newrefformat{cor}{Corollary~\ref{#1}}
\newrefformat{step}{Step~\ref{#1}}
\newrefformat{ln}{Line~\ref{#1}}
\newrefformat{eq}{Equation~\ref{#1}}
\newrefformat{pb}{Problem~\ref{#1}}
\newrefformat{it}{Item~\ref{#1}}
\newrefformat{te}{Term~\ref{#1}}
\def\Eqref Eq:#1:{\eqref{eq:#1}}
\newrefformat{Eq}{Equation~\Eqref#1:}

\newcommand{\R}{\mathbb{R}}

\newcommand{\TE}[1]{\textbf{#1}}

\NewEnviron{ResizedEq}{%
  \begin{equation}
  \adjustbox{max width=0.85\hsize}{$\begin{aligned}\BODY\end{aligned}$}%
  \end{equation}
}

\newcommand{\FPP}[2]{\frac{\partial{#1}}{\partial{#2}}}
\newcommand{\FPPR}[2]{{\partial{#1}}/{\partial{#2}}}
\newcommand{\FPPT}[2]{\frac{\partial^2{#1}}{\partial{#2}^2}}

\newcommand{\FPPTT}[3]{\frac{\partial^2{#1}}{\partial{#2}\partial{#3}}}

\newcommand{\TWO}[2]{\left(\setlength{\arraycolsep}{1pt}\begin{array}{cc}{#1}, & {#2}\end{array}\right)}
\newcommand{\TWOC}[2]{\left(\setlength{\arraycolsep}{1pt}\begin{array}{c}#1 \\ #2\end{array}\right)}

\newcommand{\THREE}[3]{\left(\setlength{\arraycolsep}{1pt}\begin{array}{ccc}{#1}, & {#2}, & {#3}\end{array}\right)}

\newcommand{\dist}{\text{dist}}
\newcommand{\gph}{\operatorname{gph}}
\newcommand{\diag}{\text{diag}}

\newcommand{\argmin}[1]{\underset{#1}{\text{argmin}}\;}

\newcommand{\ST}{\text{s.t.}}

\newcommand{\TWORCell}[2]{\begin{tabular}{@{}l@{}}#1 \\ #2\end{tabular}}

\newcommand{\proofread}[1]{}
\newif\ifArxiv
\usepackage{xcolor}
\definecolor{Blue}{rgb}{0.2, 0.2, 0.8}
\definecolor{Black}{rgb}{0.0, 0.0, 0.0}
\colorlet{TableHeadColor}{gray!28}
\colorlet{TableRowColorA}{gray!12}
\colorlet{TableRowColorB}{white}
\newcommand{\TableRowColors}{\rowcolors{2}{TableRowColorA}{TableRowColorB}}
\newcommand{\TableHeaderRow}{\rowcolor{TableHeadColor}}
\newcommand{\hldiff}[1]{\textcolor{brown}{#1}}

\makeatletter
\renewenvironment{proof}[1][\proofname]{\par
  \pushQED{\qed}%
  \normalfont \topsep6\p@\@plus6\p@\relax
  \trivlist
  \item[\hskip\labelsep
        \bfseries
    #1\@addpunct{.}]\ignorespaces
}{%
  \popQED
  \endtrivlist\@endpefalse
}
\makeatother

\usepackage{thmtools}

\def\volumeyear{2025}

\newcommand{\LControlJacobianBound}{L_{\nabla W}}
\newcommand{\LHamilton}{H^{\mathrm{c}}}
\newcommand{\LPosition}{\theta^{\mathrm{c}}}
\newcommand{\LCurvature}{L_{\nabla^2\bar U}}
\newcommand{\LControlWorkCurvature}{L_{\nabla^2\bar W}}
\newcommand{\LPotentialForce}{L_{\nabla\bar U}}
\newcommand{\LResidual}{L_r}
\newcommand{\LControlWorkForce}{L_{\nabla\bar W}}
\newcommand{\LHamBound}{\mathcal{B}_{H}}
\newcommand{\LVelBound}{\mathcal{B}_{\dot{\theta}}}
\newcommand{\LPosBound}{\mathcal{B}_{\theta}}
\newcommand{\LResidualWork}{\bar{\LResidual}}
\newcommand{\LLambdaMin}{L_\lambda}             
\newcommand{\LHessianCond}{L_B}
\newcommand{\LObjGradLip}{L_{\nabla^2 O}}

\newcommand{\PerStepObj}{\bar{\Psi}}     
\newcommand{\StatRes}{r}                 
\newcommand{\StatJac}{J}                 
\newcommand{\StackConf}{\Theta}          
\newcommand{\StackCtrl}{\Upsilon}        
\newcommand{\StatJacU}{J_\StackCtrl}     
\newcommand{\SigmaMinLICQ}{\underline{\sigma}_{\min}}
\newcommand{\Damp}{\mathfrak{D}}

\newcommand{\EpsContactPadding}{\epsilon_{\mathrm{pad}}}
\newcommand{\EpsAtan}{\epsilon_{\mathrm{atan}}}
\newcommand{\EpsFJ}{\epsilon_{\mathrm{FJ}}}
\newcommand{\EpsFJAll}{\widehat\EpsFJ}
\newcommand{\EpsDamp}{\epsilon_{\mathrm{damp}}}
\newcommand{\EpsFri}{\epsilon_{\mathrm{fri}}}
\newcommand{\EpsCRDS}{\epsilon_{\mathrm{res}}}
\newcommand{\EpsStop}{\epsilon_{\mathrm{stop}}}         
\newcommand{\EpsStopNS}{\widetilde\EpsStop}  

\newcommand{\StatResG}{r_{\mathcal{G}}}      
\newcommand{\LResidualG}{L_{r,\mathcal{G}}}    
\newcommand{\LObjGradBound}{L_{\nabla O}}    

\begin{document}

\runninghead{Pan, Wang, \& Zhu}

\title{Global Convergence of an SQP Method for Contact-Implicit Trajectory Optimization}

\author{Zherong Pan\affilnum{1}, Fuquan Wang\affilnum{2}, and Yifan Zhu\affilnum{2}}

\affiliation{\affilnum{1}Meta Reality Labs\\
\affilnum{2}Department of Computer Science, University of Illinois Chicago}

\email{zrpan@meta.com}

\begin{abstract}
Contact-Implicit Trajectory Optimization (CITO) is a powerful framework for planning motions of robots that interact with complex environments, but its convergence behavior remains difficult to characterize. Existing formulations either rely on off-the-shelf nonlinear programming solvers whose guarantees require constraint qualifications that are hard to verify for contact-rich systems, or require differentiable explicit dynamics maps that are difficult to formulate in the presence of impacts, changing contact modes, and geometric non-smoothness. This paper studies CITO for a broader class of constraint-rich dynamic systems formulated through implicit physics constraints. By exploiting maximal-coordinate structure, penalty relaxations of equality and inequality constraints, and finite-horizon Hamiltonian bounds, we prove global convergence in the numerical-optimization sense for a specific line-search Sequential Quadratic Programming (SQP) method with adaptive timestep refinement. Under the stated constraint-rich dynamic-system assumptions and fixed pre-horizon data satisfying the strict-interior predecessor condition, the method terminates finitely from any finite discretized optimized trajectory guess---a finite nonempty ordered optimized set, finite optimized configurations at its nodes, and optimized controls in the compact control set---at an $\EpsFJAll$-feasible, force-parameterized stationary point---a unit-objective, penalty-induced force-parameterized Fritz--John certificate for smooth damping, weakening to a homogeneous force-parameterized Clarke--Fritz--John certificate whose objective multiplier may vanish under nonsmooth frictional contact. In both cases the constraint forces are selected as configuration-dependent maps induced by the penalty gradient rather than optimized as independent primal variables. This is a convergence-to-stationarity result for the specified algorithm, not a guarantee of global trajectory optimality and not a guarantee for generic SQP solvers. The SQP method is applied to a smooth modified trajectory NLP whose physics constraints are the stationarity equations of a penalized, smoothed, reduced, and mollified per-step surrogate; its KKT certificate is then transferred to the original CITO formulation. The analysis further establishes a bounded optimization history: the merit-penalty parameter and optimized trajectory remain norm-bounded throughout. These results provide a convergence guarantee for a scalable gradient-based contact planner while retaining compatibility with arbitrary polyhedral contact models. Preliminary experiments on locomotion and manipulation problems support the practical applicability of the proposed solver.
\end{abstract}


\maketitle
\footnotetext[0]{\sagesf All models and code were accessed and experiments were performed by the University of Illinois Chicago authors.}

\section{Introduction}
Nonlinear dynamic systems arise throughout robotics, a central example being articulated robots interacting with their environments~\cite{todorov2012mujoco}. Planning motions for such systems is a fundamental problem in motion planning and control. Broadly, existing algorithms fall into two classes: sampling-based planners with global asymptotic guarantees and gradient-based planners that trade globality for scalability. Sampling-based planners are guaranteed, in the sampling limit, to find feasible~\cite{8584061} and globally optimal~\cite{doi:10.1177/0278364915614386} motion plans for general Lipschitz-continuous nonlinear systems. However, the number of samples required to find near-optimal solutions can grow exponentially with the dimension of the system, making these methods quickly intractable for high-dimensional dynamics. Gradient-based planners instead exploit smoothness in the dynamics, objective, and constraints to optimize a user objective under physical feasibility constraints. Although these methods generally provide only local optimality guarantees, they scale to high-dimensional systems and can perform whole-body dynamic motion planning. A prominent example is Contact-Implicit Trajectory Optimization (CITO), which jointly optimizes robot poses, controls, and contact forces through a high-dimensional trajectory optimization problem.

In contrast to the mature theory for sampling-based planning, the convergence behavior of gradient-based planners remains less well understood. Existing CITO formulations either cast the problem as a general nonlinear program (NLP) solved by off-the-shelf methods such as Sequential Quadratic Programming (SQP)~\cite{schulman2014motion,posa2014direct,winkler2018gait,manchester2019contact,landry2019bilevel,chatzinikolaidis2020contact}, or formulate it as an iterative Differential Dynamic Programming (DDP) problem~\cite{tassa2012synthesis,todorov2012mujoco,xie2017differential,kong2023hybrid,kim2025contact}. Both approaches face distinct modeling and regularity limitations. Off-the-shelf NLP solvers are not guaranteed to converge to a feasible solution unless suitable constraint qualifications (CQs) hold. Such conditions are difficult to verify for general nonlinear dynamic systems, and hard contact and collision constraints make them even more delicate. In particular, standard CQs such as the Linear Independence Constraint Qualification (LICQ) and the Mangasarian-Fromovitz Constraint Qualification (MFCQ) can fail when jointly optimizing contact points and forces induces mathematical programs with complementarity constraints (MPCCs). Weaker conditions such as MPCC-LICQ can sometimes be used, but even when verified, SQP typically provides only local convergence guarantees~\cite{doi:10.1137/S1052623403429081,byrd2010infeasibility}; it may still fail if the initial guess is not sufficiently close to a feasible solution. This contrasts sharply with sampling-based planners~\cite{8584061,doi:10.1177/0278364915614386}, whose feasibility guarantees hold from arbitrary initial samples in the asymptotic limit. We state this off-the-shelf-SQP guarantee and contrast it head-to-head with our result in~\prettyref{sec:problem}. DDP-based formulations face a different obstacle: they require the dynamics to be written as a differentiable explicit state-transition map. This structure makes the dynamics constraints regular by construction, since each transition constraint contains an identity block with respect to the next state, but such a map is difficult to define for contact-rich systems because discontinuities naturally arise from impacts, changing contact modes, and complex geometric constraints. Collision constraints between objects can also become non-differentiable under deep penetration or near singular geometric features such as medial axes. These limitations have been empirically observed and analyzed in prior work~\cite{schulman2014motion,zhang2023provably}, and they again contrast with sampling-based planners, which can handle complex environments through lightweight collision checking routines.

\TE{Main Result:}
Our goal is to strengthen the convergence analysis of CITO by studying a broader class of constraint-rich dynamic systems (CRDS): systems whose motion is governed by equality and inequality constraints, closely related to stiff dynamic systems~\cite{ashi2008numerical}, in which constraint forces can vary sharply with configuration. Articulated robots interacting with their environment form an important subclass. For this class we analyze a specific line-search SQP method with adaptive timestep refinement and prove that it converges globally in the numerical-optimization sense: given a strictly feasible startup configuration (Slater's condition), the method terminates finitely from any finite initial trajectory guess at an approximately feasible stationary trajectory of the original contact-implicit problem, in the force-parameterized Fritz--John sense made precise later. Two features distinguish this guarantee. First, no constraint qualification is assumed a priori; LICQ is instead derived from verifiable structure of the discretized dynamics, namely a sufficiently small timestep together with curvature and finite-horizon energy bounds. Second, the dynamics need not be written as an explicit differentiable transition function: the next state is defined implicitly as the solution of a system of equations, which permits differentiable contact models for robot links and environments given as arbitrary unions of convex polyhedra. We further show that the merit-penalty parameter and the optimized trajectory remain norm-bounded throughout the optimization history. The guarantee concerns convergence to stationarity rather than global optimality, and it applies to the specific method analyzed here rather than to arbitrary or off-the-shelf SQP solvers. These results are stated as~\prettyref{thm:stationarity-transfer} and~\ref{thm:nonsmooth-transfer} and proved in~\prettyref{sec:global-transfer}. We compare this guarantee at a high level with sampling-based planning and with existing gradient-based CITO formulations in~\prettyref{table:method-comparison}.

\begin{table*}[t]
\centering
\scriptsize
\setlength{\tabcolsep}{4pt}
\renewcommand{\arraystretch}{1.15}
\TableRowColors
\begin{tabular}{>{\raggedright\arraybackslash}m{0.13\textwidth}>{\raggedright\arraybackslash}m{0.13\textwidth}>{\raggedright\arraybackslash}m{0.15\textwidth}>{\raggedright\arraybackslash}m{0.28\textwidth}>{\raggedright\arraybackslash}m{0.19\textwidth}}
\toprule
\TableHeaderRow\TE{Method} & \TE{Type / solution scope} & \TE{Dynamic System} & \TE{Requirement on Contact Model} & \TE{CQ / regularity mechanism in cited guarantee} \\
\midrule
NLP~\cite{posa2014direct,schulman2014motion} & \TWORCell{Gradient-based}{(Local)} & Differentiable Implicit & \TWORCell{Differentiable distance function}{(fails at singularities)} & CQ such as LICQ/MFCQ\\
\midrule
DDP~\cite{tassa2012synthesis} & \TWORCell{Gradient-based}{(Local)} & Differentiable Explicit & \TWORCell{Differentiable transition/contact model}{(fails at singularities)} & No CQ needed\\
\midrule
Ours & \TWORCell{Gradient-based}{(Local)} & Implicit CRDS & Arbitrary unions of convex polyhedra & No CQ needed\\
\midrule
SST$^\star$~\cite{doi:10.1177/0278364915614386} & \TWORCell{Sampling-based}{(Global)} & Lipschitz Explicit & Collision-checking function & \TWORCell{Chow's condition}{Bounded second derivatives} \\
\bottomrule
\hiderowcolors
\end{tabular}
\caption{High-level comparison of representative contact-rich planning formulations. ``Local'' and ``global'' in the second column refer to solution-search scope, not to whether every listed numerical iteration has a global-convergence theorem in the same sense. The last column reports the CQ or regularity mechanism emphasized by the cited guarantee and is not an exhaustive list of assumptions. NLP-based CITO can use differentiable implicit dynamics, but standard KKT guarantees require an externally assumed CQ such as LICQ or MFCQ~\cite{gill2005snopt,biegler2009large,doi:10.1137/S1052623403429081,byrd2010infeasibility}. For DDP, the next-state identity block regularizes the explicit transition constraints, while the differentiability and local-convergence assumptions remain. For the proposed method, LICQ is derived only for the smooth modified trajectory NLP, from bounded curvature and sufficiently small timesteps enforced by adaptive refinement; no CQ is established for the original independent-force MPCC~\prettyref{pb:CITO}. The strict-interior predecessor, finite optimized mesh, fixed problem and regularization data, compact control set, and stated algorithmic safeguards remain in force. The method supports robot links and environments given as arbitrary unions of convex polyhedra through the separating-plane contact model~\cite{11266917}. In contrast, SST$^\star$ has a sampling-based global planning guarantee under Lipschitz explicit dynamics, Chow's condition, and bounded second derivatives~\cite{doi:10.1177/0278364915614386}.}
\label{table:method-comparison}
\end{table*}
Our main contributions are as follows:
\begin{enumerate}[label=(\alph*),leftmargin=*]
\item the CRDS formulation of CITO through implicit physics constraints;
\item a global-convergence guarantee, in the numerical-optimization sense, for a specific line-search SQP method with adaptive timestep refinement, certifying an approximate penalty-induced force-parameterized Fritz--John condition for the original physical formulation while deriving LICQ from the curvature and finite-horizon energy structure of the discretized dynamics rather than assuming a constraint qualification a priori;
\item a bounded-optimization-history result for the merit-penalty parameter and optimized trajectory, with common bounds over the fixed-data tolerance families of the two transfer theorems.
\end{enumerate}

The key idea behind our analysis is to expose useful structure in the physical constraints of CRDS, drawing on several recent lines of work. We develop each step for the general CRDS class and, in dedicated specialization subsections, show that articulated body dynamics is a concrete instance. The first step is \TE{Coordinate Maximization}: inspired by recent work on rigid-body simulation with maximal coordinates~\cite{brudigam2021linear,lan2022affine}, we lift the configuration space into an ambient space in which the inertial term becomes a quadratic energy, so that the kinetic energy is strongly convex with positive curvature inversely proportional to the timestep size. The second step is \TE{Penalization}: classical penalty methods~\cite{1641984,harmon2009asynchronous} convert the equality and inequality constraints into soft penalty functions, which have a physical interpretation as conservative potential energies; a variant of a recently proposed separating-plane contact model~\cite{11266917} makes the contact constraints twice differentiable for object shapes represented as unions of convex polyhedra. The third step is a \TE{Hamiltonian Upper Bound}: the total energy of the CRDS remains bounded over a finite horizon, which in turn bounds the curvature of the conservative potential energy. Choosing a sufficiently small timestep then lets the positive curvature of the kinetic energy dominate, yielding non-vanishing curvature and implying LICQ. Combined with the line-search SQP method and adaptive-refinement safeguards analyzed in this paper, this yields the global convergence guarantee described above. We also use recently developed SQP constraint-violation bounds~\cite{solodov2009global} to establish boundedness of the iterates and merit-penalty sequence throughout the specified algorithm. Finally, preliminary numerical experiments illustrate the algorithm on a range of CITO instances spanning robot locomotion and manipulation.

\section{Related Work}
\label{sec:related-work}

\subsection{Dynamic Systems and Contact Models}
\label{sub:related-dynamics}

Trajectory optimization for robotic systems depends strongly on how the dynamics are represented. Explicit formulations write the next state as a transition map, which is convenient for shooting methods, DDP, iterative Linear-Quadratic Regulator (iLQR), and model predictive control (MPC) because each rollout directly propagates states. This structure also makes the state-transition constraints regular in many direct transcriptions, since the next state appears through an identity block. However, for systems with contact, impacts, friction, closed kinematic chains, and collision constraints, constructing a globally differentiable explicit map is difficult. Implicit formulations instead represent dynamics, contact, and kinematic restrictions as equations, inequalities, complementarity conditions, or variational optimality conditions. They introduce more constraints, but they expose the physical feasibility conditions that CITO optimizes over jointly, and are therefore better suited to constraint-rich dynamic systems.

This implicit, complementarity-based view of contact has a long lineage in nonsmooth mechanics. Unilateral contact and dry friction are classically modeled through set-valued force laws and measure-differential inclusions~\cite{moreau1988unilateral,brogliato1999nonsmooth}, and implicit time-stepping schemes discretize the resulting dynamics as linear complementarity problems that stay well posed through impacts and stick-slip transitions~\cite{stewart1996implicit,stewart2000rigid}. These complementarity conditions are the direct origin of the complementarity constraints that appear in CITO. Because such constraints are difficult to solve directly, a parallel line of work reformulates MPCCs as ordinary nonlinear programs through relaxation or penalization~\cite{anitescu2000solving}. The penalty relaxation of the equality and inequality constraints we use for CRDS follows this MPCC-as-NLP tradition, but rather than treating it only as a numerical device, we combine the smooth surrogate with the constant, positive-definite mass matrix of maximal coordinates to establish the regularity our convergence analysis requires.

Maximal-coordinate models provide one route to such implicit formulations. Instead of eliminating joint constraints into a minimal coordinate chart, maximal coordinates represent bodies in an ambient space and enforce joints, loops, and other integrity constraints algebraically. \citewithauthor{brudigam2021linear} showed that this representation can be useful even for linear-quadratic optimal control, where the lifted coordinates increase dimension but simplify the treatment of constraints. Related simulation work uses the same ambient-coordinate philosophy for difficult physical systems: Affine Body Dynamics represents stiff materials through affine body variables and constraint-like potentials~\cite{lan2022affine}, while the Shape-Differentiable Robot Simulator uses convex-polyhedral shape variables and separating-plane auxiliary variables to obtain differentiable robot simulation with respect to geometry~\cite{11266917}. Our coordinate maximization step follows this line of work: it keeps the dynamics implicit and lifted, but uses the resulting constant, positive-definite mass matrix to prove regularity rather than treating maximal coordinates only as a numerical modeling device.

Contact modeling is the second major source of difficulty. Classical penalty methods are simple and often differentiable away from geometric singularities, but they allow penetration, require stiffness tuning, and can become ill-conditioned or singular under deep interpenetration~\cite{1641984,harmon2009asynchronous}. Smooth contact models used in model-based control soften complementarity or contact response so that local optimization can be applied more effectively~\cite{todorov2011convex,tassa2011stochastic,todorov2012mujoco}. These models are important precursors, but smoothing a local contact response does not by itself remove nonsmoothness caused by changing closest features, changing active contact pairs, or singular geometric configurations. More recent barrier and interior-point-inspired contact formulations, such as Incremental Potential Contact~\cite{li2020incremental} and Affine Body Dynamics~\cite{lan2022affine}, maintain nonintersection through barrier potentials and line-search safeguards. Our contact construction builds on this modern trend and on the separating-plane contact model for polyhedral geometry~\cite{11266917}: rather than assuming a globally smooth signed-distance function, it introduces auxiliary variables that keep contact constraints differentiable for links and environments given as arbitrary unions of convex polyhedra. An alternative smoothing route replaces each polyhedron by a strictly convex hull so that the proximity distance acquires continuous gradients~\cite{6710113}; our separating-plane construction instead retains the exact polyhedral geometry while keeping the contact constraints differentiable.

\subsection{Trajectory Optimization and CITO}
\label{sub:related-cito}

Sampling-based motion planners such as Stable Sparse RRT$^\star$ (SST$^\star$) provide asymptotic global optimality guarantees for general Lipschitz-continuous nonlinear systems~\cite{doi:10.1177/0278364915614386} but scale poorly with system dimension; the trajectory-optimization line of work discussed in this paper instead trades global optimality for scalability. A complementary sampling-based guarantee is probabilistic completeness: such planners find a feasible motion from arbitrary initial samples in the sampling limit~\cite{8584061}. Our result is different in kind. It establishes convergence of a specified local, gradient-based iteration, from an admissible finite discretized trajectory initialization, to an approximately feasible stationary point; it is neither probabilistic completeness nor asymptotic global optimality. Before contact-implicit formulations became common, trajectory optimization was widely used for motion-planning and control problems in which contact is absent or treated purely as a kinematic collision-avoidance constraint. Kinematic motion-planning methods such as CHOMP~\cite{ratliff2009chomp} and STOMP~\cite{kalakrishnan2011stomp} optimize smooth manipulator trajectories under workspace obstacle potentials, and Gaussian-process planners cast trajectory optimization as probabilistic inference on factor graphs~\cite{dong2016motion}. In aerial robotics, differentially flat formulations enable polynomial trajectory generation for quadrotors with no environmental contact~\cite{mellinger2011minimum}. Robust collision avoidance can also be cast as a semi-infinite program over swept volumes with subdivision-based feasibility guarantees~\cite{zhang2023provably}. These methods illustrate the scalability and convergence behavior of optimization-based planning when the dynamics are smooth and contact-free, but they sidestep the central modeling difficulty of contact-rich trajectory optimization.

CITO addresses this gap by allowing contact forces or modes to be decision variables. \citewithauthor{posa2014direct} introduced a direct transcription for rigid-body systems through contact using complementarity constraints, casting the problem as a nonlinear program. Sequential convex optimization with convex collision checking was also used to absorb collision constraints into local NLP subproblems~\cite{schulman2014motion}, while mixed-integer formulations can reason globally about kinematic contact placement or inverse kinematics without optimizing full contact dynamics~\cite{deits2014footstep,dai2019global}. Subsequent work extended the modeling and application scope of NLP-based CITO. \citewithauthor{winkler2018gait} optimized legged locomotion trajectories using phase-based end-effector parameterizations, improving practical gait generation while retaining assumptions about contact structure. Variational contact-implicit trajectory optimization was developed in~\cite{manchester2019contact,manchester2020variational}, using variational integrators to better preserve mechanical structure. \citewithauthor{landry2019bilevel} proposed a bilevel semidirect method in which lower-level contact dynamics are embedded inside the planning problem. \citewithauthor{chatzinikolaidis2020contact} used an analytically solvable contact model for locomotion on variable ground. Several of these formulations make the complementarity structure explicit by casting contact-rich planning as an MPCC, including non-prehensile planar manipulation with complementarity constraints~\cite{moura2022non}, consensus complementarity control for multi-contact model predictive control~\cite{aydinoglu2023consensus}, and simultaneous trajectory optimization and contact selection over high-fidelity geometry~\cite{zhang2025simultaneous}. These share the MPCC form we also adopt, but solve it with general nonlinear-programming or complementarity solvers whose convergence guarantees still rely on constraint qualifications. These methods show that implicit contact formulations can solve rich locomotion and manipulation problems, but the resulting NLPs inherit the regularity and smoothness issues of complementarity and geometric contact constraints.

A separate thread of research also targets global behavior for contact-rich manipulation, but in senses different from the one studied here. \citewithauthor{aceituno2020global} formulate contact-trajectory optimization as a mixed-integer convex program: taking the object motion as an input and relaxing the bilinear force-position terms with McCormick envelopes under a quasi-dynamic model, the problem is solved to global optimality by branch-and-bound, at the cost of worst-case exponential complexity in the integer variables. \citewithauthor{pang2023global} instead achieve global reach through sampling, combining an RRT-style global search with local smoothing of a convex, quasi-dynamic, differentiable contact model so that smoothed gradients can guide the search across contact modes. Both rely on a quasi-dynamic relaxation of contact and obtain globality either from combinatorial search over a convex relaxation or from randomized sampling over a smoothed model; neither establishes convergence of a continuous local solver on the full second-order dynamics. Our contribution is complementary and of a different type: we retain the implicit CRDS as the target physical model, construct from it the smooth modified trajectory NLP, and analyze a specific local, gradient-based SQP iteration on that modified problem. We make no global-optimality claim. The derived LICQ yields global convergence of the iteration to an approximate KKT point of that modified problem, and a separate transfer argument interprets that certificate as an approximately feasible, force-parameterized stationary point of the original CITO~\prettyref{pb:CITO}.

A parallel line of work keeps the rollout-based structure of DDP or iLQR. Online DDP and iLQR methods optimize trajectories through explicit dynamics and local quadratic models, and have been demonstrated on humanoid systems whose contact modes vary along the rollout when combined with smoothed contact responses~\cite{tassa2012synthesis}. Constrained DDP methods incorporate nonlinear constraints into local dynamic programming updates~\cite{xie2017differential}; hybrid iLQR and contact-implicit MPC methods then apply related ideas to contact-rich legged systems~\cite{kong2023hybrid,kim2025contact}. These methods are attractive for fast feedback control, but they typically require an explicit differentiable transition model, a fixed smooth approximation of contact, or a contact solver whose derivatives are well defined along the rollout. This is the main distinction from the setting studied here. A closely related line goes further and replaces hand-derived dynamics with locally differentiable physics simulators, obtaining search directions by locally differentiating the contact solve itself; however, such differentiability can be lost precisely at impacts and stiff contact events, where the smoothed simulation departs from the true nonsmooth dynamics~\cite{suh2022differentiable}, which again motivates a convergence guarantee that does not rely on a globally smooth differentiable rollout. Our method supplies precisely such a guarantee: it keeps the implicit CRDS as the target model and admits contact constraints for arbitrary unions of convex polyhedra through auxiliary separating planes, without requiring a differentiable rollout.

\subsection{NLP Algorithms and Convergence Conditions}
\label{sub:related-nlp}

Most direct CITO formulations are solved with general nonlinear-programming methods, especially SQP and primal-dual interior-point algorithms. SNOPT is a representative large-scale SQP solver~\cite{gill2005snopt}, and IPOPT is a representative large-scale interior-point framework~\cite{biegler2009large}. Classical SQP theory and numerical-optimization texts analyze these methods under smoothness, second-order approximation quality, and constraint qualifications such as LICQ or MFCQ~\cite{boggs1995sequential,byrd2010infeasibility,nocedal2006numerical}. These assumptions are natural for smooth NLPs, but they are difficult to verify for contact-rich trajectory optimization. Contact complementarity can violate standard CQs, and geometric collision constraints may fail to be smooth at deep penetrations, medial-axis configurations, or changes in closest features. Indeed, the sharpest guarantees available for off-the-shelf SQP are dichotomous: an exact-penalty SQP method converges either to a Karush-Kuhn-Tucker (KKT) point when a constraint qualification holds or, otherwise, only to a stationary point of the infeasibility measure~\cite{byrd2010infeasibility}, so a CQ remains an assumption to be checked rather than a property established for the problem---exactly the gap our structure-derived LICQ closes.

The optimization literature has developed weaker or specialized regularity notions for complementarity-constrained problems. MPCC formulations are especially relevant because many CITO models encode unilateral contact and friction through complementarity. \citewithauthor{doi:10.1137/S1052623403429081} analyzed interior-point methods for MPCCs using MPCC-specific regularity conditions, showing that complementarity structure can be handled more carefully than by treating the problem as an ordinary NLP. Such results clarify the limitation of simply invoking off-the-shelf SQP or interior-point convergence theory: even when MPCC-LICQ or related conditions replace ordinary LICQ, those conditions remain assumptions to be checked for the particular contact trajectory and contact geometry. General SQP theory can also weaken other hypotheses, for example by avoiding boundedness assumptions on the iterates~\cite{solodov2009global}, but it still does not prove that a CITO instance satisfies the requirements. Recent worst-case analyses of log-barrier methods for nonconvex-constrained programs sharpen these guarantees but still require a regularity condition to certify a KKT point~\cite{hinder2024worst}---exactly the structure-derived LICQ our analysis supplies.

Beyond constraint qualifications, complementarity-constrained programs also require specialized notions of stationarity, because the classical KKT conditions may be unattainable when MPCC-LICQ fails. The optimality theory for MPCCs distinguishes a hierarchy of stationarity concepts---from the weak Clarke (C-) stationarity built on the Clarke generalized gradient and normal cone~\cite{clarke1983optimization}, through Mordukhovich (M-) stationarity, up to the strong (S-) stationarity that recovers the KKT conditions of the tightened program~\cite{scheel2000mathematical,ye2005necessary,outrata1999optimality,flegel2003fritz}, with relaxation and regularization schemes known to deliver these concepts to differing degrees~\cite{hoheisel2013theoretical}. Certifying such points algorithmically typically calls for nonsmooth or semismooth generalizations of Newton and SQP methods that operate on generalized derivatives rather than classical Jacobians~\cite{qi1993nonsmooth,facchinei2003finite,curtis2012sequential}. Our nonsmooth certificate is adjacent to, but not a member of, this standard full-MPCC hierarchy. The Clarke normal in our certificate is taken to the lifted damping subgradient graph, while the equality and inequality/contact forces remain restricted to penalty-induced configuration-dependent maps. Consequently, the result is a homogeneous force-parameterized Clarke--Fritz--John certificate, not standard MPCC C-stationarity. Standard MPCC C-, M-, and S-stationarity would additionally require independent force variables and the corresponding normal cone to the complementarity set. Promoting the forces to decision variables and targeting the corresponding M- or S-stationarity notions is a natural next step.

There is also convergence analysis for online optimization and nonlinear MPC. Fast MPC and nonlinear MPC algorithms exploit repeated similar optimization problems, warm starts, and first- or second-order updates~\cite{wang2009fast,stella2017simple}. More recent analyses establish global convergence or superconvergence properties for online MPC under smoothness, sensitivity, and regularity assumptions appropriate to those control problems~\cite{na2021global,na2022superconvergence}. These results are complementary to CITO, but they do not resolve the central obstacle studied in this paper: for contact-implicit problems, the solver must cope with implicit dynamics, complementarity-like contact structure, and geometric nonsmoothness without assuming a constraint qualification a priori. Our analysis addresses this gap by proving LICQ from the curvature structure induced by maximal coordinates, penalty relaxations, and finite-horizon Hamiltonian bounds, and then applying a specific safeguarded line-search SQP framework to obtain bounded iterations that converge globally, in the numerical-optimization sense, to the stated approximate stationarity certificate.

\section{Problem Definition}
\label{sec:problem}
In this section, we present the specific class of CITO problems considered in this paper. We begin by defining the CRDS. Many problems of practical interest fall in this class: an articulated robot is a collection of rigid bodies tied together by joint equality constraints, and its interaction with the environment adds unilateral non-penetration and frictional contact inequalities. We represent the system's kinematic configuration by an ambient maximal-coordinate configuration $\theta\in\mathbb R^n$ drawn from the ambient maximal-coordinate configuration space, associated with a mass matrix $M$. We work in maximal rather than the generalized coordinates more common in robotics because the analysis needs the inertial term to be a fixed quadratic energy: in maximal coordinates the mass matrix is constant, whereas a generalized-coordinate mass matrix is configuration-dependent. Joint structure is then imposed through explicit equality constraints rather than being built into the coordinates. The system is subject to several force terms: a conservative potential force $f(\theta)=-\FPP{U(\theta)}{\theta}$, where $U(\theta)$ is the corresponding conservative potential energy; a gravitational term $g$; a control force term $c(\theta,u)$ where $u$ is the control input. We assume that the control force is integrable in the configuration variable, so it can be written as $c(\theta,u)=-\FPP{W(\theta,u)}{\theta}$ for a scalar aggregate control work potential $W(\theta,u)$. The system is also subject to hard equality constraints $h^e(\theta)=0$ and inequality constraints $h^i(\theta)\geq0$, where $h^e$ and $h^i$ are vector-valued maps. Accordingly, we introduce the constraint forces $\lambda^e$ and $\lambda^i$. We then introduce a damping force $d$ and assume that it lies in the negative subgradient set of the (possibly non-smooth) convex function $\Damp(\theta,\dot{\theta},\lambda^i)$ in $\dot{\theta}$, the dependence on $\lambda^i$ reflecting that frictional dissipation is modulated by the normal contact forces. Our formulation of damping models frictional contact under the maximal dissipation principle (detailed in~\prettyref{sec:CRDS}). These force and constraint components are illustrated on a two-link articulated example in~\prettyref{fig:crds-illustration}. Finally, we adopt a variable-timestep discretization over the finite planning horizon $H$. Let $\delta>0$ be a fixed base time scale and assign physical time
\begin{equation}
\label{eq:elapsed-time}
t_i\triangleq\delta i
\end{equation}
to every relevant real node label. Node $0$ is fixed at physical time zero and is not optimized. Its fixed pre-zero predecessor satisfies $\gamma(0)<0$. The optimized nodes form a finite nonempty totally ordered set
\begin{equation}
\mathcal I\subseteq\{i\in\mathbb R\mid 0<\delta i\leq H\}.
\end{equation}
For each $i\in\mathcal I$, $\gamma(i)$ is its immediate predecessor in $\{0\}\cup\mathcal I$; the fixed value $\gamma(0)$ completes the first two-step stencil. The local physical timestep is
\begin{equation}
\delta_i\triangleq\delta(i-\gamma(i)),\qquad i\in\{0\}\cup\mathcal I,
\end{equation}
so in particular the incoming fixed step $\delta_0=\delta(0-\gamma(0))$ is positive. Two conventions here are worth motivating. First, node $0$ and its pre-zero predecessor are fixed rather than optimized because the discrete dynamics form a two-step stencil: the residual at a node involves that node and its two predecessors, so the first optimized node already needs two configurations in hand for its stencil to be complete; $\theta_0$ additionally serves as the fixed strictly feasible configuration from which the energy bounds are later derived. Second, node labels are real rather than integer so that the adaptive refinement used later can insert a node between two existing ones, halving the local timestep, without renumbering the trajectory. For example, $\theta_i$ and $u_i$ are the kinematic state and control signal at optimized node $i$. Our discretization scheme is summarized as $<\THREE{\theta}{u}{\delta}_i>$, where angle brackets denote the column-stack of their arguments, that is, the concatenation of the listed blocks into a single vector; we use standard first- and second-order discretization for velocity and acceleration
\begin{ResizedEq}
\label{eq:discretization}
\dot{\theta}_i=&(\theta_i-\theta_{\gamma(i)})/\delta_i\\
\ddot{\theta}_i=&
\left[\theta_i-(1+\delta_i/\delta_{\gamma(i)})\theta_{\gamma(i)}+(\delta_i/\delta_{\gamma(i)})\theta_{\gamma^2(i)}\right]/\delta_i^2,
\end{ResizedEq}
Although $\delta$ is fixed, the local timestep $\delta_i$ varies with the ordered label gap $i-\gamma(i)$; this variable-step structure is central to our convergence analysis. At the first optimized node $\min\mathcal I$, the two fixed predecessor configurations in~\prettyref{eq:discretization} are $\theta_0$ and $\theta_{\gamma(0)}$. Put together, we summarize the definition of CRDS as follows.
\begin{definition}
\label{def:CRDS}
Throughout, $L_\bullet$ denotes a finite upper bound on the quantity indicated by its subscript. A CRDS is a discrete dynamic system whose complete maximal-coordinate configuration satisfies $\theta_i\in\mathbb R^n$, defined by
\begin{ResizedEq}
\label{eq:CRDS}
\begin{cases}
M\ddot{\theta}_i=
\begin{cases}
f(\theta_i)+g+c(\theta_i,u_i)+d+\\
\FPP{h^e}{\theta_i}^T\lambda_i^e+\FPP{h^i}{\theta_i}^T\lambda_i^i\\
\end{cases}\\
h^e(\theta_i)=0\\
0\leq\lambda_i^i\perp h^i(\theta_i)\geq0\\
d\in-\partial_{\dot{\theta}_i} \Damp(\theta_{\gamma(i)},\dot{\theta}_i,\lambda_{\gamma(i)}^i),
\end{cases}
\end{ResizedEq}
Define the compact convex control space $\mathcal{D}$ and hard inequality-feasible subset
\begin{equation}
\mathcal{F}^i\triangleq\{\theta\in\mathbb R^n\mid h^i(\theta)\geq0\}.
\end{equation}
The CRDS is assumed to further satisfy four conditions:
\begin{enumerate}[label=(\alph*),ref=(\alph*),leftmargin=*]
\item\label{it:crds-c2} there is an open set $\mathcal{D}_u\supset\mathcal{D}$ such that $W:\mathbb R^n\times\mathcal{D}_u\to\mathbb R$ is jointly $C^2$ in $\TWO{\theta}{u}$, while $h^e,h^i$ are $C^2$ real functions on $\mathbb R^n$;
\item\label{it:crds-potential} $U\geq0$ is proper, continuous, and extended-real-valued. $U$ is $C^2$ on $\operatorname{dom}U$;
\item\label{it:crds-damping} for every fixed predecessor configuration $\theta_\gamma\in\mathcal{F}^i$ and fixed nonnegative damping load $\lambda_\gamma^i\geq0$, $\Damp(\theta_\gamma,\cdot,\lambda_\gamma^i)$ is convex (possibly non-smooth) and $0\in\partial_{\dot{\theta}}\Damp(\theta_\gamma,0,\lambda_\gamma^i)$;
\item\label{it:crds-control} $\|c(\theta_i,u_i)\|\leq\LControlJacobianBound<\infty$ for every $(\theta_i,u_i)\in\mathbb R^n\times\mathcal{D}$, for some $\LControlJacobianBound$. 
\end{enumerate}
\end{definition}
\begin{figure}[t]
\centering
\begin{tikzpicture}[font=\small, line cap=round, line join=round, >=stealth]

\definecolor{topBoxFill}{RGB}{255,213,170}     
\definecolor{bottomBoxFill}{RGB}{255,239,178}  

\tikzset{
  capsuleOutline/.style={black, line width=1.15pt},
  forceArrow/.style={black, line width=1.0pt, ->},
}

\def\linkRadius{0.32}

\def\drawCapsule#1#2#3#4#5{
  \begin{scope}[shift={(#1,#2)}, rotate=#3]
    \fill[#4] (-#5,-\linkRadius) -- (#5,-\linkRadius)
      arc[start angle=-90, end angle=90, radius=\linkRadius]
      -- (-#5,\linkRadius)
      arc[start angle=90, end angle=270, radius=\linkRadius] -- cycle;
    \draw[capsuleOutline] (-#5,-\linkRadius) -- (#5,-\linkRadius)
      arc[start angle=-90, end angle=90, radius=\linkRadius]
      -- (-#5,\linkRadius)
      arc[start angle=90, end angle=270, radius=\linkRadius] -- cycle;
  \end{scope}
}

\useasboundingbox (-1.8,-0.9) rectangle (3.6,7.0);

\coordinate (hinge) at (0,2.7);
\coordinate (com) at (0.8,4.0);
\coordinate (topCap) at (1.6,5.3);
\coordinate (foot) at (2.0,0.32);
\coordinate (contact) at (2.0,0.0);

\draw[black, line width=1.1pt] (-0.15,6.4) -- (3.35,6.4);
\foreach \x in {0.1,0.5,...,3.3}{
  \draw[black, line width=0.5pt] (\x,6.4) -- (\x-0.3,6.7);
}

\draw[black, line width=1.0pt]
  (1.6,6.40) -- (1.6,6.30)
  -- (1.82,6.24) -- (1.38,6.15)
  -- (1.82,6.06) -- (1.38,5.97)
  -- (1.82,5.88) -- (1.38,5.79)
  -- (1.82,5.70) -- (1.38,5.61)
  -- (1.82,5.52) -- (1.38,5.44)
  -- (1.6,5.38) -- (1.6,5.30);
\node[black, anchor=west] at (1.95,5.85) {$f=-\FPP{U}{\theta_i}$};

\drawCapsule{0.8}{4.0}{58.4}{topBoxFill}{1.526}      
\drawCapsule{1.0}{1.51}{-50}{bottomBoxFill}{1.554}   

\fill[black] (com) circle (2.6pt);
\draw[forceArrow] (com) -- (0.8,3.5);
\node[black, anchor=south] at (1.0,4.0) {$g$};

\begin{scope}[shift={(hinge)}]
  \draw[forceArrow] (-25:0.85) arc[start angle=-25, end angle=35, radius=0.85];
\end{scope}
\node[black, anchor=west] at (0.98,2.48) {$c=-\FPP{W}{\theta_i}$};

\node[black, anchor=east] at (-0.7,2.7) {$h^e$};
\draw[forceArrow] (-0.58,2.74) to[out=30, in=150] (hinge);
\fill[black] (hinge) circle (2.2pt);

\draw[black, line width=1.1pt] (-0.6,0.0) -- (3.3,0.0);
\foreach \x in {-0.4,0.0,...,3.2}{
  \draw[black, line width=0.5pt] (\x,0.0) -- (\x-0.3,-0.3);
}
\fill[black] (contact) circle (2.0pt);

\draw[forceArrow] (contact) -- (2.0,0.7);
\node[black, anchor=west] at (2.2,0.7) {$h^i$};
\draw[forceArrow] (1.68,0.22) -- (0.85,0.22);
\node[black, anchor=east] at (0.8,0.22) {$d\in-\partial_{\dot{\theta}_i}\Damp$};

\end{tikzpicture}
\caption{\label{fig:crds-illustration}Illustration of the CRDS on a two-link articulated example. The peach top link and butter-yellow bottom link form the articulated body with configuration $\theta_i$. The top link is suspended from the ceiling by a spring, which supplies the conservative potential force $f=-\FPP{U}{\theta_i}$; $g$ is gravity acting at the center of mass; and the control force $c=-\FPP{W}{\theta_i}$ acts at the joint. The hinge connecting the two links enforces the equality constraint $h^e(\theta_i)=0$, while non-penetration of the foot against the ground is the inequality constraint $h^i(\theta_i)\geq0$. Finally, $d\in-\partial_{\dot{\theta}_i}\Damp$ is the frictional dissipation acting tangentially at the contact. These are exactly the terms balanced in~\prettyref{eq:CRDS}.}
\end{figure}
\begin{definition}
\label{def:smooth-CRDS}
A CRDS is smooth if, for every fixed nonnegative damping load $\lambda_\gamma^i\geq0$, $\Damp(\theta_\gamma,\cdot,\lambda_\gamma^i)$ is differentiable in velocity and $\FPPR{\Damp}{\dot{\theta}}$ is $C^1$ with respect to all arguments.
\end{definition}

Recent work~\cite{lan2022affine} represents system dynamics in maximal coordinates and uses $h^e$ to encode system integrity constraints. As a result, the mass matrix can be a constant matrix $M$, which simplifies our analysis and is assumed throughout our CRDS formulation. Our analysis further requires that the matrix $M$ is positive definite, and we define such CRDS to be full-rank.

\begin{definition}
\label{def:full-rank}
A CRDS is full-rank if the mass matrix $M$ is positive definite.
\end{definition}
We will invoke~\prettyref{def:full-rank} in the global-convergence analysis to expose the positive curvature of the kinetic energy. To define the CITO problem using CRDS as physics constraints, we further assume that a user objective $O(\{\theta_i\},\{u_i\})$ is to be minimized, subject to the following regularity requirement.
\begin{assume}
\label{ass:objective}
The user objective $O(\{\theta_i\},\{u_i\})$ is continuously differentiable, and its gradient $\nabla O$ is Lipschitz continuous (with modulus $\LObjGradLip$).
\end{assume}
\prettyref{ass:objective} is the only regularity we impose on the objective: the $C^2$ smoothness demanded by our analysis is carried entirely by the dynamics and constraint data (through~\prettyref{def:CRDS}~\ref{it:crds-c2}), not by $O$. All later invocations of the objective's smoothness---the quadratic-program (QP)/KKT gradients of~\prettyref{sec:global} and the descent estimates of its convergence proof---refer back to this assumption. We optimize over the finite ordered set $\mathcal I$ defined above and restrict every optimized control $u_i$ to the compact convex control space $\mathcal D$. The assumption on compactness almost always holds for general robotic systems, since control forces and torques must be bounded. The fixed startup data are $\theta_{\gamma(0)},\theta_0,\delta_0$, and $u_0\in\mathcal D$; these data are distinct from the optimized horizon variables, and the fixed history interval $(\gamma(0),0)$ is never subdivided. The convergence theorem permits any finite nonempty ordered initial optimized set $\mathcal I$, any finite configuration $\theta_i\in\mathbb R^n$ for every $i\in\mathcal I$, and any control $u_i\in\mathcal D$ for every $i\in\mathcal I$. The first optimized node is $\theta_{\min\mathcal I}$, and the interval $(0,\min\mathcal I)$ is the first refinable interval. Objectives, inertia, gravity, and controls depend only on the physical components $\theta_i$. This yields the following CITO formulation as a constrained optimization problem.
\begin{equation}
\begin{aligned}
\label{pb:CITO}
\argmin{\substack{\{\theta_i\},\{u_i\in\mathcal{D}\}\\\{\lambda_i^e\},\{\lambda_i^i\},i\in\mathcal{I}}}&O(\{\theta_i\},\{u_i\})\\
\ST\quad&\text{\prettyref{eq:CRDS}}\quad\forall i\in\mathcal{I},
\end{aligned}
\end{equation}
with any model-specific auxiliary variables and constraints represented abstractly through the instantiated complete CRDS configuration and maps. The concrete extended articulated-body instantiation is introduced in~\prettyref{sec:CRDS} and written explicitly in~\prettyref{pb:CRDS-ABD-extended}.
\begin{assume}[strict-interior predecessor]
\label{ass:Slater}
The fixed time-zero configuration belongs to the finite domain of the conservative potential, $\theta_0\in\operatorname{dom}U$, satisfies the equality constraints $h^e(\theta_0)=0$, and lies strictly inside every inequality domain, $h^i(\theta_0)>0$ (Slater's condition~\cite{bertsekas1997nonlinear}).
\end{assume}
\prettyref{ass:Slater} is a substantive start-up hypothesis. It is stronger than ordinary feasibility and excludes the fixed time-zero node from active unilateral contact. The strict inequality is required by the present extended-real barrier and penalty-induced predecessor-load construction. For the first optimized interval $i_1\triangleq\min\mathcal I$, $\gamma(i_1)=0$, so the predecessor configuration is the already-fixed $\theta_0$. Since $h^i(\theta_0)>0$, the exact predecessor force is well defined and complementarity gives $\lambda_0^i=0$. The normal force used to modulate damping is instead the prescribed penalty force evaluated at the same fixed $\theta_0$; it is finite on the strict-interior penalty domain and may be nonzero inside the penalty layer. Thus no additional startup state, optimization variable, or assumption is required. Anchoring the barrier at fixed node $0$ supplies a finite fixed penalty/energy base and admissible semi-implicit damping loads for the first optimized step; the feasibility and clamp-inactivity of $\theta_{\min\mathcal I}$ are then derived from the first energy step of~\prettyref{sec:Hamiltonian}, not assumed.
\begin{remark}[Initialization scope and strict-interior predecessor]
\label{rem:initialization-scope}
Because the inequality barrier is finite only where every inequality constraint has strict slack, the predecessor penalty, the penalty-induced previous-step force, and the semi-implicit dissipation load are finite and defined only for $\theta_0\in\operatorname{dom}V$. Throughout the convergence statements, the fixed startup data $\theta_{\gamma(0)},\theta_0,\delta_0,u_0$ satisfying~\prettyref{ass:Slater}, together with the regularization, clamp, and algorithm parameters, are fixed. In contrast, a finite discretized optimized trajectory guess consists of a finite nonempty ordered set $\mathcal I$, a finite $\theta_i\in\mathbb R^n$ for every $i\in\mathcal I$, and a control $u_i\in\mathcal D$ for every $i\in\mathcal I$. No modified-dynamics or original-contact feasibility is required of any optimized $\theta_i$. Arbitrary initialization does not include the fixed startup data, their contact mode, smoothing parameters, clamps, or algorithmic constants. Extending the theorem to a time-zero node with active contact would require a different start-up barrier/dissipation construction or an additional first-step argument.
\end{remark}

\subsection{Statement of Main Result}
To make this cost concrete, we recall the most representative rigorous guarantee available for off-the-shelf SQP, using a generic decision variable $x$, objective $O$, and constraint map $F$, so that the comparison with our result is direct.
\begin{framed}
\begin{theorem}[(prior work~\cite{byrd2010infeasibility}, informal)]
\label{thm:bcn-infeasibility}
Consider a general smooth NLP $\min_x O(x)$ subject to $F(x)=0$, with $F$ the constraint map (inequality constraints handled analogously) and $\|F(\cdot)\|$ the associated infeasibility measure. Apply an exact-penalty SQP method that adjusts its penalty parameter automatically. Then every limit point of the generated iterates is either
\begin{enumerate}[label=(\alph*),ref=(\alph*),leftmargin=*]
\item\label{it:bcn-kkt} a KKT point of the program, provided a constraint qualification (such as MFCQ) holds; or
\item\label{it:bcn-infeas} an infeasible stationary point, i.e., a stationary point of $\|F(\cdot)\|$ with $F\neq0$.
\end{enumerate}
\end{theorem}
\end{framed}

\begin{algorithm*}[t]
\caption{\small{CITOSolve (informal)}}
\label{alg:CITO}
\begin{algorithmic}[1]
\Require{an initial trajectory guess, with controls in the admissible set}
\Repeat
    \Repeat
        \If{\hldiff{the physics residual lies outside the working band}}
            \State \hldiff{Reduce the residual alone until it re-enters the band}
        \EndIf
        \If{\hldiff{the linearized physics constraints are ill conditioned at some node}}
            \State \hldiff{Halve that timestep}
        \Else
            \State Solve a quadratic subproblem to obtain a search direction
            \State Accept a fraction of that direction by line search on a merit function
        \EndIf
    \Until{the physics residual and the accepted step are below the stopping tolerance}
    \If{\hldiff{the energy bound fails at some node}}
        \State \hldiff{Halve that timestep}
    \EndIf
\Until{\hldiff{the energy bound holds at every node}}
\State Return the optimized trajectory and controls
\end{algorithmic}
\end{algorithm*}

In contrast to such off-the-shelf guarantees, whose convergence hinges on externally assumed constraint qualifications, the method developed in this paper assumes no CQ a priori. The main results below state that the specific line-search SQP method with adaptive timestep refinement terminates finitely and, through the stationarity transfer developed later, returns a force-parameterized approximate Fritz--John certificate for~\prettyref{pb:CITO} while maintaining the stated boundedness properties. More precisely, the direct algorithmic conclusion is the $\EpsStop$-KKT condition of the modified~\prettyref{pb:CITO-M} established in~\prettyref{thm:sqp-convergence}--\ref{thm:CITO-convergence}. The informal theorems below state the result after the separate bridge and stationarity-transfer steps have reinterpreted that certificate on the original~\prettyref{pb:CITO}. The transferred certificate is evaluated after restricting the constraint forces to the penalty-induced maps of~\prettyref{def:eps-stationary}; it is therefore not the standard stationarity condition of~\prettyref{pb:CITO} with independently variable contact forces. Throughout~\prettyref{thm:stationarity-transfer}--\ref{thm:nonsmooth-transfer}, ``force-parameterized'' has this specific restricted meaning: the constraint forces of~\prettyref{pb:CITO} are replaced by the penalty-induced maps of~\prettyref{def:eps-stationary}, so the certificate controls the physical residual and the trajectory/control stationarity for that selected force graph but provides no first-order conditions in independent force directions.
\begin{framed}
\begin{theorem}[(informal)]
\label{thm:stationarity-transfer}
Consider~\prettyref{pb:CITO} instantiated by a full-rank (\prettyref{def:full-rank}), \hldiff{smooth} (\prettyref{def:smooth-CRDS}) CRDS with controls in the compact convex control space $\mathcal{D}$, and suppose~\prettyref{ass:objective} and~\ref{ass:Slater} hold.
\begin{enumerate}[label=(\alph*),ref=(\alph*),leftmargin=*]
\item For fixed problem and pre-horizon data satisfying~\prettyref{ass:objective} and~\ref{ass:Slater}, any finite nonempty ordered initial optimized set, any finite optimized configuration guess, and any control guess in $\mathcal D$, a line-search SQP with adaptive timestep refinement terminates after finitely many iterations at a trajectory satisfying the penalty-induced force-parameterized \hldiff{$\TWO{\EpsStop}{\EpsFJAll}$-FJ-stationarity} condition of~\prettyref{pb:CITO} (\prettyref{def:eps-stationary}), for every stationarity tolerance $\EpsStop>0$, with the feasibility band $\EpsFJAll$ determined by the penalty and stopping tolerances.
\item For fixed $\EpsFJAll>0$, fixed values of all parameters on which $\EpsFJAll$ depends, and a common problem and algorithm setup including the clamps (defined in~\prettyref{sec:Hamiltonian}), the returned trajectory, prescribed constraint forces $\lambda_i^e,\lambda_i^i$, physics multipliers, and penalty weight remain bounded by finite constants uniformly as $\EpsStop\to0$.
\end{enumerate}
\end{theorem}
\end{framed}
\prettyref{thm:stationarity-transfer} is the first-stage result: it assumes the contact dynamics already smooth. Our final result removes that assumption, covering the genuinely nonsmooth frictional contact of general CRDS by resolving the damping through a normal-compatible smoothing whose scale is then sent to zero. To make the difference explicit, the two statements highlight in \hldiff{brown} the two axes along which the second generalizes the first: the dynamics class (\hldiff{smooth} in~\prettyref{thm:stationarity-transfer} versus \hldiff{nonsmooth} in~\prettyref{thm:nonsmooth-transfer}) and, correspondingly, the certified stationarity notion (\hldiff{$\TWO{\EpsStop}{\EpsFJAll}$-FJ-stationarity} sharpened to its homogeneous Clarke refinement \hldiff{C-$\TWO{\EpsStopNS}{\EpsFJAll}$-FJ-stationarity}).
\begin{framed}
\begin{theorem}[(informal)]
\label{thm:nonsmooth-transfer}
Consider~\prettyref{pb:CITO} instantiated by a full-rank (\prettyref{def:full-rank}) CRDS with possibly \hldiff{nonsmooth} frictional-contact damping and controls in the compact convex control space $\mathcal{D}$, resolved through a normal-compatible smoothing (\prettyref{def:S-D-compatible}), and suppose~\prettyref{ass:objective} and~\ref{ass:Slater} hold.
\begin{enumerate}[label=(\alph*),ref=(\alph*),leftmargin=*]
\item Under the same fixed pre-horizon data and finite discretized trajectory initialization class as~\prettyref{thm:stationarity-transfer}, the same line-search SQP with adaptive timestep refinement terminates after finitely many iterations at a trajectory satisfying the penalty-induced force-parameterized homogeneous Clarke--Fritz--John $\TWO{\EpsStopNS}{\EpsFJAll}$ certificate for~\prettyref{pb:CITO} (\prettyref{def:eps-stationary-NS}), for every stationarity tolerance $\EpsStop>0$.
\item For fixed $\EpsFJAll>0$, fixed values of all parameters on which $\EpsFJAll$ depends, and a common problem and algorithm setup including the clamps, the returned trajectory, prescribed constraint forces, and penalty weight remain bounded by finite constants uniformly as $\EpsStop\to0$.
\end{enumerate}
\end{theorem}
\end{framed}
The effective stationarity tolerance $\EpsStopNS$ aggregates the stopping tolerance $\EpsStop$ and the damping-smoothing scale, sharpening to the smooth guarantee ($\EpsStopNS=2\EpsStop$) as the smoothing is removed and vanishing as both are sent to zero. The boundedness statements concern a fixed feasibility band $\EpsFJAll$, fixed values of all parameters on which that band depends, and a common problem and algorithm setup; they make no uniform claim when that fixed setup is varied. It is instructive to contrast~\prettyref{thm:stationarity-transfer} directly with~\prettyref{thm:bcn-infeasibility}. That result certifies a feasible KKT point only when a constraint qualification (MFCQ/LICQ) holds; absent such a CQ it guarantees only convergence to an infeasible stationary point---a stationary point of the infeasibility measure $\|F(\cdot)\|$ with $F\neq0$---so feasibility is not assured a priori, and for the contact-rich MPCC~\prettyref{pb:CITO} the required CQ remains an unverifiable assumption. Our guarantee assumes no CQ a priori: LICQ is instead derived from the curvature structure of the CRDS---a sufficiently small timestep together with the finite-horizon Hamiltonian bound---and the algorithm certifies an approximately feasible, penalty-induced force-parameterized $\TWO{\EpsStop}{\EpsFJAll}$-FJ certificate for~\prettyref{pb:CITO} (\prettyref{def:eps-stationary}) with the uniform bounds of~\prettyref{thm:stationarity-transfer} as $\EpsStop\to0$. In short, where~\prettyref{thm:bcn-infeasibility} may land on a stationary point of the infeasibility measure, exploiting the structure of CRDS converts this into convergence to an approximately feasible penalty-induced force-parameterized certificate without assuming a CQ for~\prettyref{pb:CITO}. A consolidated list of all symbols used in the main text is provided in~\prettyref{table:notation}.
Throughout the paper, ``global convergence'' refers to finite termination of this specified algorithm, from the admissible initialization class of~\prettyref{rem:initialization-scope}, which permits finite optimized configurations and controls in $\mathcal D$ on any finite nonempty ordered optimized set, at an approximate stationarity certificate. It does not mean global optimality of the returned trajectory, and it does not imply an analogous guarantee for a generic SQP implementation.

\subsection{Proof Architecture and Paper Outline}
\label{sec:proof-sketch}
The solver analyzed in this paper is stated informally in~\prettyref{alg:CITO}. Stripped of the \hldiff{highlighted} lines it is a standard line-search SQP method; the highlighted lines are what this paper adds, and they are also the only ones that need justification. Every quantity named loosely there---the working band, the conditioning check, the energy check, and the stopping tolerance---is made precise in~\prettyref{sec:global}, where the algorithm is restated in full.
The proofs of~\prettyref{thm:stationarity-transfer} and~\ref{thm:nonsmooth-transfer} follow a common chain of reductions, which is most easily read backwards from the algorithm above. A line-search SQP method converges only if the constraint Jacobian is full rank at every iterate, and we are not willing to assume that. Each requirement below is what forces the next; each is discharged by one section of the paper, and the first two between them account for all four \hldiff{highlighted} steps.
\begin{enumerate}[label=(\alph*),ref=(\alph*),leftmargin=*]
\item \TE{The SQP step needs LICQ} --- this is what the \hldiff{conditioning check} and the \hldiff{timestep halving} enforce. Without a full-rank constraint Jacobian the quadratic subproblem is degenerate and no convergence guarantee applies. Rather than assume a constraint qualification, \prettyref{sec:global} derives one: the stationarity Jacobian is dominated by the kinetic term, so a sufficiently small timestep makes it full rank (\prettyref{lem:jacobian-lower-bound}); when the check fails at a node, halving that timestep is exactly the remedy the derivation prescribes. The same section also explains the \hldiff{working band}: the step and multiplier bounds behind the guarantee are stated in terms of a violation ceiling, so an arbitrary initialization must first be driven inside it.
\item \TE{That argument needs bounded curvature} --- this is what the \hldiff{energy check} enforces. The kinetic term can dominate only if the curvature of the potential is bounded. \prettyref{sec:Hamiltonian} supplies this from a finite-horizon energy bound, which simultaneously keeps the trajectory norm-bounded (\prettyref{lem:trajectory-bound}); clamping is required because the potential would otherwise develop unbounded negative curvature at contact. Verifying that bound at every node is what keeps the curvature argument valid as the iteration proceeds.
\item \TE{An energy bound needs a smooth, unconstrained problem.} Hard constraints and non-smooth friction admit no Hamiltonian to bound, so \prettyref{sec:penalty} replaces $h^e,h^i$ with smooth barriers, yielding a smooth full-rank unconstrained CRDS (\prettyref{pb:PS-CRDS}/\ref{pb:PSR-CRDS}, \prettyref{cor:SF-CRDS}, \prettyref{def:full-rank}); FJ-transfer lemmas (\prettyref{lem:eps-SR-FJ}) tie approximate minimizers of the penalized problem back to feasibility of~\prettyref{pb:CITO}.
\item \TE{Full rank also needs a constant mass matrix.} \prettyref{sec:CRDS} lifts the dynamics to maximal coordinates so that $M$ is a constant positive-definite matrix (\prettyref{lem:mass}, under the standing~\prettyref{ass:abd-mass}) and the kinetic energy is strongly convex, with curvature inversely proportional to the timestep.
\item \TE{The certificate must be transferred back.} What the loop returns is a KKT point of a surrogate, not of the original problem. \prettyref{sec:global-transfer} restricts the constraint-force variables to the penalty-induced configuration-dependent maps satisfying~\prettyref{eq:transfer-gradient-match}; on the clamp-inactive region this makes the per-step surrogate residual identical to the parameterized discrete~\prettyref{eq:CRDS} force-balance residual, so the certificate transfers to the force-parameterized first-order system of~\prettyref{pb:CITO}---\prettyref{thm:stationarity-transfer}; a graph-lifting argument on the damping subdifferential extends this to approximate Clarke--FJ stationarity for nonsmooth frictional contact---\prettyref{thm:nonsmooth-transfer}, with the articulated-body specialization completed by~\prettyref{thm:abd-C-stationarity}.
\end{enumerate}
Read instead in the order the paper presents them, these requirements become the sequence of problem transformations as summarized as
\[
\adjustbox{max width=\linewidth}{$
\begin{array}{c}
\prettyref{pb:CITO}/\ref{pb:CRDS-ABD-extended}\\[2pt]
\hphantom{{\scriptstyle\ \text{penalize, smooth, reduce}}}\Big\downarrow{\scriptstyle\ \text{penalize, smooth, reduce}}\\[2pt]
\prettyref{pb:PS-CRDS}/\ref{pb:PSR-CRDS}\\[2pt]
\hphantom{{\scriptstyle\ \text{mollify}}}\Big\downarrow{\scriptstyle\ \text{mollify}}\\[2pt]
\prettyref{pb:PS-CRDS-M}\\[2pt]
\hphantom{{\scriptstyle\ \text{stack }\nabla_{\theta_i}\PerStepObj_i=0}}\Big\downarrow{\scriptstyle\ \text{stack }\nabla_{\theta_i}\PerStepObj_i=0}\\[2pt]
\prettyref{pb:CITO-M}\\[2pt]
\hphantom{{\scriptstyle\ \text{Alg.~\ref{alg:sqp},\,\ref{alg:CITO}}}}\Big\downarrow{\scriptstyle\ \text{Alg.~\ref{alg:sqp},\,\ref{alg:CITO}}}\\[2pt]
\EpsStop\text{-KKT}(\prettyref{pb:CITO-M})\\[2pt]
\hphantom{{\scriptstyle\ \text{bridge}+\text{transfer}}}\Big\downarrow{\scriptstyle\ \text{bridge}+\text{transfer}}\\[2pt]
\text{FJ/Clarke--FJ certificate for }\prettyref{pb:CITO}
\end{array}
$}
\]

\section{\label{sec:CRDS}Coordinate Maximization:\\Articulated Body Dynamics as CRDS}

We begin our analysis by showing that articulated body dynamics (ABD) under a general class of geometric representations form a special case of CRDS. In particular, we use the recently proposed separating-plane-based contact model~\cite{11266917} to avoid the geometric singularities under deep penetrations observed in~\cite{schulman2014motion,zhang2023provably}. We further adopt the coordinate maximization technique~\cite{brudigam2021linear,lan2022affine,dai2019global}.

We further introduce equality constraints to ensure that the bodies move as the desired articulated body; these are referred to as integrity constraints. Under this representation, the mass matrix $M$ is positive definite under mild assumptions on the physical rigid bodies, which we adopt as a standing assumption throughout the articulated-body specialization.
\begin{assume}
\label{ass:abd-mass}
The articulated-body model is described by the following standing data, comprising physical mass regularity in~\ref{it:mass-volume}--\ref{it:mass-density} and a finite polyhedral contact representation in~\ref{it:contact-hull}:
\begin{enumerate}[label=(\alph*),ref=(\alph*),leftmargin=*]
\item\label{it:mass-volume} every physical rigid body $j$ occupies a bounded Lebesgue-measurable set $\Omega^j\subset\mathbb{R}^3$ in its local frame with $|\Omega^j|>0$;
\item\label{it:mass-density} the mass density $\rho^j:\Omega^j\to\mathbb{R}$ is Lebesgue measurable and bounded, with $\rho^j(x)>0$ for almost every $x\in\Omega^j$;
\item\label{it:contact-hull} the model comprises finitely many physical bodies, finitely many indexed convex hulls, and finitely many unordered pairs of hulls admitted for contact; the geometry of each body is the union of the convex hulls it carries, and every indexed convex hull is exactly the convex hull of a finite, nonempty family of fixed local-frame vertices $x^{jhk}\in\mathbb{R}^3$, where $j$ indexes the rigid body, $h$ the convex hull carried by that body, and $k$ the vertex within that hull. An admitted pair, together with a single vertex of one of its two hulls, constitutes a vertex incidence (incidence for short), and the incidences are likewise finite in number. These vertex families, the admitted pairs, and their index sets are fixed model data, independent of the physical configuration $\theta$ and of the optimization variables.
\end{enumerate}
\end{assume}
The mass-regularity clauses~\ref{it:mass-volume}--\ref{it:mass-density} apply only to physical rigid bodies; they do not apply to the auxiliary separating-plane variables introduced below, which are deliberately massless (see~\prettyref{eq:extended-config}). These two clauses supply the hypotheses of the following lemma, and they exclude degenerate physical bodies such as massless links and zero-volume point, line, or planar idealizations, unless full rank is established for those cases by a separate model-specific argument. The finite fixed-vertex clause~\ref{it:contact-hull} supplies the finite polyhedral contact representation used throughout the articulated-body development: the separating-plane construction of the normal contact model, the active-contact plane regularity analysis, and the subsequent damping-reduction, smoothing, and stationarity arguments all draw their convex hulls and vertex incidences from this fixed model data. Accordingly,~\prettyref{ass:abd-mass} is a standing model assumption: all articulated-body constructions and their applications in the main text and in the appendix proofs use this model data, whereas standalone analytic lemmas stated for arbitrary Euclidean data retain their own explicit hypotheses.
\begin{lemma}
\label{lem:mass}
Under~\prettyref{ass:abd-mass}, the rigid-body mass matrix, defined as $M=\diag(\{M^j\})$ where $M^j\in\mathbb{R}^{12\times12}$ is the mass matrix of the $j$th rigid body, satisfies $M\succ0$.
\end{lemma}
For the remainder of this section, we introduce the components needed to define the full CRDS describing articulated body dynamics. The physical configuration is
\begin{equation}
\theta\triangleq<\{R^j\},\{t^j\}>.
\end{equation}

\subsubsection{Normal Contact Model}
Each indexed convex hull used in the contact model below is the convex hull of the finite, nonempty family of fixed local-frame vertices $x^{jhk}\in\mathbb{R}^3$ furnished by clause~\ref{it:contact-hull} of~\prettyref{ass:abd-mass}. Each index ranges over a fixed finite set, so every displayed inequality indexed by $k$ or $k'$ is a finite family. We first specify the normal contact model based on the separating-plane construction of~\cite{11266917}, illustrated in~\prettyref{fig:contact-model-illustration}. For each unordered pair of convex hulls $jh$ and $j'h'$, an oriented separating plane with normal $n^{jh,j'h'}$ and offset $b^{jh,j'h'}$ induces a signed distance from every world-frame vertex $X^{j'h'k'}=R^{j'}x^{j'h'k'}+t^{j'}$ to the plane, and non-penetration is enforced by requiring the two convex hulls to lie on opposite sides of the plane
\begin{ResizedEq}
\label{eq:separating-plane}
d^{jh,j'h'k'}\triangleq&\langle n^{jh,j'h'}, X^{j'h'k'} \rangle+b^{jh,j'h'}\\
h^i_{\mathrm{sep}}(\vartheta):&
\begin{cases}
d^{j'h',jhk}-\EpsContactPadding\geq0&\forall k\\
-d^{jh,j'h'k'}-\EpsContactPadding\geq0&\forall k'\\
\end{cases},
\end{ResizedEq}
which uses a small padding $\EpsContactPadding>0$. The padding is a standard technique used by various physics engines to enforce a stricter penetration constraint.

\begin{figure}[t]
\centering
\begin{tikzpicture}[font=\small]

\definecolor{topBoxFill}{RGB}{255,213,170}     
\definecolor{bottomBoxFill}{RGB}{255,239,178}  

\def\bottomBox{
  \fill[bottomBoxFill] (-0.8,-0.8) rectangle (0.8,0.8);
  \draw[black, line width=1.2pt] (-0.8,-0.8) rectangle (0.8,0.8);
}
\def\diamondAt(#1,#2){
  \fill[topBoxFill]
    (#1,#2) ++(0,-0.8) -- ++(0.8,0.8) -- ++(-0.8,0.8) -- ++(-0.8,-0.8) -- cycle;
  \draw[black, line width=1.2pt]
    (#1,#2) ++(0,-0.8) -- ++(0.8,0.8) -- ++(-0.8,0.8) -- ++(-0.8,-0.8) -- cycle;
}

\def\dx{3.8}
\def\dy{5.0}

\def\dAx{0}      \def\dAy{2.0}
\def\sAx{0}      \def\sAy{0.8}
\def\vAx{0}      \def\vAy{1.2}   

\def\dBx{1.1}    \def\dBy{1.9}
\def\sBx{0.8}    \def\sBy{0.8}
\def\vBx{1.1}    \def\vBy{1.1}

\begin{scope}[local bounding box=PriorVE]
  \useasboundingbox (-1.9,-2.0) rectangle (1.9,3.2);
  \bottomBox
  \diamondAt(\dAx,\dAy)
  \draw[red!70!black, dashed, line width=1.0pt] (\sAx,\sAy) -- (\vAx,\vAy);
  \fill[red!70!black] (\sAx,\sAy) circle (1.4pt);
  \fill[red!70!black] (\vAx,\vAy) circle (1.4pt);
  \node[align=center] at (0,-1.4)
    {(a) Prior: vertex--edge pair};
\end{scope}

\begin{scope}[xshift=\dx cm, local bounding box=PriorVV]
  \useasboundingbox (-1.9,-2.0) rectangle (1.9,3.2);
  \bottomBox
  \diamondAt(\dBx,\dBy)
  \draw[red!70!black, dashed, line width=1.0pt] (\sBx,\sBy) -- (\vBx,\vBy);
  \fill[red!70!black] (\sBx,\sBy) circle (1.4pt);
  \fill[red!70!black] (\vBx,\vBy) circle (1.4pt);
  \node[align=center] at (0,-1.4)
    {(b) Prior: vertex--vertex pair};
\end{scope}

\begin{scope}[yshift=-\dy cm, local bounding box=SDRSVE]
  \useasboundingbox (-1.9,-2.0) rectangle (1.9,3.2);

  \def\planeY{1.0}
  \def\planeSlope{0.04}
  \draw[green!55!black, line width=1.4pt]
    (-1.8,{\planeY - 1.8*\planeSlope}) -- (1.8,{\planeY + 1.8*\planeSlope});

  \bottomBox
  \diamondAt(\dAx,\dAy)

  \fill[black] (\vAx,\vAy) circle (1.6pt);
  \node[anchor=west] at (\vAx + 0.45, \vAy + 0.15) {$x^{jhk}$};

  \fill[black] (0.8,0.8) circle (1.6pt);
  \node[anchor=west] at (0.85,0.55) {$x^{j'h'k'}$};

  \node[anchor=south, green!40!black] at (-1.3,{\planeY - 1.3*\planeSlope + 0.02})
    {$p^{jh,j'h'}$};

  \node[align=center] at (0,-1.4)
    {(c) Ours: vertex--plane};
\end{scope}

\begin{scope}[xshift=\dx cm, yshift=-\dy cm, local bounding box=SDRSVV]
  \useasboundingbox (-1.9,-2.0) rectangle (1.9,3.2);

  \def\planeIntercept{1.9}
  \draw[green!55!black, line width=1.4pt]
    (-0.5,{\planeIntercept - (-0.5)}) -- (1.8,{\planeIntercept - 1.8});

  \bottomBox
  \diamondAt(\dBx,\dBy)

  \node[align=center] at (0,-1.4)
    {(d) Ours: vertex--plane};
\end{scope}

\end{tikzpicture}
\caption{\label{fig:contact-model-illustration}Illustration of our separating-plane based 2D contact model on two convex polyhedra (a peach diamond resting on a butter-yellow square). The 3D counterpart can be derived from identical ideas. (a-b) Under prior contact models, the closest feature pair switches abruptly as the bodies slide against each other: a vertex-edge pair in (a) becomes a vertex-vertex pair in (b) after a small tangential motion, yielding a non-smooth contact distance. (c-d) \cite{11266917} introduces an auxiliary separating plane $p^{jh,j'h'}$ (green) between the two convex hulls and replaces both cases with a single signed distance from each vertex (e.g., $x^{jhk}$ on the top body, $x^{j'h'k'}$ on the bottom body) to this plane. The two distinct prior cases are unified into the same vertex-plane contact, which depends smoothly on the configuration of the bodies and the separating plane.}
\end{figure}

\subsubsection{Frictional Contact Model} To introduce frictional damping into CRDS, we use the maximal dissipation principle, which states that the friction force minimizes the rate of work done against contact subject to the Coulomb friction cone~\cite{moreau1988unilateral,goyal1991planar,stewart2000rigid}. We denote by $[\lambda^i]^{j'h',jhk}$ the constraint force applied on $x^{jhk}$ from the separating plane against $j'h'$th convex hull. Again following~\cite{11266917}, we model the separating plane $p^{jh,j'h'}$ as an auxiliary object with zero mass but potentially non-zero linear translation $\varsigma^{jh,j'h'}\in\mathbb{R}^3$ and in-plane rotation angle $\omega^{jh,j'h'}\in\mathbb{R}$.
\begin{remark}
The separating planes and their velocities are auxiliary variables, so the extended configuration is
\begin{ResizedEq}
\label{eq:extended-config}
&\vartheta\triangleq<\theta,\{\THREE{p}{\varsigma}{\omega}^{jh,j'h'}\}_{j\neq j'}>\\
&\dot\vartheta_i\triangleq\frac{\vartheta_i-\vartheta_{\gamma(i)}}{\delta_i}=\THREE{\dot\theta_i}{\dot p_i}{\varpi_i}\\
&\varpi\triangleq<\{\dot\varsigma,\dot\omega\}^{jh,j'h'}_{j\neq j'}>.
\end{ResizedEq}
Here $\dot\vartheta_i$ is the standard componentwise finite difference of the complete configuration $\vartheta_i$. Thus $\vartheta$ is the extended configuration, stacking the physical block $\theta$ together with the auxiliary separating-plane variables; it is $\vartheta$, not $\theta$, that plays the role of the generic CRDS configuration in this specialization. The separating-plane parameter $p=(n,b)$ is zero-mass and reconstructed algebraically by the conservative inner minimization, so its finite-difference block $\dot p_i$ contributes neither inertia nor explicit damping; the damping formula below consumes only the $\dot\theta_i$ and $\varpi_i$ blocks. After substituting the complete ABD configuration $\vartheta$ for the generic CRDS configuration, we intentionally reuse the generic feasible-set symbol
\begin{equation}
\label{eq:abd-extended-feasible-set}
\mathcal F^i\triangleq\{\vartheta\mid h^i(\vartheta)\geq0\}.
\end{equation}
Thus $\mathcal F^i$ is the generic CRDS feasible set evaluated in the complete full-coordinate ABD instantiation, not a second feasible-set construction. We likewise abuse notation for the inequality map: $h^i(\vartheta)$ stacks all physical, separation, and unit-normal inequalities. A physical component that was originally written $h^{i,j}(\theta)$ is interpreted as $h^{i,j}(\vartheta)\triangleq h^{i,j}(\theta)$, hence is constant in the auxiliary coordinates; the separation and unit-normal components retain their full $\vartheta$-dependence.
For each unordered convex-hull pair, $p^{jh,j'h'}$ denotes one oriented separating plane; the pair superscript is unordered, so $d^{j'h',jhk}$ and $d^{jh,j'h'k'}$ share the same plane parameters $n^{jh,j'h'},b^{jh,j'h'}$, and the two inequalities in~\prettyref{eq:separating-plane} encode the two sides of that same plane.
We emphasize that separating planes play a different role from physical configuration variables: they are not physical objects, but auxiliary variables used to express nonpenetration constraints. Therefore, separating planes should not carry mass, and their corresponding diagonal mass-matrix blocks, denoted by $M^{jh,j'h'}=0$, are set to zero as proposed in~\cite{11266917}. As a result, the extended mass matrix $\diag(\{M^j\},\{M^{jh,j'h'}\}_{j\neq j'})$ is not positive definite, while the physical rigid-body mass matrix $M$ in~\prettyref{lem:mass} remains positive definite under the standing~\prettyref{ass:abd-mass}.
\end{remark}

Given the auxiliary separating plane and its translation/rotation, frictional damping between the two convex hulls is delegated to the frictional damping between each convex hull and the separating plane. The relative velocity between $x^{jhk}$ and the separating plane is denoted by
\begin{ResizedEq}
\label{eq:tangent-velocity}
&v_\parallel^{j'h',jhk}=P^{jh,j'h'}\dot{X}^{jhk}-\\
&[\bar{n}^{jh,j'h'}\times X^{jhk}\dot{\omega}^{jh,j'h'}+\dot{\varsigma}^{jh,j'h'}],
\end{ResizedEq}
where, whenever $\vartheta\in\mathcal{F}^i$ as defined in~\prettyref{eq:abd-extended-feasible-set}, the separating-plane inequalities imply $n^{jh,j'h'}\neq0$, so we define $\bar{n}^{jh,j'h'}=n^{jh,j'h'}/\|n^{jh,j'h'}\|$ and the tangential projection matrix $P^{jh,j'h'}=I-[\bar{n}^{jh,j'h'}][\bar{n}^{jh,j'h'}]^T$. On this domain, with a fixed plane-velocity damping coefficient $\EpsDamp>0$, we define the damping term $\Damp(\vartheta,\dot\vartheta,\lambda^i)$, with the ambient CRDS velocity signature, as the sum of relative velocity norms scaled by normal force magnitudes and the friction coefficient $\eta>0$. Throughout this section and its proofs, a sum over a hull pair $jh,j'h'$ (or a joint pair $j,j'$ for $W$) is taken over unordered pairs, each counted once.
\begin{ResizedEq}
\label{eq:damping-ABD}
&\Damp(\vartheta,\dot\vartheta,\lambda^i)=\EpsDamp\sum_{jh,j'h'}\left(\|\dot{\omega}^{jh,j'h'}\|^2+\|\dot{\varsigma}^{jh,j'h'}\|^2\right)+\\
&\eta\sum_{jh,j'h'}\sum_k[\lambda^i]^{j'h',jhk}\|v_\parallel^{j'h',jhk}\|+\\
&\eta\sum_{jh,j'h'}\sum_{k'}[\lambda^i]^{jh,j'h'k'}\|v_\parallel^{jh,j'h'k'}\|.
\end{ResizedEq}
Here $[\lambda^i]^{jh,j'h'k'}$ is the normal contact force imposed by $p^{jh,j'h'}$ on $x^{j'h'k'}$ and vice versa. The first term is a quadratic damper on the separating planes' own rotational and translational velocities; since the planes are massless auxiliary variables, it regularizes their motion. The remaining two terms are the frictional dissipation of the two hulls against the plane, one per side of the pair. Each has the Coulomb form: the product of the friction coefficient $\eta$, the normal contact force, and the tangential sliding speed is precisely the rate at which a Coulomb friction force of magnitude $\eta[\lambda^i]$ opposing the sliding direction dissipates energy, so minimizing the sum over admissible velocities realizes the maximal dissipation principle stated above. The difference from a textbook Coulomb model is that friction acts between each hull and the intermediate separating plane rather than directly between the two hulls, which is what keeps the resulting constraints differentiable. Note that the damping function $\Damp$ satisfies the requirement in~\prettyref{def:CRDS}. For fixed $\vartheta_{\gamma(i)}\in\mathcal{F}^i$ and nonnegative normal contact forces, it is a non-smooth convex function of the current full velocity $\dot\vartheta_i$: it is convex in the damping-active blocks $\TWO{\dot\theta_i}{\varpi_i}$ and constant in $\dot p_i$. Moreover, $\Damp$ is a sum of nonnegative norm terms scaled by the nonnegative normal contact forces $[\lambda^i]\geq0$ plus the plane-velocity damper with positive coefficient $\EpsDamp\sum_{jh,j'h'}(\|\dot{\omega}^{jh,j'h'}\|^2+\|\dot{\varsigma}^{jh,j'h'}\|^2)$, so $\Damp\geq0$. Since the tangential relative velocities $v_\parallel$, the angular velocities $\dot{\omega}^{jh,j'h'}$, and the translational velocities $\dot{\varsigma}^{jh,j'h'}$ are all linear in $\TWO{\dot\theta_i}{\varpi_i}$, they vanish at $\dot\vartheta_i=0$, so $\Damp(\vartheta_{\gamma(i)},0,\lambda^i)=0$ attains the zero-velocity global minimum required in~\prettyref{def:CRDS}~\ref{it:crds-damping}. The non-smoothness is unavoidable because the friction force changes discontinuously at the transition between static and sliding modes. When the relative velocity $v_\parallel^{j'h',jhk}$ is non-zero, the frictional force has magnitude equal to $\eta$ times the normal contact force. Otherwise, the frictional force can take an arbitrary value in the frictional cone, which corresponds to the subgradient set of the damping function $\Damp$. Further note that in~\prettyref{def:CRDS} we discretize the damping as $\Damp(\theta_{\gamma(i)},\dot{\theta}_i,\lambda_{\gamma(i)}^i)$, which corresponds to a semi-implicit damping scheme where predecessor geometry and predecessor normal load weight the current interval velocity. In addition to the velocity dissipation term, we introduce a plane-velocity damper $\EpsDamp\sum_{jh,j'h'}(\|\dot{\omega}^{jh,j'h'}\|^2+\|\dot{\varsigma}^{jh,j'h'}\|^2)$ on both the translational and rotational velocity of each separating plane. We will show that this term is essential to ensure the regularity condition in the theoretical analysis of the next section. In practice, its effect can be controlled by an arbitrarily small positive $\EpsDamp$.

\subsubsection{Integrity Constraints}
Throughout the paper, we consider articulated bodies with two types of integrity constraints: rotation constraints and joint constraints induced by hinge or prismatic joints. For each rigid body we represent its rotation with a general $3\times3$ matrix that may deviate from a valid rotation matrix; we therefore introduce the following 9D rotation constraint
\begin{equation}
\label{eq:integrity}
h^e(\theta):(R^j)^TR^j-I=0.
\end{equation}
For both hinge and prismatic joints, we define $\TWO{R}{t}^{jj'}$ and $\TWO{R}{t}^{j'j}$ as the local transforms from the $j$th and $j'$th local frames to their link-attached constraint frames, respectively, as illustrated in~\prettyref{fig:joint-frame-illustration}. In the joint constraints below, $x,y,z$ denote the canonical orthonormal basis vectors of the constraint frame, with $z$ the hinge rotation axis (prismatic sliding axis). These transforms determine the two world-frame joint anchors $R^jt^{jj'}+t^j$ and $R^{j'}t^{j'j}+t^{j'}$, whose relative displacement is the dashed segment in the figure. First, assume that two rigid bodies, indexed by $j$ and $j'$, are connected by a hinge joint. We assume that the two rigid bodies can only undergo relative rotation about the $z$-axis in the constraint frame, with the relative rotation angle bounded by $[-L^{jj'},L^{jj'}]$. Such a joint can be expressed by the following equality and inequality constraints.
\begin{ResizedEq}
\label{eq:hinge}
h^e(\theta):&\begin{cases}
(R^jt^{jj'}+t^j)-(R^{j'}t^{j'j}+t^{j'})=0\\
(R^jR^{jj'}-R^{j'}R^{j'j})z=0\\
\end{cases}\\
h^i(\theta):&\langle R^jR^{jj'}x, R^{j'}R^{j'j}x \rangle - \cos(L^{jj'}) \geq 0.
\end{ResizedEq}

\begin{figure}[t]
\centering
\begin{tikzpicture}[font=\tiny, line cap=round, line join=round]

\definecolor{topBoxFill}{RGB}{255,213,170}     
\definecolor{bottomBoxFill}{RGB}{255,239,178}  
\definecolor{axisRed}{RGB}{190,65,55}
\definecolor{axisBlue}{RGB}{45,95,175}
\definecolor{axisGreen}{RGB}{72,145,83}

\tikzset{
  frameAxis/.style={line width=0.95pt, -stealth},
  capsuleOutline/.style={black, line width=1.15pt},
  anchorPoint/.style={circle, fill=black, inner sep=1.25pt},
}

\def\linkRadius{0.32}
\def\halfLink{1.425}

\def\drawCapsule#1#2#3#4{
  \begin{scope}[shift={(#1,#2)}, rotate=#3]
    \fill[#4] (-\halfLink,-\linkRadius) -- (\halfLink,-\linkRadius)
      arc[start angle=-90, end angle=90, radius=\linkRadius]
      -- (-\halfLink,\linkRadius)
      arc[start angle=90, end angle=270, radius=\linkRadius] -- cycle;
    \draw[capsuleOutline] (-\halfLink,-\linkRadius) -- (\halfLink,-\linkRadius)
      arc[start angle=-90, end angle=90, radius=\linkRadius]
      -- (-\halfLink,\linkRadius)
      arc[start angle=90, end angle=270, radius=\linkRadius] -- cycle;
  \end{scope}
}

\def\drawFrame#1#2#3{
  \begin{scope}[shift={(#1)}, rotate=#2]
    \draw[frameAxis, axisRed] (0,0) -- (#3,0);
    \draw[frameAxis, axisBlue] (0,0) -- (0,#3);
    \draw[frameAxis, axisGreen] (0,0) -- (-0.23,0.31);
  \end{scope}
}

\coordinate (leftCenter) at (-2.45,0.22);
\coordinate (rightCenter) at (2.45,-0.18);
\coordinate (leftTip) at (-1.067,-0.124);
\coordinate (rightTip) at (1.087,-0.573);

\drawCapsule{-2.45}{0.22}{-14}{topBoxFill}
\drawCapsule{2.45}{-0.18}{16}{bottomBoxFill}

\drawFrame{leftCenter}{-14}{0.50}
\drawFrame{rightCenter}{196}{0.50}
\drawFrame{leftTip}{-14}{0.43}
\drawFrame{rightTip}{196}{0.43}

\node[inner sep=1.2pt, anchor=south] at (-3.1,0.2) {$(R^j,t^j)$};
\node[inner sep=1.2pt, anchor=south] at (3.1,-0.15) {$(R^{j'},t^{j'})$};
\node[fill=white, inner sep=1.2pt, align=center, anchor=north] at (-1.10,-0.58)
  {$(R^jR^{jj'},$\\$R^jt^{jj'}+t^j)$};
\node[fill=white, inner sep=1.2pt, align=center, anchor=north] at (1.10,-1.03)
  {$(R^{j'}R^{j'j},$\\$R^{j'}t^{j'j}+t^{j'})$};

\node[anchorPoint] at (leftTip) {};
\node[anchorPoint] at (rightTip) {};

\draw[red!70!black, dashed, line width=1.05pt] (leftTip) -- (rightTip);

\begin{scope}[shift={(0.01,-0.35)}]
  \fill[axisGreen!18] (0,0) -- (42:1.02)
    arc[start angle=42, end angle=138, radius=1.02] -- cycle;
  \draw[axisGreen!60!black, line width=0.55pt] (42:1.02)
    arc[start angle=42, end angle=138, radius=1.02];
  \draw[axisGreen!60!black, line width=0.45pt] (0,0) -- (42:1.02);
  \draw[axisGreen!60!black, line width=0.45pt] (0,0) -- (138:1.02);
  \node[axisGreen!55!black, inner sep=0.6pt, anchor=south east] at (122:0.9)
    {$-L^{jj'}$};
  \node[axisGreen!55!black, inner sep=0.6pt, anchor=south west] at (48:1.03)
    {$L^{jj'}$};

  \draw[black!45, line width=0.45pt] (0,0) -- (66:0.88);
  \draw[black!45, line width=0.45pt] (0,0) -- (106:0.88);
  \draw[black, line width=0.95pt, -stealth] (66:0.52)
    arc[start angle=66, end angle=106, radius=0.52];
  \node[black, inner sep=1.0pt, anchor=south] at (83:0.65)
    {$\alpha^{jj'}$};
\end{scope}

\node[red!70!black, inner sep=1.0pt, anchor=north] at (0.1,-0.44)
  {$\Delta^{jj'}$};

\end{tikzpicture}
\caption{\label{fig:joint-frame-illustration}Illustration of the shared maximal-coordinate joint-frame representation used for both hinge and prismatic joints, with two capsule-shaped robot links using the same peach and butter-yellow color convention as~\prettyref{fig:contact-model-illustration}. The link-centered frames $(R^j,t^j)$ and $(R^{j'},t^{j'})$ are drawn at the centers of the two links. The local transforms $\TWO{R}{t}^{jj'}$ and $\TWO{R}{t}^{j'j}$ define two link-attached constraint frames, whose world-frame anchor positions are $R^jt^{jj'}+t^j$ and $R^{j'}t^{j'j}+t^{j'}$. The dashed red segment labeled $\Delta^{jj'}$ denotes their relative anchor displacement: hinge constraints drive this displacement to zero while constraining the relative orientation to a bounded rotation about the constraint-frame $z$ axis, whereas prismatic constraints express this displacement in the constraint frame and allow only its $z$ component to vary within $[-L^{jj'},L^{jj'}]$. The green sector illustrates the common joint-limit notation $[-L^{jj'},L^{jj'}]$ and the rotational coordinate $\alpha^{jj'}$ for the hinge case.}
\end{figure}

Now suppose the two rigid bodies are connected by a prismatic joint. We reuse the same notation for transformations from local frames to the constraint frame.
\begin{ResizedEq}
\label{eq:prismatic}
\Delta^{jj'}\triangleq&(R^jR^{jj'})^T
\left[(R^{j'}t^{j'j}+t^{j'})-(R^jt^{jj'}+t^j)\right]\\
h^e(\theta):&\begin{cases}
\langle\Delta^{jj'},x\rangle=0\\
\langle\Delta^{jj'},y\rangle=0\\
R^jR^{jj'}-R^{j'}R^{j'j}=0\\
\end{cases}\\
h^i(\theta):&\begin{cases}
L^{jj'}-\langle\Delta^{jj'},z\rangle\geq0\\
L^{jj'}+\langle\Delta^{jj'},z\rangle\geq0\\
\end{cases}.
\end{ResizedEq}
Finally, for the auxiliary separating planes, we can enforce unit-length normals through the constraint $\|n^{jh,j'h'}\|-1=0$. However, such a unit-length constraint is non-convex. Instead, we introduce the following smooth convexified unit-normal constraint
\begin{ResizedEq}
\label{eq:unit-normal}
h^i_{\mathrm{unit}}(\vartheta): 1-\|n^{jh,j'h'}\|^2\geq0.
\end{ResizedEq}
The following result shows that the convexified unit-normal constraint can equivalently enforce penetration-free constraints between convex polyhedra.
\begin{lemma}
\label{lem:separate}
The distance between two convex polyhedra, indexed by $jh$ and $j'h'$, is at least $2\EpsContactPadding$ if and only if there exists $p^{jh,j'h'}$ such that~\prettyref{eq:separating-plane} and~\ref{eq:unit-normal} hold.
\end{lemma}

\subsubsection{Control Input}
The final building block is the control input. The physical integrity and contact constraints above are still written in the actual maximal coordinates. We use the following clamp (mollification) map $\mathcal{M}$ to replace each rotation appearing in the controlled joint coordinate by a bounded surrogate and, for prismatic joints, to replace the relative anchor displacement by a bounded surrogate. Bounding these quantities is what keeps the curvature of the potential energy uniformly bounded: without clamping, the potential develops unbounded negative curvature where bodies come into contact, and neither the Hamiltonian upper bound nor the LICQ derived from it would follow. For vectors, $\|\cdot\|$ denotes the Euclidean norm; for matrices, it denotes the Frobenius norm, consistent with vectorizing the matrix. Throughout, $\mathbf 1$ denotes the all-ones element of the ambient space of its argument---a scalar, vector, or matrix as appropriate---so the $\|x\|=\infty$ branch returns an object of the same type as $x$ with norm $2L$.
\begin{ResizedEq}
\label{eq:clamp-map}
\mathcal{M}(x,L)\triangleq&
\begin{cases}
x&\|x\|\leq L\\
(1-p(t))x+p(t)2L\frac{x}{\|x\|},\quad t=\frac{\|x\|}{L}-1&L<\|x\|\leq2L\\
2L\frac{x}{\|x\|}&\|x\|>2L\\
2L\frac{\mathbf 1}{\|\mathbf 1\|}&\|x\|=\infty,
\end{cases}\\
p(t)\triangleq&
\begin{cases}
0&t\leq0,\\
35t^4-84t^5+70t^6-20t^7&0<t<1,\\
1&t\geq1.
\end{cases}
\end{ResizedEq}
We then define, for a fixed position clamp $\LPosition>0$, the clamped rotation block
\begin{equation}
\bar R^j\triangleq\mathcal M(R^j,\LPosition),
\end{equation}
which satisfies $\|\bar R^j\|\leq2\LPosition$ by~\prettyref{eq:clamp-map}. The numerical value of $\LPosition$ is not needed here: it is fixed later, once the trajectory bound is available, to exceed the norm of every configuration along the optimized trajectory by a strict margin, so the clamp is inactive there and the mollified blocks agree with the originals. This mollification does not change the articulated-body feasibility constraints; it only prevents the control-work Jacobian from depending on unbounded maximal-coordinate rotations outside the physically relevant window.

We again consider two types of joints. First, suppose the $j$th and $j'$th rigid bodies are connected by a hinge joint. The control input $u^{jj'}$ is the external torque applied to the joint angle $\alpha^{jj'}$. We use $\operatorname{atan2}_{\EpsAtan}$ to denote a smoothed version of $\operatorname{atan2}$ with smoothing width $\EpsAtan$; its explicit construction is given in the proof of~\prettyref{cor:abd-crds}. The work done by this torque is
\begin{equation}
\label{eq:hinge-control}
\begin{aligned}
\alpha^{jj'}(\theta)=&\operatorname{atan2}_{\EpsAtan}\TWOC
{\langle \bar R^jR^{jj'}y, \bar R^{j'}R^{j'j}x \rangle}
{\langle \bar R^jR^{jj'}x, \bar R^{j'}R^{j'j}x \rangle}\\
w^{jj'}(\theta,u^{jj'})=&-u^{jj'}\alpha^{jj'}(\theta).
\end{aligned}
\end{equation}
Now suppose the two rigid bodies are connected by a prismatic joint. The physical displacement $\Delta^{jj'}$ in~\prettyref{eq:prismatic} remains the integrity-constraint displacement. For control work, we first define the mollified relative anchor displacement
\begin{ResizedEq}
\label{eq:clamped-anchor-displacement}
\bar d^{jj'}\triangleq&\mathcal M\!\left((\bar R^{j'}t^{j'j}+t^{j'})-(\bar R^jt^{jj'}+t^j),\LPosition\right).
\end{ResizedEq}
By~\prettyref{eq:clamp-map}, it satisfies $\|\bar d^{jj'}\|\leq2\LPosition$. We choose the scale so that $\LPosition\geq\max_{jj'}L^{jj'}$, which makes this displacement mollification inactive for feasible prismatic joint displacements within their joint limits. The control-only prismatic displacement is then
\begin{ResizedEq}
\label{eq:clamped-prismatic-control-displacement}
\bar\Delta^{jj'}\triangleq&(\bar R^jR^{jj'})^T\bar d^{jj'}.
\end{ResizedEq}
The control input $u^{jj'}$ is the external force, and the work done by this force is
\begin{equation}
w^{jj'}(\theta,u^{jj'})=-u^{jj'}\langle \bar\Delta^{jj'},z \rangle.
\end{equation}
The corresponding control signal is stacked as $u\triangleq<\{u^{jj'}\}>$, with the per-timestep, per-joint control vector $u_i\triangleq\{u_i^{jj'}\}$. These per-joint work terms define the aggregate articulated-body control work potential
\begin{equation}
\label{eq:abd-control-work}
W(\theta,u)\triangleq\sum_{j,j'}w^{jj'}(\theta,u^{jj'}),
\end{equation}
and the control force is then defined through the duality between work and force as
\begin{equation}
\label{eq:control}
c(\theta,u)=-\FPPR{W(\theta,u)}{\theta}.
\end{equation}
This construction conforms to the definition of CRDS.
\begin{corollary}
\label{cor:abd-crds}
The articulated body dynamics under maximal coordinates, as defined by \prettyref{eq:control},\ref{eq:integrity},\ref{eq:hinge},\ref{eq:prismatic},\ref{eq:separating-plane},\ref{eq:unit-normal} with the frictional dissipation in~\prettyref{eq:tangent-velocity},\ref{eq:damping-ABD} and setting $U=0$, is a CRDS. In this specialization, the complete configuration $\vartheta$ plays the role of the generic CRDS configuration $\theta$, while $\theta$ is reserved for the reduced physical maximal-coordinate block.
\end{corollary}
The proofs of~\prettyref{lem:mass},~\ref{lem:separate}, and~\prettyref{cor:abd-crds} are collected in~\prettyref{appen:coordinate-proofs}.

\section{\label{sec:penalty}Penalization:\\CRDS as Unconstrained Optimization}
In~\prettyref{sec:CRDS}, we have seen that CRDS is a general class of dynamic systems that can describe constraint-rich articulated body dynamics. The original trajectory program~\prettyref{pb:CITO} enforces the exact discrete CRDS dynamics of~\prettyref{eq:CRDS}, together with hard equality, inequality, and complementarity structure. We do not apply trajectory-level SQP directly to this formulation. Instead, we first derive a per-step constrained variational representation of the dynamics and then absorb the hard constraints into penalty potentials. Subsequent smoothing, reduction, and mollification will turn this per-step model into the smooth residual used by the full-trajectory NLP. This section carries out the first part of that construction. By first-order optimality,~\prettyref{eq:CRDS} is represented at one timestep by the following constrained problem
\begin{ResizedEq}
\label{pb:CRDS-ABD}
\argmin{\theta_i}&
\Phi_i(\theta_i)+\delta_i\Damp(\theta_{\gamma(i)},\dot{\theta}_i,\lambda_{\gamma(i)}^i)\\
\ST\quad&h^e(\theta_i)=0,
h^i(\theta_i)\geq0.
\end{ResizedEq}
We denote the smooth objective terms by
\begin{ResizedEq}
\Phi_i(\theta_i)\triangleq
\frac{\delta_i^2}{2}\|\ddot{\theta}_i\|_M^2+U(\theta_i)-g^T\theta_i+W(\theta_i,u_i),
\end{ResizedEq}
where $\|\cdot\|_M^2\triangleq(\cdot)^TM(\cdot)$ denotes the mass-weighted norm induced by the constant positive-definite mass matrix $M$, $U(\theta_i)$ is the conservative potential from~\prettyref{def:CRDS}, and $W(\theta_i,u_i)$ is the aggregate control work potential from~\prettyref{def:CRDS}. For brevity $\Phi_i$ is written with the single argument $\theta_i$, although through $\ddot{\theta}_i$ and $W$ it also depends on $\theta_{\gamma(i)},\theta_{\gamma^2(i)}$ and the control $u_i$, exactly the dependence made explicit in the per-step objective $\PerStepObj_i$ of the mollified per-step CRDS. Next, we adopt the standard penalty method and introduce barrier-type penalty functions for equality and inequality constraints, as illustrated in~\prettyref{fig:penalty-illustration} (equality) and~\prettyref{fig:penalty-illustration-ineq} (inequality). We define the extended real function
\begin{equation}
\label{eq:Ps}
V_s(x,\EpsFJ)=
\begin{cases}
\infty & x\leq0\\
(x-\EpsFJ)^4/x^5 & x\leq\EpsFJ\\
0 & \text{otherwise}
\end{cases},
\end{equation}
with $\EpsFJ\in(0,1)$ being the barrier parameter; primes on $V_s$ denote derivatives with respect to its first argument, with $\EpsFJ$ held fixed. Because $V_s(x,\EpsFJ)=+\infty$ for $x\leq0$, finiteness of the inequality penalty requires strict slack $h^{i,j}(\theta)>0$, not merely feasibility $h^{i,j}(\theta)\geq0$. This strict finite-domain property is why~\prettyref{ass:Slater} excludes an active-contact predecessor. We then introduce the following soft penalty functions for the generic CRDS constraint maps.
\begin{equation}
\label{eq:soft-penalty}
\begin{aligned}
V^{e,j}(\theta_i,\EpsFJ)=
&V_s(\EpsFJ-h^{e,j}(\theta_i),\EpsFJ)+\\
&V_s(\EpsFJ+h^{e,j}(\theta_i),\EpsFJ),\\
V^{i,j}(\theta_i,\EpsFJ)=&V_s(h^{i,j}(\theta_i),\EpsFJ),
\end{aligned}
\end{equation}
where superscript $j$ denotes the $j$th component of each constraint map. We denote the total penalty potential by
\begin{equation}
\label{eq:penalty-potential}
V(\theta_i,\EpsFJ)=
\sum_jV^{e,j}(\theta_i,\EpsFJ)+\sum_jV^{i,j}(\theta_i,\EpsFJ),
\end{equation}
and consider the following unconstrained optimization problem.
\begin{ResizedEq}
\label{pb:P-CRDS}
\argmin{\theta_i}\Phi_i(\theta_i)+V(\theta_i,\EpsFJ)+\delta_i\Damp^V(\theta_{\gamma(i)},\dot{\theta}_i).
\end{ResizedEq}
Here the hard equality and inequality constraints are absorbed into the added penalty potential $V$, while the original conservative potential remains inside $\Phi_i$. Thus the resulting unconstrained CRDS has empty hard constraint maps and conservative potential equal to the original $U$ plus the barrier contribution represented by $V$. The term $\Damp^V$ is the penalty-induced dissipation term associated with the penalized inequality forces. It is a finite real-valued map $\Damp^V:\operatorname{dom}V\times\mathbb R^n\to\mathbb R$, whose first argument is the predecessor configuration $\theta_{\gamma(i)}$ and whose second argument is the current velocity $\dot{\theta}_i$. On this domain it is given by the natural construction $\Damp^V(\theta_{\gamma(i)},\dot{\theta}_i)=\Damp(\theta_{\gamma(i)},\dot{\theta}_i,\Lambda^i)$, where $\Lambda^i\geq0$ is the penalty-induced inequality force generated by the barrier $V$ at the previous configuration $\theta_{\gamma(i)}$. The penalty force $\Lambda^i$ is evaluated only where the predecessor penalty is finite, i.e., for $\theta_{\gamma(i)}\in\operatorname{dom}V$, and is defined componentwise from the inequality penalty derivative by
\begin{equation}
\label{eq:penalty-force-generic}
[\Lambda^i]^j=-V_s'(h^{i,j}(\theta_{\gamma(i)}),\EpsFJ)\geq0;
\end{equation}
At the first optimized node $i_1\triangleq\min\mathcal I$, this construction is evaluated at the fixed predecessor $\theta_0$. Assumption~\ref{ass:Slater} ensures $\theta_0\in\operatorname{dom}V$, so the predecessor penalty force and semi-implicit dissipation load are finite and well defined. Note that~\prettyref{pb:P-CRDS} defines a CRDS with no additional equality or inequality constraints. For every fixed predecessor configuration $\theta_{\gamma(i)}\in\operatorname{dom}V$, its dissipation $\Damp^V(\theta_{\gamma(i)},\dot{\theta}_i)=\Damp(\theta_{\gamma(i)},\dot{\theta}_i,\Lambda^i)$ inherits the properties of~\prettyref{def:CRDS}~\ref{it:crds-damping} since $\Lambda^i\geq0$; in particular, on $\operatorname{dom}V$ it is convex and nonnegative in velocity and attains its global minimum value $0$ at zero relative velocity. These properties are asserted only on $\operatorname{dom}V$, where $\Lambda^i$ and $\Damp^V$ are defined. The smoothed dissipation introduced below in~\prettyref{def:S-D} is a separately defined, everywhere-defined global extension: outside $\operatorname{dom}V$ it evaluates neither $\Lambda^i$ nor $\Damp^V$ and carries no required physical penalty-force interpretation.

\begin{figure}[t]
\centering
\begin{tikzpicture}[font=\small]

\def\penaltyR(#1){(1.6 + 0.16*sin(3*#1-25) + 0.09*cos(5*#1+55) + 0.05*sin(7*#1+10))}

\begin{scope}[local bounding box=Hard]
  \useasboundingbox (-2,-2.6) rectangle (2,2);
  \draw[black, line width=1.4pt] (0,0) circle (1.6);
  \fill[red!85!black] ({1.6*cos(55)},{1.6*sin(55)}) circle (0.14);
  \node[align=center] at (0,-2.4)
    {(a) Hard manifold constraint\\$h^{e,0}(\theta)=0$};
\end{scope}

\begin{scope}[xshift=4.1cm, local bounding box=Soft]
  \useasboundingbox (-2,-2.6) rectangle (2,2);

  \begin{scope}
    \clip (-2,-2) rectangle (2,2);
    \foreach \i in {200,199,...,1} {
      \pgfmathsetmacro{\rr}{2.85*\i/200}
      \pgfmathsetmacro{\dd}{min(1, abs(\rr-1.6)/1.25)}
      \pgfmathtruncatemacro{\pct}{round(100*\dd)}
      \fill[red!\pct!blue] (0,0) circle (\rr cm);
    }
  \end{scope}

  \draw[black, line width=1.4pt] (0,0) circle (1.6);

  \draw[green!55!black, line width=1.6pt, smooth, samples=260,
        domain=55:325]
    plot ({\penaltyR(\x)*cos(\x)}, {\penaltyR(\x)*sin(\x)});

  \pgfmathsetmacro{\rstart}{\penaltyR(55)}
  \fill[green!55!black] ({\rstart*cos(55)},{\rstart*sin(55)}) circle (0.14);

  \node[align=center] at (0,-2.4)
    {(b) Soft penalty relaxation\\$U^{e,0}(\theta,\EpsFJ)$};
\end{scope}

\end{tikzpicture}
\caption{\label{fig:penalty-illustration}Illustration of the penalty relaxation underlying our CRDS formulation. (a) The original constraint manifold $h^{e,0}(\theta)=0$ is depicted as the black circle and the configuration (red ball) is restricted to move exactly on it. Enforcing this hard constraint requires constraint qualifications (e.g., LICQ/MFCQ) that may fail for contact-rich systems, and the feasible set itself can be non-smooth. (b) Our soft penalty $V^{e,0}(\theta,\EpsFJ)$ relaxes this restriction by penalizing deviation from the manifold; the color code shows the value of the penalty in the ambient space, with blue indicating low values (near the manifold) and red indicating high values (far from it). Under this relaxation the optimized configuration (green ball) follows a trajectory (green curve) that stays close to the original manifold but is allowed to drift slightly inside or outside it. The deviation can be made arbitrarily small by reducing the barrier parameter $\EpsFJ$, while the resulting unconstrained \prettyref{pb:P-CRDS} remains smooth and differentiable.}
\end{figure}

\begin{figure}[t]
\centering
\begin{tikzpicture}[font=\small]

\def\xMin{-2}
\def\xMax{2}
\def\yMax{2.5}
\def\bandH{0.55}    

\begin{scope}[local bounding box=Hard]
  \useasboundingbox (\xMin,-1.0) rectangle (\xMax,\yMax+0.3);

  \draw[black, line width=1.4pt] (\xMin-0.15,0) -- (\xMax+0.15,0);

  \draw[red!85!black, line width=1.6pt, smooth, samples=80, domain=\xMin:0]
    plot (\x, {\yMax - (\yMax/4)*(\x+2)*(\x+2)});

  \draw[red!85!black, line width=1.6pt, smooth, samples=80, domain=0:\xMax]
    plot (\x, {\yMax - (\yMax/4)*(\x-2)*(\x-2)});

  \fill[red!85!black] (\xMin, \yMax) circle (0.14);

  \node[align=center] at (0,-0.7)
    {(a) Hard inequality constraint\\$h^{i,0}(\theta)\geq 0$};
\end{scope}

\begin{scope}[xshift=4.4cm, local bounding box=Soft]
  \useasboundingbox (\xMin,-1.0) rectangle (\xMax,\yMax+0.3);

  \begin{scope}
    \clip (\xMin-0.15,0) rectangle (\xMax+0.15,\bandH);
    \foreach \i in {0,1,...,100} {
      \pgfmathsetmacro{\yy}{\bandH*\i/100}
      \pgfmathsetmacro{\yyn}{\bandH*(\i+1)/100}
      \pgfmathtruncatemacro{\pct}{\i}
      \fill[blue!\pct!red] (\xMin-0.15, \yy) rectangle (\xMax+0.15, \yyn);
    }
  \end{scope}

  \draw[black, line width=1.4pt] (\xMin-0.15,0) -- (\xMax+0.15,0);

  \def\epsSmooth{0.15}
  \def\yHard(#1){(\yMax - (\yMax/4)*(2-abs(#1))*(2-abs(#1)))}
  \def\ySoft(#1){((\yHard(#1) + sqrt(\yHard(#1)*\yHard(#1) + 4*\epsSmooth*\epsSmooth))/2)}
  \draw[green!55!black, line width=1.6pt, smooth, samples=200,
        domain=\xMin:\xMax]
    plot (\x, {\ySoft(\x)});

  \fill[green!55!black] (\xMin, \yMax) circle (0.14);

  \node[align=center] at (0,-0.7)
    {(b) Soft penalty relaxation\\$U^{i,0}(\theta,\EpsFJ)$};
\end{scope}

\end{tikzpicture}
\caption{\label{fig:penalty-illustration-ineq}Illustration of the penalty relaxation of an inequality constraint, using a ball-on-table contact as a canonical example. (a) Under the hard inequality $h^{i,0}(\theta)\geq 0$, the table top is a sharp boundary of the feasible set: the red ball falls, touches the table exactly, and bounces back. The resulting trajectory has a non-smooth cusp at the contact instant. (b) Our soft penalty $V^{i,0}(\theta,\EpsFJ)$ replaces the rigid contact with a smoothly increasing barrier supported on a thin slab of thickness $\EpsFJ$ above the table; the color code visualizes the value of $V_s$ in this slab, with red indicating large penalty (near the table) and blue indicating vanishing penalty (at the top of the slab, where $V_s$ smoothly meets zero). The green ball never touches the table: its trajectory is a smooth quartic that dips into the colored slab and is repelled back, exactly as if a soft spring were placed between the ball and the table. Decreasing $\EpsFJ$ shrinks the slab thickness and recovers the hard constraint in the limit, while keeping the underlying \prettyref{pb:P-CRDS} unconstrained and differentiable.}
\end{figure}

\subsection{Approximate CRDS}
We introduce the following definition of the $\EpsFJ$-FJ perturbed optimality condition. This is a per-step physical condition for the constrained CRDS dynamics: it records approximate force balance, feasibility, and complementarity for the current step. It is not a stationarity condition for the full user trajectory objective $O(\{\theta_i\}_{i\in\mathcal I},\{u_i\}_{i\in\mathcal I})$; that objective enters later through the KKT system of the full-trajectory NLP and the trajectory-level transfers.
\begin{definition}[$\EpsFJ$-FJ conditions]
\label{def:eps-FJ}
A point $\theta_i^\star$ satisfies the $\EpsFJ$-FJ condition of the CRDS defined by~\prettyref{pb:CRDS-ABD} if there exist multipliers $\lambda_i^e$, $\lambda_i^i\geq0$, and $\lambda_{\gamma(i)}^i\geq0$, as well as another velocity $\dot\theta_i^{\dagger}$, such that
\begin{ResizedEq}
\label{eq:eps-fj-system}
&\EpsFJ\text{-FJ}:\\
&\begin{cases}
\begin{aligned}
\dist\!\Big(&-\nabla_{\theta_i}\Phi_i(\theta_i^\star)-\FPP{h^e}{\theta_i}(\theta_i^\star)^T\lambda_i^e+
\FPP{h^i}{\theta_i}(\theta_i^\star)^T\lambda_i^i,\\
&\partial_{\theta_i}\left(\delta_i\Damp\right)(\theta_{\gamma(i)}^\star,\dot{\theta}_i^{\dagger},\lambda_{\gamma(i)}^i)\Big)\leq\EpsFJ
\end{aligned}\\
\|\dot{\theta}_i^\star-\dot{\theta}_i^{\dagger}\|\leq\EpsFJ\\
|h^{e,j}(\theta_i^\star)|\leq\EpsFJ\quad \forall j\\
\min\!\big(h^{i,j}(\theta_i^\star),\ \lambda_i^{i,j},\ \lambda_i^{i,j}(\EpsFJ-h^{i,j}(\theta_i^\star))\big)\geq0\quad \forall j.
\end{cases}
\end{ResizedEq}
The multiplier $\lambda_i^i$ is the current constraint multiplier used for current feasibility and complementarity, while $\lambda_{\gamma(i)}^i$ is the semi-implicit dissipation multiplier induced by the previous-step penalty forces.
\end{definition}
\begin{remark}
We always use the FJ optimality condition, which is strictly weaker than the KKT condition. To ensure the KKT condition holds, MFCQ is needed, which is not assumed throughout our analysis. (See e.g., Lemma 1 of~\cite{hinder2024worst}.) The FJ condition is positively homogeneous in its multiplier tuple $(\mu_0,\lambda_i^e,\lambda_i^i)$, so it need not be written in any particular normalization: any nonzero rescaling of a valid tuple certifies the same FJ point. We state~\prettyref{eq:eps-fj-system} in the unit-objective form $\mu_0=1$, which is available here because the surrogate~\prettyref{pb:P-CRDS} is unconstrained---its penalty forces are finite at each $\EpsFJ$---so the objective multiplier can be normalized without loss. It remains an FJ rather than a KKT condition: no constraint qualification is assumed, and the constraint forces are unbounded, possibly diverging as $\EpsFJ\to0$. Moreover, because the semi-implicit damping $\Damp$ is nonsmooth in the velocity, \prettyref{eq:eps-fj-system} is not the classical smooth Fritz--John system: its force-balance line is stated in the standard nonsmooth form of approximate optimality relative to a nearby point on the (sub)differential graph, i.e.\ the damping subgradient is evaluated at a witness velocity $\dot\theta_i^\dagger$ close to $\dot\theta_i^\star$.
\end{remark}
Given the above definition, we have the following classical result when the unconstrained CRDS is solved exactly. Here $\LHamilton<\infty$ is a design clamp level chosen by the analysis rather than given problem data: it caps the admissible penalty value, and the next section shows that such a bound can indeed be arranged.
\begin{lemma}
\label{lem:eps-FJ}
Let~\prettyref{pb:CRDS-ABD} be a CRDS and $\theta_i^\star$ be an exact local minimizer of~\prettyref{pb:P-CRDS} with finite objective value. If $V(\theta_{\gamma(i)}^\star,\EpsFJ)\leq\LHamilton<\infty$, then $\theta_i^\star$ satisfies the $\EpsFJ$-FJ condition defined in~\prettyref{def:eps-FJ}.
\end{lemma}
Unfortunately,~\prettyref{lem:eps-FJ} still requires that we find the exact local minimizer of~\prettyref{pb:P-CRDS}, which is computationally infeasible. Further, the SQP convergence theory requires the constraint data to be smooth, which is violated by our non-smooth dissipation function $\Damp^V$. To address these two issues, we introduce the following smooth approximation to $\Damp^V$.
\begin{definition}
\label{def:S-D}
$\Damp_{\EpsFri,\LHamilton}^V$ is a smoothed version of $\Damp^V$ if it is an everywhere-defined finite real-valued extension $\Damp_{\EpsFri,\LHamilton}^V(\theta_{\gamma(i)},\dot{\theta}_i):\mathbb R^n\times\mathbb R^n\to\mathbb R$ such that: (i) the velocity map $\Damp_{\EpsFri,\LHamilton}^V(\theta_{\gamma(i)},\cdot)$ is convex and nonnegative; (ii) $0\in\partial_{\dot{\theta}_i}\Damp_{\EpsFri,\LHamilton}^V(\theta_{\gamma(i)},0)$; and (iii) its velocity gradient $\FPPR{\Damp_{\EpsFri,\LHamilton}^V}{\dot{\theta}_i}$ is $C^1$ in both parameters. (iv) for any $\theta_{\gamma(i)}$ with $V(\theta_{\gamma(i)},\EpsFJ)\leq\LHamilton<\infty$ and $\dot{\theta}_i^\star$, there exists a velocity $\dot{\theta}_i^{\dagger}$ such that:
\begin{enumerate}[label=(\alph*),ref=(\alph*),leftmargin=*]
\item\label{it:sd-grad-close} the smoothed and nonsmooth $\theta_i$-gradients are close,
\begin{ResizedEq}
\dist\!\Big(&\nabla_{\theta_i}\left(\delta_i\Damp_{\EpsFri,\LHamilton}^V\right)(\theta_{\gamma(i)},\dot{\theta}_i^\star),\\
&\partial_{\theta_i}\left(\delta_i\Damp^V\right)(\theta_{\gamma(i)},\dot{\theta}_i^{\dagger})\Big)\leq\delta_i\EpsFri;
\end{ResizedEq}
\item\label{it:sd-vel-close}
$\|\dot{\theta}_i^\star-\dot{\theta}_i^{\dagger}\|\leq\EpsFri$.
\end{enumerate}
\end{definition}
Next, we solve the following smoothed version of unconstrained CRDS to a gradient residual norm of $\EpsCRDS$.
\begin{ResizedEq}
\label{pb:PS-CRDS}
\argmin{\theta_i}\Phi_i(\theta_i)+V(\theta_i,\EpsFJ)+\delta_i\Damp_{\EpsFri,\LHamilton}^V(\theta_{\gamma(i)},\dot{\theta}_i).
\end{ResizedEq}
We first confirm that the smoothed problem remains inside the CRDS class.
\begin{corollary}
\label{cor:S-CRDS}
If~\prettyref{pb:CRDS-ABD} is a CRDS, $\operatorname{dom}U\cap\operatorname{dom}V\neq\emptyset$, and $\Damp_{\EpsFri,\LHamilton}^V$ is a smoothed version of $\Damp^V$ in the sense of~\prettyref{def:S-D}, then~\prettyref{pb:PS-CRDS} defines a smooth, unconstrained CRDS.
\end{corollary}
Here a smoothing in the sense of~\prettyref{def:S-D} is everywhere-defined and finite on $\mathbb R^n\times\mathbb R^n$. We then have the following result, which revises~\prettyref{lem:eps-FJ} to allow sufficient smoothness and computationally feasible error tolerance.
\begin{lemma}
\label{lem:eps-S-FJ}
Let~\prettyref{pb:CRDS-ABD} be a CRDS, let $\Damp_{\EpsFri,\LHamilton}^V$ be a smoothed version of $\Damp^V$, and let $\theta_i^\star$ be an approximate local minimizer of~\prettyref{pb:PS-CRDS} with finite objective value and gradient residual at most $\EpsCRDS$. If $V(\theta_{\gamma(i)}^\star,\EpsFJ)\leq\LHamilton<\infty$, then $\theta_i^\star$ satisfies the $\EpsFJAll$-FJ condition defined in~\prettyref{def:eps-FJ}, where
\begin{equation}
\EpsFJAll\triangleq\max\{\EpsFJ,\EpsFri,\EpsCRDS+\delta_i\EpsFri\}.
\end{equation}
\end{lemma}
The proofs of~\prettyref{lem:eps-FJ},~\ref{lem:eps-S-FJ}, and~\prettyref{cor:S-CRDS} are collected in~\prettyref{appen:penalty-proofs}.

\subsection{Articulated-Body Specialization}
The general penalization above applies directly to the CRDS maps $h^e$ and $h^i$. In the articulated-body instantiation, the optimization variable is extended from the physical configuration $\theta_i$ to $\vartheta_i$ to include separating planes and related auxiliary variables. We reuse $h^i(\vartheta_i)$ for the single stacked inequality map containing the physical, separation, and unit-normal components. Each physical component is extended trivially by $h^{i,j}(\vartheta_i)\triangleq h^{i,j}(\theta_i)$ and therefore has zero derivative in every auxiliary coordinate.
With this interpretation, the generic constrained equivalence specializes to
\begin{ResizedEq}
\label{pb:CRDS-ABD-extended}
\argmin{\vartheta_i}&
\Phi_i(\theta_i)+\delta_i\Damp(\vartheta_{\gamma(i)},\dot{\vartheta}_i,\lambda_{\gamma(i)}^i)\\
\ST\quad&h^e(\theta_i)=0,
h^i(\vartheta_i)\geq0.
\end{ResizedEq}
Here $\dot\vartheta_i=(\vartheta_i-\vartheta_{\gamma(i)})/\delta_i$ is the complete finite-difference velocity, including $\dot p_i$, while the explicit damping depends only on its damping-active subvector $\TWO{\dot\theta_i}{\varpi_i}$. This extended formulation may fail to be full-rank because $\vartheta_i$ includes massless separating-plane variables. We therefore apply two further transforms to obtain the reduced full-rank formulation analyzed in the remainder of this section.

\prettyref{pb:P-CRDS}, interpreted with the articulated-body extended variable $\vartheta_i$, is not full-rank due to the additional separating planes in the configuration vector. Fortunately, we can eliminate these separating planes using the following established result in prior works~\cite{11266917}.
\begin{corollary}[(Lemma VII.3 and Lemma VIII.1 of~\cite{11266917})]
\label{cor:convexity}
$V(\vartheta_i,\EpsFJ)$ is convex in every $p^{jh,j'h'}$ and, for every fixed nonnegative predecessor load $\lambda_{\gamma(i)}^i$, $\Damp(\vartheta_{\gamma(i)},\dot{\vartheta}_i,\lambda_{\gamma(i)}^i)$ is convex in every $\TWO{\varsigma}{\omega}_i^{jh,j'h'}$.
\end{corollary}
The added plane-velocity damper $\EpsDamp(\|\dot{\varsigma}^{jh,j'h'}\|^2+\|\dot{\omega}^{jh,j'h'}\|^2)$ is a convex quadratic in $\TWO{\varsigma}{\omega}^{jh,j'h'}$, so, together with the cited convexity of the friction part, $\Damp$ remains convex in every $\TWO{\varsigma}{\omega}^{jh,j'h'}$, and the reduction in~\prettyref{eq:reduced-penalty} stays a convex program. \prettyref{cor:convexity} allows us to eliminate all the separating planes by exact mathematical partial minimization over the separating planes and their velocities, which defines the reduced penalty value functions
\begin{ResizedEq}
\label{eq:reduced-penalty}
V_R(\theta_i,\EpsFJ)\triangleq&\min_{\{p^{jh,j'h'}\}}V(\vartheta_i,\EpsFJ)\\
\Damp_R^V(\theta_{\gamma(i)},\dot{\theta}_i)\triangleq&\min_{\varpi}\Damp^V(\vartheta_{\gamma(i)},(\dot\theta_i,\varpi)).
\end{ResizedEq}
Throughout the analysis, the minima in~\prettyref{eq:reduced-penalty} are exact mathematical definitions of value functions: auxiliary minimizers are attained and the associated stationarity and envelope identities hold without an inner-solve tolerance. This convention concerns the theoretical reduction, not a claim that a finite-precision implementation computes the minima exactly. This minimization eliminates the auxiliary damping-active velocity block $\varpi$; it does not remove the $\dot p$ block from the complete finite-difference velocity, because the explicit damping is already constant in $\dot p$. Further, although $V(\vartheta_i,\EpsFJ)$ is not strongly convex and is defined by piecewise functions, the following result shows that all discontinuous configurations are removable and that the reduced penalty value function $V_R(\theta_i,\EpsFJ)$ is globally twice differentiable; we will invoke~\prettyref{lem:twice-differentiability} in the global-convergence analysis to enable second-order arguments.
\begin{lemma}[(Lemma VII.5 of~\cite{11266917})]
\label{lem:twice-differentiability}
$V_R(\theta_i,\EpsFJ)$ is proper, continuous, and $C^2$ on its finite domain.
\end{lemma}
Finally, we smooth the reduced dissipation $\Damp_R^V$ while keeping the reduced penalty value function $V_R$ as the conservative penalty potential. Putting these pieces together yields the following reduced smoothed unconstrained problem.
\begin{ResizedEq}
\label{pb:PSR-CRDS}
\argmin{\theta_i}\Phi_i(\theta_i)+V_R(\theta_i,\EpsFJ)+\delta_i\Damp_{R,\EpsFri,\LHamilton}^V(\theta_{\gamma(i)},\dot{\theta}_i).
\end{ResizedEq}
The reduced problem defines a smooth full-rank CRDS because the auxiliary separating-plane variables have been eliminated while the physical maximal-coordinate mass matrix is retained, and the retained physical block is positive definite by~\prettyref{lem:mass} under the standing~\prettyref{ass:abd-mass}.
\begin{corollary}
\label{cor:SF-CRDS}
Under the standing~\prettyref{ass:abd-mass}, if $\Damp_{R,\EpsFri,\LHamilton}^V$ is a smoothed version of $\Damp_R^V$ in the sense of~\prettyref{def:S-D} (with $V_R$ in the role of $V$), then~\prettyref{pb:PSR-CRDS} is a smooth full-rank CRDS.
\end{corollary}
The object $\Damp_{R,\EpsFri,\LHamilton}^V$ appearing in~\prettyref{pb:PSR-CRDS} is not vacuous: such a smoothed reduced dissipation can in fact be constructed.
\begin{lemma}[(informal)]
\label{lem:smooth-DRp-Regularity}
For each fixed $\EpsFri>0$, $0<\LHamilton<\infty$, and prescribed a priori timestep floor $\delta_{\min}>0$, there exists a reduced dissipation $\Damp_{R,\EpsFri,\LHamilton}^V$ that is a smoothed version of $\Damp_R^V$ in the sense of~\prettyref{def:S-D}, with the reduced penalty potential $V_R$ playing the role of $V$.
\end{lemma}
Here $\delta_{\min}>0$ is a design constant fixed before any solve, required to satisfy $\delta_i\geq\delta_{\min}$ for every optimized timestep; a concrete admissible choice is given in~\prettyref{eq:cito-delta-min}. Together with~\prettyref{cor:SF-CRDS}, this makes~\prettyref{pb:PSR-CRDS} a bona fide smooth full-rank CRDS and discharges the smoothness hypothesis of~\prettyref{lem:eps-SR-FJ} below, under the exact-reduction convention following~\prettyref{eq:reduced-penalty}. The construction and proof are given in~\prettyref{appen:penalty-proofs}. For articulated-body contact, the attainment, conservative envelope, and plane-velocity reconstruction conditions used by the lift are established respectively by~\prettyref{lem:active-plane},~\prettyref{eq:VR-envelope-proof}, and~\prettyref{lem:varpi-reconstruction}. The following lifting result then connects approximate minimizers of the reduced smoothed problem back to perturbed FJ points of the extended constrained formulation.
\begin{lemma}
\label{lem:eps-SR-FJ}
Let $\theta_i^\star$ be an approximate local minimizer of~\prettyref{pb:PSR-CRDS} with finite objective value and gradient residual at most $\EpsCRDS$, and let $\Damp_{R,\EpsFri,\LHamilton}^V$ be a smoothed version of $\Damp_R^V$ in the sense of~\prettyref{def:S-D} (with $V_R$ in the role of $V$). Assume that $V_R(\theta_{\gamma(i)}^\star,\EpsFJ)\leq\LHamilton<\infty$. If the auxiliary minimizers in~\prettyref{eq:reduced-penalty} are attained and the reduced value functions satisfy the corresponding envelope identities at $\theta_i^\star$, then there exist separating-plane parameters inducing an extended configuration $\vartheta_i^\star$, such that $\vartheta_i^\star$ satisfies the $\EpsFJAll$-FJ condition of~\prettyref{pb:CRDS-ABD-extended}.
\end{lemma}
Here \prettyref{def:eps-FJ} is instantiated by substituting the complete extended ABD configuration $\vartheta$ for the generic CRDS configuration $\theta$. Consequently, every inequality component is evaluated as $h^{i,j}(\vartheta)$, with physical components depending only on the embedded $\theta$ block, while separation and unit-normal components depend on the relevant full coordinates, and the lifted velocity follows~\prettyref{eq:extended-config}; this substitution does not identify the full configuration $\vartheta$ with its velocity. \prettyref{lem:eps-S-FJ} and \ref{lem:eps-SR-FJ} are instrumental for our subsequent analysis because we only need to consider the unconstrained version of CRDS and the smoothness of this CRDS allows us to adopt the conventional convergence analysis technique of the SQP algorithm.

\prettyref{pb:PSR-CRDS} is the reduced, smooth per-step articulated-body dynamics surrogate. It is not yet the full-trajectory optimization problem solved by the SQP method. The next section mollifies its potential and control work, yielding the mollified per-step CRDS, whose stationarity equation is then stacked over the horizon to form the physics constraints of the full-trajectory NLP. The proofs of~\prettyref{cor:SF-CRDS} and~\prettyref{lem:eps-SR-FJ} are collected in~\prettyref{appen:penalty-proofs}.

\section{\label{sec:Hamiltonian}Hamiltonian Upper Bound}
We have seen in~\prettyref{lem:eps-S-FJ} of~\prettyref{sec:penalty} that the smoothed unconstrained CRDS~\prettyref{pb:PS-CRDS} dictates a dynamic trajectory that approximates the original CRDS~\prettyref{pb:CRDS-ABD} arbitrarily well. However, this result relies on the assumption that the penalty value $V$ be upper bounded by the design clamp $\LHamilton$, i.e., $V(\theta_{\gamma(i)}^\star,\EpsFJ)\leq\LHamilton$. In this section, we show that this bound can indeed be achieved. Our intuition is that, for the general CRDS, its smoothed version~\prettyref{pb:PS-CRDS} defines a nearly conservative dynamic system with the total energy defined by the following Hamiltonian $\mathcal{H}_i$ and its nonnegative part $\mathcal{H}_i^+$
\begin{ResizedEq}
\tilde{U}(\theta_i,\EpsFJ)\triangleq&U(\theta_i)+V(\theta_i,\EpsFJ)\\
\mathcal{H}_i(\theta_{\gamma(i)},\theta_i,\EpsFJ)\triangleq&\frac{1}{2}\|\dot{\theta}_i\|_M^2+\tilde{U}(\theta_i,\EpsFJ)-g^T\theta_i\\
\mathcal{H}_i^+(\theta_{\gamma(i)},\theta_i,\EpsFJ)\triangleq&\frac{1}{2}\|\dot{\theta}_i\|_M^2+\tilde{U}(\theta_i,\EpsFJ)\geq0.
\end{ResizedEq}
It is well known that the continuous trajectory of a conservative dynamic system preserves the Hamiltonian $\mathcal{H}$. Suppose we can upper bound its discrete counterpart $\mathcal{H}_i$ and its nonnegative part $\mathcal{H}_i^+$ in particular; we can then inherently upper bound $V$ to invoke~\prettyref{lem:eps-S-FJ}. We will show that it is indeed possible to upper bound the discrete Hamiltonian $\mathcal{H}_i$ if the negative curvature of the potential energy function $\tilde{U}$ can be lower bounded. Unfortunately, such a curvature lower bound is not possible in general. As a counterexample, consider a dynamic system involving two points, i.e., $\theta=\TWO{x}{y}$. Since the two points should not collide, we have the constraint $h^i(\theta): \|x-y\|-\EpsContactPadding\geq0$. The corresponding penalty function is $V^i(\theta,\EpsFJ)=V_s(\|x-y\|-\EpsContactPadding,\EpsFJ)$. It is easy to see that $V^i$ has unbounded negative curvature when $\|x-y\|\to0^+$. To mitigate this dilemma, we use the mollification function $\mathcal M$ introduced in the coordinate-maximization construction in~\prettyref{eq:clamp-map}, as illustrated in~\prettyref{fig:clamping}. We carefully design a mollified version of $\tilde{U}$ denoted as $\bar{U}$ that limits the support of $\tilde{U}$ in both domain and range space. Specifically, we define
\begin{ResizedEq}
\bar{U}(\theta,\EpsFJ,\LHamilton,\LPosition)\triangleq&
\mathcal{M}(\tilde{U}(\mathcal{M}(\theta,\LPosition),\EpsFJ),\LHamilton).
\end{ResizedEq}
We treat the control work symmetrically: since the control force $c=-\FPP{W}{\theta}$ enters the surrogate through the configuration argument, we mollify the control work in the same way, using only the inner (domain) clamp $\mathcal M(\cdot,\LPosition)$ since $W$ needs no range saturation,
\begin{ResizedEq}
\bar{W}(\theta,u,\LPosition)\triangleq&W(\mathcal{M}(\theta,\LPosition),u).
\end{ResizedEq}
\begin{figure}[t]
\centering
\begin{tikzpicture}[font=\small, x=0.72cm, y=0.95cm]

\def\Utilde#1{0.7 + 0.13*sin(120*(#1)) + 0.07*cos(220*(#1)+30)
  + 1.05*exp(-((#1)-1.6)*((#1)-1.6)/0.5)
  + 0.45*exp(-((#1)-3.1)*((#1)-3.1)/0.35)
  + 8*exp(-((#1)+2.2)*((#1)+2.2)/0.05)}

\def\clp#1{max(-3.6,min(3.6,#1))}

\def\Uc#1{\Utilde{\clp{#1}}}
\def\Ubar#1{2.6 - ((2.6-(\Uc{#1}))
  + sqrt((2.6-(\Uc{#1}))*(2.6-(\Uc{#1})) + 0.05))/2}

\def\Lham{2.6}   
\def\Lpos{3.6}   
\def\yTop{3.6}   

\useasboundingbox (-5.2,-1.5) rectangle (5.7,4.2);

\fill[green!55!black, opacity=0.13] (-5,0) rectangle (-\Lpos,\yTop);
\fill[green!55!black, opacity=0.13] ( \Lpos,0) rectangle ( 5.4,\yTop);

\begin{scope}
  \clip (-5,0) rectangle (5.4,\yTop);
  \draw[black, line width=2.2pt, smooth, samples=500, domain=-5:5]
    plot (\x, {\Utilde{\x}});

  \begin{scope}
    \clip (-2.5,\Lham) rectangle (-1.9,\yTop);
    \foreach \h in {-1.0,-0.7,...,4.2} {
      \draw[black, line width=0.4pt] (-2.7,\h) -- (-1.7,\h+1.0);
    }
  \end{scope}

  \draw[red!85!black, line width=1.6pt, smooth, samples=400, domain=-5:5]
    plot (\x, {\Ubar{\x}});
\end{scope}
\node[fill=white, inner sep=1pt] at (-2.2,3.05) {$\infty$};

\draw[green!55!black, line width=1.2pt, dashed]
  (-\Lpos,\Lham) -- (3.3,\Lham);
\node[green!40!black, anchor=west, fill=white, inner sep=1pt] at (2.85,\Lham)
  {$\LHamilton$};

\draw[green!55!black, ->, >=stealth] (-\Lpos,-0.75) -- (-\Lpos,-0.12);
\node[green!40!black, anchor=north] at (-\Lpos,-0.85) {$-\LPosition$};
\draw[green!55!black, ->, >=stealth] ( \Lpos,-0.75) -- ( \Lpos,-0.12);
\node[green!40!black, anchor=north] at ( \Lpos,-0.85) {$\LPosition$};

\draw[->, >=stealth, line width=1pt] (-5,0) -- (5.5,0)
  node[anchor=west] {$\theta$};
\draw[->, >=stealth, line width=1pt] (0,-0.2) -- (0,4.0)
  node[anchor=south] {$\tilde{U}(\theta,\EpsFJ)$};

\draw[black, line width=2.2pt] (0.6,3.8) -- (1.2,3.8);
\node[anchor=west] at (1.2,3.8) {original $\tilde{U}$};
\draw[red!85!black, line width=1.6pt] (0.6,3.45) -- (1.2,3.45);
\node[anchor=west] at (1.2,3.45) {mollified $\bar{U}$};

\end{tikzpicture}
\caption{\label{fig:clamping}Illustration of the mollification (clamping) used to bound the negative curvature of the potential energy. The original potential $\tilde{U}(\theta,\EpsFJ)$ (black) follows a wiggly profile but develops an unbounded spike (the hatched column marked $\infty$) wherever two bodies approach contact ($\|x-y\|\to0^+$), so its negative curvature cannot be lower bounded. The mollified potential $\bar{U}(\theta,\EpsFJ,\LHamilton,\LPosition)$ (red) coincides with $\tilde{U}$ inside the admissible window but is clamped in two complementary ways. The inner mollifier $\mathcal{M}(\theta,\LPosition)$ caps the domain: outside $|\theta|\leq\LPosition$ (shaded green regions beyond $\pm\LPosition$, the position clamp) the configuration is frozen, so $\bar{U}$ becomes constant. The outer mollifier $\mathcal{M}(\cdot,\LHamilton)$ caps the range: wherever $\tilde{U}$ would exceed $\LHamilton$ (horizontal green line at $\LHamilton$, the Hamiltonian clamp) the value is smoothly saturated, flattening the spike into a finite plateau. The transition is twice-differentiable thanks to the shared degree-seven smooth blend $p$, so $\bar{U}$ has bounded curvature (\prettyref{lem:mollification}) while remaining identical to $\tilde{U}$ on the region of interest.}
\end{figure}
The following mollification lemma formalizes the regularity and curvature properties of this construction; its regularity hypothesis holds since $\tilde{U}=U+V$ inherits these properties.
\begin{lemma}
\label{lem:mollification}
Assume $0<\LHamilton,\LPosition<\infty$ and that $\tilde{U}(\cdot,\EpsFJ)$ satisfies the potential regularity of~\prettyref{def:CRDS}~\ref{it:crds-potential}. Then:
\begin{enumerate}[label=(\alph*),ref=(\alph*),leftmargin=*]
\item\label{it:moll-identity} $\bar{U}(\theta,\EpsFJ,\LHamilton,\LPosition)=\tilde{U}(\theta,\EpsFJ)$ when $\|\theta\|\leq\LPosition$ and $\tilde{U}(\theta,\EpsFJ)\leq\LHamilton$;
\item\label{it:moll-c2} $\bar{U}$ is $C^2$ on $\mathbb R^n$;
\item\label{it:moll-curv} there exists $0<\LCurvature<\infty$ such that
\begin{ResizedEq}
-\LCurvature I\preceq&\FPPT{\bar U(\theta,\EpsFJ,\LHamilton,\LPosition)}{\theta}
\preceq\LCurvature I\qquad\forall\theta;
\end{ResizedEq}
\item\label{it:moll-bound} $0\leq\bar U(\theta,\EpsFJ,\LHamilton,\LPosition)\leq2\LHamilton$ and there exists $\LPotentialForce<\infty$ such that $\|\FPPR{\bar U(\theta,\EpsFJ,\LHamilton,\LPosition)}{\theta}\|\leq\LPotentialForce$ for all $\theta$.
\end{enumerate}
\end{lemma}
\begin{lemma}[Control-work mollification]
\label{lem:controlwork-mollification}
For any CRDS (\prettyref{def:CRDS}) and any position clamp $0<\LPosition<\infty$, the mollified control work $\bar W(\theta,u,\LPosition)=W(\mathcal M(\theta,\LPosition),u)$ satisfies, for all $\theta$ and all $u\in\mathcal{D}$:
\begin{enumerate}[label=(\alph*),ref=(\alph*),leftmargin=*]
\item\label{it:cw-identity} $\bar W(\theta,u,\LPosition)=W(\theta,u)$ when $\|\theta\|\leq\LPosition$;
\item\label{it:cw-c2} $\bar W$ is $C^2$ in $\theta$;
\item\label{it:cw-force} there exists a clamp-independent $\LControlWorkForce<\infty$ such that $\|\FPPR{\bar W(\theta,u,\LPosition)}{\theta}\|\leq\LControlWorkForce$;
\item\label{it:cw-curv} there exists $\LControlWorkCurvature<\infty$ (clamp-dependent but independent of $\theta$ and $u$) such that
\begin{ResizedEq}
-\LControlWorkCurvature I\preceq&\FPPT{\bar W(\theta,u,\LPosition)}{\theta}
\preceq\LControlWorkCurvature I;
\end{ResizedEq}
\item\label{it:cw-mixed} there exists a finite clamp-dependent constant, independent of $\theta$ and $u$, that uniformly bounds $\|\FPPTT{\bar W(\theta,u,\LPosition)}{\theta}{u}\|$.
\end{enumerate}
\end{lemma}
The proofs of~\prettyref{lem:mollification} and~\ref{lem:controlwork-mollification} are collected in~\prettyref{appen:Hamiltonian-upper-bound-proofs}.
We now define the mollified, unconstrained, smoothed per-step CRDS whose first-order stationarity equation will serve as the discrete physics residual in the full-trajectory NLP:
\begin{ResizedEq}
\label{pb:PS-CRDS-M}
\argmin{\theta_i}&\PerStepObj_i(\theta_{\gamma^2(i)},\theta_{\gamma(i)},\theta_i,u_i)\triangleq\\
&\Big[\frac{\delta_i^2}{2}\|\ddot{\theta}_i\|_M^2+\bar{U}(\theta_i,\EpsFJ,\LHamilton,\LPosition)-\\
&g^T\theta_i+\bar{W}(\theta_i,u_i,\LPosition)+\delta_i\Damp_{\EpsFri,\LHamilton}^V(\theta_{\gamma(i)},\dot{\theta}_i)\Big].
\end{ResizedEq}

\prettyref{pb:PS-CRDS-M} is therefore a per-step dynamics surrogate. The user trajectory objective is not minimized independently at each timestep through~\prettyref{pb:PS-CRDS-M}; rather, the full-trajectory optimizer will enforce $\nabla_{\theta_i}\PerStepObj_i=0$ as a physics constraint while minimizing the user objective $O(\{\theta_i\}_{i\in\mathcal I},\{u_i\}_{i\in\mathcal I})$.

Because~\prettyref{pb:PS-CRDS-M} mollifies the general smoothed CRDS~\prettyref{pb:PS-CRDS}, replacing its potential $\tilde U$ by the nonnegative $C^2$ potential $\bar U$ from~\prettyref{lem:mollification} and its control work $W$ by the $C^2$ mollified control work $\bar W$ from~\prettyref{lem:controlwork-mollification}, it is a smooth CRDS by~\prettyref{cor:S-CRDS} and~\prettyref{lem:mollification}, and it is full-rank in the sense of~\prettyref{def:full-rank} whenever the base CRDS is full-rank. With the bounded negative curvature, we can upper bound the discrete Hamiltonian of the mollified CRDS, as denoted by
\begin{ResizedEq}
&\bar{\mathcal{H}}_i(\theta_{\gamma(i)},\theta_i,\EpsFJ,\LHamilton,\LPosition)\\
\triangleq&\frac{1}{2}\|\dot{\theta}_i\|_M^2+\bar{U}(\theta_i,\EpsFJ,\LHamilton,\LPosition)-g^T\theta_i,\\
&\bar{\mathcal{H}}_i^+(\theta_{\gamma(i)},\theta_i,\EpsFJ,\LHamilton,\LPosition)\\
\triangleq&\frac{1}{2}\|\dot{\theta}_i\|_M^2+\bar{U}(\theta_i,\EpsFJ,\LHamilton,\LPosition)\geq0.
\end{ResizedEq}
We then show that due to such upper bound, the behavior of $\bar{\mathcal{H}}_i$ is locally identical to that of $\mathcal{H}_i$ due to~\prettyref{lem:mollification}. As a result, the $\EpsFJAll$-FJ solution to~\prettyref{pb:PS-CRDS-M} is also an $\EpsFJAll$-FJ solution to~\prettyref{pb:PS-CRDS}.

We recall from~\prettyref{def:full-rank} that the mass matrix is positive definite, $M\succ0$, so its smallest eigenvalue $\lambda_{\min}(M)>0$ is well defined and yields the Euclidean-$M$ norm comparison $\|v\|^2\leq\|v\|_M^2/\lambda_{\min}(M)$ for all $v$.

Before stating the bound, we recall that the mollified control force $\bar c\triangleq-\FPP{\bar W}{\theta}$ is uniformly bounded: by~\prettyref{lem:controlwork-mollification}~\ref{it:cw-force} we have $\|\bar c(\theta_i,u_i)\|=\|\FPPR{\bar W(\theta_i,u_i,\LPosition)}{\theta_i}\|\leq\LControlWorkForce$ for every configuration $\theta_i$ and every admissible control $u_i\in\mathcal{D}$. Our energy bound is anchored at the fixed strict-interior configuration $\theta_0$ of~\prettyref{ass:Slater}, introduced with the problem statement. That assumption requires the user-provided time-zero condition to strictly satisfy all constraints. Together with the fixed startup data $(\theta_{\gamma(0)},\theta_0,\delta_0)$, it guarantees the finite base velocity $\dot\theta_0=(\theta_0-\theta_{\gamma(0)})/\delta_0$ and initial Hamiltonian $\bar{\mathcal{H}}_0^+$; the one-step estimate is seeded on $(0,\min\mathcal I)$ to control the first optimized configuration $\theta_{\min\mathcal I}$. A one-step estimate, stated and proved in~\prettyref{appen:Hamiltonian-upper-bound-proofs}, bounds the increase of $\bar{\mathcal{H}}_i^+$ over a single timestep in terms of the two-sided curvature modulus $\LCurvature$, the mollified control-force bound $\LControlWorkForce$, and the gradient residual $\EpsCRDS$. This one-step estimate requires the smoothed dissipation $\Damp_{\EpsFri,\LHamilton}^V(\theta_{\gamma(i)},\cdot)$ to be convex and minimized at zero relative velocity, i.e., $0\in\partial_{\dot\theta_i}\Damp_{\EpsFri,\LHamilton}^V(\theta_{\gamma(i)},0)$. These properties hold for any CRDS by~\prettyref{def:CRDS}~\ref{it:crds-damping}; since~\prettyref{pb:PS-CRDS-M} is a smooth full-rank CRDS, its dissipation satisfies them. Accumulating this estimate over all timesteps yields a global bound that depends only on the elapsed time.
\begin{lemma}[Hamiltonian upper bound]
\label{lem:H-bound}
Take~\prettyref{ass:Slater}, and assume:
\begin{enumerate}[label=(\alph*),leftmargin=*]
\item the clamps are chosen no smaller than the fixed time-zero data, $\LPosition\geq\|\theta_0\|$ and $\LHamilton\geq\tilde U(\theta_0,\EpsFJ)$;
\item the timestep sequence satisfies
\begin{equation}
\label{eq:H-bound-timestep}
\max_{i\in\mathcal I}\delta_i<\min\left(\frac{1}{\LCurvature},\frac{\lambda_{\min}(M)}{6}\right);
\end{equation}
\item \prettyref{pb:PS-CRDS-M} is a smooth full-rank CRDS, each of whose solutions $\theta_i$ has gradient residual of norm at most $\EpsCRDS$.
\end{enumerate}
Then there exists a clamp-independent function $\LHamBound(t_i)$, finite for every finite $t_i$, such that for all $i\in\mathcal I$,
\begin{equation}
\label{eq:H-bound}
\bar{\mathcal{H}}_i^+\leq\LHamBound(t_i).
\end{equation}
\end{lemma}
The explicit closed form of $\LHamBound(t_i)$ is derived in the proof, collected in~\prettyref{appen:Hamiltonian-upper-bound-proofs}. \prettyref{lem:H-bound} formalizes the observation that the mollified surrogate~\prettyref{pb:PS-CRDS-M} is a nearly conservative system whose Hamiltonian stays constant up to a discretization error, which we upper bound. The bound function $\LHamBound(t_i)$ is clamp-independent: the two-sided curvature modulus $\LCurvature$ drops out of $\LHamBound$ because the timestep condition~\prettyref{eq:H-bound-timestep} makes $\LCurvature\delta_i^2\leq\delta_i$, while the clamps $\LHamilton$ and $\LPosition$ drop out because~\prettyref{ass:Slater} and~\prettyref{lem:mollification}~\ref{it:moll-identity} render the initial value $\bar{\mathcal{H}}_0^+=\frac{1}{2}\|\dot{\theta}_0\|_M^2+\tilde U(\theta_0,\EpsFJ)$, with $\dot\theta_0=(\theta_0-\theta_{\gamma(0)})/\delta_0$, clamp-independent. We emphasize that only the upper-bound function is clamp-independent: the validity condition~\prettyref{eq:H-bound-timestep} still depends on the two-sided curvature modulus $\LCurvature$, which in turn depends on the clamps $\LHamilton,\LPosition$ through~\prettyref{lem:mollification}.

An immediate consequence bounds the trajectory and fixes the meaning of the domain clamp $\LPosition$. Since $\bar U\geq0$ gives $\frac{1}{2}\|\dot{\theta}_i\|_M^2\leq\bar{\mathcal{H}}_i^+\leq\LHamBound(t_i)$, each velocity obeys $\|\dot{\theta}_i\|\leq\sqrt{2\LHamBound(t_i)/\lambda_{\min}(M)}$; telescoping $\theta_i=\theta_0+\sum_{k\in\mathcal I:\,k\leq i}\delta_k\dot{\theta}_k$ along the ordered predecessor chain then bounds the position on any finite horizon, which fixes the domain clamp $\LPosition$ (made precise in the constant-ordering remark below). Finally, because $V\leq\tilde U\leq\bar{\mathcal{H}}_i^+$ (immediate since $\tilde U=U+V$ with $U\geq0$) on the region where the clamps are inactive, choosing $\LHamilton>\LHamBound(H)$ over the horizon $H$ discharges the standing assumption $V(\theta_{\gamma(i)}^\star,\EpsFJ)\leq\LHamilton$ required by~\prettyref{lem:eps-S-FJ}, which is the payoff pursued in this section.
\begin{remark}
To avoid circularity, the constants are fixed in this order:
\begin{enumerate}[label=(\alph*),ref=(\alph*),leftmargin=*]
\item\label{it:order-clamps} choose gaps $\Delta_{\mathcal H}>0$ and $\Delta_\theta>0$ and set
\begin{ResizedEq}
\LHamilton\triangleq&\LHamBound(H)+\Delta_{\mathcal H}\\
\LPosition\triangleq&\|\theta_0\|+H\sqrt{2\LHamBound(H)/\lambda_{\min}(M)}+\Delta_\theta;
\end{ResizedEq}
\item let $\LCurvature$ and $\LControlWorkCurvature$ be the two-sided curvature moduli of $\bar U$ and $\bar W$ from~\prettyref{lem:mollification} and~\ref{lem:controlwork-mollification} for these clamps;
\item choose the timestep sequence satisfying $\max_{i\in\mathcal I}\delta_i<\min(1/\LCurvature,\lambda_{\min}(M)/6)$;
\item\label{it:order-floor} fix the a priori timestep floor $\delta_{\min}>0$ no larger than $\min\{\min_{i\in G_0}\delta_i^0,\tfrac12\min(\delta^\star,\delta_H)\}$---the admissible choice of~\prettyref{eq:cito-delta-min}, with $\delta^\star,\delta_H$ the singularity and Hamiltonian timestep thresholds of~\prettyref{sec:global}---and require every optimized timestep to satisfy $\delta_i\geq\delta_{\min}$;
\item\label{it:order-smoothing} choose the damping-smoothing scale as a function of $\delta_{\min}$ and the fixed clamps.
\end{enumerate}
\prettyref{it:order-clamps} is well posed because $\LHamBound(t)$---and hence both $\LHamBound(H)$ and the position majorant $\|\theta_0\|+H\sqrt{2\LHamBound(H)/\lambda_{\min}(M)}$---is independent of $\LHamilton$, $\LPosition$, $\LCurvature$, and $\LControlWorkCurvature$; the mollified control-force bound $\LControlWorkForce$ that enters $\LHamBound$ is itself clamp-independent by~\prettyref{lem:controlwork-mollification}~\ref{it:cw-force}, so this clamp choice depends on no later constant. The positive gaps ensure that every returned finite node has an open neighborhood on which the clamps are identity. Items~(d)--(e) close the ordering without circularity: the thresholds $\delta^\star,\delta_H$ entering the floor~\prettyref{eq:cito-delta-min} are built only from $\lambda_{\min}(M)$ and the clamp-dependent curvature moduli of steps~(a)--(b), never from the damping-smoothing scale, so the floor is fixed before the smoothing and the smoothing then depends only on it. Item~(c) bounds $\max_{i\in\mathcal I}\delta_i$ from above and item~(d) bounds $\min_{i\in\mathcal I}\delta_i$ from below; the two are independent and must both hold.
\end{remark}

The bound $\LHamBound(t_i)$ of~\prettyref{lem:H-bound} was deliberately made clamp-independent so that the clamps $\LHamilton,\LPosition$ could subsequently be chosen no smaller than $\LHamBound(H)$ without circularity. That design choice, however, forced two restrictive hypotheses: a small per-step gradient residual $\leq\EpsCRDS$ and the timestep condition~\prettyref{eq:H-bound-timestep} tying $\max_{i\in\mathcal I}\delta_i$ to $\LCurvature$. For the bounded-iterate claim of the abstract---that the optimized trajectory stays norm-bounded throughout the optimization---we need a complementary result that survives the regime encountered inside the line-search SQP of~\prettyref{sec:global}, where each instance of~\prettyref{pb:PS-CRDS-M} is solved inexactly (large residual) and the timestep is not tied to $\LCurvature$. The key observation is that the mollification of~\prettyref{lem:mollification}~\ref{it:moll-bound} makes the potential force globally bounded by $\LPotentialForce$; together with the mollified control force, uniformly bounded by $\LControlWorkForce$ via~\prettyref{lem:controlwork-mollification}~\ref{it:cw-force}, the total force is uniformly bounded and the discrete kinetic energy can grow at most linearly in elapsed time, for any active timestep sequence bounded below by a positive $\delta_{\min}$, with dissipation only helping. This yields velocity and position bounds that are polynomial in $t_i$ and require no timestep smallness condition and no strict feasibility (\prettyref{ass:Slater}), at the price of depending on the clamps through $\LPotentialForce$---which is legitimate here, since we are bounding the iterates rather than choosing the clamps.
\begin{lemma}[Unconditional trajectory norm bound]
\label{lem:trajectory-bound}
Assume:
\begin{enumerate}[label=(\alph*),leftmargin=*]
\item \prettyref{pb:PS-CRDS-M} is a smooth full-rank CRDS (\prettyref{def:CRDS},~\ref{def:full-rank}) with $\bar U$ satisfying~\prettyref{lem:mollification}~\ref{it:moll-bound};
\item the controls satisfy $u_i\in\mathcal D$;
\item the active timesteps satisfy $\delta_i\geq\delta_{\min}>0$;
\item each solution $\theta_i$ of~\prettyref{pb:PS-CRDS-M} has gradient residual of norm at most $\LResidual<\infty$.
\end{enumerate}
Then there exist functions $\LVelBound(t_i)$ and $\LPosBound(t_i)$, finite on every finite horizon, such that for any such timestep sequence $\{\delta_i\}$ and all $i$, the velocity and position obey
\begin{equation}
\label{eq:trajectory-velocity-bound}
\|\dot\theta_i\|_M\leq\LVelBound(t_i),\qquad\|\theta_i\|\leq\LPosBound(t_i).
\end{equation}
\end{lemma}
The explicit closed forms of $\LVelBound$ and $\LPosBound$ are derived in the proof, collected in~\prettyref{appen:Hamiltonian-upper-bound-proofs}. Unlike $\LHamBound(t_i)$, the velocity and position bounds are polynomial (linear velocity, quadratic position) rather than exponential in $t_i$, and hold without an upper timestep restriction and with an arbitrarily large but fixed finite residual ceiling $\LResidual$; the trade-off is that they are clamp-dependent through $\LPotentialForce$.
\begin{remark}
\prettyref{lem:trajectory-bound} provides the trajectory-boundedness component of the paper's bounded-optimization-history result. During SQP, every iterate keeps the per-step residual of~\prettyref{pb:PS-CRDS-M} within a common finite ceiling and a timestep schedule satisfying $\delta_i\geq\delta_{\min}>0$; the lemma then keeps the whole solution data $\{\theta_i\}$ inside a fixed ball of radius $\LPosBound(H)$ over the horizon $[0,H]$, uniformly across iterations. Together with the controls $u_i\in\mathcal D$, which are bounded by compactness of $\mathcal D$, this shows the optimized trajectory remains norm-bounded throughout the optimization process. Combined with the bounded merit-penalty sequence established in~\prettyref{thm:sqp-convergence} of the next section, this gives the paper's bounded-optimization-history result. The two bounds of this section play complementary roles: the fixed strict-interior node $0$ seeds the clamp-independent physical-energy chain of~\prettyref{lem:H-bound} and the subsequent transfer to the original constrained CRDS, whereas~\prettyref{lem:trajectory-bound} needs no strict feasibility and is instead clamp-dependent. The later should not be read as the gadget to remove~\prettyref{ass:Slater}: it bounds iterates of the globally mollified smooth problem, whereas the physical feasibility bridge still needs the finite predecessor penalty and semi-implicit load supplied by the strict-interior condition.
\end{remark}

\section{\label{sec:global}Global Convergence}
\prettyref{sec:Hamiltonian} showed that the mollified per-step \prettyref{pb:PS-CRDS-M} keeps the discrete trajectory norm-bounded (\prettyref{lem:trajectory-bound}) as long as each subproblem is solved under a common finite gradient-residual ceiling. What remains is to establish that the numerical procedure generating those iterates---a line-search SQP scheme applied to the discretized dynamics---converges globally to an approximate ($\EpsStop$-KKT) stationary solution of the modified full-trajectory \prettyref{pb:CITO-M} introduced below, which the bridge (\prettyref{rem:cito-m-bridge}) in turn relates to an $\EpsFJAll$-FJ-feasible point of the original CRDS dynamics. Consistent with the constraint-qualification-free stance of the paper (\prettyref{def:eps-FJ}), we do not assume any constraint qualification a priori. Instead, we show that the small-timestep regime already forces the stationarity Jacobian of~\prettyref{pb:PS-CRDS-M} to be well-conditioned, and that this derived well-posedness is exactly the ingredient a standard SQP method needs to converge.

The original control~\prettyref{pb:CITO} enforces the exact discrete CRDS dynamics of~\prettyref{eq:CRDS} and contains the associated force and contact variables. The SQP method is not applied directly to that MPCC. Instead, it is applied to the following modified full-trajectory NLP, obtained by replacing~\prettyref{eq:CRDS} with the stationarity equation of the per-step mollified and smoothed surrogate~\prettyref{pb:PS-CRDS-M} while retaining the control constraint $u_i\in\mathcal{D}$
\begin{ResizedEq}
\label{pb:CITO-M}
\argmin{\{\theta_i\},\{u_i\in\mathcal{D}\}}\quad&O(\{\theta_i\},\{u_i\})\\
\ST\quad&\nabla_{\theta_i}\PerStepObj_i(\theta_{\gamma^2(i)},\theta_{\gamma(i)},\theta_i,u_i)=0
\quad\forall i\in\mathcal{I}.
\end{ResizedEq}
The optimization variables are only the trajectory $(\{\theta_i\},\{u_i\})$. Unlike~\prettyref{pb:CITO}, the modified problem carries no explicit constraint multipliers $\lambda_i^e,\lambda_i^i$ and no separating-plane variables: the constraints of the original dynamics are absorbed into the reduced penalty inside the per-step objective $\PerStepObj_i$, and the separating planes are eliminated by~\prettyref{cor:SF-CRDS}. This is precisely what makes~\prettyref{pb:CITO-M} a smooth, full-rank NLP amenable to SQP. The link back to~\prettyref{pb:CITO} and its CRDS dynamics is made precise by the bridge (\prettyref{rem:cito-m-bridge}). Accordingly, in the remainder of the convergence analysis, ``the optimized problem'' and ``the SQP problem'' refer to~\prettyref{pb:CITO-M}. Statements about~\prettyref{pb:CITO} are certificate-transfer statements: the returned trajectory is lifted or interpreted as satisfying approximate feasibility or stationarity conditions of the original CITO formulation.

\subsection{\label{sec:global-licq}Curvature-Bounded Stationarity and LICQ}
Fix a smooth full-rank CRDS instance of~\prettyref{pb:PS-CRDS-M} and recall that solving it means driving the gradient of its objective to zero. We therefore make the stationarity residual the primary object of the analysis. For each optimized node $i\in\mathcal{I}$ the per-step residual is the gradient of the per-step objective $\PerStepObj_i$ of~\prettyref{pb:PS-CRDS-M}. We now introduce the stacked trajectory variables $\StackConf\triangleq(\theta_i)_{i\in\mathcal{I}}$ and $\StackCtrl\triangleq(u_i)_{i\in\mathcal{I}}$, so that the user objective reads $O(\StackConf,\StackCtrl)\equiv O(\{\theta_i\},\{u_i\})$. Both stacks contain only optimized nodes: in particular, $\StackConf$ excludes the fixed time-zero configuration $\theta_0$ and its fixed pre-zero predecessor $\theta_{\gamma(0)}$. Solver-iteration superscripts such as $\StackConf^k,\StackCtrl^k$ will denote iterates of these stacks. The per-step residual, the stacked residual, and the associated Jacobians are then
\begin{ResizedEq}
\label{eq:sqp-definition}
\StatRes_i(\theta_i)\triangleq&\nabla_{\theta_i}\PerStepObj_i,\quad
\StatRes(\StackConf)\triangleq(\StatRes_i(\theta_i))_{i\in\mathcal{I}},\\
\StatJac(\StackConf)\triangleq&\FPP{\StatRes}{\StackConf},\quad
\StatJacU(\StackConf)\triangleq\FPP{\StatRes}{\StackCtrl}.
\end{ResizedEq}
Although $\StatRes_i$ is written with the single argument $\theta_i$ for brevity, it also depends on the control $u_i$ through $W(\theta_i,u_i)$---precisely the dependence differentiated by $\StatJacU$. Since $\StatRes_i$ depends on $\StackConf$ only through $\theta_i,\theta_{\gamma(i)},\theta_{\gamma^2(i)}$ via the discretization~\prettyref{eq:discretization}, ordering the timesteps chronologically, the two Jacobians take the banded forms
\begin{ResizedEq}
\label{eq:sqp-jacobian-structure}
\StatJac=&\left(\setlength{\arraycolsep}{2pt}\begin{array}{ccccc}
\ddots & & & & \\
\StatJac_{\gamma^2(i)\gamma^4(i)} & \StatJac_{\gamma^2(i)\gamma^3(i)} & \StatJac_{\gamma^2(i)\gamma^2(i)} & & \\
& \StatJac_{\gamma(i)\gamma^3(i)} & \StatJac_{\gamma(i)\gamma^2(i)} & \StatJac_{\gamma(i)\gamma(i)} & \\
& & \StatJac_{i\gamma^2(i)} & \StatJac_{i\gamma(i)} & \StatJac_{ii}
\end{array}\right),\\
\StatJacU=&\left(\setlength{\arraycolsep}{2pt}\begin{array}{cccc}
\ddots & & & \\
& \FPP{\StatRes_{\gamma^2(i)}}{u_{\gamma^2(i)}} & & \\
& & \FPP{\StatRes_{\gamma(i)}}{u_{\gamma(i)}} & \\
& & & \FPP{\StatRes_i}{u_i}
\end{array}\right),
\end{ResizedEq}
so that each block row of $\StatJac$ has at most three nonzero blocks---the real symmetric diagonal per-step Hessian $\StatJac_{ii}\triangleq\partial\StatRes_i/\partial\theta_i$ and the two earlier-node couplings $\StatJac_{i\gamma(i)},\StatJac_{i\gamma^2(i)}$---while the control block $\StatJacU$ is block diagonal with diagonal blocks $\partial\StatRes_i/\partial u_i$. In particular $\StatJac$ is a generally nonsymmetric block lower-triangular matrix: every block row couples only to the current and two earlier nodes, so nothing sits above the diagonal. For the first optimized node $i_1\triangleq\min\mathcal I$, the stencil predecessors $\theta_0$ and $\theta_{\gamma(0)}$ are fixed data rather than decision variables, so they contribute no Jacobian columns; the $i_1$ block row therefore reduces to its diagonal block $\StatJac_{i_1i_1}$ alone, and the stacked residual and Jacobian remain square with exactly one block row and column per $i\in\mathcal{I}$. The diagonal $\StatJac_{i_1i_1}=M/\delta_{i_1}^2+\text{curvature}$ is governed by the first optimized timestep $\delta_{i_1}$ and is driven to full rank by adaptive refinement of the interval $(0,i_1)$. Consequently, once the diagonal blocks are uniformly invertible, $\StatJac$ is nonsingular, and a block-triangular inverse bound turns the per-block lower bound into a lower bound on $\sigma_{\min}(\StatJac)$ (made precise in~\prettyref{appen:global-convergence-proofs}). The residual $\StatRes_i$ is precisely the gradient residual whose norm was bounded by $\EpsCRDS$ in~\prettyref{lem:H-bound} and by $\LResidual$ in~\prettyref{lem:trajectory-bound}, so the stacked residual $\StatRes(\StackConf)=(\StatRes_i(\theta_i))_{i\in\mathcal{I}}$ is exactly the physics constraint of the modified \prettyref{pb:CITO-M}: $\StatRes(\StackConf)=0$ is feasibility of~\prettyref{pb:CITO-M} and $\|\StatRes(\StackConf)\|\leq\EpsStop$ is its approximate feasibility. A sufficiently small timestep makes each diagonal block strongly positive definite---the acceleration term $M/\delta_i^2$ dominates the bounded curvature of the potential---and hence bounds the smallest singular value of the whole Jacobian away from zero.
\begin{lemma}
\label{lem:jacobian-lower-bound}
Let~\prettyref{pb:PS-CRDS-M} be a smooth full-rank CRDS (\prettyref{def:CRDS},~\ref{def:full-rank}) with $\bar U$ satisfying~\prettyref{lem:mollification}~\ref{it:moll-curv} and $\bar W$ satisfying~\prettyref{lem:controlwork-mollification}, and let the iterates lie in the norm-bounded region guaranteed by~\prettyref{lem:trajectory-bound}. Define the nonsingularity timestep threshold
\begin{equation}
\label{eq:delta-star}
\delta^\star\triangleq\sqrt{\frac{\lambda_{\min}(M)}{\LCurvature+\LControlWorkCurvature+\LLambdaMin}}.
\end{equation}
If, for every $i\in\mathcal{I}$, the per-block singularity check $\sigma_{\min}(\StatJac_{ii})\geq\LLambdaMin>0$ holds---which is the case in particular whenever $\delta_i<\delta^\star$---then the assembled stationarity Jacobian satisfies
\begin{equation}
\sigma_{\min}(\StatJac(\StackConf))\geq\SigmaMinLICQ(\StatJac)>0.
\end{equation}
\end{lemma}
Because each $\StatJac_{ii}$ is real symmetric, $\sigma_{\min}(\StatJac_{ii})=\min_j|\lambda_j(\StatJac_{ii})|$, so the per-block check is a bound on the eigenvalue moduli of the per-step Hessian. The small-timestep condition $\delta_i<\delta^\star$ of~\prettyref{eq:delta-star} is the sufficient condition that discharges that check: it is built only from the mass spectrum $\lambda_{\min}(M)$ and the two curvature moduli $\LCurvature,\LControlWorkCurvature$ already carried by the hypotheses, it is of the same nature as~\prettyref{eq:H-bound-timestep}, and it is enforced by refining the discretization.
\begin{lemma}
\label{lem:jacobian-upper-bound}
For~\prettyref{pb:PS-CRDS-M} a smooth full-rank CRDS with $\bar U$ satisfying~\prettyref{lem:mollification}~\ref{it:moll-curv}, $\bar W$ satisfying~\prettyref{lem:controlwork-mollification}, and the iterates confined to the norm-bounded region guaranteed by~\prettyref{lem:trajectory-bound}, on a discretization scheme with $\min_{i\in\mathcal{I}}\delta_i\geq\delta_{\min}>0$, both $\sigma_{\max}(\StatJac(\StackConf))$ and $\sigma_{\max}(\StatJacU(\StackConf))$ are finite, with a bound that may depend on $\delta_{\min}$; the bound is uniform over the norm-bounded region.
\end{lemma}

Together the two lemmas bound the extreme singular values of the stationarity Jacobian,
\begin{ResizedEq}
\label{eq:sqp-conditioning}
\SigmaMinLICQ(\StatJac)\leq\sigma_{\min}(\StatJac(\StackConf))\leq\sigma_{\max}(\StatJac(\StackConf))<\infty,
\end{ResizedEq}
so that, in the small-timestep regime, $\StatJac$ is simultaneously bounded below and uniformly well-conditioned---the two ingredients that a standard SQP method requires to make progress and to keep its multipliers and search directions bounded. We stress that~\prettyref{eq:sqp-conditioning} is derived from small $\delta_i$ and bounded curvature rather than assumed: \prettyref{lem:jacobian-lower-bound} produces the LICQ modulus $\SigmaMinLICQ(\StatJac)$ as a consequence of the discretization, so the analysis never posits a constraint qualification. This is consistent with the $\EpsFJAll$-FJ, CQ-free formulation of~\prettyref{def:eps-FJ}: the conditioning of $\StatJac$ is an FJ-compatible property of the discretized dynamics, not a standing assumption on the original problem.

\subsection{\label{sec:global-sqp}Adaptive-Timestep SQP}
The certificate of~\prettyref{lem:jacobian-lower-bound} is the per-block check $\sigma_{\min}(\StatJac_{ii})\geq\LLambdaMin$, guaranteed once $\delta_i<\delta^\star$; but this admissible timestep depends on the potential curvature $\LCurvature$ and the control-work curvature $\LControlWorkCurvature$, which are not known a priori. We therefore let the solver enforce the checkable proxy adaptively: whenever a diagonal block is detected close to singular, $\sigma_{\min}(\StatJac_{ii})<\LLambdaMin$, the offending interval is bisected, which shrinks $\delta_i$ and, by~\prettyref{lem:jacobian-lower-bound}, eventually reaches a strongly positive-definite regime that guarantees the nonsingularity check; a well-conditioned indefinite block may pass before that regime is reached. This couples the discretization refinement directly to the well-posedness certificate of~\prettyref{lem:jacobian-lower-bound}. Building on this refinement, we solve the full-trajectory~\prettyref{pb:CITO-M} with a line-search SQP method.~\prettyref{pb:PS-CRDS-M} supplies the per-step stationarity residual $\StatRes_i=\nabla_{\theta_i}\PerStepObj_i$ that forms the physics constraints of~\prettyref{pb:CITO-M}. The SQP machinery is standard in its quadratic-program and merit-function components but is augmented with three problem-specific safeguards: an adaptive subdivision that enforces the singularity check, a constraint-satisfaction phase that secures the feasibility margin whenever the residual leaves the working band, and a safeguarded line search whose second branch keeps the constraint violation bounded. We develop the three ingredients in turn and then assemble them into~\prettyref{alg:sqp}, whose global convergence is~\prettyref{thm:sqp-convergence}. Throughout, every iterate is confined to the compact, norm-bounded region of~\prettyref{lem:trajectory-bound}.

\subsubsection{Stationarity Subproblem and Merit Function}
We monitor progress with the $\ell_1$ merit function
\begin{equation}
\label{eq:sqp-merit}
\phi_\varrho(\StackConf,\StackCtrl)\triangleq O(\StackConf,\StackCtrl)+\varrho\|\StatRes\|_1,
\end{equation}
whose first term $O$ drives stationarity and whose second term penalizes the residual through the penalty weight $\varrho>0$. The search direction at iterate $k$ is the pair $s^k=(s_\StackConf^k,s_\StackCtrl^k)$ of configuration and control steps solving the quadratic program
\begin{ResizedEq}
\label{pb:sqp-qp}
\text{QP}^k:
\begin{cases}
\argmin{s_\StackConf^k,\,s_\StackCtrl^k}&[s_\StackConf^k]^T\FPP{O}{\StackConf}^k+[s_\StackCtrl^k]^T\FPP{O}{\StackCtrl}^k\\
&+\tfrac12[s_\StackConf^k]^TB_\StackConf^ks_\StackConf^k+\tfrac12[s_\StackCtrl^k]^TB_\StackCtrl^ks_\StackCtrl^k\\
\ST&\StatRes^k+\StatJac^ks_\StackConf^k+\StatJacU^ks_\StackCtrl^k=0,\\
&u_i^k+s_{u_i}^k\in\mathcal{D},
\end{cases}
\end{ResizedEq}
where the objective gradients $\FPP{O}{\StackConf}^k$ and $\FPP{O}{\StackCtrl}^k$ are well defined by~\prettyref{ass:objective}, the superscript $k$ denotes evaluation at the $k$th iterate and $B_\StackConf^k,B_\StackCtrl^k$ are positive-definite Hessian approximations chosen to satisfy the uniform conditioning bound
\begin{ResizedEq}
\label{eq:sqp-hessian}
1/\LHessianCond\leq&\lambda_{\min}(B^k)\leq\lambda_{\max}(B^k)\leq\LHessianCond\\
B^k\triangleq&\diag(B_\StackConf^k,B_\StackCtrl^k),
\end{ResizedEq}
with conditioning constant $\LHessianCond\geq1$. Solving $\text{QP}^k$ returns, alongside the step $s^k$, the Lagrange multipliers $\mu_\StackConf^k$ (resp. $\mu_\StackCtrl^k$) of its physics (resp. control) constraint, and we write $\mu_{\theta_i}^k$, $\mu_{u_i}^k$, and $s_{u_i}^k$ for the timestep-$i$ sub-blocks of the stacked $\mu_\StackConf^k$, $\mu_\StackCtrl^k$, and $s_\StackCtrl^k$. These are exactly the ingredients a line-search SQP consumes; the derived conditioning~\prettyref{eq:sqp-conditioning} is what guarantees that $\text{QP}^k$ is solvable and that its step and multipliers stay bounded. The associated KKT optimality system of $\text{QP}^k$, on which the step and multiplier bounds rest, is stated where it is first used, in~\prettyref{appen:global-convergence-proofs}.

\subsubsection{Singularity Check and Adaptive Subdivision}
The certificate of~\prettyref{lem:jacobian-lower-bound} is global---it bounds $\sigma_{\min}(\StatJac)$ over the assembled Jacobian---but the solver needs a local, checkable proxy that can be evaluated node by node without assembling $\StatJac$. We use the minimum singular value of the diagonal block: the per-block test $\sigma_{\min}(\StatJac_{ii})\geq\LLambdaMin$ is exactly the hypothesis of~\prettyref{lem:jacobian-lower-bound}, which then yields $\sigma_{\min}(\StatJac)\geq\SigmaMinLICQ(\StatJac)>0$; and by~\prettyref{lem:jacobian-lower-bound} a sufficiently small timestep $\delta_i<\delta^\star$ guarantees the test, so timestep $i$ is then fine enough for the LICQ argument to apply. The current discretization is represented solely by the finite ordered optimized set $\mathcal I$; physical times are $t_i=\delta i$, and $\gamma$ is derived from the total order on $\{0\}\cup\mathcal I$. Whenever an optimized node $i$ fails the check,~\prettyref{alg:subd} inserts the midpoint label $m=(i+\gamma(i))/2$ into $\mathcal I$, initializes $\theta_m\gets(\theta_i+\theta_{\gamma(i)})/2$, and initializes $u_m\gets(u_i+u_{\gamma(i)})/2$. Convexity of $\mathcal D$ preserves control feasibility, and the insertion halves $\delta_i$, restoring the check after finitely many bisections by~\prettyref{lem:jacobian-lower-bound}. The fixed history interval $(\gamma(0),0)$ is never subdivided. The first refinable interval is $(0,\min\mathcal I)$, where the fixed values $\theta_0,u_0$ supply the predecessor data; subdivision never modifies $\theta_{\gamma(0)},\theta_0$, or $u_0$.
\begin{algorithm}[t]
\caption{Subdivide($\mathcal I,\StackConf,\StackCtrl,i$)}
\label{alg:subd}
\begin{algorithmic}[1]
\LineComment{$\gamma$ is determined by the current total order}
\State $m\gets(i+\gamma(i))/2$
\State $\theta_m\gets(\theta_i+\theta_{\gamma(i)})/2$
\State $u_m\gets(u_i+u_{\gamma(i)})/2$
\State $\mathcal I\gets\mathcal I\cup\{m\}$
\State Insert $\theta_m,u_m$ into $\StackConf,\StackCtrl$ in increasing order
\State Return $\mathcal I,\StackConf,\StackCtrl$
\end{algorithmic}
\end{algorithm}

\subsubsection{Constraint Satisfaction Phase}
The bounded-penalty guarantee of~\prettyref{thm:sqp-convergence} is established only once the residual has entered a working band $\|\StatRes\|_1<\LResidualWork$, where $\LResidualWork>0$ is an algorithm input fixed before the solve---a working residual-norm threshold---that delimits the region on which the guarantee applies: the step- and multiplier-bounds behind it are stated in terms of a violation ceiling, so an arbitrary initialization---or any subsequent iterate that leaves the band, for example after a subdivision---must first be driven into that band. Because $\|\StatRes_i\|\leq\|\StatRes\|_1$, once inside the working band $\|\StatRes\|_1<\LResidualWork$ every iterate satisfies the per-step bound $\|\StatRes_i\|\leq\LResidualWork\leq\LResidual$. During the constraint-satisfaction phase itself, the Armijo test~\prettyref{eq:sqp-constraint-linesearch} makes $\|\StatRes^k\|_2$ non-increasing. Hence norm equivalence gives the finite ceiling $\|\StatRes^k\|_1\leq\sqrt{\dim\StatRes}\,\|\StatRes^0\|_2<\infty$ throughout; therefore every iterate---in either phase---meets the finite residual-ceiling hypothesis of~\prettyref{lem:trajectory-bound} and stays inside its norm-bounded region. Because the singularity check makes $\StatJac^k$ invertible, we do this with damped Newton steps on the residual, $s_\StackConf^k=-[\StatJac^k]^{-1}\StatRes^k$, accepted by the simpler Armijo condition on $\tfrac12\|\StatRes\|^2$, i.e.,
\begin{ResizedEq}
\label{eq:sqp-constraint-linesearch}
\tfrac12\|\StatRes(\StackConf^k+[\tau^\ell]^ks_\StackConf^k)\|^2<(1-2\alpha[\tau^\ell]^k)\tfrac12\|\StatRes^k\|^2,
\end{ResizedEq}
with backtracking contraction factor $\tau\in(0,1)$ and Armijo constant $\alpha\in(0,1)$; here $\ell\in\mathbb{Z}^+$ indexes the backtracking so that $[\tau^\ell]^k$ is the accepted step size.~\prettyref{alg:constraint-solve} runs this iteration until either the working threshold $\LResidualWork$ is met or a node fails the singularity check and is returned for subdivision.
\begin{algorithm}[t]
\caption{Constraint-Solve($\mathcal I,\StackConf,\StackCtrl,\LLambdaMin,\LResidualWork$)}
\label{alg:constraint-solve}
\begin{algorithmic}[1]
\State $\StackConf^0\gets\StackConf$
\While{$k=0,1,2,\cdots$}
\For{$i\in\mathcal I$ in increasing order}
\If{$\sigma_{\min}(\StatJac_{ii}^k)<\LLambdaMin$}
\State Return $\mathcal I,\StackConf^k,i$\label{ln:sqp-constraint-singular}
\EndIf
\EndFor
\If{$\|\StatRes^k\|_1<\LResidualWork$}
\State Return $\mathcal I,\StackConf^k,\emptyset$\label{ln:sqp-constraint-return}
\EndIf
\State $s_\StackConf^k\gets-[\StatJac^k]^{-1}\StatRes^k$
\LineComment{$\tau\in(0,1),\alpha\in(0,1)$}
\State Find smallest $\ell\in\mathbb{Z}^+$ satisfying~\prettyref{eq:sqp-constraint-linesearch}
\State $\StackConf^{k+1}\gets\StackConf^k+[\tau^\ell]^ks_\StackConf^k$
\EndWhile
\end{algorithmic}
\end{algorithm}

\subsubsection{Safeguarded Inner SQP}
Once inside the working band, the inner solver takes the QP steps of~\prettyref{pb:sqp-qp} and globalizes them with the safeguarded line search
\begin{ResizedEq}
\label{eq:sqp-linesearch}
\begin{cases}
\phi_{\varrho^k}(\StackConf^k+[\tau^\ell]^ks_\StackConf^k,\StackCtrl^k+[\tau^\ell]^ks_\StackCtrl^k)\\
\qquad<\phi_{\varrho^k}^k+\alpha[\tau^\ell]^k[\bar D\phi_\varrho]^k
&\text{if }\|\StatRes^k\|_1\leq\LResidualWork\\
\begin{rcases}
\phi_{\varrho^k}(\StackConf^k+[\tau^\ell]^ks_\StackConf^k,\StackCtrl^k+[\tau^\ell]^ks_\StackCtrl^k)\\
\qquad<\phi_{\varrho^k}^k+\alpha[\tau^\ell]^k[\bar D\phi_\varrho]^k\\
\|\StatRes(\StackConf^k+[\tau^\ell]^ks_\StackConf^k,\StackCtrl^k+[\tau^\ell]^ks_\StackCtrl^k)\|_1\\
\qquad\leq\|\StatRes^k\|_1
\end{rcases}&\text{otherwise,}
\end{cases}
\end{ResizedEq}
where $[\bar D\phi_\varrho]^k\triangleq-\tfrac{1}{\LHessianCond}\|s^k\|^2-\bar\varrho\|\StatRes^k\|_1$ is a guaranteed-descent surrogate for the merit's directional derivative---negative by construction, with penalty margin $\bar\varrho>0$---and $\alpha$ enforces the Armijo condition. The first branch is the usual sufficient-decrease test, used whenever the violation is already small, $\|\StatRes^k\|_1\leq \LResidualWork$. The second branch is the problem-specific safeguard: outside the working band it additionally requires the violation not to increase, $\|\StatRes(\cdots)\|_1\leq\|\StatRes^k\|_1$. This extra branch prevents the merit's penalty term from having to chase a growing residual and is the mechanism that keeps $\|\StatRes\|_1$---and hence the penalty sequence $\{\varrho^k\}$---uniformly bounded. After each accepted step the penalty is raised just enough to dominate the current multipliers, $\varrho^k\gets\max(\varrho^k,\max_i\|\mu_{\theta_i}^k\|+\bar\varrho)$, so that the surrogate remains a strict descent direction of the merit.~\prettyref{alg:sqp-inner} iterates this until $\|s^k\|<\EpsStop/\LHessianCond$ and $\|\StatRes^k\|_1<\EpsStop$ at the stopping tolerance $\EpsStop$, or until a node fails the singularity check.
\begin{algorithm}[t]
\caption{\footnotesize{SQP-Inner($\mathcal I,\StackConf,\StackCtrl,\varrho^0,\LLambdaMin,\LResidualWork,\LHessianCond,\bar\varrho,\EpsStop$)}}
\label{alg:sqp-inner}
\begin{algorithmic}[1]
\Require{\tiny{$u_i\in\mathcal{D},\forall i,\varrho^0>0,\LLambdaMin>0,\LResidualWork>0,\LHessianCond\geq1,\bar\varrho>0,\EpsStop>0$}}
\State $\StackConf^0\gets\StackConf,\StackCtrl^0\gets\StackCtrl$
\While{$k=0,1,2,\cdots$}
\For{$i\in\mathcal I$ in increasing order}
\If{$\sigma_{\min}(\StatJac_{ii}^k)<\LLambdaMin$}
\State Return $\mathcal I,\StackConf^k,\StackCtrl^k,\varrho^k,i$\label{ln:sqp-inner-singular}
\EndIf
\EndFor
\State Choose $B_\StackConf^k,B_\StackCtrl^k$ satisfying~\prettyref{eq:sqp-hessian}
\State $s^k,\mu_\StackConf^k,\mu_\StackCtrl^k\gets\text{QP}^k$
\If{$\|s^k\|<\EpsStop/\LHessianCond$ and $\|\StatRes^k\|_1<\EpsStop$}
\State Return $\mathcal I,\StackConf^k,\StackCtrl^k,\varrho^k,\emptyset$\label{ln:sqp-inner-return}
\EndIf
\State $\varrho^k\gets\max(\varrho^k,\max_i\|\mu_{\theta_i}^k\|+\bar\varrho)$
\LineComment{$\tau\in(0,1),\alpha\in(0,1)$}
\State Find smallest $\ell\in\mathbb{Z}^+$ satisfying~\prettyref{eq:sqp-linesearch}
\State $\StackConf^{k+1}\gets\StackConf^k+[\tau^\ell]^ks_\StackConf^k$;\quad $\StackCtrl^{k+1}\gets\StackCtrl^k+[\tau^\ell]^ks_\StackCtrl^k$
\EndWhile
\end{algorithmic}
\end{algorithm}

\subsubsection{Convergence of SQP}
\prettyref{alg:sqp} assembles the three ingredients: it subdivides any node returned by the singularity check, calls the constraint-satisfaction phase whenever $\|\StatRes(\StackConf)\|_1\geq \LResidualWork$, and otherwise runs the inner SQP, returning when the inner solver meets the stopping tolerance. The method is thus parameterized by the initial finite ordered optimized set $\mathcal I$, finite optimized configuration stack $\StackConf$, and control stack $\StackCtrl$ (with $u_i\in\mathcal{D}$ for all $i$), the initial penalty $\varrho$, the per-block threshold $\LLambdaMin$, the working violation threshold $\LResidualWork$, the QP-Hessian conditioning bound $\LHessianCond$, the penalty margin $\bar\varrho$, and the stopping tolerance $\EpsStop$. The ordered set $\mathcal I$ and stacked configuration $\StackConf=(\theta_i)_{i\in\mathcal I}$ are distinct objects; physical times and predecessors are derived from the ordered labels. Constraint-Solve holds $\mathcal I$ and $\StackCtrl$ fixed while updating $\StackConf$, SQP-Inner holds $\mathcal I$ fixed while updating both stacks, and subdivision enlarges $\mathcal I$ and inserts finite interpolated values into both optimized stacks.

The initial triple $(\mathcal I^0,\StackConf^0,\StackCtrl^0)$ may consist of any finite nonempty totally ordered optimized set $\mathcal I^0\subset\{i\in\mathbb R:0<\delta i\leq H\}$, any finite optimized configuration vector, and any controls $\StackCtrl^0\in\mathcal D^{|\mathcal I^0|}$. No feasibility, small-residual, or proximity condition is imposed on the initial configurations. Since the data defining~\prettyref{pb:CITO-M} are globally defined and the damping smoothing of~\prettyref{def:S-D} is everywhere finite, every such input has a finite residual and finite residual Jacobians.

\begin{algorithm}[t]
\caption{\small{AdaptiveSQP($\mathcal I,\StackConf,\StackCtrl,\varrho,\LLambdaMin,\LResidualWork,\LHessianCond,\bar\varrho,\EpsStop$)}}
\label{alg:sqp}
\begin{algorithmic}[1]
\Require{\tiny{$u_i\in\mathcal{D},\forall i,\varrho>0,\LLambdaMin>0,\LResidualWork>0,\LHessianCond\geq1,\bar\varrho>0,\EpsStop>0$}}
\State $i\gets\emptyset$
\While{Not converged}
\If{$i\neq\emptyset$}
\State $\mathcal I,\StackConf,\StackCtrl\gets$Subdivide($\mathcal I,\StackConf,\StackCtrl,i$);\quad $i\gets\emptyset$
\EndIf
\If{$\|\StatRes(\StackConf)\|_1\geq\LResidualWork$}
\State $\mathcal I,\StackConf,i\gets$Constraint-Solve($\mathcal I,\StackConf,\StackCtrl,\LLambdaMin,\LResidualWork$)
\Else
\State{\tiny{$\mathcal I,\StackConf,\StackCtrl,\varrho,i\gets$SQP-Inner($\mathcal I,\StackConf,\StackCtrl,\varrho,\LLambdaMin,\LResidualWork,\LHessianCond,\bar\varrho,\EpsStop$)}}
\If{$i=\emptyset$}
\State Return $\mathcal I,\StackConf,\StackCtrl,\varrho$
\EndIf
\EndIf
\EndWhile
\end{algorithmic}
\end{algorithm}

Before stating the convergence theorem, we record the perturbed optimality condition that its returned point certifies. Writing $N_{\mathcal{D}}(u)$ for the normal cone of the convex control set $\mathcal{D}$ at $u\in\mathcal{D}$, we introduce the following $\EpsStop$-KKT perturbed optimality condition of the modified \prettyref{pb:CITO-M}.
\begin{definition}[$\EpsStop$-KKT condition]
\label{def:eps-KKT}
A point $(\StackConf,\StackCtrl)$ satisfies the $\EpsStop$-KKT condition of the modified \prettyref{pb:CITO-M} if there exist multipliers $\mu_\StackConf$, $\mu_\StackCtrl$ and, for every timestep $i$, an auxiliary control $u_i'$, such that
\begin{ResizedEq}
\label{eq:eps-kkt-system}
&\EpsStop\text{-KKT}:\\
&\begin{cases}
\left\|\FPP{O}{\StackConf}+[\StatJac]^T\mu_\StackConf\right\|\leq\EpsStop\\
\left\|\FPP{O}{\StackCtrl}+[\StatJacU]^T\mu_\StackConf+\mu_\StackCtrl\right\|\leq\EpsStop\\
u_i\in\mathcal{D},\ u_i'\in\mathcal{D},\ \|u_i-u_i'\|\leq\EpsStop,\ \mu_{u_i}\in N_{\mathcal{D}}(u_i')\\
\|\StatRes\|_1\leq\EpsStop.
\end{cases}
\end{ResizedEq}
The third line is the approximate (two-point) control-stationarity condition: the control multiplier $\mu_{u_i}$ is certified in the normal cone at the nearby feasible point $u_i'$ rather than at the returned control $u_i$, the two lying within $\EpsStop$ of one another. This offset arises because the SQP quadratic subproblem certifies the control multiplier in the normal cone at its own accepted iterate $u_i'=u_i+s_{u_i}$, i.e.\ $\mu_{u_i}\in N_{\mathcal{D}}(u_i+s_{u_i})$, rather than at the incoming control $u_i$. Here $\mu_\StackConf$ (resp. $\mu_\StackCtrl$) is the multiplier of the physics (resp. control) constraint, and $\mu_{u_i}$ is the timestep-$i$ sub-block of $\mu_\StackCtrl$. As $\EpsStop\to0$ the two points $u_i,u_i'$ collapse and the system reduces to the exact KKT conditions of~\prettyref{pb:CITO-M}.
\end{definition}

The convergence statement is a specialization of the classical global-convergence theory for line-search SQP~\cite{bertsekas1997nonlinear,solodov2009global} to our discretized dynamics; we phrase it in the paper's own residual and FJ currency.
\begin{theorem}
\label{thm:sqp-convergence}
Let~\prettyref{pb:PS-CRDS-M} be a smooth full-rank CRDS whose objective $O$ satisfies~\prettyref{ass:objective}. Fix any finite nonempty totally ordered initial optimized set $\mathcal I^0\subset\{i\in\mathbb R:0<\delta i\leq H\}$, any finite optimized configuration guess $\StackConf^0\in\mathbb R^{n|\mathcal I^0|}$, and any optimized control guess $\StackCtrl^0\in\mathcal D^{|\mathcal I^0|}$, and run~\prettyref{alg:sqp} with parameters meeting its input requirements. Then~\prettyref{alg:sqp} is well defined, performs only finitely many subdivisions, generates a bounded merit-penalty sequence $\{\varrho^k\}$, and terminates after finitely many iterations at a point $(\StackConf,\StackCtrl)$ satisfying the $\EpsStop$-KKT condition (\prettyref{def:eps-KKT}) of the modified \prettyref{pb:CITO-M}.
\end{theorem}
\begin{remark}[Bridge from $\EpsStop$-KKT of~\prettyref{pb:CITO-M} to approximate FJ-feasibility of~\prettyref{pb:CITO}]
\label{rem:cito-m-bridge}
The approximate stationary solution of~\prettyref{pb:CITO-M} certified by~\prettyref{thm:sqp-convergence} relates back to the original dynamics of~\prettyref{pb:CITO}. Since the returned point has $\|\StatRes\|_1\leq\EpsStop$, each per-step residual obeys $\|\StatRes_i\|_2\leq\|\StatRes\|_1\leq\EpsStop$, so it meets the gradient-residual hypothesis of~\prettyref{lem:eps-SR-FJ} with $\EpsCRDS=\EpsStop$. The remaining standing hypothesis $V_R(\theta_{\gamma(i)}^\star,\EpsFJ)\leq\LHamilton$ holds whenever the per-step Hamiltonian bound $\bar{\mathcal H}_i^+\leq\LHamBound(t_i)$ is in force: by~\prettyref{lem:H-bound} and the constant-ordering remark of~\prettyref{sec:Hamiltonian}, that bound gives $V_R\leq\tilde U\leq\bar{\mathcal H}_i^+\leq\LHamBound(H)\leq\LHamilton$ on the region where the clamps are inactive. At the first optimized node $i_1\triangleq\min\mathcal I$, the predecessor $\theta_0$ is the fixed strict-interior node of~\prettyref{ass:Slater}, with finite energy base $\bar{\mathcal H}_0^+$ seeding the chain; for every later optimized node, $\theta_{\gamma(i)}^\star$ is an optimized predecessor already covered by the Hamiltonian chain. That energy bound, however, requires the Hamiltonian timestep condition~\prettyref{eq:H-bound-timestep}, which~\prettyref{alg:sqp} does not enforce; it is secured only by the top-level~\prettyref{alg:CITO} (see~\prettyref{thm:CITO-convergence} and~\prettyref{rem:cito-ham-discharge}). Conditional on this energy bound, and since~\prettyref{pb:PS-CRDS-M} is a smooth full-rank CRDS by~\prettyref{cor:SF-CRDS} and~\prettyref{lem:mollification}, each optimized $\theta_i$---including the first optimized node $\theta_{\min\mathcal I}$---induces an extended configuration satisfying the $\EpsFJAll$-FJ condition of the original extended constrained CRDS~\prettyref{pb:CRDS-ABD-extended}. This connects the modified problem's approximate stationary solution back to the approximate FJ-feasibility of the original CRDS dynamics enforced by~\prettyref{pb:CITO}.
\end{remark}

\subsubsection{Convergence of CITOSolve}
\prettyref{thm:sqp-convergence} certifies that~\prettyref{alg:sqp} drives~\prettyref{pb:CITO-M} to an $\EpsStop$-KKT point on any fixed mesh, but its mesh control is one-sided: the singularity check of~\prettyref{ln:sqp-inner-singular} refines only until the sufficient LICQ regime $\delta_i<\delta^\star$ of~\prettyref{lem:jacobian-lower-bound} is met, and nothing in~\prettyref{alg:sqp} enforces the complementary Hamiltonian timestep condition~\prettyref{eq:H-bound-timestep} that~\prettyref{lem:H-bound} requires. Consequently~\prettyref{alg:sqp} alone does not certify the per-step energy bound $\bar{\mathcal H}_i^+\leq\LHamBound(t_i)$; and that bound is precisely what discharges the standing hypothesis $V_R(\theta_{\gamma(i)}^\star,\EpsFJ)\leq\LHamilton$ of~\prettyref{lem:eps-SR-FJ} through the bridge \prettyref{rem:cito-m-bridge}. We close this gap with a top-level solver, which we call CITOSolve, that wraps~\prettyref{alg:sqp} in an outer loop: after each SQP solve it scans the one-step Hamiltonian recursion in increasing order, subdivides the first failing interval, and restarts the outer solve until every node passes.

The test is the discrete Gr\"onwall step underlying~\prettyref{lem:H-bound}. Because~\prettyref{alg:sqp} returns $\|\StatRes\|_1\leq\EpsStop$, hence $\|\StatRes_i\|_2\leq\EpsStop$ at every node, the one-step estimate of~\prettyref{appen:Hamiltonian-upper-bound-proofs} applies with $\EpsCRDS=\EpsStop$, and a node with a small enough timestep obeys
\begin{ResizedEq}
\label{eq:cito-ham-check}
\bar{\mathcal H}_i^+\leq e^{6\delta_i/\lambda_{\min}(M)}\left(\bar{\mathcal H}_{\gamma(i)}^++\frac{\delta_i}{2}\left((\|g\|+\LControlWorkForce)^2+\EpsStop^2\right)\right).
\end{ResizedEq}
The solver checks~\prettyref{eq:cito-ham-check} for each $i\in\mathcal I$ in increasing order, records the first failure, and only after ending that scan bisects the failing interval through~\prettyref{alg:subd}; it never mutates $\mathcal I$ inside a continuing scan, and exits only when every node passes. Because the rotation constraints here are orthogonality-only with no orientation logic, the Hamiltonian check~\prettyref{eq:cito-ham-check} is the sole subdivision trigger.
\begin{repalgorithm}{alg:CITO}
\caption{\small{CITOSolve($\mathcal I,\StackConf,\StackCtrl,\varrho,\LLambdaMin,\LResidualWork,\LHessianCond,\bar\varrho,\EpsStop$)}}
\begin{algorithmic}[1]
\Require{\tiny{$u_i\in\mathcal{D},\forall i,\varrho>0,\LLambdaMin>0,\LResidualWork>0,\LHessianCond\geq1,\bar\varrho>0,\EpsStop>0$}}
\While{Not converged}
\State \small{$\mathcal I,\StackConf,\StackCtrl,\varrho\gets$AdaptiveSQP($\mathcal I,\StackConf,\StackCtrl,\varrho,\LLambdaMin,\LResidualWork,\LHessianCond,\bar\varrho,\EpsStop$)}
\State $i_{\mathrm{fail}}\gets\emptyset$
\For{$i\in\mathcal I$ in increasing order}
\If{\prettyref{eq:cito-ham-check} fails}\label{ln:cito-ham-check}
\State $i_{\mathrm{fail}}\gets i$; \textbf{break}
\EndIf
\EndFor
\If{$i_{\mathrm{fail}}=\emptyset$}
\State Return $\mathcal I,\StackConf,\StackCtrl$\label{ln:cito-return}
\EndIf
\State $\mathcal I,\StackConf,\StackCtrl\gets$Subdivide($\mathcal I,\StackConf,\StackCtrl,i_{\mathrm{fail}}$)
\EndWhile
\end{algorithmic}
\end{repalgorithm}
The initial control $\StackCtrl$ (with $u_i\in\mathcal{D}$ for all $i$) and the penalty $\varrho$ are threaded across calls---warm-starting each solve, with $\varrho$ nondecreasing---while the remaining parameters $\LLambdaMin,\LResidualWork,\LHessianCond,\bar\varrho,\EpsStop$ are exactly those consumed by~\prettyref{alg:sqp}; each outer iteration is one call to~\prettyref{alg:sqp} followed by a sweep of the check~\prettyref{eq:cito-ham-check}. Our final result certifies that this outer loop terminates and that its output carries the full-trajectory energy bound in addition to the $\EpsStop$-KKT guarantee of~\prettyref{thm:sqp-convergence}.
\begin{theorem}
\label{thm:CITO-convergence}
Fix an instance of~\prettyref{pb:PS-CRDS-M} that is a smooth full-rank CRDS whose objective $O$ satisfies~\prettyref{ass:objective}; its damping function, including any smoothing parameters used to construct it, and all other problem data remain fixed throughout the run. Fix startup data satisfying~\prettyref{ass:Slater}, including the fixed strict-interior node $\theta_0$ and its fixed pre-zero predecessor $\theta_{\gamma(0)}$, with the clamps $\LHamilton,\LPosition$ fixed by the constant-ordering remark of~\prettyref{sec:Hamiltonian}. For any finite nonempty totally ordered initial optimized set $\mathcal I^0\subset\{i\in\mathbb R:0<\delta i\leq H\}$, any finite optimized configuration guess indexed by $\mathcal I^0$, and any optimized control guess in $\mathcal D^{|\mathcal I^0|}$, run~\prettyref{alg:CITO} with parameters meeting the input requirements of~\prettyref{alg:sqp}. Then~\prettyref{alg:CITO} is well defined, performs finitely many outer iterations---each a finite call to~\prettyref{alg:sqp}---and finitely many subdivisions, generates a bounded penalty sequence, and terminates at a point $(\StackConf,\StackCtrl)$ satisfying the $\EpsStop$-KKT condition (\prettyref{def:eps-KKT}) of~\prettyref{pb:CITO-M} and, in addition, the per-step Hamiltonian bound
\begin{equation}
\label{eq:cito-ham-bound}
\bar{\mathcal H}_i^+\leq\LHamBound(t_i)\qquad\forall i\in\mathcal{I}.
\end{equation}
\end{theorem}
The proof is given in~\prettyref{appen:global-convergence-proofs}. It reuses the conclusions of~\prettyref{thm:sqp-convergence} for each SQP call, bounds the total subdivision count via the sufficient singularity-refinement threshold $\delta^\star$ and the Hamiltonian threshold $\delta_H\triangleq\min(1/\LCurvature,\lambda_{\min}(M)/6)$, and closes with a Gr\"onwall telescoping of~\prettyref{eq:cito-ham-check} to obtain~\prettyref{eq:cito-ham-bound}.
\begin{remark}[Making the~\prettyref{pb:CITO-M}-to-\prettyref{pb:CITO} feasibility bridge unconditional]
\label{rem:cito-ham-discharge}
The per-step bound~\prettyref{eq:cito-ham-bound} discharges the remaining hypothesis of~\prettyref{rem:cito-m-bridge}. With the clamps fixed by the constant-ordering remark of~\prettyref{sec:Hamiltonian} so that $\LHamilton>\LHamBound(H)$ over the horizon, the reduced potential obeys $V_R(\theta_{\gamma(i)}^\star,\EpsFJ)\leq\tilde U(\theta_{\gamma(i)}^\star,\EpsFJ)\leq\bar{\mathcal H}_{\gamma(i)}^+\leq\LHamBound(t_{\gamma(i)})\leq\LHamBound(H)\leq\LHamilton$ on the region where the clamps are inactive, so the standing hypothesis of~\prettyref{lem:eps-SR-FJ} is discharged by~\prettyref{alg:CITO} itself rather than assumed. At the first optimized node $i_1\triangleq\min\mathcal I$, this chain is seeded by the fixed strict-interior node $\theta_0$, whose finite base energy $\bar{\mathcal H}_0^+=\LHamBound(0)$ and well-defined first optimized penalty load discharge the hypothesis for $\theta_{i_1}$; for every later optimized node, the predecessor is an optimized node already covered by the chain. Combined with $\|\StatRes_i\|_2\leq\EpsStop$ from the $\EpsStop$-KKT conclusion, every optimized $\theta_i$ of the CITOSolve output---including $\theta_{\min\mathcal I}$---therefore induces an extended configuration satisfying the $\EpsFJAll$-FJ condition of~\prettyref{pb:CRDS-ABD-extended}, so the top-level solver unconditionally certifies approximate FJ-feasibility of the original CRDS.
\end{remark}

\subsection{\label{sec:global-transfer}Stationarity Transfer}
We are ready to give the formal version of~\prettyref{thm:stationarity-transfer}, previewed informally in the Problem Definition section.~\prettyref{thm:CITO-convergence} certifies that~\prettyref{alg:CITO} returns a point that is $\EpsStop$-KKT for the modified \prettyref{pb:CITO-M}, whose physics constraint is the per-step stationarity $\nabla_{\theta_i}\PerStepObj_i=0$ of the mollified surrogate~\prettyref{pb:PS-CRDS-M}. This is stationarity of a surrogate rather than of the original MPCC~\prettyref{pb:CITO}, whose physics constraint is the exact discrete CRDS dynamics~\prettyref{eq:CRDS}. The bridge remarks~\prettyref{rem:cito-m-bridge} and~\ref{rem:cito-ham-discharge} connect the returned point to the $\EpsFJAll$-FJ feasibility of the CRDS, but they stop short of certifying an approximate stationarity condition of~\prettyref{pb:CITO} itself. Closing this remaining gap is naturally a two-stage program: a first stage that assumes the CRDS underlying~\prettyref{pb:CITO} is already smooth (\prettyref{def:smooth-CRDS}), so that no damping smoothing is required and $\Damp_{\EpsFri,\LHamilton}^V=\Damp^V$; and a second stage that additionally handles the nonsmooth damping of general frictional contact. This subsection carries out both stages: the first directly, and the second through a normal-compatible smoothing (\prettyref{def:S-D-compatible}) whose fixed-smoothing stationarity certificate is normalized and lifted onto the graph of the nonsmooth damping subdifferential, yielding a homogeneous Clarke--FJ certificate. The concrete verification of that smoothing is discharged in full for the articulated-body frictional-contact model.

It is useful to separate three conclusions, organized together with the per-step physical condition and the external standard MPCC hierarchy in~\prettyref{table:certificate-scope}. First,~\prettyref{thm:CITO-convergence} gives the direct $\EpsStop$-KKT result for~\prettyref{pb:CITO-M}. Second,~\prettyref{rem:cito-m-bridge}--\ref{rem:cito-ham-discharge} provide the feasibility bridge: the Hamiltonian bound and per-step lifting lemma recover approximate physical feasibility for the original constrained CRDS. Third, the smooth and nonsmooth stationarity transfers identify the KKT residual of~\prettyref{pb:CITO-M} with the force-parameterized force-balance residual on the clamp-inactive region and transfer the adjoint relations, yielding the Fritz--John certificate statements of~\prettyref{thm:stationarity-transfer}--\ref{thm:nonsmooth-transfer} for~\prettyref{pb:CITO}.

\subsubsection{Smooth Force-Parameterized Transfer}
We first treat the case in which the CRDS underlying~\prettyref{pb:CITO} is smooth. On the clamp-inactive region, the stationarity residual of~\prettyref{pb:PS-CRDS-M} agrees with the force-balance residual obtained after substituting the penalty-induced constraint-force maps. This identifies the~\prettyref{pb:CITO-M} KKT system with the unit-objective force-parameterized certificate of~\prettyref{def:eps-stationary} and yields~\prettyref{thm:stationarity-transfer}.
The transfer rests on prescribing the MPCC constraint forces $\lambda_i^e,\lambda_i^i$ not as independent primal variables but as configuration-dependent maps $\lambda_i^e(\theta_i),\lambda_i^i(\theta_i)$ induced by the penalty gradient of $V(\cdot,\EpsFJ)$. This collapses the stationarity of~\prettyref{pb:CITO} onto the reduced trajectory space $(\StackConf,\StackCtrl)$ that~\prettyref{alg:CITO} already optimizes, where---on the region where the finite clamps $\LHamilton,\LPosition$ are inactive, so that $\bar U=\tilde U=U+V$---the per-step residual $\StatRes_i=\nabla_{\theta_i}\PerStepObj_i$ is identically the parameterized discrete force-balance residual of~\prettyref{eq:CRDS}. The $\EpsStop$-KKT residual blocks of~\prettyref{pb:CITO-M} then reduce line by line to an approximate stationarity condition of~\prettyref{pb:CITO}. We make the target condition the force-parameterized two-parameter $\TWO{\EpsStop}{\EpsFJAll}$-FJ-stationarity condition (\prettyref{def:eps-stationary})---force-parameterized because the constraint forces are prescribed as configuration-dependent maps rather than free variables: $\EpsStop$ governs the stationarity-residual and two-point-control tolerances inherited from~\prettyref{alg:CITO}, while $\EpsFJAll$ governs the physical feasibility and complementarity band. The two tolerances play different roles in the result below. Part~\ref{it:transfer-a} uses the aggregate feasibility band $\EpsFJAll=\max\{\EpsFJ,\EpsStop\}$ for $\EpsStop>0$ and the standing barrier parameter $\EpsFJ$. For the fixed-data uniform-boundedness family in part~\ref{it:transfer-b}, instead fix the barrier parameter $\EpsFJ$ and the remaining problem and algorithm data, including common clamps, and vary only $0<\EpsStop\leq\EpsFJ$; then $\EpsFJAll=\EpsFJ$ throughout the family.

\begin{remark}[Interpretation and scope of the force parameterization]
\label{rem:force-param-scope}
The certificate below is not the standard Fritz--John condition of the full MPCC~\prettyref{pb:CITO} with independent variables $(\StackConf,\StackCtrl,\lambda^e,\lambda^i)$. For each fixed barrier parameter $\EpsFJ$, the penalty construction selects configuration-dependent force maps $\Lambda_{\EpsFJ}(\StackConf)=\bigl(\lambda_i^e(\theta_i),\lambda_i^i(\theta_i)\bigr)_{i\in\mathcal{I}}$ satisfying the gradient-matching identity~\prettyref{eq:transfer-gradient-match}. \prettyref{def:eps-stationary} substitutes these maps into the force-balance, dual-feasibility, and relaxed-complementarity relations and then tests first-order stationarity only in the trajectory and control variables $(\StackConf,\StackCtrl)$; the derivatives of the force maps in~\prettyref{eq:eps-stationary-system-full} are the corresponding chain-rule terms. Thus the result certifies the physical residual, approximate feasibility, relaxed complementarity, and trajectory/control stationarity for this penalty-induced force selection. It does not certify stationarity with respect to independent variations of $\lambda^e$ or $\lambda^i$, and hence should not be identified with the standard C-, M-, S-, or KKT stationarity of the full MPCC. Promoting the forces to independent decision variables and proving one of those stronger notions is left to future work.
\end{remark}
\begin{definition}[smooth-damping penalty-induced force-parameterized $\TWO{\EpsStop}{\EpsFJAll}$-FJ certificate]
\label{def:eps-stationary}
Prescribe the current constraint forces of~\prettyref{pb:CITO} as configuration-dependent maps $\lambda_i^e(\theta_i)$ and $\lambda_i^i(\theta_i)\geq0$ obeying the gradient-matching identity, and use the predecessor evaluation $\lambda_{\gamma(i)}^i(\theta_{\gamma(i)})$ exclusively as the semi-implicit damping load. The current and predecessor evaluations have distinct roles and are not asserted to be equal. The current maps obey
\begin{ResizedEq}
\label{eq:transfer-gradient-match}
\FPP{h^e}{\theta_i}(\theta_i)^T\lambda_i^e(\theta_i)+\FPP{h^i}{\theta_i}(\theta_i)^T\lambda_i^i(\theta_i)=-\nabla_{\theta_i}V(\theta_i,\EpsFJ),
\end{ResizedEq}
and let the parameterized discrete force-balance residual of~\prettyref{eq:CRDS} be
\begin{ResizedEq}
\label{eq:transfer-residual}
F_i\triangleq&M\ddot{\theta}_i-f(\theta_i)-g-c(\theta_i,u_i)+\nabla_{\dot{\theta}_i}\Damp\!\left(\theta_{\gamma(i)},\dot{\theta}_i,\lambda_{\gamma(i)}^i(\theta_{\gamma(i)})\right)\\
&-\FPP{h^e}{\theta_i}(\theta_i)^T\lambda_i^e(\theta_i)-\FPP{h^i}{\theta_i}(\theta_i)^T\lambda_i^i(\theta_i),
\end{ResizedEq}
with $F\triangleq(F_i)_{i\in\mathcal{I}}$. A point $(\StackConf,\StackCtrl)$ satisfies the penalty-induced force-parameterized $\TWO{\EpsStop}{\EpsFJAll}$-FJ certificate for~\prettyref{pb:CITO} if there exist such force maps, physics multipliers $\mu_\StackConf$ (with timestep sub-blocks $\mu_{\theta_i}$), control multipliers $\mu_\StackCtrl$ (with sub-blocks $\mu_{u_i}$), free-sign equality-band multipliers $\mu_i^e$, nonnegative inequality multipliers $\mu_i^h,\mu_i^\lambda,\mu_i^c\geq0$, and, for every timestep $i$, an auxiliary control $u_i'$, such that
\begin{ResizedEq}
\label{eq:eps-stationary-system-full}
&\TWO{\EpsStop}{\EpsFJAll}\text{-FJ-stationarity}:\\
&\begin{cases}
\|F\|_1\leq\EpsStop\\
\begin{aligned}
&\Big\|\FPP{O}{\StackConf}+[\StatJac]^T\mu_\StackConf+\sum_i\FPP{h^e}{\theta_i}^T\mu_i^e-\sum_i\FPP{h^i}{\theta_i}^T\mu_i^h\\
&-\sum_i\FPP{\lambda_i^i}{\theta_i}^T\mu_i^\lambda-\sum_i\FPP{\big(\lambda_i^i\odot(\EpsFJAll\mathbf{1}-h^i)\big)}{\theta_i}^T\mu_i^c\Big\|\leq\EpsStop
\end{aligned}\\
\left\|\FPP{O}{\StackCtrl}+[\StatJacU]^T\mu_\StackConf+\mu_\StackCtrl\right\|\leq\EpsStop\\
u_i\in\mathcal{D},\ u_i'\in\mathcal{D},\ \|u_i-u_i'\|\leq\EpsStop,\ \mu_{u_i}\in N_{\mathcal{D}}(u_i')\\
|h^{e,j}(\theta_i)|\leq\EpsFJAll\quad\forall j\\
\begin{aligned}
&h^i(\theta_i)\geq0,\quad \lambda_i^i(\theta_i)\geq0,\\
&\lambda_i^i(\theta_i)\odot(\EpsFJAll\mathbf{1}-h^i(\theta_i))\geq0
\end{aligned}\\
\begin{aligned}
&\mu_i^h\odot h^i(\theta_i)=0,\quad\mu_i^\lambda\odot\lambda_i^i(\theta_i)=0,\\
&\mu_i^c\odot\big(\lambda_i^i(\theta_i)\odot(\EpsFJAll\mathbf{1}-h^i(\theta_i))\big)=0
\end{aligned}
\end{cases}
\end{ResizedEq}
\end{definition}
Here $\StatJac,\StatJacU$ are the Jacobians of~\prettyref{eq:sqp-definition}, coinciding with $\partial F/\partial\StackConf,\partial F/\partial\StackCtrl$ on the clamp-inactive region where $\StatRes\equiv F$, $\mu_{u_i}$ is the timestep-$i$ sub-block of $\mu_\StackCtrl$, and $\odot$ denotes the componentwise (Hadamard) product. The first condition is the residual of the exact dynamics~\prettyref{eq:CRDS} rather than the surrogate stationarity. The two Lagrangian-gradient blocks follow: the configuration block carries, besides the dynamics term $[\StatJac]^T\mu_\StackConf$, the free-sign equality-band multiplier $\mu_i^e$ and the three nonnegative multipliers $\mu_i^h,\mu_i^\lambda,\mu_i^c$ of the inequality conditions (entering with a minus sign as those conditions are $\geq0$ constraints), whereas the control block gains no such terms because neither $h^e$ nor $h^i,\lambda_i^i$ depend on $\StackCtrl$. Next is the two-point control condition. The remaining conditions add the physical feasibility that the surrogate cannot see: approximate equality feasibility (with a free-sign multiplier, hence no complementary slackness), the primal feasibility $h^i(\theta_i)\geq0$, dual feasibility $\lambda_i^i(\theta_i)\geq0$, and $\EpsFJAll$-relaxed complementarity $\lambda_i^i(\theta_i)\odot(\EpsFJAll\mathbf{1}-h^i(\theta_i))\geq0$ of~\prettyref{def:eps-FJ}, together with their complementary slackness. The terms involving $\partial\lambda_i^i(\theta_i)/\partial\theta_i$ are total-derivative terms created by substituting the prescribed force map into the relaxed dual-feasibility and complementarity relations; there is no separate stationarity equation with respect to an independent $\lambda_i^e$ or $\lambda_i^i$ variable, which is the precise difference between~\prettyref{def:eps-stationary} and a conventional first-order condition for the full MPCC~\prettyref{pb:CITO}. The objective coefficient is one because this reduced-variable certificate is inherited from the KKT condition of~\prettyref{pb:CITO-M}. Thus $\mu_0=1$ is a property of the transferred reduced-variable system, not evidence of KKT for the full MPCC: it is KKT-like only in the trajectory/control variables after the penalty-induced force maps have been substituted. No stationarity equation is certified in independent force directions, and no constraint qualification is established for~\prettyref{pb:CITO}. We retain the Fritz--John terminology to emphasize the absence of such a CQ and the possible unboundedness of the prescribed force family. The $\EpsStop$-KKT certificate of~\prettyref{thm:CITO-convergence} is instead a genuine KKT condition of the surrogate~\prettyref{pb:CITO-M}, whose LICQ is derived; its transfer to~\prettyref{pb:CITO} is what yields the FJ condition here. Sending both tolerances $\EpsStop,\EpsFJAll\to0$ collapses the two-point control, drives the parameterized force-balance residual to zero, and shrinks the penalty band onto the exact constraint manifold. This recovers exact satisfaction of the paper's force-parameterized first-order system for the selected penalty-induced force maps; it does not by itself establish C-, M-, S-, or KKT stationarity of the full independent-force MPCC.

\begin{reptheorem}{thm:stationarity-transfer}
\label{thm:stationarity-transfer-formal}
Let~\prettyref{pb:PS-CRDS-M} be a smooth full-rank CRDS whose objective $O$ satisfies~\prettyref{ass:objective}, with fixed startup data satisfying~\prettyref{ass:Slater}, including the fixed strict-interior node $\theta_0$ and its fixed pre-zero predecessor $\theta_{\gamma(0)}$, and the initialization scope of~\prettyref{rem:initialization-scope} with the clamps $\LHamilton,\LPosition$ fixed by the constant-ordering remark of~\prettyref{sec:Hamiltonian}, and assume the CRDS underlying~\prettyref{pb:CITO} is smooth (\prettyref{def:smooth-CRDS}), so that $\Damp_{\EpsFri,\LHamilton}^V=\Damp^V$. Let $(\StackConf,\StackCtrl)$ be the output of~\prettyref{alg:CITO}, and prescribe the constraint-force maps $\lambda_i^e(\theta_i),\lambda_i^i(\theta_i)$ from the penalty gradient of $V(\cdot,\EpsFJ)$. Set $\EpsFJAll\triangleq\max\{\EpsFJ,\EpsStop\}$. Then:
\begin{enumerate}[label=(\alph*),ref=(\alph*),leftmargin=*]
\item\label{it:transfer-a} The returned point $(\StackConf,\StackCtrl)$ satisfies the force-parameterized $\TWO{\EpsStop}{\EpsFJAll}$-FJ-stationarity condition of~\prettyref{pb:CITO} (\prettyref{def:eps-stationary}). This transfer holds for every $\EpsStop>0$ at the standing barrier parameter $\EpsFJ$.
\item\label{it:transfer-b} Fix $\EpsFJ$ and hold the remaining problem and algorithm data, including the clamps $\LHamilton,\LPosition$, common. Uniformly over $0<\EpsStop\leq\EpsFJ$, there are finite constants independent of $\EpsStop$ bounding the prescribed forces $\|\lambda_i^e(\theta_i)\|,\|\lambda_i^i(\theta_i)\|\leq\bar\Lambda$, the physics multipliers $\|\mu_{\theta_i}\|$ and the penalty weight $\varrho$ (both bounded by finite constants given in the proof), and the trajectory $\|\theta_i\|\leq\LPosBound(H)$, $u_i\in\mathcal{D}$.
\end{enumerate}
\end{reptheorem}
This formal version is proved in~\prettyref{appen:stationarity-proofs}. Part~\ref{it:transfer-a} reads the $\EpsStop$-KKT certificate of~\prettyref{pb:CITO-M} through the residual identity $\StatRes_i\equiv F_i$ on the clamp-inactive region after restricting the original constraint-force variables to the selected penalty-induced map. The assertion is thus relative to the penalty-induced force map selected at the fixed barrier parameter $\EpsFJ$; no stationarity conclusion is made for independent perturbations of the force variables of~\prettyref{pb:CITO}. For the fixed-$\EpsFJ$ family of part~\ref{it:transfer-b}, the Hamiltonian majorant at the endpoint $\EpsStop=\EpsFJ$ supplies common clamps, after which~\prettyref{lem:trajectory-bound}, together with the step- and multiplier-stability bounds of~\prettyref{appen:global-convergence-proofs}, and the fixed-$\EpsFJ$ level set supply the stated $\EpsStop$-independent bounds.

\subsubsection{Nonsmooth Damping-Graph Transfer}
We next remove the smooth-damping assumption. For each fixed smoothing,~\prettyref{alg:CITO} supplies the smooth $\EpsStop$-KKT certificate used above; we normalize its complete dual tuple and replace the smooth damping derivatives by nearby Clarke normals to the lifted nonsmooth damping graphs. This yields the homogeneous force-parameterized certificate of~\prettyref{def:eps-stationary-NS} and~\prettyref{thm:nonsmooth-transfer}, whose objective multiplier may vanish in the smoothing limit. To remove the smoothness assumption one replaces the genuinely nonsmooth frictional damping $\Damp^V$ by a smoothed $\Damp_{\EpsFri,\LHamilton}^V$ (\prettyref{def:S-D}) and sends the smoothing scale $\EpsFri\to0$. A unit-objective transfer in the manner of~\prettyref{thm:stationarity-transfer}\ref{it:transfer-b} would need the physics multiplier $\|\mu_{\theta_i}\|$ to stay bounded as $\EpsFri\to0$; but the physics-multiplier bound (\prettyref{appen:global-convergence-proofs}) is inversely proportional to the assembled LICQ modulus, which degrades as $\EpsFri\to0$ because the earlier-configuration (off-diagonal) dissipation coupling of the stationarity Jacobian blows up in operator norm. Rather than forcing a smoothing-uniform physics-multiplier bound, we retain a homogeneous dual certificate and transfer its smooth damping derivatives to Clarke normals of the nonsmooth damping graphs (\prettyref{def:eps-stationary-NS}). The complete tuple is normalized rather than fixing the objective coefficient, so the limiting multiplier $\mu_0$ may vanish. This construction requires no CQ for~\prettyref{pb:CITO}, but it remains force-parameterized and is not standard MPCC C-stationarity. For each fixed smoothing the smooth $\EpsStop$-KKT certificate is finite; normalizing the complete dual certificate---including its per-term smooth graph normals---by its overall scale (the Euclidean norm of the complete tuple; an overbar denotes division by this scale, so the normalized tuple lies on the unit sphere) places the objective coefficient $\mu_0\in(0,1]$ and automatically yields the normalized bounds $\|\bar\mu_{\theta_i}\|\leq1$ and $\bar\mu_{\theta_i}^T K_i\bar\mu_{\theta_i}\leq|\mathcal{T}|/\delta_{\min}$, which are exactly the fixed constants that the normal-compatibility condition consumes. The smoothing therefore needs to satisfy only~\prettyref{def:S-D-compatible}, as we state next and prove in~\prettyref{appen:stationarity-proofs}.
\begin{definition}
\label{def:S-D-compatible}
Assume the damping energy decomposes as a finite sum $\Damp^V=\sum_{m\in\mathcal{T}}\Damp^{V,m}$ of convex, possibly nonsmooth terms, smoothed termwise so that $\Damp_{\EpsFri,\LHamilton}^V=\sum_{m\in\mathcal{T}}\Damp_{\EpsFri,\LHamilton}^{V,m}$ with each $\Damp_{\EpsFri,\LHamilton}^{V,m}$ a~\prettyref{def:S-D} smoothing of $\Damp^{V,m}$ and the cardinality $|\mathcal{T}|$ independent of $\EpsFri,\EpsStop$. With the earlier-configuration mixed-derivative block $\Gamma_i$ and the diagonal damping block $K_i$,
\begin{ResizedEq}
\label{eq:S-D-compatible-blocks}
\Gamma_i\triangleq&\FPPTT{\Damp_{\EpsFri,\LHamilton}^V(\theta_{\gamma(i)},\dot{\theta}_i)}{\dot{\theta}_i}{\theta_{\gamma(i)}}\\
K_i\triangleq&\delta_i\FPPT{\Damp_{\EpsFri,\LHamilton}^V(\theta_{\gamma(i)},\dot{\theta}_i)}{\theta_i}\succeq0,
\end{ResizedEq}
the per-term blocks $\Gamma_i^m,K_i^m$ defined by the same second derivatives of $\Damp_{\EpsFri,\LHamilton}^{V,m}$ ($\dot{\theta}_i$ held fixed in $\Gamma_i$; so $\Gamma_i=\sum_{m\in\mathcal{T}}\Gamma_i^m$ and $K_i=\sum_{m\in\mathcal{T}}K_i^m$ with each $K_i^m\succeq0$, whence $\xi^T K_i^m\xi\leq\xi^T K_i\xi$), and the fixed, $\EpsFri$-independent energy constant
\begin{ResizedEq}
C_E^\star\triangleq\frac{|\mathcal{T}|}{\delta_{\min}},
\end{ResizedEq}
$\delta_{\min}>0$ the a priori timestep floor of~\prettyref{eq:cito-delta-min}, a smoothed dissipation $\Damp_{\EpsFri,\LHamilton}^V$ satisfying~\prettyref{def:S-D} is normal-compatible if there is a modulus $\omega_N(\EpsFri)\downarrow0$ as $\EpsFri\downarrow0$, uniform over $i\in\mathcal{I}$, over $m\in\mathcal{T}$, and over every previous configuration with $V(\theta_{\gamma(i)},\EpsFJ)\leq\LHamilton<\infty$, such that whenever the tested direction obeys the normalized bounds $\|\xi\|\leq1$ and $\xi^T K_i\xi\leq C_E^\star$, each term $m\in\mathcal{T}$ admits, at its~\prettyref{def:S-D} nearby point $(\theta_{\gamma(i)},\dot{\theta}_i^{m\dagger},q_i^{m\dagger})$ with $q_i^{m\dagger}\in\partial_{\dot{\theta}_i}\Damp^{V,m}(\theta_{\gamma(i)},\dot{\theta}_i^{m\dagger})$, a Clarke normal
\begin{ResizedEq}
&n_i^{m\sharp}\in N^C_{\gph\partial_{\dot{\theta}}\Damp^{V,m}}(\theta_{\gamma(i)},\dot{\theta}_i^{m\dagger},q_i^{m\dagger})\\
&\|n_i^{m\sharp}-((\Gamma_i^m)^T\xi,\delta_i K_i^m\xi,-\xi)\|\leq\omega_N(\EpsFri).
\end{ResizedEq}
Here $\gph$ denotes the graph of a set-valued map, and $N^C_\cdot(z)$ the Clarke normal cone to a closed set $\cdot$ at $z\in S$~\cite{clarke1983optimization}.
\end{definition}
To state the nonsmooth transfer, lift each damping term's force to a free variable: writing $q_i^m$ for a stand-in of the velocity gradient $\nabla_{\dot{\theta}_i}\Damp^{V,m}$ (a lifted damping force, not a generalized coordinate) and $q_i\triangleq\sum_{m\in\mathcal{T}}q_i^m$, replace the aggregate gradient $\nabla_{\dot{\theta}_i}\Damp^V=\sum_{m\in\mathcal{T}}\nabla_{\dot{\theta}_i}\Damp^{V,m}$ in~\prettyref{eq:transfer-residual} by $\sum_{m\in\mathcal{T}}q_i^m$,
\begin{ResizedEq}
\label{eq:transfer-residual-lifted}
&\widehat F_i(\StackConf,\StackCtrl,(q_i^m)_{m\in\mathcal{T}})\\
\triangleq&M\ddot{\theta}_i-f(\theta_i)-g-c(\theta_i,u_i)+\sum_{m\in\mathcal{T}}q_i^m\\
&-\FPP{h^e}{\theta_i}(\theta_i)^T\lambda_i^e(\theta_i)-\FPP{h^i}{\theta_i}(\theta_i)^T\lambda_i^i(\theta_i),
\end{ResizedEq}
with $\widehat F\triangleq(\widehat F_i)_{i\in\mathcal{I}}$, and introduce the graph-argument map $Q_i$ and its pullback
\begin{ResizedEq}
\label{eq:transfer-Qmap}
Q_i(\StackConf,q_i)\triangleq&\THREE{\theta_{\gamma(i)}}{\frac{\theta_i-\theta_{\gamma(i)}}{\delta_i}}{q_i}\\
[\nabla Q_i]^T(a_i,b_i,c_i)=&\THREE{\tfrac{b_i}{\delta_i}}{a_i-\tfrac{b_i}{\delta_i}}{c_i},
\end{ResizedEq}
the pullback acting on the $(\theta_i,\theta_{\gamma(i)},q_i)$ blocks for a covector $(a_i,b_i,c_i)$ in the $(\theta_{\gamma(i)},\dot{\theta}_i,q_i)$ graph coordinates; the same map $Q_i$ is applied per term, to $(\StackConf,q_i^m)$, with $[\nabla Q_i]^T$ unchanged.
\begin{definition}[penalty-induced force-parameterized homogeneous Clarke--Fritz--John $\TWO{\EpsStop}{\EpsFJAll}$ certificate]
\label{def:eps-stationary-NS}
Prescribe the constraint-force maps as the free-sign equality-force map $\lambda_i^e(\theta_i)$ and the nonnegative inequality-force map $\lambda_i^i(\theta_i)\geq0$ as in~\prettyref{def:eps-stationary}, and let $\widehat F,Q_i$ be as in~\prettyref{eq:transfer-residual-lifted} and~\ref{eq:transfer-Qmap}. A point $(\StackConf,\StackCtrl)$ satisfies the penalty-induced force-parameterized homogeneous Clarke (C-)$\TWO{\EpsStop}{\EpsFJAll}$-FJ certificate for~\prettyref{pb:CITO} if there exist an objective multiplier $\mu_0\geq0$, physics multipliers $\mu_\StackConf$ (with timestep sub-blocks $\mu_{\theta_i}$), control multipliers $\mu_\StackCtrl$ (with sub-blocks $\mu_{u_i}$), free-sign equality-band multipliers $\mu_i^e$, nonnegative inequality multipliers $\mu_i^h,\mu_i^\lambda,\mu_i^c\geq0$, auxiliary controls $u_i'$, and, for every timestep $i$ and every term $m\in\mathcal{T}$, a nearby velocity $\dot{\theta}_i^{m\sharp}$, a damping force $q_i^{m\sharp}\in\partial_{\dot{\theta}_i}\Damp^{V,m}(\theta_{\gamma(i)},\dot{\theta}_i^{m\sharp})$, and a Clarke graph normal $n_i^{m\sharp}\in N^C_{\gph\partial_{\dot{\theta}}\Damp^{V,m}}(\theta_{\gamma(i)},\dot{\theta}_i^{m\sharp},q_i^{m\sharp})$, such that
\begin{ResizedEq}
&\TWO{\EpsStop}{\EpsFJAll}\text{ homogeneous Clarke--FJ certificate}:\\
&\begin{cases}
\|\widehat F(\StackConf,\StackCtrl,(q^{m\sharp}))\|_1\leq\EpsStop,\quad\|\dot{\theta}_i^{m\sharp}-\dot{\theta}_i\|\leq\EpsStop\ \ \forall m\\
\begin{aligned}
&\Big\|\nabla_{(\StackConf,q)}\Big[\mu_0 O+\sum_i\langle\mu_{\theta_i},\widehat F_i\rangle+\sum_i\langle\mu_i^e,h^e(\theta_i)\rangle-\sum_i\langle\mu_i^h,h^i(\theta_i)\rangle\\
&-\sum_i\langle\mu_i^\lambda,\lambda_i^i(\theta_i)\rangle-\sum_i\langle\mu_i^c,\lambda_i^i(\theta_i)\odot(\EpsFJAll\mathbf{1}-h^i(\theta_i))\rangle\Big]\\
&+\sum_i\sum_{m\in\mathcal{T}}[\nabla Q_i]^T n_i^{m\sharp}\Big\|\leq\EpsStop
\end{aligned}\\
\left\|\mu_0\FPP{O}{\StackCtrl}+[\StatJacU]^T\mu_\StackConf+\mu_\StackCtrl\right\|\leq\EpsStop\\
u_i\in\mathcal{D},\ u_i'\in\mathcal{D},\ \|u_i-u_i'\|\leq\EpsStop,\ \mu_{u_i}\in N_{\mathcal{D}}(u_i')\\
|h^{e,j}(\theta_i)|\leq\EpsFJAll\quad\forall j\\
\begin{aligned}
&h^i(\theta_i)\geq0,\quad \lambda_i^i(\theta_i)\geq0,\\
&\lambda_i^i(\theta_i)\odot(\EpsFJAll\mathbf{1}-h^i(\theta_i))\geq0
\end{aligned}\\
\begin{aligned}
&\mu_i^h\odot h^i(\theta_i)=0,\quad\mu_i^\lambda\odot\lambda_i^i(\theta_i)=0,\\
&\mu_i^c\odot\big(\lambda_i^i(\theta_i)\odot(\EpsFJAll\mathbf{1}-h^i(\theta_i))\big)=0
\end{aligned}
\end{cases}
\end{ResizedEq}
together with the dual-nontriviality normalization
\begin{ResizedEq}
\label{eq:eps-stationary-NS-normalization}
1=&\mu_0^2+\|\mu_\StackConf\|^2+\|\mu_\StackCtrl\|^2+\\
&\sum_i\left(\|\mu_i^e\|^2+\|\mu_i^h\|^2+\|\mu_i^\lambda\|^2+\|\mu_i^c\|^2\right)+\\
&\sum_i\sum_{m\in\mathcal{T}}\|n_i^{m\sharp}\|^2.
\end{ResizedEq}
\end{definition}
Here $\nabla_{(\StackConf,q)}$ is evaluated at $q=(q^{m\sharp})$, and each $q_i^m$-block of the bracket is $\mu_{\theta_i}$ because $\FPP{\widehat F_i}{q_i^m}=I$ and $\mu_0 O$ is $q$-independent. The equality-band and inequality multipliers $\mu_i^e,\mu_i^h,\mu_i^\lambda,\mu_i^c$ enter the configuration block and satisfy dual feasibility and complementary slackness exactly as in~\prettyref{def:eps-stationary}; the objective multiplier $\mu_0$ and the per-term Clarke graph normals $n_i^{m\sharp}$ of the nonsmooth damping subdifferentials are the genuinely new ingredients. When the damping is smooth, each $\partial_{\dot{\theta}}\Damp^{V,m}$ is single-valued, $\dot{\theta}_i^{m\sharp}=\dot{\theta}_i$, $q_i^{m\sharp}=\nabla_{\dot{\theta}_i}\Damp^{V,m}$, and the per-term graph normals $n_i^{m\sharp}=((\Gamma_i^m)^T\mu_{\theta_i},\delta_i K_i^m\mu_{\theta_i},-\mu_{\theta_i})$ make each $q_i^m$-block vanish and sum over $m$ (using $\sum_{m}\Gamma_i^m=\Gamma_i$, $\sum_{m}K_i^m=K_i$) to reproduce the smooth graph normal $(\Gamma_i^T\mu_{\theta_i},\delta_i K_i\mu_{\theta_i},\cdot)$, reducing the line to $\|\mu_0\FPP{O}{\StackConf}+[\StatJac]^T\mu_\StackConf\|\leq\EpsStop$; in the nonabnormal case $\mu_0>0$, dividing through by $\mu_0$ recovers the configuration block of~\prettyref{def:eps-stationary}, so the whole condition collapses to~\prettyref{def:eps-stationary} up to the positive rescaling by $\mu_0$. This is the paper's homogeneous Clarke--Fritz--John certificate for the penalty-induced force-parameterized formulation. The Clarke normals concern the lifted nonsmooth damping graphs. Because the equality and inequality/contact forces remain prescribed maps, the conclusion is not standard MPCC C-stationarity of~\prettyref{pb:CITO} with independent force variables. As in the FJ-homogeneity remark following~\prettyref{def:eps-FJ}, the whole multiplier tuple $(\mu_0,\mu_\StackConf,\mu_\StackCtrl,\mu_i^e,\mu_i^h,\mu_i^\lambda,\mu_i^c,(n_i^{m\sharp}))$ is determined only up to positive scaling and is fixed here by the unit-sphere normalization~\prettyref{eq:eps-stationary-NS-normalization} rather than by $\mu_0=1$; the objective coefficient $\mu_0$ may vanish in the smoothing limit, so unlike~\prettyref{def:eps-stationary} the condition is not asserted in unit-objective form. When $\mu_0>0$, the homogeneous tuple can be divided by $\mu_0$ to obtain a unit-objective form. When $\mu_0=0$, the certificate is abnormal: its nontriviality is carried entirely by the constraint, control-set, complementarity-band, and damping-graph multipliers, and it gives no first-order information about descent of the user objective. The Clarke cone $N^C$ additionally relaxes the smooth damping-graph stationarity of~\prettyref{def:eps-stationary} to its nonsmooth analog. Replacing only these damping-graph Clarke normals by limiting normals would sharpen the damping component of the force-parameterized certificate, but would not establish standard MPCC M-stationarity without independent force variables and the appropriate normal cone to the complementarity set.
\begin{reptheorem}{thm:nonsmooth-transfer}
\label{thm:nonsmooth-transfer-formal}
Take the hypotheses of~\prettyref{thm:CITO-convergence} (including the fixed strict-interior node $\theta_0$ of~\prettyref{ass:Slater} and its fixed pre-zero predecessor $\theta_{\gamma(0)}$), fix $\EpsFJ$, and use a smoothed dissipation $\Damp_{\EpsFri,\LHamilton}^V$ that is normal-compatible (\prettyref{def:S-D-compatible}), let $(\StackConf,\StackCtrl)$ be the output of~\prettyref{alg:CITO}, prescribe the force-parameterized maps $\lambda_i^e(\theta_i),\lambda_i^i(\theta_i)$, and set $\EpsFJAll\triangleq\max\{\EpsFJ,\EpsStop\}$. Write $N^\star\triangleq|\mathcal{T}|(|G_0|+N_{\mathrm{sub}})$ for the fixed, $\EpsFri,\EpsStop$-independent number of timestep-term pairs, and assume $\EpsFri$ is small enough that $\sqrt{N^\star}\,\omega_N(\EpsFri)\leq1/2$. Then:
\begin{enumerate}[label=(\alph*),ref=(\alph*),leftmargin=*]
\item\label{it:ns-transfer-a} The returned point satisfies the penalty-induced force-parameterized homogeneous Clarke--Fritz--John $\TWO{\EpsStopNS}{\EpsFJAll}$ certificate for~\prettyref{pb:CITO} (\prettyref{def:eps-stationary-NS}), with
\begin{ResizedEq}
\label{eq:ns-transfer-band}
\EpsStopNS\triangleq&\max\big\{2\EpsStop+2C_N\sqrt{N^\star}\,\omega_N(\EpsFri),\ \EpsStop+C_F N^\star\EpsFri\big\},
\end{ResizedEq}
for constants $C_N,C_F$ independent of $\EpsFri,\EpsStop$.
\item\label{it:ns-transfer-b} Fix $\EpsFJ$ and $\EpsFri>0$, and hold the remaining problem and algorithm data, including the clamps $\LHamilton,\LPosition$, common. Uniformly over $0<\EpsStop\leq\EpsFJ$, there are finite constants independent of $\EpsStop$ bounding the trajectory ($\|\theta_i\|\leq\LPosBound(H)$, $u_i\in\mathcal{D}$), the prescribed forces ($\|\lambda_i^e(\theta_i)\|,\|\lambda_i^i(\theta_i)\|\leq\bar\Lambda$), and the penalty weight ($\varrho$, bounded by a finite constant given in the proof).
\end{enumerate}
\end{reptheorem}
This final version is proved in~\prettyref{appen:stationarity-proofs}. Part~\ref{it:ns-transfer-a} normalizes the fixed-smoothing $\EpsStop$-KKT certificate of~\prettyref{pb:CITO-M} by its complete-certificate scale: this places the objective coefficient $\mu_0\in(0,1]$ and, because the normalized per-term smooth graph normals $((\Gamma_i^m)^T\bar\mu_{\theta_i},\delta_i K_i^m\bar\mu_{\theta_i},-\bar\mu_{\theta_i})$ are part of the unit-sphere certificate, automatically yields $\|\bar\mu_{\theta_i}\|\leq1$ and $\bar\mu_{\theta_i}^T K_i\bar\mu_{\theta_i}\leq C_E^\star=|\mathcal{T}|/\delta_{\min}$; these are exactly the fixed constants that~\prettyref{def:S-D-compatible} consumes term by term, replacing each normalized smooth graph normal by a genuine Clarke normal at its~\prettyref{def:S-D} nearby graph point, and a final conic rescaling $\beta\leq2$ restores exact unit-sphere nontriviality; the dynamics line alone obeys the sharper band $\EpsStop+C_F N^\star\EpsFri$. Part~\ref{it:ns-transfer-b} fixes both $\EpsFJ$ and $\EpsFri$, so the smoothed problem and its Jacobian/SQP constants are common while only $\EpsStop\in(0,\EpsFJ]$ varies; it then collects the same $\EpsStop$-independent primal bounds as the smooth transfer. As in~\prettyref{thm:stationarity-transfer}, the conclusion is relative to the prescribed penalty-induced force maps: the nonsmooth graph lift changes the damping part of the certificate but does not promote the constraint forces to independent optimization variables.
\begin{remark}
\label{rem:ns-transfer-concrete}
\prettyref{thm:nonsmooth-transfer} is a general-CRDS statement, holding for any smoothing that satisfies the single normal-compatibility condition~\prettyref{def:S-D-compatible}. Verifying~\prettyref{def:S-D-compatible} for a concrete smoothing of the frictional-contact damping~\prettyref{eq:damping-ABD}---through the three-region (sticking/transition/sliding) geometry of a weighted-norm profile, a small-weight split accommodating vanishing penalty-induced contact forces, and a vertexwise Schur reduction of the plane-velocity block---is carried out unconditionally for the articulated-body model in the specialization below (\prettyref{lem:abd-normal-compatible}), so for that model the nonsmooth transfer is gap-free.
\end{remark}

\subsubsection{Articulated Body Specialization}
We now instantiate the general theory on articulated-body dynamics and discharge, for that model, the concrete smoothing left open by~\prettyref{rem:ns-transfer-concrete}. The frictional-contact damping~\prettyref{eq:damping-ABD} is a penalty-induced dissipation $\Damp^V$ whose friction part decomposes into a finite set $\mathcal{T}$ of per-vertex--plane terms $\eta[\Lambda^i]^m\|v_\parallel^m\|$---one for each vertex-plane friction incidence over the convex-hull pairs $jh,j'h'$, indexed by $m\in\mathcal{T}$ exactly as required by the finite-sum hypothesis of~\prettyref{def:S-D-compatible}---plus the smooth plane-velocity damper $\EpsDamp\|\varpi\|^2$ in the separating-plane velocities $\varpi$. The concrete smoothing, built in~\prettyref{appen:penalty-proofs}, is a pseudo-Huber replacement of each Euclidean friction-velocity norm $\|v_\parallel^m\|$ at a smoothing scale controlled by $\EpsFri$, retaining the plane-velocity damper unchanged. This same smoothing was already shown $\Damp_R^V$-close in the sense of~\prettyref{def:S-D} by~\prettyref{lem:smooth-DRp-Regularity}; we now show it is moreover normal-compatible.
\begin{lemma}[(informal)]
\label{lem:abd-compatible}
For each fixed $\EpsFri>0$, $0<\LHamilton<\infty$, and prescribed a priori timestep floor $\delta_{\min}>0$, there exists a smoothing of the reduced articulated-body damping in the sense of~\prettyref{def:S-D} that is moreover Schur-lifted per-vertex normal-compatible in the damping-active physical-plus-plane velocity $\TWO{\dot\theta}{\varpi}$, with the plane-velocity damper $\EpsDamp\|\varpi\|^2$ retained exactly and externally.
\end{lemma}
The formal restatement and its proof---the trajectory-localized per-vertex normal compatibility through the three-region sticking/transition/sliding analysis---are recorded in~\prettyref{appen:abd-specialization}. The modulus is uniform on the compact configuration and bounded full-velocity region containing every returned CITO trajectory, rather than on the generally unbounded bare reduced-energy sublevel; this introduces no new physical assumption because~\prettyref{thm:abd-C-stationarity} concerns only those returned trajectories. The ABD case is handled separately because the physically meaningful conclusion is per contact vertex in the unreduced physical-plus-plane velocities, whereas the reduced dissipation is a partial minimization over the plane velocities that does not decompose vertexwise; the proof therefore performs its own reduced-to-full certificate lift, with the plane-velocity damper kept exact and external.
\begin{theorem}[articulated-body per-vertex Clarke--Fritz--John certificate]
\label{thm:abd-C-stationarity}
Under the standing~\prettyref{ass:abd-mass}, take all hypotheses, fixed startup and initialization data, clamp choices, and algorithm-input conditions of~\prettyref{thm:CITO-convergence}, and fix $\EpsFJ$ and $\EpsFri>0$. For the finite initial mesh, form the smoothing-independent thresholds $\delta^\star$ and $\delta_H$ of~\prettyref{thm:CITO-convergence}; set prospectively $\delta_{\min}\triangleq\min\{\min_{i\in G_0}\delta_i^0,\min(\delta^\star,\delta_H)/2\}>0$. Apply~\prettyref{lem:smooth-DRp-Regularity} with this prescribed floor to construct the reduced articulated-body damping, and run~\prettyref{alg:CITO} on~\prettyref{pb:PSR-CRDS} with the normal-compatible pseudo-Huber smoothing constructed in~\prettyref{appen:abd-specialization}, prescribe the force-parameterized maps $\lambda_i^e(\theta_i),\lambda_i^i(\theta_i)$, and set $\EpsFJAll\triangleq\max\{\EpsFJ,\EpsStop\}$. Writing $N^\star\triangleq|\mathcal{T}|(|G_0|+N_{\mathrm{sub}})$ and assuming $\EpsFri$ small enough that $\sqrt{N^\star}\,\omega_N(\EpsFri)\leq1/2$, the returned physical trajectory admits attained separating planes $p_i^\star$, unique minimizing plane velocities $\varpi_i^\star$, and auxiliary coordinates reconstructed forward from their fixed node-$0$ base values along the predecessor chain, yielding temporally consistent $(\varsigma_i^\star,\omega_i^\star)$ and full configurations $\vartheta_i^\star$ for every $i\in\mathcal I$. The stacked extended trajectory $(\vartheta_i^\star)_{i\in\mathcal I}$ satisfies auxiliary feasibility and exact plane-configuration stationarity, and satisfies the full physical-plus-plane, penalty-induced force-parameterized homogeneous C-$\TWO{\EpsStopNS}{\EpsFJAll}$-FJ-stationarity condition on the extended formulation~\prettyref{pb:CRDS-ABD-extended}---in the damping-active physical-plus-plane velocity $\TWO{\dot\theta}{\varpi}$, each vertex term $m\in\mathcal{T}$ represented by an exact Clarke normal to the corresponding damping-active subgradient graph and the plane-velocity damper $\EpsDamp\|\varpi\|^2$ kept exact and external---with
\begin{ResizedEq}
\label{eq:abd-tildeeps}
\EpsStopNS\triangleq&\max\big\{2\EpsStop+2C_N\sqrt{N^\star}\,\omega_N(\EpsFri),\ \EpsStop+C_F N^\star\EpsFri\big\},
\end{ResizedEq}
for constants $C_N,C_F$ independent of $\EpsFri,\EpsStop$.
\end{theorem}
This theorem is proved in~\prettyref{appen:abd-specialization}. For articulated-body dynamics it discharges the closeness~\prettyref{lem:smooth-DRp-Regularity} and the per-vertex normal-compatibility condition on the compact region containing all returned CITO configurations and reconstructed full velocities, then lifts the reduced certificate to the full physical-plus-plane formulation and normalizes the complete full certificate there. The reconstructed separating planes and plane velocities are exact mathematical minimizers of the partial minimizations in~\prettyref{eq:reduced-penalty}; plane-velocity uniqueness follows from~\prettyref{lem:varpi-reconstruction} and the retained $\EpsDamp\|\varpi\|^2$ term. Consequently, exact auxiliary stationarity, the envelope identities, and the Schur lift introduce no additional theoretical tolerance. Finite-precision inner-solve error is outside the present theorem and is not implicitly absorbed into $\EpsStop$ or $\EpsStopNS$. The reconstructed variables enlarge the damping-active physical system, but the equality and inequality/contact forces remain the prescribed penalty-induced maps.

\subsection{Hierarchy and Interpretation of the Certificates}
The paper uses several related but nonidentical first-order certificates. They differ along two axes: the optimization problem and variable space on which they are stated, and the smooth or nonsmooth geometry used in the transfer. The direct algorithmic result is an $\EpsStop$-KKT point of the smooth modified~\prettyref{pb:CITO-M}. The per-step condition of~\prettyref{def:eps-FJ} records approximate physical feasibility and does not include stationarity of the user trajectory objective. \prettyref{thm:stationarity-transfer} transfers the KKT system of~\prettyref{pb:CITO-M} to a unit-objective, penalty-induced force-parameterized certificate for~\prettyref{pb:CITO}. \prettyref{thm:nonsmooth-transfer} replaces smooth damping derivatives by Clarke normals to lifted damping graphs and normalizes the complete dual tuple, so its objective multiplier $\mu_0$ may vanish. The adjective ``Clarke'' in this last statement refers to the damping-graph normal; it is not standard MPCC C-stationarity of~\prettyref{pb:CITO}. Standard MPCC C-, M-, S-, and KKT stationarity would require independent force variables and is not established here. \prettyref{table:certificate-scope} summarizes these distinctions.

\begin{table*}[t]
\centering
\scriptsize
\setlength{\tabcolsep}{3.5pt}
\renewcommand{\arraystretch}{1.18}
\TableRowColors
\begin{tabular}{>{\raggedright\arraybackslash}m{0.13\textwidth}>{\raggedright\arraybackslash}m{0.15\textwidth}>{\raggedright\arraybackslash}m{0.18\textwidth}>{\raggedright\arraybackslash}m{0.17\textwidth}>{\raggedright\arraybackslash}m{0.14\textwidth}>{\raggedright\arraybackslash}m{0.15\textwidth}}
\toprule
\TableHeaderRow\textbf{Result} & \textbf{Problem / variables} & \textbf{Force treatment} & \textbf{First-order notion} & \textbf{Objective multiplier} & \textbf{Interpretation and limitation}\\
\midrule
\prettyref{def:eps-FJ}; \prettyref{lem:eps-FJ},~\ref{lem:eps-S-FJ} & Per-step constrained CRDS (\prettyref{pb:CRDS-ABD}/\ref{pb:PS-CRDS}) & Existential or penalty-induced forces used to recover physical relations & $\EpsFJ$-FJ physical condition & No trajectory-objective multiplier & Approximate force balance, feasibility, and complementarity at one step; not stationarity of the user trajectory objective.\\
\prettyref{def:eps-KKT}; \prettyref{thm:sqp-convergence},~\ref{thm:CITO-convergence} & Modified smooth trajectory NLP,~\prettyref{pb:CITO-M}, in $(\StackConf,\StackCtrl)$ & No independent constraint-force variables; forces are absorbed in the surrogate & Genuine $\EpsStop$-KKT of~\prettyref{pb:CITO-M} & Normalized to one & Direct algorithmic convergence result; LICQ is derived for~\prettyref{pb:CITO-M}, not for the original MPCC.\\
\prettyref{def:eps-stationary}; \prettyref{thm:stationarity-transfer} & \prettyref{pb:CITO} restricted to $(\lambda^e,\lambda^i)=\Lambda_{\EpsFJ}(\StackConf)$ & Penalty-induced configuration-dependent force maps & Unit-objective force-parameterized FJ certificate & $\mu_0=1$ & KKT-like only in the reduced trajectory/control variables; not KKT or standard MPCC stationarity of full~\prettyref{pb:CITO}.\\
\prettyref{def:eps-stationary-NS}; \prettyref{thm:nonsmooth-transfer},~\ref{thm:abd-C-stationarity} & Same force-restricted~\prettyref{pb:CITO}, with lifted nonsmooth damping variables & Same prescribed force maps; nearby damping forces and graph normals are lifted & Homogeneous force-parameterized Clarke--FJ certificate & $\mu_0\geq0$ and may vanish & ``Clarke'' refers to damping-graph normals. If $\mu_0=0$, the certificate is abnormal and contains no objective first-order information.\\
Standard MPCC C/M/S/KKT (not proved) & Full~\prettyref{pb:CITO} in $(\StackConf,\StackCtrl,\lambda^e,\lambda^i)$ & Constraint forces vary independently & Standard MPCC C-, M-, S-, or KKT stationarity & Usually unit-objective under a suitable CQ; homogeneous FJ variants otherwise & Would test independent force directions and complementarity geometry; outside the present result.\\
\bottomrule
\hiderowcolors
\end{tabular}
\caption{Hierarchy and scope of the first-order certificates. The rows are not a single implication chain because they concern different problems, variable spaces, and smoothness settings. The direct solver result is KKT for~\prettyref{pb:CITO-M}; the conclusions for~\prettyref{pb:CITO} are transferred, penalty-induced force-parameterized certificates.}
\label{table:certificate-scope}
\end{table*}

Across the two transferred certificates, the paper certifies the force-balance residual for the selected penalty-induced forces, approximate equality feasibility, nonnegative selected inequality/contact forces, relaxed complementarity, and trajectory/control adjoint conditions including the force-map chain rules. It does not certify stationarity under independent force perturbations, KKT of the original MPCC, standard MPCC C-, M-, or S-stationarity, independent optimality of the contact-force variables, or a constraint qualification for~\prettyref{pb:CITO}.

\begin{figure}[t]
\centering
\begin{tikzpicture}
\node[anchor=north west,inner sep=0] at (0.10\linewidth,0)
  {\includegraphics[width=0.90\linewidth]{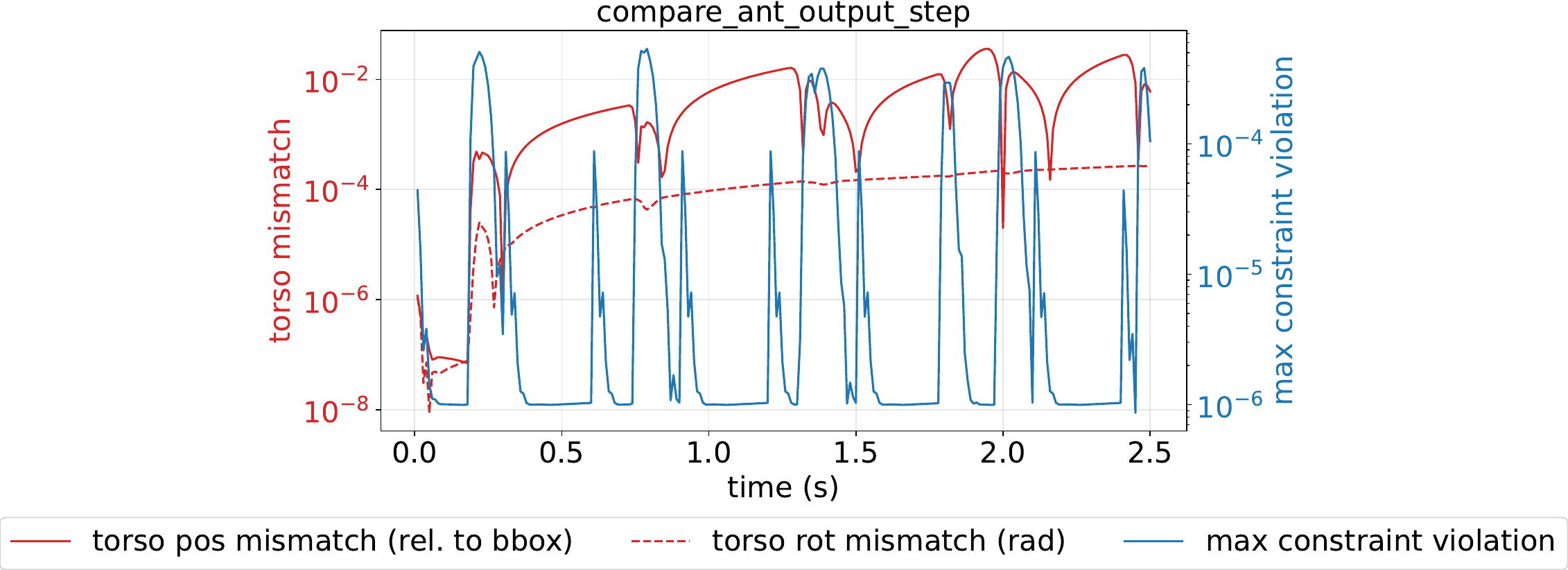}};
\node[anchor=west,inner sep=0] at (0,-0.148\linewidth)
  {\includegraphics[trim={560 172 611 250},clip,frame,width=0.24\linewidth]{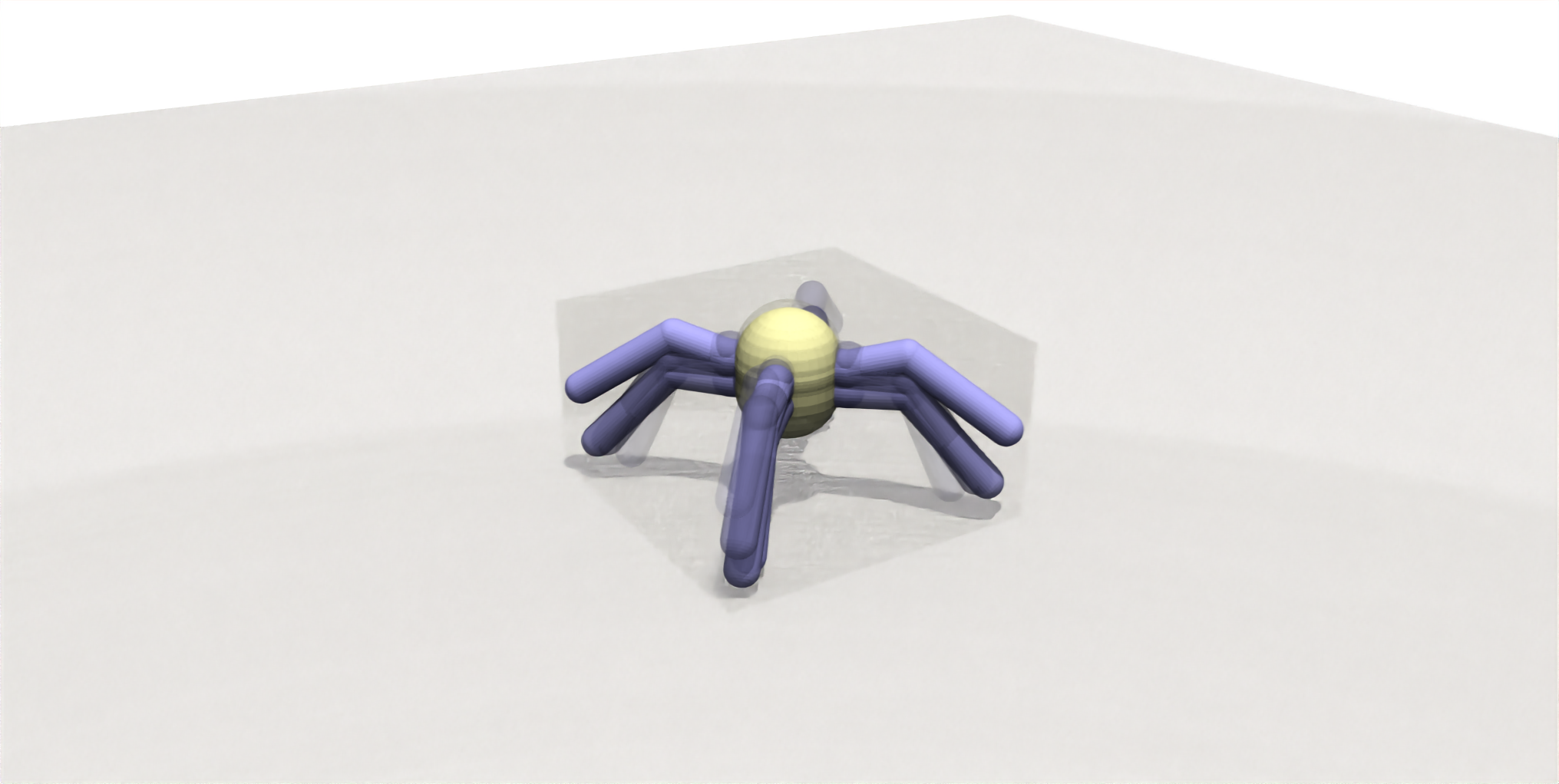}};
\end{tikzpicture}
\caption{\label{fig:discrepancy-ant}Left: Maximal- versus minimal-coordinate simulator accuracy on the 18-DOF ant trotting gait. Right: per-frame discrepancy between the two simulators over time, all on a logarithmic scale.}
\end{figure}

\begin{figure}[t]
\centering
\begin{tikzpicture}
\node[anchor=north west,inner sep=0] at (0.10\linewidth,0)
  {\includegraphics[width=0.90\linewidth]{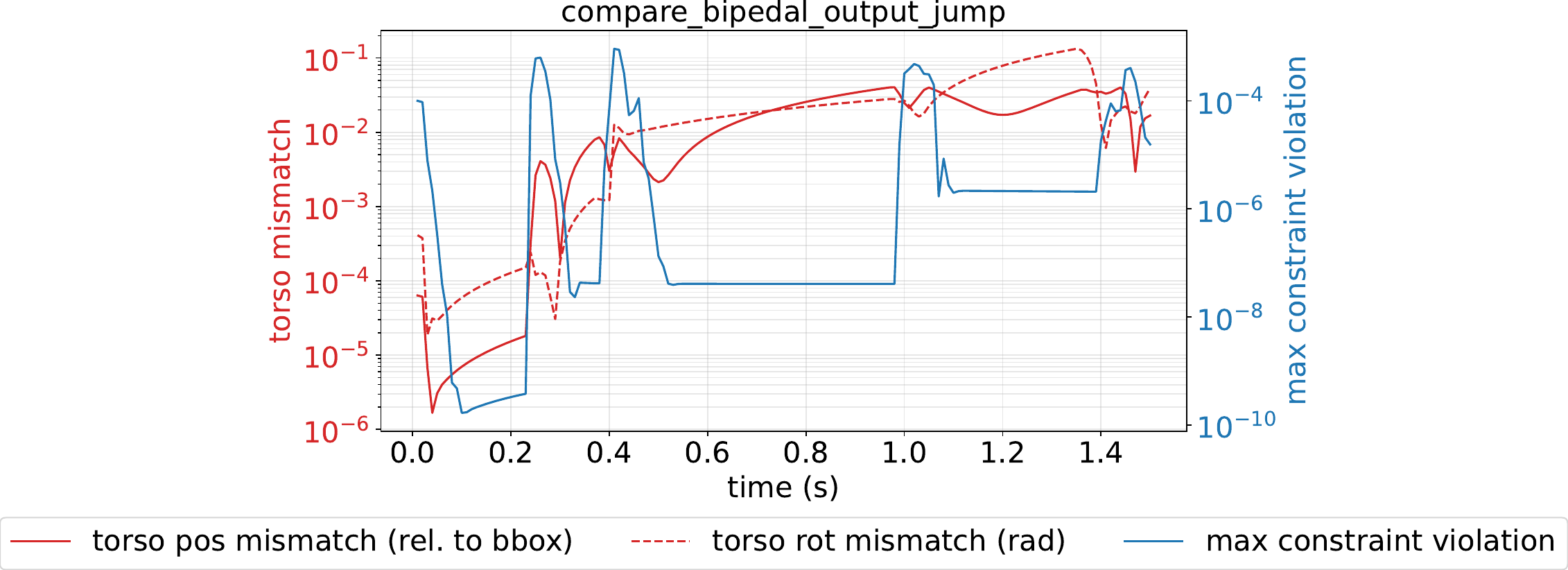}};
\node[anchor=west,inner sep=0] at (0,-0.148\linewidth)
  {\includegraphics[trim={620 162 511 190},clip,frame,width=0.24\linewidth]{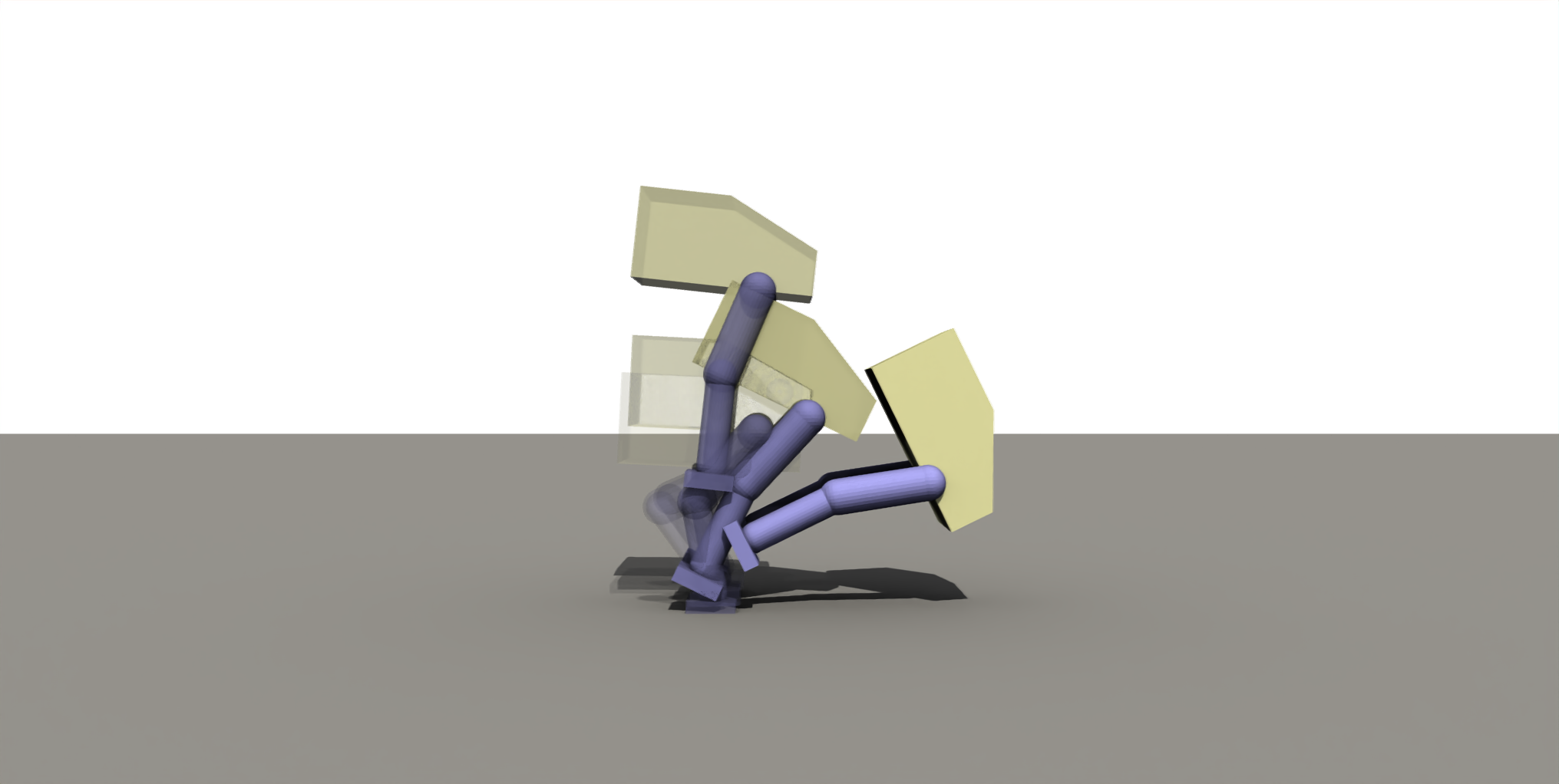}};
\end{tikzpicture}
\caption{\label{fig:discrepancy-bipedal}Left: Maximal- versus minimal-coordinate simulator accuracy on the 14-DOF bipedal jump. Right: per-frame discrepancy between the two simulators over time, all on a logarithmic scale.}
\end{figure}

\section{\label{sec:experiments}Numerical Experiments}
We conduct preliminary numerical experiments using a modified version of the open-source code released along with~\cite{11266917}, to which we make two changes. First, we replace the minimal-coordinates formulation with the maximal-coordinates one. Second, we modify their contact model to incorporate all the energy mollifications discussed in~\prettyref{sec:Hamiltonian}. On top of this, we implement a customized line-search SQP algorithm following~\prettyref{sec:global}. In the notation of the theory, the implementation applies the SQP framework to the articulated-body instance of~\prettyref{pb:CITO-M}. The reduced articulated-body surrogate~\prettyref{pb:PSR-CRDS}, after the mollification of~\prettyref{pb:PS-CRDS-M}, supplies the per-step stationarity residual used as the physics constraint. The original~\prettyref{pb:CITO} is not solved directly; its force-parameterized stationarity certificate is reconstructed numerically after convergence. The value functions in~\prettyref{eq:reduced-penalty} are defined theoretically by exact mathematical partial minimization, whereas the implementation necessarily approximates those inner minimizers in finite precision. We follow~\cite{11266917} and use mixed precision: all computations are performed under 64-bit double precision, except for the separating-plane and plane-velocity inner solves, which use 128-bit quadmath~\cite{bailey2005high}. Across all such inner solves in the reported experiments, the maximum final gradient norm with respect to the inner variables is below $10^{-14}$, with approximately $20$ inner iterations per solve on average. This rapid convergence reflects the convexity of the inner problems and our use of a Newton solver. Thus the experiments approximate the exact-reduction method to high numerical accuracy but are not a literal exact-arithmetic realization of the theorem. We use several implementation techniques to improve practical convergence speed. First, we attempt a standard second-order correction (SOC)~\cite{nocedal2006numerical} after an unsuccessful trial of the baseline line search. The correction proposes only an alternative trial point and is accepted only if it satisfies the same control-admissibility and safeguarded merit and constraint-violation tests as~\prettyref{eq:sqp-linesearch}. If the corrected trial fails, it is discarded and backtracking continues with the original QP direction. Thus SOC never bypasses the safeguards of the analyzed method. Because the corrected displacement is not literally the scaled QP direction in~\prettyref{alg:sqp-inner}, we treat SOC as an optional practical acceleration; the formal convergence theorem remains stated for the baseline update. Second, the Jacobian matrix $J^k$ is narrow-banded and lower triangular, which allows the QP subproblem to be reduced efficiently using a fast block-tridiagonal sparse linear solver such as~\cite{frison2020hpipm}. Specifically, the full rank of $J^k$ allows us to assemble a reduced system in the control step $s_\StackCtrl^k$, which we denote by QP-assembly. An active-set routine then solves the same convex QP~\prettyref{pb:sqp-qp} for $s_\StackCtrl^k$, after which we recover $s_\StackConf^k$ and the associated QP multipliers; we denote this stage by QP-solve. The active-set routine changes only how the QP is solved, and the convergence analysis does not prescribe a particular QP solver. Finally, we use a damped L-BFGS update~\cite{nocedal2006numerical} to propose the QP Hessian approximation. Before every QP solve, we apply curvature clamping so that the matrix supplied to~\prettyref{pb:sqp-qp} satisfies the uniform spectral bounds in~\prettyref{eq:sqp-hessian}. Thus the quasi-Newton update changes the practical approximation quality but not the Hessian-conditioning hypothesis used in the convergence proof. We use multi-threading to assemble the system matrices at each timestep in parallel.

\subsection{CITO Benchmarks}
We evaluate our method on a suite of trajectory optimization tasks, three of which concern robot locomotion gait discovery and the other two focus on table-top manipulation. We first relate these benchmarks to the theorem's initialization assumptions. The convergence theorem permits any finite optimized-node set $\mathcal I$, finite configurations $\theta_i$ for $i\in\mathcal I$, and admissible controls $u_i\in\mathcal D$. In every benchmark, the fixed time-zero robot configuration $\theta_0$ is initialized floating above the ground, with minimum ground-contact inequality slack equal to $\EpsFJ>0$ under the barrier model. Thus $\theta_0$ is strictly interior and satisfies~\prettyref{ass:Slater}; the first optimized node $\min\mathcal I$ may subsequently enter the $\EpsFJ$ contact layer. This fixed time-zero configuration should be distinguished from the dynamically infeasible optimized-trajectory guesses used to warm-start later receding-horizon solves. Those guesses are finite at every node of the finite experimental set $\mathcal I$ and are therefore covered by~\prettyref{thm:sqp-convergence}; the constraint-satisfaction phase drives their finite residual into the working band. For all experiments, we adopt a receding-horizon control paradigm as in model predictive controllers (MPC): we repeatedly optimize a trajectory of horizon $H$ but only commit the control signals of the first $H/2$ steps. We then extend the horizon by padding the remaining $H/2$ steps, using the last control signal and robot configuration to form the initial guess for the next optimization window of horizon $H$. This process is repeated until the task is completed. For all experiments, we use $H=2s$, an initial local timestep $\delta_i=0.025s$ for every $i\in\mathcal I$, $\EpsStop=0.001$, $\bar\varrho=1$, $\LLambdaMin=0.01$, and $\LResidualWork=1$. Every reported MPC solve was continued until its final inner SQP call satisfied both strict stopping tests in~\prettyref{alg:sqp-inner}, $\|\StatRes^k\|_1<\EpsStop$ and $\|s^k\|<\EpsStop/\LHessianCond$. No reported solve was terminated by a separate wall-clock or iteration-count rule. The adaptive minimum-singular-value check $\sigma_{\min}(\StatJac_{ii})\geq\LLambdaMin$ of~\prettyref{alg:sqp-inner}--\prettyref{alg:sqp} and the Hamiltonian check of~\prettyref{alg:CITO} remain enabled in all experiments. Every reported benchmark passes both checks throughout the solve, so the numbers of singularity-triggered and Hamiltonian-triggered subdivisions are both zero. Adaptive refinement is therefore an inactive safeguard on these instances rather than a mechanism omitted from the implementation. Other parameters, including $\EpsFri$ and $\EpsFJ$, are fine-tuned per experiment. We select them so that the simulated motion stays physically faithful: the maximal-coordinate constraint violations remain small and the resulting trajectory stays close to a minimal-coordinate simulation of the same system, both of which are quantified in~\prettyref{sec:max-vs-min}. Consequently, we do not drive $\EpsFri\downarrow0$ or $\EpsFJ\downarrow0$ toward the true (stiff and non-smooth) dynamics, since this vanishing-regularization regime lies outside the theory's uniform-boundedness and conditioning guarantees (\prettyref{thm:stationarity-transfer} and~\ref{thm:nonsmooth-transfer}). All experiments are performed on a server with a 166-core AMD EPYC-Genoa CPU and 128GB memory.

\begin{figure}[t]
\centering
\setlength{\fboxsep}{0pt}%
\fbox{\includegraphics[trim={740 606 484 110},clip,width=\linewidth]{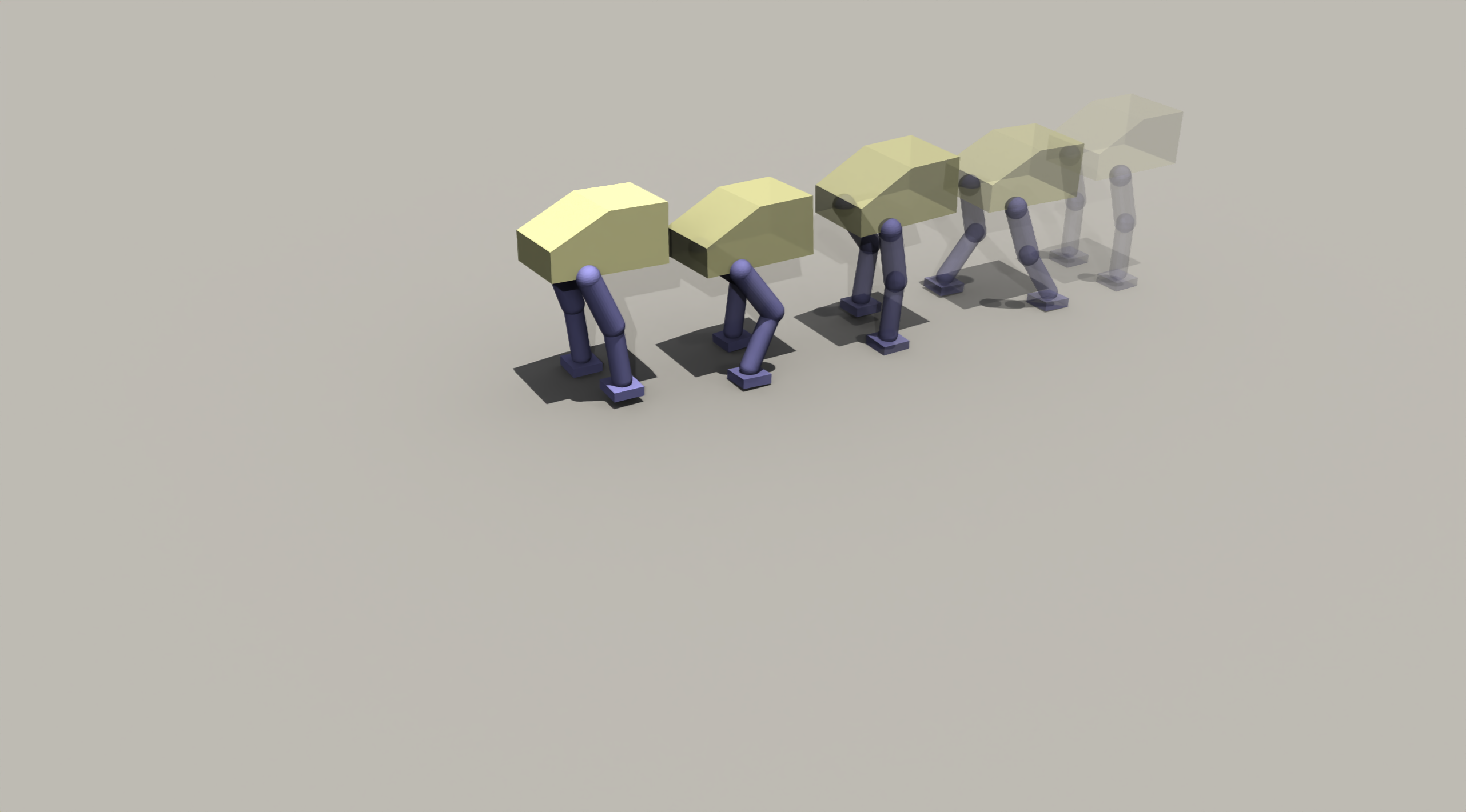}}
\caption{\label{fig:bipedal}Accumulated overlay of the generated locomotion trajectory for the 14-DOF bipedal robot on flat ground. Starting with zero predecessor velocity and trivial optimized controls, while the fixed predecessor is placed $\EpsFJ$ above the ground, the trajectory optimizer discovers a forward walking gait over the long horizon. Poses along the trajectory are superimposed and fade from opaque (early timesteps) to transparent (late timesteps), visualizing the forward progress from left to right.}
\end{figure}

\begin{figure}[t]
\centering
\setlength{\fboxsep}{0pt}%
\fbox{\includegraphics[trim={300 516 564 180},clip,width=\linewidth]{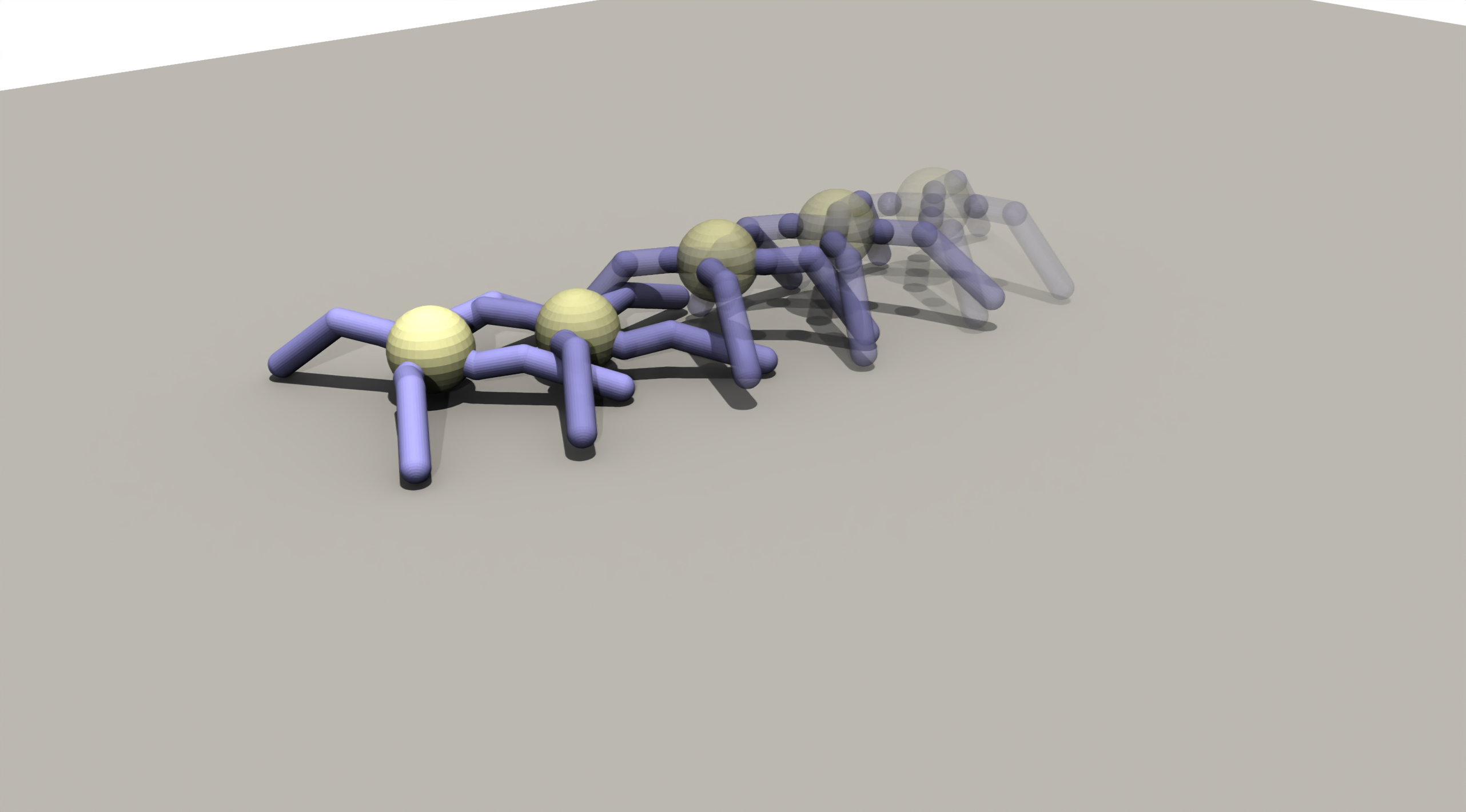}}
\caption{\label{fig:ant}Accumulated overlay of the generated locomotion trajectory for the 18-DOF quadruped ant robot on flat ground. Starting with zero predecessor velocity and trivial optimized controls, while the fixed predecessor is placed $\EpsFJ$ above the ground, the trajectory optimizer discovers a forward gait over the long horizon. Poses along the trajectory are superimposed and fade from opaque (early timesteps) to transparent (late timesteps), visualizing the forward progress from left to right.}
\end{figure}

\begin{figure*}[t]
\centering
\setlength{\fboxsep}{0pt}%
\setlength{\tabcolsep}{2pt}
\resizebox{\textwidth}{!}{%
\begin{tabular}{ccc}
\fbox{\includegraphics[trim={900 201 574 140},clip,height=4cm]{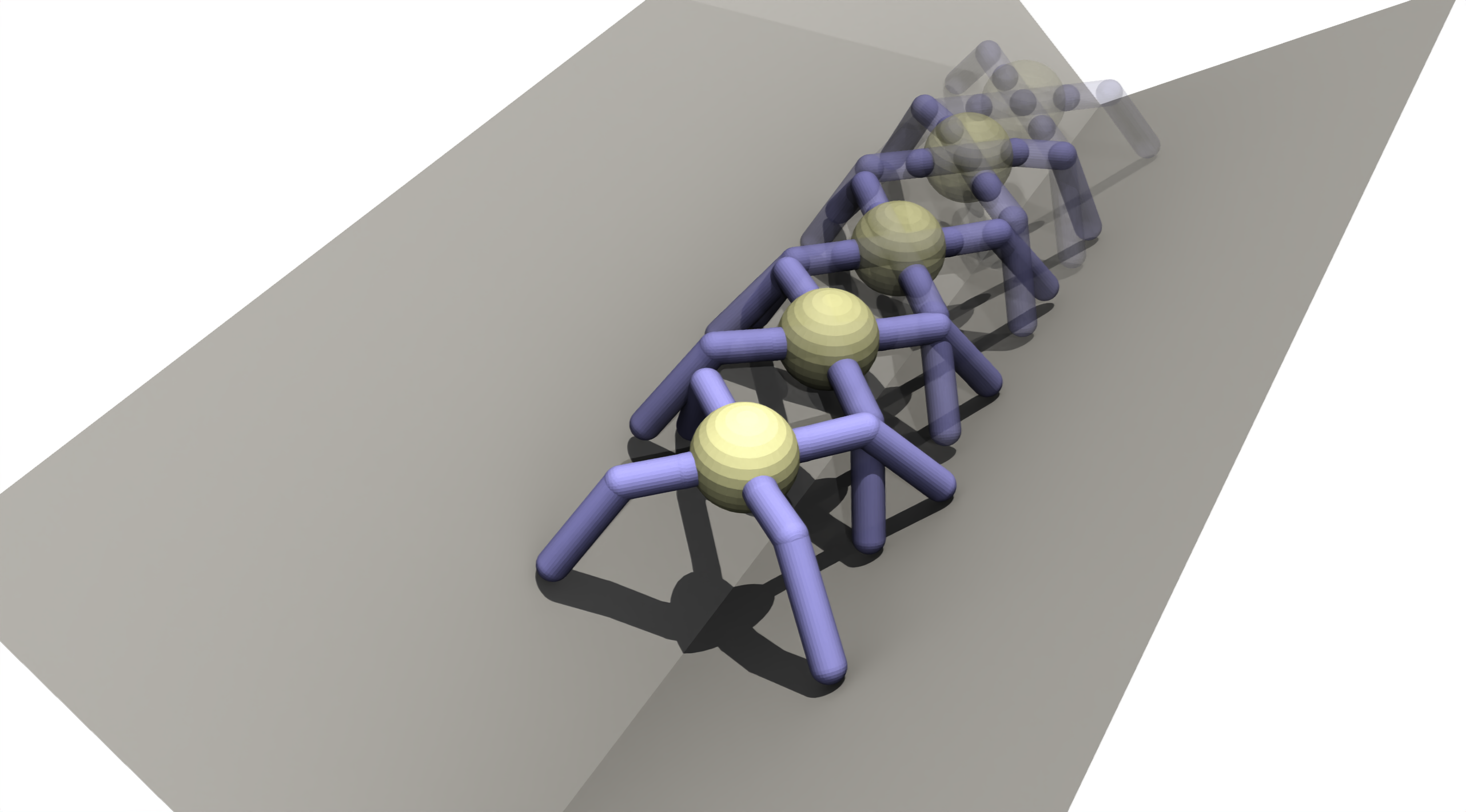}} &
\fbox{\includegraphics[trim={690 464 714 412},clip,height=4cm]{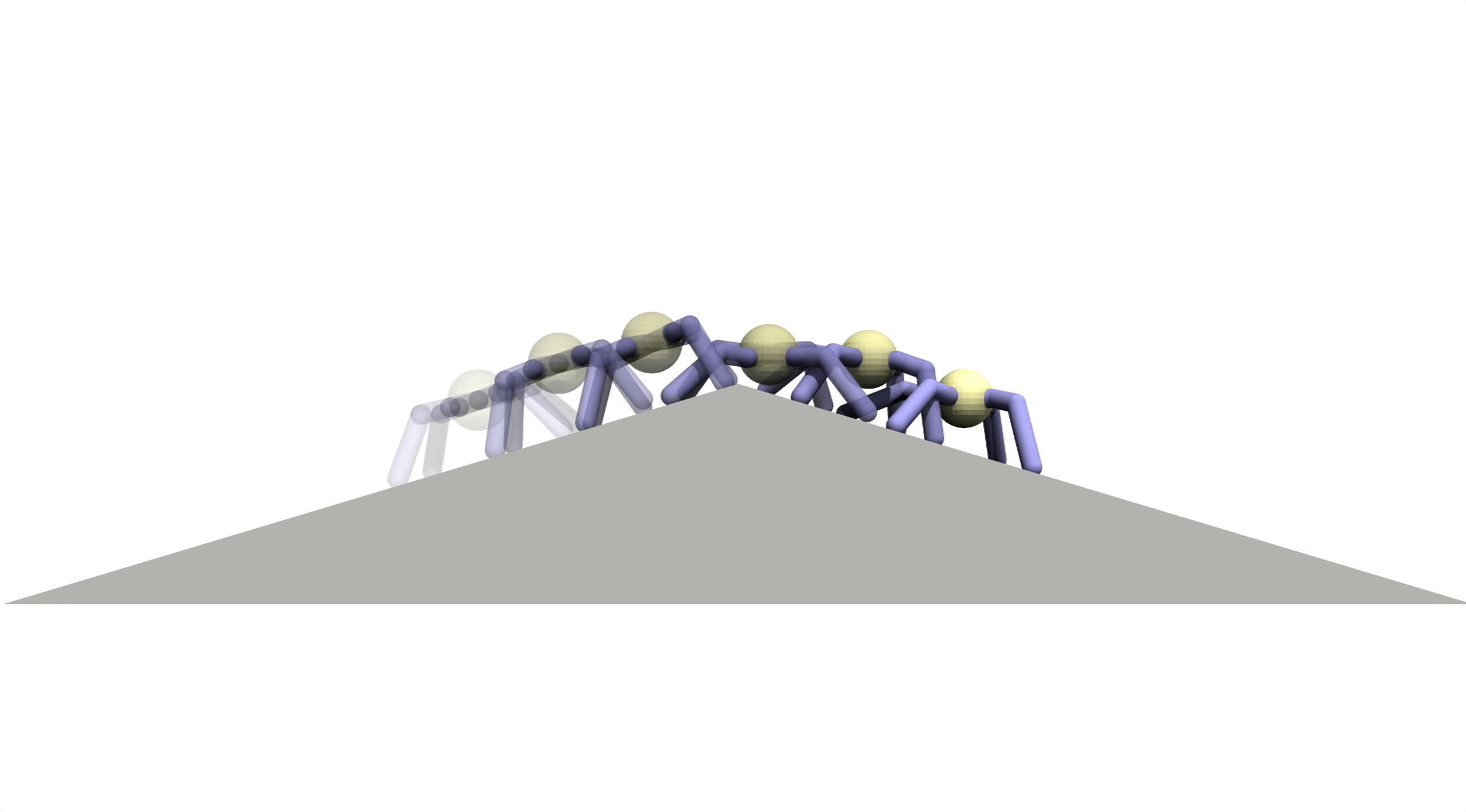}} &
\fbox{\includegraphics[trim={720 183 354 655},clip,height=4cm]{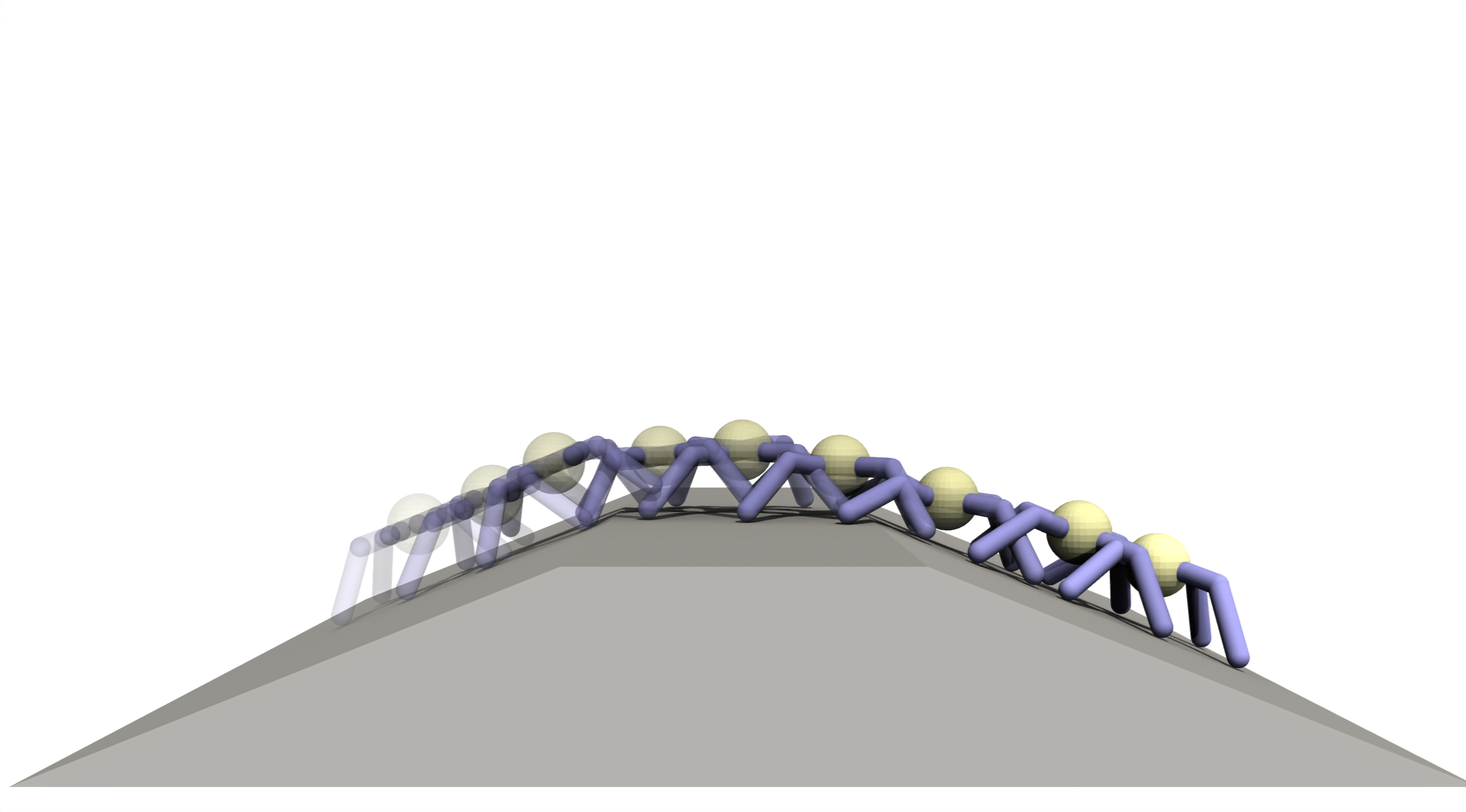}} \\
\end{tabular}%
}
\caption{\label{fig:ant-terrain}The 18-DOF quadruped ant robot traversing three different kinds of terrain: an inclined ramp (left), a sharp ridge (middle), and a smooth hill (right). For each terrain, poses along the generated trajectory are superimposed and fade from opaque (early timesteps) to transparent (late timesteps), visualizing the traversal.}
\end{figure*}

Our first locomotion task has a 14-DOF ($7\times12$-DOF under maximal coordinates) bipedal robot move forward on flat ground along the x-axis at a constant speed of $0.3m/s$. Our objective function is defined as
\begin{ResizedEq}
O(\{\theta_i\},\{u_i\})=
&O^\text{foot}(\{\theta_i\})+O^\text{torso}(\{\theta_i\})+\\
&5O^\text{target}(\{\theta_i\},0.3)+0.01O^\text{reg}(\{u_i\})\\
O^\text{foot}=&\int_{t\in[0,H]}\sum_{j=0}^1\|n^{\text{foot},j}(\theta(t))-z\|^2dt\\
O^\text{torso}=&
\int_{t\in[0,H]}\|z^\text{torso}(\theta(t))-z^\text{torso}(\theta_0)\|^2dt+\\
&\int_{t\in[0,H]}\|R^\text{torso}(\theta(t))-R^\text{torso}(\theta_0)\|^2dt\\
O^\text{target}(\{\theta_i\},\dot{x}^{\text{torso},\star})=&\int_{t\in[0,H]}\|\dot{x}^\text{torso}(\theta(t))-\dot{x}^{\text{torso},\star}\|^2dt\\
O^\text{reg}=&\int_{t\in[0,H]}\|u(t)\|^2dt.
\end{ResizedEq}
Note that all objective functions are expressed as integrals over the trajectory, where the configuration and control histories are discretized as piecewise-constant functions of time ($\theta(t)$ and $u(t)$). Therefore, the objective is invariant to trajectory subdivision. Here $O^\text{foot}$ requires that the robot's foot normal ($n^{\text{foot},j}$) is aligned with the z-axis. $O^\text{torso}$ requires the robot torso to stay upright with minimal rotation ($R^\text{torso}$) and remain at the same altitude ($z^\text{torso}$). $O^\text{target}$ is the target velocity tracking term and $O^\text{reg}$ is a small control regularization. For this task, the optimized controls and base velocity $\dot\theta_0=(\theta_0-\theta_{\gamma(0)})/\delta_0$ are initialized trivially, while the fixed time-zero robot configuration $\theta_0$ is placed $\EpsFJ$ above the ground and hence is not in active support contact. Our trajectory optimizer automatically discovers a locomotion gait over the long horizon, with frequent contact-mode changes against the ground, as illustrated in~\prettyref{fig:bipedal}. Our second task is likewise a gait search problem, but applied to an 18-DOF ($9\times12$-DOF under maximal coordinates) quadruped ant robot, with the generated trajectory illustrated in~\prettyref{fig:ant}. Since the ant robot is passively stable and has longer feet, we allow it to move faster at a speed of $0.6m/s$, using the following simplified objective
\begin{ResizedEq}
O(\{\theta_i\},\{u_i\})=&
5O^\text{target}(\{\theta_i\},0.6)+
0.01O^\text{reg}(\{u_i\}).
\end{ResizedEq}
Our third task has the ant robot traverse three different kinds of terrain, an inclined ramp, a sharp ridge, and a smooth hill, as shown in~\prettyref{fig:ant-terrain}, where we use the same objective function configuration and parameters as the case with flat ground. We also tried more complex terrain shapes; unfortunately, our solver becomes stuck in local minima and cannot make meaningful progress. This failure is consistent with the local solution quality certified by the theory. The global-convergence result concerns finite termination of the specified SQP method, from any finite discretized optimized trajectory guess within the admissible initialization class of~\prettyref{rem:initialization-scope}, at an approximately feasible force-parameterized stationary trajectory; it does not ensure a globally optimal, or even task-successful, trajectory.

\begin{figure*}[t]
\centering
\setlength{\fboxsep}{0pt}%
\setlength{\tabcolsep}{0.5pt}
\renewcommand{\arraystretch}{0.5}
\resizebox{\textwidth}{!}{%
\begin{tabular}{ccccc}
\includegraphics[trim={650 50 650 100},clip,width=0.24\textwidth,frame]{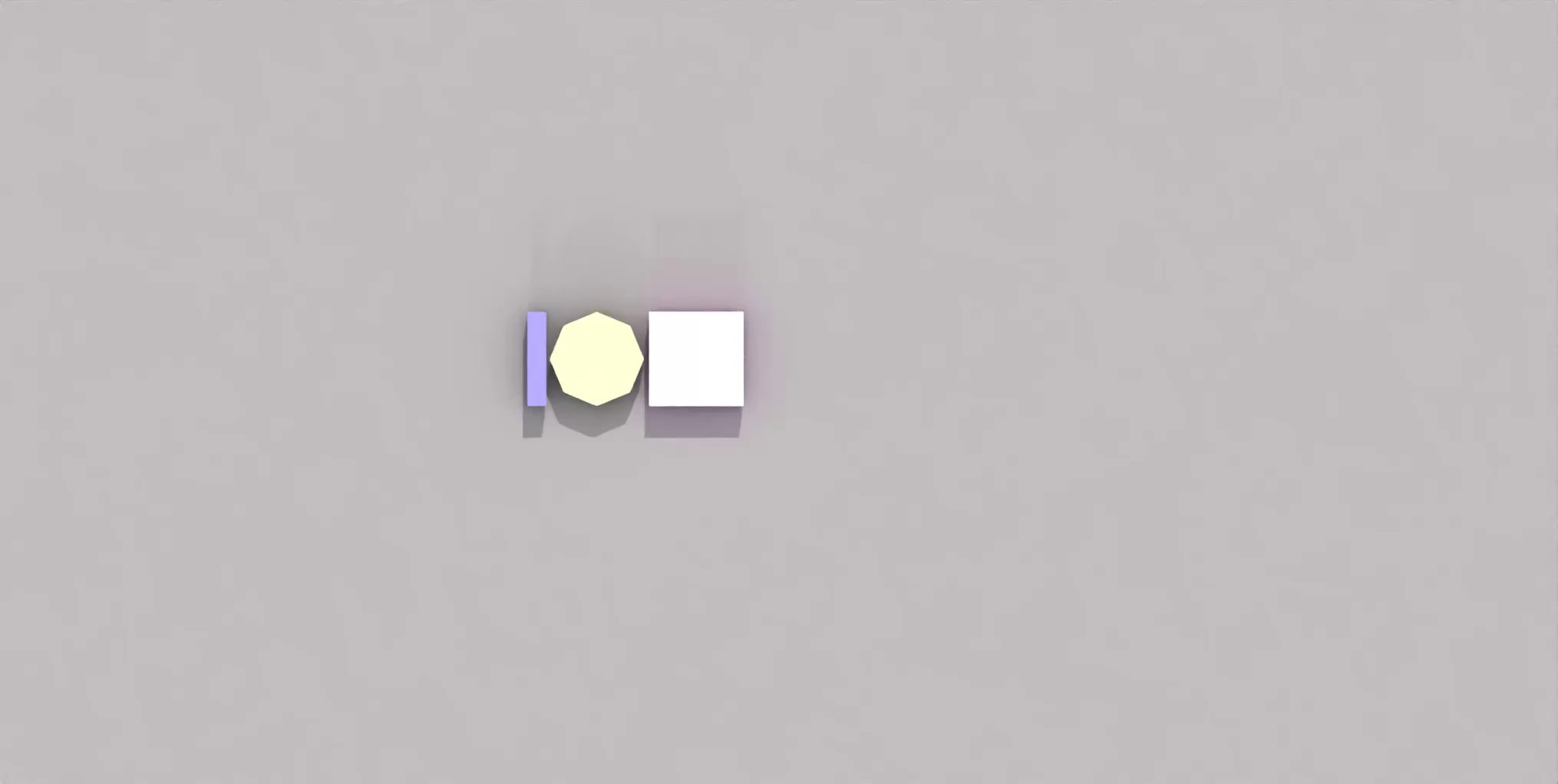} & \includegraphics[trim={650 50 650 100},clip,width=0.24\textwidth,frame]{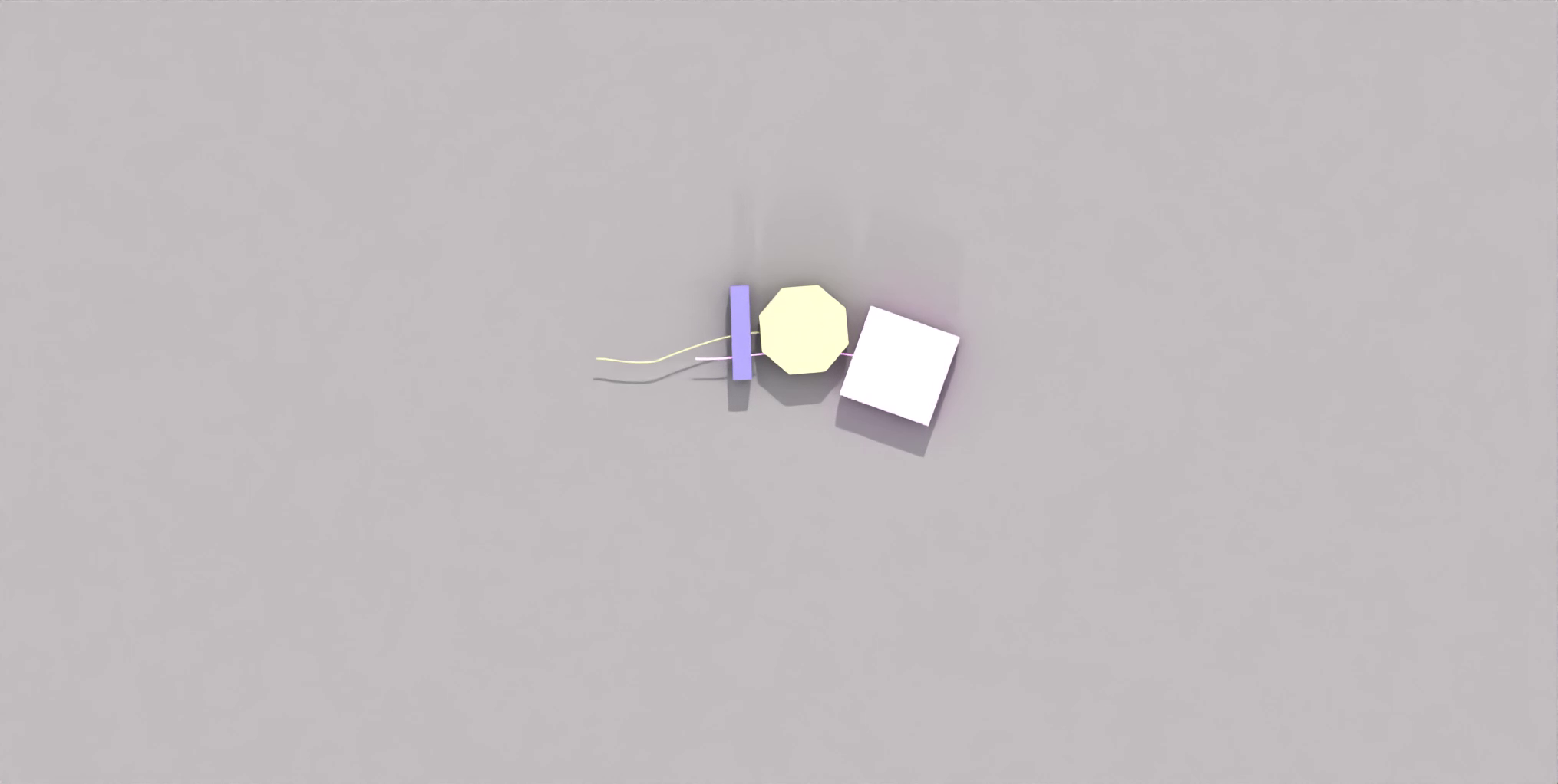} & \includegraphics[trim={650 50 650 100},clip,width=0.24\textwidth,frame]{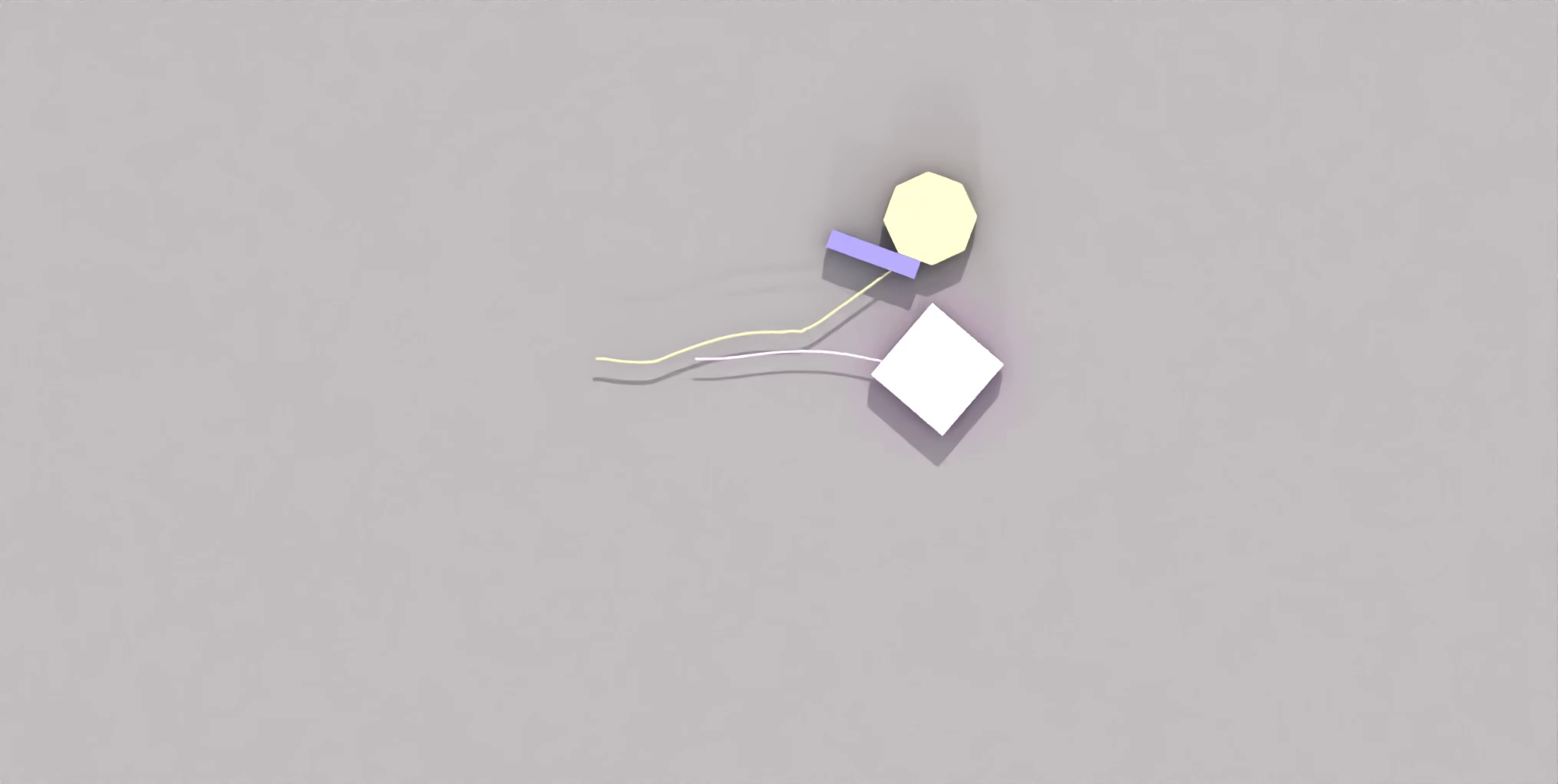} & \includegraphics[trim={650 50 650 100},clip,width=0.24\textwidth,frame]{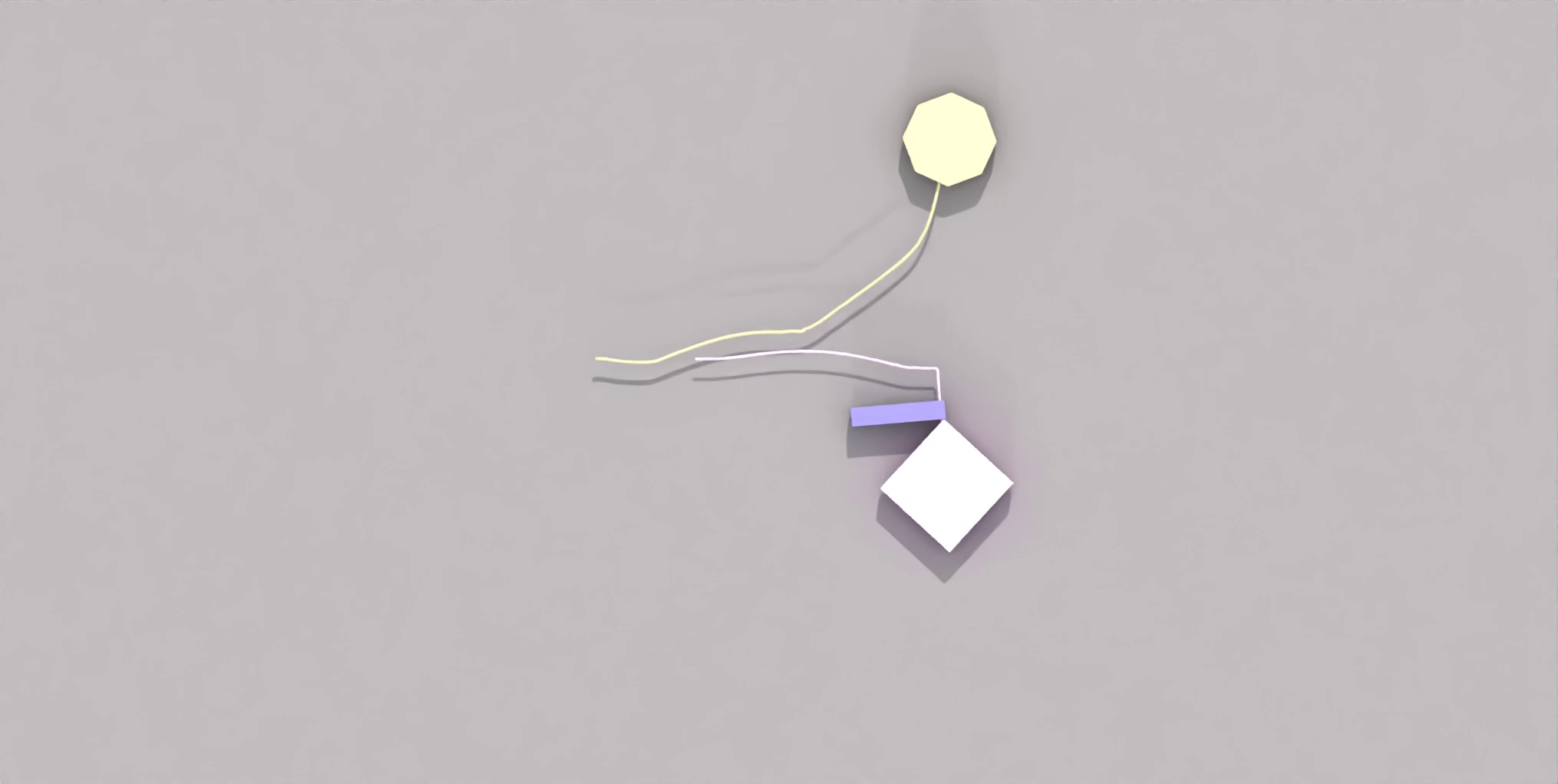} & \includegraphics[trim={650 50 650 100},clip,width=0.24\textwidth,frame]{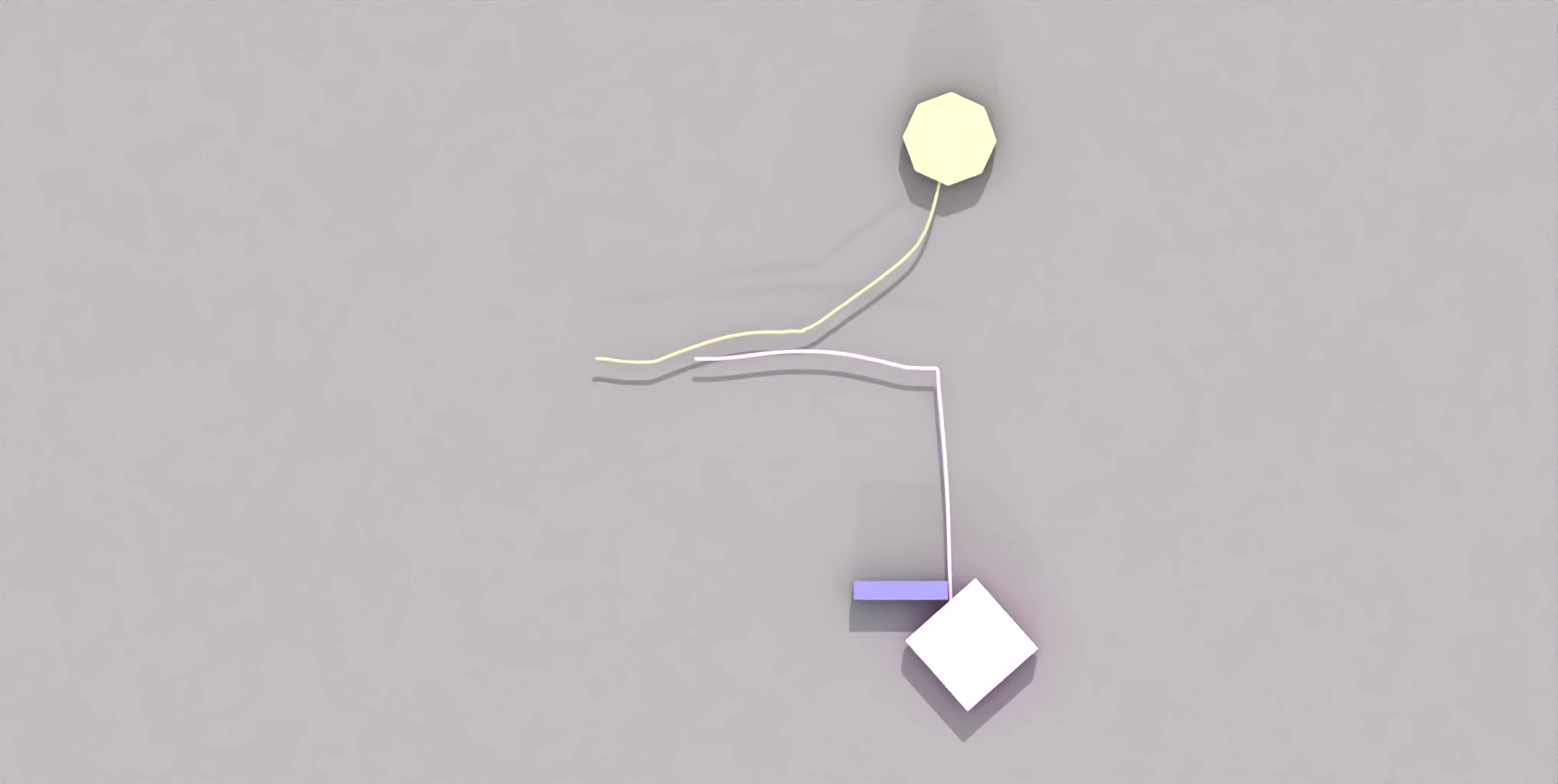} \\
\end{tabular}%
}
\caption{\label{fig:push}Frames of the generated table-top pushing trajectory. The controllable purple pusher must separate the two initially stacked objects---a butter-yellow octagon and a white square---and drive them to two target positions that are far apart. The thin curves trace the world-frame center-of-mass trajectories of the octagon (yellow) and the square (white) up to each frame, showing that the pusher first nudges the octagon upward and then reroutes to push the square downward.}
\end{figure*}

\begin{figure*}[t]
\centering
\setlength{\fboxsep}{0pt}%
\setlength{\tabcolsep}{0.5pt}
\renewcommand{\arraystretch}{0.5}
\resizebox{\textwidth}{!}{%
\begin{tabular}{ccccc}
\fbox{\includegraphics[trim={1090 465 1020 520},clip,width=0.19\textwidth]{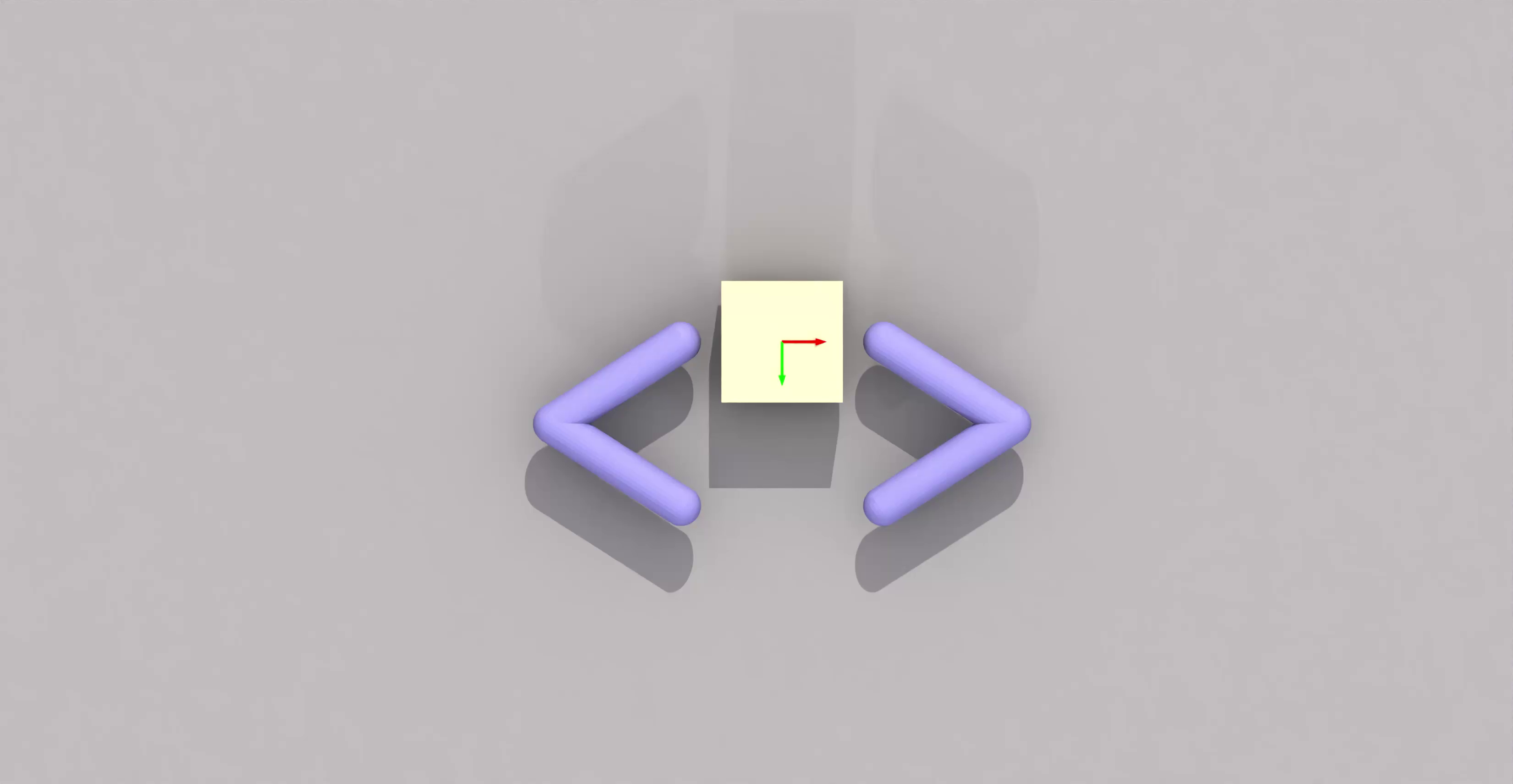}} & \fbox{\includegraphics[trim={1090 465 1020 520},clip,width=0.19\textwidth]{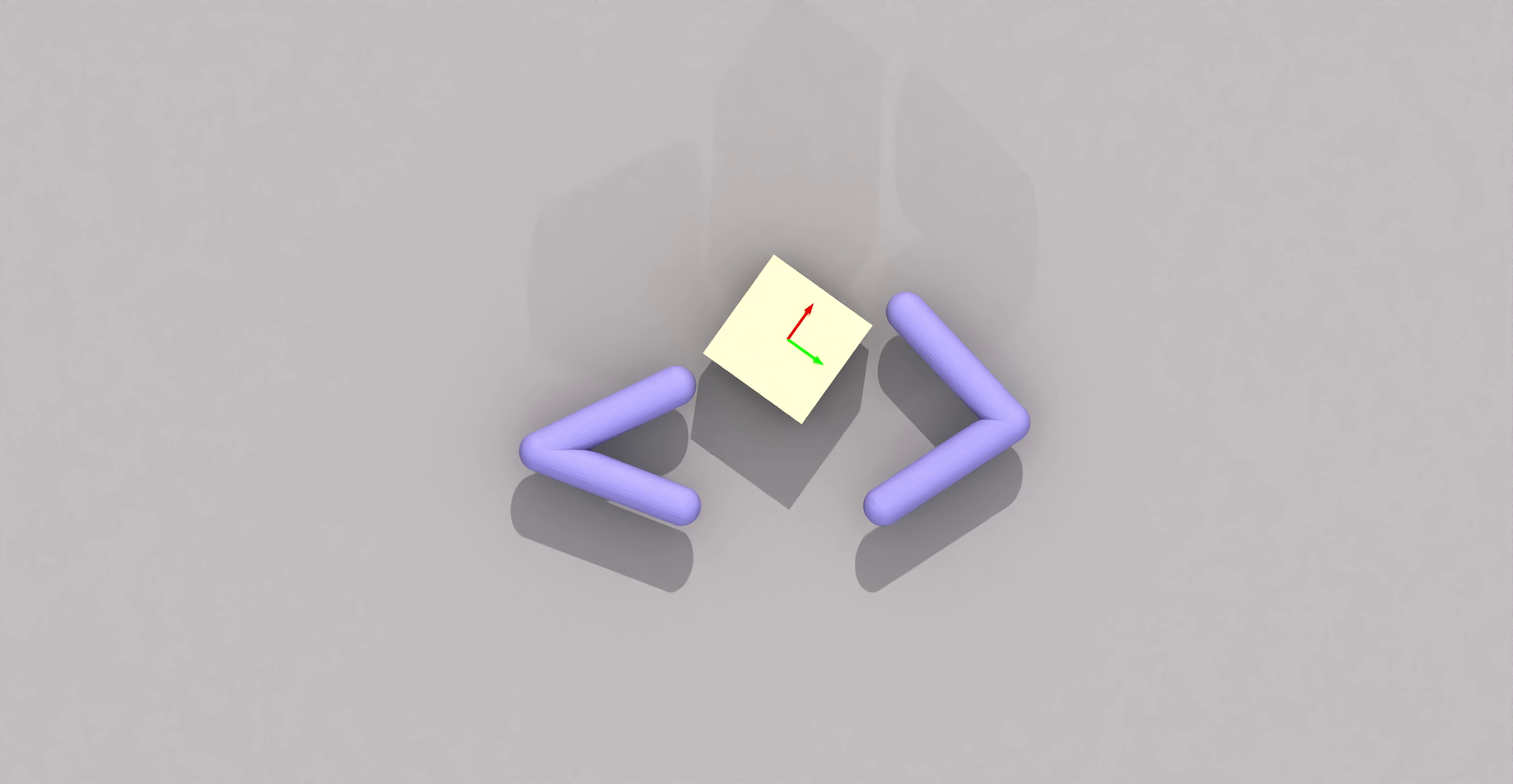}} & \fbox{\includegraphics[trim={1090 465 1020 520},clip,width=0.19\textwidth]{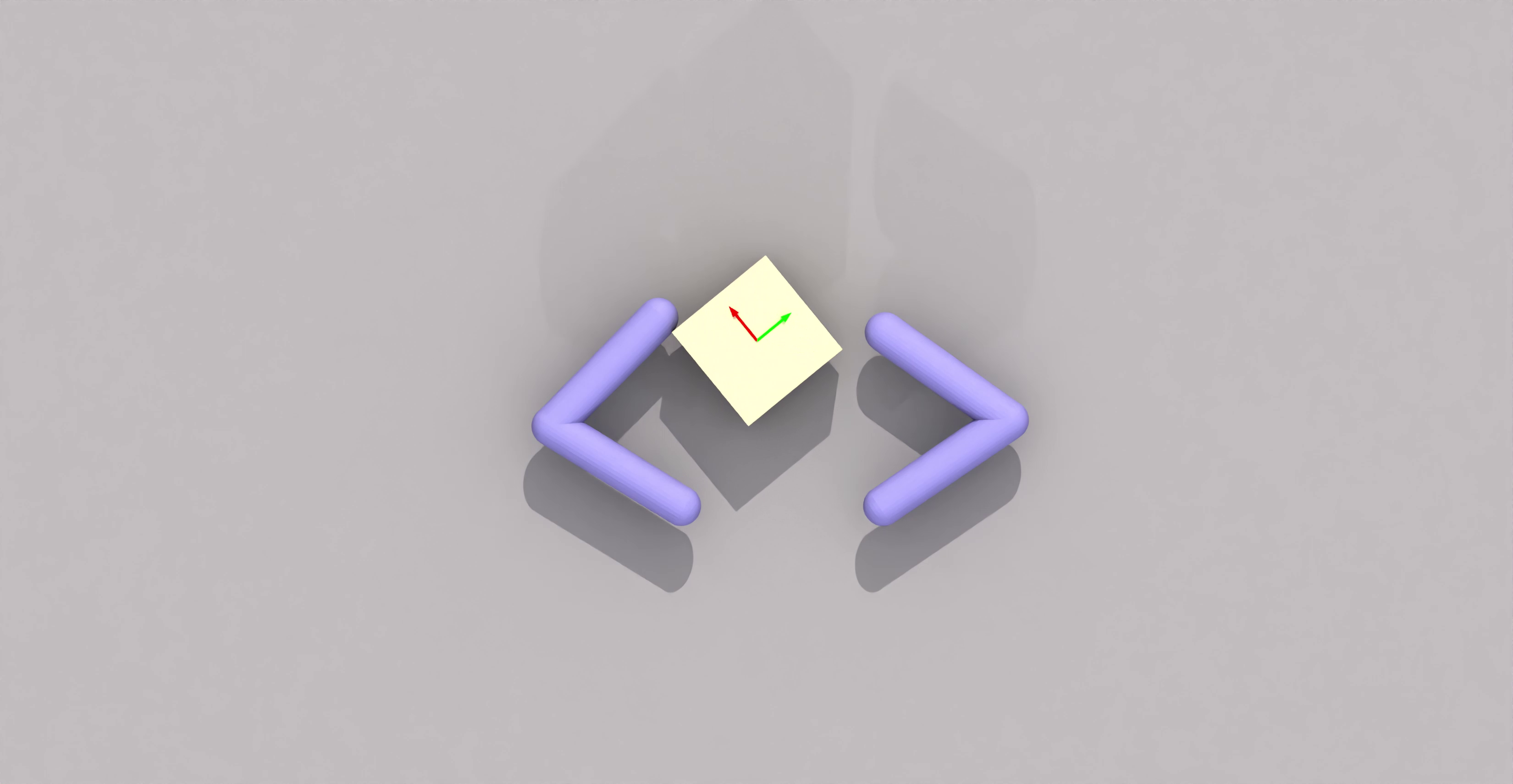}} & \fbox{\includegraphics[trim={1090 465 1020 520},clip,width=0.19\textwidth]{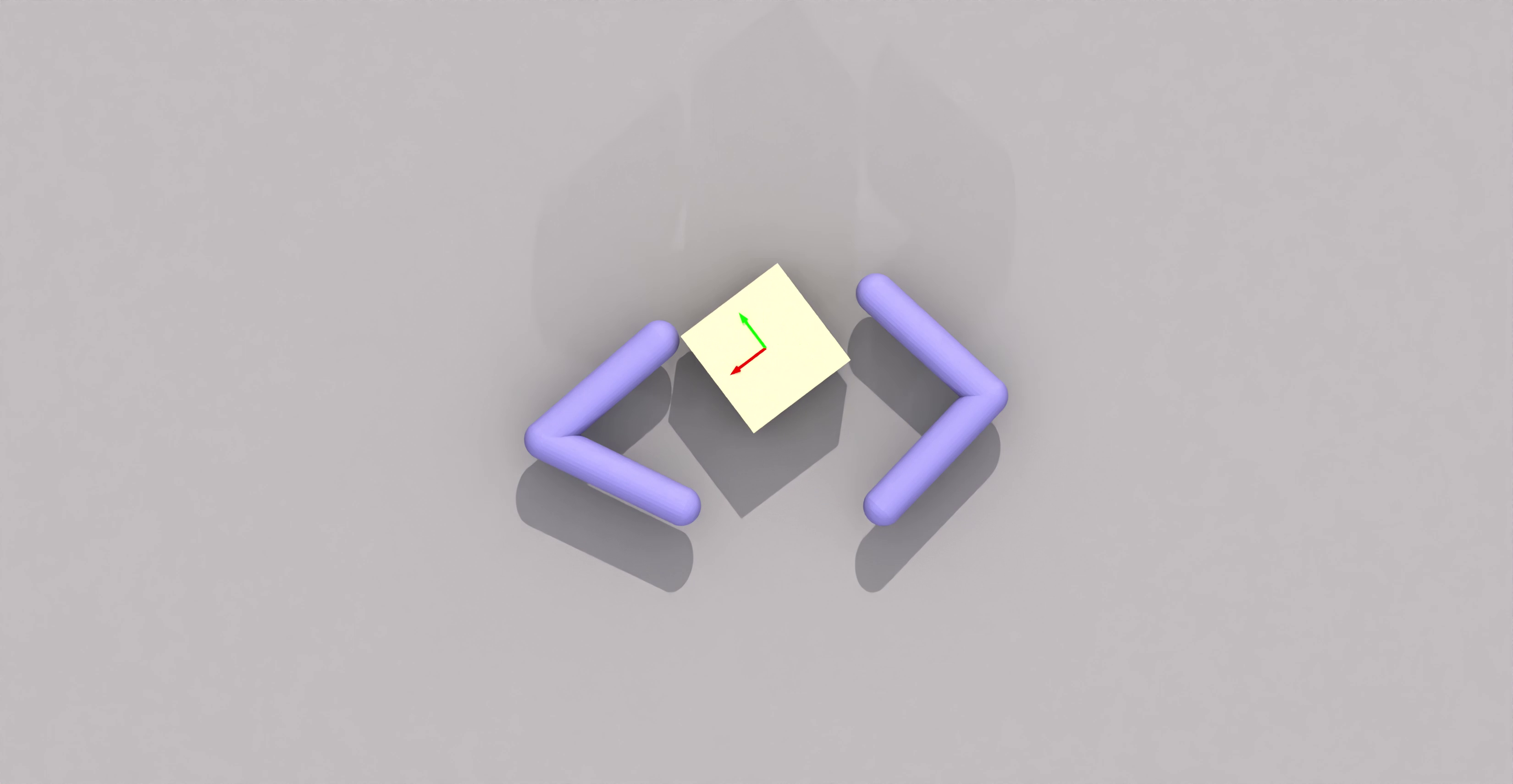}} & \fbox{\includegraphics[trim={1090 465 1020 520},clip,width=0.19\textwidth]{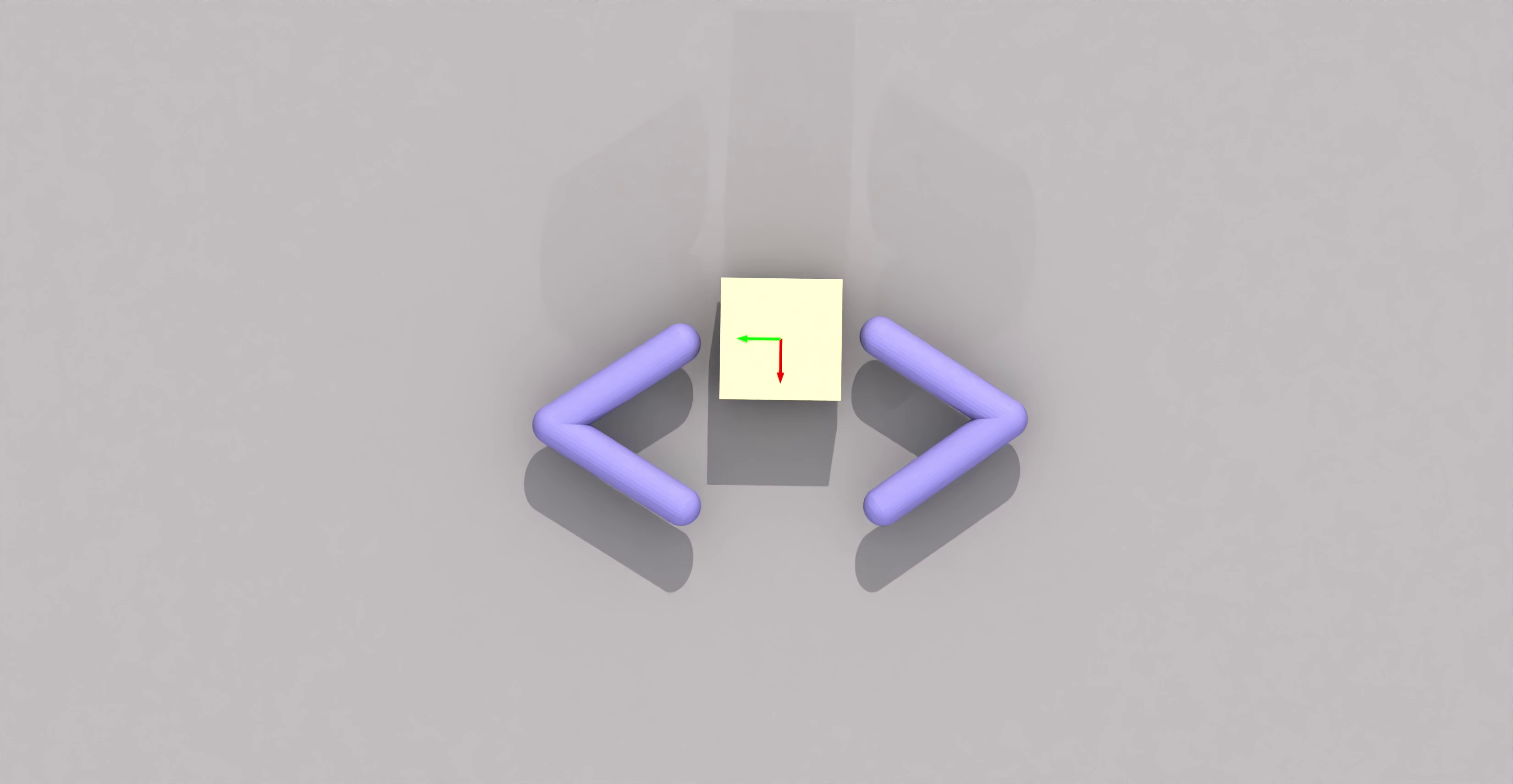}} \\
\end{tabular}%
}
\caption{\label{fig:rotate}Frames of the generated in-hand rotation trajectory, in which a two-finger gripper---each finger actuated by two hinge joints---reorients a cubic object by $270^\circ$. The red and green arrows drawn on the object denote its body-frame axes, making the accumulated rotation visible across frames as the fingers repeatedly regrasp and turn the cube.}
\end{figure*}

Besides locomotion tasks, we evaluate our method on two table-top manipulation scenarios. Our first scenario consists of three objects, where the purple pusher's orientation and position can be controlled while the two other objects are unactuated. Each object has 3 planar DOFs, totaling 9 DOFs under minimal coordinates ($3\times4$ DOFs under maximal coordinates). The three objects are initially stacked together, and the pusher is tasked with moving the two other objects to two target positions ($t^{0,\star}$ and $t^{1,\star}$) that are far apart. For this example, we use the following cost as the objective function
\begin{ResizedEq}
O(\{\theta_i\},\{u_i\})=&
O^\text{terminal}(\{\theta_i\})+
0.01O^\text{reg}(\{u_i\})\\
O^\text{terminal}(\{\theta_i\})=&\sum_{j=0}^1\|t^j(\theta(H))-t^{j,\star}\|^2,
\end{ResizedEq}
where we penalize the terminal position discrepancy. Several frames of the generated pushing trajectory are visualized in~\prettyref{fig:push}, along with the pushed objects' trajectories. From such trivial initialization, our solver automatically finds the switching contact modes for the pusher to first move the butter-yellow octagon upwards and then move the white square downwards.
Our second scenario has a robot gripper rotate a cubic object by $270^\circ$ to achieve the target orientation $\log R^{0,\star}$. The gripper consists of two fingers, each with two actuated hinge joints. The entire system has 7 DOFs ($5\times4$ DOFs under maximal coordinates). Frames of the generated manipulation trajectory are visualized in~\prettyref{fig:rotate}. Note that the position and orientation of the object are both movable instead of fixed. Therefore, the gripper needs to change contact mode to rotate the object while preventing the object from being pushed away from the reachable set. To this end, we use the following objective function that explicitly requires the cubic object to be within reach
\begin{ResizedEq}
O(\{\theta_i\},\{u_i\})=
&O^{\text{terminal},\text{rotate}}(\{\theta_i\})+200O^\text{reach}(\{\theta_i\})+\\
&0.01O^\text{reg}(\{u_i\})\\
O^{\text{terminal},\text{rotate}}(\{\theta_i\})=&\|\log R^0(\theta(H))-\log R^{0,\star}\|^2\\
O^\text{reach}(\{\theta_i\})=&\sum_{j=0}^1\|t^{\text{fingertip},j}(\theta(H))-t^0(\theta(H))\|^2,
\end{ResizedEq}
where $O^\text{reach}$ is a strong regularization that requires the two fingertips ($t^{\text{fingertip},j}$) to be as close as possible to the position of the object ($t^0$). Again, our solver automatically searches for the trajectory that involves several contact mode changes to regrasp the cube.

\subsection{Maximal vs. Minimal Coordinates}
\label{sec:max-vs-min}
Our theoretical analysis relies heavily on the dynamic model being formulated under full maximal coordinates. Theoretically, this poses no problem because the two formulations coincide when the constraint violation is exactly zero. In practice, however, exactly zero constraint violation is not attainable. Fortunately, several prior works~\cite{lan2022affine,10.1145/3811276} have shown that using full maximal coordinates does not degrade physical accuracy in many complex practical setups. For self-containment, we reevaluate the maximal-versus-minimal simulator discrepancy under our default parameter setup on the two locomotion CITO setups introduced above. In both settings, we simulate the robot under both maximal and minimal coordinates and then profile the discrepancy in the robot's torso translation (relative to the robot length scale) and torso rotation (measured in radians), as well as the maximal-coordinate integrity and joint constraint violations (\prettyref{eq:integrity},~\ref{eq:hinge}, and~\ref{eq:prismatic}) at every frame. As shown in~\prettyref{fig:discrepancy-ant}, our first scenario has an 18-DOF ant robot perform a trotting gait; here the constraint violation stays well below $10^{-3}$, the torso's relative translation error stays within $3.6\%$, and the rotation error stays within $0.1^\circ$. Our second scenario, shown in~\prettyref{fig:discrepancy-bipedal}, simulates a jumping trajectory for the 14-DOF bipedal robot, which is more challenging due to the rapid contact-mode change. Here the constraint violation remains well below $10^{-3}$, the torso's relative translation error stays below $4\%$, and the rotation error stays within $7.6^\circ$.

\subsection{Performance Analysis}
\begin{figure}[t]
\centering
\includegraphics[width=.9\linewidth]{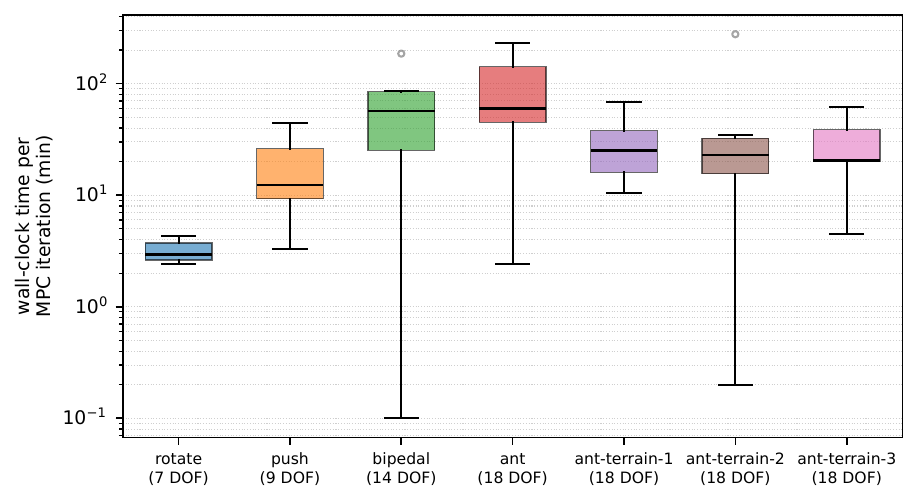}
\caption{\label{fig:cost-dof}Distribution of the wall-clock time per MPC iteration for each benchmark, ordered by the number of robot DOFs (in minimal coordinates). Each box summarizes the per-MPC-iteration times of one benchmark.}
\end{figure}

\begin{figure*}[ht]
\centering
\setlength{\tabcolsep}{0pt}
\begin{tabular}{cc}
\toprule
Bipedal Locomotion & Ant Locomotion\\
\midrule
\includegraphics[trim={0 90 0 0},clip,width=.45\textwidth]{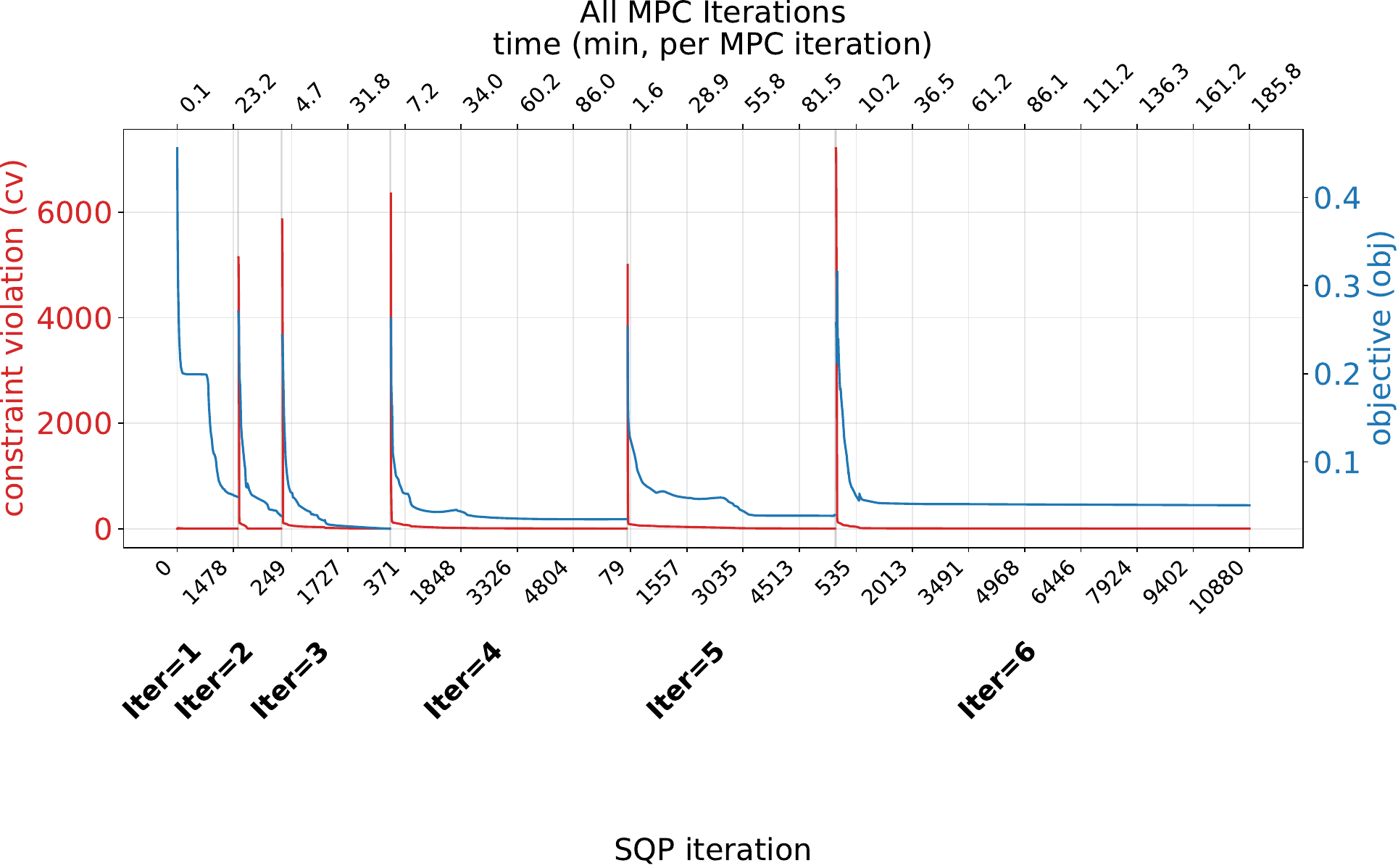} &
\includegraphics[trim={0 90 0 0},clip,width=.45\textwidth]{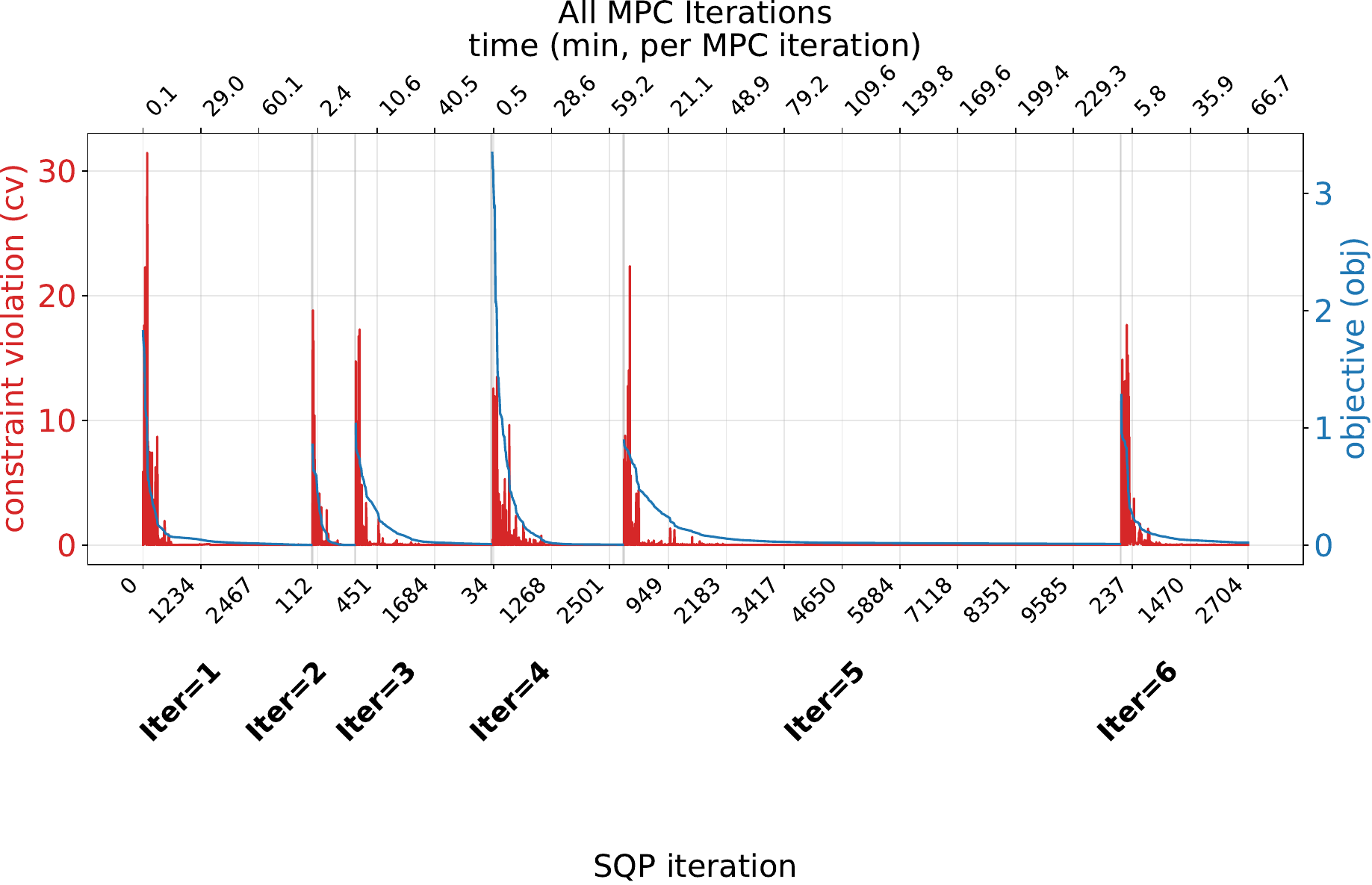} \\
\toprule
Ant-Terrain-1 & Ant-Terrain-2\\
\midrule
\includegraphics[trim={0 90 0 0},clip,width=.45\textwidth]{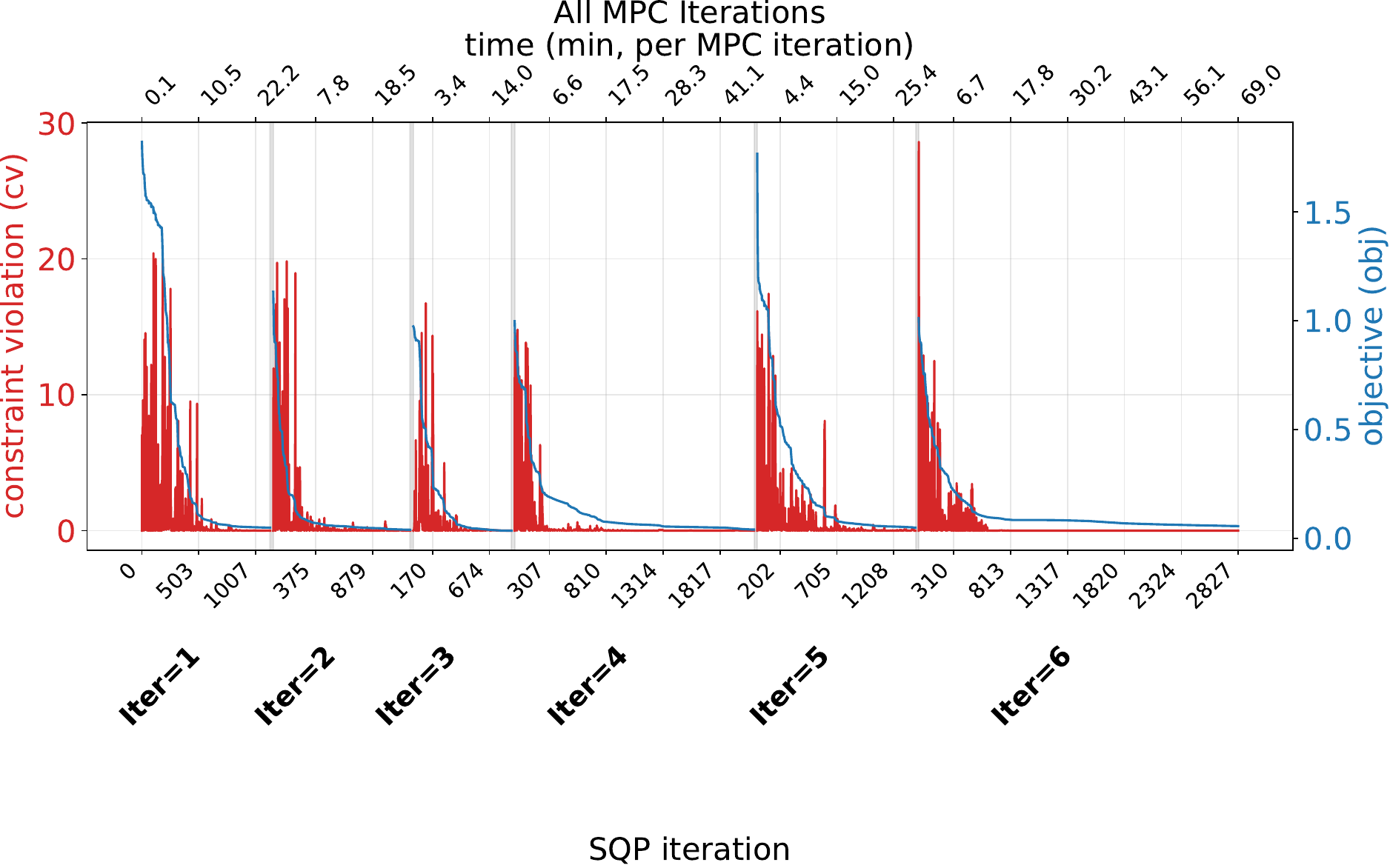} &
\includegraphics[trim={0 90 0 0},clip,width=.45\textwidth]{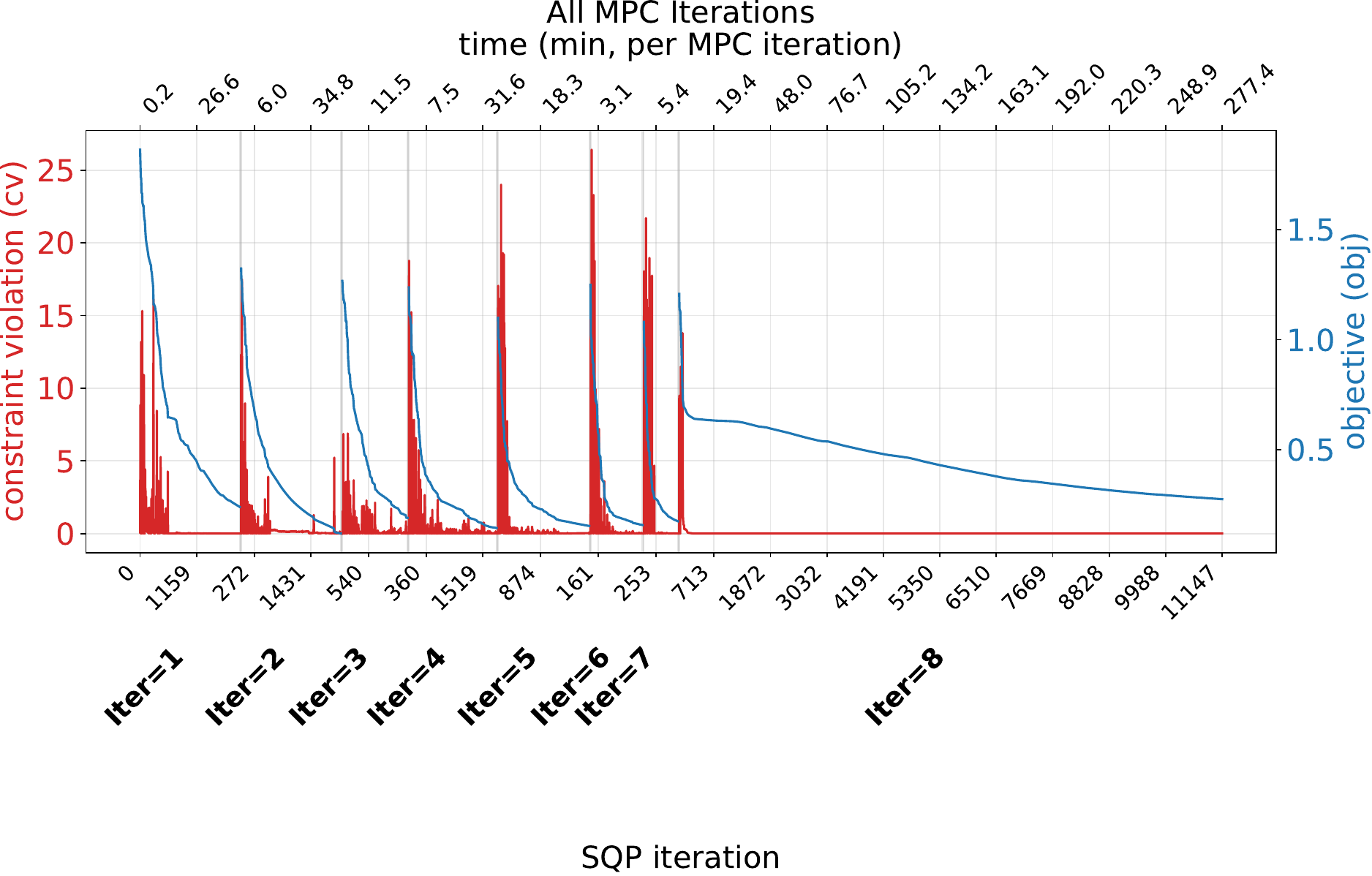} \\
\toprule
Ant-Terrain-3 & Push\\
\midrule
\includegraphics[trim={0 80 0 0},clip,width=.45\textwidth]{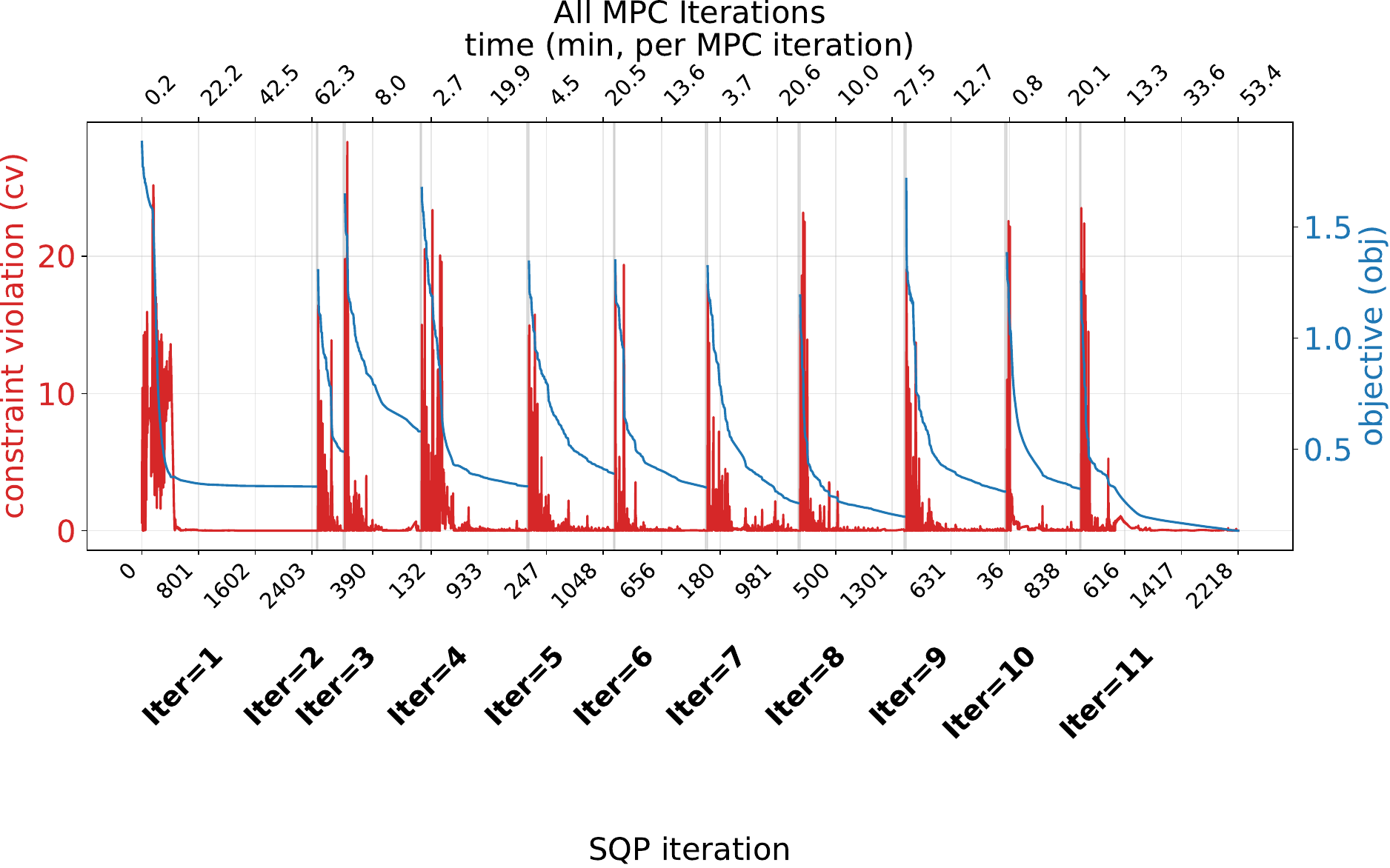} &
\includegraphics[trim={0 90 0 0},clip,width=.45\textwidth]{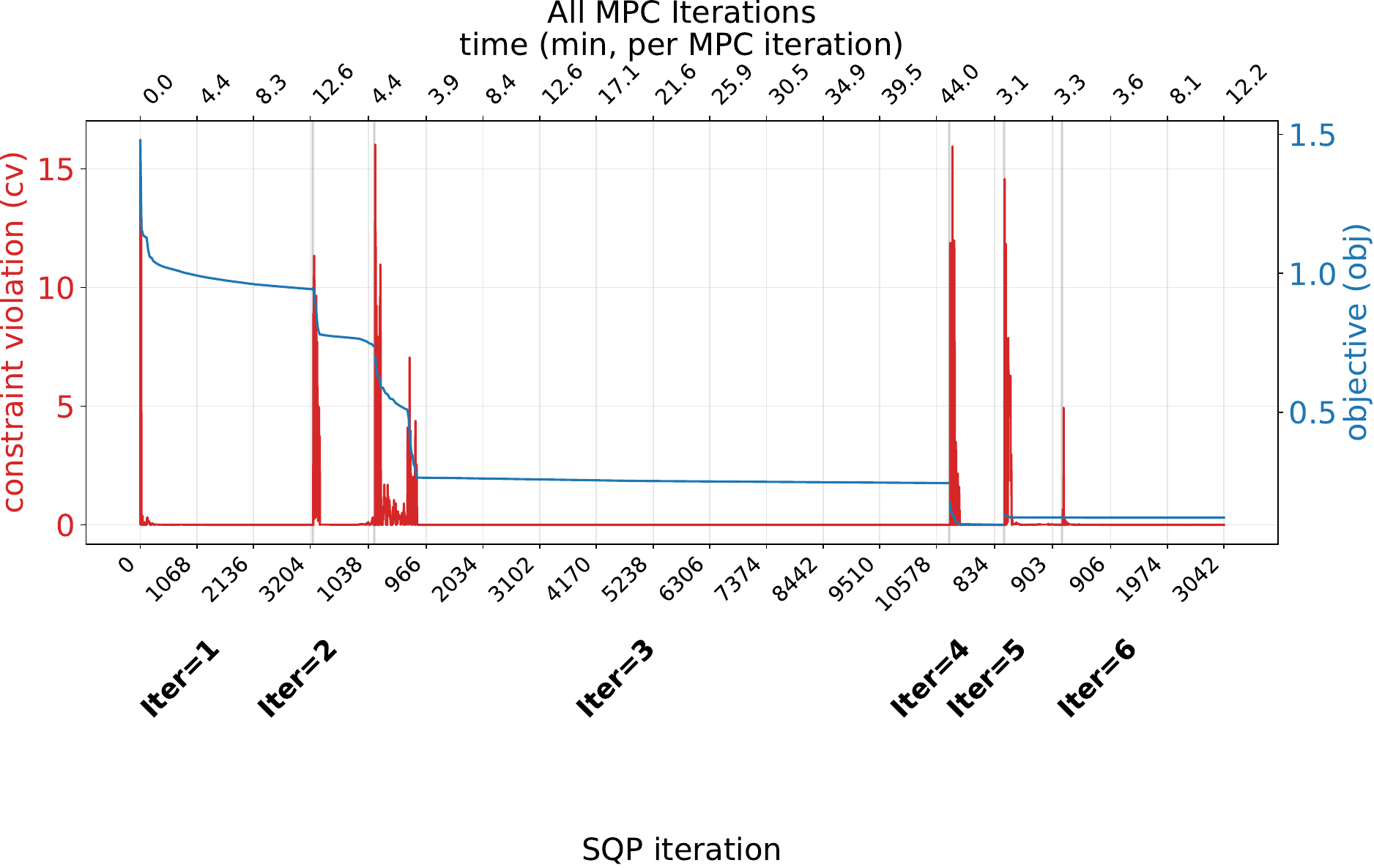} \\
\toprule
\multicolumn{2}{c}{Rotate}\\
\midrule
\multicolumn{2}{c}{\includegraphics[trim={0 90 0 0},clip,width=.45\textwidth]{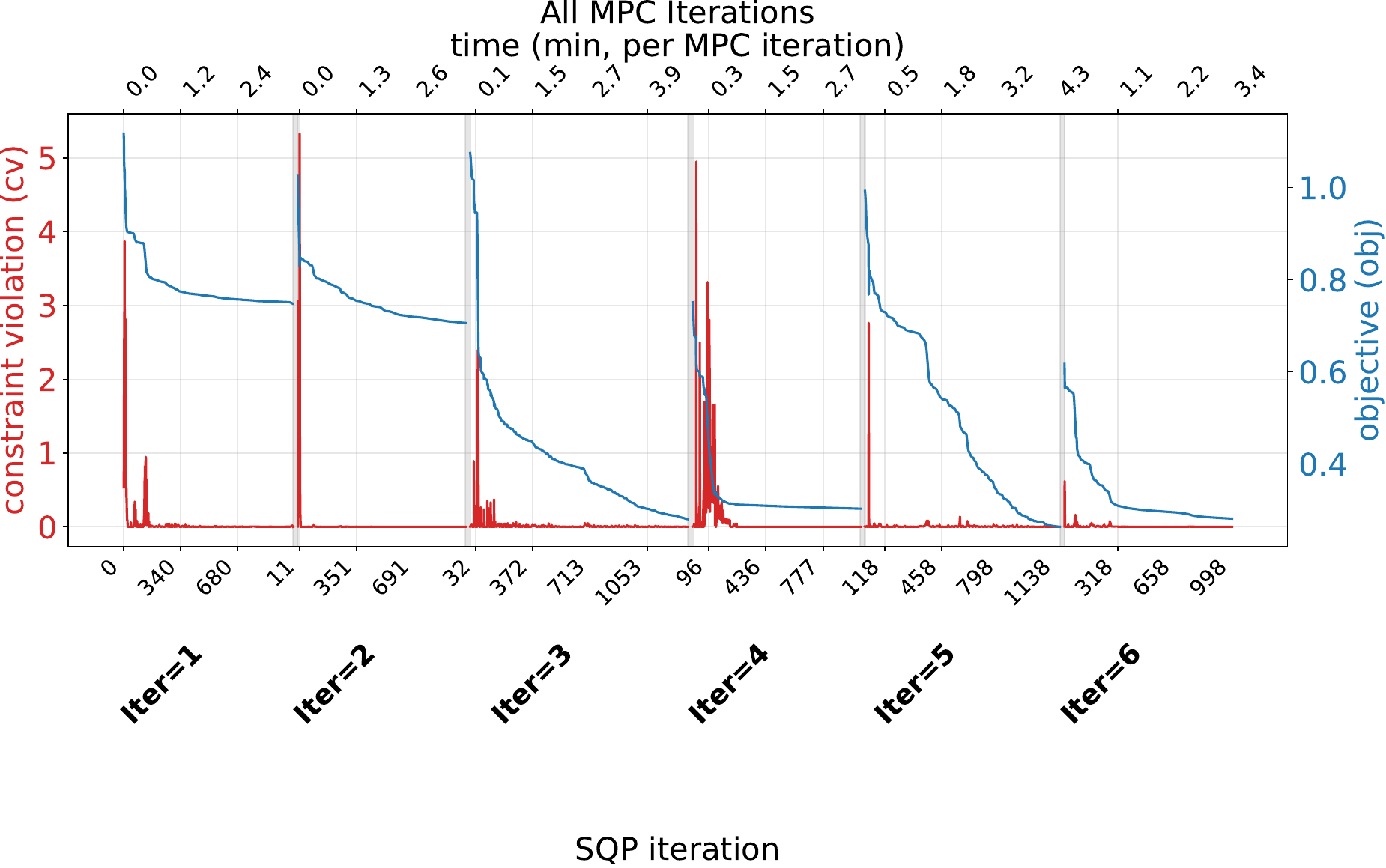}}\\
\bottomrule
\end{tabular}
\caption{\label{fig:performance}We visualize the physics-constraint violation $\|\StatRes^k\|_1$ and objective function at every trajectory optimization iteration. Each iteration generates an $H$ second sub-trajectory and we then extrapolate the last frame to build the initial guess for the next iteration. The extrapolated configuration vector is finite but need not satisfy the modified dynamics or original contact constraints; the constraint-satisfaction phase therefore produces the observed initial spike and subsequent decrease in violation.}
\end{figure*}

We summarize the residual and objective histories of our method in~\prettyref{fig:performance}. In all the scenarios, we observe consistent behavior in the solver's convergence history. From the rough initial guess, it quickly drives the physics-constraint residual down to our prescribed stopping threshold, after which it makes local trajectory adjustments that further reduce the objective function value. In~\prettyref{fig:cost-dof}, we plot the distribution of average wall-clock time per MPC solve over the horizon $H=2s$. Reaching the prescribed stopping threshold takes between $0.1$--$2\,\mathrm{hr}$ per MPC solve in these experiments. We also observe that the average cost does not increase significantly with the number of DOFs. Faster implementations and evaluation-complexity bounds for reaching the prescribed stopping test remain directions for future work.
\begin{figure}[t]
\centering
\includegraphics[width=.9\linewidth]{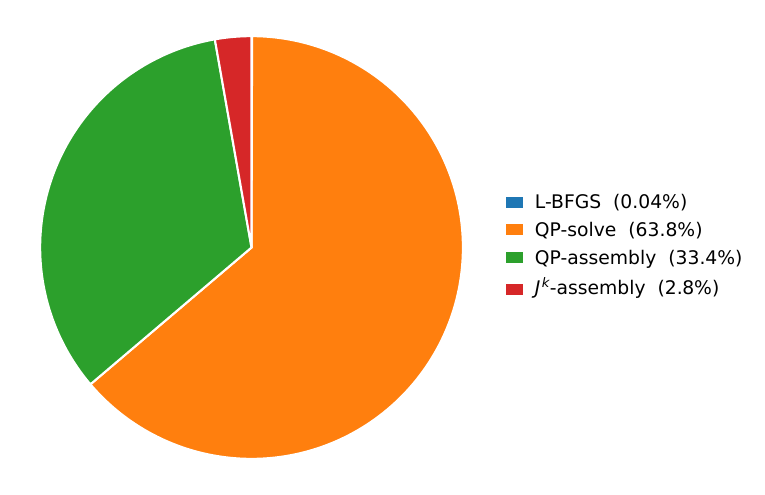}
\caption{\label{fig:perf-breakdown}Breakdown of the total wall-clock time of a representative MPC run (the Ant-Terrain-2 benchmark) into its four major components: the L-BFGS Hessian update, the active-set QP solve (QP-solve), the reduced-system assembly within each QP solve (QP-assembly), and the physics constraint Jacobian assembly ($J^k$-assembly).}
\end{figure}
Finally, to identify the computational bottleneck of our solver, we further profile a representative MPC run on the Ant-Terrain-2 benchmark and break down its total wall-clock time into the four major components, as reported in~\prettyref{fig:perf-breakdown}. The QP subproblem dominates the overall cost: the active-set QP solve alone accounts for $63.8\%$ of the runtime, and together with the reduced-system assembly it consumes over $97\%$ of the total. In comparison, the physics constraint Jacobian assembly ($J^k$-assembly) takes only $2.8\%$, and the L-BFGS Hessian update is negligible. This indicates that the efficiency of our method is primarily governed by the QP stage, and points to faster QP solvers and reduced-system assembly as the most promising directions for future acceleration.

\subsection{Compare with Off-the-shelf SQP}
\begin{figure*}[ht]
\centering
\setlength{\tabcolsep}{0pt}
\begin{tabular}{cc}
\toprule
Bipedal Locomotion & Ant Locomotion\\
\midrule
\includegraphics[width=.45\textwidth]{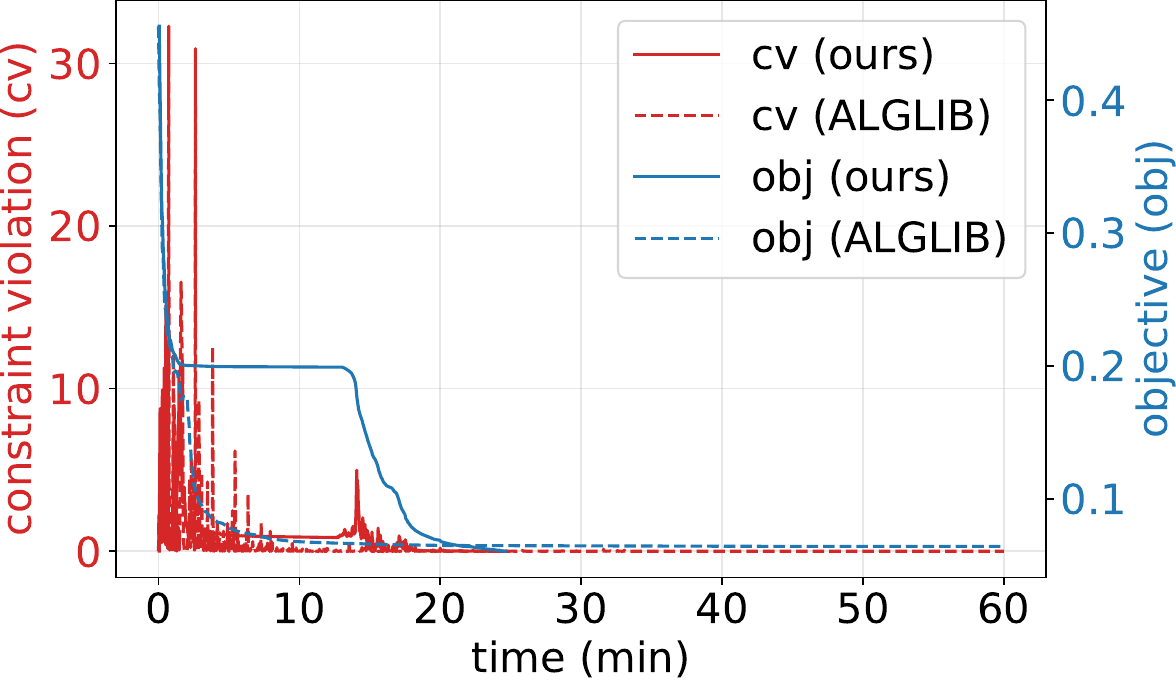} &
\includegraphics[width=.45\textwidth]{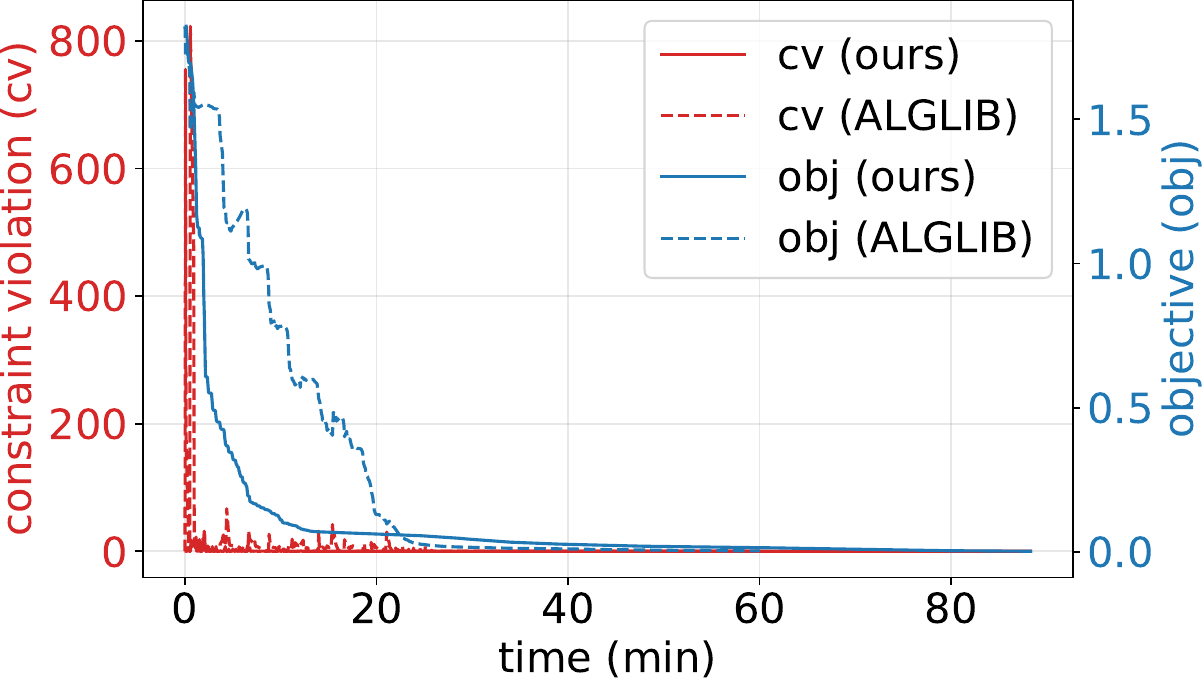} \\
\toprule
Ant-Terrain-1 & Ant-Terrain-2\\
\midrule
\includegraphics[width=.45\textwidth]{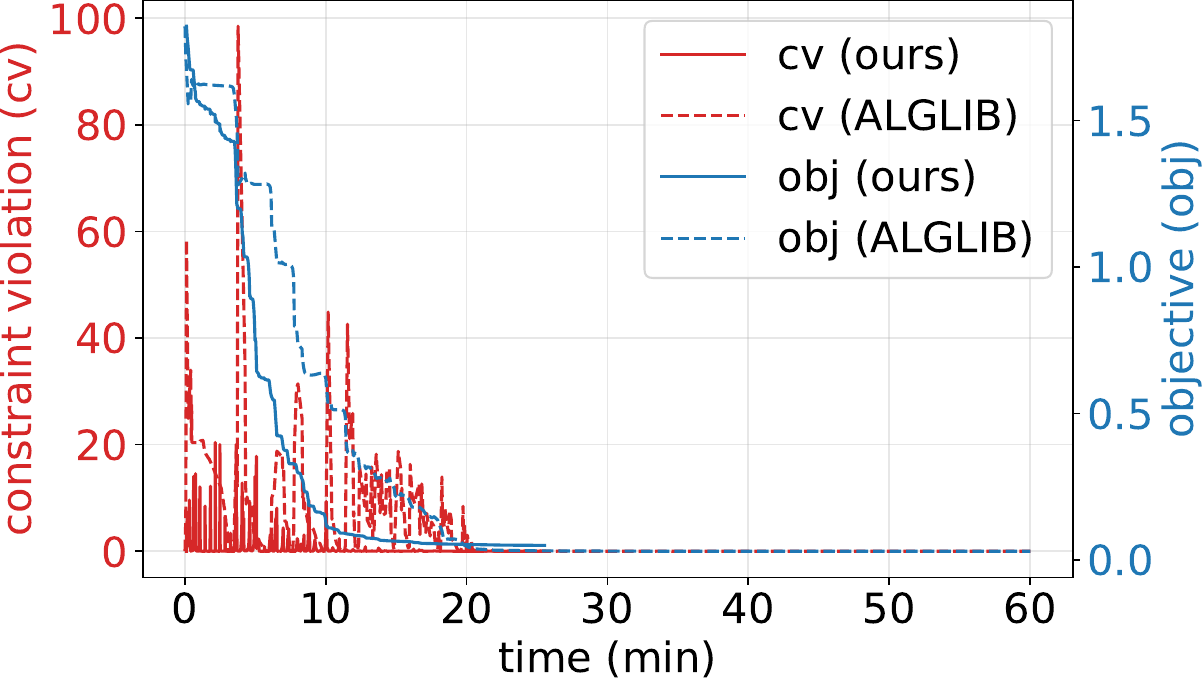} &
\includegraphics[width=.45\textwidth]{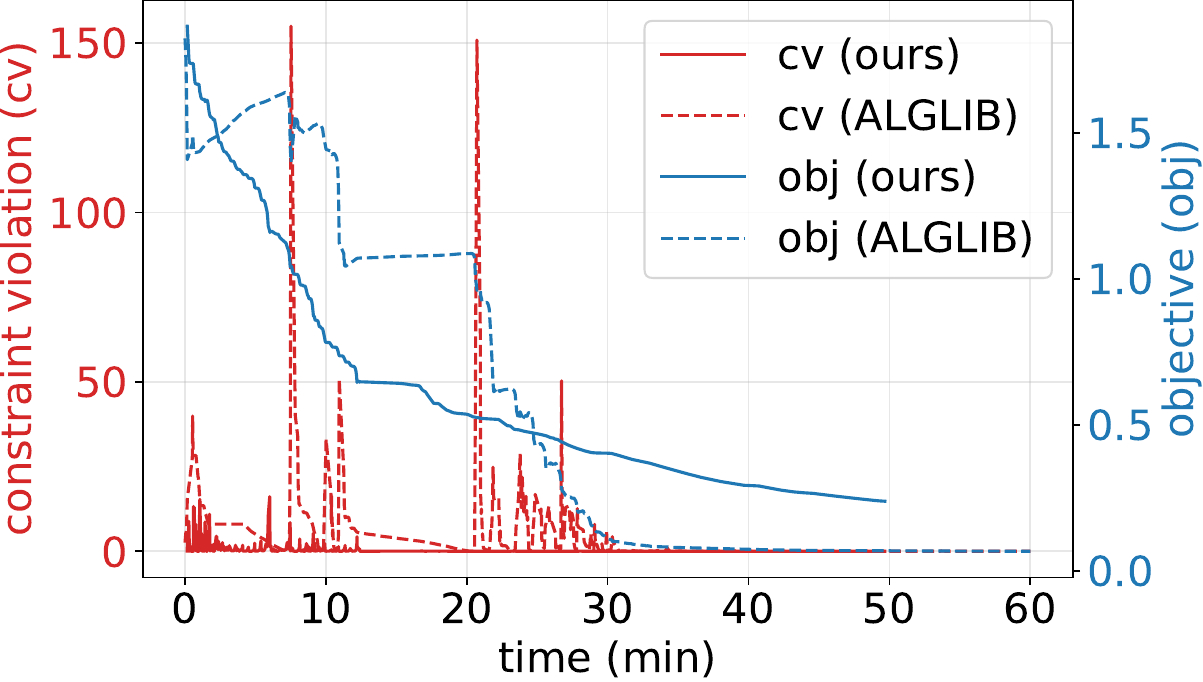} \\
\toprule
Ant-Terrain-3 & Push\\
\midrule
\includegraphics[width=.45\textwidth]{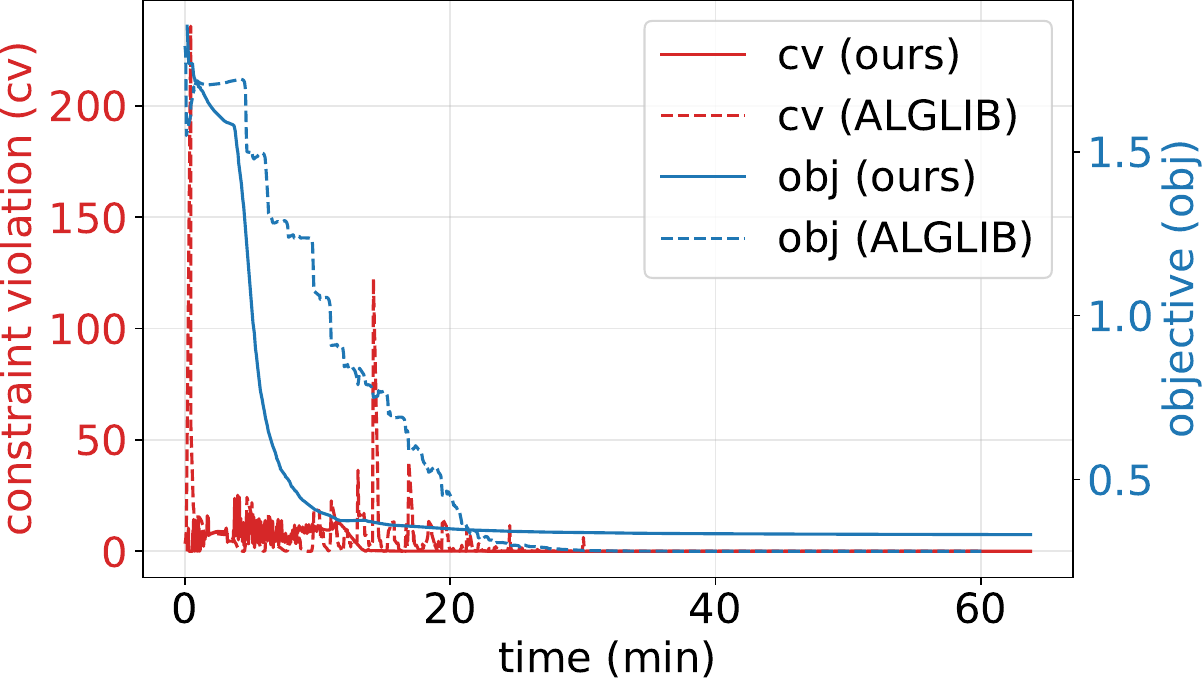} &
\includegraphics[width=.45\textwidth]{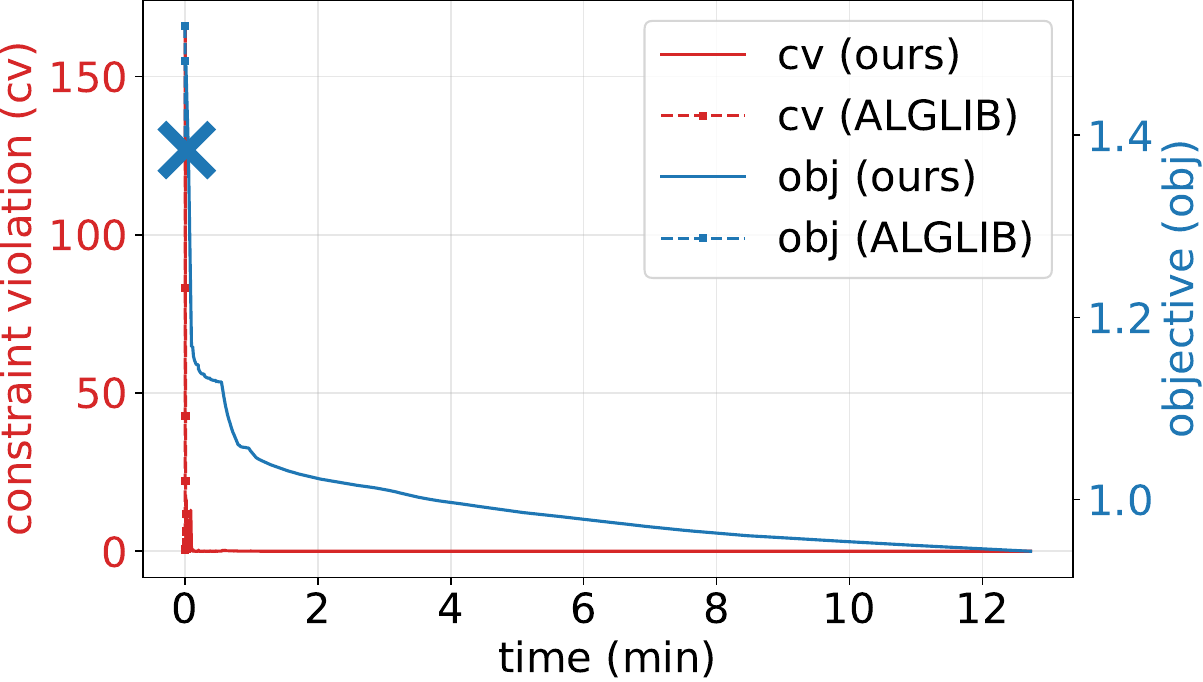} \\
\toprule
\multicolumn{2}{c}{Rotate}\\
\midrule
\multicolumn{2}{c}{\includegraphics[width=.45\textwidth]{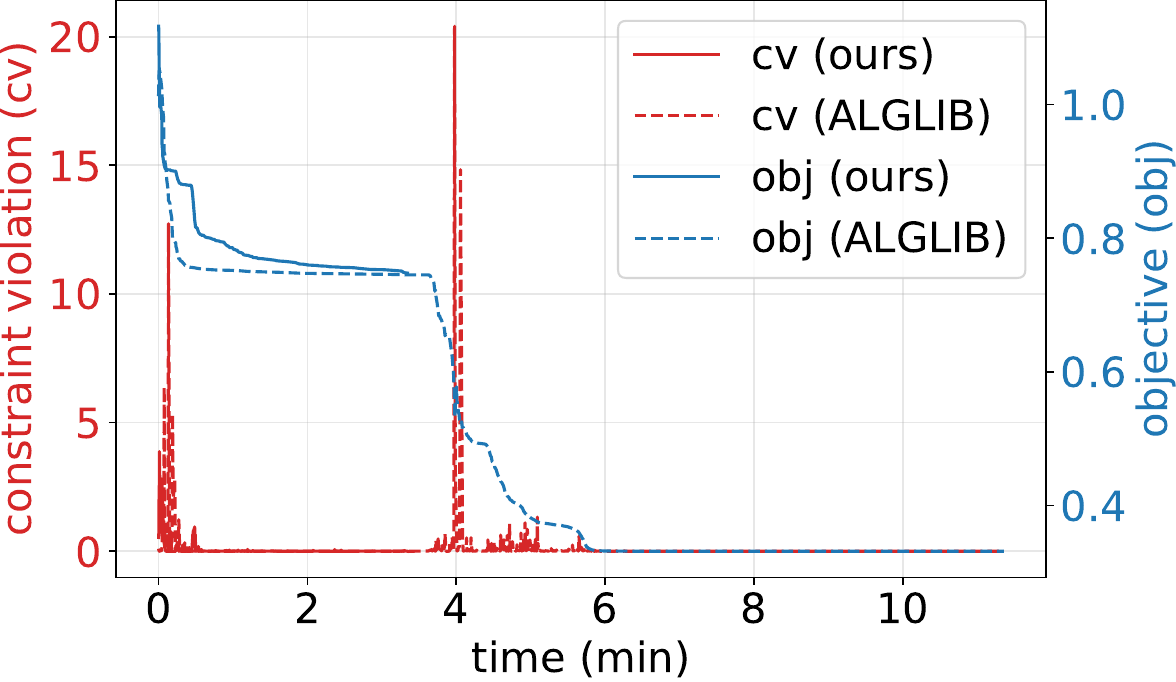}}\\
\bottomrule
\end{tabular}
\caption{\label{fig:alglib-compare}Wall-clock comparison between our customized SQP solver and the filter-SQP solver in ALGLIB on the first MPC solve of each benchmark. Both solvers start from the same initial trajectory. The red curves show constraint violation (CV, left axis), and the blue curves show objective value (right axis); solid curves denote our method and dashed curves denote ALGLIB. ALGLIB does not enforce a prescribed upper bound on constraint violation at its intermediate iterates and consequently exhibits frequent violation spikes. The blue cross in the Push benchmark marks the point at which ALGLIB fails. The comparison is limited to the first MPC solve because the solvers generally return different local solutions, which would give subsequent MPC windows different warm starts.}
\end{figure*}

We additionally compare our customized solver with the large-scale sparse filter-SQP algorithm in ALGLIB MinNLC, version 4.08.0~\cite{alglib408}. Both solvers are applied to the same CITO formulation, so this comparison isolates the effect of the solver. The comparison is restricted to the first MPC solve of each benchmark, for which both algorithms start from the same initial trajectory. Because the problem is nonconvex, the two algorithms can converge to different local solutions during this first solve. The next MPC window would then be initialized from different trajectories, so a subsequent comparison would conflate solver behavior with the effect of different warm starts and would no longer be meaningful. The first-window residual and objective histories are reported in~\prettyref{fig:alglib-compare}.
Unlike our safeguarded line search, ALGLIB does not provide a mechanism for imposing a prescribed upper bound on the nonlinear-constraint violation of every accepted iterate; this type of iterate-wise cap is also absent from most off-the-shelf SQP implementations. Consequently, the ALGLIB histories exhibit numerous constraint-violation spikes. In the Push benchmark, one such excursion leads to solver failure, marked by the blue cross in~\prettyref{fig:alglib-compare}. Recall that no singularity- or Hamiltonian-triggered subdivision occurs in our experiments. Together, these observations indicate that, although subdivision is inactive on these instances, explicitly bounding the constraint violation during optimization is important for robustness.
ALGLIB reaches a lower final objective in some benchmarks, whereas our solver reaches a lower one in others. This variation is expected because neither method guarantees global optimality. In wall-clock time, ALGLIB is generally slower at reducing the objective. Our implementation exploits the block-tridiagonal structure of the trajectory QP with a specialized linear-system solver, whereas ALGLIB uses a more general sparse matrix solver; the latter supports a broader class of nonlinear programs but cannot exploit this problem-specific structure as efficiently.

\subsection{Compare with Alternative CITO Formulations}
\begin{figure}[t]
\centering
\setlength{\tabcolsep}{0pt}
\begin{tabular}{c}
\toprule
Bipedal Locomotion\\
\midrule
\includegraphics[width=\linewidth]{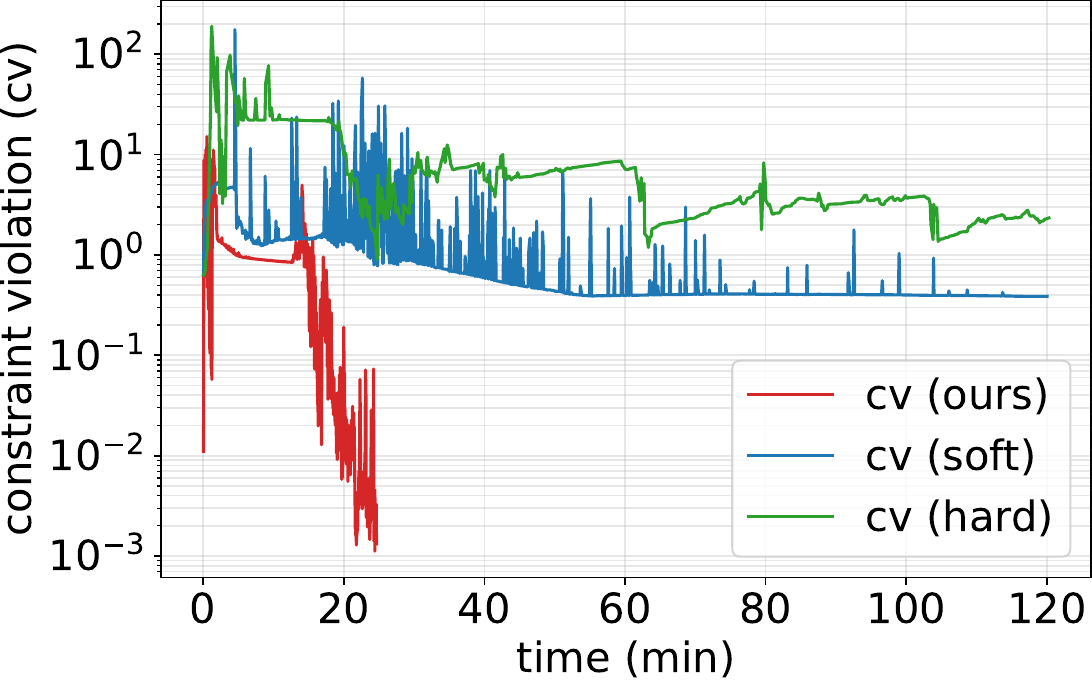}\\
\toprule
Ant Locomotion\\
\midrule
\includegraphics[width=\linewidth]{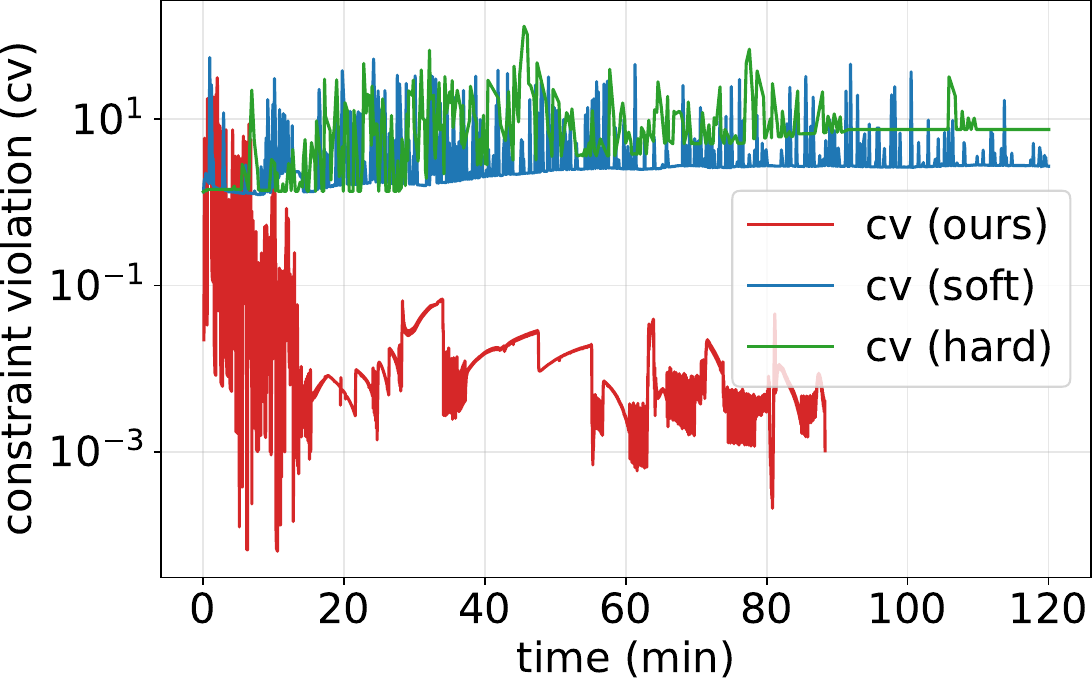}\\
\bottomrule
\end{tabular}
\caption{\label{fig:direct-compare}Reduction of the physics-constraint violation over wall-clock time for the three contact formulations on the bipedal locomotion benchmark (top) and the ant locomotion benchmark (bottom). All three formulations solve the same first MPC window of the same benchmark, start from the identical trivial initialization with the robot standing still, and are run for a budget of $2\,\mathrm{hr}$; the constraint violation is plotted on a logarithmic axis. The curve labeled ``cv (ours)'' is the penalized, force-parameterized formulation solved by the safeguarded line-search SQP method of~\prettyref{sec:global}, ``cv (soft)'' is the soft-penalty contact-invariant formulation of~\cite{mordatch2012discovery}, and ``cv (hard)'' is the hard-complementarity MPCC formulation of~\cite{posa2014direct}.} 
\end{figure}

The comparison above varies the solver while holding our formulation fixed. We now do the converse, holding the benchmark and the initialization fixed while varying the way the contact model itself is encoded, and compare against the two dominant choices in the CITO literature. All alternative formulations are solved with the same ALGLIB SQP solver used in the previous subsection, so this comparison isolates the effect of the formulation. The first is the direct method of~\cite{posa2014direct}, which transcribes the contact model as hard complementarity constraints and solves the resulting MPCC with SQP; we refer to it as the hard-constraint formulation. The second is contact-invariant optimization~\cite{mordatch2012discovery}, which replaces the complementarity conditions with soft penalty terms, a contact-invariant cost together with a physics-violation cost, and minimizes their weighted sum; we refer to it as the soft-penalty formulation. As stated in~\cite{mordatch2012discovery}, physical consistency is enforced through a soft cost rather than a hard constraint precisely so that the optimizer is allowed to leave the constraint manifold during the search. We evaluate all three formulations on the bipedal and ant locomotion benchmarks and focus on how effectively each of them reduces the physics-constraint violation. All three solve the same first MPC window from the same trivial initialization used throughout our experiments, in which the robot stands still, and the two baselines are given a wall-clock budget of $2\,\mathrm{hr}$, considerably more than our solver needs to reach its stopping threshold on these two benchmarks. The resulting violation histories are reported in~\prettyref{fig:direct-compare}.
The hard-constraint formulation does not find a feasible trajectory on either benchmark: its constraint violation remains high throughout the entire budget. This outcome is consistent with what is claimed in~\cite{posa2014direct}, where the reported results are obtained from case-by-case initial guesses, such as a user-specified intermediate grasp state interpolated to the terminal state for manipulation, or an initial trajectory already carrying a nominal mode sequence for gait synthesis. Imposing complementarity as a hard constraint requires every accepted iterate to remain close to the constraint manifold, so the search has to be started sufficiently close to a feasible trajectory; a standing-still initialization does not meet that requirement.
The soft-penalty formulation behaves differently and attains a smaller constraint violation than the hard-constraint formulation over time on both benchmarks. Because the penalty admits temporary excursions away from the constraint manifold, the optimizer can explore contact-mode changes that are unreachable for the hard-constraint formulation from the same starting point. This is the mechanism through which~\cite{mordatch2012discovery} reports behaviors synthesized from uninformative initializations such as ours. A smaller violation does not by itself imply that a feasible or physically meaningful motion has been found. On the bipedal benchmark, the soft-penalty result is visually plausible; a numerical residual remains, which is expected and acceptable for a soft-penalty formulation that does not claim exact constraint satisfaction. On the ant benchmark, by contrast, the returned motion exhibits non-physical striding of the robot feet against the ground. Our formulation reaches the prescribed stopping threshold on both benchmarks. It retains the ability of the soft-penalty formulation to start from a trivial initialization, since our penalized surrogate likewise tolerates intermediate iterates that are not exactly physical, while the safeguarded line search additionally caps the violation of every accepted iterate and~\prettyref{thm:sqp-convergence} certifies finite termination at the prescribed residual level.

\section{\label{sec:conclusion}Conclusion and Future Work}
This paper developed, to our knowledge, the first global-convergence analysis of a specific SQP framework tailored to contact-implicit trajectory optimization. We posed CITO as~\prettyref{pb:CITO} over a broad class of constraint-rich dynamic systems with implicit physics constraints, and instantiated the class using articulated-body dynamics with contact between arbitrary unions of convex polyhedra. The analyzed method is the modified line-search SQP scheme of~\prettyref{alg:sqp} wrapped by the adaptive top-level solver of~\prettyref{alg:CITO}. These algorithms optimize the modified full-trajectory~\prettyref{pb:CITO-M}, whose physics constraints are the stationarity equations of the mollified per-step surrogate~\prettyref{pb:PS-CRDS-M}; they do not apply SQP directly to the original MPCC~\prettyref{pb:CITO}. The approximate force-parameterized certificate for~\prettyref{pb:CITO} is obtained through the subsequent feasibility bridge and stationarity transfer. Under the stated assumptions, for fixed startup data with $\theta_0$ satisfying the strict-interior condition and from any finite discretized optimized trajectory guess whose controls lie in the compact control set, this specified algorithm terminates after finitely many iterations and subdivisions at an $\EpsStop$-KKT point of~\prettyref{pb:CITO-M}, which the bridge and transfer results convert into an approximately feasible, force-parameterized stationary certificate for~\prettyref{pb:CITO} without assuming a constraint qualification a priori. This is global convergence in the numerical-optimization sense, not global optimality of the returned trajectory. The mechanism is that LICQ is derived rather than posited: maximal-coordinate lifting makes the discrete kinetic energy strongly convex with curvature of order $M/\delta_i^2$, the penalty relaxation together with the finite-horizon Hamiltonian bound (\prettyref{lem:H-bound}) caps the potential curvature, and a sufficiently small timestep---enforced by the adaptive subdivision of~\prettyref{alg:CITO}---makes the positive-definite mass matrix dominate, conditioning the stationarity Jacobian (\prettyref{lem:jacobian-lower-bound},~\ref{lem:jacobian-upper-bound}). The same analysis keeps both the merit-penalty parameter and the optimized trajectory norm-bounded throughout the optimization. These are exact-arithmetic boundedness conclusions under the stated fixed-data assumptions, not a roundoff- or backward-error guarantee for the floating-point implementation. The guarantee comes in two forms: for smooth damping the transfer attains the unit-objective, penalty-induced force-parameterized condition of~\prettyref{thm:stationarity-transfer}---``almost-KKT'' only in the reduced trajectory/control sense, not KKT of the full independent-force MPCC (\prettyref{def:eps-stationary}, objective multiplier $\mu_0=1$), while for genuinely nonsmooth frictional contact it weakens to the homogeneous penalty-induced force-parameterized Clarke--Fritz--John certificate of~\prettyref{thm:nonsmooth-transfer}. Here the Clarke normals belong to the lifted damping graphs, and the normalized objective multiplier $\mu_0$ may vanish; in that abnormal case the certificate contains no first-order objective information. In both cases the force selection is inherited from the penalty gradient; the certificate therefore establishes physical residual, approximate feasibility, relaxed complementarity, and trajectory/control stationarity for that selected force graph, but not independent optimality of the force variables.

It is worth being precise about what this result does, and does not, establish. First, global convergence here is in the numerical-optimization sense---finite termination of the specified algorithm, from any finite discretized optimized trajectory guess under the fixed pre-horizon and other admissibility assumptions, at an approximate stationarity certificate---and not global optimality: gradient-based information certifies only local stationarity properties of the returned trajectory. Second, and more subtly, the result is existential over algorithms rather than over solvers in general. What we prove is that one particular modified SQP, equipped with the three problem-specific safeguards (adaptive subdivision, a constraint-satisfaction pre-phase, and a safeguarded line search), converges globally for CITO; it does not imply that an off-the-shelf or generic SQP method enjoys the same guarantee. This distinction is not academic. As the representative infeasibility result~\prettyref{thm:bcn-infeasibility} shows, absent a CQ, a general exact-penalty SQP may converge only to an infeasible stationary point of the constraint-violation measure $\|F(\cdot)\|$; the classical failure examples of interior-point and SQP methods on problems lacking a CQ~\cite{wachter2000failure} and the delicate local behavior of SQP applied directly to complementarity constraints~\cite{fletcher2006local} illustrate the same hazard from the MPCC side. The three algorithmic safeguards, together with the derived LICQ modulus, are precisely what convert this potential failure into the approximate feasibility and stationarity guarantee established for the specified method under the stated conditions. Finally, feasibility is only approximate: the certified point lies within an $\EpsFJAll$-band of the exact constraint manifold, with the band tied to the barrier parameter $\EpsFJ$. The uniform-boundedness results hold for fixed positive primitive regularization parameters and make no uniform claim in the stiff or vanishing-smoothing limits. A further limitation is the strict-interior requirement at the fixed time-zero configuration $\theta_0$: the present extended-real barrier needs positive inequality slack there so that the first optimized penalty-induced force and semi-implicit damping load are finite. Allowing boundary contact at node $0$, with $h^i(\theta_0)=0$, would require a boundary-finite barrier, an independently initialized first-step contact force, or a separate startup contact solve and transfer argument. The articulated-body certificate also uses exact mathematical partial minimization over the auxiliary separating-plane variables and velocities: the present theorem does not propagate finite-precision inner-solve error into its final stationarity tolerance. These limitations indicate the following directions for potential future work.

\begin{enumerate}[label=(\alph*),ref=(\alph*),leftmargin=*]
\item\label{it:fw-complexity} A complexity theory to accompany the convergence theory. Our results establish finite termination but provide no bound on the number of iterations or subdivisions required. The machinery for an upper bound seems within reach: an evaluation-complexity argument in the style of adaptive-regularization and log-barrier analyses~\cite{cartis2014complexity,hinder2024worst} would translate the per-successful-step decrease of the $\ell_1$ merit function into an iteration count for reaching an $\EpsStop$-KKT point, with constants governed by the derived LICQ modulus $\SigmaMinLICQ(\StatJac)$ and hence by the curvature ratio $M/\delta_i^2$; the resulting bound would depend explicitly on the timestep floor $\delta_{\min}$ and the horizon $H$. The far more interesting question is how such a bound scales with the horizon, and here we conjecture that the dependence is exponential. The intuition is combinatorial: a contact-rich trajectory of horizon $H$ with a fixed number of candidate contacts admits on the order of $2^{(\#\text{contacts}\cdot|\mathcal I|)}$ discrete contact-mode sequences, and selecting among them inherits the hardness of the underlying linear-complementarity and MPCC subproblems~\cite{chung1989np,scheel2000mathematical} as well as the PSPACE-hardness of motion planning through contact~\cite{reif1979complexity}; the shortest-path-through-convex-sets and mixed-integer footstep formulations~\cite{marcucci2024shortest,deits2014footstep} make the same mode-selection combinatorics explicit. Establishing a matching lower bound is the harder and more valuable half: it would require a hard-instance construction---a family of CITO problems in which any method of our class is provably forced to explore $\Omega(2^{cH})$ contact-mode branches---in the spirit of the lower-bound constructions for first-order optimization~\cite{carmon2020lower}. Together, matching upper and lower bounds would characterize the intrinsic difficulty of contact planning and delineate when a globally convergent method can also be efficient.
\item\label{it:fw-nonsmooth-sqp} Closing the gap between the smooth and nonsmooth guarantees. In the smooth case the transfer attains the unit-objective form $\mu_0=1$, which is ``almost-KKT'' only in the reduced trajectory/control sense after the penalty-induced force maps have been substituted. In the nonsmooth case the complete certificate is homogeneous and $\mu_0$ may vanish because the physics multiplier is not proved uniformly bounded as $\EpsFri\to0$---the assembled LICQ modulus degrades through the off-diagonal dissipation coupling of the stationarity Jacobian. A promising remedy is to regularize the damping by its Moreau--Yosida envelope~\cite{rockafellar1998variational,lemarechal1997practical}, which may restore a uniform multiplier bound, prove $\mu_0>0$ uniformly, and thereby recover a unit-objective force-parameterized conclusion. The obstacle is that the Moreau envelope is only $C^{1,1}$ rather than $C^2$, so its second-order structure falls outside the classical twice-differentiable SQP theory used here. Extending the analysis would require a nonsmooth SQP framework built on generalized (semismooth) derivatives~\cite{qi1993nonsmooth,facchinei2003finite,curtis2012sequential,clarke1983optimization}, an interesting undertaking in its own right that would also broaden the applicability of the convergence result beyond the smoothing route we take.
\item\label{it:fw-inexact-reduction} Extending the transfer to inexact partial minimization. A finite-precision implementation only approximates the auxiliary minimizers defining $V_R$ and $\Damp_R^V$, so its auxiliary stationarity and reduced-value envelope identities carry nonzero errors that are not represented by the present $\EpsStop$ or $\EpsStopNS$ bounds. An inexact-reduction theory should quantify these residuals, propagate the resulting envelope and Schur-lift errors through the full certificate, and expose their contribution to the final stationarity tolerance.
\item\label{it:fw-force-param} Sharpening the force parameterization at the heart of the stationarity transfer. Throughout, the constraint forces $\lambda_i^e,\lambda_i^i$ are not free decision variables but prescribed configuration-dependent maps induced by the penalty gradient of $V(\cdot,\EpsFJ)$; this penalty-induced graph restriction is what collapses the stationarity transfer of~\prettyref{pb:CITO} onto the reduced trajectory space that~\prettyref{alg:CITO} optimizes, but it does not test independent force directions, so the certified forces are consistent with the penalty model rather than independently optimal. There are two logically separate strengthening directions. First, replacing the Clarke normals to the damping graphs in~\prettyref{def:eps-stationary-NS} by limiting (Mordukhovich) damping-graph normals would sharpen only the nonsmooth damping component of the present force-parameterized certificate. By itself, this would not establish standard MPCC M-stationarity, because the equality and inequality/contact forces would still be prescribed maps and the normal cone to the complementarity set would still not be tested in the independent force directions. Second, promoting those forces to genuine decision variables and deriving the appropriate limiting-normal conditions for the full complementarity system could yield standard MPCC M-stationarity; stronger regularity or a tightened-program argument would be needed to reach S-stationarity~\cite{scheel2000mathematical,ye2005necessary,outrata1999optimality,hoheisel2013theoretical,flegel2003fritz}. Pursuing these directions would connect the constructive, algorithmically derived guarantees of this paper to the sharpest stationarity notions available for MPCCs.
\end{enumerate}

\bibliographystyle{SageH}
\bibliography{references}

@phdthesis{ashi2008numerical,
  title  = {Numerical methods for stiff systems},
  author = {Ashi, Hala},
  year   = {2008},
  school = {University of Nottingham}
}

@inproceedings{todorov2012mujoco,
  title        = {Mujoco: A physics engine for model-based control},
  author       = {Todorov, Emanuel and Erez, Tom and Tassa, Yuval},
  booktitle    = {2012 IEEE/RSJ international conference on intelligent robots and systems},
  pages        = {5026--5033},
  year         = {2012},
  organization = {IEEE}
}

@inproceedings{1641984,
  author    = {Yamane, K. and Nakamura, Y.},
  booktitle = {Proceedings 2006 IEEE International Conference on Robotics and Automation, 2006. ICRA 2006.},
  title     = {Stable penalty-based model of frictional contacts},
  year      = {2006},
  volume    = {},
  number    = {},
  pages     = {1904-1909},
  doi       = {10.1109/ROBOT.2006.1641984}
}

@article{harmon2009asynchronous,
  author     = {Harmon, David and Vouga, Etienne and Smith, Breannan and Tamstorf, Rasmus and Grinspun, Eitan},
  title      = {Asynchronous contact mechanics},
  year       = {2012},
  issue_date = {April 2012},
  publisher  = {Association for Computing Machinery},
  address    = {New York, NY, USA},
  volume     = {55},
  number     = {4},
  issn       = {0001-0782},
  url        = {https://doi.org/10.1145/2133806.2133828},
  doi        = {10.1145/2133806.2133828},
  journal    = {Commun. ACM},
  month      = {apr},
  pages      = {102–109},
  numpages   = {8}
}

@inproceedings{dong2016motion,
  title     = {Motion planning as probabilistic inference using Gaussian processes and factor graphs.},
  author    = {Dong, Jing and Mukadam, Mustafa and Dellaert, Frank and Boots, Byron},
  booktitle = {Robotics: Science and Systems},
  volume    = {12},
  number    = {4},
  year      = {2016}
}

@article{schulman2014motion,
  title     = {Motion planning with sequential convex optimization and convex collision checking},
  author    = {Schulman, John and Duan, Yan and Ho, Jonathan and Lee, Alex and Awwal, Ibrahim and Bradlow, Henry and Pan, Jia and Patil, Sachin and Goldberg, Ken and Abbeel, Pieter},
  journal   = {The International Journal of Robotics Research},
  volume    = {33},
  number    = {9},
  pages     = {1251--1270},
  year      = {2014},
  publisher = {SAGE Publications Sage UK: London, England}
}

@article{posa2014direct,
  title     = {A direct method for trajectory optimization of rigid bodies through contact},
  author    = {Posa, Michael and Cantu, Cecilia and Tedrake, Russ},
  journal   = {The International Journal of Robotics Research},
  volume    = {33},
  number    = {1},
  pages     = {69--81},
  year      = {2014},
  publisher = {Sage Publications Sage UK: London, England}
}

@article{bertsekas1997nonlinear,
  title     = {Nonlinear programming},
  author    = {Bertsekas, Dimitri P},
  journal   = {Journal of the Operational Research Society},
  volume    = {48},
  number    = {3},
  pages     = {334--334},
  year      = {1997},
  publisher = {Taylor \& Francis}
}

@inproceedings{brudigam2021linear,
  title        = {Linear-quadratic optimal control in maximal coordinates},
  author       = {Br{\"u}digam, Jan and Manchester, Zachary},
  booktitle    = {2021 IEEE International Conference on Robotics and Automation (ICRA)},
  pages        = {9775--9781},
  year         = {2021},
  organization = {IEEE}
}

@article{biegler2009large,
  title     = {Large-scale nonlinear programming using IPOPT: An integrating framework for enterprise-wide dynamic optimization},
  author    = {Biegler, Lorenz T and Zavala, Victor M},
  journal   = {Computers \& Chemical Engineering},
  volume    = {33},
  number    = {3},
  pages     = {575--582},
  year      = {2009},
  publisher = {Elsevier}
}

@article{gill2005snopt,
  title     = {SNOPT: An SQP algorithm for large-scale constrained optimization},
  author    = {Gill, Philip E and Murray, Walter and Saunders, Michael A},
  journal   = {SIAM review},
  volume    = {47},
  number    = {1},
  pages     = {99--131},
  year      = {2005},
  publisher = {SIAM}
}

@article{anitescu2000solving,
  title   = {On solving mathematical programs with complementarity constraints as nonlinear programs},
  author  = {Anitescu, Mihai},
  journal = {Preprint ANL/MCS-P864-1200, Argonne National Laboratory, Argonne, IL},
  volume  = {3},
  year    = {2000}
}

@article{dai2019global,
  title     = {Global inverse kinematics via mixed-integer convex optimization},
  author    = {Dai, Hongkai and Izatt, Gregory and Tedrake, Russ},
  journal   = {The International Journal of Robotics Research},
  volume    = {38},
  number    = {12-13},
  pages     = {1420--1441},
  year      = {2019},
  publisher = {SAGE Publications Sage UK: London, England}
}

@inproceedings{suh2022differentiable,
  title        = {Do differentiable simulators give better policy gradients?},
  author       = {Suh, Hyung Ju and Simchowitz, Max and Zhang, Kaiqing and Tedrake, Russ},
  booktitle    = {International Conference on Machine Learning},
  pages        = {20668--20696},
  year         = {2022},
  organization = {PMLR}
}

@article{winkler2018gait,
  title     = {Gait and trajectory optimization for legged systems through phase-based end-effector parameterization},
  author    = {Winkler, Alexander W and Bellicoso, C Dario and Hutter, Marco and Buchli, Jonas},
  journal   = {IEEE Robotics and Automation Letters},
  volume    = {3},
  number    = {3},
  pages     = {1560--1567},
  year      = {2018},
  publisher = {IEEE}
}

@article{manchester2019contact,
  title     = {Contact-implicit trajectory optimization using variational integrators},
  author    = {Manchester, Zachary and Doshi, Neel and Wood, Robert J and Kuindersma, Scott},
  journal   = {The International Journal of Robotics Research},
  volume    = {38},
  number    = {12-13},
  pages     = {1463--1476},
  year      = {2019},
  publisher = {SAGE Publications Sage UK: London, England}
}

@article{solodov2009global,
  title     = {Global convergence of an SQP method without boundedness assumptions on any of the iterative sequences},
  author    = {Solodov, Mikhail V},
  journal   = {Mathematical programming},
  volume    = {118},
  number    = {1},
  pages     = {1--12},
  year      = {2009},
  publisher = {Springer}
}

@book{brogliato1999nonsmooth,
  title     = {Nonsmooth mechanics},
  author    = {Brogliato, Bernard and Brogliato, B},
  volume    = {3},
  year      = {1999},
  publisher = {Springer}
}

@article{zhang2023provably,
  title   = {Provably Robust Semi-Infinite Program Under Collision Constraints via Subdivision},
  author  = {Zhang, Duo and Liang, Chen and Gao, Xifeng and Wu, Kui and Pan, Zherong},
  journal = {arXiv preprint arXiv:2302.01135},
  year    = {2023}
}

@article{6710113,
  author   = {Escande, Adrien and Miossec, Sylvain and Benallegue, Mehdi and Kheddar, Abderrahmane},
  journal  = {IEEE Transactions on Robotics},
  title    = {A Strictly Convex Hull for Computing Proximity Distances With Continuous Gradients},
  year     = {2014},
  volume   = {30},
  number   = {3},
  pages    = {666-678},
  doi      = {10.1109/TRO.2013.2296332}
}

@inproceedings{xie2017differential,
  title        = {Differential dynamic programming with nonlinear constraints},
  author       = {Xie, Zhaoming and Liu, C Karen and Hauser, Kris},
  booktitle    = {2017 IEEE International Conference on Robotics and Automation (ICRA)},
  pages        = {695--702},
  year         = {2017},
  organization = {IEEE}
}

@article{chatzinikolaidis2020contact,
  title     = {Contact-implicit trajectory optimization using an analytically solvable contact model for locomotion on variable ground},
  author    = {Chatzinikolaidis, Iordanis and You, Yangwei and Li, Zhibin},
  journal   = {IEEE Robotics and Automation Letters},
  volume    = {5},
  number    = {4},
  pages     = {6357--6364},
  year      = {2020},
  publisher = {IEEE}
}

@inproceedings{landry2019bilevel,
  title        = {Bilevel optimization for planning through contact: A semidirect method},
  author       = {Landry, Benoit and Lorenzetti, Joseph and Manchester, Zachary and Pavone, Marco},
  booktitle    = {The International Symposium of Robotics Research},
  pages        = {789--804},
  year         = {2019},
  organization = {Springer}
}

@article{stewart1996implicit,
  title     = {An implicit time-stepping scheme for rigid body dynamics with inelastic collisions and coulomb friction},
  author    = {Stewart, David E and Trinkle, Jeffrey C},
  journal   = {International Journal for Numerical Methods in Engineering},
  volume    = {39},
  number    = {15},
  pages     = {2673--2691},
  year      = {1996},
  publisher = {Wiley Online Library}
}

@inproceedings{manchester2020variational,
  title        = {Variational contact-implicit trajectory optimization},
  author       = {Manchester, Zachary and Kuindersma, Scott},
  booktitle    = {Robotics Research: The 18th International Symposium ISRR},
  pages        = {985--1000},
  year         = {2020},
  organization = {Springer}
}

@inproceedings{tassa2012synthesis,
  title        = {Synthesis and stabilization of complex behaviors through online trajectory optimization},
  author       = {Tassa, Yuval and Erez, Tom and Todorov, Emanuel},
  booktitle    = {2012 IEEE/RSJ International Conference on Intelligent Robots and Systems},
  pages        = {4906--4913},
  year         = {2012},
  organization = {IEEE}
}

@inproceedings{todorov2011convex,
  title        = {A convex, smooth and invertible contact model for trajectory optimization},
  author       = {Todorov, Emanuel},
  booktitle    = {2011 IEEE International Conference on Robotics and Automation},
  pages        = {1071--1076},
  year         = {2011},
  organization = {IEEE}
}

@article{na2021global,
  title   = {Global convergence of online optimization for nonlinear model predictive control},
  author  = {Na, Sen},
  journal = {Advances in Neural Information Processing Systems},
  volume  = {34},
  pages   = {12441--12453},
  year    = {2021}
}

@article{na2022superconvergence,
  title     = {Superconvergence of online optimization for model predictive control},
  author    = {Na, Sen and Anitescu, Mihai},
  journal   = {IEEE Transactions on Automatic Control},
  volume    = {68},
  number    = {3},
  pages     = {1383--1398},
  year      = {2022},
  publisher = {IEEE}
}

@article{pang2023global,
  title     = {Global planning for contact-rich manipulation via local smoothing of quasi-dynamic contact models},
  author    = {Pang, Tao and Suh, HJ Terry and Yang, Lujie and Tedrake, Russ},
  journal   = {IEEE Transactions on Robotics},
  year      = {2023},
  publisher = {IEEE}
}

@inproceedings{tassa2011stochastic,
  author    = {Y. Tassa AND E. Todorov},
  title     = {Stochastic Complementarity for Local Control of Discontinuous Dynamics},
  booktitle = {Proceedings of Robotics: Science and Systems},
  year      = {2010},
  address   = {Zaragoza, Spain},
  month     = {June},
  doi       = {10.15607/RSS.2010.VI.022}
}

@article{wang2009fast,
  title     = {Fast model predictive control using online optimization},
  author    = {Wang, Yang and Boyd, Stephen},
  journal   = {IEEE Transactions on control systems technology},
  volume    = {18},
  number    = {2},
  pages     = {267--278},
  year      = {2009},
  publisher = {IEEE}
}

@inproceedings{deits2014footstep,
  title        = {Footstep planning on uneven terrain with mixed-integer convex optimization},
  author       = {Deits, Robin and Tedrake, Russ},
  booktitle    = {2014 IEEE-RAS international conference on humanoid robots},
  pages        = {279--286},
  year         = {2014},
  organization = {IEEE}
}

@inproceedings{aceituno2020global,
  author    = {Bernardo Aceituno-Cabezas AND Alberto Rodriguez},
  title     = {{A Global Quasi-Dynamic Model for Contact-Trajectory Optimization in Manipulation}},
  booktitle = {Proceedings of Robotics: Science and Systems},
  year      = {2020},
  address   = {Corvalis, Oregon, USA},
  month     = {July},
  doi       = {10.15607/RSS.2020.XVI.047}
}

@article{aydinoglu2023consensus,
  title   = {Consensus complementarity control for multi-contact mpc},
  author  = {Aydinoglu, Alp and Wei, Adam and Huang, Wei-Cheng and Posa, Michael},
  journal = {arXiv preprint arXiv:2304.11259},
  year    = {2023}
}

@inproceedings{stella2017simple,
  title        = {A simple and efficient algorithm for nonlinear model predictive control},
  author       = {Stella, Lorenzo and Themelis, Andreas and Sopasakis, Pantelis and Patrinos, Panagiotis},
  booktitle    = {2017 IEEE 56th Annual Conference on Decision and Control (CDC)},
  pages        = {1939--1944},
  year         = {2017},
  organization = {IEEE}
}

@article{8584061,
  author   = {Kleinbort, Michal and Solovey, Kiril and Littlefield, Zakary and Bekris, Kostas E. and Halperin, Dan},
  journal  = {IEEE Robotics and Automation Letters},
  title    = {Probabilistic Completeness of RRT for Geometric and Kinodynamic Planning With Forward Propagation},
  year     = {2019},
  volume   = {4},
  number   = {2},
  pages    = {i-vii},
  doi      = {10.1109/LRA.2018.2888947}
}

@article{doi:10.1177/0278364915614386,
  author  = {Yanbo Li and Zakary Littlefield and Kostas E. Bekris},
  title   = {Asymptotically optimal sampling-based kinodynamic planning},
  journal = {The International Journal of Robotics Research},
  volume  = {35},
  number  = {5},
  pages   = {528-564},
  year    = {2016},
  doi     = {10.1177/0278364915614386},
  url     = {https://doi.org/10.1177/0278364915614386},
  eprint  = {https://doi.org/10.1177/0278364915614386}
}

@article{kong2023hybrid,
  title     = {Hybrid iLQR model predictive control for contact implicit stabilization on legged robots},
  author    = {Kong, Nathan J and Li, Chuanzheng and Council, George and Johnson, Aaron M},
  journal   = {IEEE Transactions on Robotics},
  volume    = {39},
  number    = {6},
  pages     = {4712--4727},
  year      = {2023},
  publisher = {IEEE}
}

@article{kim2025contact,
  title     = {Contact-implicit model predictive control: Controlling diverse quadruped motions without pre-planned contact modes or trajectories},
  author    = {Kim, Gijeong and Kang, Dongyun and Kim, Joon-Ha and Hong, Seungwoo and Park, Hae-Won},
  journal   = {The International Journal of Robotics Research},
  volume    = {44},
  number    = {3},
  pages     = {486--510},
  year      = {2025},
  publisher = {Sage Publications Sage UK: London, England}
}

@article{doi:10.1137/S1052623403429081,
  author  = {Raghunathan, Arvind U. and Biegler, Lorenz T.},
  title   = {An Interior Point Method for Mathematical Programs with Complementarity Constraints (MPCCs)},
  journal = {SIAM Journal on Optimization},
  volume  = {15},
  number  = {3},
  pages   = {720-750},
  year    = {2005},
  doi     = {10.1137/S1052623403429081},
  url     = {https://doi.org/10.1137/S1052623403429081},
  eprint  = {https://doi.org/10.1137/S1052623403429081}
}

@article{11266917,
  author   = {Ye, Xiaohan and Gao, Xifeng and Wu, Kui and Pan, Zherong and Komura, Taku},
  journal  = {IEEE Transactions on Robotics},
  title    = {SDRS: Shape-Differentiable Robot Simulator},
  year     = {2026},
  volume   = {42},
  number   = {},
  pages    = {1309-1329},
  doi      = {10.1109/TRO.2025.3636344}
}

@article{li2020incremental,
  title   = {Incremental Potential Contact: Intersection- and Inversion-free, Large-Deformation Dynamics},
  author  = {Li, Minchen and Ferguson, Zachary and Schneider, Teseo and Langlois, Timothy and Zorin, Denis and Panozzo, Daniele and Jiang, Chenfanfu and Kaufman, Danny M.},
  journal = {ACM Transactions on Graphics},
  volume  = {39},
  number  = {4},
  year    = {2020},
  doi     = {10.1145/3386569.3392425}
}

@article{lan2022affine,
  title   = {Affine Body Dynamics: Fast, Stable and Intersection-free Simulation of Stiff Materials},
  author  = {Lan, Lei and Kaufman, Danny M. and Li, Minchen and Jiang, Chenfanfu and Yang, Yin},
  journal = {ACM Transactions on Graphics},
  volume  = {41},
  number  = {4},
  pages   = {1--14},
  year    = {2022},
  doi     = {10.1145/3528223.3530064}
}

@article{boggs1995sequential,
  title   = {Sequential Quadratic Programming},
  author  = {Boggs, Paul T. and Tolle, Jon W.},
  journal = {Acta Numerica},
  volume  = {4},
  pages   = {1--51},
  year    = {1995},
  doi     = {10.1017/S0962492900002518}
}

@article{byrd2010infeasibility,
  title   = {Infeasibility Detection and {SQP} Methods for Nonlinear Optimization},
  author  = {Byrd, Richard H. and Curtis, Frank E. and Nocedal, Jorge},
  journal = {SIAM Journal on Optimization},
  volume  = {20},
  number  = {5},
  pages   = {2281--2299},
  year    = {2010},
  doi     = {10.1137/080738222}
}

@book{nocedal2006numerical,
  title     = {Numerical Optimization},
  author    = {Nocedal, Jorge and Wright, Stephen J.},
  edition   = {2},
  year      = {2006},
  publisher = {Springer},
  address   = {New York},
  doi       = {10.1007/978-0-387-40065-5}
}

@article{zhang2025simultaneous,
  title     = {Simultaneous trajectory optimization and contact selection for contact-rich manipulation with high-fidelity geometry},
  author    = {Zhang, Mengchao and Jha, Devesh K and Raghunathan, Arvind U and Hauser, Kris},
  journal   = {IEEE Transactions on Robotics},
  year      = {2025},
  publisher = {IEEE}
}

@inproceedings{moura2022non,
  title        = {Non-prehensile planar manipulation via trajectory optimization with complementarity constraints},
  author       = {Moura, Joao and Stouraitis, Theodoros and Vijayakumar, Sethu},
  booktitle    = {2022 International Conference on Robotics and Automation (ICRA)},
  pages        = {970--976},
  year         = {2022},
  organization = {IEEE}
}

@article{hinder2024worst,
  title     = {Worst-case iteration bounds for log barrier methods on problems with nonconvex constraints},
  author    = {Hinder, Oliver and Ye, Yinyu},
  journal   = {Mathematics of Operations Research},
  volume    = {49},
  number    = {4},
  pages     = {2402--2424},
  year      = {2024},
  publisher = {INFORMS}
}

@inproceedings{mellinger2011minimum,
  title        = {Minimum snap trajectory generation and control for quadrotors},
  author       = {Mellinger, Daniel and Kumar, Vijay},
  booktitle    = {2011 IEEE International Conference on Robotics and Automation},
  pages        = {2520--2525},
  year         = {2011},
  organization = {IEEE}
}

@inproceedings{ratliff2009chomp,
  title        = {{CHOMP}: Gradient optimization techniques for efficient motion planning},
  author       = {Ratliff, Nathan and Zucker, Matt and Bagnell, J. Andrew and Srinivasa, Siddhartha},
  booktitle    = {2009 IEEE International Conference on Robotics and Automation},
  pages        = {489--494},
  year         = {2009},
  organization = {IEEE}
}

@inproceedings{kalakrishnan2011stomp,
  title        = {{STOMP}: Stochastic trajectory optimization for motion planning},
  author       = {Kalakrishnan, Mrinal and Chitta, Sachin and Theodorou, Evangelos and Pastor, Peter and Schaal, Stefan},
  booktitle    = {2011 IEEE International Conference on Robotics and Automation},
  pages        = {4569--4574},
  year         = {2011},
  organization = {IEEE}
}

@incollection{moreau1988unilateral,
  title     = {Unilateral contact and dry friction in finite freedom dynamics},
  author    = {Moreau, Jean Jacques},
  booktitle = {Nonsmooth Mechanics and Applications},
  editor    = {Moreau, J. J. and Panagiotopoulos, P. D.},
  pages     = {1--82},
  year      = {1988},
  publisher = {Springer},
  address   = {Vienna},
  series    = {CISM Courses and Lectures}
}

@article{goyal1991planar,
  title     = {Planar sliding with dry friction Part 1. Limit surface and moment function},
  author    = {Goyal, Suresh and Ruina, Andy and Papadopoulos, Jim},
  journal   = {Wear},
  volume    = {143},
  number    = {2},
  pages     = {307--330},
  year      = {1991},
  publisher = {Elsevier}
}

@article{stewart2000rigid,
  title     = {Rigid-body dynamics with friction and impact},
  author    = {Stewart, David E},
  journal   = {SIAM Review},
  volume    = {42},
  number    = {1},
  pages     = {3--39},
  year      = {2000},
  publisher = {SIAM}
}

@book{bauschke2017convex,
  author    = {Bauschke, Heinz H. and Combettes, Patrick L.},
  title     = {Convex Analysis and Monotone Operator Theory in Hilbert Spaces},
  edition   = {2},
  series    = {CMS Books in Mathematics},
  publisher = {Springer},
  address   = {Cham},
  year      = {2017},
  doi       = {10.1007/978-3-319-48311-5}
}

@article{wachter2000failure,
  title     = {Failure of global convergence for a class of interior point methods for nonlinear programming},
  author    = {W{\"a}chter, Andreas and Biegler, Lorenz T.},
  journal   = {Mathematical Programming},
  volume    = {88},
  number    = {3},
  pages     = {565--574},
  year      = {2000},
  publisher = {Springer}
}

@article{fletcher2006local,
  title     = {Local convergence of {SQP} methods for mathematical programs with equilibrium constraints},
  author    = {Fletcher, Roger and Leyffer, Sven and Ralph, Daniel and Scholtes, Stefan},
  journal   = {SIAM Journal on Optimization},
  volume    = {17},
  number    = {1},
  pages     = {259--286},
  year      = {2006},
  publisher = {SIAM}
}

@article{cartis2014complexity,
  title     = {On the complexity of finding first-order critical points in constrained nonlinear optimization},
  author    = {Cartis, Coralia and Gould, Nicholas I. M. and Toint, Philippe L.},
  journal   = {Mathematical Programming},
  volume    = {144},
  number    = {1--2},
  pages     = {93--106},
  year      = {2014},
  publisher = {Springer}
}

@article{carmon2020lower,
  title     = {Lower bounds for finding stationary points {I}},
  author    = {Carmon, Yair and Duchi, John C. and Hinder, Oliver and Sidford, Aaron},
  journal   = {Mathematical Programming},
  volume    = {184},
  number    = {1--2},
  pages     = {71--120},
  year      = {2020},
  publisher = {Springer}
}

@article{chung1989np,
  title     = {{NP}-completeness of the linear complementarity problem},
  author    = {Chung, Soo Y.},
  journal   = {Journal of Optimization Theory and Applications},
  volume    = {60},
  number    = {3},
  pages     = {393--399},
  year      = {1989},
  publisher = {Springer}
}

@article{marcucci2024shortest,
  title     = {Shortest paths in graphs of convex sets},
  author    = {Marcucci, Tobia and Umenberger, Jack and Parrilo, Pablo and Tedrake, Russ},
  journal   = {SIAM Journal on Optimization},
  volume    = {34},
  number    = {1},
  pages     = {507--532},
  year      = {2024},
  publisher = {SIAM}
}

@inproceedings{reif1979complexity,
  title        = {Complexity of the mover's problem and generalizations},
  author       = {Reif, John H.},
  booktitle    = {20th Annual Symposium on Foundations of Computer Science (FOCS)},
  pages        = {421--427},
  year         = {1979},
  organization = {IEEE}
}

@book{rockafellar1998variational,
  title     = {Variational Analysis},
  author    = {Rockafellar, R. Tyrrell and Wets, Roger J.-B.},
  series    = {Grundlehren der mathematischen Wissenschaften},
  volume    = {317},
  publisher = {Springer},
  address   = {Berlin},
  year      = {1998},
  doi       = {10.1007/978-3-642-02431-3}
}

@article{lemarechal1997practical,
  title     = {Practical aspects of the {Moreau--Yosida} regularization: Theoretical preliminaries},
  author    = {Lemar{\'e}chal, Claude and Sagastiz{\'a}bal, Claudia},
  journal   = {SIAM Journal on Optimization},
  volume    = {7},
  number    = {2},
  pages     = {367--385},
  year      = {1997},
  publisher = {SIAM}
}

@article{qi1993nonsmooth,
  title     = {A nonsmooth version of {Newton's} method},
  author    = {Qi, Liqun and Sun, Jie},
  journal   = {Mathematical Programming},
  volume    = {58},
  number    = {1--3},
  pages     = {353--367},
  year      = {1993},
  publisher = {Springer}
}

@book{facchinei2003finite,
  title     = {Finite-Dimensional Variational Inequalities and Complementarity Problems},
  author    = {Facchinei, Francisco and Pang, Jong-Shi},
  series    = {Springer Series in Operations Research},
  publisher = {Springer},
  address   = {New York},
  year      = {2003},
  doi       = {10.1007/b97543}
}

@article{curtis2012sequential,
  title     = {A sequential quadratic programming algorithm for nonconvex, nonsmooth constrained optimization},
  author    = {Curtis, Frank E. and Overton, Michael L.},
  journal   = {SIAM Journal on Optimization},
  volume    = {22},
  number    = {2},
  pages     = {474--500},
  year      = {2012},
  publisher = {SIAM}
}

@book{clarke1983optimization,
  title     = {Optimization and Nonsmooth Analysis},
  author    = {Clarke, Frank H.},
  publisher = {Wiley},
  address   = {New York},
  year      = {1983}
}

@article{scheel2000mathematical,
  title     = {Mathematical programs with complementarity constraints: Stationarity, optimality, and sensitivity},
  author    = {Scheel, Holger and Scholtes, Stefan},
  journal   = {Mathematics of Operations Research},
  volume    = {25},
  number    = {1},
  pages     = {1--22},
  year      = {2000},
  publisher = {INFORMS}
}

@article{ye2005necessary,
  title     = {Necessary and sufficient optimality conditions for mathematical programs with equilibrium constraints},
  author    = {Ye, Jane J.},
  journal   = {Journal of Mathematical Analysis and Applications},
  volume    = {307},
  number    = {1},
  pages     = {350--369},
  year      = {2005},
  publisher = {Elsevier}
}

@article{outrata1999optimality,
  title     = {Optimality conditions for a class of mathematical programs with equilibrium constraints},
  author    = {Outrata, Ji{\v{r}}{\'\i} V.},
  journal   = {Mathematics of Operations Research},
  volume    = {24},
  number    = {3},
  pages     = {627--644},
  year      = {1999},
  publisher = {INFORMS}
}

@article{hoheisel2013theoretical,
  title     = {Theoretical and numerical comparison of relaxation methods for mathematical programs with complementarity constraints},
  author    = {Hoheisel, Tim and Kanzow, Christian and Schwartz, Alexandra},
  journal   = {Mathematical Programming},
  volume    = {137},
  number    = {1--2},
  pages     = {257--288},
  year      = {2013},
  publisher = {Springer}
}

@article{flegel2003fritz,
  title     = {A {Fritz John} approach to first order optimality conditions for mathematical programs with equilibrium constraints},
  author    = {Flegel, Michael L. and Kanzow, Christian},
  journal   = {Optimization},
  volume    = {52},
  number    = {3},
  pages     = {277--286},
  year      = {2003},
  publisher = {Taylor \& Francis}
}

@article{bailey2005high,
  title     = {High-precision floating-point arithmetic in scientific computation},
  author    = {Bailey, David H},
  journal   = {Computing in science \& engineering},
  volume    = {7},
  number    = {3},
  pages     = {54--61},
  year      = {2005},
  publisher = {IEEE}
}

@article{frison2020hpipm,
  title     = {HPIPM: a high-performance quadratic programming framework for model predictive control},
  author    = {Frison, Gianluca and Diehl, Moritz},
  journal   = {IFAC-PapersOnLine},
  volume    = {53},
  number    = {2},
  pages     = {6563--6569},
  year      = {2020},
  publisher = {Elsevier}
}

@article{10.1145/3811276,
  author     = {He, Zhiyong and Guo, Dewen and Guo, Minghao and Zhao, Yili and Matusik, Wojciech and Su, Hao and Jiang, Chenfanfu and Chen, Peter Yichen and Yang, Yin},
  title      = {M-ABD: Scalable, Efficient, and Robust Multi-Affine-Body Dynamics},
  year       = {2026},
  issue_date = {July 2026},
  publisher  = {Association for Computing Machinery},
  address    = {New York, NY, USA},
  volume     = {45},
  number     = {4},
  issn       = {0730-0301},
  url        = {https://doi.org/10.1145/3811276},
  doi        = {10.1145/3811276},
  journal    = {ACM Trans. Graph.},
  month      = jul,
  articleno  = {99},
  numpages   = {19}
}

@misc{alglib408,
  author       = {{ALGLIB Project}},
  title        = {{ALGLIB} 4.08.0},
  year         = {2026},
  howpublished = {Computer software},
  url          = {https://www.alglib.net/},
  note         = {MinNLC large-scale sparse filter-based SQP solver; C++ source generated 2026-06-08}
}

@article{mordatch2012discovery,
  title={Discovery of complex behaviors through contact-invariant optimization},
  author={Mordatch, Igor and Todorov, Emanuel and Popovi{\'c}, Zoran},
  journal={ACM Transactions on Graphics (ToG)},
  volume={31},
  number={4},
  pages={1--8},
  year={2012},
  publisher={ACM New York, NY, USA}
}
\clearpage
\newsavebox{\NotationBox}%
\newsavebox{\NotationBoxL}%
\newsavebox{\NotationBoxR}%
\newsavebox{\NotationCapBox}%
\newlength{\NotationAvail}%
\newcommand{\NotationHeadColor}{gray!28}
\newcommand{\NotationRowColorA}{gray!12}
\newcommand{\NotationRowColorB}{white}
\newcommand{\NotationCaption}{Notation table. (Index convention: the subscript slot is reserved for the timestep index (e.g.\ $\theta_i$, $\delta_i$, $\theta_{\gamma(i)}$); every other index is a superscript (e.g.\ $R^j$, $x^{jhk}$, $p^{jh,j'h'}$). With a matrix argument, $\lambda$ denotes an eigenvalue (e.g.\ $\lambda_{\min}(M)$), whereas $\lambda^e$ and $\lambda^i$ denote equality and inequality constraint forces.)}
%
\newcommand{\SecIndices}{%
\rowcolor{\NotationHeadColor}\multicolumn{2}{@{}l}{\textbf{Indices \& discretization}}\\
\midrule
$\delta>0,\ H$ & fixed base time scale / planning-horizon duration\\
$i\in\mathbb R$ & node label; optimized labels satisfy $0<\delta i\le H$\\
$\mathcal I$ & finite nonempty totally ordered set of optimized labels\\
$0,\ \gamma(0)<0$ & fixed time-zero node / its fixed pre-zero predecessor\\
$\gamma(i),\ \gamma^2(i)$ & previous / second-previous node labels\\
$t_i=\delta i$ & physical time at node $i$\\
$\delta_i=\delta(i-\gamma(i)),\ \delta_0>0$ & local timestep / fixed incoming timestep\\
$\delta_{\min}$ & a priori timestep floor; every optimized local timestep satisfies $\delta_i\geq\delta_{\min}$\\
$\delta^\star,\ \delta_H$ & nonsingularity / Hamiltonian timestep thresholds\\
$j,h,k$ (superscr.) & rigid-body / convex-hull / vertex index\\
$jj',\ jh,j'h'$ (superscr.) & joint-pair / convex-hull-pair index\\
}%
\newcommand{\SecCRDS}{%
\rowcolor{\NotationHeadColor}\multicolumn{2}{@{}l}{\textbf{CRDS dynamics}}\\
\midrule
$\theta,\ \theta_i$ & complete generic CRDS configuration\\
$\dot{\theta}_i,\ \ddot{\theta}_i$ & discrete velocity, acceleration\\
$M$ & constant mass matrix; positive definite in the full-rank case\\
$u,\ u_i$ & control input / current-step control input\\
$\mathcal{D}$ & compact convex admissible control space\\
$\mathcal{D}_u$ & open control-space neighborhood used for smooth extension\\
$U(\theta)$ & conservative potential energy\\
$f(\theta)$ & force induced by the conservative potential\\
$g$ & gravity term\\
$c(\theta,u),\ \bar c(\theta,u)$ & control force / mollified control force\\
$W(\theta,u)$ & potential representing control work\\
$d$ & damping force\\
$\Damp(\theta,\dot{\theta},\lambda^i)$ & damping / dissipation function\\
$h^e(\theta),\ h^i(\theta)$ & equality / inequality constraint maps\\
$\lambda_i^e,\ \lambda_i^i$ & equality / inequality constraint forces; free primal variables in~\prettyref{pb:CITO}, but prescribed penalty-induced maps $\lambda_i^e(\theta_i),\lambda_i^i(\theta_i)$ in~\prettyref{def:eps-stationary} and~\ref{def:eps-stationary-NS}\\
$\mathcal{F}^i$ & generic feasible set; also used for the ABD inequality-feasible set\\
$\partial_{\dot{\theta}}$ & subdifferential in the velocity argument\\
$O(\{\theta_i\},\{u_i\})$ & user-specified trajectory objective; see~\prettyref{ass:objective}\\
$\LControlJacobianBound$ & uniform control-force magnitude bound\\
}%
\newcommand{\SecArticulated}{%
\rowcolor{\NotationHeadColor}\multicolumn{2}{@{}l}{\textbf{Articulated body (maximal coord.)}}\\
\midrule
$x^{jhk}$ & local vertex position indexed by body, convex hull, and vertex\\
$R^j,\ t^j$ & body rotation / translation\\
$X^{jhk}$ & world-frame vertex and its time derivatives\\
$M^j$ & per-body mass matrix\\
$\Omega^j,\ \rho^j(x)$ & body volume in local frame; mass density\\
$I$ & identity matrix\\
}%
\newcommand{\SecContact}{%
\rowcolor{\NotationHeadColor}\multicolumn{2}{@{}l}{\textbf{Contact \& friction model}}\\
\midrule
$\vartheta,\ \vartheta_i$ & complete extended ABD configuration\\
$\dot\vartheta$ & complete ABD finite-difference velocity, including the plane-parameter block $\dot p$\\
$p^{jh,j'h'}$ & separating plane between a convex-hull pair\\
$n^{jh,j'h'},\ \bar{n}^{jh,j'h'}$ & plane normal, unit normal\\
$b^{jh,j'h'}$ & plane offset\\
$d^{jh,j'h'k'}$ & signed vertex-to-plane distance\\
$h^i_{\mathrm{sep}},\ h^i_{\mathrm{unit}}$ & separating-plane / unit-normal maps\\
$P^{jh,j'h'}$ & tangential projection associated with the separating plane\\
$\varsigma^{jh,j'h'},\ \omega^{jh,j'h'}$ & auxiliary plane translation / in-plane angle\\
$M^{jh,j'h'}$ & massless separating-plane block\\
$[\lambda^i]^{j'h',jhk}$ & normal contact force (plane on vertex)\\
$v_\parallel^{j'h',jhk}$ & tangential relative velocity\\
$\eta$ & friction coefficient\\
$\EpsContactPadding$ & contact padding\\
$\EpsDamp$ & plane-velocity damper coefficient\\
}%
\newcommand{\SecJoints}{%
\rowcolor{\NotationHeadColor}\multicolumn{2}{@{}l}{\textbf{Integrity constraints \& joints}}\\
\midrule
$R^{jj'},\ t^{jj'}$ & local transform to the joint constraint frame\\
$\Delta^{jj'}$ & relative anchor displacement (prismatic)\\
$L^{jj'}$ & joint limit (angle / displacement bound)\\
$\alpha^{jj'}$ & hinge joint angle\\
$x,\ y,\ z$ & constraint-frame axes (joint constraints)\\
}%
\newcommand{\SecControlMoll}{%
\rowcolor{\NotationHeadColor}\multicolumn{2}{@{}l}{\textbf{Control mollification}}\\
\midrule
$\mathcal{M}(x,L)$ & clamp / mollification map\\
$p(t)$ & shared thrice-continuously-differentiable smoothstep\\
$\bar{R}^j$ & clamped body-rotation block\\
$\operatorname{atan2}_{\EpsAtan}$ & smoothed two-argument arctangent with prescribed width\\
$u^{jj'},\ w^{jj'}$ & per-joint control (torque/force), its work\\
$\bar{d}^{jj'}$ & mollified relative anchor displacement\\
$\bar{\Delta}^{jj'}$ & control-only prismatic displacement\\
$\LPosition$ & position clamp\\
}%
\newcommand{\SecPenalization}{%
\rowcolor{\NotationHeadColor}\multicolumn{2}{@{}l}{\textbf{Penalization}}\\
\midrule
$\Phi_i(\theta_i)$ & smooth objective terms\\
$\|\cdot\|_M$ & mass-weighted norm\\
$V_s(x,\EpsFJ)$ & scalar barrier penalty\\
$V^{e,j},\ V^{i,j}$ & per-component eq.\ / ineq.\ penalties\\
$V(\theta_i,\EpsFJ)$ & total penalty potential\\
$\Damp^V$ & penalty-induced dissipation\\
$\Damp^{V,m},\ \mathcal{T}$ & individual penalty-damping term / finite damping-term index set\\
$\Lambda^i$ & penalty-induced inequality force\\
}%
\newcommand{\SecSmoothing}{%
\rowcolor{\NotationHeadColor}\multicolumn{2}{@{}l}{\textbf{Smoothing, certificate \& reduction}}\\
\midrule
$\theta_i^\star,\ \theta_i^{\dagger}$ & FJ stationary / auxiliary point\\
$\lambda_{\gamma(i)}^i$ & semi-implicit dissipation multiplier\\
$\dist(v,S),\ \partial\varphi$ & point-to-set distance; convex subdifferential\\
$\Damp_{\EpsFri,\LHamilton}^V$ & smoothed dissipation\\
$\bar\Lambda$ & constraint-force ceiling\\
$\varpi$ & stacked separating-plane velocity, distinct from plane configuration variables\\
$K_i,\ K_i^m$ & aggregate / per-term diagonal smoothed-damping block\\
$\Gamma_i,\ \Gamma_i^m$ & aggregate / per-term earlier-configuration mixed damping block\\
$C_E^\star$ & fixed energy constant\\
$\xi$ & normalized tested direction in the normal-compatibility condition ($\|\xi\|\le1,\ \xi^T K_i\xi\le C_E^\star$)\\
$\omega_N$ & normal-compatibility modulus vanishing with friction smoothing\\
$N^C_{\gph\partial_{\dot{\theta}}\Damp^{V,m}}$ & Clarke normal cone to the per-term damping graph\\
$n_i^{m\sharp}$ & per-term Clarke graph normal\\
$\EpsFJAll$ & FJ-stationarity feasibility band (chosen per result)\\
$\dot\vartheta_i$ & current-step complete ABD finite-difference velocity $<\dot\theta_i,\dot p_i,\varpi_i>$\\
$h^i(\vartheta)$ & stacked ABD inequality map with trivially extended physical components\\
$V_R(\theta_i,\EpsFJ)$ & reduced penalty value function\\
$\Damp_R^V,\ \Damp_{R,\EpsFri,\LHamilton}^V$ & reduced physical dissipation / globally extended smoothed aggregate\\
$\EpsFJ,\ \EpsFri,\ \EpsCRDS$ & penalty/FJ, friction-smoothing, residual tol.\\
}%
\newcommand{\SecHamiltonian}{%
\rowcolor{\NotationHeadColor}\multicolumn{2}{@{}l}{\textbf{Hamiltonian \& trajectory bounds}}\\
\midrule
$\mathcal{H}_i,\ \mathcal{H}_i^+$ & discrete Hamiltonian, nonnegative part\\
$\tilde{U}$ & total conservative and penalty potential\\
$\widetilde U_R$ & reduced total conservative and penalty potential $U+V_R$\\
$\bar{U}$ & position- and energy-mollified total potential\\
$\PerStepObj_i$ & per-step objective of~\prettyref{pb:PS-CRDS-M}\\
$\bar{\mathcal{H}}_i,\ \bar{\mathcal{H}}_i^+$ & mollified Hamiltonian, nonnegative part\\
$\lambda_{\min}(M)$ & smallest mass-matrix eigenvalue\\
$\theta_{\gamma(0)},\ \theta_0,\ \delta_0,\ u_0$ & fixed startup data; $\theta_0$ satisfies~\prettyref{ass:Slater}, $\delta_0>0$, and $u_0\in\mathcal D$\\
$i_1$ & first optimized node, $i_1\triangleq\min\mathcal I$\\
$\theta_{\min\mathcal I}$ & first optimized configuration\\
$\LHamBound(t_i),\ \LVelBound(t_i),\ \LPosBound(t_i)$ & Hamiltonian, velocity, position upper bounds\\
$\LHamilton$ & Hamiltonian clamp\\
$\Delta_{\mathcal H},\ \Delta_\theta$ & strict safety gaps used to select the clamps\\
$\bar W(\theta,u,\LPosition)$ & position-mollified control work\\
$\LCurvature$ & uniform two-sided curvature bound for the mollified potential\\
$\LControlWorkCurvature$ & uniform two-sided curvature bound for mollified control work\\
$\LPotentialForce$ & potential-force bound\\
$\LControlWorkForce$ & uniform mollified control-force magnitude bound\\
$\LResidual$ & per-step stationarity-residual ceiling\\
$\LObjGradLip$ & objective-gradient Lipschitz modulus (\prettyref{ass:objective})\\
}%
\newcommand{\SecSQP}{%
\rowcolor{\NotationHeadColor}\multicolumn{2}{@{}l}{\textbf{SQP convergence}}\\
\midrule
$\StackConf=(\theta_i)_{i\in\mathcal I},\ \StackCtrl=(u_i)_{i\in\mathcal I}$ & optimized configuration / control arrays on the current $\mathcal I$\\
$(\vartheta_i^\star)_{i\in\mathcal I}$ & reconstructed stacked extended ABD trajectory\\
$\StatRes(\StackConf),\ \StatRes_i(\theta_i)$ & stationarity residual (stacked / per-step)\\
$F,\ F_i$ & stacked / per-step force-balance residual\\
$\widehat{F},\ \widehat{F}_i$ & lifted stacked / per-step force-balance residual\\
$q_i,\ q_i^m$ & aggregate / per-term lifted damping force (not a generalized coordinate)\\
$Q_i$ & damping-graph argument map used in aggregate and per-term forms\\
$\StatJac(\StackConf),\ \StatJac_{ii}$ & generally nonsymmetric assembled stationarity Jacobian / its real symmetric timestep-diagonal Hessian block\\
$\StatJacU(\StackConf)$ & stationarity-residual Jacobian with respect to controls\\
$N_{\mathcal{D}}(u)$ & control-set normal cone at the specified control\\
$\odot$ & componentwise (Hadamard) product\\
$\LLambdaMin$ & threshold on $\sigma_{\min}(\StatJac_{ii})=\min_j|\lambda_j(\StatJac_{ii})|$\\
$\SigmaMinLICQ(\StatJac)$ & assembled LICQ lower-bound modulus\\
}%
\newcommand{\SecAdaptive}{%
\rowcolor{\NotationHeadColor}\multicolumn{2}{@{}l}{\textbf{Adaptive-timestep SQP}}\\
\midrule
$s$ & QP search step with configuration and control blocks\\
$B$ & block-structured QP Hessian approximation\\
$\mu_0$ & Fritz--John objective multiplier; fixed to $1$ in the unit-objective smooth certificate of~\prettyref{def:eps-stationary}, but only constrained by $\mu_0\geq0$ in the homogeneous nonsmooth certificate of~\prettyref{def:eps-stationary-NS}, where $\mu_0=0$ denotes an abnormal certificate with no first-order objective contribution\\
$\mu_\StackConf,\ \mu_\StackCtrl$ & physics-residual / control-set constraint multipliers\\
$\mu_i^e,\mu_i^h,\mu_i^\lambda,\mu_i^c$ & eq.-band / feasibility / dual-feas.\ / complementarity FJ mult.\\
$\phi_\varrho$ & one-norm merit function\\
$[\bar D\phi_\varrho]$ & guaranteed-descent surrogate for the merit function\\
$\varrho,\ \varrho^0,\ \bar\varrho$ & merit penalty weight / initial value / update margin\\
$\LHessianCond$ & uniform QP-Hessian conditioning constant\\
$\LResidualWork$ & working one-norm residual threshold\\
$\EpsStop$ & stopping tolerance for residual and step tests\\
$u_i'$ & nearby feasible control used by the two-point stationarity certificate\\
$G_0,\ N_{\mathrm{sub}}$ & initial optimized interval right endpoints / number of inserted optimized nodes\\
$N^\star$ & fixed number of timestep-damping-term pairs\\
$C_N,\ C_F$ & smoothing- and stopping-independent stationarity-band constants\\
$\EpsStopNS$ & nonsmooth stationarity tolerance; see~\prettyref{def:eps-stationary-NS}\\
}%
%
\newcommand{\NotationFontSize}{\fontsize{6.4}{7.2}\selectfont}%
\newcommand{\NotationStretch}{0.84}%
\newcommand{\NotationPanelSetup}{%
  \NotationFontSize
  \setlength{\tabcolsep}{0pt}%
  \renewcommand{\arraystretch}{\NotationStretch}%
  \rowcolors{2}{\NotationRowColorA}{\NotationRowColorB}%
}%
%
\newsavebox{\ScratchBox}%
\newcommand{\MeasureSec}[2]{%
  \begin{lrbox}{\ScratchBox}%
  \begin{minipage}[t]{0.495\textwidth}\NotationPanelSetup
  \begin{tabularx}{\linewidth}{@{}l@{\hspace{0.6em}}X@{}}\toprule #2\bottomrule\end{tabularx}%
  \end{minipage}%
  \end{lrbox}%
  \typeout{SECH #1=\the\dimexpr\ht\ScratchBox+\dp\ScratchBox\relax}%
}%
%
\begin{table*}[p]
\centering
\caption{\NotationCaption}
\label{table:notation}
\scriptsize
%
\sbox{\NotationBoxL}{%
\begin{minipage}[t]{0.495\textwidth}\NotationPanelSetup
\renewcommand{\tabularxcolumn}[1]{m{#1}}%
\begin{tabularx}{\linewidth}{@{}>{\raggedright\arraybackslash}m{0.43\linewidth}@{\hspace{0.6em}}>{\raggedright\arraybackslash}X@{}}
\toprule
\SecIndices
\midrule
\SecCRDS
\midrule
\SecJoints
\midrule
\SecControlMoll
\midrule
\SecPenalization
\midrule
\SecHamiltonian
\bottomrule
\end{tabularx}
\end{minipage}%
}%
\sbox{\NotationBoxR}{%
\begin{minipage}[t]{0.495\textwidth}\NotationPanelSetup
\renewcommand{\tabularxcolumn}[1]{m{#1}}%
\begin{tabularx}{\linewidth}{@{}>{\raggedright\arraybackslash}m{0.43\linewidth}@{\hspace{0.6em}}>{\raggedright\arraybackslash}X@{}}
\toprule
\SecContact
\midrule
\SecArticulated
\midrule
\SecSmoothing
\midrule
\SecSQP
\midrule
\SecAdaptive
\bottomrule
\end{tabularx}
\end{minipage}%
}%
\sbox{\NotationBox}{\usebox{\NotationBoxL}\hspace{0.01\textwidth}\usebox{\NotationBoxR}}%
%
\sbox{\NotationCapBox}{\normalsize\parbox{\textwidth}{\NotationCaption}}%
\setlength{\NotationAvail}{\dimexpr\textheight-\ht\NotationCapBox-\dp\NotationCapBox-2\baselineskip\relax}%
%
\adjustbox{width=\textwidth,max totalheight=\NotationAvail,center}{\usebox{\NotationBox}}%
\end{table*}

\clearpage
\noindent The following appendices collect the proofs of the results stated
in the main text. A consolidated list of the symbols used in the analysis
is provided in~\prettyref{table:notation}; the benchmark-specific objective
terms of the numerical experiments are defined in place. Individual proofs may
additionally introduce their own local symbols (proof-specific auxiliary
quantities); these are defined where they appear and are intentionally
omitted from~\prettyref{table:notation}. Throughout the appendices, the damping/dissipation function is written $\Damp$ (with variants $\Damp^V$, $\Damp_{\EpsFri,\LHamilton}^V$, $\Damp_R^V$, $\Damp_{R,\EpsFri,\LHamilton}^V$); a prefix $D$ applied to a map (e.g.\ $Dq$, $Dp^\star$, $DG$, $D\varpi^\star$, $D\mathcal{M}$) instead denotes that map's total derivative/Jacobian.

\appendix
\section{Proofs for Coordinate Maximization}
\label{appen:coordinate-proofs}

\begin{proof}[Proof of~\prettyref{lem:mass}]
It suffices to prove the claim for a single rigid body, so we drop the superscript $j$. Let $q=(\operatorname{vec}(\dot R),\dot t)\in\mathbb{R}^{12}$ and define
\begin{equation}
A(x)=\TWO{x^T\otimes I_3}{I_3},
\end{equation}
so that $A(x)\in\mathbb{R}^{3\times12}$ and $\dot R x+\dot t=A(x)q$, using the standard vectorization identity $\dot R x=(x^T\otimes I_3)\operatorname{vec}(\dot R)$. The kinetic energy of the body is
\begin{equation}
K(q)=\frac{1}{2}\int_\Omega \rho(x)\|A(x)q\|^2dx=\frac{1}{2}q^TMq,
\end{equation}
where
\begin{equation}
M=\int_\Omega \rho(x)A(x)^TA(x)dx
\end{equation}
is the Hessian of $K$. Since $\Omega$ is bounded, there is a finite $r_\Omega$ such that $\|x\|\leq r_\Omega$ for every $x\in\Omega$, and its Lebesgue measure is finite. Bounded measurability of $\rho$ then gives finite mass. The entries of $A(x)^TA(x)$ are polynomials of degree at most two in $x$, so boundedness of $\Omega$ also gives a finite second moment. Hence the integral defining $M$ is finite entrywise, and $M\in\mathbb{R}^{12\times12}$. Therefore
\begin{equation}
q^TMq=\int_\Omega \rho(x)\|\dot R x+\dot t\|^2dx\geq0.
\end{equation}
If $q^TMq=0$, then $\dot R x+\dot t=0$ for almost all $x\in\Omega$, because $\rho(x)>0$ almost everywhere. A nonzero affine map from $\mathbb{R}^3$ to $\mathbb{R}^3$ has a zero set of Lebesgue measure zero. Since $|\Omega|>0$, the affine map $x\mapsto \dot R x+\dot t$ must be identically zero, which implies $\dot R=0$ and $\dot t=0$. Hence $q=0$, so the single-body mass block is positive definite. Applying this argument to every rigid body gives $M^j\succ0$ for every $j$. Since $M=\diag(\{M^j\})$ is block diagonal with positive-definite blocks, the rigid-body mass matrix satisfies $M\succ0$.
\end{proof}

\begin{proof}[Proof of~\prettyref{lem:separate}]
Let $y^k=X^{jhk}=R^jx^{jhk}+t^j$, $z^{k'}=X^{j'h'k'}=R^{j'}x^{j'h'k'}+t^{j'}$, and set $A=\operatorname{conv}\{y^k\}_k$, $B=\operatorname{conv}\{z^{k'}\}_{k'}$.
First suppose the oriented separating plane $p^{jh,j'h'}=\TWO{n^{jh,j'h'}}{b^{jh,j'h'}}$ satisfies~\prettyref{eq:separating-plane} and~\ref{eq:unit-normal}, with $A$ on the positive side and $B$ on the negative side. Since the plane inequalities are affine, for every $a\in A$ and $q\in B$, we have
\begin{equation}
\begin{aligned}
\langle n^{jh,j'h'},a\rangle+b^{jh,j'h'}\geq&\EpsContactPadding,\\
\langle n^{jh,j'h'},q\rangle+b^{jh,j'h'}\leq&-\EpsContactPadding.
\end{aligned}
\end{equation}
Therefore,
\begin{equation}
\|a-q\|\geq\langle n^{jh,j'h'},a-q\rangle\geq2\EpsContactPadding,
\end{equation}
where the first inequality uses $\|n^{jh,j'h'}\|\leq1$, which follows from~\prettyref{eq:unit-normal}. Taking the infimum over $A\times B$ gives $\dist(A,B)\geq2\EpsContactPadding$.

Conversely, suppose $\dist(A,B)\geq2\EpsContactPadding$. Compactness gives closest points $a^\star\in A$ and $q^\star\in B$, with $a^\star\neq q^\star$. Define
\begin{equation}
n=\frac{a^\star-q^\star}{\|a^\star-q^\star\|},
b_0=-\frac{1}{2}\langle n,a^\star+q^\star\rangle.
\end{equation}
The closest-point (metric projection) optimality conditions---the variational inequalities for projection onto the convex sets $A$ and $B$---imply
\begin{equation}
\begin{aligned}
\langle n,a-a^\star\rangle\geq&0 &&\forall a\in A,\\
\langle n,q-q^\star\rangle\leq&0 &&\forall q\in B.
\end{aligned}
\end{equation}
Thus, for every $a\in A$ and $q\in B$,
\begin{equation}
\begin{aligned}
\langle n,a\rangle+b_0\geq&\frac{1}{2}\|a^\star-q^\star\|\geq \EpsContactPadding,\\
\langle n,q\rangle+b_0\leq&-\frac{1}{2}\|a^\star-q^\star\|\leq-\EpsContactPadding.
\end{aligned}
\end{equation}
These inequalities hold at all vertices, so~\prettyref{eq:separating-plane} holds with $p^{jh,j'h'}=\TWO{n}{b_0}$. Since $\|n\|=1$,~\prettyref{eq:unit-normal} holds as well.
\end{proof}

\begin{proof}[Proof of~\prettyref{cor:abd-crds}]
We verify conditions~\ref{it:crds-c2}, \ref{it:crds-potential}, \ref{it:crds-damping}, and \ref{it:crds-control} of~\prettyref{def:CRDS} for the articulated-body construction, with the extended configuration variable $\vartheta$ in~\prettyref{eq:extended-config} playing the role of the CRDS configuration.

Condition~\ref{it:crds-c2}: The constraint maps in~\prettyref{eq:integrity},\ref{eq:hinge},\ref{eq:prismatic},\ref{eq:separating-plane},\ref{eq:unit-normal} are finite-dimensional real-valued maps whose components are polynomial or affine functions of the maximal-coordinate variables and separating-plane variables. Hence $h^e$ and $h^i$ are $C^2$ real functions of $\vartheta_i$.

For the hinge-angle map in~\prettyref{eq:hinge-control}, we use the following smoothed version of $\operatorname{atan2}$. Let $r=\sqrt{x^2+y^2}$, choose fixed branch-cut smoothing constants $0<\xi_0<\xi_1<1$, and, only when $r>\EpsAtan$, define $\xi=(r+x)/(2r)$. Using the shared polynomial $p$ from~\prettyref{eq:clamp-map}, set
\begin{ResizedEq}
&\operatorname{atan2}_{\EpsAtan}(y,x)=\\
&\begin{cases}
0,&r\leq\EpsAtan,\\
0,&r>\EpsAtan\ \text{and}\ \xi\leq\xi_0,\\
p\!\left(\frac{r-\EpsAtan}{\EpsAtan}\right)
p\!\left(\frac{\xi-\xi_0}{\xi_1-\xi_0}\right)
\operatorname{atan2}(y,x),&r>\EpsAtan\ \text{and}\ \xi>\xi_0.
\end{cases}
\end{ResizedEq}
The second branch is used only where the selected branch of $\operatorname{atan2}$ is smooth. The cutoff $p$ is $C^3$ and has zero first, second, and third derivatives at $0$ and $1$. Therefore the radial cutoff gives a $C^3$ extension through $r=\EpsAtan$, while the branch-cut cutoff glues to the zero branch before the discontinuity of $\operatorname{atan2}$ is reached. Away from these cutoff boundaries the expression is a composition of $C^3$ functions. Thus $\operatorname{atan2}_{\EpsAtan}$ is $C^3$, and the hinge work is $C^2$ because its two scalar arguments are $C^2$ functions of the clamped rotations. The prismatic contribution is $C^2$, and finite sums preserve $C^2$ regularity, so the aggregate control work potential $W$ in~\prettyref{eq:abd-control-work} is a $C^2$ real function; being linear in the controls $u$ with these $C^2$ configuration coefficients, it is jointly $C^2$ in $(\theta_i,u_i)$, verifying~\ref{it:crds-c2}. The force map in~\prettyref{eq:control} is then obtained as $c=-\FPPR{W}{\vartheta_i}$.

Condition~\ref{it:crds-potential}: The articulated-body specialization sets the conservative potential to $U(\vartheta_i)\equiv0$. This maps the configuration domain into $[0,+\infty]$ and is proper because it is finite everywhere. It is continuous as an extended-real-valued function, and its finite domain is the full configuration domain, where it is trivially $C^2$. Its conservative force is $f=-\FPPR{U}{\vartheta_i}=0$, agreeing with the articulated-body construction.

Condition~\ref{it:crds-damping}: Fix $\vartheta_\gamma\in\mathcal{F}^i$, where the overloaded full-coordinate feasible set $\mathcal{F}^i$ is defined in~\prettyref{eq:abd-extended-feasible-set}, and fix a neutral nonnegative damping load $\lambda_\gamma^i\geq0$, as required by~\prettyref{def:CRDS}. On this domain, every separating-plane normal appearing in~\prettyref{eq:damping-ABD} is nonzero: by the orientation convention above, the two sides of~\prettyref{eq:separating-plane} use the same oriented plane, so otherwise the padded separation inequalities would require the same affine offset to be both at least $\EpsContactPadding$ and at most $-\EpsContactPadding$. Hence the normalized normal $\bar n$ and projection $P$ are well-defined fixed finite coefficients when $\Damp(\vartheta_\gamma,\cdot,\lambda_\gamma^i)$ is viewed as a function of $\dot{\vartheta}_i$. Under the semi-implicit convention in~\prettyref{eq:CRDS}, all non-velocity quantities in~\prettyref{eq:tangent-velocity} and~\ref{eq:damping-ABD}, including $\bar n$, $P$, and the positions $X^{jhk}$, are evaluated at the fixed predecessor configuration $\vartheta_\gamma$. Each tangential relative velocity in~\prettyref{eq:tangent-velocity} is therefore linear in the fixed projection $\TWO{\dot\theta_i}{\varpi_i}$ of the full velocity and independent of $\dot p_i$. Composing a convex function with this fixed linear projection preserves convexity. The terms weighted by $\eta\lambda_\gamma^i$ are nonnegative multiples of norms of linear functions, and the plane-velocity damper $\EpsDamp(\|\dot\omega^{jh,j'h'}\|^2+\|\dot{\varsigma}^{jh,j'h'}\|^2)$ is convex because $\EpsDamp>0$. Finite sums preserve convexity, so $\Damp(\vartheta_\gamma,\cdot,\lambda_\gamma^i)$ is convex in the complete $\dot\vartheta_i$ on the domain where it is defined. Its ambient subgradient is the damping-active subgradient lifted with an exact zero $\dot p_i$ block, and the damping force satisfies $d\in-\partial_{\dot\vartheta_i}\Damp$. Moreover, each tangential relative velocity $v_\parallel$ in~\prettyref{eq:tangent-velocity}, the angular velocity $\dot\omega^{jh,j'h'}$, and the translational velocity $\dot{\varsigma}^{jh,j'h'}$ is linear in $\TWO{\dot\theta_i}{\varpi_i}$, hence has $0$ in its ambient subdifferential at $\dot\vartheta_i=0$; since $\Damp(\vartheta_\gamma,\cdot,\lambda_\gamma^i)$ is convex in $\dot{\vartheta}_i$ and $0\in\partial_{\dot{\vartheta}_i}\Damp(\vartheta_\gamma,0,\lambda_\gamma^i)$, zero relative velocity is a global minimizer. This verifies the minimizer condition required in~\prettyref{def:CRDS}~\ref{it:crds-damping}. It is not smooth in general, because the norm terms are non-differentiable at zero tangential relative velocity when the corresponding contact weight is positive.

Condition~\ref{it:crds-control}: We first make explicit what the derivative with respect to the control means. For each controlled joint $a=jj'$, write its scalar controlled coordinate as $\phi^a(\theta)$, with $\phi^a=\alpha^{jj'}$ for a hinge joint and $\phi^a=q^{jj'}$ for a prismatic joint. Both control-work definitions have the form
\begin{ResizedEq}
w^a(\theta,u^a)=&-u^a\phi^a(\theta)\\
W(\theta,u)=&-\sum_a u^a\phi^a(\theta).
\end{ResizedEq}
The force map is $c(\vartheta_i,u_i)=-\FPPR{W}{\vartheta_i}$. Since $W$ depends only on the physical coordinates inside $\vartheta_i$, the entries of $c$ corresponding to the auxiliary separating-plane variables are zero. For any physical coordinate component $\vartheta_i^\beta$,
\begin{ResizedEq}
c_\beta(\vartheta_i,u_i)=&\sum_a u^a\FPP{\phi^a(\theta_i)}{\vartheta_i^\beta}\\
\FPP{c_\beta(\vartheta_i,u_i)}{u^a}=&\FPP{\phi^a(\theta_i)}{\vartheta_i^\beta}.
\end{ResizedEq}
Thus the $a$th column of $\FPPR{c(\vartheta_i,u_i)}{u_i}$ is exactly the extended configuration gradient $\nabla_{\vartheta_i}\phi^a$, with zeros in the auxiliary-variable rows. It remains only to bound the configuration gradient of each controlled coordinate globally.

For a hinge joint, $\phi^a=\alpha^{jj'}$ in~\prettyref{eq:hinge-control}. Define
\begin{ResizedEq}
Y=\langle \bar R^jR^{jj'}y,\bar R^{j'}R^{j'j}x\rangle,
X=\langle \bar R^jR^{jj'}x,\bar R^{j'}R^{j'j}x\rangle.
\end{ResizedEq}
Then $\alpha^{jj'}=\operatorname{atan2}_{\EpsAtan}(Y,X)$. By~\prettyref{eq:clamp-map}, $\|\bar R^j\|\leq2\LPosition$. The Jacobian of $R^j\mapsto\bar R^j=\mathcal M(R^j,\LPosition)$ is uniformly bounded because $\mathcal M$ is a fixed $C^2$ radial map that is the identity on the inner region, has a compact transition region, and is radial with constant norm outside the transition. Hence $Y,X$ and their first derivatives with respect to $R^j,R^{j'}$ are bounded by constants depending only on $\LPosition$ and the fixed joint-frame data. The derivatives of $\operatorname{atan2}_{\EpsAtan}$ are globally bounded by the smoothing construction above: the radial cutoff $p\!\left(\frac{r-\EpsAtan}{\EpsAtan}\right)$ forces $\operatorname{atan2}_{\EpsAtan}\equiv0$ whenever $r=\sqrt{Y^2+X^2}\leq\EpsAtan$, so the raw gradient $\nabla\operatorname{atan2}$ (of norm $1/r$) is only ever evaluated where $r>\EpsAtan$ and is thus capped by $1/\EpsAtan$, while the two smoothstep-derivative terms each carry a single $1/\EpsAtan$ factor (using $\|\nabla r\|\leq1$ and $\|\nabla \xi\|\leq1/(2r)$, with the branch-cut constants $\xi_0,\xi_1$ kept fixed); hence $\|\nabla\operatorname{atan2}_{\EpsAtan}\|=\mathcal{O}(1/\EpsAtan)$ with no residual singularity, so the chain rule
\begin{ResizedEq}
\nabla\alpha^{jj'}=&
(\partial_Y\operatorname{atan2}_{\EpsAtan})(Y,X)\nabla Y+\\
&(\partial_X\operatorname{atan2}_{\EpsAtan})(Y,X)\nabla X
\end{ResizedEq}
shows that every hinge column is uniformly bounded by a constant depending only on $\LPosition$, $\EpsAtan$, and the fixed joint-frame data.

For a prismatic joint, $\phi^a=q^{jj'}$, where
\begin{equation}
q^{jj'}=\langle\bar\Delta^{jj'},z\rangle=\langle\bar d^{jj'},a^{jj'}\rangle,
a^{jj'}\triangleq\bar R^jR^{jj'}z,
\end{equation}
and $\bar d^{jj'}$ is defined in~\prettyref{eq:clamped-anchor-displacement}. By~\prettyref{eq:clamp-map}, $\|a^{jj'}\|\leq2\LPosition\|R^{jj'}z\|$ and $\|\bar d^{jj'}\|\leq2\LPosition$. For any physical coordinate component $\xi$, the product rule gives
\begin{equation}
\FPPR{q^{jj'}}{\xi}=\left\langle\FPPR{\bar d^{jj'}}{\xi},a^{jj'}\right\rangle+
\left\langle\bar d^{jj'},\FPPR{a^{jj'}}{\xi}\right\rangle.
\end{equation}
If $\xi$ is a translation component, then $\partial a^{jj'}/\partial\xi=0$, while $\partial\bar d^{jj'}/\partial\xi$ is the uniformly bounded Jacobian of $\mathcal M(\cdot,\LPosition)$ applied to a fixed vector with entries $\pm1$ or $0$. Therefore translation derivatives are globally bounded. If $\xi$ is a rotation component, then $\partial a^{jj'}/\partial\xi$ is the uniformly bounded Jacobian of the clamped rotation map applied to the fixed vector $R^{jj'}z$, so the second product-rule term is bounded by $\|\bar d^{jj'}\|$ times a fixed constant. The first product-rule term is also bounded: $\bar d^{jj'}=\mathcal M(s^{jj'},\LPosition)$ with
\begin{equation}
s^{jj'}=(\bar R^{j'}t^{j'j}+t^{j'})-(\bar R^jt^{jj'}+t^j),
\end{equation}
so $\partial\bar d^{jj'}/\partial\xi=D\mathcal M(s^{jj'},\LPosition)\partial s^{jj'}/\partial\xi$; the first factor is uniformly bounded and the second factor contains only uniformly bounded derivatives of clamped rotations applied to the fixed anchors $t^{jj'},t^{j'j}$. Consequently every derivative of $q^{jj'}$ with respect to the maximal-coordinate variables is uniformly bounded by a constant depending only on $\LPosition$ and fixed joint-frame data.

Every controlled joint therefore contributes a uniformly bounded controlled-coordinate gradient $\nabla_{\vartheta_i}\phi^a$ globally, with a bound depending only on $\LPosition$, $\EpsAtan$, and the fixed joint-frame data.

It remains to bound the control force itself for admissible controls $u_i\in\mathcal{D}$. From~\prettyref{eq:control} we have $c(\vartheta_i,u_i)=\sum_a u^a\nabla_{\vartheta_i}\phi^a$, so the triangle inequality gives $\|c(\vartheta_i,u_i)\|\leq\sum_a|u^a|\,\|\nabla_{\vartheta_i}\phi^a\|$. Since $\mathcal{D}$ is compact, $U_{\max}\triangleq\max_{u\in\mathcal{D}}\|u\|<\infty$, and the controlled-coordinate gradients $\nabla_{\vartheta_i}\phi^a$ are globally bounded by the hinge and prismatic estimates established above; hence $\|c(\vartheta_i,u_i)\|\leq U_{\max}\sum_a\|\nabla_{\vartheta_i}\phi^a\|\triangleq C_c<\infty$, which is finite because the articulated body has finitely many controlled joints. Choosing $\LControlJacobianBound\geq C_c$ yields $\|c(\vartheta_i,u_i)\|\leq\LControlJacobianBound$ for every configuration $\vartheta_i$ and every admissible control $u_i\in\mathcal{D}$, verifying~\prettyref{def:CRDS}~\ref{it:crds-control}.

The articulated-body construction therefore defines a CRDS. Finally, the extended configuration includes massless separating-plane variables, whose mass blocks are zero by construction; therefore the extended mass matrix is singular and the extended CRDS is not full-rank in the sense of~\prettyref{def:full-rank}. The norm terms in the frictional dissipation also show that, in general, the extended CRDS is not smooth in the sense of~\prettyref{def:smooth-CRDS}.
\end{proof}

\section{Proofs for Penalization}
\label{appen:penalty-proofs}

\subsection{Approximate CRDS}
\begin{proof}[Proof of~\prettyref{lem:eps-FJ}]
Finite objective value in~\prettyref{pb:P-CRDS} implies that $V(\theta_i^\star,\EpsFJ)<\infty$. Since $V_s(x,\EpsFJ)=\infty$ when $x\leq0$, finiteness of the equality penalties gives
\begin{equation}
\EpsFJ-h^{e,j}(\theta_i^\star)>0\quad\text{and}\quad\EpsFJ+h^{e,j}(\theta_i^\star)>0
\end{equation}
for every equality component $j$. Hence $|h^{e,j}(\theta_i^\star)|<\EpsFJ$. Finiteness of the inequality penalties similarly gives $h^{i,j}(\theta_i^\star)>0$ for every inequality component $j$.

We now construct the current constraint multipliers from the current penalty gradients and the semi-implicit dissipation multiplier from the previous-step penalty forces. Define
\begin{ResizedEq}
\lambda_i^{e,j}\triangleq
&-\FPP{V_s}{x}(\EpsFJ-h^{e,j}(\theta_i^\star),\EpsFJ)+\\
&\FPP{V_s}{x}(\EpsFJ+h^{e,j}(\theta_i^\star),\EpsFJ)
\end{ResizedEq}
and
\begin{equation}
\lambda_i^{i,j}\triangleq-\FPP{V_s}{x}(h^{i,j}(\theta_i^\star),\EpsFJ).
\end{equation}
For $x>0$, $V_s(\cdot,\EpsFJ)$ is nonincreasing on its nonzero branch and has zero derivative on its zero branch, so $\lambda_i^{i,j}\geq0$. The hypothesis $V(\theta_{\gamma(i)}^\star,\EpsFJ)\leq\LHamilton<\infty$ and nonnegativity of the penalty terms imply that every previous-step inequality term $V_s(h^{i,j}(\theta_{\gamma(i)}^\star),\EpsFJ)$ is finite. Since $V_s(x,\EpsFJ)=+\infty$ for $x\leq0$, it follows that $h^{i,j}(\theta_{\gamma(i)}^\star)>0$ for every $j$, and hence
\begin{equation}
\theta_{\gamma(i)}^\star\in\operatorname{int}\mathcal F^i\subseteq\mathcal F^i.
\end{equation}
Define the semi-implicit dissipation multiplier componentwise from the previous-step penalty force~\prettyref{eq:penalty-force-generic} evaluated at $\theta_{\gamma(i)}^\star$,
\begin{equation}
\lambda_{\gamma(i)}^{i,j}\triangleq-V_s'(h^{i,j}(\theta_{\gamma(i)}^\star),\EpsFJ);
\end{equation}
since its argument is positive, this derivative is finite, and since $V_s'\leq0$ on its nonzero branch and vanishes elsewhere, $\lambda_{\gamma(i)}^i$ is finite and nonnegative. The chain rule gives
\begin{equation}
\nabla_{\theta_i}V(\theta_i^\star,\EpsFJ)
=\FPP{h^e}{\theta_i}(\theta_i^\star)^T\lambda_i^e-
\FPP{h^i}{\theta_i}(\theta_i^\star)^T\lambda_i^i.
\end{equation}

Finite objective value in~\prettyref{pb:P-CRDS} also implies $U(\theta_i^\star)<+\infty$: all the other terms in the objective are finite real-valued or nonnegative extended-real terms at a point where their sum is finite. Hence $\theta_i^\star\in\operatorname{dom}U$. By~\prettyref{def:CRDS}, $U:\mathbb R^n\to[0,+\infty]$ is continuous as an extended-real-valued map. Since $[0,+\infty)$ is open in the extended-real codomain, $\operatorname{dom}U=U^{-1}([0,+\infty))$ is open. The strict inequalities established above likewise put $\theta_i^\star$ in $\operatorname{dom}V$; because $h^e$ and $h^i$ are continuous and these inequalities are strict, $\operatorname{dom}V$ is open around $\theta_i^\star$. Therefore $\theta_i^\star$ lies in the open set $\operatorname{dom}U\cap\operatorname{dom}V$, on which $\Phi_i+V(\cdot,\EpsFJ)$ is $C^1$: $U$ is $C^2$ on its finite domain, $V$ is $C^2$ on its finite domain, and the remaining terms are smooth by~\prettyref{def:CRDS}.

Hold $\theta_{\gamma(i)}^\star$ and the positive timestep $\delta_i$ fixed. Since $\theta_{\gamma(i)}^\star\in\mathcal F^i$ as established above,~\prettyref{def:CRDS}~\ref{it:crds-damping} says that $\Damp^V(\theta_{\gamma(i)}^\star,\cdot)$ is convex in velocity. The map $\theta_i\mapsto\dot\theta_i=(\theta_i-\theta_{\gamma(i)}^\star)/\delta_i$ is affine. Thus the composition of the damping with this affine map, multiplied by $\delta_i>0$, is convex in $\theta_i$. Applying the smooth-plus-convex Fermat rule at the interior local minimizer $\theta_i^\star$ gives
\begin{ResizedEq}
0\in\nabla_{\theta_i}\Phi_i(\theta_i^\star)+\nabla_{\theta_i}V(\theta_i^\star,\EpsFJ)+\partial_{\theta_i}\left(\delta_i\Damp^V\right)(\theta_{\gamma(i)}^\star,\dot{\theta}_i^\star).
\end{ResizedEq}
Set $\theta_i^{\dagger}=\theta_i^\star$. Then $\|\dot{\theta}_i^\star-\dot{\theta}_i^{\dagger}\|=0\leq\EpsFJ$, so the closeness condition in~\prettyref{eq:eps-fj-system} holds. With this choice, $\Damp^V(\theta_{\gamma(i)}^\star,\dot\theta_i)=\Damp(\theta_{\gamma(i)}^\star,\dot\theta_i,\lambda_{\gamma(i)}^i)$ by definition of the semi-implicit penalty-induced dissipation. Substituting the multiplier identity above gives the stationarity line in~\prettyref{eq:eps-fj-system} with $\lambda_i^i$ in the current inequality-gradient term and $\lambda_{\gamma(i)}^i$ in the dissipation subgradient. Finally, we verify the combined inequality line in~\prettyref{eq:eps-fj-system}. We have $h^{i,j}(\theta_i^\star)>0$ and $\lambda_i^{i,j}\geq0$, both shown above. For the product term: if $h^{i,j}(\theta_i^\star)\geq\EpsFJ$, then $V_s'(h^{i,j}(\theta_i^\star),\EpsFJ)=0$ (the derivative vanishes on the zero branch $x>\EpsFJ$ and, by the order-four zero at $x=\EpsFJ$, also at the boundary), so $\lambda_i^{i,j}=0$ and $\lambda_i^{i,j}(\EpsFJ-h^{i,j}(\theta_i^\star))=0$; if $h^{i,j}(\theta_i^\star)<\EpsFJ$, then $\EpsFJ-h^{i,j}(\theta_i^\star)>0$ and $\lambda_i^{i,j}\geq0$, so $\lambda_i^{i,j}(\EpsFJ-h^{i,j}(\theta_i^\star))\geq0$. Hence $\min\!\big(h^{i,j}(\theta_i^\star),\lambda_i^{i,j},\lambda_i^{i,j}(\EpsFJ-h^{i,j}(\theta_i^\star))\big)\geq0$ for all $j$. Thus $\theta_i^\star$ satisfies the $\EpsFJ$-FJ condition.
\end{proof}

\begin{proof}[Proof of~\prettyref{cor:S-CRDS}]
We verify conditions~\ref{it:crds-c2}, \ref{it:crds-potential}, \ref{it:crds-damping}, and \ref{it:crds-control} of~\prettyref{def:CRDS}.

Condition~\ref{it:crds-c2}: The hard constraint maps of~\prettyref{pb:PS-CRDS} are empty. Empty maps are closed, differentiable, real-valued maps on their domains, and the corresponding multiplier terms disappear from the equations. The remaining object governed by condition~\ref{it:crds-c2} is the control-work potential $W$, which is inherited unchanged from the original CRDS and hence remains a $C^2$ real function; the regularity of the induced control map $c=-\partial W/\partial\theta_i$ is treated under condition~\ref{it:crds-control} below.

Condition~\ref{it:crds-potential}: Define the conservative potential of the unconstrained system by
\begin{equation}
\widetilde U(\theta_i,\EpsFJ)\triangleq U(\theta_i)+V(\theta_i,\EpsFJ).
\end{equation}
The inherited conservative potential $U(\theta_i)$ is already contained in $\Phi_i$, while $V(\theta_i,\EpsFJ)$ is the added barrier potential. Since the original system is a CRDS, $U$ satisfies the potential condition of~\prettyref{def:CRDS}: it is nonnegative, proper, extended-real continuous, and $C^2$ on $\operatorname{dom}U$.

The scalar penalty $V_s(\cdot,\EpsFJ)$ is nonnegative and extended-real continuous as a map into $[0,+\infty]$. Its finite domain is the open interval $(0,+\infty)$, and it is $C^2$ there: the rational branch is $C^2$ on $(0,\EpsFJ)$, the zero branch is $C^2$ on $(\EpsFJ,+\infty)$, and the two branches have matching value, first derivative, and second derivative at $x=\EpsFJ$, because the rational branch carries a factor $(x-\EpsFJ)^4$ with a zero of order four at $x=\EpsFJ$, so its value, first, and second derivatives all vanish there and match the zero branch. Because $h^e$ and $h^i$ are $C^2$ real functions, the composed penalties $V^{e,j}$ and $V^{i,j}$ are extended-real continuous and $C^2$ on their finite domains. The finite domain of $V$ is defined by strict inequalities such as $|h^{e,j}(\theta_i)|<\EpsFJ$ and $h^{i,j}(\theta_i)>0$, hence $\operatorname{dom}V$ is open. Therefore
\begin{equation}
\operatorname{dom}\widetilde U=\operatorname{dom}U\cap\operatorname{dom}V
\end{equation}
is open. The nonempty-intersection assumption in~\prettyref{cor:S-CRDS} makes $\widetilde U$ proper, and the preceding continuity and differentiability properties show that $\widetilde U$ is nonnegative, extended-real continuous, and $C^2$ on this finite domain. Its conservative force is
\begin{ResizedEq}
\widetilde f(\theta_i)
=&-\FPP{\widetilde U(\theta_i,\EpsFJ)}{\theta_i}\\
=&-\FPP{U(\theta_i)}{\theta_i}-\FPP{V(\theta_i,\EpsFJ)}{\theta_i},
\end{ResizedEq}
which is exactly the original conservative force plus the differentiable penalty force obtained by absorbing the hard equality and inequality constraints into $V$.

Condition~\ref{it:crds-damping}: The unconstrained \prettyref{pb:PS-CRDS} has empty hard inequality maps, so under the unchanged definition of $\mathcal F^i$ its hard feasible set is $\mathcal F^i=\mathbb R^n$. By the strengthened~\prettyref{def:S-D}, the smoothed dissipation $\Damp_{\EpsFri,\LHamilton}^V$ is an everywhere-defined finite extension on all of $\mathbb R^n\times\mathbb R^n$, so it supplies a smoothing on this full set: for every fixed predecessor configuration $\theta_{\gamma(i)}\in\mathbb R^n$, $\Damp_{\EpsFri,\LHamilton}^V(\theta_{\gamma(i)},\cdot)$ is everywhere evaluable, convex and nonnegative in $\dot{\theta}_i$, and globally minimized at zero relative velocity, i.e. $0\in\partial_{\dot{\theta}_i}\Damp_{\EpsFri,\LHamilton}^V(\theta_{\gamma(i)},0)$. These convexity, nonnegativity, and zero-velocity global-minimization properties therefore satisfy the CRDS damping condition~\prettyref{def:CRDS}~\ref{it:crds-damping} for every required predecessor configuration, on the full domain required by~\prettyref{def:CRDS}, not merely on $\operatorname{dom}\widetilde U$. The penalty force $\Lambda^i$ and the raw nonsmooth dissipation $\Damp^V$ are never evaluated outside $\operatorname{dom}V$: they enter only through the approximation clauses of~\prettyref{def:S-D}, which are restricted to the finite penalty sublevel $V(\theta_{\gamma(i)},\EpsFJ)\leq\LHamilton<\infty\subseteq\operatorname{dom}V$.

Condition~\ref{it:crds-control}: The control map is unchanged, so the original control bounds are inherited: the same constant $\LControlJacobianBound>0$ satisfies $\|c(\theta_i,u_i)\|\leq\LControlJacobianBound$ for every configuration $\theta_i$ and every admissible control $u_i\in\mathcal{D}$.

Finally, by the strengthened~\prettyref{def:S-D}, $\Damp_{\EpsFri,\LHamilton}^V$ is differentiable in $\dot{\theta}_i$ and its velocity gradient $\FPPR{\Damp_{\EpsFri,\LHamilton}^V}{\dot{\theta}_i}$ is globally jointly $C^1$ with respect to all arguments appearing in the discrete equations. Hence the CRDS is smooth in the sense of~\prettyref{def:smooth-CRDS}.
\end{proof}

\begin{proof}[Proof of~\prettyref{lem:eps-S-FJ}]
Set
\begin{equation}
\EpsFJAll\triangleq\max\{\EpsFJ,
\EpsFri,
\EpsCRDS+
\delta_i\EpsFri\}.
\end{equation}
We verify the lines of~\prettyref{eq:eps-fj-system} with this tolerance, constructing the current constraint multipliers and the semi-implicit dissipation multiplier separately.

Feasibility: Finite objective value in~\prettyref{pb:PS-CRDS} implies that $V(\theta_i^\star,\EpsFJ)<\infty$. By the definition of the equality and inequality penalties, exactly as in the proof of~\prettyref{lem:eps-FJ}, every equality component satisfies
\begin{equation}
|h^{e,j}(\theta_i^\star)|<\EpsFJ\leq\EpsFJAll,
\end{equation}
and every inequality component satisfies
\begin{equation}
h^{i,j}(\theta_i^\star)>0.
\end{equation}
This gives the equality-feasibility line $|h^{e,j}(\theta_i^\star)|\leq\EpsFJAll$ of~\prettyref{eq:eps-fj-system} and the primal-feasibility term $h^{i,j}(\theta_i^\star)>0$ of the combined inequality line, which is concluded below at tolerance $\EpsFJAll$.

Penalty multipliers and complementarity: Define the penalty-induced multipliers by
\begin{ResizedEq}
\lambda_i^{e,j}\triangleq
&-\FPP{V_s}{x}(\EpsFJ-h^{e,j}(\theta_i^\star),\EpsFJ)+\\
&\FPP{V_s}{x}(\EpsFJ+h^{e,j}(\theta_i^\star),\EpsFJ)
\end{ResizedEq}
and
\begin{equation}
\lambda_i^{i,j}\triangleq-\FPP{V_s}{x}(h^{i,j}(\theta_i^\star),\EpsFJ).
\end{equation}
The inequality penalty is nonincreasing in its constraint argument ($V_s'\leq0$ on its nonzero branch and $V_s'=0$ elsewhere), hence $\lambda_i^{i,j}\geq0$. Define $\lambda_{\gamma(i)}^i\triangleq\Lambda^i$, where $\Lambda^i$ is the semi-implicit penalty force from~\prettyref{eq:penalty-force-generic} evaluated at $\theta_{\gamma(i)}^\star$. Thus $\lambda_{\gamma(i)}^i\geq0$. These are distinct multipliers: $\lambda_i^i$ comes from current penalty gradients, while $\lambda_{\gamma(i)}^i$ comes from previous-step penalty forces. The chain rule gives
\begin{equation}
\nabla_{\theta_i}V(\theta_i^\star,
\EpsFJ)=
\FPP{h^e}{\theta_i}(\theta_i^\star)^T\lambda_i^e-
\FPP{h^i}{\theta_i}(\theta_i^\star)^T\lambda_i^i.
\end{equation}
We now conclude the combined inequality line of~\prettyref{eq:eps-fj-system} at tolerance $\EpsFJAll$. We have $h^{i,j}(\theta_i^\star)>0$ and $\lambda_i^{i,j}\geq0$, both shown above. For the product term: if $h^{i,j}(\theta_i^\star)\geq\EpsFJAll\ (\geq\EpsFJ)$, then $V_s'(h^{i,j}(\theta_i^\star),\EpsFJ)=0$, so $\lambda_i^{i,j}=0$ and $\lambda_i^{i,j}(\EpsFJAll-h^{i,j}(\theta_i^\star))=0$; if $h^{i,j}(\theta_i^\star)<\EpsFJAll$, then $\EpsFJAll-h^{i,j}(\theta_i^\star)>0$ and $\lambda_i^{i,j}\geq0$, so $\lambda_i^{i,j}(\EpsFJAll-h^{i,j}(\theta_i^\star))\geq0$. Hence $\min\!\big(h^{i,j}(\theta_i^\star),\lambda_i^{i,j},\lambda_i^{i,j}(\EpsFJAll-h^{i,j}(\theta_i^\star))\big)\geq0$ for all $j$.

Velocity perturbation: In the one-step problem, $\theta_{\gamma(i)}^\star$ and $\delta_i$ are fixed, and $\dot\theta_i=(\theta_i-\theta_{\gamma(i)}^\star)/\delta_i$. Since $V(\theta_{\gamma(i)}^\star,\EpsFJ)\leq\LHamilton<\infty$, the closeness clauses~\prettyref{it:sd-grad-close} and~\ref{it:sd-vel-close} of~\prettyref{def:S-D}, applied to $\Damp_{\EpsFri,\LHamilton}^V$, give a velocity $\dot{\theta}_i^{\dagger}$ such that
\begin{ResizedEq}
\dist\!\left(
\FPP{\delta_i\Damp_{\EpsFri,\LHamilton}^V(\theta_{\gamma(i)}^\star,
\dot{\theta}_i^\star)}{\theta_i^\star},
\partial_{\theta_i}\left(\delta_i\Damp^V\right)(\theta_{\gamma(i)}^\star,
\dot{\theta}_i^{\dagger})\right)
\leq\delta_i\EpsFri,
\end{ResizedEq}
and
\begin{equation}
\|\dot{\theta}_i^\star-
\dot{\theta}_i^{\dagger}\|\leq\EpsFri\leq\EpsFJAll.
\end{equation}
Choose $\dot\theta_i^{\dagger}=(\theta_i^{\dagger}-\theta_{\gamma(i)}^\star)/\delta_i$ with $\theta_i^{\dagger}$ to be the current-step point that realizes this velocity. This proves the velocity-closeness line.

Stationarity: The gradient-residual assumption for~\prettyref{pb:PS-CRDS} gives
\begin{ResizedEq}
\left\|\nabla_{\theta_i}\Phi_i(\theta_i^\star)+
\nabla_{\theta_i}V(\theta_i^\star,\EpsFJ)+
\FPP{\delta_i\Damp_{\EpsFri,\LHamilton}^V(\theta_{\gamma(i)}^\star,
\dot{\theta}_i^\star)}{\theta_i^\star}\right\|
\leq\EpsCRDS.
\end{ResizedEq}
The smoothed dissipation $\Damp_{\EpsFri,\LHamilton}^V$ approximates the penalized nonsmooth dissipation $\Damp^V$. By the semi-implicit definition, $\Damp^V(\theta_{\gamma(i)}^\star,\dot\theta_i)$ is the CRDS dissipation $\Damp(\theta_{\gamma(i)}^\star,
\dot\theta_i,
\lambda_{\gamma(i)}^i)$ appearing in~\prettyref{eq:eps-fj-system}. Therefore the preceding distance bound gives
\begin{ResizedEq}
\dist\!\Big(&\FPP{\delta_i\Damp_{\EpsFri,\LHamilton}^V(\theta_{\gamma(i)}^\star,
\dot{\theta}_i^\star)}{\theta_i^\star},\\
&\partial_{\theta_i}\left(\delta_i\Damp\right)(\theta_{\gamma(i)}^\star,
\dot{\theta}_i^{\dagger},
\lambda_{\gamma(i)}^i)\Big)
\leq\delta_i\EpsFri.
\end{ResizedEq}
Combining this inequality with the gradient-residual bound and the penalty-gradient identity yields
\begin{ResizedEq}
\dist\!\Big(&-\nabla_{\theta_i}\Phi_i(\theta_i^\star)-
\FPP{h^e}{\theta_i}(\theta_i^\star)^T\lambda_i^e+
\FPP{h^i}{\theta_i}(\theta_i^\star)^T\lambda_i^i,\\
&\partial_{\theta_i}\left(\delta_i\Damp\right)(\theta_{\gamma(i)}^\star,
\dot{\theta}_i^{\dagger},
\lambda_{\gamma(i)}^i)\Big)
\leq\EpsCRDS+
\delta_i\EpsFri\leq\EpsFJAll.
\end{ResizedEq}
This is the stationarity line. Hence all lines of~\prettyref{eq:eps-fj-system} hold with tolerance $\EpsFJAll$, so $\theta_i^\star$ satisfies the $\EpsFJAll$-FJ condition of~\prettyref{def:eps-FJ} with current penalty-gradient multipliers and the separate previous-step dissipation multiplier.
\end{proof}

\subsection{Articulated-Body Specialization}
In the articulated-body setting, the generic penalty-induced inequality force~\prettyref{eq:penalty-force-generic} is instantiated through the separating-plane constraints of~\prettyref{eq:separating-plane}: the per-pair force $\Lambda^i$ is assembled componentwise from the $jh$-side and $j'h'$-side padded signed distances, which enter with opposite signs, weighted by the plane-normal length $\|n^{jh,j'h'}\|$,
\begin{ResizedEq}
\label{eq:penalty-force}
&[\Lambda^i]^{j'h',jhk}=-V_s'(d^{j'h',jhk}-\EpsContactPadding,\EpsFJ)\|n^{jh,j'h'}\|,\\
&[\Lambda^i]^{jh,j'h'k'}=-V_s'(-d^{jh,j'h'k'}-\EpsContactPadding,\EpsFJ)\|n^{jh,j'h'}\|.
\end{ResizedEq}
\begin{lemma}[Penalty-derivative flatness at the contact cutoff]
\label{lem:Vs-flatness}
Consider the quartic penalty $V_s(x,\EpsFJ)=(x-\EpsFJ)^4/x^5$ from~\prettyref{eq:Ps} restricted to its cutoff branch $0<x\leq\EpsFJ$, and write $d_+=(\EpsFJ-x)_+$. Then there exist constants $0<c=c(\EpsFJ)\leq C=C(\EpsFJ)$ and a threshold $d^\star=d^\star(\EpsFJ)>0$ such that, whenever $0<d_+\leq d^\star$,
\begin{ResizedEq}
&c\,d_+^3\leq|V_s'(x,\EpsFJ)|\leq C\,d_+^3\\
&c\,d_+^2\leq V_s''(x,\EpsFJ)\leq C\,d_+^2\\
&|V_s'''(x,\EpsFJ)|\leq C\,d_+.
\end{ResizedEq}
On the closed cutoff branch $0<x\leq\EpsFJ$, one has $V_s'(x,\EpsFJ)\leq0$ and $V_s''(x,\EpsFJ)\geq0$. Both inequalities are strict on the active branch $0<x<\EpsFJ$, and the first three derivatives vanish at $x=\EpsFJ$.
\end{lemma}
\begin{proof}
On the active branch $d_+=\EpsFJ-x>0$, direct differentiation of $V_s=d_+^4x^{-5}$ gives
\begin{ResizedEq}
&V_s'=-4\,d_+^3x^{-5}-5\,d_+^4x^{-6}\\
&V_s''=12\,d_+^2x^{-5}+40\,d_+^3x^{-6}+30\,d_+^4x^{-7}\\
&V_s'''=-24\,d_+x^{-5}-180\,d_+^2x^{-6}-360\,d_+^3x^{-7}-210\,d_+^4x^{-8}.
\end{ResizedEq}
Set $d^\star=\min\{\EpsFJ/2,1\}$. If $0<d_+\leq d^\star$, then $\EpsFJ/2\leq x<\EpsFJ$ and $d_+\leq1$. Factoring $d_+^3$, $d_+^2$, and $d_+$, respectively, in the displayed formulas, and bounding every remaining nonnegative term using $x^{-q}\leq(\EpsFJ/2)^{-q}$ and $d_+^q\leq1$, gives all three upper bounds for a constant $C=C(\EpsFJ)$. The first term in each of the first two formulas, together with $x^{-5}\geq\EpsFJ^{-5}$, gives
\begin{ResizedEq}
&|V_s'|\geq4\EpsFJ^{-5}d_+^3,
&V_s''\geq12\EpsFJ^{-5}d_+^2,
\end{ResizedEq}
so the lower bounds hold, for example, with $c=4\EpsFJ^{-5}$. Every term in $V_s'$ is negative and every term in $V_s''$ is positive when $0<x<\EpsFJ$, proving the strict active-branch signs. At $x=\EpsFJ$ one has $d_+=0$, so the displayed formulas extend with $V_s'=V_s''=V_s'''=0$, proving the closed-branch claims.
\end{proof}

\begin{lemma}[Active-contact plane regularity]
\label{lem:active-plane}
Fix a convex-hull pair $(jh,j'h')$ and define the per-pair normal contact potential and its reduced value
\begin{ResizedEq}
&V^{jh,j'h'}(\theta_{\gamma(i)},p^{jh,j'h'})=\\
&\begin{cases}
V_s(1-\|n^{jh,j'h'}\|^2,\EpsFJ)+\\
\sum_kV_s(d^{j'h',jhk}-\EpsContactPadding,\EpsFJ)+\\
\sum_{k'}V_s(-d^{jh,j'h'k'}-\EpsContactPadding,\EpsFJ)
\end{cases}\\
&V_R^{jh,j'h'}(\theta_{\gamma(i)})=
\inf_{p^{jh,j'h'}}V^{jh,j'h'}(\theta_{\gamma(i)},p^{jh,j'h'}),
\end{ResizedEq}
where $p^{jh,j'h'}=(n^{jh,j'h'},b^{jh,j'h'})$ parameterizes the separating plane of~\prettyref{eq:separating-plane} and $d^{jh,j'h'k'}=\langle n^{jh,j'h'},X^{j'h'k'}\rangle+b^{jh,j'h'}$. The companion family $d^{j'h',jhk}=\langle n^{jh,j'h'},X^{jhk}\rangle+b^{jh,j'h'}$ is defined analogously; both families are signed distances from the same oriented separating plane $p^{jh,j'h'}=(n^{jh,j'h'},b^{jh,j'h'})$, entering with opposite signs on the two hulls---the $jh$-side vertices through $d^{j'h',jhk}-\EpsContactPadding>0$ and the $j'h'$-side vertices through $-d^{jh,j'h'k'}-\EpsContactPadding>0$. Choose $0<\EpsFJ<1$ and any finite level $L<\infty$, and assume the reduced level-set bound $V_R(\theta_{\gamma(i)},\EpsFJ)\leq L$. Since each separating-plane variable $p^{jh,j'h'}$ enters only its own pair summand $V^{jh,j'h'}$ while the remaining nonnegative integrity and joint penalty terms, denoted $V_{\mathrm{rest}}$, depend on no separating-plane variable, the minimization separates exactly across pairs: $V_R=V_{\mathrm{rest}}+\sum_{jh,j'h'}V_R^{jh,j'h'}$. In particular, $V_R\geq\sum_{jh,j'h'}V_R^{jh,j'h'}$ and every per-pair term is nonnegative, so $V_R^{jh,j'h'}(\theta_{\gamma(i)})\leq L$. Then the following hold, where~\prettyref{it:active-plane-normal}, \ref{it:active-plane-both}, and~\ref{it:active-plane-convex} additionally assume $V_R^{jh,j'h'}(\theta_{\gamma(i)})>0$:
\begin{enumerate}[label=(\alph*),ref=(\alph*),leftmargin=*]
\item\label{it:active-plane-attain} the infimum defining $V_R^{jh,j'h'}(\theta_{\gamma(i)})$ is attained at a finite minimizer $p^\star=(n^\star,b^\star)$---so the infimum is in fact a minimum---at which every $V_s$ argument is strictly positive; in particular $\|n^\star\|<1$;
\item\label{it:active-plane-normal} the normal-length penalty is active, $V_s(1-\|n^\star\|^2,\EpsFJ)>0$, hence $\|n^\star\|>\sqrt{1-\EpsFJ}$;
\item\label{it:active-plane-both} both hulls have active-slab vertices: there exist $k$ with $[\Lambda^i]^{j'h',jhk}>0$ and $k'$ with $[\Lambda^i]^{jh,j'h'k'}>0$;
\item\label{it:active-plane-convex} $V^{jh,j'h'}(\theta_{\gamma(i)},\cdot)$ is locally strongly convex at $p^\star$, so $p^\star(\theta_{\gamma(i)})$ is the unique global minimizer and is $C^1$ in $\theta_{\gamma(i)}$; the reduced contact force $\Lambda^i(\theta_{\gamma(i)})$ is therefore unambiguously defined on the active side, and consequently $\bar{n}^{jh,j'h'}$, $P^{jh,j'h'}=I-\bar{n}^{jh,j'h'}[\bar{n}^{jh,j'h'}]^T$, every $d^{jh,j'h'k'}$, and every force
\begin{ResizedEq}
[\Lambda^i]^{jh,j'h'k'}=-V_s'(-d^{jh,j'h'k'}-\EpsContactPadding,\EpsFJ)\|n^{jh,j'h'}\|
\end{ResizedEq}
are $C^1$ functions of $\theta_{\gamma(i)}$ on the active side $V_R^{jh,j'h'}(\theta_{\gamma(i)})>0$ (the $C^1$ matching across the zero-force boundary $V_R^{jh,j'h'}=0$ is established separately in~\prettyref{lem:boundary-flatness}).
\end{enumerate}
\end{lemma}
\begin{proof}[Proof of~\prettyref{lem:active-plane}]
We first prove~\prettyref{it:active-plane-attain}. By the reduced level-set bound and nonnegativity of the per-pair terms, $V_R^{jh,j'h'}(\theta_{\gamma(i)})\leq L<\infty$. Fix any $C>L$. For any $p=(n,b)$ with $V^{jh,j'h'}(\theta_{\gamma(i)},p)\leq C$, finiteness of every summand---since $V_s(x,\EpsFJ)=+\infty$ for $x\leq0$ and finite for $x>0$---forces all $V_s$ arguments strictly positive: the normal-length argument satisfies $1-\|n\|^2>0$, so $\|n\|<1$; and, for every vertex, $d^{j'h',jhk}-\EpsContactPadding>0$ and $-d^{jh,j'h'k'}-\EpsContactPadding>0$. By clause~\ref{it:contact-hull} of~\prettyref{ass:abd-mass}, each hull is the convex hull of a finite, nonempty family of fixed local-frame vertices, so each side contributes at least one vertex whose world-frame coordinate $X^{jhk}=R^jx^{jhk}+t^j$ is finite; since $d^{jh,j'h'k'}=\langle n,X^{j'h'k'}\rangle+b$ is affine in $b$ with unit coefficient and $\|n\|<1$, the $jh$-vertex inequality bounds $b$ from below and the $j'h'$-vertex inequality bounds $b$ from above. Hence the sublevel set $\{p:V^{jh,j'h'}(\theta_{\gamma(i)},p)\leq C\}$ is bounded; it is closed because $V^{jh,j'h'}$, a sum of lower semicontinuous extended-real functions $V_s$, is lower semicontinuous. A bounded closed set in finite dimensions is compact, and this set is nonempty since $C>L\geq V_R^{jh,j'h'}(\theta_{\gamma(i)})$; so by Weierstrass $V^{jh,j'h'}(\theta_{\gamma(i)},\cdot)$ attains its minimum over this set at a finite $p^\star=(n^\star,b^\star)$, and this minimum equals the infimum $V_R^{jh,j'h'}(\theta_{\gamma(i)})$. Since $V^{jh,j'h'}(\theta_{\gamma(i)},p^\star)=V_R^{jh,j'h'}(\theta_{\gamma(i)})\leq L<\infty$ is finite, all its $V_s$ arguments are strictly positive, and in particular $\|n^\star\|<1$. The remaining items assume $V_R^{jh,j'h'}(\theta_{\gamma(i)})>0$ and are established at any such finite minimizer $p^\star$.

Since $V_R^{jh,j'h'}(\theta_{\gamma(i)})>0$, at least one penalty term is active at $p^\star$, meaning its $V_s$ argument lies in $(0,\EpsFJ)$; a $jh$-hull vertex $X^{jhk}$ is active when $d^{j'h',jhk}-\EpsContactPadding\in(0,\EpsFJ)$, and a $j'h'$-hull vertex $X^{j'h'k'}$ is active when $-d^{jh,j'h'k'}-\EpsContactPadding\in(0,\EpsFJ)$.

We next prove~\prettyref{it:active-plane-normal}. Write $p_\alpha=\alpha\,p^\star$ and use the radial stationarity of the minimizer,
\begin{ResizedEq}
0=\frac{d}{d\alpha}V^{jh,j'h'}(\theta_{\gamma(i)},p_\alpha)\Big|_{\alpha=1}.
\end{ResizedEq}
Increasing $\alpha$ scales every signed distance $d^{jh,j'h'k'}\to\alpha\,d^{jh,j'h'k'}$---measured from the same oriented plane $p^{jh,j'h'}$, so both hulls' active margins move consistently---which strictly decreases each active vertex penalty---the $jh$-vertex arguments $d^{j'h',jhk}-\EpsContactPadding$ grow more positive and the $j'h'$-vertex arguments $-d^{jh,j'h'k'}-\EpsContactPadding$ grow more positive as well, both receding from the cutoff on which $V_s'<0$ (explicitly, writing $s\in(0,\EpsFJ)$ for the corresponding active $V_s$ argument, $\partial_\alpha(d^{j'h',jhk}-\EpsContactPadding)=d^{j'h',jhk}=s+\EpsContactPadding\geq\EpsContactPadding>0$ on the $jh$-side and $\partial_\alpha(-d^{jh,j'h'k'}-\EpsContactPadding)=-d^{jh,j'h'k'}=s+\EpsContactPadding\geq\EpsContactPadding>0$ on the $j'h'$-side)---while it raises $\|n^\star\|$. Suppose the normal-length penalty were inactive, i.e. $1-\|n^\star\|^2\geq\EpsFJ$, so $V_s=V_s'=0$ there and the normal-length term contributes nothing to the derivative. Then the derivative above would consist solely of the strictly negative vertex contributions (some vertex is active because $V_R^{jh,j'h'}>0$), yielding $\tfrac{d}{d\alpha}V^{jh,j'h'}<0$ strictly, contradicting stationarity. Hence $V_s(1-\|n^\star\|^2,\EpsFJ)>0$, so $1-\|n^\star\|^2<\EpsFJ$ and hence $\|n^\star\|>\sqrt{1-\EpsFJ}$ for $0<\EpsFJ<1$.

Next we prove~\prettyref{it:active-plane-both}. Since the normal-length term is active by~\prettyref{it:active-plane-normal}, its radial derivative along $p_\alpha=\alpha\,p^\star$ is strictly positive; radial stationarity therefore requires at least one active vertex penalty term to supply the balancing negative contribution. To see that the active vertices must occur on both sides, use the offset stationarity $\FPP{V^{jh,j'h'}}{b^{jh,j'h'}}=0$: since $\FPP{d^{jh,j'h'k'}}{b^{jh,j'h'}}=1$ and $V_s'<0$ on the cutoff branch, each active $jh$-vertex term contributes $V_s'(d^{j'h',jhk}-\EpsContactPadding,\EpsFJ)\cdot(+1)\leq0$, whereas each active $j'h'$-vertex term contributes $V_s'(-d^{jh,j'h'k'}-\EpsContactPadding,\EpsFJ)\cdot(-1)\geq0$. If all active vertices lay on a single side, $\FPP{V^{jh,j'h'}}{b^{jh,j'h'}}$ would be strictly single-signed, contradicting stationarity. Hence both sides carry an active vertex, which by~\prettyref{eq:penalty-force} means $[\Lambda^i]^{j'h',jhk}>0$ for some $k$ and $[\Lambda^i]^{jh,j'h'k'}>0$ for some $k'$.

Finally we prove~\prettyref{it:active-plane-convex}. The Hessian $V_{pp}=\FPPT{V^{jh,j'h'}}{p^{jh,j'h'}}$ decomposes as the normal-length Hessian on the $n^{jh,j'h'}$-block plus the positive-semidefinite vertex Gram term $\sum V_s''(s)\,(\partial_p s)(\partial_p s)^T$, each $V_s''\geq0$ by~\prettyref{lem:Vs-flatness}. Writing $q=1-\|n^{jh,j'h'}\|^2$, the normal-length Hessian is
\begin{ResizedEq}
&\nabla_n^2 V_s(q,\EpsFJ)\\
=&4V_s''(q,\EpsFJ)\,n^{jh,j'h'}[n^{jh,j'h'}]^T-2V_s'(q,\EpsFJ)I.
\end{ResizedEq}
By~\prettyref{it:active-plane-normal} the normal-length penalty is active, so $V_s'(q,\EpsFJ)<0$ and $-2V_s'(q,\EpsFJ)I\succ0$; hence the $n^{jh,j'h'}$-block is strictly positive definite. For any nonzero $(\delta n,\delta b)$: if $\delta n\neq0$, the $n^{jh,j'h'}$-block yields a strictly positive quadratic form; if $\delta n=0$,~\prettyref{it:active-plane-both} supplies an active vertex with $V_s''>0$ and $\partial_b s=\pm1$, so $\sum V_s''(s)(\partial_b s)^2\,\delta b^2>0$. Thus $V_{pp}\succ0$, i.e. $V^{jh,j'h'}(\theta_{\gamma(i)},\cdot)$ is locally strongly convex at $p^\star$ and the minimizer is locally unique. Because $V^{jh,j'h'}(\theta_{\gamma(i)},\cdot)$ is convex in $p^{jh,j'h'}$ by~\prettyref{cor:convexity}, its set of global minimizers is convex; the strict local minimizer $p^\star$ is an isolated point of that set, so the set is the singleton $\{p^\star\}$ and $p^\star$ is the unique global minimizer. Hence the reduced force $\Lambda^i(\theta_{\gamma(i)})$ is unambiguously defined on the active side. Before invoking the implicit function theorem we verify its hypotheses on an open neighborhood of $(\theta_{\gamma(i)},p^\star)$ contained in the active region: by clause~\ref{it:contact-hull} of~\prettyref{ass:abd-mass} the local vertices $x^{jhk}$ are fixed model data, so each world-frame vertex $X^{jhk}=R^jx^{jhk}+t^j$ is affine in the maximal coordinates; the penalty arguments $1-\|n^{jh,j'h'}\|^2$ and $\pm d-\EpsContactPadding$ are affine or smooth in $(\theta_{\gamma(i)},p^{jh,j'h'})$; the finite sum over the fixed vertex families preserves this regularity; every active $V_s$ argument stays strictly positive on this neighborhood, where $V_s\in C^3$; and the established $V_{pp}\succ0$ makes the derivative of the stationarity map in $p^{jh,j'h'}$ nonsingular. Hence the stationarity map $\nabla_{p^{jh,j'h'}}V^{jh,j'h'}(\theta_{\gamma(i)},p^{jh,j'h'})$ is jointly $C^1$ there. Applying the implicit function theorem to the stationarity condition $0=\FPP{V^{jh,j'h'}}{p^{jh,j'h'}}(\theta_{\gamma(i)},p^\star)$ shows that $p^\star(\theta_{\gamma(i)})$ is $C^1$, and after shrinking the neighborhood the resulting local branch remains within the active region. Finally, since $(x-\EpsFJ)^4$ vanishes to fourth order at the cutoff, $V_s\in C^3$ there and $V_s'\in C^1$; combined with $\|n^\star\|\geq\sqrt{1-\EpsFJ}>0$ from~\prettyref{it:active-plane-normal}, the maps $\bar{n}^{jh,j'h'}=n^{jh,j'h'}/\|n^{jh,j'h'}\|$, $P^{jh,j'h'}=I-\bar{n}^{jh,j'h'}[\bar{n}^{jh,j'h'}]^T$, each $d^{jh,j'h'k'}$, and each force $[\Lambda^i]^{jh,j'h'k'}=-V_s'(-d^{jh,j'h'k'}-\EpsContactPadding,\EpsFJ)\|n^{jh,j'h'}\|$ inherit $C^1$ dependence on $\theta_{\gamma(i)}$ by the chain rule.
\end{proof}

\begin{lemma}[Boundary flatness of the reduced dissipation gradient]
\label{lem:boundary-flatness}
Fix a convex-hull pair $(jh,j'h')$, choose $0<\EpsFJ<1$, and fix any smoothing scale $\mu>0$. Let a finite proof-local level $C$ and a boundary point $\theta_{\gamma(i)}\in\Omega_C\triangleq\{\theta:V_R(\theta,\EpsFJ)\leq C\}$ satisfy $V_R^{jh,j'h'}(\theta_{\gamma(i)})=0$; the level $C$ serves only to localize the compactness and derivative estimates below and neither determines nor alters the fixed scale $\mu$. Consider active-side points $\theta_{\gamma(i)}'$ with $V_R^{jh,j'h'}(\theta_{\gamma(i)}')>0$ and $r\triangleq\|\theta_{\gamma(i)}'-\theta_{\gamma(i)}\|\to0$; let $p^\star(\theta_{\gamma(i)}')$ be the unique active-side separating-plane minimizer of~\prettyref{lem:active-plane}, $[\Lambda^i]^{jh,j'h'k'}$ the induced penalty forces, $\varpi^\star(\theta_{\gamma(i)}',\dot{\theta}_i)$ the active-side plane-velocity minimizer, and $G\triangleq\FPP{\widehat\Damp_{R,\mu}^{V,jh,j'h'}}{\dot{\theta}_i}$. Let $\rho$ be the maximal active slab-penetration depth over the vertices of both hulls, as defined in the proof. Then there exist local constants $c_p,C_p>0$ such that, along the active side sufficiently near the boundary,
\begin{ResizedEq}
\FPPT{V^{jh,j'h'}}{p^{jh,j'h'}}(\theta_{\gamma(i)}',p^\star(\theta_{\gamma(i)}'))
\succeq c_p\rho^3I,
&\qquad \|Dp^\star(\theta_{\gamma(i)}')\|\leq C_p\rho^{-1}.
\end{ResizedEq}
Moreover, uniformly on every bounded set of $\dot{\theta}_i$, $G(\theta_{\gamma(i)}',\cdot)\to0$ and $DG(\theta_{\gamma(i)}',\cdot)\to0$ as $r\to0$. Consequently $G$, extended by zero on the inactive side $V_R^{jh,j'h'}=0$, is $C^1$ across the boundary.
\end{lemma}
\begin{proof}[Proof of~\prettyref{lem:boundary-flatness}]
We first note that, at the boundary point $\theta_{\gamma(i)}$,
\begin{ResizedEq}
\FPP{\widehat\Damp_{R,\mu}^{V,jh,j'h'}}{\dot{\theta}_i}=0,\quad
\FPPT{\widehat\Damp_{R,\mu}^{V,jh,j'h'}}{\dot{\theta}_i}=0.
\end{ResizedEq}
We next obtain a value estimate without using sensitivity of the minimizing plane velocity. Consider the non-reduced version and its gradient
\begin{ResizedEq}
&\FPP{\widehat\Damp_{\mu}^{V,jh,j'h'}}{\dot{\theta}_i}=\\
&\begin{cases}
\eta\sum_k[\Lambda^i]^{j'h',jhk}\nabla\psi_\mu(v_\parallel^{j'h',jhk})^TP^{jh,j'h'}
\FPP{\dot{X}^{jhk}}{\dot{\theta}_i}+\\
\eta\sum_{k'}[\Lambda^i]^{jh,j'h'k'}\nabla\psi_\mu(v_\parallel^{jh,j'h'k'})^TP^{jh,j'h'}
\FPP{\dot{X}^{j'h'k'}}{\dot{\theta}_i}.
\end{cases}
\end{ResizedEq}
Since $\|\nabla\psi_\mu\|\leq1$, $\|P^{jh,j'h'}\|\leq1$, and the velocity Jacobian has bounded operator norm, we immediately have for any $\theta_{\gamma(i)}'$ and $\varpi$ that $\widehat\Damp_{\mu}^{V,jh,j'h'}$ is Lipschitz continuous in $\dot{\theta}_i$ with Lipschitz modulus being $L(\theta_{\gamma(i)}')=C_1\eta(\sum_k[\Lambda^i]^{j'h',jhk}+\sum_{k'}[\Lambda^i]^{jh,j'h'k'})$ with $C_1$ being a constant upper bound. The fixed coefficient $\eta$ only rescales the local constants and does not change any of the asymptotic rates below. Since minimizing with respect to $\varpi$ preserves the Lipschitz constant, we know that $\widehat\Damp_{R,\mu}^{V,jh,j'h'}$ is also Lipschitz continuous with modulus $L(\theta_{\gamma(i)}')$, which implies
\begin{ResizedEq}
\left\|\FPP{\widehat\Damp_{R,\mu}^{V,jh,j'h'}}{\dot{\theta}_i}(\theta_{\gamma(i)}',\cdot)\right\|
\leq L(\theta_{\gamma(i)}').
\end{ResizedEq}
If we can show that $L(\theta_{\gamma(i)}')=\mathcal{O}(\|\theta_{\gamma(i)}'-\theta_{\gamma(i)}\|^3)$ near $\theta_{\gamma(i)}$, then we have
\begin{ResizedEq}
\left\|\FPP{\widehat\Damp_{R,\mu}^{V,jh,j'h'}}{\dot{\theta}_i}\right\|\leq C_1\|\theta_{\gamma(i)}'-\theta_{\gamma(i)}\|^3,
\end{ResizedEq}
and
\begin{ResizedEq}
&\frac{\left\|\FPP{\widehat\Damp_{R,\mu}^{V,jh,j'h'}}{\dot{\theta}_i}(\theta_{\gamma(i)}+\Delta\theta_{\gamma(i)},\dot{\theta}_i+\Delta\dot{\theta}_i)\right\|}{\|\Delta\theta_{\gamma(i)}\|+\|\Delta\dot{\theta}_i\|}\\
\leq&\frac{C_1\|\Delta\theta_{\gamma(i)}\|^3}{\|\Delta\theta_{\gamma(i)}\|+\|\Delta\dot{\theta}_i\|}\to0,
\end{ResizedEq}
which establishes the value-level flatness. Writing $G\triangleq\FPP{\widehat\Damp_{R,\mu}^{V,jh,j'h'}}{\dot{\theta}_i}$ for the reduced velocity gradient, this yields $\|G(\theta_{\gamma(i)}',\cdot)\|=\mathcal{O}(\|\theta_{\gamma(i)}'-\theta_{\gamma(i)}\|^3)$; hence $G(\theta_{\gamma(i)}',\cdot)\to0$ as $\theta_{\gamma(i)}'\to\theta_{\gamma(i)}$, so $G=0$ at the boundary point and the mixed derivative $\FPPTT{\widehat\Damp_{R,\mu}^{V,jh,j'h'}}{\theta_{\gamma(i)}}{\dot{\theta}_i}$ vanishes there. The full $C^1$ regularity of $G$ is upgraded from this value-level flatness to a derivative-level estimate in the boundary-matching step below.

We only need to show that every $[\Lambda^i]^{jh,j'h'k'}=\mathcal{O}(\|\theta_{\gamma(i)}'-\theta_{\gamma(i)}\|^3)$ on one side of the convex hull and the other side is symmetric. We denote the penetration depth into the active slab as $d^{jh,j'h'k'}_+=\max(\EpsFJ+d^{jh,j'h'k'}+\EpsContactPadding,0)$ and the potential defined on the vertex is $V_s(-d^{jh,j'h'k'}-\EpsContactPadding,\EpsFJ)$. By~\prettyref{lem:Vs-flatness} and the boundedness of the normal vector, we have
\begin{ResizedEq}
&[\Lambda^i]^{jh,j'h'k'}\\
\leq&-V_s'(-d^{jh,j'h'k'}-\EpsContactPadding,\EpsFJ)=\mathcal{O}([d^{jh,j'h'k'}_+]^3).
\end{ResizedEq}
So we only need to show that $d^{jh,j'h'k'}_+(\theta_{\gamma(i)}')=\mathcal{O}(\|\theta_{\gamma(i)}'-\theta_{\gamma(i)}\|)$. Here we face a difficulty similar to the one above, i.e. the optimal separating plane $p^{jh,j'h'}$ might not be unique. We choose a version of the separating plane, denoted by $p_0^{jh,j'h'}$ that achieves $V^{jh,j'h'}(\theta_{\gamma(i)},p^{jh,j'h'})=0$. We fix that separating plane and perturb $\theta_{\gamma(i)}$ to $\theta_{\gamma(i)}'$ to have
\begin{ResizedEq}
&V_R^{jh,j'h'}(\theta_{\gamma(i)}')\\
\leq&V^{jh,j'h'}(\theta_{\gamma(i)}',p_0^{jh,j'h'})=\mathcal{O}(\|\theta_{\gamma(i)}'-\theta_{\gamma(i)}\|^4),
\end{ResizedEq}
where the first inequality is due to the optimality of separating plane at $\theta_{\gamma(i)}'$ and the second bound is due to the definition of $V_s$. Now we focus on one particular term $V_s(-d^{jh,j'h'k'}-\EpsContactPadding,\EpsFJ)$ to derive our desired bound
\begin{ResizedEq}
&\frac{(d^{jh,j'h'k'}_+(\theta_{\gamma(i)}'))^4}{\EpsFJ^5}\\
\leq&\frac{(d^{jh,j'h'k'}_+(\theta_{\gamma(i)}'))^4}{(-d^{jh,j'h'k'}-\EpsContactPadding)^5}\\
=&V_s(-d^{jh,j'h'k'}-\EpsContactPadding,\EpsFJ)\\
\leq&V_R^{jh,j'h'}(\theta_{\gamma(i)}')
=\mathcal{O}(\|\theta_{\gamma(i)}'-\theta_{\gamma(i)}\|^4).
\end{ResizedEq}
From the above inequality we derive our desired result $d^{jh,j'h'k'}_+(\theta_{\gamma(i)}')=\mathcal{O}(\|\theta_{\gamma(i)}'-\theta_{\gamma(i)}\|)$. It remains to upgrade this value-level flatness to a derivative-level estimate, from which the full $C^1$ regularity of $G$ will follow.

Before applying the active-plane regularity result to the approaching points, we place them in one finite total level set. Write $\theta_{\mathrm b}\triangleq\theta_{\gamma(i)}$. Since $V_R(\theta_{\mathrm b},\EpsFJ)\leq C<\infty$, every nonnegative per-pair reduced value is finite and bounded by the total value. For every contact pair $a$, choose at $\theta_{\mathrm b}$ a finite minimizing plane $p_a^0$ using the attainment conclusion established above; that conclusion also applies when the corresponding per-pair reduced value is zero. Finiteness of the per-pair objective at $p_a^0$ forces its normal-length and vertex-slab barrier arguments to be strictly positive, because $V_s(x,\EpsFJ)=+\infty$ for $x\leq0$. Likewise, finiteness of the plane-independent part of $V_R(\theta_{\mathrm b},\EpsFJ)$ forces all equality, physical-inequality, integrity, and joint barrier arguments to be strictly positive at $\theta_{\mathrm b}$.

By the standing finite contact model, there are finitely many contact pairs and vertices with fixed local-frame vertex families; the equality, physical-inequality, integrity, and joint maps also have finitely many components in this finite articulated-body model. With every $p_a^0$ frozen, all of the corresponding barrier arguments are continuous functions of $\theta$. Hence their strict positivity persists simultaneously on one neighbourhood $U$ of $\theta_{\mathrm b}$. Define the unreduced frozen-plane comparison
\begin{equation}
\mathfrak{F}(\theta)\triangleq V(\theta,\{p_a^0\}_a,\EpsFJ),
\end{equation}
where this notation includes every pairwise contact penalty and every plane-independent penalty term. On $U$ all barrier arguments remain positive, so the explicit finite sum defining $\mathfrak{F}$ is finite and continuous. Separability and optimality of the plane minimizations give
\begin{ResizedEq}
V_R(\theta,\EpsFJ)\leq \mathfrak{F}(\theta),
&\qquad \mathfrak{F}(\theta_{\mathrm b})=V_R(\theta_{\mathrm b},\EpsFJ)\leq C.
\end{ResizedEq}
Set $C'\triangleq C+1$. By continuity of $\mathfrak{F}$, after shrinking $U$ we have $\mathfrak{F}(\theta)<C'$ throughout $U$, and therefore $V_R(\theta,\EpsFJ)<C'$ there. In particular, every sufficiently close active-side point lies in this single finite total level set, so the active-plane regularity result applies uniformly with $L=C'$. All local constants below are understood on the common neighbourhood $U\subseteq\Omega_{C'}$.

We work on the active side $V_R^{jh,j'h'}>0$, where the preceding application of~\prettyref{lem:active-plane} provides a locally unique and differentiable optimal plane $p^{jh,j'h'}(\theta_{\gamma(i)})$, abbreviated $p^\star$, and the active-side plane velocity $\varpi^\star(\theta_{\gamma(i)},\dot{\theta}_i)$ is likewise locally unique and differentiable, since the smoothed objective $\widehat\Damp_{\mu}^{V,jh,j'h'}$ is strongly convex in $\varpi$---the plane-velocity damper $\EpsDamp\|\varpi\|^2$ is strongly convex in the entire $\varpi$ (\prettyref{cor:convexity}), uniformly and independent of the vanishing forces, so the first-order condition $\FPP{\widehat\Damp_{\mu}^{V,jh,j'h'}}{\varpi}=0$ defines a unique minimizer $\varpi^\star$ that is differentiable by the implicit function theorem, with the $C^1$ plane supplied by~\prettyref{lem:active-plane}. Throughout we write $\rho\triangleq\max\{\max_k d^{j'h',jhk}_+,\ \max_{k'}d^{jh,j'h'k'}_+\}$ for the maximal slab-penetration depth over the active vertex penalties on both hulls---each depth being $\EpsFJ$ minus the corresponding $V_s$ argument clamped at $0$, with $d^{jh,j'h'k'}_+=\max(\EpsFJ+d^{jh,j'h'k'}+\EpsContactPadding,0)$ on the $j'h'$-side and the symmetric $d^{j'h',jhk}_+=\max(\EpsFJ-(d^{j'h',jhk}-\EpsContactPadding),0)$ on the $jh$-side---together with $r\triangleq\|\theta_{\gamma(i)}'-\theta_{\gamma(i)}\|$, so that the value-level estimate derived above, which holds on both sides, reads $\rho=\mathcal{O}(r)$. Because $\rho$ ranges over both hulls, the phrases ``the deepest active vertex'' and ``the other side is symmetric'' below refer to active vertices on either hull. Throughout this proof $\mu>0$ is held fixed (the smoothing parameter of~\prettyref{lem:smooth-DRp}), and all curvature and derivative constants are local in the following sense: they are taken in a neighbourhood of the boundary point $\theta_{\gamma(i)}$ and are uniform only for $\dot{\theta}_i$ ranging over a bounded set. In particular, on such bounded velocity sets the radial curvature $\psi_\mu''(v)=\mu^2/s^3$ (with $s=\sqrt{\|v\|^2+\mu^2}$) is bounded below by a positive constant, a fact used in Step~4 below.

Step 1: uniform coercivity of the plane Hessian. Abbreviate $V_{pp}\triangleq\FPPT{V^{jh,j'h'}}{p^{jh,j'h'}}$. The entire plane-sensitivity estimate rests on a uniform lower bound for $V_{pp}$, which we isolate and prove as a standalone claim.

\TE{Claim.} Locally near the boundary, for every perturbation $(\delta n^{jh,j'h'},\delta b^{jh,j'h'})$,
\begin{ResizedEq}
&\TWO{\delta n^{jh,j'h'}}{\delta b^{jh,j'h'}}^{T}V_{pp}\TWO{\delta n^{jh,j'h'}}{\delta b^{jh,j'h'}}\\
\geq&c\,\rho^3\left(\|\delta n^{jh,j'h'}\|^2+[\delta b^{jh,j'h'}]^2\right),
\end{ResizedEq}
i.e. $V_{pp}\succeq c\,\rho^3 I$; consequently $\|V_{pp}^{-1}\|=\mathcal{O}(\rho^{-3})$. We deliberately claim only the $\rho^3$ floor, which holds for arbitrary polyhedral contacts, rather than the stronger $\rho^2$ floor that would require the active vertices to span all tangential normal directions.

Proof of Claim. Index the active slab penalties on both hulls by $a\in\mathcal{A}_{\rm slab}$, each with $V_s$ argument $s_a$ and vertex coordinate $X_a$: on the $j'h'$-side $s_a=-d^{jh,j'h'k'}-\EpsContactPadding$ and on the $jh$-side $s_a=d^{j'h',jhk}-\EpsContactPadding$. Each $s_a$ is affine in $p=(n^{jh,j'h'},b^{jh,j'h'})$, with $\partial s_a/\partial b^{jh,j'h'}=\pm1$ and $\partial s_a/\partial n^{jh,j'h'}=\pm X_a$ (sign $-$ on the $j'h'$-side, $+$ on the $jh$-side); the sign is immaterial below since each vertex enters $V_{pp}$ through the rank-one $(\partial_p s_a)(\partial_p s_a)^T$. Hence
\begin{ResizedEq}
V_{pp}=\sum_{a\in\mathcal{A}_{\rm slab}}V_s''(s_a)\,(\partial_p s_a)(\partial_p s_a)^T+V_{nn},
\end{ResizedEq}
where the sum runs over all active slab arguments on both hulls, $V_{nn}=\FPPT{V_s(1-\|n^{jh,j'h'}\|^2,\EpsFJ)}{p^{jh,j'h'}}$ is the normal-length Hessian, supported on the $n^{jh,j'h'}$-block, and each $V_s''(s_a)\geq0$ by~\prettyref{lem:Vs-flatness}, so the vertex Gram term is positive semidefinite. Along the offset direction, the $b^{jh,j'h'}$-$b^{jh,j'h'}$ entry is $\sum_{a}V_s''(s_a)\geq c\,\rho^2$, contributed by the deepest active vertex on either hull through the lower half of~\prettyref{lem:Vs-flatness}; this is a clean, geometry-free floor.

We next certify the normal-length curvature scale. By~\prettyref{lem:active-plane} the normal-length penalty is active at the optimal plane. Write
\begin{ResizedEq}
q_n\triangleq1-\|n^{jh,j'h'}\|^2,
&\rho_n\triangleq\EpsFJ-q_n>0.
\end{ResizedEq}
We first localize $\rho_n$ without invoking any cutoff asymptotics for the normal term. Since $\rho=\mathcal{O}(r)$, we have $\rho\to0$; after restricting to a sufficiently small active-side tail, every active slab depth is at most $\rho\leq d^\star$. Thus~\prettyref{lem:Vs-flatness} applies to the slab arguments and, because the vertex set is finite and $\EpsContactPadding\leq\partial_\alpha s_a=s_a+\EpsContactPadding\leq\EpsFJ+\EpsContactPadding$, their total radial contribution satisfies
\begin{ResizedEq}
\sum_{a\in\mathcal{A}_{\rm slab}}V_s'(s_a,\EpsFJ)\,\partial_\alpha s_a=\mathcal{O}(\rho^3)\to0.
\end{ResizedEq}
On the other hand, radial stationarity along $p_\alpha=\alpha p^\star$ gives
\begin{ResizedEq}
0=&-2\|n^{jh,j'h'}\|^2V_s'(q_n,\EpsFJ)\\
&+\sum_{a\in\mathcal{A}_{\rm slab}}V_s'(s_a,\EpsFJ)\,\partial_\alpha s_a.
\end{ResizedEq}
The bound $\|n^{jh,j'h'}\|^2\geq1-\EpsFJ>0$ therefore implies $|V_s'(q_n,\EpsFJ)|\to0$, without any local estimate for the normal argument.

It remains to identify the limiting normal argument. The normal penalty $V_s(q_n,\EpsFJ)$ is a nonnegative summand of the per-pair minimized value, and the comparison estimate above gives $V_R^{jh,j'h'}(\theta_{\gamma(i)}')=\mathcal{O}(r^4)$, so in particular these values are locally bounded. Since $V_s(q,\EpsFJ)\to+\infty$ as $q\downarrow0$, the sequence $q_n$ stays uniformly away from zero. Any accumulation point in a compact subinterval of $(0,\EpsFJ)$ is impossible because $V_s'(q,\EpsFJ)$ is continuous and strictly negative there on the active branch, whereas $|V_s'(q_n,\EpsFJ)|\to0$. As activity gives $0<q_n<\EpsFJ$, it follows that $q_n\to\EpsFJ$, hence $\rho_n\to0$.

We may now restrict to a smaller tail on which $0<\rho_n\leq d^\star$. Only at this point do the local estimates of~\prettyref{lem:Vs-flatness} give
\begin{ResizedEq}
|V_s'(q_n,\EpsFJ)|\asymp\rho_n^3,
&V_s''(q_n,\EpsFJ)\asymp\rho_n^2.
\end{ResizedEq}
The slab radial total is in fact $\Theta(\rho^3)$: all of its terms have the same negative sign, the deepest active vertex has depth $\rho$ and contributes at least a fixed multiple of $\rho^3$ because $\partial_\alpha s_a\geq\EpsContactPadding>0$, while finiteness of the vertex set and the uniform upper bound on $\partial_\alpha s_a$ give the matching $\mathcal{O}(\rho^3)$ bound. Since $\|n^{jh,j'h'}\|$ is bounded above and bounded away from zero, balancing the normal and slab radial contributions in the stationarity identity yields $\rho_n^3\asymp\rho^3$, and therefore $\rho_n\asymp\rho$.

With $\rho_n\asymp\rho$ in hand we examine the anisotropy of the normal-length Hessian. Writing $q=1-\|n^{jh,j'h'}\|^2$,
\begin{ResizedEq}
&\nabla_n^2 V_s(q,\EpsFJ)\\
=&4V_s''(q,\EpsFJ)\,n^{jh,j'h'}[n^{jh,j'h'}]^T-2V_s'(q,\EpsFJ)I,
\end{ResizedEq}
whose two eigenvalue scales differ sharply: the radial direction $n^{jh,j'h'}$ carries curvature $\asymp V_s''(q)\|n^{jh,j'h'}\|^2\asymp\rho_n^2\asymp\rho^2$, whereas the tangential normal directions carry only the isotropic floor $-2V_s'(q)\asymp\rho_n^3\asymp\rho^3$ (recall $V_s'<0$ on the cutoff branch). Hence, without an extra geometric spanning assumption on the active vertices, $V_{nn}$ guarantees an isotropic floor $V_{nn}\succeq c\,\rho^3 I$ on the $n^{jh,j'h'}$-block.

It remains to combine the offset and normal-length floors over an arbitrary direction $(\delta n^{jh,j'h'},\delta b^{jh,j'h'})$; here the possible cancellation in each active vertex term between $\delta b^{jh,j'h'}$ and $\langle X_a,\delta n^{jh,j'h'}\rangle$ must be handled explicitly. Let $R\triangleq\max_{a\in\mathcal{A}_{\rm slab}}\|X_a\|$ bound the vertex coordinates and fix a split constant $C_{\rm split}>2R$. The vertex Gram term contributes the nonnegative amount $\sum_{a\in\mathcal{A}_{\rm slab}}V_s''(s_a)(\langle X_a,\delta n^{jh,j'h'}\rangle+\delta b^{jh,j'h'})^2$. If $|\delta b^{jh,j'h'}|\leq C_{\rm split}\|\delta n^{jh,j'h'}\|$, drop this nonnegative term and invoke the $n^{jh,j'h'}$-block floor $V_{nn}\succeq c\,\rho^3 I$
\begin{ResizedEq}
&\TWO{\delta n^{jh,j'h'}}{\delta b^{jh,j'h'}}^{T}V_{pp}\TWO{\delta n^{jh,j'h'}}{\delta b^{jh,j'h'}}\\
\geq&c\,\rho^3\|\delta n^{jh,j'h'}\|^2
\geq\frac{c\,\rho^3}{1+C_{\rm split}^2}\left(\|\delta n^{jh,j'h'}\|^2+[\delta b^{jh,j'h'}]^2\right),
\end{ResizedEq}
the final step using $[\delta b^{jh,j'h'}]^2\leq C_{\rm split}^2\|\delta n^{jh,j'h'}\|^2$. If instead $|\delta b^{jh,j'h'}|>C_{\rm split}\|\delta n^{jh,j'h'}\|$, retain only the deepest active vertex $a$ (on whichever hull attains the depth $\rho$); since $|\langle X_a,\delta n^{jh,j'h'}\rangle|\leq R\|\delta n^{jh,j'h'}\|<\tfrac12|\delta b^{jh,j'h'}|$, the cross term cannot cancel the offset, so $|\langle X_a,\delta n^{jh,j'h'}\rangle+\delta b^{jh,j'h'}|\geq\tfrac12|\delta b^{jh,j'h'}|$ and, using $V_s''(s_a)\geq c\,\rho^2$,
\begin{ResizedEq}
&\TWO{\delta n^{jh,j'h'}}{\delta b^{jh,j'h'}}^{T}V_{pp}\TWO{\delta n^{jh,j'h'}}{\delta b^{jh,j'h'}}\\
\geq&\tfrac{c}{4}\,\rho^2[\delta b^{jh,j'h'}]^2
\geq\tfrac{c}{4}\,\rho^3\,\frac{C_{\rm split}^2}{1+C_{\rm split}^2}\left(\|\delta n^{jh,j'h'}\|^2+[\delta b^{jh,j'h'}]^2\right),
\end{ResizedEq}
using $\|\delta n^{jh,j'h'}\|^2<C_{\rm split}^{-2}[\delta b^{jh,j'h'}]^2$ and $\rho^2\geq\rho^3$ near the boundary. In both cases, after renaming constants, the quadratic form is at least $c\,\rho^3(\|\delta n^{jh,j'h'}\|^2+[\delta b^{jh,j'h'}]^2)$, so $V_{pp}\succeq c\,\rho^3 I$ and $\|V_{pp}^{-1}\|=\mathcal{O}(\rho^{-3})$. This is the quantitative form of the active-side nonsingularity established in~\prettyref{lem:active-plane}---the ``uniformly nonsingular normalized system'' underlying the plane-sensitivity bound---made explicit via the two-sided $V_s''$ scale, and completes the proof of the Claim.

Step 2: bounded plane singularity. Differentiating the plane first-order optimality condition $0=\FPP{V^{jh,j'h'}}{p^{jh,j'h'}}(\theta_{\gamma(i)},p^\star)$ in $\theta_{\gamma(i)}$ and solving yields
\begin{ResizedEq}
Dp^\star=-\left(\FPPT{V^{jh,j'h'}}{p^{jh,j'h'}}\right)^{-1}
\FPPTT{V^{jh,j'h'}}{p^{jh,j'h'}}{\theta_{\gamma(i)}}=-V_{pp}^{-1}V_{p\theta},
\end{ResizedEq}
where $V_{p\theta}\triangleq\FPPTT{V^{jh,j'h'}}{p^{jh,j'h'}}{\theta_{\gamma(i)}}$ is a partial derivative at fixed $p$. We now expand this block. For an active slab term $V_s(s_a,\EpsFJ)$,
\begin{ResizedEq}
\partial_{p\theta}^2[V_s(s_a,\EpsFJ)]
=V_s''(s_a,\EpsFJ)(\partial_p s_a)(\partial_\theta s_a)^T
+V_s'(s_a,\EpsFJ)\partial_{p\theta}^2s_a.
\end{ResizedEq}
At fixed $p$, the signed-distance argument is affine in the world vertex $X_a(\theta)$ and in $(n,b)$. In a common neighbourhood of the boundary point, rigid-body kinematics are $C^2$ and the vertices and plane normal are bounded; hence $\partial_p s_a$, $\partial_\theta s_a$, and $\partial_{p\theta}^2s_a$ are uniformly bounded. Therefore the two displayed terms are respectively $\mathcal{O}(\rho^2)$ and $\mathcal{O}(\rho^3)$ by~\prettyref{lem:Vs-flatness}. The normal-length penalty is independent of $\theta$, inactive terms have zero derivatives, and the finite active sum is uniform, so $\|V_{p\theta}\|=\mathcal{O}(\rho^2)$. For the inverse Hessian, the Step~1 Claim gives $\|V_{pp}^{-1}\|=\mathcal{O}(\rho^{-3})$. Combining the two estimates as a norm product,
\begin{ResizedEq}
\|Dp^\star\|\leq\|V_{pp}^{-1}\|\,\|V_{p\theta}\|
=\mathcal{O}(\rho^{-3})\cdot \mathcal{O}(\rho^2)=\mathcal{O}(\rho^{-1}).
\end{ResizedEq}
This plane sensitivity is mildly singular rather than bounded, but the downstream steps show it is amply sufficient, because every mildly singular factor is ultimately multiplied by an active force $[\Lambda^i]_a=\mathcal{O}(\rho^3)$.

Step 3: derivative-level force estimate. For each active vertex, write $[\Lambda^i]^{jh,j'h'k'}=-V_s'(s_{k'},\EpsFJ)\|n^{jh,j'h'}\|$ with $s_{k'}(\theta_{\gamma(i)})=-d^{jh,j'h'k'}(\theta_{\gamma(i)},p^\star(\theta_{\gamma(i)}))-\EpsContactPadding$. By Step 2, the total derivative $Ds_{k'}=\FPP{s_{k'}}{\theta_{\gamma(i)}}+\FPP{s_{k'}}{p^{jh,j'h'}}Dp^\star=\mathcal{O}(\rho^{-1})$ (the plane-sensitivity factor $Dp^\star=\mathcal{O}(\rho^{-1})$ dominates), and $\|n^{jh,j'h'}\|$ is bounded away from $0$ by~\prettyref{lem:active-plane} with $D\|n^{jh,j'h'}\|=\mathcal{O}(\|Dp^\star\|)=\mathcal{O}(\rho^{-1})$. Differentiating the product $[\Lambda^i]^{jh,j'h'k'}=-V_s'(s_{k'},\EpsFJ)\|n^{jh,j'h'}\|$ and using~\prettyref{lem:Vs-flatness},
\begin{ResizedEq}
&\|D[\Lambda^i]^{jh,j'h'k'}\|\\
\leq&|V_s''(s_{k'},\EpsFJ)|\,\|Ds_{k'}\|\,\|n^{jh,j'h'}\|
+|V_s'(s_{k'},\EpsFJ)|\,\mathcal{O}(\rho^{-1})\\
\leq&C\,\rho^2\cdot\rho^{-1}+C\,\rho^3\cdot\rho^{-1}=C\rho+C\rho^2
=\mathcal{O}(\rho)=\mathcal{O}(r),
\end{ResizedEq}
using $\rho=\mathcal{O}(r)$ from above. This is the derivative-level replacement for the value-level estimate $[\Lambda^i]^{jh,j'h'k'}=\mathcal{O}(r^3)$.

Step 4: velocity-gradient derivative. On the active side, the envelope identity $\FPP{\widehat\Damp_{\mu}^{V,jh,j'h'}}{\varpi}=0$ at $\varpi^\star$ lets us evaluate $G$ from the non-reduced gradient displayed above at $\varpi=\varpi^\star$; abbreviating the active summand index by $a$, this reads $G=\eta\sum_a[\Lambda^i]_a\,\nabla\psi_\mu(v_a)^TP^{jh,j'h'}\FPP{\dot{X}_a}{\dot{\theta}_i}$. Differentiating $G$ jointly in $(\theta_{\gamma(i)},\dot{\theta}_i)$ splits into
\begin{ResizedEq}
DG=&\underbrace{\eta\sum_aD[\Lambda^i]_a\left[\nabla\psi_\mu(v_a)^TP^{jh,j'h'}A_a\right]}_{\text{(I)}}\\
&+\underbrace{\eta\sum_a[\Lambda^i]_a\,D\left[\nabla\psi_\mu(v_a)^TP^{jh,j'h'}A_a\right]}_{\text{(II)}}
\end{ResizedEq}
where $A_a\triangleq\FPP{\dot{X}_a}{\dot{\theta}_i}$. For term (I), $D[\Lambda^i]_a=\mathcal{O}(\rho)=\mathcal{O}(r)$ by Step 3 and the bracket is bounded ($\|\nabla\psi_\mu\|\leq1$, $\|P^{jh,j'h'}\|\leq1$, and $A_a$ has bounded operator norm, exactly the bounds used above), so $\text{(I)}=\mathcal{O}(r)$. For term (II), $[\Lambda^i]_a=\mathcal{O}(\rho^3)=\mathcal{O}(r^3)$; the differentiated bracket $D[\nabla\psi_\mu(v_a)^TP^{jh,j'h'}A_a]$ now carries mildly singular factors---$DP^{jh,j'h'}=D(I-\bar{n}^{jh,j'h'}[\bar{n}^{jh,j'h'}]^T)$ and $D\bar{n}^{jh,j'h'}$ scale as $\mathcal{O}(\|Dp^\star\|)=\mathcal{O}(\rho^{-1})$ by Step 2, while $A_a=\FPP{\dot{X}_a}{\dot{\theta}_i}$ remains purely kinematic---together with a dependence on the inner plane velocity $\varpi^\star$. To bound $D\varpi^\star$ we differentiate the first-order optimality condition $\FPP{\widehat\Damp_{\mu}^{V,jh,j'h'}}{\varpi}(\theta_{\gamma(i)},\dot{\theta}_i,\varpi^\star)=0$ jointly in $(\theta_{\gamma(i)},\dot{\theta}_i)$; on the active side, where $\varpi^\star$ is unique and differentiable as established above, this yields
\begin{ResizedEq}
D\varpi^\star=-\left(\FPPT{\widehat\Damp_{\mu}^{V,jh,j'h'}}{\varpi}\right)^{-1}
\FPPTT{\widehat\Damp_{\mu}^{V,jh,j'h'}}{\varpi}{(\theta_{\gamma(i)},\dot{\theta}_i)}.
\end{ResizedEq}
Before bounding the $\varpi$-Hessian we record that, for $\dot{\theta}_i$ in a bounded set, the auxiliary minimizer $\varpi^\star$ stays bounded: comparing with $\varpi=0$ gives $\widehat\Damp_{\mu}^{V,jh,j'h'}(\theta_{\gamma(i)},\dot{\theta}_i,\varpi^\star)\leq C\sum_a[\Lambda^i]_a$, and the coercive plane-velocity damper $\EpsDamp\|\varpi^\star\|^2\leq C\sum_a[\Lambda^i]_a$ bounds the entire $\varpi^\star$ directly, independent of the vanishing forces. Hence each $v_a(\theta_{\gamma(i)},\dot{\theta}_i,\varpi^\star)$ remains in a bounded set, so $\psi_\mu''(v_a)=\mu^2/s^3\geq c_\mu>0$ there, which justifies the lower bound on $\nabla^2\psi_\mu$ used below. We analyze the $\varpi$-Hessian through its block structure. By~\prettyref{eq:tangent-velocity}, $v_a=P^{jh,j'h'}\dot{X}_a-[\bar{n}^{jh,j'h'}\times X_a\,\dot{\omega}^{jh,j'h'}+\dot{\varsigma}^{jh,j'h'}]$, so the velocity Jacobian with respect to the relevant pair block of $\varpi$ is $J_a\triangleq\partial v_a/\partial\varpi=[\,{-}I\ \ {-}(\bar{n}^{jh,j'h'}\times X_a)\,]$; in particular $\partial v_a/\partial\dot{\varsigma}^{jh,j'h'}=-I$ is isotropic. Writing $\Psi\triangleq \widehat\Damp_{\mu}^{V,jh,j'h'}$,
\begin{ResizedEq}
\Psi_{\varpi\varpi}=\eta\sum_a[\Lambda^i]_a\,J_a^T\nabla^2\psi_\mu(v_a)\,J_a
+2\EpsDamp\,I_{\varpi},
\end{ResizedEq}
where $\nabla^2\psi_\mu(v)=\tfrac1s(I-vv^T/s^2)$ with $s=\sqrt{\|v\|^2+\mu^2}$, whose smallest eigenvalue is the radial one $\psi_\mu''(v)=\mu^2/s^3>0$. Since the friction Gram term $\sum_a[\Lambda^i]_a\,J_a^T\nabla^2\psi_\mu(v_a)\,J_a$ is positive semidefinite, the damper supplies the uniform floor $\Psi_{\varpi\varpi}\succeq2\EpsDamp\,I_\varpi$ on the entire plane velocity, giving directly $\|\Psi_{\varpi\varpi}^{-1}\|\leq1/(2\EpsDamp)=\mathcal{O}(1)$; both plane-velocity components are uniformly damped, so the inverse Hessian is bounded independently of the vanishing forces.

We next derive the mixed blocks at fixed $\varpi$, before differentiating the optimizer. The auxiliary first-order condition is
\begin{ResizedEq}
\Psi_\varpi=\eta\sum_a[\Lambda^i]_aJ_a^T\nabla\psi_\mu(v_a)+2\EpsDamp\varpi.
\end{ResizedEq}
For a configuration direction, the product and chain rules give
\begin{ResizedEq}
\Psi_{\varpi\theta}=&\eta\sum_a(D_\theta[\Lambda^i]_a)J_a^T\nabla\psi_\mu(v_a)\\
&+\eta\sum_a[\Lambda^i]_a(D_\theta J_a)^T\nabla\psi_\mu(v_a)\\
&+\eta\sum_a[\Lambda^i]_aJ_a^T\nabla^2\psi_\mu(v_a)D_\theta v_a.
\end{ResizedEq}
Here all derivatives are partial derivatives after substituting $p=p^\star(\theta)$ but holding $\varpi$ fixed. Step~3 gives $D_\theta[\Lambda^i]_a=\mathcal{O}(\rho)$. The formulas for $P$, $\bar n$, $J_a$, and $v_a$ are smooth in the bounded kinematic variables and in $p$ away from $\|n\|=0$; their total configuration derivatives are therefore $\mathcal{O}(1+\|Dp^\star\|)=\mathcal{O}(\rho^{-1})$. Since $\nabla\psi_\mu$ and $\nabla^2\psi_\mu$ are bounded for fixed $\mu>0$, while $[\Lambda^i]_a=\mathcal{O}(\rho^3)$, the three sums are respectively $\mathcal{O}(\rho)$, $\mathcal{O}(\rho^2)$, and $\mathcal{O}(\rho^2)$. Thus $\|\Psi_{\varpi\theta}\|=\mathcal{O}(\rho)$.

For a physical-velocity direction, $[\Lambda^i]_a$, $J_a$, and $P$ are independent of $\dot\theta_i$ at fixed configuration and fixed $\varpi$, while $D_{\dot\theta_i}v_a=P A_a$ is uniformly bounded. Hence
\begin{ResizedEq}
\Psi_{\varpi\dot\theta_i}
=\eta\sum_a[\Lambda^i]_aJ_a^T\nabla^2\psi_\mu(v_a)P A_a
=\mathcal{O}(\rho^3).
\end{ResizedEq}
All constants above are common on a sufficiently small configuration neighbourhood and for $\dot\theta_i$ in the fixed bounded set: the preceding coercivity comparison bounds $\varpi^\star$ uniformly there, and fixed $\mu>0$ bounds the pseudo-Huber derivatives on the resulting velocity set. Substituting these mixed-block estimates into the displayed identity,
\begin{ResizedEq}
&\|D_{\theta_{\gamma(i)}}\varpi^\star\|=\mathcal{O}(1)\cdot \mathcal{O}(\rho)=\mathcal{O}(\rho),\\
&\|D_{\dot{\theta}_i}\varpi^\star\|=\mathcal{O}(1)\cdot \mathcal{O}(\rho^3)=\mathcal{O}(\rho^3).
\end{ResizedEq}
Intuitively, $\Psi_{\varpi\varpi}^{-1}=\mathcal{O}(1)$ while $\Psi_{\varpi\theta_{\gamma(i)}}=\mathcal{O}(\rho)$ and $\Psi_{\varpi\dot{\theta}_i}=\mathcal{O}(\rho^3)$, now backed by the uniform damper floor rather than asserted. Hence every occurrence of $D\varpi^\star$ in $DG$ is multiplied by $[\Lambda^i]_a=\mathcal{O}(\rho^3)$: $[\Lambda^i]_a\,D_{\theta_{\gamma(i)}}\varpi^\star=\mathcal{O}(\rho^3)\cdot \mathcal{O}(\rho)=\mathcal{O}(\rho^4)=\mathcal{O}(r^4)$; likewise the plane-sensitivity factors $DP^{jh,j'h'},D\bar{n}^{jh,j'h'}=\mathcal{O}(\rho^{-1})$ enter multiplied by $[\Lambda^i]_a=\mathcal{O}(\rho^3)$, contributing $\mathcal{O}(\rho^2)$, which is subdominant. Hence $\text{(II)}=\mathcal{O}(r)$. Combining the two terms,
\begin{ResizedEq}
\|DG(\theta_{\gamma(i)}',\dot{\theta}_i')\|\leq C\,r\to0
\end{ResizedEq}
from the active side.

Step 5: boundary matching. On the inactive side $V_R^{jh,j'h'}=0$ the induced forces vanish, so $\widehat\Damp_{R,\mu}^{V,jh,j'h'}$ is identically zero in $\dot{\theta}_i$; hence $G=0$ there and its $\dot{\theta}_i$-derivatives vanish, while the value-level part above gives the mixed derivative $=0$ as well. On the active side, the value-level flatness gives $G(\theta_{\gamma(i)}',\cdot)\to0$ and Step 4 gives $DG(\theta_{\gamma(i)}',\cdot)\to0$, uniformly on every bounded set of $\dot{\theta}_i$, as $r\to0$. Hence $DG$ is continuous across the boundary and matches its zero boundary value, so $G=\FPP{\widehat\Damp_{R,\mu}^{V,jh,j'h'}}{\dot{\theta}_i}$, extended by zero on the inactive side, is $C^1$ across the boundary. Combined with the active-side differentiability of $G$ established above, this proves the lemma.
\end{proof}

\begin{lemma}[Smoothed reduced dissipation: dissipation and smoothness]
\label{lem:smooth-DRp}
Fix any smoothing scale $\mu>0$ and let $\widehat\Damp_{R,\mu}^V$ be the reduced pseudo-Huber dissipation constructed in the proof from the fixed-scale formula depending only on the scalar $\mu$. For every finite level $C$, on the reduced sublevel set $\Omega_C=\{\theta_{\gamma(i)}:V_R(\theta_{\gamma(i)},\EpsFJ)\leq C\}$ the single function $\widehat\Damp_{R,\mu}^V$ is well-defined by finite attained auxiliary plane-velocity minimizers, satisfies the dissipation condition~\prettyref{def:CRDS}~\ref{it:crds-damping} (convexity in $\dot{\theta}_i$ and zero relative velocity as a global minimizer), and has velocity gradient $\FPPR{\widehat\Damp_{R,\mu}^V}{\dot{\theta}_i}$ jointly $C^1$ locally around every point of $\Omega_C$, as required by~\prettyref{def:smooth-CRDS}. The level $C$ is an arbitrary proof-local localization device: it does not enter the formula for $\widehat\Damp_{R,\mu}^V$, which uses the same $\mu$ for every $C$; the local derivative constants may depend on the base point, on $C$, and on the fixed $\mu$. Since every finite-domain configuration lies in some $\Omega_C$, the single fixed-scale function is thereby well-defined and locally smooth throughout $\operatorname{dom}V_R$.
\end{lemma}
\begin{proof}[Proof of~\prettyref{lem:smooth-DRp}]
We constructively define, for the fixed scalar $\mu>0$, the $\mu$-smoothed reduced dissipation $\widehat\Damp_{R,\mu}^V$ of $\Damp_R^V$. We define the pseudo-Huber smooth norm function template: $\psi_\mu(x)=\sqrt{\|x\|^2+\mu^2}-\mu$. We then replace every velocity norm with $\psi_\mu$, while retaining the plane-velocity damper $\EpsDamp(\|\dot{\omega}^{jh,j'h'}\|^2+\|\dot{\varsigma}^{jh,j'h'}\|^2)$ already present in $\Damp_R^V$. Put together, we define
\begin{ResizedEq}
&\widehat\Damp_{R,\mu}^{V,jh,j'h'}(\theta_{\gamma(i)},\dot{\theta}_i)=\\
&\min_{\varpi}\widehat\Damp_{\mu}^{V,jh,j'h'}(\theta_{\gamma(i)},\dot{\theta}_i,\varpi)\\
&\widehat\Damp_{\mu}^{V,jh,j'h'}(\theta_{\gamma(i)},\dot{\theta}_i,\varpi)=\\
&\begin{cases}
\eta\sum_k[\Lambda^i]^{j'h',jhk}\psi_\mu(v_\parallel^{j'h',jhk})+\\
\eta\sum_{k'}[\Lambda^i]^{jh,j'h'k'}\psi_\mu(v_\parallel^{jh,j'h'k'})+\\
\EpsDamp\|\varpi\|^2.
\end{cases}\\
&\widehat\Damp_{R,\mu}^V(\theta_{\gamma(i)},\dot{\theta}_i)=
\sum_{jh,j'h'}\widehat\Damp_{R,\mu}^{V,jh,j'h'}(\theta_{\gamma(i)},\dot{\theta}_i).
\end{ResizedEq}
We will show that for any $\mu>0$, our $\widehat\Damp_{R,\mu}^V$ satisfies the damping condition~\prettyref{def:CRDS}~\ref{it:crds-damping} and the smoothness condition~\prettyref{def:smooth-CRDS} (the closeness clauses of~\prettyref{def:S-D} are established separately in~\prettyref{lem:smooth-DRp-Regularity}). The convexity is due to $\widehat\Damp_{\mu}^{V,jh,j'h'}$ being jointly convex in both $\dot{\theta}_i,\varpi$ and partial minimization with respect to $\varpi$ preserves convexity. For the zero-velocity minimizer condition, each smoothed norm $\psi_\mu(x)=\sqrt{\|x\|^2+\mu^2}-\mu$ is smooth with $\nabla\psi_\mu(0)=0$, and the plane-velocity damper $\EpsDamp\|\varpi\|^2$ has vanishing gradient at $\varpi=0$; since every $v_\parallel$ is linear in $(\dot{\theta}_i,\varpi)$ and the friction weights $\eta[\Lambda^i]\geq0$, the convex $\widehat\Damp_{\mu}^{V,jh,j'h'}(\theta_{\gamma(i)},\cdot,\cdot)$ has $0\in\partial_{(\dot{\theta}_i,\varpi)}\widehat\Damp_{\mu}^{V,jh,j'h'}(\theta_{\gamma(i)},0,0)$. Since $\varpi$ is unconstrained, the partial minimization inherits this stationarity, so $0\in\partial_{\dot{\theta}_i}\widehat\Damp_{R,\mu}^V(\theta_{\gamma(i)},0)$, i.e. zero relative velocity is a global minimizer of the convex $\widehat\Damp_{R,\mu}^V(\theta_{\gamma(i)},\cdot)$. We only need to prove the smoothness of each term $\widehat\Damp_{R,\mu}^{V,jh,j'h'}$. All curvature and derivative estimates below are local: their constants may depend on the bounded neighborhood of $(\theta_{\gamma(i)},\dot{\theta}_i)$ under consideration---in particular, the lower bounds on $\nabla^2\psi_\mu$ hold on the regular set of bounded velocities---which suffices for the local $C^1$ smoothness required by~\prettyref{def:smooth-CRDS}.

Fix an arbitrary finite level $C$ and work on $\Omega_C=\{\theta_{\gamma(i)}:V_R(\theta_{\gamma(i)},\EpsFJ)\leq C\}$. The attainment, active, inactive, and boundary arguments below are carried out on this fixed $\Omega_C$ with $\mu$ held unchanged; the level $C$ localizes these estimates only and does not enter the formula for $\widehat\Damp_{R,\mu}^V$. Because $C$ is arbitrary and $\mu$ is fixed once and for all, establishing joint $C^1$ regularity on each such $\Omega_C$ certifies the single fixed-scale function at every configuration it covers; as every finite-domain configuration lies in some $\Omega_C$, this yields well-definedness and local smoothness throughout $\operatorname{dom}V_R$.

We first verify that the inner minimization over $\varpi$ is attained on $\Omega_C$, so that $\widehat\Damp_{R,\mu}^{V,jh,j'h'}$ is well-defined pointwise. Fix $(\theta_{\gamma(i)},\dot{\theta}_i)$ with $\theta_{\gamma(i)}\in\Omega_C$. When the pair carries an active force, the optimal separating plane $p^\star$ it uses is itself attained and finite by~\prettyref{it:active-plane-attain}, since its level-set hypothesis $V_R^{jh,j'h'}(\theta_{\gamma(i)})\leq C$ holds on $\Omega_C$. The objective $\widehat\Damp_{\mu}^{V,jh,j'h'}(\theta_{\gamma(i)},\dot{\theta}_i,\cdot)$ is continuous, convex, and nonnegative in $\varpi$ by~\prettyref{cor:convexity}, and the plane-velocity damper $\EpsDamp\|\varpi\|^2$ is coercive in the entire $\varpi$. Hence $\widehat\Damp_{\mu}^{V,jh,j'h'}$ is coercive in $\varpi$ for every configuration, independent of whether any force is active, so a finite minimizer $\varpi^\star$ is always attained. If no force is active, $\widehat\Damp_{R,\mu}^{V,jh,j'h'}\equiv0$ and the canonical minimizer $\varpi^\star=0$ is attained. In either case $\widehat\Damp_{R,\mu}^V$ is well-defined by a finite attained minimizer on $\Omega_C$.

The remaining smoothness claim we establish by an exhaustive partition of the configurations $\theta_{\gamma(i)}\in\Omega_C$ into three cases for the pair $(jh,j'h')$.

Case I: We then consider $\theta_{\gamma(i)}$ such that $[\Lambda^i]^{j'h',jhk}=0$ for all $k$ and $[\Lambda^i]^{jh,j'h'k'}=0$ for all $k'$ in an open neighborhood of $\theta_{\gamma(i)}$. Then all forces vanish and the objective reduces to the plane-velocity damper alone, $\widehat\Damp_{\mu}^{V,jh,j'h'}(\theta_{\gamma(i)},\dot{\theta}_i,\varpi)=\EpsDamp\|\varpi\|^2$, which is not literally zero for arbitrary $\varpi$; however its minimum over $\varpi$ is attained at $\varpi=0$ with value $0$, so the reduced value $\widehat\Damp_{R,\mu}^{V,jh,j'h'}(\theta_{\gamma(i)},\dot{\theta}_i)=0$ whenever $\theta_{\gamma(i)}$ lies in the same open neighborhood. The smoothness in both $\theta_{\gamma(i)}$ and $\dot{\theta}_i$ is trivial.

Case II: We then consider $\theta_{\gamma(i)}$ such that $[\Lambda^i]^{j'h',jhk}>0$ for some $k$, whence $V_R^{jh,j'h'}(\theta_{\gamma(i)})>0$. By~\prettyref{lem:active-plane}, the other convex hull also carries an active-slab vertex, i.e. $[\Lambda^i]^{jh,j'h'k'}>0$ for some $k'$; the optimal separating plane $p^{jh,j'h'}(\theta_{\gamma(i)})$ is locally unique and $C^1$ with $\|n^{jh,j'h'}\|\geq\sqrt{1-\EpsFJ}$; and $d^{jh,j'h'k'}$, $n^{jh,j'h'}$, and $[\Lambda^i]^{jh,j'h'k'}$ are $C^1$ functions of $\theta_{\gamma(i)}$.

Now we turn to analyzing the frictional damping. Because of the additional plane-velocity damper, $\widehat\Damp_{\mu}^{V,jh,j'h'}$ is jointly convex in $(\dot{\theta}_i,\varpi)$ and, in fact, strongly convex in $\varpi$ uniformly: the regularizer $\EpsDamp\|\varpi\|^2$ is strongly convex in the entire $\varpi$ with the positive lower bound $2\EpsDamp$ independent of the vanishing forces, while the smoothed friction terms $\psi_\mu(v_\parallel^\cdot)$ add a further positive-semidefinite curvature.

Define the velocity-stationarity map
\begin{ResizedEq}
\Xi(\theta_{\gamma(i)},\dot{\theta}_i,\varpi)\triangleq\FPP{\widehat\Damp_{\mu}^{V,jh,j'h'}}{\varpi}.
\end{ResizedEq}
On the active side $\Xi$ is $C^1$ jointly in $(\theta_{\gamma(i)},\dot{\theta}_i,\varpi)$: by~\prettyref{lem:active-plane} the plane $p^\star$ and hence $n^{jh,j'h'},\bar{n}^{jh,j'h'},P^{jh,j'h'}$, the signed distances $d^{jh,j'h'k'}$, and the forces $[\Lambda^i]^{jh,j'h'k'}$ are $C^1$ in $\theta_{\gamma(i)}$; each relative velocity $v_\parallel^\cdot$ of~\prettyref{eq:tangent-velocity} is affine in $(\dot{\theta}_i,\varpi)$ with these $C^1$ coefficients; and $\psi_\mu$ is $C^\infty$. Hence both first-order gradient maps $\FPP{\widehat\Damp_{\mu}^{V,jh,j'h'}}{\varpi}$ and $\FPP{\widehat\Damp_{\mu}^{V,jh,j'h'}}{\dot{\theta}_i}$ are $C^1$ jointly in $(\theta_{\gamma(i)},\dot{\theta}_i,\varpi)$; in particular $\Xi$ is $C^1$. Moreover
\begin{ResizedEq}
\FPP{\Xi}{\varpi}=\FPPT{\widehat\Damp_{\mu}^{V,jh,j'h'}}{\varpi}\succ0
\end{ResizedEq}
by the local strong convexity in $\varpi$, so $\partial_\varpi \Xi$ is nonsingular. Applying the $C^1$ implicit function theorem to $\Xi(\theta_{\gamma(i)},\dot{\theta}_i,\varpi)=0$, the minimizer map $\varpi^\star(\theta_{\gamma(i)},\dot{\theta}_i)$ is locally unique and $C^1$ in $(\theta_{\gamma(i)},\dot{\theta}_i)$.

Finally, by the envelope identity---$\FPP{\widehat\Damp_{\mu}^{V,jh,j'h'}}{\varpi}=\Xi=0$ at $\varpi^\star$---
\begin{ResizedEq}
&\FPP{\widehat\Damp_{R,\mu}^{V,jh,j'h'}}{\dot{\theta}_i}(\theta_{\gamma(i)},\dot{\theta}_i)=\\
&\FPP{\widehat\Damp_{\mu}^{V,jh,j'h'}}{\dot{\theta}_i}(\theta_{\gamma(i)},\dot{\theta}_i,\varpi^\star(\theta_{\gamma(i)},\dot{\theta}_i)),
\end{ResizedEq}
the $\FPP{\widehat\Damp_{\mu}^{V,jh,j'h'}}{\varpi}\,D\varpi^\star$ term dropping because $\Xi=0$. The right-hand side composes the $C^1$ partial gradient $\FPP{\widehat\Damp_{\mu}^{V,jh,j'h'}}{\dot{\theta}_i}$ with the $C^1$ map $(\theta_{\gamma(i)},\dot{\theta}_i)\mapsto(\theta_{\gamma(i)},\dot{\theta}_i,\varpi^\star(\theta_{\gamma(i)},\dot{\theta}_i))$; a composition of $C^1$ maps is $C^1$, so $\FPP{\widehat\Damp_{R,\mu}^{V,jh,j'h'}}{\dot{\theta}_i}$ is $C^1$ jointly in $(\theta_{\gamma(i)},\dot{\theta}_i)$ on the active side.

Case III: Let $\theta_{\gamma(i)}\in\Omega_C$ satisfy $V_R^{jh,j'h'}(\theta_{\gamma(i)})=0$, approached by active configurations $\theta_{\gamma(i)}'$ with $V_R^{jh,j'h'}>0$, and write $G\triangleq\FPP{\widehat\Damp_{R,\mu}^{V,jh,j'h'}}{\dot{\theta}_i}$. On the inactive side the induced forces vanish, so the pairwise contribution to $\widehat\Damp_{R,\mu}^V$ and its $\dot{\theta}_i$-gradient are identically zero. On the active side,~\prettyref{lem:boundary-flatness} gives $G(\theta_{\gamma(i)}',\cdot)\to0$ and $DG(\theta_{\gamma(i)}',\cdot)\to0$ as $\theta_{\gamma(i)}'\to\theta_{\gamma(i)}$; hence $G$, extended by zero on the inactive side, is $C^1$ across the boundary. This $C^1$ boundary matching rests on the strengthened boundary-flatness estimate of the revised~\prettyref{lem:boundary-flatness}---specifically its Step~1 coercivity Claim $V_{pp}\succeq c\,\rho^3 I$ and the resulting plane-sensitivity bound $\|Dp^\star\|=\mathcal{O}(\rho^{-1})$, which drive the derivative-level flatness $DG\to0$. This is the required $C^1$-matching in Case III.

Every configuration $\theta_{\gamma(i)}\in\Omega_C$ falls into exactly one of Case~I (no active force in a neighborhood), Case~II ($V_R^{jh,j'h'}(\theta_{\gamma(i)})>0$), or Case~III (the boundary $V_R^{jh,j'h'}(\theta_{\gamma(i)})=0$ approached by positive points), and the derivatives match continuously across the transitions; hence $\FPP{\widehat\Damp_{R,\mu}^V}{\dot{\theta}_i}$ is $C^1$ jointly in $(\theta_{\gamma(i)},\dot{\theta}_i)$ locally around every point with $\theta_{\gamma(i)}\in\Omega_C$, which, together with the arbitrariness of the finite level $C$ and the fixed scale $\mu$, is exactly the local smoothness of the single fixed-scale dissipation $\widehat\Damp_{R,\mu}^V$ required by~\prettyref{def:smooth-CRDS} at every point of every $\Omega_C$, and the proof is complete.
\end{proof}

\begin{replemma}{lem:smooth-DRp-Regularity}
\label{lem:smooth-DRp-Regularity-formal}
Fix $\EpsFri>0$, a design/closeness level $0<\LHamilton<\infty$, and a prescribed a priori timestep floor $\delta_{\min}>0$. There is a single common smoothing scale $\mu_{\EpsFri,\LHamilton}>0$ and a globally defined dissipation $\Damp_{R,\EpsFri,\LHamilton}^V:\mathbb R^n\times\mathbb R^n\to\mathbb R$ obtained by applying the decreasing mollifier $\mathfrak M$ to the local core $\widehat\Damp_{R,\mu_{\EpsFri,\LHamilton}}^V$ of~\prettyref{lem:smooth-DRp}. The global dissipation satisfies the damping condition of~\prettyref{def:CRDS} at every configuration, and its velocity gradient is jointly $C^1$ on $\mathbb R^n\times\mathbb R^n$ as required by~\prettyref{def:smooth-CRDS}. It is also a smoothed version of $\Damp_R^V$ in the sense of~\prettyref{def:S-D}, with $V_R$ playing the role of $V$: for any $\theta_{\gamma(i)}^\star$ with $V_R(\theta_{\gamma(i)}^\star,\EpsFJ)\leq\LHamilton<\infty$ and any $\dot{\theta}_i^\star$, there exists $\dot{\theta}_i^{\dagger}$ satisfying~\prettyref{it:sd-grad-close} and~\ref{it:sd-vel-close}.
\end{replemma}
\begin{proof}[Proof of~\prettyref{lem:smooth-DRp-Regularity}]
Fix $\EpsFri>0$, the design/closeness level $0<\LHamilton<\infty$, and the prescribed a priori floor $\delta_{\min}>0$ once and for all, and work on the fixed sublevel $\Omega_{\LHamilton}=\{\theta:V_R(\theta,\EpsFJ)\leq\LHamilton\}$. We use the fixed-scale family $\widehat\Damp_{R,\mu}^V$ of~\prettyref{lem:smooth-DRp}, whose formula depends only on the scalar $\mu>0$, and select a single common scale $\mu=\mu_{\EpsFri,\LHamilton}$ at the end of the proof. We keep the fixed friction coefficient $\eta$ explicit, so the plane-velocity damper remains unscaled. We introduce the following set of notations to index the set of vertices in the active slab
\begin{ResizedEq}
\mathcal{A}^+=&\bigcup_{jh,j'h'}\left\{\TWO{[\Lambda^i]^{j'h',jhk}>0}{v_\parallel^{j'h',jhk}}\right\}\\
\mathcal{A}^-=&\bigcup_{jh,j'h'}\left\{\TWO{[\Lambda^i]^{jh,j'h'k'}>0}{v_\parallel^{jh,j'h'k'}}\right\}\\
\TWO{\lambda_r}{v_r}\in\mathcal{A}=&\mathcal{A}^+\cup\mathcal{A}^-\\
v_r=&A_r\dot{\theta}_i-B_r\varpi.
\end{ResizedEq}
Note that we use aggregated notation over all convex hull pairs. Similarly, $\varpi$ is the concatenated separating plane velocities. Here $A_r$ and $B_r$ are the Jacobian of relative planar velocity with respect to the rigid body velocity and the plane velocity, respectively. We can then write
\begin{ResizedEq}
\Damp^V(\theta_{\gamma(i)}^\star,\dot{\theta}_i,\varpi)=&\eta\sum_{\mathcal{A}}\lambda_r\|v_r\|+\EpsDamp\sum_{jh,j'h'}\|\varpi^{jh,j'h'}\|^2\\
\widehat\Damp_\mu^V(\theta_{\gamma(i)}^\star,\dot{\theta}_i,\varpi)=&\eta\sum_{\mathcal{A}}\lambda_r\psi_\mu(v_r)+\EpsDamp\sum_{jh,j'h'}\|\varpi^{jh,j'h'}\|^2\\
\Damp_R^V(\theta_{\gamma(i)}^\star,\dot{\theta}_i)=&\min_{\bar{\varpi}}\Damp^V(\theta_{\gamma(i)}^\star,\dot{\theta}_i,\bar{\varpi})\\
\widehat\Damp_{R,\mu}^V(\theta_{\gamma(i)}^\star,\dot{\theta}_i)=&\min_{\varpi}\widehat\Damp_\mu^V(\theta_{\gamma(i)}^\star,\dot{\theta}_i,\varpi).
\end{ResizedEq}
Because $\|x\|-\mu\leq\psi_\mu(x)\leq\|x\|$, we have the following bound between $\Damp^V$ and $\widehat\Damp_\mu^V$ for any fixed $\varpi$
\begin{ResizedEq}
\label{eq:DV-bound}
    &\Damp^V(\theta_{\gamma(i)}^\star,\dot{\theta}_i^\star,\varpi)-\mu\bar{\Lambda}\\
\leq&\Damp^V(\theta_{\gamma(i)}^\star,\dot{\theta}_i^\star,\varpi)-\mu\eta\sum_\mathcal{A}\lambda_r\\
\leq&\widehat\Damp_\mu^V(\theta_{\gamma(i)}^\star,\dot{\theta}_i^\star,\varpi)\\
\leq&\Damp^V(\theta_{\gamma(i)}^\star,\dot{\theta}_i^\star,\varpi).
\end{ResizedEq}
In the above inequality we use a finite aggregate load ceiling on the fixed sublevel $\Omega_{\LHamilton}$, derived as follows. For $\theta_{\gamma(i)}^\star\in\Omega_{\LHamilton}$ and every active vertex $r$, the barrier value bounds the corresponding penalty argument: $V_s(s_r,\EpsFJ)\leq V_R^{jh,j'h'}(\theta_{\gamma(i)}^\star)\leq V_R(\theta_{\gamma(i)}^\star,\EpsFJ)\leq\LHamilton<\infty$, so $s_r$ stays in a compact band bounded away from the penalty singularity $s=0$, on which $|V_s'(\cdot,\EpsFJ)|$ is finite by~\prettyref{lem:Vs-flatness}; hence each contact load $\lambda_r=|V_s'(s_r,\EpsFJ)|\,\|n^{jh,j'h'}\|\leq|V_s'(s_r,\EpsFJ)|$ (using $\|n^{jh,j'h'}\|\leq1$) is bounded by a finite constant depending only on $\LHamilton$. Because the number of convex-hull pairs $(jh,j'h')$ and of active vertices per pair is finite, this finite incidence bounds their aggregate, giving the weighted contact load $\eta\sum_\mathcal{A}\lambda_r\leq\bar{\Lambda}$ with $\bar{\Lambda}=\bar{\Lambda}(\LHamilton)<\infty$ a finite constant depending on $\LHamilton$ and the fixed $\eta$. Minimizing over $\varpi$ on the left-hand side of the last inequality of~\prettyref{eq:DV-bound}, we have
\begin{ResizedEq}
&\widehat\Damp_{R,\mu}^V(\theta_{\gamma(i)}^\star,\dot{\theta}_i^\star)
=\min_{\varpi}\widehat\Damp_\mu^V(\theta_{\gamma(i)}^\star,\dot{\theta}_i^\star,\varpi)\\
\leq&\widehat\Damp_\mu^V(\theta_{\gamma(i)}^\star,\dot{\theta}_i^\star,\varpi)
\leq \Damp^V(\theta_{\gamma(i)}^\star,\dot{\theta}_i^\star,\varpi).
\end{ResizedEq}
Taking the infimum over $\varpi$ on the right-hand side, we have $\widehat\Damp_{R,\mu}^V\leq \Damp_R^V$. Following a similar argument to minimize $\Damp^V$, we have $\Damp_R^V-\mu\bar{\Lambda}\leq \widehat\Damp_{R,\mu}^V$. In summary, we have $|\Damp_R^V-\widehat\Damp_{R,\mu}^V|\leq\mu\bar{\Lambda}$. Now since $\widehat\Damp_{R,\mu}^V$ is convex and differentiable by~\prettyref{lem:smooth-DRp}, we have
\begin{ResizedEq}
    &\Damp_R^V(\theta_{\gamma(i)}^\star,\dot{\theta}_i)\\
\geq&\widehat\Damp_{R,\mu}^V(\theta_{\gamma(i)}^\star,\dot{\theta}_i)\\
\geq&\widehat\Damp_{R,\mu}^V(\theta_{\gamma(i)}^\star,\dot{\theta}_i')+
\langle\partial_{\dot{\theta}}\widehat\Damp_{R,\mu}^V(\theta_{\gamma(i)}^\star,\dot{\theta}_i'),\dot{\theta}_i-\dot{\theta}_i'\rangle\\
\geq&\Damp_{R}^V(\theta_{\gamma(i)}^\star,\dot{\theta}_i')+
\langle\partial_{\dot{\theta}_i}\widehat\Damp_{R,\mu}^V(\theta_{\gamma(i)}^\star,\dot{\theta}_i'),\dot{\theta}_i-\dot{\theta}_i'\rangle-\mu\bar{\Lambda},
\end{ResizedEq}
which essentially means $\nabla_{\dot{\theta}_i}\widehat\Damp_{R,\mu}^V(\theta_{\gamma(i)}^\star,\dot{\theta}_i')\in\partial_{\mu\bar{\Lambda}}[\Damp_R^V(\theta_{\gamma(i)}^\star,\cdot)](\dot{\theta}_i')$. Before invoking the perturbation argument we dispose of the trivial case $\bar{\Lambda}=0$: then $\eta\sum_\mathcal{A}\lambda_r\leq\bar{\Lambda}=0$ makes the weighted frictional terms vanish, so $\Damp_R^V(\theta_{\gamma(i)}^\star,\cdot)\equiv \widehat\Damp_{R,\mu}^V(\theta_{\gamma(i)}^\star,\cdot)\equiv0$ and both are smooth with vanishing gradients; the required closeness conditions~\prettyref{it:sd-grad-close} and~\ref{it:sd-vel-close} then hold immediately with $\dot{\theta}_i^{\dagger}=\dot{\theta}_i^\star$ (both left-hand sides in the displayed system below are zero), and no division by $\bar{\Lambda}$ is needed. Assume henceforth $\bar{\Lambda}>0$. The function $\Damp_R^V$ is proper, closed, and convex, so by the Br{\o}ndsted-Rockafellar~\cite{bauschke2017convex} theorem for closed convex functions, for any fixed $\EpsFri$, there exists $\dot{\theta}_i^{\dagger}$ such that
\begin{ResizedEq}
\begin{cases}
\dist\!\Big(\nabla_{\theta_i}\left(\delta_i\widehat\Damp_{R,\mu}^V\right)(\theta_{\gamma(i)}^\star,\dot{\theta}_i^\star),
\partial_{\theta_i}\left(\delta_i\Damp_R^V\right)(\theta_{\gamma(i)}^\star,\dot{\theta}_i^{\dagger})\Big)\leq\delta_i\EpsFri\\
\|\dot{\theta}_i^{\dagger}-\dot{\theta}_i^\star\|\leq\EpsFri.
\end{cases}
\end{ResizedEq}
In the remaining case $\bar{\Lambda}>0$ we make the parameterization explicit. The Br{\o}ndsted-Rockafellar theorem is applied in the velocity variable $\dot{\theta}_i$ at fixed $\theta_{\gamma(i)}$. Since $\dot{\theta}_i=(\theta_i-\theta_{\gamma(i)})/\delta_i$, the configuration-gradient scaling in~\prettyref{def:S-D} satisfies, by the chain rule, $\nabla_{\theta_i}(\delta_i \Damp_R^V)=\nabla_{\dot{\theta}_i}\Damp_R^V$ (and likewise for $\widehat\Damp_{R,\mu}^V$), so the $\delta_i$ cancels and the displayed configuration gradients are exactly the velocity gradients. Applied to the $\mu\bar{\Lambda}$-subgradient relation above with velocity-perturbation radius $\EpsFri$, the theorem yields the point $\dot{\theta}_i^{\dagger}$ of the displayed system satisfying $\|\dot{\theta}_i^{\dagger}-\dot{\theta}_i^\star\|\leq\EpsFri$ (its second line) and a gradient discrepancy of at most $\mu\bar{\Lambda}/\EpsFri$, which is exactly the discrepancy appearing in~\prettyref{it:sd-grad-close} of~\prettyref{def:S-D}. We now fix the single common scale for the pair $(\EpsFri,\LHamilton)$. The floor $\delta_{\min}$ is design data fixed before the first smooth solve, independently of the realized mesh and of the smoothing scale. In the adaptive CITO application it is computed algebraically from the finite initial mesh and the smoothing-independent singularity and Hamiltonian threshold formulas before invoking this lemma; the later convergence proof verifies that subdivision respects this declared floor. Choose the single common scale
\begin{ResizedEq}
0<\mu_{\EpsFri,\LHamilton}\leq\frac{\delta_{\min}\EpsFri^2}{\bar{\Lambda}},
\end{ResizedEq}
and, in the degenerate case $\bar{\Lambda}=0$ handled above, any fixed $\mu_{\EpsFri,\LHamilton}>0$, since the approximation discrepancy is then identically zero. For every timestep $\delta_i\geq\delta_{\min}$ this makes the gradient discrepancy $\mu_{\EpsFri,\LHamilton}\bar{\Lambda}/\EpsFri\leq\delta_{\min}\EpsFri\leq\delta_i\EpsFri$, exactly as required by~\prettyref{it:sd-grad-close}, so the one common scale validates every timestep simultaneously. We accordingly fix the local core $\widehat\Damp_{R,\mu_{\EpsFri,\LHamilton}}^V$ uniformly for all $\delta_i\geq\delta_{\min}$. It remains to globalize this local physical core. Using the shared polynomial $p$ of~\prettyref{eq:clamp-map}, define the decreasing damping mollifier
\begin{equation}
\label{eq:reduced-damping-cutoff}
\mathfrak M(s,L)\triangleq
\begin{cases}
1-p(s/L-1)&0\leq s<+\infty,\\
0&s=+\infty.
\end{cases}
\end{equation}
Since $p$ is nondecreasing with values in $[0,1]$, $\mathfrak M$ is nonnegative and nonincreasing. It equals one for $s\leq L$, equals zero for $s\geq2L$ and at $s=+\infty$, and its first three scalar derivatives vanish at the finite transition endpoints. Define the final aggregate reduced damping by
\begin{ResizedEq}
\label{eq:global-reduced-damping}
&\Damp_{R,\EpsFri,\LHamilton}^V(\theta,\dot\theta)\\
\triangleq&
\begin{cases}
\mathfrak M(V_R(\theta,\EpsFJ),\LHamilton)\widehat\Damp_{R,\mu_{\EpsFri,\LHamilton}}^V(\theta,\dot\theta),&V_R(\theta,\EpsFJ)<2\LHamilton,\\
0,&V_R(\theta,\EpsFJ)\geq2\LHamilton\text{ or }V_R(\theta,\EpsFJ)=+\infty.
\end{cases}
\end{ResizedEq}
The first branch is contained in $\operatorname{dom}V_R$, so the local core and its separating-plane construction are never evaluated outside their finite domain. This cutoff globalizes only the aggregate reduced damping; it does not redefine $V_R$, its conservative force, or any physical per-contact load, and no per-contact physical interpretation is assigned to the transition or zero branch.

We verify the global damping contract. Evaluability follows directly from the guarded branches of~\prettyref{eq:global-reduced-damping}. At fixed $\theta$, the finite branch is a nonnegative, velocity-independent scalar multiple of the convex local core, so it is convex in $\dot\theta$ and retains zero velocity as a global minimizer; the zero branch has the same properties trivially.

Write $G=\FPPR{\Damp_{R,\EpsFri,\LHamilton}^V}{\dot\theta}$ and $\widehat G=\FPPR{\widehat\Damp_{R,\mu_{\EpsFri,\LHamilton}}^V}{\dot\theta}$. On the open set $V_R<2\LHamilton$,~\prettyref{lem:twice-differentiability} and~\ref{lem:smooth-DRp} give
\begin{equation}
G(\theta,\dot\theta)=\mathfrak M(V_R(\theta,\EpsFJ),\LHamilton)\widehat G(\theta,\dot\theta),
\end{equation}
which is jointly $C^1$ by the chain and product rules. At a finite point with $V_R=2\LHamilton$, choose a finite proof-local sublevel around it. The local $C^1$ regularity of $\widehat G$ makes $\widehat G$ and its first derivatives locally bounded there, while $\mathfrak M(s,\LHamilton)\to0$ and $\partial_s\mathfrak M(s,\LHamilton)\to0$ as $s\uparrow2\LHamilton$. Hence both $G$ and its first derivatives tend to zero and match the identically zero outer branch. Finally, if $V_R(\theta,\EpsFJ)=+\infty$, extended-real continuity of $V_R$ makes $\{V_R>2\LHamilton\}$ an open neighborhood of $\theta$, on which the damping and $G$ are identically zero. These three regions prove that $G$ is jointly $C^1$ globally. In particular the continuous blocks $\Damp_{\dot\theta\dot\theta}$ and $\Damp_{\theta\dot\theta}$ are bounded on compact regions; no scalar $\Damp_{\theta\theta}$ derivative is required by~\prettyref{def:smooth-CRDS}.

On $\Omega_{\LHamilton}$ the mollifier $\mathfrak M(V_R,\LHamilton)$ is identically one, so the final damping agrees exactly with the local core as a function on that sublevel. Therefore the Br{\o}ndsted-Rockafellar argument above, including both closeness clauses, is unchanged for every $\theta_{\gamma(i)}^\star$ with $V_R(\theta_{\gamma(i)}^\star,\EpsFJ)\leq\LHamilton$. This proves all claims.
\end{proof}

\begin{proof}[Proof of~\prettyref{cor:SF-CRDS}]
We verify conditions~\ref{it:crds-c2}, \ref{it:crds-potential}, \ref{it:crds-damping}, and \ref{it:crds-control} of~\prettyref{def:CRDS} for the reduced smoothed \prettyref{pb:PSR-CRDS}, using the assumed~\prettyref{def:S-D} properties of $\Damp_{R,\EpsFri,\LHamilton}^V$ throughout.

Condition~\ref{it:crds-c2}: The hard constraint maps $h^e$ and $h^i$ of~\prettyref{pb:PSR-CRDS} are empty because equality constraints, physical inequality constraints, and auxiliary separating-plane inequalities are absorbed into the reduced penalty potential $V_R$. Empty constraint maps are vacuously $C^2$ real maps, and the corresponding multiplier terms disappear. The aggregate control work potential $W$ is inherited from the articulated-body construction of~\prettyref{cor:abd-crds}. Since the reduction eliminates only auxiliary separating-plane variables, $W$ remains a $C^2$ real function of the physical configuration $\theta_i$.

Condition~\ref{it:crds-potential}: Define the conservative potential of the reduced unconstrained system by
\begin{equation}
\widetilde U_R(\theta_i)\triangleq U(\theta_i)+V_R(\theta_i,\EpsFJ).
\end{equation}
The original potential $U$ satisfies the CRDS potential condition: it is nonnegative, proper, extended-real continuous, and $C^2$ on $\operatorname{dom}U$. By~\prettyref{lem:twice-differentiability}, the reduced penalty $V_R$ is proper, extended-real continuous, and $C^2$ on its finite domain; moreover $V_R\geq0$ because each scalar penalty $V_s\geq0$ in~\prettyref{eq:Ps} and nonnegativity is preserved under the sum and the min in~\prettyref{eq:reduced-penalty}. Therefore
\begin{equation}
\operatorname{dom}\widetilde U_R=\operatorname{dom}U\cap\operatorname{dom}V_R.
\end{equation}
Note that the finite domains overlap, $\operatorname{dom}U\cap\operatorname{dom}V_R\neq\emptyset$ (automatic in the articulated-body specialization with $U=0$ of~\prettyref{cor:abd-crds}); hence $\widetilde U_R$ is nonnegative, proper, extended-real continuous, and $C^2$ on its nonempty finite domain $\operatorname{dom}U\cap\operatorname{dom}V_R$. Its conservative force on this finite domain is
\begin{ResizedEq}
\widetilde f_R(\theta_i)=&-\FPPR{\widetilde U_R(\theta_i)}{\theta_i}\\
=&-\FPPR{U(\theta_i)}{\theta_i}-\FPPR{V_R(\theta_i,\EpsFJ)}{\theta_i},
\end{ResizedEq}
which is exactly the original conservative force plus the reduced penalty force induced after eliminating the separating planes.

Condition~\ref{it:crds-damping}: The reduced problem has empty hard inequality maps, so under the unchanged definition of $\mathcal F^i$ its hard feasible set is $\mathcal F^i=\mathbb R^n$. The assumed~\prettyref{def:S-D} contract applies at every predecessor configuration in this set: $\Damp_{R,\EpsFri,\LHamilton}^V(\theta_{\gamma(i)},\cdot)$ is everywhere evaluable, convex in $\dot\theta_i$, and globally minimized at zero velocity. Thus the damping condition holds on the full domain required by~\prettyref{def:CRDS}, not merely on $\operatorname{dom}\widetilde U_R$.

Condition~\ref{it:crds-control}: The reduction eliminates only auxiliary separating-plane variables and does not alter the joint-control construction. Thus the control map and the compact control domain are inherited from~\prettyref{cor:abd-crds}, and eliminating separating-plane variables does not change the bounded control force. The same finite constant $\LControlJacobianBound>0$ satisfies $\|c(\theta_i,u_i)\|\leq\LControlJacobianBound$ for every configuration $\theta_i$ and every admissible control $u_i\in\mathcal{D}$.

The system defined by~\prettyref{pb:PSR-CRDS} is therefore a CRDS. It is smooth globally in the sense of~\prettyref{def:smooth-CRDS}, because~\prettyref{def:S-D} requires $\FPPR{\Damp_{R,\EpsFri,\LHamilton}^V}{\dot\theta_i}$ to be jointly $C^1$ on $\mathbb R^n\times\mathbb R^n$. Finally, after eliminating the massless separating-plane variables, the complete configuration variable is the physical maximal-coordinate vector $\theta_i$; by~\prettyref{lem:mass}, whose positive-volume and positive-density hypotheses are supplied by the standing~\prettyref{ass:abd-mass}, its mass matrix satisfies $M\succ0$. Hence the reduced smoothed CRDS is full-rank in the sense of~\prettyref{def:full-rank}.
\end{proof}

\begin{lemma}[Marginal subdifferential lifting]
\label{lem:marginal-subdiff-lifting}
Let $\mathfrak{F}(x,z)$ be a proper convex function and let $f(x)=\min_z\mathfrak{F}(x,z)$. Suppose the minimum is attained at $\bar z$, so $f(x)=\mathfrak{F}(x,\bar z)<\infty$. Then, for every $g\in\partial f(x)$,
\begin{equation}
(g,0)\in\partial \mathfrak{F}(x,\bar z).
\end{equation}
\end{lemma}
In words, the reduced subgradient lifts to an unreduced subgradient whose eliminated-variable component is zero.
\begin{proof}[Proof of~\prettyref{lem:marginal-subdiff-lifting}]
Since $g\in\partial f(x)$,
\begin{equation}
f(y)\geq f(x)+\langle g,y-x\rangle.
\end{equation}
Also, $\mathfrak{F}(y,w)\geq f(y)$ for all $y,w$ by the definition of $f$ as a minimum over the eliminated variable. Using the attainment identity $f(x)=\mathfrak{F}(x,\bar z)$ gives
\begin{ResizedEq}
\mathfrak{F}(y,w)\geq \mathfrak{F}(x,\bar z)+\langle g,y-x\rangle+
\langle 0,w-\bar z\rangle.
\end{ResizedEq}
This is exactly the subgradient inequality for $(g,0)\in\partial \mathfrak{F}(x,\bar z)$.
\end{proof}

\begin{proof}[Proof of~\prettyref{lem:eps-SR-FJ}]
For this proof, let $P_D$ denote the fixed damping projection $\dot\vartheta^D=P_D\dot\vartheta=<\dot\theta,\varpi>$, with adjoint $P_D^T\TWO{a}{b}=\THREE{a}{0_{\dot p}}{b}$. We verify each line of the $\EpsFJAll$-FJ condition in~\prettyref{def:eps-FJ} (recall $\EpsFJAll=\max\{\EpsFJ,\EpsFri,\EpsCRDS+\delta_i\EpsFri\}$ from~\prettyref{lem:eps-S-FJ}), applied to the extended constrained \prettyref{pb:CRDS-ABD-extended}, by constructing the current constraint multipliers and the semi-implicit dissipation multiplier separately.

First choose the separating-plane parameters that lift the reduced penalty at the current and previous timesteps. By the attainment assumption in~\prettyref{eq:reduced-penalty}, there exist $p_i^{\star,jh,j'h'}$ and $p_{\gamma(i)}^{\star,jh,j'h'}$ attaining the two reduced penalty values. Consistent with~\prettyref{eq:extended-config}, the lifted variables also contain the auxiliary separating-plane translation and rotation coordinates
\begin{ResizedEq}
&\vartheta_i^\star\triangleq<\theta_i^\star,
\left\{\THREE{p_i^\star}{\varsigma_i^\star}{\omega_i^\star}^{jh,j'h'}\right\}_{j\neq j'}>\\
&\vartheta_{\gamma(i)}^\star\triangleq<\theta_{\gamma(i)}^\star,
\left\{\THREE{p_{\gamma(i)}^\star}{\varsigma_{\gamma(i)}^\star}{\omega_{\gamma(i)}^\star}^{jh,j'h'}\right\}_{j\neq j'}>.
\end{ResizedEq}
Here the penalty lift fixes only $p_i^\star$ and $p_{\gamma(i)}^\star$; the coordinates $\varsigma_i^\star,\omega_i^\star$ will be specified below so their discrete velocity realizes the auxiliary velocity used by the dissipation reduction. Since $V$ is independent of these auxiliary translation and rotation coordinates,
\begin{ResizedEq}
\label{eq:lift-VR-min}
V_R(\theta_i^\star,\EpsFJ)=&V(\vartheta_i^\star,\EpsFJ)\\
V_R(\theta_{\gamma(i)}^\star,\EpsFJ)=&V(\vartheta_{\gamma(i)}^\star,\EpsFJ).
\end{ResizedEq}
The proof uses the following two envelope identities for the reduced penalty value function at this minimizer
\begin{ResizedEq}
\label{eq:VR-envelope-proof}
&\FPPR{V_R(\theta_i^\star,\EpsFJ)}{\theta_i}
=\FPPR{V(\vartheta_i^\star,\EpsFJ)}{\theta_i}\\
&0=\FPPR{V(\vartheta_i^\star,\EpsFJ)}{p_i^{jh,j'h'}}.
\end{ResizedEq}
The first identity lifts the physical reduced penalty gradient to the unreduced variables; the second is the auxiliary stationarity condition for the current separating-plane minimization.

Define $\lambda_{\gamma(i)}^i\triangleq\Lambda^i(\vartheta_{\gamma(i)}^\star)$ from the unreduced penalty forces in~\prettyref{eq:penalty-force} evaluated at the previous lifted configuration. This nonnegative multiplier is the one used by the semi-implicit dissipation term.

We first verify feasibility. Since $\theta_i^\star$ has finite objective value in~\prettyref{pb:PSR-CRDS},
\begin{equation}
V_R(\theta_i^\star,\EpsFJ)<\infty.
\end{equation}
By~\prettyref{eq:lift-VR-min}, $V(\vartheta_i^\star,\EpsFJ)<\infty$. As in the proof of~\prettyref{lem:eps-FJ}, finiteness of the equality penalties gives, for every equality component $j$,
\begin{equation}
\EpsFJ-h^{e,j}(\theta_i^\star)>0,
\qquad
\EpsFJ+h^{e,j}(\theta_i^\star)>0,
\end{equation}
and hence
\begin{equation}
|h^{e,j}(\theta_i^\star)|<\EpsFJ\leq\EpsFJAll.
\end{equation}
Finiteness of the inequality penalties gives, for every component of the single extended map $h^i(\vartheta_i)$ (with physical components extended trivially from $\theta_i$),
\begin{equation}
h^{i,j}(\vartheta_i^\star)>0.
\end{equation}
Thus the equality-feasibility line $|h^{e,j}(\theta_i^\star)|\leq\EpsFJAll$ of~\prettyref{eq:eps-fj-system} holds, and the primal-feasibility term $h^{i,j}(\vartheta_i^\star)>0$ of the combined inequality line holds for the lifted extended point $\vartheta_i^\star$; the combined line is concluded below at tolerance $\EpsFJAll$.

We next construct the multipliers from the unreduced penalty gradients at the minimizing separating planes. Define
\begin{ResizedEq}
\lambda_i^{e,j}\triangleq
&-\FPP{V_s}{x}(\EpsFJ-h^{e,j}(\theta_i^\star),\EpsFJ)+\\
&\FPP{V_s}{x}(\EpsFJ+h^{e,j}(\theta_i^\star),\EpsFJ),
\end{ResizedEq}
and
\begin{equation}
\lambda_i^{i,j}\triangleq-\FPP{V_s}{x}(h^{i,j}(\vartheta_i^\star),\EpsFJ).
\end{equation}
Since $V_s(\cdot,\EpsFJ)$ is nonincreasing on its nonzero finite branch and constant on its zero branch,
\begin{equation}
\lambda_i^{i,j}\geq0.
\end{equation}
The nonnegativity of $\lambda_{\gamma(i)}^i$ was already established by its construction from the previous-step unreduced penalty forces. The chain rule and~\prettyref{eq:VR-envelope-proof} yield the physical reduced-gradient identity
\begin{ResizedEq}
\FPPR{V_R(\theta_i^\star,\EpsFJ)}{\theta_i}
=&\FPP{h^e}{\theta_i}(\theta_i^\star)^T\lambda_i^e-
\FPP{h^i}{\theta_i}(\vartheta_i^\star)^T\lambda_i^i,
\end{ResizedEq}
and the auxiliary identity
\begin{equation}
0=\FPP{h^i}{p_i^{jh,j'h'}}(\vartheta_i^\star)^T\lambda_i^i,
\qquad\forall jh,j'h'.
\end{equation}
The equality constraints do not depend on the auxiliary separating-plane variables in~\prettyref{pb:CRDS-ABD-extended}.

We now conclude the combined inequality line of~\prettyref{eq:eps-fj-system} at tolerance $\EpsFJAll$ for every physical and auxiliary inequality component. We have $h^{i,j}(\vartheta_i^\star)>0$ and $\lambda_i^{i,j}\geq0$, both shown above. For the product term: if
\begin{equation}
h^{i,j}(\vartheta_i^\star)\geq\EpsFJAll,
\end{equation}
then $h^{i,j}(\vartheta_i^\star)\geq\EpsFJ$, so $V_s$ is locally constant at $h^{i,j}(\vartheta_i^\star)$ and
\begin{equation}
\lambda_i^{i,j}=0,
\end{equation}
whence $\lambda_i^{i,j}(\EpsFJAll-h^{i,j}(\vartheta_i^\star))=0$; otherwise $\EpsFJAll-h^{i,j}(\vartheta_i^\star)>0$ and $\lambda_i^{i,j}\geq0$, so $\lambda_i^{i,j}(\EpsFJAll-h^{i,j}(\vartheta_i^\star))\geq0$. Hence $\min\!\big(h^{i,j}(\vartheta_i^\star),\lambda_i^{i,j},\lambda_i^{i,j}(\EpsFJAll-h^{i,j}(\vartheta_i^\star))\big)\geq0$ for every physical and auxiliary inequality component.

We now verify velocity closeness and stationarity. The gradient residual for~\prettyref{pb:PSR-CRDS} gives
\begin{ResizedEq}
\left\|\nabla_{\theta_i}\Phi_i(\theta_i^\star)+\FPPR{V_R(\theta_i^\star,\EpsFJ)}{\theta_i}+
\FPP{\delta_i\Damp_{R,\EpsFri,\LHamilton}^V(\theta_{\gamma(i)}^\star,\dot\theta_i^\star)}{\theta_i^\star}\right\|
\leq\EpsCRDS.
\end{ResizedEq}
Since $\Damp_{R,\EpsFri,\LHamilton}^V$ is a smoothed version of $\Damp_R^V$ in the sense of~\prettyref{def:S-D} (with $V_R$ in the role of $V$) and $V_R(\theta_{\gamma(i)}^\star,\EpsFJ)\leq\LHamilton<\infty$, its closeness clauses~\prettyref{it:sd-grad-close} and~\ref{it:sd-vel-close} yield a velocity $\dot\theta_i^{\dagger}$ and a reduced nonsmooth subgradient $g_R^{\dagger}\in\partial_{\theta_i}\left(\delta_i\Damp_R^V\right)(\theta_{\gamma(i)}^\star,\dot\theta_i^{\dagger})$ such that
\begin{equation}
\left\|
\FPP{\delta_i\Damp_{R,\EpsFri,\LHamilton}^V(\theta_{\gamma(i)}^\star,\dot\theta_i^\star)}{\theta_i^\star}
-g_R^{\dagger}\right\|\leq\delta_i\EpsFri,
\end{equation}
and
\begin{equation}
\|\dot\theta_i^\star-\dot\theta_i^{\dagger}\|\leq\EpsFri\leq\EpsFJAll.
\end{equation}

Now choose the nonsmooth auxiliary velocity minimizer
\begin{equation}
\label{eq:nonsmooth-varpi-min}
\bar\varpi_i^{\dagger}\in\argmin{\varpi}
\Damp^V(\vartheta_{\gamma(i)}^\star,\dot\theta_i^{\dagger},\varpi).
\end{equation}
This choice does not alter the separating-plane parameters in $\vartheta_i^\star$; it only fixes the eliminated auxiliary velocity used by the dissipation reduction. Writing $\bar\varpi_i^{\dagger}$ componentwise as $(\bar{\dot\varsigma}_i^{\dagger},\bar{\dot\omega}_i^{\dagger})$, now specify the previously free current auxiliary coordinates in $\vartheta_i^\star$ by
\begin{ResizedEq}
\varsigma_i^{\star,jh,j'h'}=&\varsigma_{\gamma(i)}^{\star,jh,j'h'}+\delta_i\bar{\dot\varsigma}_i^{\dagger,jh,j'h'},\\
\omega_i^{\star,jh,j'h'}=&\omega_{\gamma(i)}^{\star,jh,j'h'}+\delta_i\bar{\dot\omega}_i^{\dagger,jh,j'h'}.
\end{ResizedEq}
This is harmless because the conservative penalty $V$ depends on the separating-plane parameters $p_i^\star$, whereas the eliminated auxiliary velocities enter only the dissipation value functions. Choose $\vartheta_i^{\dagger}$ using the same current and predecessor plane parameters $p_i^\star,p_{\gamma(i)}^\star$ as $\vartheta_i^\star$, with physical velocity $\dot\theta_i^{\dagger}$ and the same auxiliary velocity $\bar\varpi_i^{\dagger}$. Thus the starred and dagger points have the identical finite-difference block $\dot p_i^\star=(p_i^\star-p_{\gamma(i)}^\star)/\delta_i$ and the identical $\varpi$ block. Then
\begin{equation}
\dot\vartheta_i^\star-\dot\vartheta_i^{\dagger}
=\THREE{\dot\theta_i^\star-\dot\theta_i^{\dagger}}{0_{\dot p}}{0_{\varpi}},
\end{equation}
and therefore
\begin{equation}
\|\dot\vartheta_i^\star-\dot\vartheta_i^{\dagger}\|
=\|\dot\theta_i^\star-\dot\theta_i^{\dagger}\|
\leq\EpsFri\leq\EpsFJAll.
\end{equation}
This proves the full extended velocity-closeness line. The chosen nonsmooth auxiliary velocity induces
\begin{equation}
\Damp_R^V(\theta_{\gamma(i)}^\star,\dot\theta_i^{\dagger})
=\Damp^V(\vartheta_{\gamma(i)}^\star,\dot\theta_i^{\dagger},\bar\varpi_i^{\dagger}).
\end{equation}
By~\prettyref{lem:marginal-subdiff-lifting}, applied to the convex marginal function defining $\Damp_R^V$, the reduced nonsmooth subgradient lifts to the unreduced nonsmooth subgradient used by the extended stationarity condition
\begin{ResizedEq}
\label{eq:lifted-D-subgradient}
&g_R^{\dagger}
\in\partial_{\theta_i}\left(\delta_i\Damp_R^V\right)(\theta_{\gamma(i)}^\star,\dot\theta_i^{\dagger})\\
\Longrightarrow&
(g_R^{\dagger},0_{\varpi})
\in\partial_{\theta_i,\varpi}\left(\delta_i\Damp^V\right)(\vartheta_{\gamma(i)}^\star,\dot\theta_i^{\dagger},\bar\varpi_i^{\dagger}).
\end{ResizedEq}
Here $(g_R^{\dagger},0_{\varpi})$ is a covector only in the damping-active variables $(\theta_i,\varpi)$; the zero $\varpi$ component is exactly the first-order optimality condition for the nonsmooth inner minimization in~\prettyref{eq:nonsmooth-varpi-min}.

Substituting the physical reduced-penalty identity above into the gradient residual and using the triangle inequality with the $\delta_i\EpsFri$ approximation of $g_R^{\dagger}$ gives
\begin{ResizedEq}
\dist\!\Big(&-\nabla_{\theta_i}\Phi_i(\theta_i^\star)-\FPP{h^e}{\theta_i}(\theta_i^\star)^T\lambda_i^e+
\FPP{h^i}{\theta_i}(\vartheta_i^\star)^T\lambda_i^i,\\
&\partial_{\theta_i}\left(\delta_i\Damp\right)(\vartheta_{\gamma(i)}^\star,\dot\vartheta_i^{\dagger},\lambda_{\gamma(i)}^i)\Big)
\leq\EpsCRDS+\delta_i\EpsFri\leq\EpsFJAll.
\end{ResizedEq}
This is the physical part of the stationarity line. In the ambient complete-velocity space the corresponding damping covector is $P_D^T(g_R^{\dagger},0_{\varpi})=\THREE{g_R^{\dagger}}{0_{\dot p}}{0_{\varpi}}$, so the $\dot p$ damping-force block is exactly zero.

It remains to check the eliminated auxiliary components. From~\prettyref{eq:VR-envelope-proof}, current-plane optimality gives
\begin{equation}
0=\FPP{V(\vartheta_i^\star,\EpsFJ)}{p_i^{jh,j'h'}}
=\FPP{h^i}{p_i^{jh,j'h'}}(\vartheta_i^\star)^T\lambda_i^i.
\end{equation}
From~\prettyref{eq:lifted-D-subgradient}, nonsmooth velocity optimality gives
\begin{equation}
0\in\partial_{\varpi}\Damp^V(\vartheta_{\gamma(i)}^\star,\dot\theta_i^{\dagger},\bar\varpi_i^{\dagger}).
\end{equation}
Therefore the damping-active auxiliary stationarity residual is exactly zero, and the ambient $\dot p$ damping block is independently zero by the $P_D^T$ lift. Combining these zero auxiliary components with the physical residual above proves the full extended stationarity line in~\prettyref{eq:eps-fj-system} for~\prettyref{pb:CRDS-ABD-extended}.

All conditions in~\prettyref{eq:eps-fj-system} hold with tolerance $\EpsFJAll$ for the lifted extended configuration $\vartheta_i^\star$. This is exactly the multiplier split required by~\prettyref{def:eps-FJ}: the current multipliers are induced by unreduced penalty gradients at current minimizing separating-plane parameters, while the dissipation multiplier is induced by previous-step unreduced penalty forces.
\end{proof}

\section{Proofs for Hamiltonian Upper Bound}
\label{appen:Hamiltonian-upper-bound-proofs}

\begin{proof}[Proof of~\prettyref{lem:mollification}]
Fix $\EpsFJ,\LHamilton,\LPosition$ and write
\begin{equation}
q(\theta)=\mathcal M(\theta,\LPosition),\qquad \phi(s)=\mathcal M(s,\LHamilton),
\end{equation}
and define
\begin{ResizedEq}
K_{2\LHamilton,2\LPosition}\triangleq&
\{z:\|z\|\leq2\LPosition,\tilde U(z,\EpsFJ)\leq2\LHamilton\}.
\end{ResizedEq}
Since $\tilde U(\cdot,\EpsFJ)$ is extended-real continuous, its $2\LHamilton$-sublevel set is closed as the preimage of a closed subset of $[0,+\infty]$.
Therefore $K_{2\LHamilton,2\LPosition}$ is closed and bounded, hence compact.
Every point of $K_{2\LHamilton,2\LPosition}$ has finite $\tilde U$ value, so $K_{2\LHamilton,2\LPosition}\subset\operatorname{dom}\tilde U$.
Since $[0,+\infty)$ is open in $[0,+\infty]$, extended-real continuity implies
\begin{equation}
\operatorname{dom}\tilde U=\tilde U^{-1}([0,+\infty))
\end{equation}
is open.
Thus the compact set $K_{2\LHamilton,2\LPosition}$ is compactly contained in the $C^2$ domain of $\tilde U$.
Now
\begin{equation}
\bar U(\theta,\EpsFJ,\LHamilton,\LPosition)
=\phi(\tilde U(q(\theta),\EpsFJ)).
\end{equation}
For the transition branch of $\mathcal M$, write its radial output as
\begin{equation}
h(r)=r+(2L-r)p(r/L-1),\qquad L<r<2L.
\end{equation}
The shared polynomial in~\prettyref{eq:clamp-map} satisfies $p(0)=0$, $p(1)=1$, and its first three derivatives vanish at both endpoints. Direct differentiation therefore gives
\begin{ResizedEq}
&h(L)=L,\quad h'(L)=1,\quad h''(L)=0,\\
&h(2L)=2L,\quad h'(2L)=0,\quad h''(2L)=0.
\end{ResizedEq}
Therefore the transition branch of $\mathcal M$ matches both the identity branch at radius $L$ and the saturated branch at radius $2L$ to second order.
Thus $q$ and $\phi$ are $C^2$ maps with globally bounded first and second derivatives (a $C^2$ radial profile with second-order endpoint matching extends to a $C^2$ map, the transition occurring at a radius bounded away from the origin).
Also, $\|q(\theta)\|\leq2\LPosition$.

If $\|\theta\|\leq\LPosition$, then $q(\theta)=\theta$.
If also $\tilde U(\theta,\EpsFJ)\leq\LHamilton$, then
$\phi(\tilde U(\theta,\EpsFJ))=\tilde U(\theta,\EpsFJ)$.
This proves~\ref{it:moll-identity}.

For~\ref{it:moll-c2}, first consider the open set where $\tilde U(q(\theta),\EpsFJ)<2\LHamilton$.
There $q(\theta)\in K_{2\LHamilton,2\LPosition}$, so $\tilde U$ is $C^2$ at $q(\theta)$.
Hence the composition with $q$ and $\phi$ is $C^2$.
On the open set where $\tilde U(q(\theta),\EpsFJ)>2\LHamilton$, the scalar clamp is saturated and
$\bar U(\theta,\EpsFJ,\LHamilton,\LPosition)=2\LHamilton$.
Hence $\bar U$ is $C^2$ there.
At the transition level $\tilde U(q(\theta),\EpsFJ)=2\LHamilton$, the scalar clamp satisfies
$\phi'(2\LHamilton)=\phi''(2\LHamilton)=0$.
Thus the first and second derivatives of the composition match the constant branch.
The same second-order matching at $\LHamilton$ gives smooth gluing to the identity branch.
Hence $\bar U$ is $C^2$ on $\mathbb R^n$.

For~\ref{it:moll-curv}, let $v$ be any unit vector. Wherever the scalar clamp is not saturated, the chain rule gives
\begin{ResizedEq}
&v^T\FPPT{\bar U(\theta,\EpsFJ,\LHamilton,\LPosition)}{\theta}v\\
=&\phi''(\tilde U(q,\EpsFJ))
\langle\FPPR{\tilde U(q,\EpsFJ)}{q},Dq(\theta)v\rangle^2+\\
&\phi'(\tilde U(q,\EpsFJ))
(Dq(\theta)v)^T\FPPT{\tilde U(q,\EpsFJ)}{q}(Dq(\theta)v)+\\
&\phi'(\tilde U(q,\EpsFJ))
\left\langle\FPPR{\tilde U(q,\EpsFJ)}{q},D^2q(\theta)[v,v]\right\rangle,
\end{ResizedEq}
where $q=q(\theta)$.
This nonsaturated region has $q(\theta)\in K_{2\LHamilton,2\LPosition}$.
Since this set is compact and $\tilde U$ is $C^2$ on a neighborhood of it,
there exist finite constants $C_1,C_2$ satisfying
\begin{ResizedEq}
\left\|\FPPR{\tilde U(z,\EpsFJ)}{z}\right\|\leq C_1,
\left\|\FPPT{\tilde U(z,\EpsFJ)}{z}\right\|\leq C_2
\forall z\in K_{2\LHamilton,2\LPosition}.
\end{ResizedEq}
The global $C^2$ bounds for the clamps give finite constants $C_3,C_4,C_5,C_6$ such that
\begin{equation}
\|Dq\|\leq C_3,
\|D^2q\|\leq C_4,
|\phi'|\leq C_5,
|\phi''|\leq C_6.
\end{equation}
Thus, in the nonsaturated region, applying the triangle inequality together with the Cauchy-Schwarz and submultiplicative bounds to the three terms above yields the two-sided magnitude estimate
\begin{ResizedEq}
&\left|v^T\FPPT{\bar U(\theta,\EpsFJ,\LHamilton,\LPosition)}{\theta}v\right|\\
\leq& C_6C_1^2C_3^2+C_5C_2C_3^2+C_5C_1C_4.
\end{ResizedEq}
In the saturated region the Hessian is zero, and the same bound holds by continuity on the transition.
Define
\begin{ResizedEq}
\LCurvature\triangleq 1+C_6C_1^2C_3^2+C_5C_2C_3^2+C_5C_1C_4.
\end{ResizedEq}
Then $|v^T\FPPT{\bar U}{\theta}v|\leq\LCurvature$ for every unit $v$ and every $\theta$.
This is equivalent to the claimed two-sided matrix bound $-\LCurvature I\preceq\FPPT{\bar U}{\theta}\preceq\LCurvature I$.

For~\ref{it:moll-bound}, the outer clamp $\phi=\mathcal M(\cdot,\LHamilton)$ takes values in $[0,2\LHamilton]$ by~\prettyref{eq:clamp-map} (the radial profile satisfies $p(t)\in[0,1]$ on $[0,1]$, so the transition radius $h(r)=(1-p(t))r+p(t)2\LHamilton$ lies in $[\LHamilton,2\LHamilton]$ and the saturated branch caps at $2\LHamilton$), and $\tilde U\geq0$ is nonnegative, so $0\leq\bar U=\phi(\tilde U(q(\theta),\EpsFJ))\leq2\LHamilton$ for all $\theta$.
For the gradient, the chain rule gives, wherever the scalar clamp is not saturated,
\begin{ResizedEq}
\FPP{\bar U(\theta,\EpsFJ,\LHamilton,\LPosition)}{\theta}
=\phi'(\tilde U(q,\EpsFJ))\,Dq(\theta)^T\FPP{\tilde U(q,\EpsFJ)}{q},
\end{ResizedEq}
with $q=q(\theta)\in K_{2\LHamilton,2\LPosition}$, while in the saturated region $\phi'=0$ makes the gradient vanish.
Using $|\phi'|\leq C_5$, $\|Dq\|\leq C_3$, and $\|\FPPR{\tilde U(z,\EpsFJ)}{z}\|\leq C_1$ on $K_{2\LHamilton,2\LPosition}$, we obtain $\|\FPPR{\bar U}{\theta}\|\leq C_1C_3C_5$ in the nonsaturated region, and the same bound holds by continuity on the transition and trivially (value $0$) in the saturated region.
Define $\LPotentialForce\triangleq C_1C_3C_5$; then $\|\FPPR{\bar U}{\theta}\|\leq\LPotentialForce$ for every $\theta$, proving~\ref{it:moll-bound}.
\end{proof}

\begin{proof}[Proof of~\prettyref{lem:controlwork-mollification}]
Write $q(\theta)=\mathcal M(\theta,\LPosition)$. By the mollifier construction in~\prettyref{eq:clamp-map}, $q$ is $C^2$ with globally bounded first and second derivatives, $\|Dq\|\leq C_3$ with $C_3$ a clamp-independent absolute constant and $\|D^2q\|\leq C_4$, and $\|q(\theta)\|\leq2\LPosition$. By~\prettyref{def:CRDS}~\ref{it:crds-c2}, $W$ is jointly $C^2$ on the open set $\mathbb R^n\times\mathcal{D}_u$ containing $\mathbb R^n\times\mathcal D$, so $\FPPR{W}{\theta}$, $\FPPT{W}{\theta}$, and $\FPPTT{W}{\theta}{u}$ are continuous there; by~\prettyref{def:CRDS}~\ref{it:crds-control}, $\|c(\theta,u)\|=\|\FPPR{W(\theta,u)}{\theta}\|\leq\LControlJacobianBound$ for all $(\theta,u)\in\mathbb R^n\times\mathcal D$. Fix $u\in\mathcal D$.

Property~\ref{it:cw-identity} is immediate: if $\|\theta\|\leq\LPosition$ then $q(\theta)=\theta$, so $\bar W(\theta,u,\LPosition)=W(\theta,u)$.
For~\ref{it:cw-c2}, $q$ is $C^2$ and $W$ is $C^2$ in $\theta$, so the composition $\bar W$ is $C^2$ in $\theta$.
The chain rule gives, for every unit vector $v$ and every $\theta$,
\begin{ResizedEq}
\FPP{\bar W(\theta,u,\LPosition)}{\theta}=&Dq(\theta)^T\FPP{W(q(\theta),u)}{q},\\
v^T\FPPT{\bar W(\theta,u,\LPosition)}{\theta}v=&(Dq(\theta)v)^T\FPPT{W(q(\theta),u)}{q}(Dq(\theta)v)\\
&+\left\langle\FPP{W(q(\theta),u)}{q},D^2q(\theta)[v,v]\right\rangle.
\end{ResizedEq}
Since $q(\theta)\in\mathbb R^n$ and $u\in\mathcal D$, the global CRDS control-force bound gives $\|\FPPR{W(q(\theta),u)}{q}\|=\|c(q(\theta),u)\|\leq\LControlJacobianBound$.
With $\|Dq\|\leq C_3$, we obtain, for every $\theta$ and every $u\in\mathcal D$,
\begin{ResizedEq}
\left\|\FPP{\bar W(\theta,u,\LPosition)}{\theta}\right\|\leq C_3\LControlJacobianBound\triangleq\LControlWorkForce,
\end{ResizedEq}
which proves~\ref{it:cw-force}; since $C_3$ is an absolute mollifier constant and $\LControlJacobianBound$ is clamp-independent, so is $\LControlWorkForce$.
For~\ref{it:cw-curv}, the set $\{q:\|q\|\leq2\LPosition\}\times\mathcal D$ is compact and lies in the open $C^2$ domain $\mathbb R^n\times\mathcal{D}_u$. Hence the joint continuity of $\FPPT{W}{\theta}$ yields a finite constant
\begin{ResizedEq}
C_W\triangleq\sup_{\|q\|\leq2\LPosition,\,u\in\mathcal D}\left\|\FPPT{W(q,u)}{q}\right\|<\infty.
\end{ResizedEq}
Using $\|Dq\|\leq C_3$, $\|D^2q\|\leq C_4$, and $\|\FPPR{W(q(\theta),u)}{q}\|\leq\LControlJacobianBound$ in the Hessian identity above, the triangle and submultiplicative inequalities give
\begin{ResizedEq}
\left|v^T\FPPT{\bar W(\theta,u,\LPosition)}{\theta}v\right|\leq C_3^2C_W+\LControlJacobianBound C_4\triangleq\LControlWorkCurvature
\end{ResizedEq}
for every unit $v$, every $\theta$, and every $u\in\mathcal D$, which is equivalent to the two-sided bound of~\ref{it:cw-curv}; it is clamp-dependent (through $C_W$, whose region $\|q\|\leq2\LPosition$ grows with $\LPosition$, and through $C_4$), but independent of $\theta$ and $u$.
Finally, the joint continuity of $\FPPTT{W}{\theta}{u}$ on the same compact set gives $C_{W2}\triangleq\sup_{\|q\|\leq2\LPosition,\,u\in\mathcal D}\|\FPPTT{W(q,u)}{q}{u}\|<\infty$. The chain-rule identity $\FPPTT{\bar W(\theta,u,\LPosition)}{\theta}{u}=Dq(\theta)^T\FPPTT{W(q(\theta),u)}{q}{u}$ then yields the uniform clamp-dependent bound $\|\FPPTT{\bar W}{\theta}{u}\|\leq C_3C_{W2}$ for all $\theta$ and $u$, proving~\ref{it:cw-mixed}.
\end{proof}

\begin{lemma}[One-step Hamiltonian bound]
\label{lem:H-one-step}
Let~\prettyref{pb:PS-CRDS-M} define a smooth full-rank CRDS (\prettyref{def:CRDS},~\ref{def:full-rank}) and let $\bar U$ be $\LCurvature$-weakly convex as in~\prettyref{lem:mollification}~\ref{it:moll-curv}. If $\theta_i$ solves~\prettyref{pb:PS-CRDS-M} to a gradient residual of norm at most $\EpsCRDS$, then
\begin{ResizedEq}
\label{eq:H-one-step}
&\bar{\mathcal{H}}_i^+-\bar{\mathcal{H}}_{\gamma(i)}^+\\
\leq&\frac{\LCurvature\delta_i^2+2\delta_i}{2\lambda_{\min}(M)}\|\dot{\theta}_i\|_M^2
+\frac{\delta_i}{2}\left((\|g\|+\LControlWorkForce)^2+\EpsCRDS^2\right).
\end{ResizedEq}
\end{lemma}
\begin{proof}[Proof of~\prettyref{lem:H-one-step}]
Since~\prettyref{pb:PS-CRDS-M} is a smooth CRDS (\prettyref{def:smooth-CRDS}), and using the variational-integrator identity $\ddot\theta_i=(\dot\theta_i-\dot\theta_{\gamma(i)})/\delta_i$ implied by~\prettyref{eq:discretization} (so that the kinetic term contributes $\nabla_{\theta_i}[\tfrac{\delta_i^2}{2}\|\ddot\theta_i\|_M^2]=M(\dot\theta_i-\dot\theta_{\gamma(i)})/\delta_i$), the gradient-residual stationarity of $\theta_i$ reads
\begin{ResizedEq}
\label{eq:one-step-stationarity}
M\frac{\dot\theta_i-\dot\theta_{\gamma(i)}}{\delta_i}
=&-\FPP{\bar U}{\theta_i}+g+\bar c(\theta_i,u_i)\\
&-\partial_{\dot\theta_i}\Damp_{\EpsFri,\LHamilton}^V(\theta_{\gamma(i)},\dot\theta_i)+\rho_i,
\end{ResizedEq}
with $\|\rho_i\|\leq\EpsCRDS$, where we have used $\FPPR{\bar W}{\theta_i}=-\bar c(\theta_i,u_i)$ with $\bar c\triangleq-\FPP{\bar W}{\theta}$ the mollified control force (analogous to $c=-\FPP{W}{\theta}$ in~\prettyref{eq:control}), the identities $\dot\theta_i=(\theta_i-\theta_{\gamma(i)})/\delta_i$ and $\delta_i\partial_{\theta_i}\Damp_{\EpsFri,\LHamilton}^V=\partial_{\dot\theta_i}\Damp_{\EpsFri,\LHamilton}^V$, and abbreviate $\bar U(\theta)\triangleq\bar U(\theta,\EpsFJ,\LHamilton,\LPosition)$.
Denote the right-hand side by $s_i$, so that $\dot\theta_i-\dot\theta_{\gamma(i)}=\delta_iM^{-1}s_i$.
Completing the square, the kinetic difference is
\begin{ResizedEq}
&\frac12\|\dot\theta_i\|_M^2-\frac12\|\dot\theta_{\gamma(i)}\|_M^2\\
=&\frac12\|\dot\theta_i\|_M^2-\frac12\|\dot\theta_i-\delta_iM^{-1}s_i\|_M^2\\
=&\delta_i\dot\theta_i^Ts_i-\frac{\delta_i^2}{2}s_i^TM^{-1}s_i\\
=&(\theta_i-\theta_{\gamma(i)})^Ts_i-\frac{\delta_i^2}{2}s_i^TM^{-1}s_i,
\end{ResizedEq}
using $\delta_i\dot\theta_i=\theta_i-\theta_{\gamma(i)}$.
Substituting $s_i$ and collecting terms,
\begin{ResizedEq}
&\bar{\mathcal H}_i^+-\bar{\mathcal H}_{\gamma(i)}^+\\
=&\bar U(\theta_i)-\bar U(\theta_{\gamma(i)})+\frac12\|\dot\theta_i\|_M^2-\frac12\|\dot\theta_{\gamma(i)}\|_M^2\\
=&\left[\bar U(\theta_i)-\bar U(\theta_{\gamma(i)})-(\theta_i-\theta_{\gamma(i)})^T\FPP{\bar U}{\theta_i}\right]\\
&+(\theta_i-\theta_{\gamma(i)})^T(g+\bar c(\theta_i,u_i))\\
&-(\theta_i-\theta_{\gamma(i)})^T\partial_{\dot\theta_i}\Damp_{\EpsFri,\LHamilton}^V\\
&+(\theta_i-\theta_{\gamma(i)})^T\rho_i-\frac{\delta_i^2}{2}s_i^TM^{-1}s_i.
\end{ResizedEq}
We bound the five terms in turn.
For the first, the $\LCurvature$-weak convexity of $\bar U$ (\prettyref{lem:mollification}~\ref{it:moll-curv}) gives
\begin{ResizedEq}
&\bar U(\theta_i)-\bar U(\theta_{\gamma(i)})-(\theta_i-\theta_{\gamma(i)})^T\FPP{\bar U}{\theta_i}\\
\leq&\frac{\LCurvature}{2}\|\theta_i-\theta_{\gamma(i)}\|^2.
\end{ResizedEq}
For the dissipation term,~\prettyref{def:CRDS}~\ref{it:crds-damping} yields
\begin{ResizedEq}
&\dot\theta_i^T\partial_{\dot\theta_i}\Damp_{\EpsFri,\LHamilton}^V(\theta_{\gamma(i)},\dot\theta_i)\\
\geq&\Damp_{\EpsFri,\LHamilton}^V(\theta_{\gamma(i)},\dot\theta_i)-\Damp_{\EpsFri,\LHamilton}^V(\theta_{\gamma(i)},0)\geq0,
\end{ResizedEq}
so $-(\theta_i-\theta_{\gamma(i)})^T\partial_{\dot\theta_i}\Damp_{\EpsFri,\LHamilton}^V=-\delta_i\dot\theta_i^T\partial_{\dot\theta_i}\Damp_{\EpsFri,\LHamilton}^V\leq0$: the dissipation term is dropped.
The negative square $-\frac{\delta_i^2}{2}s_i^TM^{-1}s_i\leq0$ because $M\succ0$.
For the two remaining linear terms, Young's inequality with $\delta_i\dot\theta_i=\theta_i-\theta_{\gamma(i)}$ gives
\begin{ResizedEq}
&(\theta_i-\theta_{\gamma(i)})^T(g+\bar c(\theta_i,u_i))\\
\leq&\frac{\delta_i}{2}\|\dot\theta_i\|^2+\frac{\delta_i}{2}\|g+\bar c(\theta_i,u_i)\|^2,\\
&(\theta_i-\theta_{\gamma(i)})^T\rho_i\\
\leq&\frac{\delta_i}{2}\|\dot\theta_i\|^2+\frac{\delta_i}{2}\EpsCRDS^2.
\end{ResizedEq}
Combining the five bounds,
\begin{ResizedEq}
&\bar{\mathcal H}_i^+-\bar{\mathcal H}_{\gamma(i)}^+\\
\leq&\frac{\LCurvature\delta_i^2+2\delta_i}{2}\|\dot\theta_i\|^2
+\frac{\delta_i}{2}\left(\|g+\bar c(\theta_i,u_i)\|^2+\EpsCRDS^2\right).
\end{ResizedEq}
Finally, $\|g+\bar c(\theta_i,u_i)\|\leq\|g\|+\LControlWorkForce$ by~\prettyref{lem:controlwork-mollification}~\ref{it:cw-force}, and $\|\dot\theta_i\|^2\leq\|\dot\theta_i\|_M^2/\lambda_{\min}(M)$, which give~\prettyref{eq:H-one-step}.
\end{proof}

\begin{proof}[Proof of~\prettyref{lem:H-bound}]
The timestep condition~\prettyref{eq:H-bound-timestep} implies $\delta_i<1/\LCurvature$, hence $\LCurvature\delta_i^2<\delta_i$ and $\LCurvature\delta_i^2+2\delta_i\leq3\delta_i$.
Since $\bar U\geq0$, we have $\frac12\|\dot\theta_i\|_M^2\leq\bar{\mathcal H}_i^+$, so~\prettyref{lem:H-one-step} gives
\begin{ResizedEq}
&\bar{\mathcal H}_i^+-\bar{\mathcal H}_{\gamma(i)}^+\\
\leq&\frac{3\delta_i}{2\lambda_{\min}(M)}\|\dot\theta_i\|_M^2
+\frac{\delta_i}{2}\left((\|g\|+\LControlWorkForce)^2+\EpsCRDS^2\right)\\
\leq&\frac{3\delta_i}{\lambda_{\min}(M)}\bar{\mathcal H}_i^+
+\frac{\delta_i}{2}\left((\|g\|+\LControlWorkForce)^2+\EpsCRDS^2\right).
\end{ResizedEq}
Rearranging and using $e^{-2x}\leq1-x$ for $x=3\delta_i/\lambda_{\min}(M)\in(0,1/2)$, which is valid because $\delta_i<\lambda_{\min}(M)/6$ by~\prettyref{eq:H-bound-timestep}, we obtain the discrete Gr\"onwall recursion
\begin{ResizedEq}
&e^{-6\delta_i/\lambda_{\min}(M)}\bar{\mathcal H}_i^+\\
\leq&\left[1-\frac{3\delta_i}{\lambda_{\min}(M)}\right]\bar{\mathcal H}_i^+\\
\leq&\bar{\mathcal H}_{\gamma(i)}^++\frac{\delta_i}{2}\left((\|g\|+\LControlWorkForce)^2+\EpsCRDS^2\right).
\end{ResizedEq}
Telescoping from fixed node $0$ to $i$ along the ordered predecessor chain and using $t_i=\sum_{j\in\mathcal I:\,j\leq i}\delta_j$,
\begin{ResizedEq}
\bar{\mathcal H}_i^+
\leq&e^{6t_i/\lambda_{\min}(M)}\bar{\mathcal H}_0^+\\
&+\frac{1}{2}\left((\|g\|+\LControlWorkForce)^2+\EpsCRDS^2\right)\\
&\quad\times\sum_{k\in\mathcal I:\,k\leq i}\delta_ke^{6(t_i-t_{\gamma(k)})/\lambda_{\min}(M)}\\
\leq&e^{6t_i/\lambda_{\min}(M)}\Big[\bar{\mathcal H}_0^+
+\frac{t_i}{2}\left((\|g\|+\LControlWorkForce)^2+\EpsCRDS^2\right)\Big]\\
\triangleq&\LHamBound(t_i),
\end{ResizedEq}
which defines the explicit (exponential) closed form of $\LHamBound(t_i)$ and proves~\prettyref{eq:H-bound}.
It remains to verify clamp independence.
By~\prettyref{ass:Slater}, $\theta_0$ is strictly feasible, so $\tilde U(\theta_0,\EpsFJ)<\infty$; with the clamp choices $\LPosition\geq\|\theta_0\|$ and $\LHamilton\geq\tilde U(\theta_0,\EpsFJ)$,~\prettyref{lem:mollification}~\ref{it:moll-identity} gives $\bar U(\theta_0,\EpsFJ,\LHamilton,\LPosition)=\tilde U(\theta_0,\EpsFJ)$.
Hence $\bar{\mathcal H}_0^+=\frac12\|\dot\theta_0\|_M^2+\tilde U(\theta_0,\EpsFJ)$, with $\dot\theta_0=(\theta_0-\theta_{\gamma(0)})/\delta_0$, depends only on the fixed startup data $(\theta_{\gamma(0)},\theta_0,\delta_0)$ and the CRDS data $M$ and $\tilde U$, and in particular not on the clamps $\LHamilton,\LPosition$.
Since $\LCurvature$ was eliminated through~\prettyref{eq:H-bound-timestep}, $\LHamBound(t_i)$ is independent of $\LHamilton$, $\LPosition$, and $\LCurvature$; the only control-work constant it retains is the mollified control-force bound $\LControlWorkForce$, which is itself clamp-independent by~\prettyref{lem:controlwork-mollification}~\ref{it:cw-force}.
\end{proof}

\begin{proof}[Proof of~\prettyref{lem:trajectory-bound}]
Since~\prettyref{pb:PS-CRDS-M} is a smooth CRDS, the gradient-residual stationarity of $\theta_i$ is again~\prettyref{eq:one-step-stationarity}, now with $\|\rho_i\|\leq\LResidual$. Denoting its right-hand side by $s_i$, we have $M(\dot\theta_i-\dot\theta_{\gamma(i)})=\delta_i s_i$.
Take the standard inner product with $\dot\theta_i$. The left-hand side is
\begin{ResizedEq}
&\dot\theta_i^TM(\dot\theta_i-\dot\theta_{\gamma(i)})\\
=&\|\dot\theta_i\|_M^2-\dot\theta_i^TM\dot\theta_{\gamma(i)}\\
\geq&\|\dot\theta_i\|_M^2-\|\dot\theta_i\|_M\|\dot\theta_{\gamma(i)}\|_M,
\end{ResizedEq}
by the Cauchy-Schwarz inequality in the $M$-inner product.
For the right-hand side, the dissipation contributes $-\dot\theta_i^T\partial_{\dot\theta_i}\Damp_{\EpsFri,\LHamilton}^V(\theta_{\gamma(i)},\dot\theta_i)\leq0$, exactly as in~\prettyref{eq:H-one-step}: by the global convexity and zero-velocity minimization in~\prettyref{def:S-D}, $\dot\theta_i^T\partial_{\dot\theta_i}\Damp_{\EpsFri,\LHamilton}^V\geq \Damp_{\EpsFri,\LHamilton}^V(\theta_{\gamma(i)},\dot\theta_i)-\Damp_{\EpsFri,\LHamilton}^V(\theta_{\gamma(i)},0)\geq0$. Dropping this term,
\begin{ResizedEq}
\dot\theta_i^Ts_i
\leq&\dot\theta_i^T\left(-\FPP{\bar U}{\theta_i}+g+\bar c(\theta_i,u_i)+\rho_i\right)\\
\leq&\|\dot\theta_i\|\left(\LPotentialForce+\|g\|+\LControlWorkForce+\LResidual\right)
=\|\dot\theta_i\|\,C,
\end{ResizedEq}
where $C\triangleq\LPotentialForce+\|g\|+\LControlWorkForce+\LResidual$ is the uniform total-force bound and we used~\prettyref{lem:mollification}~\ref{it:moll-bound} for $\|\FPPR{\bar U}{\theta_i}\|\leq\LPotentialForce$,~\prettyref{lem:controlwork-mollification}~\ref{it:cw-force} for $\|\bar c(\theta_i,u_i)\|\leq\LControlWorkForce$, and $\|\rho_i\|\leq\LResidual$.
With the norm comparison $\|\dot\theta_i\|\leq\|\dot\theta_i\|_M/\sqrt{\lambda_{\min}(M)}$, the identity $M(\dot\theta_i-\dot\theta_{\gamma(i)})=\delta_i s_i$ therefore yields
\begin{ResizedEq}
\|\dot\theta_i\|_M^2-\|\dot\theta_i\|_M\|\dot\theta_{\gamma(i)}\|_M
\leq\delta_i\frac{C}{\sqrt{\lambda_{\min}(M)}}\|\dot\theta_i\|_M.
\end{ResizedEq}
If $\|\dot\theta_i\|_M=0$ the one-step velocity increment below is trivial; otherwise divide by $\|\dot\theta_i\|_M>0$ to obtain
\begin{ResizedEq}
\label{eq:velocity-increment}
\|\dot\theta_i\|_M\leq\|\dot\theta_{\gamma(i)}\|_M+\delta_i\frac{C}{\sqrt{\lambda_{\min}(M)}}.
\end{ResizedEq}
Telescoping~\prettyref{eq:velocity-increment} along the ordered $\gamma$-chain from fixed node $0$ to $i$ and using $\sum_{k\in\mathcal I:\,k\leq i}\delta_k=t_i$ from~\prettyref{eq:elapsed-time} gives the velocity bound: we may take
\begin{ResizedEq}
\|\dot\theta_i\|_M\leq\LVelBound(t_i),\qquad\LVelBound(t_i)=\|\dot\theta_0\|_M+\frac{C}{\sqrt{\lambda_{\min}(M)}}\,t_i,
\end{ResizedEq}
which is the velocity bound function $\LVelBound$.
For the position, telescoping $\theta_i=\theta_0+\sum_{k\in\mathcal I:\,k\leq i}\delta_k\dot\theta_k$ along the ordered predecessor chain and applying the norm comparison to the velocity bound,
\begin{ResizedEq}
\|\theta_i\|
\leq&\|\theta_0\|+\sum_{k\in\mathcal I:\,k\leq i}\delta_k\|\dot\theta_k\|\\
\leq&\|\theta_0\|+\sum_{k\in\mathcal I:\,k\leq i}\delta_k\left(\frac{\|\dot\theta_0\|_M}{\sqrt{\lambda_{\min}(M)}}+\frac{C}{\lambda_{\min}(M)}\,t_k\right).
\end{ResizedEq}
Using $\sum_{k\in\mathcal I:\,k\leq i}\delta_k=t_i$ and $\sum_{k\in\mathcal I:\,k\leq i}\delta_kt_k\leq t_i\sum_{k\in\mathcal I:\,k\leq i}\delta_k=t_i^2$ (since $t_k\leq t_i$) gives the position bound: we may take
\begin{ResizedEq}
\|\theta_i\|\leq&\LPosBound(t_i)\\
\LPosBound(t_i)=&\|\theta_0\|+\frac{t_i}{\sqrt{\lambda_{\min}(M)}}\|\dot\theta_0\|_M+\frac{C}{\lambda_{\min}(M)}\,t_i^2,
\end{ResizedEq}
the position bound function $\LPosBound$. Together these establish~\prettyref{eq:trajectory-velocity-bound}.
No upper-timestep smallness condition was used. The finite-difference setup requires each individual timestep to be positive, while the stronger uniform numerical floor $\delta_i\geq\delta_{\min}>0$ is retained for the downstream fixed-discretization and SQP applications. The initial value entered only through the fixed startup data $\dot\theta_0,\theta_0$ (with $\dot\theta_0=(\theta_0-\theta_{\gamma(0)})/\delta_0$), so~\prettyref{ass:Slater} is not needed; the bound holds for any active timestep sequence satisfying that floor.
\end{proof}

\section{Proofs for Global Convergence}
\label{appen:global-convergence-proofs}

\begin{proof}[Proof of~\prettyref{lem:jacobian-lower-bound}]
Differentiating~\prettyref{pb:PS-CRDS-M} with $\partial\ddot{\theta}_i/\partial\theta_i=I/\delta_i^2$, the per-step Hessian $\StatJac_{ii}=\FPPT{\PerStepObj_i}{\theta_i}$ expands into
\begin{ResizedEq}
\label{eq:sqp-diagonal-expand}
\StatJac_{ii}=&\frac{M}{\delta_i^2}+\FPPT{\bar U(\theta_i,\EpsFJ,\LHamilton,\LPosition)}{\theta_i}+\\
&\FPPT{\bar W(\theta_i,u_i,\LPosition)}{\theta_i}+\delta_i\FPPT{\Damp_{\EpsFri,\LHamilton}^V(\theta_{\gamma(i)},\dot{\theta}_i)}{\theta_i}.
\end{ResizedEq}
We treat the three non-kinetic terms.
First, the mollified potential is $\LCurvature$-weakly convex, $\FPPT{\bar U}{\theta_i}\succeq-\LCurvature I$, by~\prettyref{lem:mollification}~\ref{it:moll-curv}. Second, $\Damp_{\EpsFri,\LHamilton}^V(\theta_{\gamma(i)},\cdot)$ is convex and, through the affine map $\dot{\theta}_i=(\theta_i-\theta_{\gamma(i)})/\delta_i$, so is $\theta_i\mapsto \Damp_{\EpsFri,\LHamilton}^V(\theta_{\gamma(i)},\dot{\theta}_i)$; hence $\delta_i\FPPT{\Damp_{\EpsFri,\LHamilton}^V}{\theta_i}\succeq0$ by~\prettyref{def:CRDS}~\ref{it:crds-damping} and~\prettyref{def:S-D}, so this term contributes no negative curvature. Third, the mollified control work $\bar W$ is $C^2$ in $\theta_i$ by~\prettyref{lem:controlwork-mollification}~\ref{it:cw-c2}, and by~\prettyref{lem:controlwork-mollification}~\ref{it:cw-curv} its configuration Hessian is bounded globally and independently of the region: there is a region-independent (clamp-dependent) constant $\LControlWorkCurvature<\infty$ with
\begin{ResizedEq}
\label{eq:sqp-cw}
-\LControlWorkCurvature I\preceq\FPPT{\bar W(\theta_i,u_i,\LPosition)}{\theta_i}\preceq\LControlWorkCurvature I,
\end{ResizedEq}
valid for every $\theta_i$ and every $u_i\in\mathcal{D}$, clamp-dependent but region-independent by~\prettyref{lem:controlwork-mollification}~\ref{it:cw-curv}; in particular $\FPPT{\bar W}{\theta_i}\succeq-\LControlWorkCurvature I$ everywhere, so that $\LControlWorkCurvature$---and thus the timestep threshold $\delta^\star$ used later---does not depend on the region radius. Consequently the combined non-kinetic Hessian is bounded below by $-(\LCurvature+\LControlWorkCurvature)I$, giving the unconditional curvature identity
\begin{ResizedEq}
\label{eq:sqp-diagonal-lower}
\StatJac_{ii}\succeq\frac{M}{\delta_i^2}-(\LCurvature+\LControlWorkCurvature) I\succeq\left[\frac{\lambda_{\min}(M)}{\delta_i^2}-(\LCurvature+\LControlWorkCurvature)\right]I,
\end{ResizedEq}
using $M\succeq\lambda_{\min}(M)I$ from~\prettyref{def:full-rank}, whose combined curvature modulus sums the potential curvature $\LCurvature$ of~\prettyref{lem:mollification}~\ref{it:moll-curv} and the region-independent control-work Hessian bound $\LControlWorkCurvature$ of~\prettyref{eq:sqp-cw}; it holds for every $\delta_i$ with no timestep restriction.

It remains to bound $\sigma_{\min}(\StatJac)$ away from zero uniformly. Order the optimized nodes chronologically, so that, as noted in~\prettyref{sec:global-licq}, $\StatJac$ is block lower-triangular. Split it as $\StatJac=J_\text{diag}+J_\text{off}$, where the block diagonal $J_\text{diag}\triangleq\diag(\StatJac_{ii})_{i\in\mathcal{I}}$ collects the per-step Hessians and the strictly block-lower-triangular $J_\text{off}$ collects only those earlier-node couplings $\StatJac_{i\gamma(i)},\StatJac_{i\gamma^2(i)}$ whose column indices belong to $\mathcal I$. In particular, for $i_1\triangleq\min\mathcal I$, both stencil predecessors are fixed ($0$ and $\gamma(0)$), so the $i_1$-row has only the diagonal optimized block $\StatJac_{i_1i_1}$. Under the singularity-check hypothesis $\sigma_{\min}(\StatJac_{ii})\geq\LLambdaMin$ for every $i\in\mathcal{I}$, each diagonal block is invertible with $\|\StatJac_{ii}^{-1}\|=1/\sigma_{\min}(\StatJac_{ii})\leq1/\LLambdaMin$. Hence $J_\text{diag}$ is invertible and $\|J_\text{diag}^{-1}\|=\max_i\|\StatJac_{ii}^{-1}\|\leq1/\LLambdaMin$. Since $\mathcal{I}$ is finite and every off-diagonal block of $J_\text{off}$ (the couplings $\StatJac_{i\gamma(i)},\StatJac_{i\gamma^2(i)}$ displayed in~\prettyref{eq:sqp-jacobian-structure}) is a continuous function of the $C^2$ Hessian blocks on the compact region of~\prettyref{lem:trajectory-bound}, each such block has finite norm; taking the maximum over the finitely many blocks indexed by $\mathcal{I}$, we have $\|J_\text{off}\|\leq C_2<\infty$ for a finite proof-local constant $C_2$ (its explicit values are recorded later).

Set $N\triangleq J_\text{diag}^{-1}J_\text{off}$. Since $J_\text{diag}^{-1}$ is block diagonal and $J_\text{off}$ strictly block lower-triangular, $N$ is strictly block lower-triangular, hence nilpotent with $N^{|\mathcal{I}|}=0$. Therefore $\StatJac=J_\text{diag}(I+N)$ is invertible, with the terminating Neumann series $(I+N)^{-1}=\sum_{k=0}^{|\mathcal{I}|-1}(-N)^k$; in particular $\StatJac(\StackConf)$ is nonsingular. Using $\|N\|\leq\|J_\text{diag}^{-1}\|\,\|J_\text{off}\|\leq C_2/\LLambdaMin$,
\begin{ResizedEq}
\|(I+N)^{-1}\|\leq\sum_{k=0}^{|\mathcal{I}|-1}\left(\frac{C_2}{\LLambdaMin}\right)^k\triangleq C_3,
\end{ResizedEq}
so $\|\StatJac^{-1}\|=\|(I+N)^{-1}J_\text{diag}^{-1}\|\leq C_3/\LLambdaMin$, and hence
\begin{ResizedEq}
\sigma_{\min}(\StatJac)=\frac{1}{\|\StatJac^{-1}\|}\geq\frac{\LLambdaMin}{C_3}\triangleq\SigmaMinLICQ(\StatJac)>0
\end{ResizedEq}
Finally, the check is discharged by a sufficiently small timestep, which is the sufficiency claim of the statement: if $\delta_i<\delta^\star$, the threshold of~\prettyref{eq:delta-star}, then $\lambda_{\min}(M)/\delta_i^2-(\LCurvature+\LControlWorkCurvature)>\LLambdaMin$, so~\prettyref{eq:sqp-diagonal-lower} gives $\StatJac_{ii}\succeq\LLambdaMin I$. Thus the block is positive definite and $\sigma_{\min}(\StatJac_{ii})\geq\LLambdaMin$, discharging the hypothesis; this is the threshold enforced by the adaptive subdivision of~\prettyref{sec:global-sqp}.
\end{proof}

\begin{proof}[Proof of~\prettyref{lem:jacobian-upper-bound}]
We use the pessimistic block-sum estimate, in which each of the $|\mathcal{I}|$ optimized block rows of $\StatJac$ has at most three nonzero blocks (with fixed-predecessor columns omitted, so the first optimized row has only its diagonal block),
\begin{ResizedEq}
\label{eq:sqp-block-sum}
&\sigma_{\max}(\StatJac)=\|\StatJac\|\\
\leq&\sum_{i\in\mathcal{I}}\sum_{j\in\{\gamma^2(i),\gamma(i),i\}\cap\mathcal I}\|\StatJac_{ij}\|.
\end{ResizedEq}
For the diagonal block, by~\prettyref{lem:jacobian-lower-bound}, specifically the diagonal expansion~\prettyref{eq:sqp-diagonal-expand}, we have $\|M/\delta_i^2\|=\lambda_{\max}(M)/\delta_i^2$; the mollified-potential Hessian obeys $\|\FPPT{\bar U}{\theta_i}\|\leq\LCurvature$ by the two-sided bound of~\prettyref{lem:mollification}~\ref{it:moll-curv}; the mollified control-work Hessian obeys $\|\FPPT{\bar W}{\theta_i}\|\leq\LControlWorkCurvature$ by the two-sided bound of~\prettyref{lem:controlwork-mollification}~\ref{it:cw-curv} (equivalently~\prettyref{eq:sqp-cw}); and the convex-dissipation Hessians are bounded on the norm-bounded region $\mathcal{R}$ of~\prettyref{lem:trajectory-bound}, which we make explicit through the proof-local constants
\begin{ResizedEq}
C_4\triangleq&\sup_{i\in\mathcal{I},\,\mathcal{R}}\left\|\delta_i\FPPT{\Damp_{\EpsFri,\LHamilton}^V}{\theta_i}\right\|,\\
C_5\triangleq&\sup_{i\in\mathcal{I},\,\mathcal{R}}\left\|\delta_i\FPPTT{\Damp_{\EpsFri,\LHamilton}^V}{\theta_i}{\theta_{\gamma(i)}}\right\|,
\end{ResizedEq}
all finite for a fixed discretization scheme by continuity on the compact region, where $C_4$ and $C_5$ absorb the $\delta_i$ factors of the dissipation Hessians and are therefore allowed to depend on $\delta_{\min}$, while $\LControlWorkCurvature$ is the region-independent control-work curvature bound of~\prettyref{eq:sqp-cw} (equivalently~\prettyref{lem:controlwork-mollification}~\ref{it:cw-curv}). Collecting the four diagonal terms of~\prettyref{eq:sqp-diagonal-expand},
\begin{ResizedEq}
\|\StatJac_{ii}\|\leq\frac{\lambda_{\max}(M)}{\delta_i^2}+\LCurvature+\LControlWorkCurvature+C_4<\infty.
\end{ResizedEq}
For the predecessor derivatives, only the kinetic term $\tfrac{\delta_i^2}{2}\|\ddot{\theta}_i\|_M^2$ and the dissipation couple $\theta_i$ to earlier timesteps; differentiating the discretization~\prettyref{eq:discretization} gives
\begin{ResizedEq}
\StatJac_{i\gamma(i)}=&-\frac{(1+\delta_i/\delta_{\gamma(i)})M}{\delta_i^2}+\delta_i\FPPTT{\Damp_{\EpsFri,\LHamilton}^V}{\theta_i}{\theta_{\gamma(i)}},\\
\StatJac_{i\gamma^2(i)}=&\frac{M}{\delta_i\delta_{\gamma(i)}},
\end{ResizedEq}
These derivatives are optimized off-diagonal blocks only when their column labels belong to $\mathcal I$; for $i_1=\min\mathcal I$, the fixed stencil $\theta_{\gamma(0)},\theta_0,\theta_{i_1}$ contributes only $\StatJac_{i_1i_1}$ to the optimized Jacobian. With $C_5$ bounding the cross dissipation,
\begin{ResizedEq}
\|\StatJac_{i\gamma(i)}\|\leq&\left(1+\frac{\delta_i}{\delta_{\gamma(i)}}\right)\frac{\lambda_{\max}(M)}{\delta_i^2}+C_5\\
\|\StatJac_{i\gamma^2(i)}\|\leq&\frac{\lambda_{\max}(M)}{\delta_i\delta_{\gamma(i)}}.
\end{ResizedEq}
On a fixed discretization scheme every $\delta_i\geq\delta_{\min}>0$ and each of the block norms displayed above is finite (with a bound that may depend on $\delta_{\min}$). Since $\mathcal{I}$ collects finitely many timesteps, the sum~\prettyref{eq:sqp-block-sum} has at most $3|\mathcal{I}|$ finite terms, so
\begin{ResizedEq}
&\sigma_{\max}(\StatJac)\leq|\mathcal{I}|\times\\
&\max_{i\in\mathcal{I}}\sum_{j\in\{\gamma^2(i),\gamma(i),i\}\cap\mathcal I}\|\StatJac_{ij}\|<\infty,
\end{ResizedEq}
finite for any fixed discretization scheme with $\delta_{\min}>0$, with a bound that may depend on $\delta_{\min}$.
For the control block, $\StatJacU$ is block diagonal with diagonal blocks
\begin{ResizedEq}
\frac{\partial\StatRes_i}{\partial u_i}=\FPPTT{\bar W(\theta_i,u_i,\LPosition)}{\theta_i}{u_i}=-\FPP{\bar c(\theta_i,u_i)}{u_i},
\end{ResizedEq}
using $\bar c=-\FPPR{\bar W}{\theta_i}$, the mollified control force. By~\prettyref{lem:controlwork-mollification}~\ref{it:cw-mixed} the mixed second derivative $\FPPTT{\bar W(\theta_i,u_i,\LPosition)}{\theta_i}{u_i}$ is finite, and by continuity it is bounded on the compact region $\mathcal{R}$ of~\prettyref{lem:trajectory-bound} times the compact control set $\mathcal{D}$; denote this region-dependent proof-local bound by $C_u\triangleq\sup_{\mathcal{R},\,u_i\in\mathcal{D}}\|\FPPTT{\bar W}{\theta_i}{u_i}\|<\infty$. Each such block satisfies $\|\partial\StatRes_i/\partial u_i\|\leq C_u$, so the block-diagonal operator norm attains the same bound, $\sigma_{\max}(\StatJacU)\leq C_u<\infty$, with no dependence on $\delta_i$.
\end{proof}

\subsection{Adaptive-Timestep SQP}
We now collect the proof-local machinery behind~\prettyref{thm:sqp-convergence}. The adaptive-timestep SQP method it analyzes---the $\ell_1$ merit~\prettyref{eq:sqp-merit}, the quadratic program~\prettyref{pb:sqp-qp} with conditioning bound~\prettyref{eq:sqp-hessian}, the safeguarded line searches~\prettyref{eq:sqp-linesearch} and~\ref{eq:sqp-constraint-linesearch}, and the four algorithms~\prettyref{alg:subd}, \ref{alg:constraint-solve}, \ref{alg:sqp-inner}, and \ref{alg:sqp}---is defined in~\prettyref{sec:global-sqp}, whose notation we reuse throughout ($s=(s_\StackConf,s_\StackCtrl)$, $B=(B_\StackConf,B_\StackCtrl)$, $\varrho$, $\mu_\StackConf,\mu_\StackCtrl$, $\tau$, $\alpha$, $\phi_\varrho$, $\LHessianCond$, $\LResidualWork$, $\bar\varrho$, $\EpsStop$; the timestep-subscripted $\mu_{\theta_i}$ and $s_{u_i}$ are the timestep-$i$ sub-blocks of the stacked $\mu_\StackConf$ and $s_\StackCtrl$). The proofs below additionally introduce a few genuinely proof-local constants, defined where they first appear and excluded from~\prettyref{table:notation}: the violation-dependent bounds $c_O(\cdot),c_s(\cdot),c_\mu(\cdot)$ and the uniform-continuity moduli $\omega_\StatJac,\omega_r,\omega_\sigma,\omega_\Lambda$. The development follows the standard line-search SQP theory~\cite{bertsekas1997nonlinear,solodov2009global}, specialized to the residual $\StatRes$ and Jacobians $\StatJac,\StatJacU$ of~\prettyref{eq:sqp-definition}.

Throughout, every iterate is confined to the compact, norm-bounded region of~\prettyref{lem:trajectory-bound}, and the optimized node set $\mathcal I$ of~\prettyref{sec:global-sqp} is refined by bisection through~\prettyref{alg:subd}. The evolving algorithmic state is the explicit tuple $(\mathcal I,\StackConf,\StackCtrl,\varrho)$; routines return the updated components specified in~\prettyref{sec:global-sqp}, while the predecessor map is derived from the total order on $\{0\}\cup\mathcal I$ rather than carried as mutable state. On this fixed region, the objective regularity of~\prettyref{ass:objective} together with the twice-differentiable constraint and dynamics data of~\prettyref{def:CRDS}~\ref{it:crds-c2} yield four facts used repeatedly.
\begin{enumerate}[label=(\alph*),ref=(\alph*),leftmargin=*]
\item\label{it:fact-Obounded} By~\prettyref{ass:objective} the user objective $O$ is continuous (indeed $C^1$), so on the compact region it is bounded below; shifting $O$ by the harmless constant $-\inf O$ (which alters neither the minimizer nor the iterates) we may assume $O\geq0$.
\item\label{it:fact-gradO-lip} By~\prettyref{ass:objective} $\nabla O$ is Lipschitz with modulus $\LObjGradLip<\infty$; and since $O$ is $C^1$, $\nabla O$ is continuous, hence bounded on the compact region of~\prettyref{lem:trajectory-bound}, $\|\nabla O\|\leq\LObjGradBound<\infty$.
\item\label{it:fact-J-uc} Because $\PerStepObj$ is only $C^2$, the stationarity Jacobian $\StatJac=\nabla\StatRes$ is continuous but need not be Lipschitz; on the compact region it is uniformly continuous, with modulus $\omega_\StatJac$, i.e. $\|\StatJac(\StackConf)-\StatJac(\StackConf')\|\leq\omega_\StatJac(\|\StackConf-\StackConf'\|)$ with $\omega_\StatJac(t)\to0$ as $t\to0^+$. The same argument applies to the full residual Jacobian $\nabla_{(\StackConf,\StackCtrl)}\StatRes=(\StatJac,\StatJacU)$, which---being continuous on the compact region since $\PerStepObj$ is $C^2$---is uniformly continuous with modulus $\omega_r$; this is the modulus that governs any expansion of $\StatRes$ stepping in both the configuration and control blocks.
\item\label{it:fact-sigmin-uc} Likewise $\sigma_{\min}(\StatJac(\cdot))$ is uniformly continuous with modulus $\omega_\sigma$.
The only consequence of~\prettyref{it:fact-J-uc} and~\ref{it:fact-sigmin-uc} is that residual-based step-acceptance thresholds below are stated as existence results ($\bar\tau>0$) rather than closed forms; the $O$-based thresholds remain explicit through $\LObjGradLip$. No a-priori iterate bound is needed, since~\prettyref{lem:trajectory-bound} already fixes the region.
\item\label{it:fact-singcheck} Singularity-Check Conditioning: This is exactly~\prettyref{lem:jacobian-lower-bound}. Whenever the per-block check $\sigma_{\min}(\StatJac_{ii}^k)\geq\LLambdaMin$ passes for every $i$, that lemma yields $\sigma_{\min}(\StatJac^k)\geq\SigmaMinLICQ(\StatJac)=\LLambdaMin/C_3>0$. The derived modulus depends on more than the threshold $\LLambdaMin$: through the terminating Neumann series $C_3=\sum_{k=0}^{|\mathcal{I}|-1}(C_2/\LLambdaMin)^k$ it also involves the region-dependent off-diagonal coupling bound $C_2$ (set by the dissipation Hessian couplings on the compact region of~\prettyref{lem:trajectory-bound}) and the mesh size $|\mathcal{I}|$ (which increases under subdivision). All three are fixed once the mesh is fixed, so $\SigmaMinLICQ(\StatJac)$ is uniform over the run on that fixed mesh. This is the form of the lower bound invoked below.
\end{enumerate}

We first establish that the constraint-satisfaction phase is well posed and finite. Both results hold whenever the singularity check of~\prettyref{ln:sqp-constraint-singular} passes, so that $\sigma_{\min}(\StatJac)\geq\SigmaMinLICQ(\StatJac)>0$ by~\prettyref{lem:jacobian-lower-bound}.
\begin{lemma}[Constraint-phase line-search well-posedness]
\label{lem:sqp-constraint-linesearch}
Consider one fixed invocation of~\prettyref{alg:constraint-solve}. Suppose the singularity check $\sigma_{\min}(\StatJac_{ii}^k)\geq\LLambdaMin$ of~\prettyref{ln:sqp-constraint-singular} holds for every $i$, and the iterates lie in the region of~\prettyref{lem:trajectory-bound}. At every iteration at which the algorithm has not returned at~\prettyref{ln:sqp-constraint-return}, i.e. $\|\StatRes^k\|_1\geq\LResidualWork$, there exists $\bar\tau>0$, uniform over all such iterations of this invocation, such that every $[\tau^\ell]^k<\bar\tau$ satisfies the constraint line search~\prettyref{eq:sqp-constraint-linesearch}.
\end{lemma}
\begin{proof}[Proof of~\prettyref{lem:sqp-constraint-linesearch}]
Because the singularity check passes,~\prettyref{lem:jacobian-lower-bound} gives $\sigma_{\min}(\StatJac^k)\geq\SigmaMinLICQ(\StatJac)>0$, so the Newton step $s_\StackConf^k=-[\StatJac^k]^{-1}\StatRes^k$ is well defined with
\begin{ResizedEq}
\label{eq:sqp-constraint-step}
\|s_\StackConf^k\|\leq&\|[\StatJac^k]^{-1}\|\,\|\StatRes^k\|
\leq\frac{\|\StatRes^k\|}{\SigmaMinLICQ(\StatJac)}\leq\frac{\|\StatRes^0\|}{\SigmaMinLICQ(\StatJac)},
\end{ResizedEq}
the last step by induction: each previously accepted instance of~\prettyref{eq:sqp-constraint-linesearch} makes $\Lambda^k=\tfrac12\|\StatRes^k\|^2$ strictly decrease, so $\|\StatRes^k\|\leq\|\StatRes^0\|$. Since $\|s_\StackConf^k\|$ is bounded by~\prettyref{eq:sqp-constraint-step}, every trial point $\StackConf^k+[\tau^\ell]^ks_\StackConf^k$ with $[\tau^\ell]^k$ small enough stays inside a fixed compact neighborhood $\mathcal{R}'$ of the region of~\prettyref{lem:trajectory-bound} (that region enlarged by a fixed radius), on which the uniform-continuity moduli invoked below are valid; the descent property itself will follow from the Taylor expansion of $\Lambda$, not from invertibility of $\StatJac$ along the accepted segment.
Now write $\Lambda(\StackConf)\triangleq\tfrac12\|\StatRes(\StackConf)\|^2$, whose gradient $\nabla\Lambda=[\StatJac]^T\StatRes$ is, on the compact set $\mathcal{R}'$, a product of the uniformly continuous $\StatJac$ (fact~\ref{it:fact-J-uc}) and the continuously differentiable $\StatRes$, hence itself uniformly continuous; denote its modulus by $\omega_\Lambda$. For any step $p=[\tau^\ell]^ks_\StackConf^k$, we have
\begin{ResizedEq}
\Lambda(\StackConf^k+p)\leq\Lambda^k+\nabla\Lambda(\StackConf^k)^Tp+\|p\|\,\omega_\Lambda(\|p\|).
\end{ResizedEq}
Since $\nabla\Lambda(\StackConf^k)^Ts_\StackConf^k=[\StatRes^k]^T\StatJac^k(-[\StatJac^k]^{-1}\StatRes^k)=-\|\StatRes^k\|^2=-2\Lambda^k$, we have
\begin{ResizedEq}
&\Lambda(\StackConf^k+[\tau^\ell]^ks_\StackConf^k)\\
\leq&(1-2[\tau^\ell]^k)\Lambda^k+[\tau^\ell]^k\|s_\StackConf^k\|\,\omega_\Lambda([\tau^\ell]^k\|s_\StackConf^k\|),
\end{ResizedEq}
so the Armijo condition~\prettyref{eq:sqp-constraint-linesearch}, i.e. $\Lambda(\StackConf^k+[\tau^\ell]^ks_\StackConf^k)<(1-2\alpha[\tau^\ell]^k)\Lambda^k$, holds as soon as
\begin{ResizedEq}
\|s_\StackConf^k\|\,\omega_\Lambda([\tau^\ell]^k\|s_\StackConf^k\|)<2(1-\alpha)\Lambda^k.
\end{ResizedEq}
By the hypothesis $\|\StatRes^k\|_1\geq\LResidualWork$ of the statement, Cauchy--Schwarz gives $\Lambda^k=\tfrac12\|\StatRes^k\|^2\geq\tfrac12\|\StatRes^k\|_1^2/(\dim\StatRes)\geq\tfrac12\LResidualWork^2/(\dim\StatRes)>0$, while by~\prettyref{eq:sqp-constraint-step} $\|s_\StackConf^k\|$ is bounded and $\omega_\Lambda([\tau^\ell]^k\|s_\StackConf^k\|)\to0$ as $[\tau^\ell]^k\to0^+$. Hence the displayed strict inequality holds for all sufficiently small $[\tau^\ell]^k$, which also keep the trial point inside $\mathcal{R}'$. Since the mesh, controls, initial residual, and Jacobian modulus are fixed throughout this invocation, every bound is uniform over $\mathcal{R}'$ and all its non-returning iterations, yielding a common $\bar\tau>0$.
\end{proof}

\begin{lemma}[Finite termination of the constraint-satisfaction phase]
\label{lem:sqp-constraint-termination}
Under the standing regularity facts,~\prettyref{alg:constraint-solve} is well defined and terminates after finitely many iterations, returning either a node index for subdivision (\prettyref{ln:sqp-constraint-singular}) or a point with $\|\StatRes\|_1<\LResidualWork$ (\prettyref{ln:sqp-constraint-return}).
\end{lemma}
\begin{proof}[Proof of~\prettyref{lem:sqp-constraint-termination}]
Consider one invocation of~\prettyref{alg:constraint-solve}. Its finite mesh and controls remain fixed until the routine returns: subdivision occurs only after the caller receives a node index. Thus $\delta_{\min}>0$ is fixed and the input controls remain in $\mathcal D$.

Well-posedness: If the singularity check~\prettyref{ln:sqp-constraint-singular} fails, the algorithm returns immediately. Otherwise $\StatJac^k$ is invertible, so the Newton step exists. Since the residual return test~\prettyref{ln:sqp-constraint-return} precedes the step, whenever the iteration reaches the line search we have $\|\StatRes^k\|_1\geq\LResidualWork$. The accepted Armijo inequality makes $\Lambda^k\triangleq\tfrac12\|\StatRes^k\|_2^2$ non-increasing, and therefore
\begin{ResizedEq}
\|\StatRes^k\|_1\leq\sqrt{\dim\StatRes}\,\|\StatRes^k\|_2
\leq\sqrt{\dim\StatRes}\,\|\StatRes^0\|_2\triangleq R_{\rm cs}<\infty.
\end{ResizedEq}
With this common residual ceiling, the fixed positive timestep floor, and fixed admissible controls,~\prettyref{lem:trajectory-bound} places every accepted iterate in one common compact region. Hence~\prettyref{lem:sqp-constraint-linesearch} applies uniformly and supplies an acceptable trial at every non-returning iteration. Thus each iteration is well defined.

Finite termination: Suppose for contradiction that the routine runs indefinitely, so neither the singularity exit nor the residual exit is taken. Then the singularity check passes and $\|\StatRes^k\|_1\geq\LResidualWork$ for every $k$. Consequently
\begin{ResizedEq}
\Lambda^k\geq\frac{\|\StatRes^k\|_1^2}{2\dim\StatRes}
\geq\frac{\LResidualWork^2}{2\dim\StatRes}>0.
\end{ResizedEq}
Write $\beta_k\in\{1,\tau,\tau^2,\ldots\}$ for the first accepted geometric-backtracking trial. Telescoping the accepted decreases gives
\begin{ResizedEq}
2\alpha\sum_{m=0}^k\beta_m\Lambda^m
\leq\Lambda^0-\Lambda^{k+1}\leq\Lambda^0,
\end{ResizedEq}
so $\sum_m\beta_m\leq(\dim\StatRes)\Lambda^0/(\alpha\LResidualWork^2)<\infty$.

On the common compact region,~\prettyref{lem:sqp-constraint-linesearch} gives one $\bar\tau>0$ such that every trial $0<\beta<\bar\tau$ passes. Replace $\bar\tau$ by $\min\{1,\bar\tau\}$. If the initial trial is accepted, then $\beta_k=1\geq\bar\tau$. Otherwise its rejected predecessor $\beta_k/\tau$ cannot be smaller than $\bar\tau$, since every such trial would pass; hence $\beta_k\geq\tau\bar\tau$. Therefore every accepted step satisfies $\beta_k\geq\tau\bar\tau>0$, forcing $\sum_m\beta_m=\infty$, a contradiction. The routine must return after finitely many iterations, through one of its two exits.
\end{proof}

The following three results are the standard SQP support lemmas, restated for our residual system; each holds whenever the singularity check of~\prettyref{ln:sqp-inner-singular} passes, so that LICQ with modulus $\SigmaMinLICQ(\StatJac)$ is in force by~\prettyref{lem:jacobian-lower-bound}.
\begin{lemma}[QP subproblem solvability and merit descent]
\label{lem:sqp-qp}
Suppose $u_i^k\in\mathcal{D}$ for all $i$, the singularity check passes so that $\sigma_{\min}(\StatJac^k)\geq\SigmaMinLICQ(\StatJac)>0$ (\prettyref{lem:jacobian-lower-bound}), the QP Hessian obeys~\prettyref{eq:sqp-hessian}, and $\varrho^k\geq\max_i\|\mu_{\theta_i}^k\|+\bar\varrho$. Then $\text{QP}^k$ is convex and feasible with a unique solution, and its solution $s^k$ is a descent direction of the merit~\prettyref{eq:sqp-merit} with directional derivative bounded by $[D\phi_\varrho]^k\leq[\bar D\phi_\varrho]^k=-\tfrac{1}{\LHessianCond}\|s^k\|^2-\bar\varrho\|\StatRes^k\|_1\leq0$, with strict inequality unless $s^k=0$ and $\StatRes^k=0$.
\end{lemma}
\begin{proof}[Proof of~\prettyref{lem:sqp-qp}]
By~\prettyref{eq:sqp-hessian} the objective of~\prettyref{pb:sqp-qp} is strongly convex, so over the nonempty closed convex feasible set (nonemptiness shown next) a minimizer exists and is unique. Feasibility follows because $\mathcal{D}$ is nonempty---pick any $u_i^k+s_{u_i}^k\in\mathcal{D}$---and $\StatJac^k$ is nonsingular under LICQ (\prettyref{lem:jacobian-lower-bound}), so the equality constraint is solvable for $s_\StackConf^k$. Because $\text{QP}^k$ is convex with a nonsingular equality-constraint Jacobian, its unique solution $s^k$ is characterized by the associated KKT system
\begin{ResizedEq}
\label{eq:sqp-kktqp}
\text{KKT-QP}^k:
\begin{cases}
0=\FPP{O}{\StackConf}^k+B_\StackConf^ks_\StackConf^k+[\StatJac^k]^T\mu_\StackConf^k\\
0=\FPP{O}{\StackCtrl}^k+B_\StackCtrl^ks_\StackCtrl^k+[\StatJacU^k]^T\mu_\StackConf^k+\mu_\StackCtrl^k\\
0=\StatRes^k+\StatJac^ks_\StackConf^k+\StatJacU^ks_\StackCtrl^k\\
u_i^k+s_{u_i}^k\in\mathcal{D},\quad\mu_{u_i}^k\in N_{\mathcal{D}}(u_i^k+s_{u_i}^k),
\end{cases}
\end{ResizedEq}
where $\mu_\StackConf^k$ (resp. $\mu_\StackCtrl^k$) is the multiplier of the physics (resp. control) constraint and $\mu_{\theta_i}^k$, $\mu_{u_i}^k$ are the timestep-$i$ sub-blocks of the stacked $\mu_\StackConf^k$, $\mu_\StackCtrl^k$ introduced in~\prettyref{sec:global-sqp}. We compute the directional derivative of the merit~\prettyref{eq:sqp-merit} along $s^k$. Splitting the residual entries into $I_{\neq0}\triangleq\{(i,j):[\StatRes_i^k]_j\neq0\}$ and its complement $I_{=0}$, the directional derivative of $\varrho\|\StatRes\|_1$ is
\begin{ResizedEq}
&\varrho^k\!\!\sum_{(i,j)\in I_{\neq0}}\!\!\frac{[\StatRes_i^k]_j\,[\StatJac^ks_\StackConf^k+\StatJacU^ks_\StackCtrl^k]_{ij}}{|[\StatRes_i^k]_j|}\\
&+\varrho^k\!\!\sum_{(i,j)\in I_{=0}}\!\!\big|[\StatJac^ks_\StackConf^k+\StatJacU^ks_\StackCtrl^k]_{ij}\big|.
\end{ResizedEq}
The third line of~\prettyref{eq:sqp-kktqp} gives $\StatJac^ks_\StackConf^k+\StatJacU^ks_\StackCtrl^k=-\StatRes^k$ entrywise, so each bracket equals $-[\StatRes_i^k]_j$; the $I_{=0}$ terms then vanish and the $I_{\neq0}$ terms collapse to $-\varrho^k\|\StatRes^k\|_1$. Hence
\begin{ResizedEq}
[D\phi_\varrho]^k=[s_\StackConf^k]^T\FPP{O}{\StackConf}^k+[s_\StackCtrl^k]^T\FPP{O}{\StackCtrl}^k-\varrho^k\|\StatRes^k\|_1.
\end{ResizedEq}
Substituting $\FPPR{O}{\StackConf}^k=-B_\StackConf^ks_\StackConf^k-[\StatJac^k]^T\mu_\StackConf^k$ and $\FPPR{O}{\StackCtrl}^k=-B_\StackCtrl^ks_\StackCtrl^k-[\StatJacU^k]^T\mu_\StackConf^k-\mu_\StackCtrl^k$ from the first two lines of~\prettyref{eq:sqp-kktqp}, and using $\StatJac^ks_\StackConf^k+\StatJacU^ks_\StackCtrl^k=-\StatRes^k$ once more,
\begin{ResizedEq}
[D\phi_\varrho]^k=&-[s_\StackConf^k]^TB_\StackConf^ks_\StackConf^k-[s_\StackCtrl^k]^TB_\StackCtrl^ks_\StackCtrl^k\\
&+[\StatRes^k]^T\mu_\StackConf^k-[s_\StackCtrl^k]^T\mu_\StackCtrl^k-\varrho^k\|\StatRes^k\|_1.
\end{ResizedEq}
Because $\mu_{u_i}^k\in N_{\mathcal{D}}(u_i^k+s_{u_i}^k)$ lies in the normal cone of the convex set $\mathcal{D}$ at $u_i^k+s_{u_i}^k$, and $u_i^k\in\mathcal{D}$ is a feasible comparison point, we have $\langle\mu_{u_i}^k,\,u_i^k-(u_i^k+s_{u_i}^k)\rangle\leq0$, i.e.\ $\langle\mu_{u_i}^k,s_{u_i}^k\rangle\geq0$; summing over $i$ gives $[s_\StackCtrl^k]^T\mu_\StackCtrl^k\geq0$. Moreover $[\StatRes^k]^T\mu_\StackConf^k\leq\sum_i\|\mu_{\theta_i}^k\|\,\|\StatRes_i^k\|_1$ by the per-block Cauchy--Schwarz inequality together with $\|\StatRes_i^k\|_2\leq\|\StatRes_i^k\|_1$. Therefore
\begin{ResizedEq}
[D\phi_\varrho]^k\leq&-[s_\StackConf^k]^TB_\StackConf^ks_\StackConf^k-[s_\StackCtrl^k]^TB_\StackCtrl^ks_\StackCtrl^k\\
&-\sum_{i}(\varrho^k-\|\mu_{\theta_i}^k\|)\|\StatRes_i^k\|_1\leq[\bar D\phi_\varrho]^k,
\end{ResizedEq}
where the last step uses $\lambda_{\min}(B^k)\geq1/\LHessianCond$ from~\prettyref{eq:sqp-hessian} and $\varrho^k-\|\mu_{\theta_i}^k\|\geq\bar\varrho$. Since $B^k\succ0$ with $\lambda_{\min}(B^k)\geq1/\LHessianCond$, the quadratic term vanishes iff $s^k=0$, while $\bar\varrho\|\StatRes^k\|_1=0$ iff $\StatRes^k=0$; hence $[\bar D\phi_\varrho]^k\leq0$ with equality iff $s^k=0$ and $\StatRes^k=0$, and $[\bar D\phi_\varrho]^k<0$ otherwise.
\end{proof}
\begin{lemma}[QP step and multiplier bounds]
\label{lem:sqp-multiplier}
Suppose $u_i^k\in\mathcal{D}$ for all $i$, the singularity check passes so that $\sigma_{\min}(\StatJac^k)\geq\SigmaMinLICQ(\StatJac)>0$ (\prettyref{lem:jacobian-lower-bound}), the QP Hessian satisfies~\prettyref{eq:sqp-hessian}, and $\|\StatRes^k\|_1\leq R$ for a residual ceiling $R>0$. Then $\|s^k\|\leq c_s(R)$ and $\|\mu_\StackConf^k\|\leq c_\mu(R)$, where, for the fixed mesh and compact region, $c_s(R)$ and $c_\mu(R)$ are finite and nondecreasing in $R$ and also depend on the fixed LICQ modulus $\SigmaMinLICQ(\StatJac)$, objective-gradient bound $\LObjGradBound$, and QP-conditioning constant $\LHessianCond$ appearing in their formulas, but are independent of $\EpsStop$.
\end{lemma}
\begin{proof}[Proof of~\prettyref{lem:sqp-multiplier}]
Step bound: Construct a feasible reference step $s^\star=(s_\StackConf^\star,s_\StackCtrl^\star)$ with $s_\StackCtrl^\star=0$ (which trivially meets the last line of~\prettyref{eq:sqp-kktqp}, the KKT characterization established in~\prettyref{lem:sqp-qp}, as $u_i^k\in\mathcal{D}$) and $s_\StackConf^\star=-[\StatJac^k]^{-1}\StatRes^k$, so that, by~\prettyref{lem:jacobian-lower-bound} and $\|\StatRes^k\|_2\leq\|\StatRes^k\|_1\leq R$,
\begin{ResizedEq}
\|s^\star\|\leq\frac{\|\StatRes^k\|}{\SigmaMinLICQ(\StatJac)}\leq\frac{R}{\SigmaMinLICQ(\StatJac)}.
\end{ResizedEq}
By fact~\ref{it:fact-gradO-lip} and~\prettyref{eq:sqp-hessian}, the QP objective at $s^\star$ is bounded above by
\begin{ResizedEq}
&\LObjGradBound\|s^\star\|+\frac{\LHessianCond}{2}\|s^\star\|^2\\
\leq&\LObjGradBound\frac{R}{\SigmaMinLICQ(\StatJac)}+\frac{\LHessianCond}{2}\frac{R^2}{\SigmaMinLICQ(\StatJac)^2}
\triangleq c_O(R).
\end{ResizedEq}
Conversely, for any feasible $s$ with $\|s\|>\LHessianCond \LObjGradBound+\sqrt{\LHessianCond^2\LObjGradBound^2+2\LHessianCond c_O(R)}\triangleq c_s(R)$, the QP objective is at least $\tfrac{1}{2\LHessianCond}\|s\|^2-\LObjGradBound\|s\|>c_O(R)$, so no such $s$ can be optimal. Hence $\|s^k\|\leq c_s(R)$.
Multiplier bound: The first line of~\prettyref{eq:sqp-kktqp} gives $[\StatJac^k]^T\mu_\StackConf^k=-\FPPR{O}{\StackConf}^k-B_\StackConf^ks_\StackConf^k$; using $\sigma_{\min}([\StatJac^k]^T)=\sigma_{\min}(\StatJac^k)\geq\SigmaMinLICQ(\StatJac)$,
\begin{ResizedEq}
\|\mu_\StackConf^k\|\leq&\frac{\|\FPPR{O}{\StackConf}^k\|+\LHessianCond\|s_\StackConf^k\|}{\SigmaMinLICQ(\StatJac)}\\
\leq&\frac{\LObjGradBound+\LHessianCond c_s(R)}{\SigmaMinLICQ(\StatJac)}\triangleq c_\mu(R),
\end{ResizedEq}
so that $\|\mu_\StackConf^k\|\leq c_\mu(R)$. For the fixed mesh and compact region, both bounds are finite and nondecreasing in $R$; their other dependencies are the fixed quantities $\SigmaMinLICQ(\StatJac)$, $\LObjGradBound$, and $\LHessianCond$ displayed above, and they are independent of $\EpsStop$.
\end{proof}
\begin{lemma}[Uniform residual, step, and penalty bounds]
\label{lem:sqp-stability}
Consider a run segment of~\prettyref{alg:sqp-inner} before any singularity return, so every iteration that reaches the QP has passed the check of~\prettyref{ln:sqp-inner-singular} at its current iterate and satisfies $\sigma_{\min}(\StatJac^k)\geq\SigmaMinLICQ(\StatJac)>0$ (\prettyref{lem:jacobian-lower-bound}); suppose also that the QP Hessian satisfies~\prettyref{eq:sqp-hessian}. Suppose further that the initialization obeys $u_i^0\in\mathcal{D}$ for all $i$ and $\|\StatRes^0\|_1\leq\LResidualWork$, with $\StackConf^0$ in the configuration region of~\prettyref{lem:trajectory-bound} induced by ceiling $\LResidualWork$. Let
\begin{ResizedEq}
\label{eq:sqp-stability-trial-region}
\mathcal K'\triangleq\left[\mathcal R_\theta(\LResidualWork)\oplus\mathcal B_\theta(c_s(\LResidualWork))\right]\times\mathcal D^{|\mathcal I|},
\end{ResizedEq}
where $\mathcal R_\theta(\LResidualWork)$ is that configuration region, $\mathcal B_\theta(\rho)$ is the closed configuration-space ball of radius $\rho$, and $\oplus$ denotes the Minkowski sum. Define the per-step residual ceiling by
\begin{ResizedEq}
\label{eq:sqp-stability-consistency}
\LResidual\triangleq\LResidualWork+(\dim\StatRes)\left[\sigma_{\max}(\StatJac)+\sigma_{\max}(\StatJacU)\right]c_s(\LResidualWork),
\end{ResizedEq}
where $\sigma_{\max}(\StatJac),\sigma_{\max}(\StatJacU)$ are the operator-norm bounds of~\prettyref{lem:jacobian-upper-bound} on the fixed compact joint trial region $\mathcal K'$ (so $\LResidualWork\leq\LResidual$ automatically). Then~\prettyref{alg:sqp-inner} generates iterates with
\begin{enumerate}[label=(\alph*),ref=(\alph*),leftmargin=*]
\item\label{it:stab-feas} $u_i^k\in\mathcal{D}$ for all $i$ and all $k$;
\item\label{it:stab-viol} a uniform violation bound
\begin{ResizedEq}
\|\StatRes^k\|_1\leq\LResidual\quad\text{for all }k,
\end{ResizedEq}
so that every iterate remains in the region of~\prettyref{lem:trajectory-bound}; and hence
\item\label{it:stab-bounds} uniformly bounded steps, multipliers, and penalty,
\begin{ResizedEq}
&\|s^k\|\leq c_s(\LResidual)\\
&\|\mu_\StackConf^k\|\leq c_\mu(\LResidual)\\
&\varrho^k\leq\varrho_{\max}\triangleq\max\{\varrho^0,c_\mu(\LResidual)+\bar\varrho\},
\end{ResizedEq}
where $\varrho^k$ denotes the post-update penalty used by the $k$th merit line search; this post-update sequence is nondecreasing across accepted line searches.
\end{enumerate}
\end{lemma}
\begin{proof}[Proof of~\prettyref{lem:sqp-stability}]
Penalty-state convention and control feasibility~\ref{it:stab-feas}. Let $\varrho_k^-$ denote the penalty entering iteration $k$ and
\begin{ResizedEq}
\varrho_k^+\triangleq\max\{\varrho_k^-,\max_i\|\mu_{\theta_i}^k\|+\bar\varrho\}
\end{ResizedEq}
its post-update value used by that iteration's merit line search. The in-place assignment in~\prettyref{alg:sqp-inner} persists across loop iterations, so after an accepted step $\varrho_{k+1}^-=\varrho_k^+$; the algorithm's post-update notation $\varrho^k$ denotes $\varrho_k^+$. The controls stay feasible by induction: if $u_i^k\in\mathcal{D}$, the last line of~\prettyref{eq:sqp-kktqp}, from~\prettyref{lem:sqp-qp}, gives $u_i^k+s_{u_i}^k\in\mathcal{D}$, and $u_i^{k+1}=u_i^k+[\tau^\ell]^ks_{u_i}^k$ is a convex combination of the two, hence belongs to $\mathcal{D}$ by convexity; the base case is $u_i^0\in\mathcal{D}$. Consequently~\prettyref{lem:sqp-multiplier} is applicable at every iterate whose violation is controlled and that reaches the QP.

Uniform violation bound~\ref{it:stab-viol}. We prove by induction that $\|\StatRes^k\|_1\leq\LResidual$ for every iterate generated before a singularity return. The base case follows from $\|\StatRes^0\|_1\leq\LResidualWork\leq\LResidual$. For the inductive step, assume $\|\StatRes^k\|_1\leq\LResidual$. Then $\|\StatRes_i^k\|_2\leq\LResidual$ for every $i$, so~\prettyref{lem:trajectory-bound} places the current configuration in its compact region for ceiling $\LResidual$. If iteration $k$ reaches the QP, the algorithmic check at~\prettyref{ln:sqp-inner-singular} has passed at this current iterate, so~\prettyref{lem:jacobian-lower-bound} supplies the LICQ modulus $\SigmaMinLICQ(\StatJac)>0$; together with~\prettyref{eq:sqp-hessian}, this makes~\prettyref{lem:sqp-multiplier} applicable. No singularity check is required at intermediate trial points, where only the upper Jacobian bounds on $\mathcal K'$ are used.

The violation can increase only on an inside-band step. Indeed, when $\|\StatRes^k\|_1>\LResidualWork$, the second branch of the safeguarded line search~\prettyref{eq:sqp-linesearch} enforces $\|\StatRes^{k+1}\|_1\leq\|\StatRes^k\|_1\leq\LResidual$. Every residual-increasing step therefore starts with $\|\StatRes^k\|_1\leq\LResidualWork$, so its configuration lies in $\mathcal R_\theta(\LResidualWork)$ and~\prettyref{lem:sqp-multiplier} with $R=\LResidualWork$ gives $\|s^k\|\leq c_s(\LResidualWork)$. Writing the accepted joint update as $[\tau^\ell]^ks^k$, its configuration segment lies in $\mathcal R_\theta(\LResidualWork)\oplus\mathcal B_\theta(c_s(\LResidualWork))$ because $\|s_\StackConf^k\|\leq\|s^k\|$. Its control segment lies in $\mathcal D^{|\mathcal I|}$ because each point is a convex combination of $u_i^k$ and the QP-feasible endpoint $u_i^k+s_{u_i}^k$. Hence the complete joint segment lies in the compact set $\mathcal K'$ of~\prettyref{eq:sqp-stability-trial-region}.

We now justify the Jacobian bounds on this enlarged trial region directly. On the fixed finite mesh every timestep is positive, so the discrete velocities and accelerations entering $\StatRes$ are affine functions of the stacked configurations. The set $\mathcal R_\theta(\LResidualWork)$ and the closed ball $\mathcal B_\theta(c_s(\LResidualWork))$ are compact; their Minkowski sum is compact, and so is $\mathcal D^{|\mathcal I|}$. Thus $\mathcal K'$ is compact. The mollified potential $\bar U$ is globally $C^2$; on the open control domain containing $\mathcal D$, the mollified control work has continuous configuration Hessian and mixed configuration--control derivative; and the everywhere-defined smoothed dissipation has a velocity gradient that is jointly $C^1$. Consequently the full residual Jacobian $\nabla_{(\StackConf,\StackCtrl)}\StatRes=(\StatJac,\StatJacU)$ is defined and continuous on an open neighborhood of $\mathcal K'$. In~\prettyref{eq:sqp-stability-consistency}, define the proof-local bounds
\begin{ResizedEq}
\sigma_{\max}(\StatJac)\triangleq&\sup_{(\StackConf,\StackCtrl)\in\mathcal K'}\|\StatJac\|\\
\sigma_{\max}(\StatJacU)\triangleq&\sup_{(\StackConf,\StackCtrl)\in\mathcal K'}\|\StatJacU\|,
\end{ResizedEq}
which are finite by compactness and continuity, and are the operator-norm bounds of~\prettyref{lem:jacobian-upper-bound} for $\StatJac$ and $\StatJacU$ recorded at the result level.

The mean-value formula for $\StatRes$ along this segment, with these full-residual-Jacobian bounds on $\mathcal K'$, gives
\begin{ResizedEq}
&\|\StatRes^{k+1}\|_1\\
\leq&\sum_{i,j}\left(|[\StatRes_i^k]_j|+\left|\left[\left({\textstyle\int_0^1}\nabla_{(\StackConf,\StackCtrl)}\StatRes\,dt\right)[\tau^\ell]^ks^k\right]_{ij}\right|\right)\\
\leq&\|\StatRes^k\|_1+(\dim\StatRes)\left[\sigma_{\max}(\StatJac)+\sigma_{\max}(\StatJacU)\right]\|s^k\|\\
\leq&\LResidualWork+(\dim\StatRes)\left[\sigma_{\max}(\StatJac)+\sigma_{\max}(\StatJacU)\right]c_s(\LResidualWork)=\LResidual.
\end{ResizedEq}
Thus either safeguard branch gives $\|\StatRes^{k+1}\|_1\leq\LResidual$, closing the induction. Every generated iterate therefore remains in the compact trajectory region associated with ceiling $\LResidual$.

Bounded steps, multipliers, and penalty~\ref{it:stab-bounds}. Applying~\prettyref{lem:sqp-multiplier} with $R=\LResidual$ yields $\|s^k\|\leq c_s(\LResidual)$ and $\|\mu_\StackConf^k\|\leq c_\mu(\LResidual)$. Moreover,
\begin{ResizedEq}
\varrho_k^+=\max\{\varrho_k^-,\max_i\|\mu_{\theta_i}^k\|+\bar\varrho\},
\qquad
\varrho_{k+1}^-=\varrho_k^+.
\end{ResizedEq}
Since $\max_i\|\mu_{\theta_i}^k\|\leq\|\mu_\StackConf^k\|\leq c_\mu(\LResidual)$, induction from $\varrho_0^-=\varrho^0$ gives $\varrho_k^-,\varrho_k^+\leq\varrho_{\max}\triangleq\max\{\varrho^0,c_\mu(\LResidual)+\bar\varrho\}$ for every $k$. The carried-forward update also gives $\varrho_k^+\leq\varrho_{k+1}^+$ whenever both line searches occur, so the post-update sequence denoted $\varrho^k$ is nondecreasing and uniformly bounded.
\end{proof}
\begin{lemma}[Main SQP line-search well-posedness]
\label{lem:sqp-linesearch}
Under the hypotheses of~\prettyref{lem:sqp-stability}, with the iterates confined to the norm-bounded region of~\prettyref{lem:trajectory-bound}, suppose the singularity check of~\prettyref{ln:sqp-inner-singular} passes, the QP Hessian satisfies~\prettyref{eq:sqp-hessian}, the penalty has been updated so that the descent condition of~\prettyref{lem:sqp-qp} holds ($\varrho^k\geq\max_i\|\mu_{\theta_i}^k\|+\bar\varrho$), and the stopping test of~\prettyref{ln:sqp-inner-return} has not triggered, so that $\|s^k\|\geq\EpsStop/\LHessianCond$ or $\|\StatRes^k\|_1\geq\EpsStop$. Then there exists $\bar\tau>0$, uniform over all such iterations, such that every step $[\tau^\ell]^k<\bar\tau$ satisfies the safeguarded line search~\prettyref{eq:sqp-linesearch}; hence the backtracking accepts a step bounded below.
\end{lemma}
\begin{proof}[Proof of~\prettyref{lem:sqp-linesearch}]
Abbreviate $\beta\triangleq[\tau^\ell]^k$, set $C_s\triangleq c_s(\LResidual)$, and introduce the proof-local trial region
\begin{ResizedEq}
\mathcal K_{\rm ls}\triangleq\left[\mathcal R_\theta(\LResidual)\oplus\mathcal B_\theta(C_s)\right]\times\mathcal D^{|\mathcal I|}.
\end{ResizedEq}
By~\prettyref{lem:sqp-stability}, the current configuration belongs to $\mathcal R_\theta(\LResidual)$, the complete step satisfies $\|s^k\|\leq C_s$, and the current controls belong to $\mathcal D^{|\mathcal I|}$. By~\prettyref{lem:sqp-qp}, the QP-feasible endpoint also belongs to $\mathcal D^{|\mathcal I|}$, so convexity of $\mathcal D$ keeps its entire control segment feasible. Consequently, for every $t,\beta\in[0,1]$, every joint segment point $z^k+t\beta s^k$, where $z\triangleq(\StackConf,\StackCtrl)$, lies in the compact set $\mathcal K_{\rm ls}$.

By the same data regularity used in~\prettyref{lem:sqp-stability}, the full residual Jacobian $J_r\triangleq\nabla_{(\StackConf,\StackCtrl)}\StatRes=(\StatJac,\StatJacU)$ is defined and continuous on an open neighborhood of $\mathcal K_{\rm ls}$. On this trial region define the nondecreasing supremum modulus
\begin{ResizedEq}
\omega_r^{\rm ls}(a)\triangleq\sup\left\{\|J_r(z')-J_r(z)\|:\ z,z'\in\mathcal K_{\rm ls},\ \|z'-z\|\leq a\right\}.
\end{ResizedEq}
Compactness and continuity give $\omega_r^{\rm ls}(a)\to0$ as $a\downarrow0$. The integral mean-value formula and the QP identity $J_r(z^k)s^k=-\StatRes^k$ yield
\begin{ResizedEq}
\StatRes(z^k+\beta s^k)=(1-\beta)\StatRes^k+r_k(\beta),
\end{ResizedEq}
where
\begin{ResizedEq}
\|r_k(\beta)\|\leq\beta\|s^k\|\omega_r^{\rm ls}(\beta\|s^k\|).
\end{ResizedEq}
Thus
\begin{ResizedEq}
\|\StatRes(z^k+\beta s^k)\|_1\leq(1-\beta)\|\StatRes^k\|_1+(\dim\StatRes)\beta\|s^k\|\omega_r^{\rm ls}(\beta\|s^k\|).
\end{ResizedEq}

Since $\nabla O$ is Lipschitz with modulus $\LObjGradLip$ (fact~\ref{it:fact-gradO-lip}), combining its descent estimate with the accepted QP descent surrogate and the residual estimate gives
\begin{ResizedEq}
\label{eq:sqp-merit-decrease}
\phi_{\varrho^k}(z^k+\beta s^k)\leq&\phi_{\varrho^k}^k+\beta[\bar D\phi_\varrho]^k+\beta A(\beta),\\
A(\beta)\triangleq&\frac{\LObjGradLip}{2}\beta C_s^2+\varrho_{\max}(\dim\StatRes)C_s\omega_r^{\rm ls}(\beta C_s).
\end{ResizedEq}
The function $A$ is proof-local and satisfies $A(\beta)\to0$ as $\beta\downarrow0$.

Case I ($\|\StatRes^k\|_1\leq\LResidualWork$): because the stopping test of~\prettyref{ln:sqp-inner-return} has not triggered,
\begin{ResizedEq}
-[\bar D\phi_\varrho]^k=\frac{1}{\LHessianCond}\|s^k\|^2+\bar\varrho\|\StatRes^k\|_1\geq m_I\triangleq\min\left\{\frac{\EpsStop^2}{\LHessianCond^3},\bar\varrho\EpsStop\right\}>0.
\end{ResizedEq}
Since $A(\beta)\to0$, there is a uniform $\bar\tau_I>0$ such that $A(\beta)<(1-\alpha)m_I$ whenever $0<\beta<\bar\tau_I$. For such $\beta$,~\prettyref{eq:sqp-merit-decrease} gives explicitly
\begin{ResizedEq}
&\phi_{\varrho^k}(z^k+\beta s^k)\\
<&\phi_{\varrho^k}^k+\beta[\bar D\phi_\varrho]^k+(1-\alpha)\beta(-[\bar D\phi_\varrho]^k)=\phi_{\varrho^k}^k+\alpha\beta[\bar D\phi_\varrho]^k,
\end{ResizedEq}
which is the strict first branch of~\prettyref{eq:sqp-linesearch}.

Case II ($\|\StatRes^k\|_1>\LResidualWork$): here $-[\bar D\phi_\varrho]^k\geq\bar\varrho\LResidualWork$. Hence there is a uniform $\bar\tau_{II}^{\phi}>0$ such that $A(\beta)<(1-\alpha)\bar\varrho\LResidualWork$ for $0<\beta<\bar\tau_{II}^{\phi}$, and the same calculation proves the strict merit Armijo inequality. Separately, because $\omega_r^{\rm ls}(\beta C_s)\to0$, there is a uniform $\bar\tau_{II}^{r}>0$ such that
\begin{ResizedEq}
(\dim\StatRes)C_s\omega_r^{\rm ls}(\beta C_s)\leq\LResidualWork
\end{ResizedEq}
for $0<\beta<\bar\tau_{II}^{r}$. The residual estimate then gives
\begin{ResizedEq}
\|\StatRes(z^k+\beta s^k)\|_1\leq(1-\beta)\|\StatRes^k\|_1+\beta\LResidualWork\leq\|\StatRes^k\|_1,
\end{ResizedEq}
which is the required non-strict safeguard.

Taking $\bar\tau\triangleq\min\{1,\bar\tau_I,\bar\tau_{II}^{\phi},\bar\tau_{II}^{r}\}>0$, every trial step $0<\beta<\bar\tau$ passes the applicable branch. Under the first-passing backtracking convention, the accepted trial therefore has size at least $\tau\bar\tau>0$, as in~\prettyref{lem:sqp-constraint-termination}.
\end{proof}

\begin{proof}[Proof of~\prettyref{thm:sqp-convergence}]
Throughout, the input requirements of~\prettyref{alg:sqp}---feasible initial controls $u_i^0\in\mathcal{D}$, a positive initial penalty $\varrho^0>0$, a positive penalty margin $\bar\varrho>0$, a positive stopping tolerance $\EpsStop>0$, and QP Hessians satisfying~\prettyref{eq:sqp-hessian} with $\LHessianCond\geq1$ (enforced by~\prettyref{alg:sqp-inner} when it chooses $B_\StackConf^k,B_\StackCtrl^k$)---are in force; we invoke these below without further comment.

Finitely many subdivisions: The bisection count is bounded globally and independently of any trajectory region. Recall from~\prettyref{eq:delta-star} the uniform timestep threshold $\delta^\star$ of~\prettyref{lem:jacobian-lower-bound},
where $\LCurvature$ is the global two-sided potential-curvature modulus of~\prettyref{lem:mollification} and $\LControlWorkCurvature$ is the region-independent control-work curvature bound~\prettyref{eq:sqp-cw}, valid globally over configurations and controls in $\mathcal{D}$. The unconditional diagonal estimate~\prettyref{eq:sqp-diagonal-lower} of~\prettyref{lem:jacobian-lower-bound} therefore gives, at every configuration and independently of any residual ceiling or compact region, $\StatJac_{ii}\succeq\LLambdaMin I$ whenever $\delta_i\leq\delta^\star$. This stronger positive-definite bound implies $\sigma_{\min}(\StatJac_{ii})\geq\LLambdaMin$, so such an interval passes the singularity check and is never subdivided. Only intervals with $\delta_i>\delta^\star$ can fail the check, and each bisection halves such an offending interval. Let $G_0\triangleq\mathcal I^0$ be the finite set of optimized right endpoints of the initial mesh over the fixed planning horizon $[0,H]$, with $i\in G_0$ indexing the interval $(\gamma(i),i)$. The first refinable interval is $(0,\min\mathcal I^0)$; the fixed history interval $(\gamma(0),0)$ is not indexed by $G_0$ and is never subdivided. Each initial right endpoint $i\in G_0$ spawns a binary subdivision tree of depth at most $\left\lceil\log_2(\delta_i^0/\delta^\star)\right\rceil$, so the total number of bisections over the whole run is bounded by
\begin{ResizedEq}
N_{\mathrm{sub}}\leq\sum_{i\in G_0}\left(2^{\max\left\{0,\left\lceil\log_2(\delta_i^0/\delta^\star)\right\rceil\right\}}-1\right)<\infty,
\end{ResizedEq}
with $\delta_i^0$ the length of initial interval $i\in G_0$; this count uses only the finite initial mesh and the region-independent threshold $\delta^\star$. Moreover, every active timestep is bounded below throughout the run by the positive floor $\delta_{\mathrm{floor}}\triangleq\min\{\delta^\star/2,\min_{i\in G_0}\delta_i^0\}>0$: a newly created child has a parent of length greater than $\delta^\star$ and hence has length greater than $\delta^\star/2$, while an interval never subdivided retains its initial length. After at most $N_{\mathrm{sub}}$ bisections no interval can fail the singularity check again, and the run decomposes into at most $N_{\mathrm{sub}}+1$ mesh episodes, namely maximal stretches of computation on a single fixed mesh, each begun either on the initial mesh or immediately after a subdivision and ended either by the next subdivision or by the algorithm's return.

Well-posedness: Each inner iteration either fails the singularity check of~\prettyref{ln:sqp-inner-singular} and returns an index for subdivision, or has LICQ in force, in which case $\text{QP}^k$ is solvable (\prettyref{lem:sqp-qp}) and a line-search step exists (\prettyref{lem:sqp-linesearch}); the constraint-satisfaction phase is finite by~\prettyref{lem:sqp-constraint-termination}, and~\prettyref{alg:subd} restores the singularity check in finitely many bisections through~\prettyref{eq:sqp-diagonal-lower}, so~\prettyref{alg:sqp} is well defined.

Proof constants and compact region: We next induct over the at most $N_{\mathrm{sub}}+1$ mesh episodes to establish finite episode-initial residuals and per-episode compactness before defining a common whole-run ceiling. By the stated initialization conditions, $\mathcal I^0$ is finite, the explicit initial arrays contain finitely many finite configuration values, every timestep is positive, and all controls lie in $\mathcal D$. The globally mollified potential and control-work terms and the everywhere-finite damping smoothing therefore make every per-step residual and residual derivative finite. Hence the initial stacked residual is finite, although it may be arbitrarily large. Suppose an episode on a fixed reached optimized set $\mathcal I_{\mathcal G}$ begins from an explicit tuple $(\mathcal I_{\mathcal G},\StackConf,\StackCtrl,\varrho)$ with finite residual. During its constraint-satisfaction phase,~\prettyref{eq:sqp-constraint-linesearch} keeps $\Lambda=\tfrac12\|\StatRes\|^2$ non-increasing, so that phase has the finite residual ceiling $\sqrt{\dim\StatResG}\,\|\StatResG^0\|_2$; together with the already-established floor $\delta_{\mathrm{floor}}>0$,~\prettyref{lem:trajectory-bound} confines its iterates to a compact region. If the phase reaches the working band and enters the inner SQP, let $\LResidualG$ denote the fixed-mesh ceiling defined by~\prettyref{eq:sqp-stability-consistency}. This is a genuine explicit definition, not a self-referential assumption: (i) the mesh $\mathcal{G}$ is fixed, so $\dim\StatResG$ and all finite-dimensional mesh data are fixed; (ii) by the residual-grows-only-inside-band mechanism of~\prettyref{lem:sqp-stability}, the step bound $c_s(\LResidualWork)$ and the objective-gradient bound $\LObjGradBound$ entering~\prettyref{eq:sqp-stability-consistency} are evaluated through~\prettyref{lem:sqp-multiplier} only at inside-band iterates, which lie in the band region $\mathcal{R}(\LResidualWork)$ of~\prettyref{lem:trajectory-bound}, determined by $\LResidualWork$ alone; in particular the fixed-mesh LICQ modulus $\SigmaMinLICQ(\StatJac)$, through its off-diagonal bound $C_2$, and hence $c_s(\LResidualWork)$ and $\LObjGradBound$ are finite on that region by~\prettyref{ass:objective}; and (iii) the enlarged region $\mathcal{R}'=\mathcal{R}(\LResidualWork)\oplus\mathcal{B}(c_s(\LResidualWork))$ is built from the already-fixed $c_s(\LResidualWork)$, so the excursion operator norms $\sigma_{\max}(\StatJac)$ and $\sigma_{\max}(\StatJacU)$ on its right-hand side are finite before $\LResidualG$ is defined. Thus the right-hand side of~\prettyref{eq:sqp-stability-consistency} depends only on the working band, the fixed mesh, and the global curvature moduli and defines a finite $\LResidualG$ without reference to $\mathcal{R}(\LResidualG)$; the per-mesh definitional argument is therefore non-circular. By~\prettyref{lem:sqp-stability}, the inner-SQP iterates of this episode satisfy $\|\StatRes^k\|_1\leq\LResidualG$, and~\prettyref{lem:trajectory-bound} again confines them to an episode-specific compact region. If the episode ends by subdivision, the neighboring states used by~\prettyref{alg:subd} belong to this preceding episode's compact region and its controls belong to $\mathcal{D}$; the inserted time is the midpoint of a positive interval, the inserted states and controls are finite convex interpolants of finite neighboring values, and convexity preserves membership of the new control in $\mathcal D$. Since the new timesteps remain at least $\delta_{\mathrm{floor}}>0$, continuity of the residual map implies that the next episode's initialization residual is finite. This proves the induction, without presupposing a common compact region, over all reached episodes. There are only finitely many such episodes, so we may now define
\begin{ResizedEq}
R_{\mathrm{pre}}\triangleq\max\Big\{&\max_{\mathcal{G}\text{ a mesh reached}}\sqrt{\dim\StatResG}\,\|\StatResG^0\|_2,\\
&\max_{\substack{\mathcal{G}\text{ a mesh reached whose episode}\\\text{enters inner SQP}}}\LResidualG\Big\}<\infty.
\end{ResizedEq}
Here $\StatResG^0$ is the residual at the point initializing the constraint-satisfaction phase on mesh $\mathcal{G}$. If no reached episode enters inner SQP, the second maximum is simply omitted. The maximum is finite because at most $N_{\mathrm{sub}}+1$ meshes are reached, every such initialization residual is finite by the episode induction, and every $\LResidualG$ defined for an episode entering inner SQP is finite by the preceding fixed-mesh definition. The first maximum covers every constraint-satisfaction phase, including one that exits with a singularity index, while the second covers exactly the episodes that generate inner-SQP iterates. Thus every iterate of the whole run satisfies $\|\StatRes^k\|_1\leq R_{\mathrm{pre}}$. Applying~\prettyref{lem:trajectory-bound} with this whole-run ceiling and $\delta_{\mathrm{floor}}>0$ confines all iterates to a single compact region $\mathcal{R}_{\mathrm{pre}}$. Consequently all region-dependent constants used later at general iterates---the Jacobian operator norms $\sigma_{\max}(\StatJac),\sigma_{\max}(\StatJacU)$ of~\prettyref{lem:jacobian-upper-bound}, the fixed-mesh step and multiplier bounds $c_s(\LResidualG),c_\mu(\LResidualG)$, and the off-diagonal bounds entering the fixed-mesh LICQ moduli---are finite on $\mathcal{R}_{\mathrm{pre}}$; none of them was used to obtain $\delta^\star$ or $N_{\mathrm{sub}}$, and none re-enters the already-closed definition~\prettyref{eq:sqp-stability-consistency} of $\LResidualG$.

Bounded Penalty: Within any mesh episode, run the constraint-satisfaction phase (\prettyref{alg:constraint-solve}) from that episode's starting point. Since~\prettyref{eq:sqp-constraint-linesearch} keeps $\Lambda$ non-increasing, every pre-phase iterate satisfies $\|\StatRes_i^k\|_2\leq\|\StatRes^k\|_2\leq\|\StatResG^0\|_2$, hence $\|\StatRes^k\|_1\leq\sqrt{\dim\StatResG}\,\|\StatResG^0\|_2\leq R_{\mathrm{pre}}$ and the iterate lies in $\mathcal{R}_{\mathrm{pre}}$, while the controls stay at the feasible initial $u_i^0\in\mathcal{D}$ because~\prettyref{alg:constraint-solve} updates only $\StackConf$; thus~\prettyref{lem:sqp-constraint-linesearch} and~\ref{lem:sqp-constraint-termination} apply and the phase terminates finitely through one of its two exits.

If~\prettyref{lem:sqp-constraint-termination} returns a singularity index, the episode ends without entering inner SQP. The signatures of~\prettyref{alg:constraint-solve} and~\ref{alg:subd} contain no penalty argument or return value, and neither algorithm changes $\varrho$; hence the incoming penalty is carried unchanged to the next episode. If instead the constraint phase reaches the working band, take the explicitly returned updated pair $(\mathcal I,\StackConf)$ together with the unchanged $\StackCtrl$ and carried penalty $\varrho$ as the inner-SQP initialization. The initialization hypotheses of~\prettyref{lem:sqp-stability}---$u_i^0\in\mathcal{D}$ for all $i$ and $\|\StatRes^0\|_1\leq\LResidualWork$---then hold, so~\prettyref{lem:sqp-stability} applies and supplies, once and for all on this fixed mesh: control feasibility $u_i^k\in\mathcal{D}$ for all $i,k$; the uniform violation bound $\|\StatRes^k\|_1\leq\LResidualG$, so that~\prettyref{lem:trajectory-bound} keeps applying with $R_{\mathrm{pre}}$; and the uniformly bounded steps, multipliers, and penalty $\|s^k\|\leq c_s(\LResidualG)$, $\|\mu_\StackConf^k\|\leq c_\mu(\LResidualG)$, and $\{\varrho^k\}$ bounded by the fixed-mesh $\varrho_{\max}$. This bootstrap follows from the line-search safeguard of~\prettyref{eq:sqp-linesearch}, independently of the conclusions of the present theorem, so no circularity arises.

Finite induction over the reached episodes now bounds the whole-run penalty by the maximum of the initial $\varrho$ and, for every reached episode that enters inner SQP, its fixed-mesh candidate $c_\mu(\LResidualG)+\bar\varrho$. If no episode enters inner SQP, this maximum reduces to the initial penalty. Singularity-only episodes add no candidate because they leave $\varrho$ unchanged, and there are only finitely many inner-SQP episodes, so this maximum is finite.

Finite termination: We show each mesh episode is finite and then induct over the finitely many episodes. Fix a mesh episode. Its constraint-satisfaction phase terminates in finitely many iterations by~\prettyref{lem:sqp-constraint-termination} (Bounded Penalty), returning either a singularity index---which ends the episode with one subdivision---or a point in the working band. In the latter case the episode proceeds to the inner SQP on the fixed mesh, which either returns a singularity index in finitely many iterations---again ending the episode with one subdivision---or produces no singularity return; it remains to show that an inner SQP running on a fixed mesh with no singularity return terminates in finitely many iterations. Suppose not, so that on this fixed mesh the inner SQP neither returns a singularity index nor meets the stopping test. Because the penalty update is monotone, $\varrho^{k-1}\leq\varrho^k$, and $O\geq0$ (fact~\ref{it:fact-Obounded}), the accepted decrease $\phi_{\varrho^k}^{k+1}\leq\phi_{\varrho^k}^k+\alpha[\tau^\ell]^k[\bar D\phi_\varrho]^k$ telescopes with the ratio chain
\begin{ResizedEq}
\frac{\varrho^0}{\varrho^k}\phi_{\varrho^k}^{k+1}\leq\phi_{\varrho^0}^0+\alpha\sum_{m=0}^k\frac{\varrho^0}{\varrho^m}[\tau^\ell]^m[\bar D\phi_\varrho]^m,
\end{ResizedEq}
where each step uses $(\varrho^{m-1}/\varrho^m)\phi_{\varrho^m}^m\leq\phi_{\varrho^{m-1}}^m$, valid because $O\geq0$ and $\varrho^{m-1}/\varrho^m\leq1$. Since $\phi_{\varrho^k}^{k+1}\geq0$ and $\varrho^m\leq\varrho_{\max}$,
\begin{ResizedEq}
&-\frac{\varrho^0\alpha}{\varrho_{\max}}\sum_{m=0}^\infty[\tau^\ell]^m[\bar D\phi_\varrho]^m\\
\leq&-\alpha\sum_{m=0}^\infty\frac{\varrho^0}{\varrho^m}[\tau^\ell]^m[\bar D\phi_\varrho]^m\leq\phi_{\varrho^0}^0<\infty,
\end{ResizedEq}
so $\sum_k[\tau^\ell]^k(-[\bar D\phi_\varrho]^k)<\infty$. Because this inner SQP is assumed not to terminate, the stopping test of~\prettyref{ln:sqp-inner-return} never triggers, so~\prettyref{lem:sqp-linesearch} applies at every inner iteration and, since every iterate stays in the region of~\prettyref{lem:trajectory-bound}, bounds the accepted step sizes below by a uniform $\bar\tau>0$. Hence $\bar\tau\sum_k(-[\bar D\phi_\varrho]^k)\leq\sum_k[\tau^\ell]^k(-[\bar D\phi_\varrho]^k)<\infty$, and since every term $-[\bar D\phi_\varrho]^k\geq0$ we get $-[\bar D\phi_\varrho]^k\to0$. Because $-[\bar D\phi_\varrho]^k=\tfrac{1}{\LHessianCond}\|s^k\|^2+\bar\varrho\|\StatRes^k\|_1$ with both terms nonnegative (the penalty margin $\bar\varrho>0$ being an input requirement of~\prettyref{alg:sqp}), both $\|s^k\|\to0$ and $\|\StatRes^k\|_1\to0$. But then, for the positive stopping tolerance $\EpsStop$ (likewise required by~\prettyref{alg:sqp}), the stopping test $\|s^k\|<\EpsStop/\LHessianCond$ and $\|\StatRes^k\|_1<\EpsStop$ of~\prettyref{ln:sqp-inner-return} is eventually satisfied, contradicting non-termination. Hence an inner SQP on a fixed mesh with no singularity return terminates after finitely many iterations. In every case a mesh episode is finite; since at most $N_{\mathrm{sub}}$ subdivisions occur there are at most $N_{\mathrm{sub}}+1$ episodes, and~\prettyref{alg:sqp} either returns within some episode or, having exhausted its subdivisions, terminates within the final episode. Inducting over the at most $N_{\mathrm{sub}}+1$ mesh episodes, the whole run terminates after finitely many iterations.

Stationarity: The routines pass their state explicitly: accepted updates of~\prettyref{alg:constraint-solve} and~\ref{alg:sqp-inner} update the arrays on fixed $\mathcal I$, while~\prettyref{alg:subd} returns the refined tuple $(\mathcal I,\StackConf,\StackCtrl)$. Thus the point returned through~\prettyref{ln:sqp-inner-return} (and by~\prettyref{alg:sqp}) is the current $(\StackConf^k,\StackCtrl^k)$ indexed by the returned $\mathcal I$, and the $\EpsStop$-KKT verification below is carried out at that returned point. At return through~\prettyref{ln:sqp-inner-return}, $\|s^k\|<\EpsStop/\LHessianCond$ and $\|\StatRes^k\|_1<\EpsStop$; we verify the returned point $(\StackConf^k,\StackCtrl^k)$ satisfies the $\EpsStop$-KKT condition of~\prettyref{def:eps-KKT} with the multipliers $\mu_\StackConf^k,\mu_\StackCtrl^k$ of~\prettyref{eq:sqp-kktqp}. From the first line of~\prettyref{eq:sqp-kktqp}, $\|\FPPR{O}{\StackConf}^k+[\StatJac^k]^T\mu_\StackConf^k\|=\|B_\StackConf^ks_\StackConf^k\|\leq\LHessianCond\|s_\StackConf^k\|\leq\LHessianCond\|s^k\|\leq\EpsStop$, and identically from its second line $\|\FPPR{O}{\StackCtrl}^k+[\StatJacU^k]^T\mu_\StackConf^k+\mu_\StackCtrl^k\|=\|B_\StackCtrl^ks_\StackCtrl^k\|\leq\LHessianCond\|s^k\|\leq\EpsStop$. For the two-point control line, set $u_i'\triangleq u_i^k+s_{u_i}^k$: the fourth line of~\prettyref{eq:sqp-kktqp} gives $u_i'\in\mathcal{D}$ and $\mu_{u_i}^k\in N_{\mathcal{D}}(u_i')$, the returned control satisfies $u_i^k\in\mathcal{D}$ by control feasibility~\ref{it:stab-feas} of~\prettyref{lem:sqp-stability}, and $\|u_i^k-u_i'\|=\|s_{u_i}^k\|\leq\|s^k\|\leq\EpsStop/\LHessianCond\leq\EpsStop$ since $\LHessianCond\geq1$. Finally $\|\StatRes^k\|_1<\EpsStop$ by the stopping test. Hence the returned point satisfies the $\EpsStop$-KKT condition of~\prettyref{def:eps-KKT}.
\end{proof}

\begin{proof}[Proof of~\prettyref{thm:CITO-convergence}]
This is the guarantee for the outer loop~\prettyref{alg:CITO} that wraps~\prettyref{alg:sqp} with the Hamiltonian check~\prettyref{eq:cito-ham-check}. Throughout we abbreviate the discretization-error constant
\begin{ResizedEq}
C\triangleq(\|g\|+\LControlWorkForce)^2+\EpsStop^2,
\end{ResizedEq}
and write the two timestep thresholds that govern the two subdivision mechanisms in play: the sufficient singularity-refinement threshold $\delta^\star$ of~\prettyref{eq:delta-star} (used by~\prettyref{alg:sqp} internally) and the Hamiltonian threshold
\begin{ResizedEq}
\delta_H\triangleq\min\left(\frac{1}{\LCurvature},\frac{\lambda_{\min}(M)}{6}\right),
\end{ResizedEq}
which is exactly the timestep condition~\prettyref{eq:H-bound-timestep} of~\prettyref{lem:H-bound}. Both $\delta^\star$ and $\delta_H$ are built from $\lambda_{\min}(M)$, the two-sided potential-curvature modulus $\LCurvature$ of~\prettyref{lem:mollification}, the region-independent control-work Hessian bound $\LControlWorkCurvature$ of~\prettyref{eq:sqp-cw}, and the algorithm threshold $\LLambdaMin$; none depends on the residual ceiling or the current mesh. Before the first smooth solve, define the design floor
\begin{equation}
\label{eq:cito-delta-min}
\delta_{\min}\triangleq\min\!\left\{\min_{i\in G_0}\delta_i^0,\frac{\min(\delta^\star,\delta_H)}{2}\right\}>0.
\end{equation}
This algebraic quantity uses only the finite initial mesh and the two threshold formulas, is the admissible choice of the a priori timestep floor $\delta_{\min}$ prescribed by~\prettyref{lem:smooth-DRp-Regularity}, and will be shown below to bound every generated timestep from below. The damping function is already fixed by the hypothesis that~\prettyref{pb:PS-CRDS-M} is a smooth CRDS; this theorem selects no damping or smoothing parameter. The thresholds require no damping-curvature upper bound: in the diagonal Jacobian bound~\prettyref{eq:sqp-diagonal-lower} the damping Hessian is positive semidefinite and is discarded (\prettyref{lem:jacobian-lower-bound}), while the Hamiltonian threshold $\delta_H$ uses damping only through convexity and minimality at zero velocity (\prettyref{lem:H-bound}). Any constants used to upper-bound derivatives on compact sets may depend on the one fixed damping function, but they do not feed back into $\delta^\star$, $\delta_H$, or $\delta_{\min}$.

The three input requirements added over~\prettyref{thm:sqp-convergence}---\prettyref{ass:Slater} and the positive-gap clamp definitions of the constant-ordering remark of~\prettyref{sec:Hamiltonian}---are in force throughout; they are used only in Part~3. Thus~\prettyref{ass:Slater} is not used to prove finite termination and $\EpsStop$-KKT of each fixed smooth~\prettyref{pb:CITO-M} SQP call in Part~1. It enters in Part~3 to make the predecessor penalty energy and penalty-induced load finite, seed the Hamiltonian recursion, establish clamp inactivity, and complete the transfer to the original CRDS. We argue in three parts.

Step 1: each SQP call is finite, bounded-penalty, and $\EpsStop$-KKT. Every outer iteration of~\prettyref{alg:CITO} calls~\prettyref{alg:sqp} on the current mesh with parameters that meet its input requirements (the same $\LLambdaMin,\LResidualWork,\LHessianCond,\bar\varrho,\EpsStop$ throughout, controls $u_i\in\mathcal{D}$, and a positive penalty $\varrho$ threaded from the previous call). The current mesh is finite: the initial mesh is finite and each outer iteration inserts finitely many midpoints through~\prettyref{alg:subd}. Since~\prettyref{pb:PS-CRDS-M} is a smooth full-rank CRDS with objective satisfying~\prettyref{ass:objective}, the hypotheses of~\prettyref{thm:sqp-convergence} hold for this call, so it is well defined, performs finitely many (internal singularity-driven) subdivisions, generates a bounded merit-penalty sequence, and terminates finitely at a point satisfying the $\EpsStop$-KKT condition of~\prettyref{def:eps-KKT}; in particular the returned point has $\|\StatRes\|_1\leq\EpsStop$. Applied to the final outer iteration---the one that returns at~\prettyref{ln:cito-return}---this gives the $\EpsStop$-KKT conclusion of the present theorem verbatim, including $\|\StatRes\|_1\leq\EpsStop$.

Step 2: finitely many outer subdivisions, hence a bounded penalty. We first show any node whose interval satisfies $\delta_i\leq\delta_H$ passes the check~\prettyref{eq:cito-ham-check}, so it is never subdivided by the outer loop. Fix the point returned by the SQP call of the current outer iteration. By Part~1 it has $\|\StatRes\|_1\leq\EpsStop$, so each node obeys $\|\StatRes_i\|_2\leq\|\StatRes\|_1\leq\EpsStop$, and~\prettyref{lem:H-one-step} applies at that node with $\EpsCRDS=\EpsStop$
\begin{ResizedEq}
\label{eq:cito-onestep-raw}
\bar{\mathcal H}_i^+-\bar{\mathcal H}_{\gamma(i)}^+
\leq\frac{\LCurvature\delta_i^2+2\delta_i}{2\lambda_{\min}(M)}\|\dot\theta_i\|_M^2+\frac{\delta_i}{2}C.
\end{ResizedEq}
Suppose $\delta_i\leq\delta_H$. Then $\delta_i\leq1/\LCurvature$ gives $\LCurvature\delta_i^2\leq\delta_i$, hence $\LCurvature\delta_i^2+2\delta_i\leq3\delta_i$; and $\bar U\geq0$ gives $\tfrac12\|\dot\theta_i\|_M^2\leq\bar{\mathcal H}_i^+$. Substituting both into~\prettyref{eq:cito-onestep-raw},
\begin{ResizedEq}
\bar{\mathcal H}_i^+-\bar{\mathcal H}_{\gamma(i)}^+
\leq&\frac{3\delta_i}{2\lambda_{\min}(M)}\|\dot\theta_i\|_M^2+\frac{\delta_i}{2}C\\
\leq&\frac{3\delta_i}{\lambda_{\min}(M)}\bar{\mathcal H}_i^++\frac{\delta_i}{2}C.
\end{ResizedEq}
Rearranging, $\left[1-\tfrac{3\delta_i}{\lambda_{\min}(M)}\right]\bar{\mathcal H}_i^+\leq\bar{\mathcal H}_{\gamma(i)}^++\tfrac{\delta_i}{2}C$. Writing $x\triangleq3\delta_i/\lambda_{\min}(M)$, the bound $\delta_i\leq\delta_H\leq\lambda_{\min}(M)/6$ gives $x\in(0,1/2]$, on which $e^{-2x}\leq1-x$; hence
\begin{ResizedEq}
\label{eq:cito-gronwall-step}
e^{-6\delta_i/\lambda_{\min}(M)}\bar{\mathcal H}_i^+
\leq&\left[1-\frac{3\delta_i}{\lambda_{\min}(M)}\right]\bar{\mathcal H}_i^+\\
\leq&\bar{\mathcal H}_{\gamma(i)}^++\frac{\delta_i}{2}C,
\end{ResizedEq}
which is exactly~\prettyref{eq:cito-ham-check} after multiplying through by $e^{6\delta_i/\lambda_{\min}(M)}$. Thus every node with $\delta_i\leq\delta_H$ passes the check of~\prettyref{ln:cito-ham-check} and is not subdivided by the outer loop.

Only nodes with $\delta_i>\min(\delta^\star,\delta_H)$ can therefore trigger a subdivision, either by the SQP-internal singularity check (which never fires below $\delta^\star$, by~\prettyref{eq:sqp-diagonal-lower}) or by the outer Hamiltonian check (which never fires at or below $\delta_H$, by the previous paragraph). Each subdivision halves the offending interval, and once an interval's length has dropped to $\min(\delta^\star,\delta_H)$ neither it nor any descendant can fail either check again. Letting $G_0\triangleq\mathcal I^0$ be the finite set of optimized right endpoints of the initial mesh over the horizon $[0,H]$, with $i\in G_0$ indexing $(\gamma(i),i)$ and first refinable interval $(0,\min\mathcal I^0)$, each initial interval spawns a binary subdivision tree of depth at most $\left\lceil\log_2\!\left(\delta_i^0/\min(\delta^\star,\delta_H)\right)\right\rceil$, so the total number of subdivisions over the entire run---SQP-internal and outer-loop combined---is bounded by
\begin{ResizedEq}
\label{eq:cito-Nsub}
N_{\mathrm{sub}}\leq\sum_{i\in G_0}\left(2^{\max\left\{0,\left\lceil\log_2\!\left(\frac{\delta_i^0}{\min(\delta^\star,\delta_H)}\right)\right\rceil\right\}}-1\right)<\infty.
\end{ResizedEq}
In particular the outer loop can subdivide only finitely often, so after finitely many outer iterations no node fails the check and~\prettyref{alg:CITO} returns at~\prettyref{ln:cito-return}; each outer iteration is one finite SQP call by Part~1, so the whole run is finite. The proof-level floor~\prettyref{eq:cito-delta-min} was declared before the first solve, and the same binary-tree argument shows that every generated mesh satisfies $\delta_i\geq\delta_{\min}$. This lower bound supplies the positive mesh floor required by~\prettyref{lem:trajectory-bound} and~\ref{lem:jacobian-upper-bound}; no smoothing construction is performed in this general theorem. The final mesh is fixed with this positive lower bound (an unsubdivided initial interval retains its length $\delta_i^0$, while any bisected interval has length exceeding $\min(\delta^\star,\delta_H)/2$), which is exactly the $\delta_{\min}>0$ that~\prettyref{lem:trajectory-bound} and~\ref{lem:jacobian-upper-bound} require. The penalty is bounded across the run: within each SQP call it is bounded by $\varrho_{\max}$ of~\prettyref{lem:sqp-stability}, and since only finitely many meshes arise and $\varrho$ is threaded nondecreasingly, the overall penalty is capped by the maximum of these finitely many per-mesh bounds.

Step 3: the Hamiltonian bound at exit. At return through~\prettyref{ln:cito-return} no node fails the check of~\prettyref{ln:cito-ham-check}, so~\prettyref{eq:cito-ham-check} holds at every active node $i$
\begin{ResizedEq}
\label{eq:cito-exit-check}
&\bar{\mathcal H}_i^+\leq e^{6\delta_i/\lambda_{\min}(M)}\left(\bar{\mathcal H}_{\gamma(i)}^++\frac{\delta_i}{2}C\right)\\
\Rightarrow&
e^{-6\delta_i/\lambda_{\min}(M)}\bar{\mathcal H}_i^+\leq\bar{\mathcal H}_{\gamma(i)}^++\frac{\delta_i}{2}C.
\end{ResizedEq}
(The passage of the check is the statement that~\prettyref{eq:cito-exit-check} holds; the derivation~\prettyref{eq:cito-gronwall-step} above shows this is precisely the discrete Gr\"onwall step of~\prettyref{lem:H-bound} with $\EpsCRDS=\EpsStop$, and that any node with $\delta_i\leq\delta_H$ satisfies it automatically---which is what forced the finiteness in Part~2.) We telescope~\prettyref{eq:cito-exit-check} along the predecessor chain from fixed node $0$ to $i$. In its equivalent multiplied form the step at optimized node $k$ reads $\bar{\mathcal H}_k^+\leq e^{6\delta_k/\lambda_{\min}(M)}(\bar{\mathcal H}_{\gamma(k)}^++\tfrac{\delta_k}{2}C)$. Unrolling this recursion, the growth factors from $k$ through $i$ multiply to $\prod_{\substack{j\in\mathcal I\\k\leq j\leq i}}e^{6\delta_j/\lambda_{\min}(M)}=e^{6(t_i-t_{\gamma(k)})/\lambda_{\min}(M)}$ because $\sum_{\substack{j\in\mathcal I\\k\leq j\leq i}}\delta_j=t_i-t_{\gamma(k)}$, so the chain gives, term by term,
\begin{ResizedEq}
\bar{\mathcal H}_i^+
\leq&e^{6t_i/\lambda_{\min}(M)}\bar{\mathcal H}_0^++\\
&\frac{C}{2}\sum_{\substack{k\in\mathcal I\\k\leq i}}\delta_k\,e^{6(t_i-t_{\gamma(k)})/\lambda_{\min}(M)},
\end{ResizedEq}
where the predecessor-chain identity $t_i=\sum_{\substack{j\in\mathcal I\\j\leq i}}\delta_j$ follows from $t_i=\delta i$ and $\delta_j=\delta(j-\gamma(j))$, and $\bar{\mathcal H}_0^+$ is the fixed time-zero energy. Bounding each factor $e^{6(t_i-t_{\gamma(k)})/\lambda_{\min}(M)}\leq e^{6t_i/\lambda_{\min}(M)}$ (since $t_{\gamma(k)}\geq0$) and using $\sum_{\substack{k\in\mathcal I\\k\leq i}}\delta_k=t_i$,
\begin{ResizedEq}
\label{eq:cito-closedform}
\bar{\mathcal H}_i^+\leq e^{6t_i/\lambda_{\min}(M)}\left[\bar{\mathcal H}_0^++\frac{t_i}{2}C\right]=\LHamBound(t_i),
\end{ResizedEq}
which is the closed form of $\LHamBound$ from the proof of~\prettyref{lem:H-bound}. This proves~\prettyref{eq:cito-ham-bound}. Together with the position-majorant telescoping of~\prettyref{lem:H-bound}, the positive-gap clamp definitions give
\begin{equation}
\bar{\mathcal H}_i^+\leq\LHamBound(H)<\LHamilton,\qquad\|\theta_i\|<\LPosition.
\end{equation}
The inverse clamp implications follow from~\prettyref{eq:clamp-map}: a scalar clamped value below $\LHamilton$ must be on the identity branch, and a radial clamped norm below $\LPosition$ has the original configuration on the radial identity branch. Since every returned node is finite and the inequalities are strict, continuity supplies an open neighborhood of every node and predecessor on which both clamps are identity. On this neighborhood the position and energy clamps are identity, so the clamped potential and control work agree as functions with the fixed unclamped smooth CRDS data. The damping function is the same fixed function throughout the theorem. Thus all residual and Jacobian identities transfer by functional equality on an open neighborhood. In the separate reduced articulated-body instantiation, the additional identity $\mathfrak M(V_R,\LHamilton)\equiv1$ follows from $V_R\leq\widetilde U_R<\LHamilton$ and~\prettyref{eq:global-reduced-damping}, discharging the reduced-model hypothesis used by~\prettyref{rem:cito-ham-discharge}.

For completeness we confirm~\prettyref{eq:cito-closedform} by induction on the $\gamma$-chain based at node $0$, checking that the closed form $\LHamBound(t_i)=e^{6t_i/\lambda_{\min}(M)}[\bar{\mathcal H}_0^++\tfrac{t_i}{2}C]$ satisfies the same recursion~\prettyref{eq:cito-exit-check}. The base case is $\bar{\mathcal H}_0^+=\LHamBound(0)$: by~\prettyref{ass:Slater}, the strict-interior fixed configuration $\theta_0\in\operatorname{dom}U$ gives $U(\theta_0)<\infty$, while strict feasibility makes the penalty contribution $V(\theta_0,\EpsFJ)$ finite; hence $\tilde U(\theta_0,\EpsFJ)<\infty$. With the clamp choices $\LPosition\geq\|\theta_0\|$, $\LHamilton\geq\tilde U(\theta_0,\EpsFJ)$ (implied by the constant-ordering remark),~\prettyref{lem:mollification}~\ref{it:moll-identity} gives $\bar U(\theta_0,\EpsFJ,\LHamilton,\LPosition)=\tilde U(\theta_0,\EpsFJ)$, hence $\bar{\mathcal H}_0^+=\tfrac12\|\dot\theta_0\|_M^2+\tilde U(\theta_0,\EpsFJ)=\LHamBound(0)$, where $\dot\theta_0=(\theta_0-\theta_{\gamma(0)})/\delta_0$. For the inductive step, we verify
\begin{ResizedEq}
\label{eq:cito-induction-recursion}
\LHamBound(t_i)\geq e^{6\delta_i/\lambda_{\min}(M)}\left(\LHamBound(t_{\gamma(i)})+\frac{\delta_i}{2}C\right).
\end{ResizedEq}
Using $t_i=t_{\gamma(i)}+\delta_i$, so that $e^{6\delta_i/\lambda_{\min}(M)}e^{6t_{\gamma(i)}/\lambda_{\min}(M)}=e^{6t_i/\lambda_{\min}(M)}$, the right-hand side of~\prettyref{eq:cito-induction-recursion} equals
\begin{ResizedEq}
e^{6t_i/\lambda_{\min}(M)}\left[\bar{\mathcal H}_0^++\frac{t_{\gamma(i)}}{2}C\right]+e^{6\delta_i/\lambda_{\min}(M)}\frac{\delta_i}{2}C,
\end{ResizedEq}
so the difference $\LHamBound(t_i)$ minus this quantity is
\begin{ResizedEq}
&e^{6t_i/\lambda_{\min}(M)}\frac{t_i-t_{\gamma(i)}}{2}C-e^{6\delta_i/\lambda_{\min}(M)}\frac{\delta_i}{2}C\\
=&\frac{\delta_i}{2}C\left(e^{6t_i/\lambda_{\min}(M)}-e^{6\delta_i/\lambda_{\min}(M)}\right)\geq0,
\end{ResizedEq}
using $t_i-t_{\gamma(i)}=\delta_i$ and $t_i\geq\delta_i$ (so the exponential factors are ordered). This proves~\prettyref{eq:cito-induction-recursion}. Combining the exit check~\prettyref{eq:cito-exit-check}, the induction hypothesis $\bar{\mathcal H}_{\gamma(i)}^+\leq\LHamBound(t_{\gamma(i)})$, and~\prettyref{eq:cito-induction-recursion},
\begin{ResizedEq}
&\bar{\mathcal H}_i^+\leq e^{6\delta_i/\lambda_{\min}(M)}\left(\bar{\mathcal H}_{\gamma(i)}^++\frac{\delta_i}{2}C\right)\\
\leq&e^{6\delta_i/\lambda_{\min}(M)}\left(\LHamBound(t_{\gamma(i)})+\frac{\delta_i}{2}C\right)
\leq\LHamBound(t_i),
\end{ResizedEq}
closing the induction. Hence $\bar{\mathcal H}_i^+\leq\LHamBound(t_i)$ for all $i$, which is~\prettyref{eq:cito-ham-bound}.
\end{proof}

\subsection{Stationarity Transfer}
\label{appen:stationarity-proofs}
Recall the target condition that~\prettyref{thm:stationarity-transfer} certifies (\prettyref{def:eps-stationary}): it is the $\EpsStop$-KKT condition~\prettyref{def:eps-KKT} of the modified \prettyref{pb:CITO-M} reread as an approximate stationarity condition of the original mathematical program with complementarity constraints~\prettyref{pb:CITO}, once the constraint forces are prescribed as configuration-dependent maps and the physical feasibility band is added back. The smooth transfer proof of~\prettyref{thm:stationarity-transfer} assumes the CRDS underlying~\prettyref{pb:CITO} is smooth (\prettyref{def:smooth-CRDS}), so no damping smoothing is needed and $\Damp_{\EpsFri,\LHamilton}^V=\Damp^V$; the nonsmooth transfer of~\prettyref{thm:nonsmooth-transfer} removes this assumption below. Throughout we work with the general CRDS formulation (no separating planes, so $V_R=V$ and $\Damp_R^V=\Damp^V$).

\begin{proof}[Proof of~\prettyref{thm:stationarity-transfer}]
Write $(\StackConf,\StackCtrl)$ for the point returned by~\prettyref{alg:CITO}. By~\prettyref{thm:CITO-convergence}, it admits a finite $\EpsStop$-KKT certificate in the sense of~\prettyref{def:eps-KKT} for~\prettyref{pb:CITO-M}; denote its physics and control multipliers by $\mu_\StackConf,\mu_\StackCtrl$ and its auxiliary controls by $u_i'$. The same theorem supplies the per-step Hamiltonian bound~\prettyref{eq:cito-ham-bound}.

Step~1: multiplier selection and reduction to the succinct system. From the $\EpsStop$-KKT certificate of~\prettyref{pb:CITO-M} inherit the physics and control multipliers $\mu_\StackConf,\mu_\StackCtrl$ (with timestep sub-blocks $\mu_{\theta_i},\mu_{u_i}$) and the auxiliary control $u_i'$. For the remaining multipliers of the full-form system~\prettyref{eq:eps-stationary-system-full} choose
\begin{ResizedEq}
\label{eq:transfer-zero-mult}
\mu_i^e=\mu_i^h=\mu_i^\lambda=\mu_i^c=0\quad\forall i\in\mathcal I.
\end{ResizedEq}
This is admissible because the constraint forces $\lambda_i^e,\lambda_i^i$ are prescribed configuration-dependent maps rather than free primal variables, so no separate stationarity equation with respect to them survives to force these multipliers nonzero. With this choice the full-form system collapses: the complementary-slackness conditions hold trivially, each of $\mu_i^h,\mu_i^\lambda,\mu_i^c$ being zero; the configuration block loses its $\mu_i^e,\mu_i^h,\mu_i^\lambda,\mu_i^c$ terms, leaving $\|\FPP{O}{\StackConf}+[\StatJac]^T\mu_\StackConf\|\leq\EpsStop$; the control block is unchanged; and the three feasibility conditions read componentwise as $h^{i,j}(\theta_i)\geq0$, $\lambda_i^{i,j}(\theta_i)\geq0$, and $\lambda_i^{i,j}(\theta_i)(\EpsFJAll-h^{i,j}(\theta_i))\geq0$, equivalently the single combined line $\min(h^{i,j}(\theta_i),\lambda_i^{i,j}(\theta_i),\lambda_i^{i,j}(\theta_i)(\EpsFJAll-h^{i,j}(\theta_i)))\geq0$. It therefore remains to verify the succinct system
\begin{ResizedEq}
\label{eq:eps-stationary-system}
&\TWO{\EpsStop}{\EpsFJAll}\text{-FJ-stationarity (succinct)}:\\
&\begin{cases}
\|F\|_1\leq\EpsStop\\
\left\|\FPP{O}{\StackConf}+[\StatJac]^T\mu_\StackConf\right\|\leq\EpsStop\\
\left\|\FPP{O}{\StackCtrl}+[\StatJacU]^T\mu_\StackConf+\mu_\StackCtrl\right\|\leq\EpsStop\\
u_i\in\mathcal{D},\ u_i'\in\mathcal{D},\ \|u_i-u_i'\|\leq\EpsStop,\ \mu_{u_i}\in N_{\mathcal{D}}(u_i')\\
|h^{e,j}(\theta_i)|\leq\EpsFJAll\quad\forall j\\
\min\!\big(h^{i,j}(\theta_i),\ \lambda_i^{i,j}(\theta_i),\ \lambda_i^{i,j}(\theta_i)(\EpsFJAll-h^{i,j}(\theta_i))\big)\geq0\quad\forall j,
\end{cases}
\end{ResizedEq}
whose six lines are established in Steps~2-8 below. Every node-indexed line is quantified over optimized nodes $i\in\mathcal I$ only.

Step~2: clamp inactivity and the residual identity. By~\prettyref{eq:cito-ham-bound}, $\bar{\mathcal H}_i^+\leq\LHamBound(t_i)$ at every node. With the positive-gap clamp definitions of the constant-ordering remark in~\prettyref{sec:Hamiltonian}, exactly as in~\prettyref{rem:cito-ham-discharge}, every returned finite node satisfies $\|\theta_i\|<\LPosition$ and $\tilde U(\theta_i,\EpsFJ)\leq\bar{\mathcal H}_i^+\leq\LHamBound(H)<\LHamilton$. This covers every optimized node $i\in\mathcal I$, including the first optimized node $i_1\triangleq\min\mathcal I$: its clamp inactivity and residual identity are derived from the predecessor-chain energy bound based at fixed node $0$. Its stencil predecessors are the fixed strict-interior configuration $\theta_0$ and the fixed pre-zero configuration $\theta_{\gamma(0)}$, which supply the semi-implicit damping and predecessor penalty force. The inverse radial and scalar clamp implications therefore put each node in an open neighborhood on which both clamps are identity. Hence~\prettyref{lem:mollification}~\ref{it:moll-identity} gives $\bar U=\tilde U=U+V$ as a function on that neighborhood, while~\prettyref{lem:controlwork-mollification}~\ref{it:cw-identity} gives $\bar W=W$ (hence $\bar c=c$) there; the first-stage smoothness assumption gives $\Damp_{\EpsFri,\LHamilton}^V=\Damp^V$ on the same neighborhood. On this neighborhood, define the predecessor damping-load map by evaluating the prescribed inequality-force map at the predecessor configuration. By the construction of~\prettyref{eq:penalty-force-generic},
\begin{equation}
\label{eq:transfer-damping-identity}
\Damp^V(\theta_{\gamma(i)},\dot\theta_i)=\Damp\!\left(\theta_{\gamma(i)},\dot\theta_i,\lambda_{\gamma(i)}^i(\theta_{\gamma(i)})\right).
\end{equation}
This predecessor evaluation parameterizes damping, whereas $\lambda_i^i(\theta_i)$ remains the current force used in force balance, feasibility, and complementarity; no equality between the two force values is asserted. There the per-step objective of~\prettyref{pb:PS-CRDS-M} reads
\begin{ResizedEq}
\PerStepObj_i=&\frac{\delta_i^2}{2}\|\ddot{\theta}_i\|_M^2+U(\theta_i)+V(\theta_i,\EpsFJ)\\
&-g^T\theta_i+W(\theta_i,u_i)+\delta_i\Damp^V(\theta_{\gamma(i)},\dot{\theta}_i).
\end{ResizedEq}
Differentiating in $\theta_i$ with $\partial\ddot{\theta}_i/\partial\theta_i=I/\delta_i^2$ and $\partial\dot{\theta}_i/\partial\theta_i=I/\delta_i$ from~\prettyref{eq:discretization}, and using $f=-\nabla_{\theta_i}U$ (\prettyref{def:CRDS}~\ref{it:crds-potential}) and $c=-\nabla_{\theta_i}W$ (\prettyref{eq:control}),
\begin{ResizedEq}
\label{eq:transfer-stat-expand}
\StatRes_i=\nabla_{\theta_i}\PerStepObj_i=&M\ddot{\theta}_i-f(\theta_i)-g-c(\theta_i,u_i)+\\
&\nabla_{\dot{\theta}_i}\Damp^V(\theta_{\gamma(i)},\dot{\theta}_i)+\nabla_{\theta_i}V(\theta_i,\EpsFJ).
\end{ResizedEq}
Using~\prettyref{eq:transfer-damping-identity} and substituting the gradient-matching identity~\prettyref{eq:transfer-gradient-match} for $\nabla_{\theta_i}V$ (established in Step~3 below) turns~\prettyref{eq:transfer-stat-expand} into the force-balance residual~\prettyref{eq:transfer-residual}, so $\StatRes_i\equiv F_i$ as an identity of functions on the open neighborhood. Because the identity holds on an open set, the Jacobians coincide: $\StatJac=\partial\StatRes/\partial\StackConf=\partial F/\partial\StackConf$ and $\StatJacU=\partial\StatRes/\partial\StackCtrl=\partial F/\partial\StackCtrl$.

Step~3: prescribed forces and the gradient-matching identity. Following the penalty-force construction of~\prettyref{lem:eps-FJ} and~\prettyref{eq:penalty-force-generic}, prescribe the maps componentwise from the penalty derivative,
\begin{ResizedEq}
\label{eq:transfer-force-maps}
\lambda_i^{i,j}(\theta_i)\triangleq&-V_s'(h^{i,j}(\theta_i),\EpsFJ)\geq0,\\
\lambda_i^{e,j}(\theta_i)\triangleq&V_s'(\EpsFJ-h^{e,j}(\theta_i),\EpsFJ)-V_s'(\EpsFJ+h^{e,j}(\theta_i),\EpsFJ),
\end{ResizedEq}
the second being the two-sided combination for the equalities. Nonnegativity of $\lambda_i^{i,j}$ follows from $V_s'\leq0$ on its active branch (\prettyref{lem:Vs-flatness}). Differentiating the penalty potential~\prettyref{eq:penalty-potential} through~\prettyref{eq:soft-penalty} by the chain rule gives, term by term,
\begin{ResizedEq}
\nabla_{\theta_i}V(\theta_i,\EpsFJ)=-\FPP{h^e}{\theta_i}(\theta_i)^T\lambda_i^e(\theta_i)-\FPP{h^i}{\theta_i}(\theta_i)^T\lambda_i^i(\theta_i),
\end{ResizedEq}
which is exactly the gradient-matching identity~\prettyref{eq:transfer-gradient-match} posited by~\prettyref{def:eps-stationary}.

Step~4: stationarity blocks. Because $(\StackConf,\StackCtrl)$ is $\EpsStop$-KKT for~\prettyref{pb:CITO-M}, the first two lines of~\prettyref{eq:eps-kkt-system} read $\|\FPP{O}{\StackConf}+[\StatJac]^T\mu_\StackConf\|\leq\EpsStop$ and $\|\FPP{O}{\StackCtrl}+[\StatJacU]^T\mu_\StackConf+\mu_\StackCtrl\|\leq\EpsStop$. By Step~2 the Jacobians $\StatJac,\StatJacU$ appearing there are precisely $\partial F/\partial\StackConf,\partial F/\partial\StackCtrl$, so these two inequalities are lines~2-3 of~\prettyref{eq:eps-stationary-system} verbatim, with the inherited multipliers $\mu_\StackConf,\mu_\StackCtrl$.

Step~5: residual and control lines. The fourth line of~\prettyref{eq:eps-kkt-system} gives $\|\StatRes\|_1\leq\EpsStop$; by the identity $\StatRes_i\equiv F_i$ of Step~2 this is $\|F\|_1\leq\EpsStop$, which is line~1 of~\prettyref{eq:eps-stationary-system}. The third line of~\prettyref{eq:eps-kkt-system}---$u_i\in\mathcal{D}$, $u_i'\in\mathcal{D}$, $\|u_i-u_i'\|\leq\EpsStop$, and $\mu_{u_i}\in N_{\mathcal{D}}(u_i')$---is line~4 of~\prettyref{eq:eps-stationary-system}.

Step~6: equality band. By Step~2 the returned point has $V(\theta_i,\EpsFJ)\leq\tilde U(\theta_i,\EpsFJ)\leq\LHamilton<\infty$, so every equality penalty term $V_s(\EpsFJ\pm h^{e,j}(\theta_i),\EpsFJ)$ in~\prettyref{eq:soft-penalty} is finite. By the support of $V_s$ in~\prettyref{eq:Ps}, finiteness forces both arguments positive, $\EpsFJ\pm h^{e,j}(\theta_i)>0$, i.e. $|h^{e,j}(\theta_i)|<\EpsFJ\leq\EpsFJAll$, which is line~5.

Step~7: inequality band. Finiteness of $V_s(h^{i,j}(\theta_i),\EpsFJ)$ likewise forces $h^{i,j}(\theta_i)>0$, and the prescribed force obeys $\lambda_i^{i,j}(\theta_i)=-V_s'(h^{i,j}(\theta_i),\EpsFJ)\geq0$ by Step~3. On a component with $h^{i,j}(\theta_i)\geq\EpsFJ$ the derivative $V_s'$ vanishes (\prettyref{eq:Ps}), so $\lambda_i^{i,j}(\theta_i)=0$; on a component with $0<h^{i,j}(\theta_i)<\EpsFJ$ both factors of $\lambda_i^{i,j}(\theta_i)(\EpsFJ-h^{i,j}(\theta_i))$ are nonnegative. Hence $\lambda_i^{i,j}(\theta_i)(\EpsFJ-h^{i,j}(\theta_i))\geq0$ in all cases, and relaxing $\EpsFJ\leq\EpsFJAll$,
\begin{ResizedEq}
&\lambda_i^{i,j}(\theta_i)(\EpsFJAll-h^{i,j}(\theta_i))\\
=&\lambda_i^{i,j}(\theta_i)(\EpsFJ-h^{i,j}(\theta_i))+(\EpsFJAll-\EpsFJ)\lambda_i^{i,j}(\theta_i)\geq0,
\end{ResizedEq}
which, with $h^{i,j}(\theta_i)>0$ and $\lambda_i^{i,j}(\theta_i)\geq0$, is line~6. This verifies every line of the succinct system~\prettyref{eq:eps-stationary-system}, which by the multiplier selection of Step~1 certifies the full-form condition~\prettyref{eq:eps-stationary-system-full}. Every line holds for arbitrary $\EpsStop,\EpsFJAll>0$, proving part~\ref{it:transfer-a}.

Step~8: uniform boundedness over the fixed-data tolerance family. Fix $\EpsFJ$, fix the remaining problem and algorithm data and initial data, and use common clamps $\LHamilton,\LPosition$ while $0<\EpsStop\leq\EpsFJ$ varies. On this family $\EpsFJAll=\max\{\EpsFJ,\EpsStop\}=\EpsFJ$. Moreover, the one-step Hamiltonian majorant in~\prettyref{eq:cito-ham-check} is monotone in its only $\EpsStop$-dependent term, $\EpsStop^2$. Its value at the fixed endpoint $\EpsStop=\EpsFJ$ therefore dominates every member of the family. The constant-ordering construction of~\prettyref{sec:Hamiltonian} may consequently be applied once at that endpoint to choose common finite $\LHamilton$ and $\LPosition$ for all $0<\EpsStop\leq\EpsFJ$.

We now make the remaining constants common before choosing $\EpsStop$. With the problem data, algorithm parameters, barrier parameter, and clamps fixed, the singularity and Hamiltonian thresholds $\delta^\star$ and $\delta_H$ in the proof of~\prettyref{thm:CITO-convergence} are independent of $\EpsStop$; set
\begin{ResizedEq}
\delta_c\triangleq\min\{\delta^\star,\delta_H\}.
\end{ResizedEq}
For these common clamps the a priori timestep floor $\delta_{\min}$ of~\prettyref{eq:cito-delta-min} equals $\min\{\min_{i\in G_0}\delta_i^0,\delta_c/2\}>0$. For each initial interval $i\in G_0$, let $\mathcal T_i$ be its exact midpoint-subdivision tree truncated at depth $d_i\triangleq\max\{0,\lceil\log_2(\delta_i^0/\delta_c)\rceil\}$. Let $\mathfrak G$ contain the meshes whose restriction to each initial interval is the leaf frontier of a rooted pruned subtree of $\mathcal T_i$. Every $\mathcal T_i$ is finite and $G_0$ is finite, hence $\mathfrak G$ is finite; this uses the bounded-depth exact midpoint trees, not merely a bound on mesh cardinality. The subdivision estimate~\prettyref{eq:cito-Nsub} shows that every mesh reached by every tolerance-indexed run belongs to $\mathfrak G$: once an interval has length at most $\delta_c$, neither the singularity check nor the Hamiltonian check can subdivide it. Consequently every such mesh has timestep at least $\delta_{\min}$ and at most $|G_0|+N_{\mathrm{sub}}$ active intervals, with both bounds fixed before $\EpsStop$; this floor is a lower bound only and need not be attained, since a bisected interval has length strictly greater than $\delta_c/2$.

Order $\mathfrak G$ by refinement depth and construct compact initialization and all-iterate enclosures inductively. The root initialization is fixed. Suppose the initialization enclosure $K_{\mathrm{init},\mathcal G}$ is compact. Residual continuity gives $R_{\mathrm{init},\mathcal G}\triangleq\sup_{z\in K_{\mathrm{init},\mathcal G}}\|\StatResG(z)\|_2<\infty$. Since the constraint phase makes $\|\StatResG\|_2$ nonincreasing, throughout that phase
\begin{ResizedEq}
R_{\mathrm{cs},\mathcal G}\triangleq\sqrt{\dim\StatResG}\,R_{\mathrm{init},\mathcal G}
\end{ResizedEq}
is a residual one-norm ceiling. Together with~\prettyref{lem:trajectory-bound} and the common timestep floor, it gives a compact constraint-phase enclosure.

If the working band is reached, its terminal point lies in the corresponding compact working-band region. On that fixed mesh, bounded off-diagonal Jacobian blocks and~\prettyref{lem:jacobian-lower-bound} give a positive LICQ modulus; then~\prettyref{lem:sqp-multiplier} at $\LResidualWork$ gives the fixed step bound used to enlarge the working-band region. This compact trial enlargement is precisely the region on which~\prettyref{eq:sqp-stability-consistency} takes finite Jacobian upper bounds and defines the fixed-mesh inner-SQP ceiling $\LResidualG$. Thus the construction is noncircular: the working-band and trial regions are fixed before $\LResidualG$ is defined. The fixed-mesh SQP stability construction and~\prettyref{lem:trajectory-bound} then give a compact inner-SQP enclosure; its finite union with the constraint-phase enclosure is a compact all-iterate enclosure $K_{\mathrm{all},\mathcal G}$.

The exact midpoint initialization map of~\prettyref{alg:subd} is continuous in both states and controls. Its image of $K_{\mathrm{all},\mathcal G}$ is compact, and convexity of the compact set $\mathcal D$ keeps the interpolated controls feasible. Taking the finite union over intervals that can be subdivided gives compact initialization enclosures for all child meshes. This completes the induction without assuming continuity of the complete variable-iteration SQP output map.

Since $\mathfrak G$ is finite, the genuine all-phase residual ceiling
\begin{ResizedEq}
R^\star\triangleq\max_{\mathcal G\in\mathfrak G}\max\{R_{\mathrm{cs},\mathcal G},\LResidualG\}<\infty
\end{ResizedEq}
is fixed before $\EpsStop$. On every compact $K_{\mathrm{all},\mathcal G}$, the fixed-mesh LICQ modulus has a positive lower bound and all other bounds entering~\prettyref{lem:sqp-multiplier} are finite. Taking the minimum of these positive LICQ moduli and the maxima of the upper bounds over $\mathfrak G$, or equivalently taking the maximum of the resulting meshwise multiplier bounds, defines $C_\mu^\star<\infty$ such that
\begin{ResizedEq}
\|\mu_{\theta_i}\|\leq C_\mu^\star
\end{ResizedEq}
for every QP multiplier, in particular every multiplier in a returned certificate. The actual update is the maximum update $\varrho^+=\max\{\varrho^-,\max_i\|\mu_{\theta_i}\|+\bar\varrho\}$, not an additive increment. Therefore
\begin{ResizedEq}
\label{eq:transfer-common-penalty}
\varrho\leq\varrho_{\max}\triangleq\max\{\varrho^0,C_\mu^\star+\bar\varrho\}.
\end{ResizedEq}
Because the updated value is threaded, rather than reset or incremented, across internal mesh episodes and outer calls to~\prettyref{alg:CITO}, the same cap holds regardless of their number.

Applying~\prettyref{lem:trajectory-bound} with $R^\star$, the common timestep floor $\delta_{\min}$, the fixed horizon, and the common CRDS data gives one finite trajectory bound $\|\theta_i\|\leq\LPosBound(H)$; control iterates remain in the common compact set $\mathcal D$. Finally, because $\EpsFJ$ is fixed, the common level set $V(\theta_i,\EpsFJ)\leq\LHamilton$ stays a positive distance from each penalty singularity: $V_s(x,\EpsFJ)\to\infty$ as $x\to0^+$, and the penalty-derivative bound of~\prettyref{lem:Vs-flatness} is finite on the resulting compact argument band. Thus~\prettyref{eq:transfer-force-maps} gives a finite common force ceiling $\bar\Lambda$ for all current evaluations $\lambda_i^e(\theta_i),\lambda_i^i(\theta_i)$ and all predecessor damping-load evaluations $\lambda_{\gamma(i)}^i(\theta_{\gamma(i)})$. For optimized predecessors this follows from the same compact penalty sublevel; for the first optimized node it follows from the single fixed strict-interior configuration $\theta_0$. The trajectory, physics multipliers, penalty, and prescribed forces are therefore bounded independently of $\EpsStop$ on $0<\EpsStop\leq\EpsFJ$, proving part~\ref{it:transfer-b}.
\end{proof}

We now assemble the conditional certificate-normalization machinery that the nonsmooth transfer uses in place of any smoothing-uniform multiplier bound. Fix $\EpsFri>0$ and suppose that $(\StackConf,\StackCtrl)$ admits a finite $\EpsStop$-KKT certificate of the corresponding smooth modified \prettyref{pb:CITO-M} in the sense of~\prettyref{def:eps-KKT}. Denote its physics and control multipliers by $\mu_\StackConf,\mu_\StackCtrl$ (with timestep sub-blocks $\mu_{\theta_i},\mu_{u_i}$) and its auxiliary controls by $u_i'$. Separately, for the later homogeneous transfer select the equality-band and inequality multipliers to be zero; these multipliers are not part of~\prettyref{def:eps-KKT}. All second derivatives below exist because~\prettyref{def:S-D} inherits~\prettyref{def:smooth-CRDS}. Assume in addition the level-set and clamp-inactivity properties established for the eventual algorithmic output by~\prettyref{eq:cito-ham-bound} and~\prettyref{rem:cito-ham-discharge}, so every $\theta_i$ lies in the norm-bounded region of~\prettyref{lem:trajectory-bound}. On the accepted mesh, every timestep is bounded below and the cardinality is bounded above,
\begin{ResizedEq}
\label{eq:norm-mesh}
\delta_i\geq\delta_{\min}>0,\qquad |\mathcal{I}|\leq|G_0|+N_{\mathrm{sub}},
\end{ResizedEq}
with $\delta_{\min}$ of~\prettyref{eq:cito-delta-min} and $N_{\mathrm{sub}}$ of~\prettyref{eq:cito-Nsub}, both built from the thresholds $\delta^\star,\delta_H$. These thresholds depend only on the fixed nonsmoothing problem data, initial mesh, clamps, and algorithmic thresholds: the diagonal damping curvature $K_i\succeq0$ only helps the singularity check and cannot force additional subdivision, while the Hamiltonian threshold uses damping only through convexity and minimality at zero velocity. Hence, for this common nonsmoothing data, both bounds are independent of $\EpsFri$ and $\EpsStop$.

For a direction $\xi$ in the timestep-$i$ multiplier space and each term $m\in\mathcal{T}$, the smooth graph $\{q=\nabla_{\dot{\theta}}\Damp_{\EpsFri,\LHamilton}^{V,m}(\theta_{\gamma(i)},\dot{\theta}_i)\}$ is the zero set of $q-\nabla_{\dot{\theta}}\Damp_{\EpsFri,\LHamilton}^{V,m}$; multiplying the transposed Jacobian of this defining map by $-\xi$ and using $\nabla^2_{\dot{\theta}\dot{\theta}}\Damp_{\EpsFri,\LHamilton}^{V,m}=\delta_i K_i^m$ (from $\dot{\theta}_i=(\theta_i-\theta_{\gamma(i)})/\delta_i$ and~\prettyref{eq:S-D-compatible-blocks}) gives the smooth per-term graph normal generated by $\xi$,
\begin{ResizedEq}
\label{eq:ns-smooth-normal}
n_i^{m,\EpsFri}(\xi)\triangleq((\Gamma_i^m)^T\xi,\ \delta_i K_i^m\xi,\ -\xi).
\end{ResizedEq}
Its pullback under $Q_i$ by~\prettyref{eq:transfer-Qmap} is
\begin{ResizedEq}
\label{eq:ns-pullback}
[\nabla Q_i]^T n_i^{m,\EpsFri}(\xi):
\begin{cases}
K_i^m\xi&\theta_i\text{-block}\\
(\Gamma_i^m)^T\xi-K_i^m\xi&\theta_{\gamma(i)}\text{-block}\\
-\xi&q_i^m\text{-block},
\end{cases}
\end{ResizedEq}
and, writing $A_{ji}$ for the kinetic part of the predecessor Jacobian block, the chain rule through $\partial\dot{\theta}_j/\partial\theta_j=I/\delta_j$, $\partial\dot{\theta}_j/\partial\theta_{\gamma(j)}=-I/\delta_j$, together with $\nabla^2_{\dot{\theta}\dot{\theta}}\Damp_{\EpsFri,\LHamilton}^V=\delta_j K_j$, gives the damping split
\begin{ResizedEq}
\label{eq:um-splitGK}
\delta_j\FPPTT{\Damp_{\EpsFri,\LHamilton}^V}{\theta_j}{\theta_{\gamma(j)}}=\Gamma_j-K_j\ \Rightarrow\ \StatJac_{ji}=A_{ji}+\Gamma_j-K_j
\end{ResizedEq}
for a successor $\gamma(j)=i$ (\prettyref{eq:sqp-diagonal-expand}). Consequently, summing~\prettyref{eq:ns-pullback} over $m$ with $\xi=\mu_{\theta_i}$ and over $i$ reproduces the damping content of the adjoint $[\StatJac]^T\mu_\StackConf$ exactly: the $\theta_i$-block $\sum_m K_i^m\mu_{\theta_i}=K_i\mu_{\theta_i}$ matches the diagonal damping term of $\StatJac_{ii}$ (\prettyref{eq:sqp-diagonal-expand}), and the $\theta_{\gamma(j)}$-blocks $\sum_m((\Gamma_j^m)^T-K_j^m)\mu_{\theta_j}=(\Gamma_j-K_j)^T\mu_{\theta_j}$ furnish the off-diagonal damping term of $\StatJac_{ji}$ (second-neighbor successors carry only kinetic terms, absent from $Q_j$'s $\theta_i$-argument).
\begin{definition}
\label{def:S}
For the finite certificate multipliers $\mu_\StackConf,\mu_\StackCtrl$ fixed above, form the per-term smooth graph normals $n_i^{m,\EpsFri}(\mu_{\theta_i})$ of~\prettyref{eq:ns-smooth-normal}, all finite for fixed $\EpsFri$ since the smoothed problem is finite-dimensional and smooth, and set
\begin{ResizedEq}
\label{eq:S-scale}
S_{\EpsFri}^2\triangleq1+\|\mu_\StackConf\|^2+\|\mu_\StackCtrl\|^2+\sum_{i\in\mathcal{I}}\sum_{m\in\mathcal{T}}\|n_i^{m,\EpsFri}(\mu_{\theta_i})\|^2.
\end{ResizedEq}
Define the normalized multipliers and normals
\begin{ResizedEq}
\label{eq:normalized-mult}
&\bar\mu_0\triangleq S_{\EpsFri}^{-1},\quad \bar\mu_\StackConf\triangleq S_{\EpsFri}^{-1}\mu_\StackConf,\quad \bar\mu_\StackCtrl\triangleq S_{\EpsFri}^{-1}\mu_\StackCtrl,\\
&\bar n_i^{m,\EpsFri}\triangleq S_{\EpsFri}^{-1}n_i^{m,\EpsFri}(\mu_{\theta_i}),
\end{ResizedEq}
with timestep sub-blocks $\bar\mu_{\theta_i}=S_{\EpsFri}^{-1}\mu_{\theta_i}$, $\bar\mu_{u_i}=S_{\EpsFri}^{-1}\mu_{u_i}$.
\end{definition}
\begin{lemma}[Certificate normalization to the unit sphere]
\label{lem:normalized-certificate}
Suppose $(\StackConf,\StackCtrl)$ admits the finite $\EpsStop$-KKT certificate of~\prettyref{def:eps-KKT} used in~\prettyref{def:S}, with multipliers $\mu_\StackConf,\mu_\StackCtrl$ and auxiliary controls $u_i'$. Then the normalized quantities satisfy
\begin{ResizedEq}
\label{eq:norm-unit-sphere}
&\bar\mu_0^2+\|\bar\mu_\StackConf\|^2+\|\bar\mu_\StackCtrl\|^2+\sum_{i,m}\|\bar n_i^{m,\EpsFri}\|^2=1\\
&0<\bar\mu_0\leq1,
\end{ResizedEq}
and $\bar n_i^{m,\EpsFri}=n_i^{m,\EpsFri}(\bar\mu_{\theta_i})$. Moreover, after replacing the objective gradient by $\bar\mu_0\nabla O$, the configuration and control stationarity residuals of the $\EpsStop$-KKT certificate are bounded by $\EpsStop/S_{\EpsFri}\leq\EpsStop$
\begin{ResizedEq}
\label{eq:norm-residuals}
&\left\|\bar\mu_0\FPP{O}{\StackConf}+[\StatJac]^T\bar\mu_\StackConf\right\|\leq\frac{\EpsStop}{S_{\EpsFri}}\\
&\left\|\bar\mu_0\FPP{O}{\StackCtrl}+[\StatJacU]^T\bar\mu_\StackConf+\bar\mu_\StackCtrl\right\|\leq\frac{\EpsStop}{S_{\EpsFri}}.
\end{ResizedEq}
The scaled control multiplier stays in the cone $\bar\mu_{u_i}\in N_{\mathcal{D}}(u_i')$.
\end{lemma}
\begin{proof}[Proof of~\prettyref{lem:normalized-certificate}]
\prettyref{eq:norm-unit-sphere} is immediate from~\prettyref{eq:S-scale} and~\ref{eq:normalized-mult}; $S_{\EpsFri}<\infty$ gives $\bar\mu_0>0$ and $S_{\EpsFri}\geq1$ gives $\bar\mu_0\leq1$. The identity $\bar n_i^{m,\EpsFri}=n_i^{m,\EpsFri}(\bar\mu_{\theta_i})$ follows because each block of~\prettyref{eq:ns-smooth-normal} is linear in $\xi$, so with $\bar\mu_{\theta_i}=S_{\EpsFri}^{-1}\mu_{\theta_i}$ we have $\bar n_i^{m,\EpsFri}=S_{\EpsFri}^{-1}n_i^{m,\EpsFri}(\mu_{\theta_i})=n_i^{m,\EpsFri}(\bar\mu_{\theta_i})$. By the assumed certificate, the first two lines of~\prettyref{eq:eps-kkt-system} give the stationarity inequalities. Dividing those inequalities---$\|\FPP{O}{\StackConf}+[\StatJac]^T\mu_\StackConf\|\leq\EpsStop$ and $\|\FPP{O}{\StackCtrl}+[\StatJacU]^T\mu_\StackConf+\mu_\StackCtrl\|\leq\EpsStop$---by $S_{\EpsFri}$ scales the objective coefficient to $\bar\mu_0=1/S_{\EpsFri}$, every multiplier to its bar-version, and the residual from $\EpsStop$ to $\EpsStop/S_{\EpsFri}$, giving~\prettyref{eq:norm-residuals}. Finally $N_{\mathcal{D}}(u_i')$ is a cone and $S_{\EpsFri}^{-1}>0$.
\end{proof}
\begin{lemma}[Energy bounds on the normalized multipliers]
\label{lem:energy}
Let the normalized quantities satisfy the unit-sphere identity~\prettyref{eq:norm-unit-sphere} and the normalized-normal identity $\bar n_i^{m,\EpsFri}=n_i^{m,\EpsFri}(\bar\mu_{\theta_i})$. Then, for every $i\in\mathcal{I}$ and $m\in\mathcal{T}$,
\begin{ResizedEq}
&\|\bar\mu_{\theta_i}\|\leq1\\
&\|K_i^m\bar\mu_{\theta_i}\|\leq\frac1{\delta_i}\leq\frac1{\delta_{\min}}\\
&\bar\mu_{\theta_i}^T K_i\bar\mu_{\theta_i}\leq C_E^\star=\frac{|\mathcal{T}|}{\delta_{\min}}.
\end{ResizedEq}
All three constants are independent of $\EpsFri,\EpsStop$ and of the unnormalized multipliers.
\end{lemma}
\begin{proof}[Proof of~\prettyref{lem:energy}]
The timestep-$i$ block norm is bounded by the stacked norm $\|\bar\mu_\StackConf\|\leq1$ (\prettyref{eq:norm-unit-sphere}), giving $\|\bar\mu_{\theta_i}\|\leq1$. By~\prettyref{eq:norm-unit-sphere} and $\bar n_i^{m,\EpsFri}=n_i^{m,\EpsFri}(\bar\mu_{\theta_i})$, the velocity block $\delta_i K_i^m\bar\mu_{\theta_i}$ obeys $\delta_i\|K_i^m\bar\mu_{\theta_i}\|\leq\|\bar n_i^{m,\EpsFri}\|\leq1$; divide by $\delta_i$ and use~\prettyref{eq:norm-mesh}. Finally, with $K_i=\sum_m K_i^m$, $K_i^m\succeq0$, and Cauchy-Schwarz,
\begin{ResizedEq}
\bar\mu_{\theta_i}^T K_i\bar\mu_{\theta_i}
=&\sum_{m\in\mathcal{T}}\bar\mu_{\theta_i}^T K_i^m\bar\mu_{\theta_i}\\
\leq&\sum_{m\in\mathcal{T}}\|\bar\mu_{\theta_i}\|\,\|K_i^m\bar\mu_{\theta_i}\|
\leq\frac{|\mathcal{T}|}{\delta_{\min}}.
\end{ResizedEq}
\end{proof}
\begin{corollary}
\label{cor:apply-8b}
The normalized direction $\xi=\bar\mu_{\theta_i}$ satisfies the hypotheses of~\prettyref{def:S-D-compatible}, namely $\|\xi\|\leq1$ and $\xi^T K_i\xi\leq C_E^\star=|\mathcal{T}|/\delta_{\min}$, with constants common to all $i\in\mathcal{I}$, all $m\in\mathcal{T}$, and all smoothing parameters.
\end{corollary}
\begin{proof}[Proof of~\prettyref{cor:apply-8b}]
Immediate from~\prettyref{lem:energy}.
\end{proof}
The final two lemmas restore exact dual nontriviality after the smooth graph normals are replaced by nearby Clarke normals. Here $G_0=\mathcal I^0$ consists of the initial optimized interval right endpoints and $N_{\mathrm{sub}}$ counts inserted optimized nodes. Write $N^\star\triangleq|\mathcal{T}|(|G_0|+N_{\mathrm{sub}})$, a fixed bound on the number of optimized-node--term pairs.
\begin{lemma}
\label{lem:nonzero}
Let $\bar Z\triangleq(\bar\mu_0,\bar\mu_\StackConf,\bar\mu_\StackCtrl,(\bar n_i^{m,\EpsFri}))$ satisfy $\|\bar Z\|=1$. Suppose that for each $i,m$ a Clarke normal $n_i^{m\sharp}$ satisfies $\|n_i^{m\sharp}-\bar n_i^{m,\EpsFri}\|\leq\omega_N(\EpsFri)$, and set $Z^\sharp\triangleq(\bar\mu_0,\bar\mu_\StackConf,\bar\mu_\StackCtrl,(n_i^{m\sharp}))$. Then $\|Z^\sharp\|\geq1-\sqrt{N^\star}\,\omega_N(\EpsFri)$; in particular if $\sqrt{N^\star}\,\omega_N(\EpsFri)\leq1/2$ then $\|Z^\sharp\|\geq1/2$.
\end{lemma}
\begin{proof}[Proof of~\prettyref{lem:nonzero}]
By hypothesis $\|\bar Z\|=1$. The stacked difference obeys $\|Z^\sharp-\bar Z\|\leq(\sum_{i,m}\omega_N(\EpsFri)^2)^{1/2}\leq\sqrt{N^\star}\,\omega_N(\EpsFri)$ by~\prettyref{eq:norm-mesh} and $|\mathcal{T}|$ terms; the reverse triangle inequality gives the claim.
\end{proof}
\begin{lemma}
\label{lem:final-beta}
Let $Z^\sharp$ be the nonzero replaced homogeneous tuple of~\prettyref{lem:nonzero}, and suppose the band multipliers vanish, $\mu_i^e=\mu_i^h=\mu_i^\lambda=\mu_i^c=0$. Under $\sqrt{N^\star}\,\omega_N(\EpsFri)\leq1/2$, set $\beta\triangleq\|Z^\sharp\|^{-1}\leq2$ and multiply the objective multiplier, all physics, control, and band multipliers, and all graph normals by $\beta$. Then the nontriviality normalization~\prettyref{eq:eps-stationary-NS-normalization} holds exactly, every cone and sign condition is preserved, and every dual stationarity residual is multiplied by at most two.
\end{lemma}
\begin{proof}[Proof of~\prettyref{lem:final-beta}]
Clarke normal cones and the convex control cone $N_{\mathcal{D}}(u_i')$ are cones; the nonnegativity of $\mu_0,\mu_i^h,\mu_i^\lambda,\mu_i^c$ and all complementary-slackness equalities are preserved by multiplication by $\beta>0$. By~\prettyref{lem:nonzero}, $\beta=\|Z^\sharp\|^{-1}\leq2$, giving unit norm and the factor-two residual bound.
\end{proof}

\begin{proof}[Proof of~\prettyref{thm:nonsmooth-transfer}]
We follow the smooth proof of~\prettyref{thm:stationarity-transfer}, but normalize the fixed-smoothing certificate before transferring the damping graph normals; the control, equality, and inequality lines are unchanged. By~\prettyref{thm:CITO-convergence}, the returned point admits a finite $\EpsStop$-KKT certificate in the sense of~\prettyref{def:eps-KKT}; denote its physics and control multipliers by $\mu_\StackConf,\mu_\StackCtrl$ and its auxiliary controls by $u_i'$.

Step~1: multiplier selection and reduction to the succinct system. From the $\EpsStop$-KKT certificate of~\prettyref{pb:CITO-M} inherit $\mu_\StackConf,\mu_\StackCtrl$ (with timestep sub-blocks $\mu_{\theta_i},\mu_{u_i}$) and $u_i'$, and select $\mu_i^e=\mu_i^h=\mu_i^\lambda=\mu_i^c=0$ as in~\prettyref{eq:transfer-zero-mult}, admissible for the same reason as in Step~1 of the proof of~\prettyref{thm:stationarity-transfer}: the constraint forces $\lambda_i^e,\lambda_i^i$ are prescribed configuration-dependent maps, so no separate stationarity equation forces these multipliers nonzero. It remains to verify the succinct nonsmooth system
\begin{ResizedEq}
\label{eq:eps-stationary-NS-system-succinct}
&\text{C-}\TWO{\EpsStop}{\EpsFJAll}\text{-FJ-stationarity (succinct)}:\\
&\begin{cases}
\|\widehat F(\StackConf,\StackCtrl,(q^{m\sharp}))\|_1\leq\EpsStop,\quad\|\dot{\theta}_i^{m\sharp}-\dot{\theta}_i\|\leq\EpsStop\ \ \forall m\\
\left\|\nabla_{(\StackConf,q)}\big[\mu_0 O+\sum_i\langle\mu_{\theta_i},\widehat F_i\rangle\big]+\sum_i\sum_{m}[\nabla Q_i]^T n_i^{m\sharp}\right\|\leq\EpsStop\\
\left\|\mu_0\FPP{O}{\StackCtrl}+[\StatJacU]^T\mu_\StackConf+\mu_\StackCtrl\right\|\leq\EpsStop\\
u_i\in\mathcal{D},\ u_i'\in\mathcal{D},\ \|u_i-u_i'\|\leq\EpsStop,\ \mu_{u_i}\in N_{\mathcal{D}}(u_i')\\
|h^{e,j}(\theta_i)|\leq\EpsFJAll\quad\forall j\\
\min\!\big(h^{i,j}(\theta_i),\ \lambda_i^{i,j}(\theta_i),\ \lambda_i^{i,j}(\theta_i)(\EpsFJAll-h^{i,j}(\theta_i))\big)\geq0\quad\forall j,
\end{cases}
\end{ResizedEq}
together with the nontriviality normalization~\prettyref{eq:eps-stationary-NS-normalization}, at the multiplier tuple constructed below; every node-indexed line and sum is over optimized $i\in\mathcal I$; the two stationarity lines and the dynamics line will hold at the decorated tolerance $\EpsStopNS$ of~\prettyref{eq:ns-transfer-band}, the remaining lines at $(\EpsStop,\EpsFJAll)$.

Step~2: clamp inactivity and the lifted residual identity. Exactly as in Step~2 of the proof of~\prettyref{thm:stationarity-transfer}, every returned $\theta_i$ lies in the interior of the clamp-inactive region, where $\bar U=\tilde U=U+V$, $\bar W=W$, and the gradient-matching identity~\prettyref{eq:transfer-gradient-match} holds. With the smoothed damping in force the per-step residual~\prettyref{eq:transfer-stat-expand} reads
\begin{ResizedEq}
&\StatRes_i=\widehat F_i(\StackConf,\StackCtrl,(q_i^{m,\EpsFri})_{m\in\mathcal{T}})\\
&q_i^{m,\EpsFri}\triangleq\nabla_{\dot{\theta}_i}\Damp_{\EpsFri,\LHamilton}^{V,m}(\theta_{\gamma(i)},\dot{\theta}_i),
\end{ResizedEq}
with $\sum_{m\in\mathcal{T}}q_i^{m,\EpsFri}=\nabla_{\dot{\theta}_i}\Damp_{\EpsFri,\LHamilton}^V(\theta_{\gamma(i)},\dot{\theta}_i)$, i.e.\ the lifted residual~\prettyref{eq:transfer-residual-lifted} at the smoothed per-term forces.

Step~3: normalization and exact lifting of the configuration line. Apply~\prettyref{lem:normalized-certificate} to this finite certificate and normalize it by the scale $S_{\EpsFri}$ of~\prettyref{def:S}, producing $\bar\mu_0,\bar\mu_\StackConf,\bar\mu_\StackCtrl$ and the normalized smooth graph normals $\bar n_i^{m,\EpsFri}=n_i^{m,\EpsFri}(\bar\mu_{\theta_i})$ (\prettyref{lem:normalized-certificate}). By the exact lifting~\prettyref{eq:ns-pullback} and~\ref{eq:um-splitGK}, summing the pullbacks of $n_i^{m,\EpsFri}(\bar\mu_{\theta_i})$ over $i,m$ reproduces the damping content of $[\StatJac]^T\bar\mu_\StackConf$; the kinetic, potential, control-force, and constraint-force terms of $[\StatJac]^T\bar\mu_\StackConf$ are the $q$-independent part of $\nabla_{(\StackConf,q)}\sum_i\langle\bar\mu_{\theta_i},\widehat F_i\rangle$ (since $\widehat F_i$ agrees with $F_i$ except in the damping term), and $\FPP{\widehat F_i}{q_i^m}=I$ makes each $q_i^m$-block $\bar\mu_{\theta_i}$, cancelled by the $q_i^m$-block $-\bar\mu_{\theta_i}$ of $[\nabla Q_i]^T n_i^{m,\EpsFri}(\bar\mu_{\theta_i})$. Hence the normalized configuration line~\prettyref{eq:norm-residuals} is the exact rewriting
\begin{ResizedEq}
\label{eq:ns-lifted-kkt}
&\Big\|\nabla_{(\StackConf,q)}\big[\bar\mu_0 O+\textstyle\sum_i\langle\bar\mu_{\theta_i},\widehat F_i\rangle\big]+\textstyle\sum_i\sum_{m}[\nabla Q_i]^T n_i^{m,\EpsFri}(\bar\mu_{\theta_i})\Big\|\\
\leq&\frac{\EpsStop}{S_{\EpsFri}}\leq\EpsStop,
\end{ResizedEq}
evaluated at $q=(q^{m,\EpsFri})$; no approximation has been made.

Step~4: fixed constants and the nearby graph points. By~\prettyref{lem:energy} and~\prettyref{cor:apply-8b}, $\xi=\bar\mu_{\theta_i}$ obeys $\|\xi\|\leq1$ and $\xi^T K_i\xi\leq C_E^\star=|\mathcal{T}|/\delta_{\min}$, exactly the hypotheses of~\prettyref{def:S-D-compatible}. That condition supplies, for each $m\in\mathcal{T}$, a~\prettyref{def:S-D} nearby point $(\theta_{\gamma(i)},\dot{\theta}_i^{m\dagger},q_i^{m\dagger})$ with $q_i^{m\dagger}\in\partial_{\dot{\theta}_i}\Damp^{V,m}(\theta_{\gamma(i)},\dot{\theta}_i^{m\dagger})$ and a Clarke normal $n_i^{m\sharp}\in N^C_{\gph\partial_{\dot{\theta}}\Damp^{V,m}}(\theta_{\gamma(i)},\dot{\theta}_i^{m\dagger},q_i^{m\dagger})$ with
\begin{ResizedEq}
\label{eq:ns-normal-close}
\|n_i^{m\sharp}-n_i^{m,\EpsFri}(\bar\mu_{\theta_i})\|=\|n_i^{m\sharp}-\bar n_i^{m,\EpsFri}\|\leq\omega_N(\EpsFri).
\end{ResizedEq}
Set $\dot{\theta}_i^{m\sharp}=\dot{\theta}_i^{m\dagger}$, $q_i^{m\sharp}=q_i^{m\dagger}$. Since $\nabla_{\theta_i}(\delta_i\Damp_{\EpsFri,\LHamilton}^{V,m})=q_i^{m,\EpsFri}$ and $\partial_{\theta_i}(\delta_i\Damp^{V,m})(\theta_{\gamma(i)},\dot{\theta}_i^{m\dagger})=\partial_{\dot{\theta}_i}\Damp^{V,m}(\theta_{\gamma(i)},\dot{\theta}_i^{m\dagger})$,~\prettyref{def:S-D}\ref{it:sd-grad-close} and~\ref{it:sd-vel-close} give
\begin{ResizedEq}
\label{eq:ns-force-close}
&\|q_i^{m\sharp}-q_i^{m,\EpsFri}\|\leq\delta_i\EpsFri\\
&\|\dot{\theta}_i^{m\sharp}-\dot{\theta}_i\|\leq\EpsFri.
\end{ResizedEq}

Step~5: dynamics line. Since $\StatRes_i=\widehat F_i(\StackConf,\StackCtrl,(q_i^{m,\EpsFri}))$ with each $q_i^m$ entering with identity coefficient, the residual line $\|\StatRes\|_1\leq\EpsStop$ of~\prettyref{def:eps-KKT} and~\prettyref{eq:ns-force-close} give
\begin{ResizedEq}
\label{eq:ns-dyn-band}
&\|\widehat F(\StackConf,\StackCtrl,(q^{m\sharp}))\|_1\\
\leq&\|\StatRes\|_1+\sum_{i\in\mathcal{I}}\sum_{m\in\mathcal{T}}\|q_i^{m\sharp}-q_i^{m,\EpsFri}\|_1\\
\leq&\EpsStop+C_F N^\star\EpsFri,
\end{ResizedEq}
where $\|q_i^{m\sharp}-q_i^{m,\EpsFri}\|_1\leq\sqrt{\dim\theta_i}\,\delta_i\EpsFri$, $N^\star=|\mathcal{T}|(|G_0|+N_{\mathrm{sub}})$ (governed by the subdivision threshold $\delta_c$, not by $H/\delta_{\min}$, so it does not inflate as the floor $\delta_{\min}$ shrinks), and $C_F\triangleq\max\{1,\sqrt{\max_i\dim\theta_i}\,\delta_{\max}\}$ is $\EpsFri,\EpsStop$-independent by~\prettyref{eq:norm-mesh}. The nearby-velocity residual $\|\dot{\theta}_i^{m\sharp}-\dot{\theta}_i\|\leq\EpsFri$ holds by~\prettyref{eq:ns-force-close}. Both are primal quantities, unaffected by the dual normalization.

Step~6: configuration C-stationarity line and final renormalization. Replacing $n_i^{m,\EpsFri}(\bar\mu_{\theta_i})$ by $n_i^{m\sharp}$ in~\prettyref{eq:ns-lifted-kkt} changes it by $\sum_i\sum_m[\nabla Q_i]^T(n_i^{m\sharp}-\bar n_i^{m,\EpsFri})$. Each block $\theta_i$ appears in the pullbacks of $Q_i$ and of the maps $Q_j$ with $\gamma(j)=i$, over all $m\in\mathcal{T}$---a finite set by the banded mesh structure~\prettyref{eq:sqp-jacobian-structure} and $\delta_i\geq\delta_{\min}>0$---so the assembled pullback operator obeys
\begin{ResizedEq}
&\Big\|\sum_i\sum_m[\nabla Q_i]^T z_i^m\Big\|\leq C_N\Big(\sum_i\sum_m\|z_i^m\|^2\Big)^{1/2}\\
&C_N=C_N(\delta_{\min},\text{incidence}),
\end{ResizedEq}
uniformly in $\EpsFri,\EpsStop$. With $z_i^m=n_i^{m\sharp}-\bar n_i^{m,\EpsFri}$ and~\prettyref{eq:ns-normal-close},
\begin{ResizedEq}
\label{eq:ns-stat-band}
&\Big\|\nabla_{(\StackConf,q)}\big[\bar\mu_0 O+\textstyle\sum_i\langle\bar\mu_{\theta_i},\widehat F_i\rangle\big]+\textstyle\sum_i\sum_m[\nabla Q_i]^T n_i^{m\sharp}\Big\|\\
\leq&\EpsStop+C_N\sqrt{N^\star}\,\omega_N(\EpsFri),
\end{ResizedEq}
each $q_i^m$-block reducing only to $\|\bar\mu_{\theta_i}+[n_i^{m\sharp}]_q\|\leq\omega_N(\EpsFri)$, with $[n_i^{m\sharp}]_q$ the damping-force block of $n_i^{m\sharp}$, counted in the same estimate. By~\prettyref{lem:nonzero} the replaced dual vector $Z^\sharp=(\bar\mu_0,\bar\mu_\StackConf,\bar\mu_\StackCtrl,(n_i^{m\sharp}))$ has $\|Z^\sharp\|\geq1/2$ under the hypothesis $\sqrt{N^\star}\,\omega_N(\EpsFri)\leq1/2$; applying the conic normalization $\beta=\|Z^\sharp\|^{-1}\leq2$ of~\prettyref{lem:final-beta}---multiplying $\bar\mu_0,\bar\mu_\StackConf,\bar\mu_\StackCtrl$, the zero band multipliers, and all graph normals by $\beta$---restores~\prettyref{eq:eps-stationary-NS-normalization} exactly and multiplies the dual residuals~\prettyref{eq:norm-residuals} and~\ref{eq:ns-stat-band} by at most two. Write the final certificate $\mu_0\triangleq\beta\bar\mu_0$, $\mu_\StackConf\triangleq\beta\bar\mu_\StackConf$, $\mu_\StackCtrl\triangleq\beta\bar\mu_\StackCtrl$, $n_i^{m\sharp}\mapsto\beta n_i^{m\sharp}$; the configuration line then holds at $2\EpsStop+2C_N\sqrt{N^\star}\,\omega_N(\EpsFri)$.

Step~7: control line and bands. The damping graph adds no control dependence, so the control-stationarity line $\|\mu_0\FPP{O}{\StackCtrl}+[\StatJacU]^T\mu_\StackConf+\mu_\StackCtrl\|$ (at most $2\EpsStop/S_{\EpsFri}\leq2\EpsStop$ by~\prettyref{eq:norm-residuals} and $\beta\leq2$), the two-point control condition, and the equality, inequality, and complementarity bands are established verbatim as in the proof of~\prettyref{thm:stationarity-transfer} (Steps~4-7 there), holding at tolerance $\EpsStop$ and band $\EpsFJAll$; the scaled control multiplier keeps $\mu_{u_i}\in N_{\mathcal{D}}(u_i')$ by~\prettyref{lem:normalized-certificate} and~\ref{lem:final-beta}. Collecting the bands, each line of~\prettyref{eq:eps-stationary-NS-system-succinct} holds at its own tolerance under the single scalar $\EpsStopNS$: the configuration line at $2\EpsStop+2C_N\sqrt{N^\star}\,\omega_N(\EpsFri)$, which dominates the control line's $2\EpsStop$; and the dynamics line~\prettyref{eq:ns-dyn-band} at $\EpsStop+C_F N^\star\EpsFri$, which dominates the nearby-velocity residual $\EpsFri$ since $C_F N^\star\geq1$. These lines are decoupled---no aggregate constraint sums their residuals---so the smallest common tolerance is the maximum of the two bands, namely $\EpsStopNS$ of~\prettyref{eq:ns-transfer-band}. Every line then holds at $\TWO{\EpsStopNS}{\EpsFJAll}$ together with the normalization~\prettyref{eq:eps-stationary-NS-normalization}, which by the multiplier selection of Step~1 certifies~\prettyref{def:eps-stationary-NS}, proving part~\ref{it:ns-transfer-a}; as $\EpsFri\to0$ both $\omega_N(\EpsFri)\to0$ and $\EpsFri\to0$, so $\EpsStopNS\to2\EpsStop$, and jointly $\EpsStop,\EpsFri\to0$ gives $\EpsStopNS\to0$.

Step~8: uniform boundedness over the fixed-data tolerance family. Fix $\EpsFJ$ and $\EpsFri>0$, fix the remaining problem and algorithm data and initial data, and use common clamps $\LHamilton,\LPosition$ while $0<\EpsStop\leq\EpsFJ$ varies. As in Step~8 of the smooth proof, $\EpsFJAll=\EpsFJ$ on this family and the Hamiltonian majorant at the endpoint $\EpsStop=\EpsFJ$ dominates the family, so the constant-ordering construction supplies one common clamp choice. Because $\EpsFri$ is fixed in part~\ref{it:ns-transfer-b}, the smoothed dissipation and its relevant derivatives, the stationarity-Jacobian constants, and all smoothing-dependent SQP constants are fixed problem data while only $\EpsStop$ varies.

The remaining construction in Step~8 of the smooth proof therefore applies verbatim to this fixed smoothed CRDS. The runs share the common initial optimized right-endpoint set $G_0=\mathcal I^0$ and lie in the finite family $\mathfrak G$ of reachable midpoint-refined meshes constructed there, with the common positive timestep floor $\delta_{\min}$. The same finite-mesh induction gives compact enclosures for every restoration, trial, and inner-SQP iterate and the common residual ceiling $R^\star$. On those enclosures the meshwise LICQ bounds and~\prettyref{lem:sqp-multiplier} give the common multiplier ceiling $C_\mu^\star$ used in the maximum update, and threading that update across mesh episodes and outer calls gives the common penalty cap $\varrho_{\max}$ of~\prettyref{eq:transfer-common-penalty}. Applying~\prettyref{lem:trajectory-bound} with $R^\star$ and $\delta_{\min}$ gives the common trajectory bound, while the fixed-$\EpsFJ$ level set gives the finite prescribed-force ceiling $\|\lambda_i^e(\theta_i)\|,\|\lambda_i^i(\theta_i)\|\leq\bar\Lambda$, exactly as in the smooth proof.

The homogeneous dual certificate is normalized to the unit sphere~\prettyref{eq:eps-stationary-NS-normalization} by construction. Part~\ref{it:ns-transfer-b} asserts no bound on the corresponding unnormalized nonsmooth dual tuple; in particular, the per-term Clarke normals' velocity blocks may be large at sticking. Hence the trajectory, prescribed-force, and penalty bounds are uniform in $\EpsStop$ on $0<\EpsStop\leq\EpsFJ$ for the fixed $\EpsFJ,\EpsFri$ setup, proving part~\ref{it:ns-transfer-b}.
\end{proof}

\subsection{Articulated Body Specialization}
\label{appen:abd-specialization}
This subsection discharges, for articulated-body dynamics, the two results claimed in the main body---the formal restatement of~\prettyref{lem:abd-compatible} and~\prettyref{thm:abd-C-stationarity}---together with all supporting lemmas. Throughout, the smoothing is the pseudo-Huber construction of~\prettyref{lem:smooth-DRp}: every Euclidean velocity norm $\|v_\parallel^m\|$ of the frictional damping~\prettyref{eq:damping-ABD} is replaced by $\psi_\mu(v_\parallel^m)$, the plane-velocity damper $\EpsDamp\|\varpi\|^2$ is retained, and the plane velocities $\varpi$ are eliminated by partial minimization. We first record the primitive weighted-norm geometry (\prettyref{lem:pseudo-huber-derivatives},~\ref{lem:weighted-norm-clarke}) and its weight-lifted master graph (\prettyref{lem:master-weighted-norm-clarke}), then the core three-region normal transfer with its weight-normal bookkeeping (\prettyref{lem:three-region-transfer},~\ref{lem:tau-bookkeeping}) and its exact return to full velocities (\prettyref{lem:vertex-preimage},~\ref{lem:exact-pullback},~\ref{lem:normal-pullback}), resting on the per-vertex activation estimates of~\prettyref{lem:pervertex-activation} and the smooth-normal identification of~\prettyref{lem:smooth-normal-identification}, then the per-vertex normal item (\prettyref{lem:abd-normal-compatible}), assembling~\prettyref{lem:abd-compatible}, and finally the exact reduced-to-full Schur lift---linear, hence commuting with scalar multiplication only and not with the certificate norm, so the full certificate scale is formed on the full tuple before normalization rather than inherited from~\prettyref{def:S}---and the concrete Clarke-stationarity theorem (\prettyref{lem:varpi-reconstruction} through~\ref{thm:abd-C-stationarity}). The pseudo-Huber profile $\psi_\mu$ acts on each vertex relative velocity $v_\parallel^m$. We collect its derivatives and their closeness to the Euclidean-norm subdifferential, then the Clarke normal cone of the weighted-norm subdifferential graph, both used pointwise in the three-region transfer.

To begin with, define the proof-local returned-configuration region and a fixed buffer
\begin{ResizedEq}
\label{eq:abd-compact-regions}
&\mathcal K_H\triangleq\{x:V_R(x,\EpsFJ)\leq\LHamilton,\ \|x\|\leq\LPosBound(H)\}\\
&\widehat{\mathcal K}_H\triangleq\{x:V_R(x,\EpsFJ)\leq\LHamilton+\Delta_H,\ \|x\|\leq\LPosBound(H)+\Delta_x\},
\end{ResizedEq}
where the fixed margins $\Delta_H,\Delta_x>0$ are chosen before the smoothing constants. We record directly why the displayed sets are compact, rather than appealing to dissipation regularity. Let $x_k\to x$ lie in a fixed norm ball with $V_R(x_k,\EpsFJ)\leq C$, and choose minimizing separating planes $p_k$. Finiteness of the normal-length penalty gives $\|n_k\|<1$, while one vertex from each nonempty hull and the two oppositely oriented finite slab penalties bound every offset $b_k$ from below and above, uniformly for bounded $x_k$. Thus $(p_k)$ is bounded. Pass first to a subsequence realizing $\liminf_k V_R(x_k,\EpsFJ)$ and then to a further subsequence on which $p_k\to p$. Extended-real lower semicontinuity of the unreduced scalar-penalty sum gives
\begin{ResizedEq}
V_R(x,\EpsFJ)
\leq&V(x,p,\EpsFJ)\\
\leq&\liminf_k V(x_k,p_k,\EpsFJ)\\
=&\liminf_k V_R(x_k,\EpsFJ)\leq C.
\end{ResizedEq}
Hence every finite sublevel of $V_R$ intersected with a norm ball is closed and bounded, and therefore compact in finite dimensions; in particular both sets in~\prettyref{eq:abd-compact-regions} are compact. This argument also supplies the extended-real lower semicontinuity needed below and is consistent with the finite-domain regularity of~\prettyref{lem:twice-differentiability}. All boundary neighborhoods below are taken inside the fixed buffer $\widehat{\mathcal K}_H$. Constants may depend on the fixed horizon, $\EpsFJ$, clamps, buffer, geometry, and $\EpsDamp$, but not on $\EpsFri$ or $\EpsStop$ over the fixed tolerance regime. Membership of the returned trajectory in $\mathcal K_H$ is established only in the proof of~\prettyref{thm:abd-C-stationarity}, after the reduced solve.

The per-vertex normal transfer of~\prettyref{appen:abd-specialization} consumes the boundary asymptotics of~\prettyref{lem:Vs-flatness} and~\ref{lem:boundary-flatness} in a per-vertex form: it needs, term by term, the size of the penalty-induced friction weight, of the derivative of that weight, and of the derivative of the relative-velocity map, all measured against the vertex penetration depth. We record these consequences once here. Crucially, the plane-derivative $D_{\theta_{\gamma(i)}}C^m$ is not uniformly bounded across the activation boundary---it inherits the mildly singular $\mathcal{O}(\rho^{-1})$ scaling of $Dp^\star$ from~\prettyref{lem:boundary-flatness}---whereas its product with the vanishing weight is bounded; the transfer will always pair the singular factor with the weight.
\begin{lemma}[Per-vertex activation and coefficient estimates]
\label{lem:pervertex-activation}
Fix a vertex term $m$, an active-slab vertex penalty on one hull of a convex-hull pair $(jh,j'h')$, and abbreviate $x\triangleq\theta_{\gamma(i)}$, the penalty-induced friction weight $\alpha^m\triangleq\eta[\Lambda^i]^m$ of~\prettyref{eq:penalty-force}, the relative-velocity map $C^m$ defined by $v_\parallel^m=C^m(x)(\dot{\theta}_i,\varpi)$ through~\prettyref{eq:tangent-velocity}, and the slab-penetration depth $d^m\triangleq d^{jh,j'h'k'}_+\in[0,\rho]$ of that vertex, with $\rho$ the maximal active depth. Work on the active side $V_R^{jh,j'h'}(x)>0$ with $x\in\widehat{\mathcal K}_H$ and $\dot{\theta}_i$ ranging over a bounded set. Then the estimates of~\prettyref{lem:Vs-flatness} and~\ref{lem:boundary-flatness} imply, with constants $0<c\leq C$ uniform over the finite vertex set, timesteps, and $x\in\widehat{\mathcal K}_H$,
\begin{ResizedEq}
&c\,(d^m)^3\leq\alpha^m\leq C\,(d^m)^3\\
&\|C^m(x)\|\leq C\\
&\|D_xC^m\|\leq C(1+\rho^{-1})\\
&\|D_x\alpha^m\|\leq C\,d^m\leq C\,(\alpha^m)^{1/3}\\
&\alpha^m\,\|D_xC^m\|\leq C\,(\alpha^m)^{2/3};
\end{ResizedEq}
consequently, for every threshold $t\in(0,1]$,
\begin{ResizedEq}
\label{eq:pv-threshold}
\alpha^m\geq t\quad\Longrightarrow\quad\|D_xC^m\|\leq C\,t^{-1/3}.
\end{ResizedEq}
\end{lemma}
\begin{proof}[Proof of~\prettyref{lem:pervertex-activation}]
Write $s_m$ for the active $V_s$ argument of the vertex term $m$, so its slab-penetration depth is $d^m=(\EpsFJ-s_m)_+$ and, by~\prettyref{eq:penalty-force}, $\alpha^m=-\eta\,V_s'(s_m,\EpsFJ)\,\|n^\star\|$. The two-sided cubic bound $c\,(d^m)^3\leq|V_s'(s_m,\EpsFJ)|\leq C\,(d^m)^3$ of~\prettyref{lem:Vs-flatness}, together with $\sqrt{1-\EpsFJ}\leq\|n^\star\|\leq1$ from~\prettyref{lem:active-plane}\ref{it:active-plane-normal}, gives the first line.

The map $C^m$ is assembled from the rigid-body kinematics of~\prettyref{eq:tangent-velocity} and the normalized separating-plane data $\bar{n}^{jh,j'h'},P^{jh,j'h'}=I-\bar{n}^{jh,j'h'}[\bar{n}^{jh,j'h'}]^T$, which are $C^1$ on the active side by~\prettyref{lem:active-plane}\ref{it:active-plane-convex}. Bounded $x$ and the finite local vertex set uniformly bound every world vertex $X^m(x)$ and its rigid-kinematic derivative. Write the convex-hull-pair potential of~\prettyref{lem:active-plane} as $V^{\mathrm{pair}}(x,p)$, where $p=(n,b)$, and define the zero-cost configuration-plane lift
\begin{ResizedEq}
\mathcal Z=\{(x,p):x\in\widehat{\mathcal K}_H,\ V^{\mathrm{pair}}(x,p)=0\}.
\end{ResizedEq}
The plane variable is uniformly bounded on this lift, and more generally on every fixed finite sublevel of $V^{\mathrm{pair}}$ over $x\in\widehat{\mathcal K}_H$. Indeed, finiteness of the normal-length penalty gives $\|n\|<1$. Choose one vertex from each nonempty hull. Finiteness of their two oppositely oriented scalar penalties gives
\begin{ResizedEq}
&\langle n,X^{jhk}(x)\rangle+b>\EpsContactPadding\\
&\langle n,X^{j'h'k'}(x)\rangle+b<-\EpsContactPadding.
\end{ResizedEq}
Since the world vertices and $n$ are uniformly bounded, the first inequality bounds $b$ from below and the second bounds it from above, uniformly in $x$. This also bounds all active minimizing planes on the fixed finite pair-potential level inherited from the buffered Hamiltonian sublevel. The potential $V^{\mathrm{pair}}$ is nonnegative and extended-real lower semicontinuous (and continuous on its finite domain) by the scalar-penalty construction, so its zero sublevel is closed. Thus $\mathcal Z$ is a closed subset of the product of the compact set $\widehat{\mathcal K}_H$ and a fixed finite-dimensional ball, and is compact.

We next show that the active boundary region is uniformly controlled by this compact lift. Let $x_k\in\widehat{\mathcal K}_H$ be active with maximal penetration $\rho_k\downarrow0$, and put $p_k=p^\star(x_k)=(n_k,b_k)$. The minimizing pair values remain on the fixed finite level just described, hence $(p_k)$ is bounded. After passage to a subsequence, compactness of $\widehat{\mathcal K}_H$ and boundedness of $(p_k)$ give $(x_k,p_k)\to(\bar x,\bar p)$. Every active scalar vertex penalty tends to zero by~\prettyref{lem:Vs-flatness}, because each vertex depth is at most $\rho_k$. It remains to control the normal-length term without presupposing that the limit lies in $\mathcal Z$. Define
\begin{ResizedEq}
q_{n,k}\triangleq1-\|n_k\|^2,
\qquad
\rho_{n,k}\triangleq\EpsFJ-q_{n,k}.
\end{ResizedEq}
The nonnegative summand $V_s(q_{n,k},\EpsFJ)$ is bounded by the fixed finite pair-potential level and diverges as $q_{n,k}\downarrow0$, so $q_{n,k}$ stays uniformly away from zero. For each slab vertex $a$, write $s_{a,k}$ for its $V_s$ argument and $d_{a,k}=(\EpsFJ-s_{a,k})_+\leq\rho_k$. Radial stationarity of the minimizing plane along $p_\alpha=\alpha p_k$ gives
\begin{ResizedEq}
0=&-2\|n_k\|^2V_s'(q_{n,k},\EpsFJ)
+\sum_a V_s'(s_{a,k},\EpsFJ)
\left.\frac{\partial s_{a,k}(\alpha)}{\partial\alpha}\right|_{\alpha=1}.
\end{ResizedEq}
The sum is finite by the standing finite-vertex model. Moreover, boundedness of $p_k$ and of the world vertices on $\widehat{\mathcal K}_H$, together with the fixed contact padding, uniformly bounds every displayed radial derivative. By~\prettyref{lem:Vs-flatness}, each active summand is $\mathcal O(d_{a,k}^3)=\mathcal O(\rho_k^3)$, while each inactive summand vanishes. Since $\|n_k\|^2\geq1-\EpsFJ>0$ by~\prettyref{lem:active-plane}\ref{it:active-plane-normal}, radial stationarity therefore yields
\begin{ResizedEq}
|V_s'(q_{n,k},\EpsFJ)|=\mathcal O(\rho_k^3)\longrightarrow0.
\end{ResizedEq}
The derivative $V_s'(\cdot,\EpsFJ)$ is continuous and strictly negative on compact subsets of $(0,\EpsFJ)$. Because the active normal-length argument $q_{n,k}<\EpsFJ$ stays uniformly away from zero, the preceding limit forces $q_{n,k}\to\EpsFJ$, equivalently $\rho_{n,k}\to0$. Hence $V_s(q_{n,k},\EpsFJ)\to0$. Together with the vanishing slab penalties this proves $V^{\mathrm{pair}}(x_k,p_k)\to0$. Lower semicontinuity and nonnegativity give $V^{\mathrm{pair}}(\bar x,\bar p)=0$, so $(\bar x,\bar p)\in\mathcal Z$. It follows that
\begin{ResizedEq}
\operatorname{dist}((x_k,p_k),\mathcal Z)\longrightarrow0:
\end{ResizedEq}
otherwise a subsequence staying a fixed positive distance from $\mathcal Z$ would, by the same boundedness argument, have a further subsequence converging to a point of $\mathcal Z$.

Fix $(\bar x,\bar p)\in\mathcal Z$. The proof of~\prettyref{lem:boundary-flatness}, applied with the zero-cost comparison plane $\bar p$, supplies a neighborhood in configuration-plane space whose active minimizing pairs satisfy the plane coercivity and sensitivity estimates. The constants are locally uniform in the pair $(x,p)$: $\EpsFJ$, the contact padding, and the hull geometry are fixed; the incidence and vertex sets are finite; the active normal has the positive lower bound from~\prettyref{lem:active-plane}\ref{it:active-plane-normal}; world vertices and their rigid-kinematic derivatives are uniformly bounded on $\widehat{\mathcal K}_H$; and the scalar-penalty derivative estimates near the fixed cutoff are explicit and uniform. If $\mathcal Z$ is nonempty, these pair neighborhoods cover the compact set $\mathcal Z$. Extracting a finite subcover, taking the minimum positive coercivity constant and the maximum derivative constant, and shrinking to a common positive penetration radius gives $\rho_0>0$ and common constants for all active points with $0<\rho\leq\rho_0$. If $\mathcal Z$ is empty, the preceding subsequence argument instead shows directly that $\rho$ has a positive lower bound on the active minimizing-plane graph; choose $\rho_0$ smaller than that bound, so the away-from-boundary compactness argument below applies to the entire graph. Finally, in the nonempty case, if no common $\rho_0$ existed, there would be active points with $\rho_k\downarrow0$ outside the covered active neighborhoods, contradicting the distance-to-$\mathcal Z$ conclusion above. Consequently, whenever the boundary regime is nonempty, the proof of~\prettyref{lem:boundary-flatness} gives uniformly
\begin{ResizedEq}
\|Dp^\star\|\leq C\rho^{-1},\qquad 0<\rho\leq\rho_0.
\end{ResizedEq}
Differentiating $C^m$ and using $D\bar n,DP=\mathcal{O}(\|Dp^\star\|)$ now gives $\|C^m\|\leq C$ and $\|D_xC^m\|\leq C(1+\rho^{-1})$ in this boundary region; the additive $1$ absorbs the purely kinematic part. In particular, bounded world positions control the term $(D_x\bar n)\times X^m$, which is exactly why the position ball in~\prettyref{eq:abd-compact-regions} is necessary.

It remains to make the constants uniform away from that boundary region. Define the active minimizing-plane remainder graph
\begin{ResizedEq}
\mathcal G_{\rho_0}=\{(x,p^\star(x)):x\in\widehat{\mathcal K}_H,\ x\text{ active},\ \rho(x)\geq\rho_0\}.
\end{ResizedEq}
It is bounded by compactness of $\widehat{\mathcal K}_H$ and the finite-level plane bound above. To prove it closed, consider a convergent graph sequence $(x_k,p^\star(x_k))\to(\bar x,\bar p)$. The fixed pair-potential bound and lower semicontinuity keep $(\bar x,\bar p)$ in the finite domain. The finitely many signed-distance depth functions are continuous there, so $\rho(\bar x,\bar p)\geq\rho_0$ and the limiting configuration remains active. For any $q$ with $V^{\mathrm{pair}}(\bar x,q)<\infty$, all scalar-barrier arguments at $(\bar x,q)$ are strictly positive. The same fixed $q$ therefore remains in the finite domain for all sufficiently large $k$, and continuity there together with the minimizing inequalities gives
\begin{ResizedEq}
V^{\mathrm{pair}}(\bar x,\bar p)
\leq&\liminf_k V^{\mathrm{pair}}(x_k,p^\star(x_k))\\
\leq&\limsup_k V^{\mathrm{pair}}(x_k,q)\\
=&V^{\mathrm{pair}}(\bar x,q).
\end{ResizedEq}
Thus $\bar p$ is a minimizer at $\bar x$, and uniqueness from~\prettyref{lem:active-plane}\ref{it:active-plane-convex} identifies $\bar p=p^\star(\bar x)$. Hence $\mathcal G_{\rho_0}$ is closed and bounded in finite dimensions and is compact.

On $\mathcal G_{\rho_0}$ the plane Hessian is positive definite by~\prettyref{lem:active-plane}\ref{it:active-plane-convex}. Its smallest eigenvalue therefore has a positive minimum, while the continuous implicit-function expression $Dp^\star=-V_{pp}^{-1}V_{p x}$ is uniformly bounded. The same compactness and the explicit cutoff limits in~\prettyref{lem:Vs-flatness} uniformly bound $p^\star,Dp^\star,C^m,D_xC^m,\alpha^m,D_x\alpha^m$ and preserve the two-sided cubic ratio $\alpha^m/(d^m)^3$ for active $m$ (with its positive cutoff limit). Enlarging $c,C$ once combines these remainder bounds with the finite-cover boundary bounds.

For the third line differentiate $\alpha^m=-\eta\,V_s'(s_m,\EpsFJ)\,\|n^\star\|$
\begin{ResizedEq}
\|D_x\alpha^m\|\leq\eta\,|V_s''(s_m)|\,\|Ds_m\|\,\|n^\star\|+\eta\,|V_s'(s_m)|\,\big\|D\|n^\star\|\big\|.
\end{ResizedEq}
By~\prettyref{lem:Vs-flatness}, $|V_s''(s_m)|\leq C\,(d^m)^2$ and $|V_s'(s_m)|\leq C\,(d^m)^3$; by Step~2-3 of~\prettyref{lem:boundary-flatness}, $\|Ds_m\|+\big\|D\|n^\star\|\big\|\leq C(1+\rho^{-1})$. Hence
\begin{ResizedEq}
\|D_x\alpha^m\|\leq C\Big((d^m)^2+\frac{(d^m)^2}{\rho}+\frac{(d^m)^3}{\rho}\Big)\leq C\,d^m,
\end{ResizedEq}
using $d^m\leq\rho$ and $d^m\leq1$ locally; the bound $C\,d^m\leq C\,(\alpha^m)^{1/3}$ then follows from the first line.

For the fourth line, the first two lines give $\alpha^m\|D_xC^m\|\leq C(d^m)^3(1+\rho^{-1})\leq C\big((d^m)^3+(d^m)^2\big)\leq C(d^m)^2\leq C(\alpha^m)^{2/3}$, using $d^m\leq\rho$, $d^m\leq1$, and the first line. Finally, if $\alpha^m\geq t$ then the first line gives $\rho\geq d^m\geq c\,t^{1/3}$, and substitution into the derivative bound yields~\prettyref{eq:pv-threshold}.
\end{proof}

\begin{lemma}[Pseudo-Huber derivative and error bounds]
\label{lem:pseudo-huber-derivatives}
For $v\in\R^n$ and $\mu>0$ write $s\triangleq\sqrt{\|v\|^2+\mu^2}$, and for $v\neq0$ write $e\triangleq v/\|v\|$. Then $\psi_\mu(v)=\sqrt{\|v\|^2+\mu^2}-\mu$ is $C^\infty$ with
\begin{ResizedEq}
\label{eq:ph-grad-hess}
\nabla\psi_\mu(v)=&\frac{v}{s}\\
\nabla^2\psi_\mu(v)=&\frac1s\left(I-\frac{vv^T}{s^2}\right)=\frac1s P_e+\frac{\mu^2}{s^3}ee^T,
\end{ResizedEq}
the second identity holding for $v\neq0$ with $P_e\triangleq I-ee^T$, and $\nabla^2\psi_\mu(0)=I/\mu$. Consequently $\nabla^2\psi_\mu(v)\succeq0$ with, for $n\geq2$, tangential eigenvalue $1/s$ (multiplicity $n-1$), and radial eigenvalue $\mu^2/s^3$; when $n=1$, the tangential multiplicity is zero. Moreover, for $\|v\|>0$,
\begin{ResizedEq}
&\|\nabla\psi_\mu(v)\|\leq1\\
&\left\|\nabla\psi_\mu(v)-\frac{v}{\|v\|}\right\|\leq\frac{\mu^2}{2\|v\|^2}\\
&\left\|\nabla^2\psi_\mu(v)-\frac{P_e}{\|v\|}\right\|_{\mathrm{op}}\leq\frac{\mu^2}{\|v\|^3}.
\end{ResizedEq}
\end{lemma}
\begin{proof}[Proof of~\prettyref{lem:pseudo-huber-derivatives}]
Differentiating $\psi_\mu(v)=s-\mu$ with $s=\sqrt{\|v\|^2+\mu^2}$ gives $\nabla s=v/s$, hence $\nabla\psi_\mu(v)=v/s$; differentiating once more, $\nabla^2\psi_\mu(v)=\frac1s I-\frac{v}{s^2}(\nabla s)^T=\frac1s I-\frac{vv^T}{s^3}=\frac1s(I-vv^T/s^2)$, which at $v=0$ ($s=\mu$) is $I/\mu$. For $v\neq0$, $vv^T=\|v\|^2ee^T$ and $I-\|v\|^2ee^T/s^2=P_e+(1-\|v\|^2/s^2)ee^T=P_e+(\mu^2/s^2)ee^T$, giving the split~\prettyref{eq:ph-grad-hess}; its eigenvalues are $1/s$ on $e^\perp$ and $\frac1s\cdot\frac{\mu^2}{s^2}=\mu^2/s^3$ along $e$, both nonnegative. The gradient bound is $\|\nabla\psi_\mu(v)\|=\|v\|/s\leq1$. For the gradient error, with $r=\|v\|$,
\begin{ResizedEq}
\left\|\frac{v}{s}-\frac{v}{r}\right\|=r\Big(\frac1r-\frac1s\Big)=\frac{s-r}{s}=\frac{\mu^2}{s(s+r)}\leq\frac{\mu^2}{2r^2},
\end{ResizedEq}
using $s\geq r$ and $s+r\geq2r$. For the Hessian error, subtract $P_e/r$ from~\prettyref{eq:ph-grad-hess}: the tangential part contributes $(\frac1s-\frac1r)P_e$ with $\big|\frac1s-\frac1r\big|=\frac{s-r}{rs}=\frac{\mu^2}{rs(s+r)}\leq\frac{\mu^2}{2r^3}$, and the radial part contributes $\frac{\mu^2}{s^3}ee^T$ with norm $\frac{\mu^2}{s^3}\leq\frac{\mu^2}{r^3}$. Since $P_e$ and $ee^T$ have orthogonal ranges $e^\perp$ and $\operatorname{span}\{e\}$, respectively, the induced Euclidean operator norm of their weighted sum equals the maximum of the tangential and radial block norms, and is therefore at most $\max\{\mu^2/(2r^3),\mu^2/r^3\}=\mu^2/r^3$.
\end{proof}
\begin{lemma}[Clarke normal cone of the weighted-norm graph]
\label{lem:weighted-norm-clarke}
Fix a weight $\alpha\geq0$ and let $\mathcal{G}_\alpha\triangleq\gph(\alpha\partial\|\cdot\|)\subseteq\R^n\times\R^n$, the graph of the set-valued map $v\mapsto\alpha\partial\|v\|$, so $\mathcal{G}_\alpha=\{(0,q):\|q\|\leq\alpha\}\cup\{(v,\alpha v/\|v\|):v\neq0\}$. Its Clarke normal cone $N^C_{\mathcal{G}_\alpha}$ is, for $\alpha>0$,
\begin{enumerate}[label=(\alph*),ref=(\alph*),leftmargin=*]
\item\label{it:wnc-strict} strict stick, $v=0$ and $\|q\|<\alpha$: $\ N^C_{\mathcal{G}_\alpha}(0,q)=\R^n\times\{0\}$;
\item\label{it:wnc-boundary} stick boundary, $v=0$ and $q=\alpha e$ with $\|e\|=1$: $\ N^C_{\mathcal{G}_\alpha}(0,\alpha e)=\R^n\times\mathrm{span}\{e\}$;
\item\label{it:wnc-slide} sliding, $v=re$ with $r>0$ and $q=\alpha e$: $\ N^C_{\mathcal{G}_\alpha}(re,\alpha e)=\{(\tfrac{\alpha}{r}P_e z,-z):z\in\R^n\}$.
\end{enumerate}
For $\alpha=0$, $\mathcal{G}_0=\R^n\times\{0\}$ and $N^C_{\mathcal{G}_0}(v,0)=\{0\}\times\R^n$.
\end{lemma}
\begin{proof}[Proof of~\prettyref{lem:weighted-norm-clarke}]
The convex subdifferential of $v\mapsto\alpha\|v\|$ is $\{q:\|q\|\leq\alpha\}$ at $v=0$ and the singleton $\{\alpha v/\|v\|\}$ for $v\neq0$, giving the stated $\mathcal{G}_\alpha$. We compute the Clarke normal cone (the closed convex hull of limiting normals) at each point.

\ref{it:wnc-strict} Near $(0,q)$ with $\|q\|<\alpha$ the graph is locally the affine set $\{0\}\times\R^n$ (only the $v=0$ branch enters, as the sliding branch carries $\|q'\|=\alpha$), and $q$ is interior to the force ball, so the normal cone of this manifold is $\R^n\times\{0\}$.

\ref{it:wnc-slide} For $v\neq0$ the graph is the smooth manifold $\{(v,\alpha v/\|v\|)\}$; its tangent space at $v=re$ is the range of $(I,\tfrac{\alpha}{r}P_e)$ (using $D(v/\|v\|)=P_e/\|v\|$). A covector $(a,-z)$ is normal iff $a=(\tfrac{\alpha}{r}P_e)^Tz=\tfrac{\alpha}{r}P_e z$ ($P_e$ symmetric); ranging $z$ gives the stated cone, and regular, limiting, and Clarke cones coincide on a smooth manifold.

\ref{it:wnc-boundary} Put $p=(0,\alpha e)$. Directions through the stick-ball branch are $\{(0,w):e^Tw\leq0\}$, while directions through sliding points $(re',\alpha e')\to p$ are $\{(\rho e,w):\rho\geq0,\ w\perp e\}$. Indeed, at a tangent scale $t_k\downarrow0$, finite convergence along the sliding branch forces $r_k/t_k\to\rho\geq0$ and $\alpha(e_k-e)/t_k\to w\perp e$, whence $r_ke_k/t_k\to\rho e$; conversely, both displayed families are realized directly. Thus
\begin{ResizedEq}
T_{\mathcal G_\alpha}(p)=\{(0,w):e^Tw\leq0\}\cup\{(\rho e,w):\rho\geq0,\ w\perp e\}.
\end{ResizedEq}
Polarizing this Bouligand tangent cone gives the regular normal
\begin{ResizedEq}
\widehat N_{\mathcal G_\alpha}(p)=\{(a,\lambda e):a^Te\leq0,\ \lambda\geq0\}.
\end{ResizedEq}
We now include limiting normals from all nearby strata. Strict-stick points contribute $\R^n\times\{0\}$. Nearby boundary points $(0,\alpha e')$, $e'\to e$, contribute $\{(a,\lambda e):a^Te\leq0,\ \lambda\geq0\}$ by the preceding regular-normal formula with $e'$ in place of $e$. At sliding points $(re',\alpha e')$, regular normals are $(\tfrac{\alpha}{r}P_{e'}z,-z)$. If such normals converge finitely to $(a,b)$ as $r\downarrow0$ and $e'\to e$, then $P_{e'}z=\mathcal{O}(r)$, so $b\in\mathrm{span}\{e\}$, and the first block's orthogonality to $e'$ gives $a\perp e$. Conversely, every $(a,\lambda e)$ with $a\perp e$ and $\lambda\in\R$ is obtained by taking $e'=e$ and $z=-\lambda e+(r/\alpha)a$. Hence every limiting normal has force block in $\mathrm{span}\{e\}$. Strict-stick limits supply every $(a,0)$, while boundary and sliding limits supply $(0,\lambda e)$ with both signs of $\lambda$; their closed convex hull is exactly $\R^n\times\mathrm{span}\{e\}$.

For $\alpha=0$, $\alpha\partial\|v\|=\{0\}$, so $\mathcal{G}_0=\R^n\times\{0\}$ is a subspace and its normal cone is $\{0\}\times\R^n$.
\end{proof}

A fixed-weight normal from~\prettyref{lem:weighted-norm-clarke} carries no covector dual to the weight $\alpha$, yet in the articulated model $\alpha=\eta[\Lambda^i]^m(\theta_{\gamma(i)})$ varies with the previous configuration, so the eventual full normal must include a weight-normal component. We record it by promoting the weight to a graph coordinate; the resulting master graph is used in the weight-normal bookkeeping of~\prettyref{lem:tau-bookkeeping} and the pullback of~\prettyref{lem:exact-pullback}.
\begin{lemma}[Clarke normal cone of the master weighted-norm graph]
\label{lem:master-weighted-norm-clarke}
Let $\widehat{\mathcal{G}}\triangleq\{(\alpha,v,q):\alpha\geq0,\ q\in\alpha\partial\|v\|\}\subseteq\R\times\R^n\times\R^n$ be the master weighted-norm graph, whose slice at fixed $\alpha$ is the graph $\mathcal{G}_\alpha$ of~\prettyref{lem:weighted-norm-clarke}. For $\alpha>0$ its Clarke normal cone $N^C_{\widehat{\mathcal{G}}}$ is
\begin{enumerate}[label=(\alph*),ref=(\alph*),leftmargin=*]
\item\label{it:mwnc-strict} strict stick, $v=0$ and $\|q\|<\alpha$: $\ N^C_{\widehat{\mathcal{G}}}(\alpha,0,q)=\{0\}\times\R^n\times\{0\}$;
\item\label{it:mwnc-boundary} stick boundary, $v=0$ and $q=\alpha e$ with $\|e\|=1$: $\ N^C_{\widehat{\mathcal{G}}}(\alpha,0,\alpha e)=\{(t,a,-t e):t\in\R,\ a\in\R^n\}$;
\item\label{it:mwnc-slide} sliding, $v=re$ with $r>0$ and $q=\alpha e$: $\ N^C_{\widehat{\mathcal{G}}}(\alpha,re,\alpha e)=\{(e^Tz,\tfrac{\alpha}{r}P_e z,-z):z\in\R^n\}$.
\end{enumerate}
For $\alpha=0$ two cases must be distinguished, according to whether the velocity is nonzero or the point is fully degenerate:
\begin{enumerate}[label=(\alph*),ref=(\alph*),leftmargin=*,resume]
\item\label{it:mwnc-zero-v} zero weight with nonzero velocity, $v\neq0$ and $e=v/\|v\|$: $\ N^C_{\widehat{\mathcal{G}}}(0,v,0)=\{(\tau,0,b):\tau+b^Te\leq0\}$;
\item\label{it:mwnc-origin} fully degenerate apex, $v=0$: $\ N^C_{\widehat{\mathcal{G}}}(0,0,0)=\R\times\R^n\times\R^n$.
\end{enumerate}
For $\alpha>0$, and also for $\alpha=0$ with $v\neq0$, dropping the first (weight) coordinate returns the corresponding fixed-weight cone of~\prettyref{lem:weighted-norm-clarke} exactly; at the fully degenerate apex $(0,0,0)$ this projection is strictly larger than the fixed-weight cone, as the master graph is more singular there.
\end{lemma}
\begin{proof}[Proof of~\prettyref{lem:master-weighted-norm-clarke}]
At a strict-stick point the force ball is inactive and $\widehat{\mathcal{G}}$ is locally the manifold $\{v=0\}$ in $(\alpha,v,q)$-space, whose normal cone is $\{0\}\times\R^n\times\{0\}$; this is~\ref{it:mwnc-strict}.

At a sliding point $\widehat{\mathcal{G}}$ is the smooth manifold $q=\alpha v/\|v\|$; differentiating gives the tangent variations $\delta q=e\,\delta\alpha+\tfrac{\alpha}{r}P_e\,\delta v$ (using $D(v/\|v\|)=P_e/\|v\|$). A covector $(\tau,a,b)$ is normal iff it annihilates every tangent $(\delta\alpha,\delta v,\delta q)$, i.e. $\tau+b^Te=0$ and $a+\tfrac{\alpha}{r}P_e b=0$ ($P_e$ symmetric); setting $b=-z$ yields~\ref{it:mwnc-slide}, and on a smooth manifold the regular, limiting, and Clarke cones coincide.

At the kink $P=(\alpha,0,\alpha e)$, directions through the stick branch are $\{(s,0,w):e^Tw\leq s\}$, while directions through the sliding branch are $\{(s,\rho e,se+w_\perp):s\in\R,\ \rho\geq0,\ w_\perp\perp e\}$. Polarizing their union gives the regular normal
\begin{ResizedEq}
\widehat N_{\widehat{\mathcal G}}(P)=\{(t,a,-te):t\leq0,\ a^Te\leq0\}.
\end{ResizedEq}
We again include all nearby strata. Strict-stick points contribute $(0,a,0)$ for arbitrary $a\in\R^n$. Nearby boundary points $(\alpha',0,\alpha'e')\to P$ contribute $\{(t,a,-te):t\leq0,\ a^Te\leq0\}$ by the preceding regular-normal formula. Sliding points contribute $(e'^Tz,\tfrac{\alpha'}{r}P_{e'}z,-z)$. A finite sliding limit has force block $-z\to-te$ for some $t\in\R$ and velocity block orthogonal to $e$; conversely every $(t,a,-te)$ with $a\perp e$ is obtained by taking $e'=e$ and $z=te+(r/\alpha)a$. Thus every limiting normal satisfies the coupling between its weight and force blocks, and the strict-stick family supplies arbitrary velocity blocks. The closed convex hull of the strict-stick, boundary-to-boundary, and sliding limiting normals is therefore exactly $\{(t,a,-te):t\in\R,\ a\in\R^n\}$, giving~\ref{it:mwnc-boundary}.

For $\alpha=0$ the weight is one-sided at the nonnegativity floor, and the geometry differs between $v\neq0$ and the apex $v=0$.

If $v\neq0$, set $e=v/\|v\|$; in a neighborhood of $(0,v,0)$ the velocity stays nonzero, so $\widehat{\mathcal{G}}$ is the image of the manifold-with-boundary parameterization $\mathfrak{F}(\alpha,w)=(\alpha,w,\alpha w/\|w\|)$ with $\alpha\geq0$. At $\alpha=0$ the $w$-derivative of the third block carries the factor $\alpha$ and hence vanishes, so the Clarke tangent cone is $T^C_{\widehat{\mathcal{G}}}(0,v,0)=\{(s,w,se):s\geq0,\,w\in\R^n\}$. A covector $(\tau,a,b)$ lies in its polar iff $\tau s+a^Tw+b^T(se)\leq0$ for all $s\geq0$ and $w\in\R^n$; the free $w$ forces $a=0$ and the one-sided scalar $s\geq0$ gives $\tau+b^Te\leq0$, which is~\ref{it:mwnc-zero-v}.

At the apex we show $T^C_{\widehat{\mathcal{G}}}(0,0,0)=\{0\}$, whose polar is the full space $\R\times\R^n\times\R^n$, giving~\ref{it:mwnc-origin}. Let $d=(\dot\alpha,\dot v,\dot q)$ be a Clarke tangent direction. A Clarke tangent must be realizable along every sequence of graph points approaching the origin and every positive step sequence decreasing to zero. Along strict-stick points $(\alpha_k,0,0)$ with $\alpha_k\downarrow0$, choose a step scale $s_k=o(\alpha_k)$. Since the realizing directions are bounded, the realizing endpoints remain in the strict interior of the force ball, where the graph is locally $\{v=0\}$; hence $\dot v=0$.

Now fix a unit vector $e$ and set $r_k=1/k$, $\alpha_k=1/k^2$, and $x_k=(\alpha_k,r_ke,\alpha_ke)$. Choose $t_k=1/k^3$, so that $t_k=o(\alpha_k)$ and $t_k=o(r_k)$. For realizing directions $d_k=(\dot\alpha_k,\dot v_k,\dot q_k)\to d$, the endpoints $x_k+t_kd_k$ have positive weight and nonzero velocity for all large $k$, and therefore lie on the smooth sliding branch $q=\alpha v/\|v\|$. Taylor expansion at $x_k$, followed by division by $t_k$, gives
\begin{ResizedEq}
&\dot q_k=e\,\dot\alpha_k+\frac{\alpha_k}{r_k}P_e\dot v_k+\varepsilon_k\\
&\|\varepsilon_k\|\leq C\left(\frac{t_k}{r_k}+\frac{\alpha_k t_k}{r_k^2}\right)\longrightarrow0,
\end{ResizedEq}
where boundedness of the realizing directions is used. Since $\alpha_k/r_k=1/k\to0$, passing to the limit yields $\dot q=e\,\dot\alpha$. Repeating the argument with $-e$ gives $\dot q=-e\,\dot\alpha$; hence $\dot\alpha=0$ and $\dot q=0$ for every $n\geq1$. Together with $\dot v=0$, this proves $d=0$ and $T^C_{\widehat{\mathcal{G}}}(0,0,0)=\{0\}$.

Projecting out the weight component maps the stated cone onto the corresponding cone of~\prettyref{lem:weighted-norm-clarke} for $\alpha>0$ and for $\alpha=0$ with $v\neq0$; at the apex, however, the projection of~\ref{it:mwnc-origin} is $\R^n\times\R^n$, strictly larger than $N^C_{\mathcal{G}_0}(0,0)=\{0\}\times\R^n$, so the projection identity fails there.
\end{proof}

The next lemma is the analytic heart of the specialization: it transfers a smooth pseudo-Huber graph normal to a nearby Clarke normal of the exact weighted-norm graph, with a modulus vanishing as $\mu\downarrow0$. It is stated in the primitive relative-velocity coordinate $v$ and the multiplier direction $z$; the return to full articulated velocities follows in~\prettyref{lem:vertex-preimage} and~\ref{lem:normal-pullback}.
\begin{lemma}[Three-region smooth-to-Clarke normal transfer]
\label{lem:three-region-transfer}
Fix finite constants $C_z,C_E,\bar\Lambda$. Let $0\leq\alpha\leq\bar\Lambda$ and $v,z\in\R^n$ satisfy
\begin{ResizedEq}
\label{eq:tr-hyp}
&\|z\|\leq C_z\\
&\alpha\,z^T\nabla^2\psi_\mu(v)\,z\leq C_E,
\end{ResizedEq}
and set the scales $\alpha_\mu=\mu^{1/5}$, $r_c=\mu^{4/5}$, $r_s=\mu^{2/5}$. Then for all sufficiently small $\mu$ there exist a graph point $(v^\sharp,q^\sharp)\in\gph(\alpha\partial\|\cdot\|)$ and a Clarke normal $n^\sharp\in N^C_{\gph(\alpha\partial\|\cdot\|)}(v^\sharp,q^\sharp)$ with
\begin{ResizedEq}
\label{eq:tr-bounds}
&\|v^\sharp-v\|+\|q^\sharp-\alpha\nabla\psi_\mu(v)\|\leq\omega_g(\mu)\\
&\|n^\sharp-(\alpha\nabla^2\psi_\mu(v)z,-z)\|\leq\omega_n(\mu),
\end{ResizedEq}
where $\omega_g(\mu)=C\mu^{1/5}$ and $\omega_n(\mu)=C\mu^{1/10}$ for a constant $C=C(C_z,C_E,\bar\Lambda)$, uniformly over all $(\alpha,v,z)$ obeying~\prettyref{eq:tr-hyp}.
\end{lemma}
\begin{proof}[Proof of~\prettyref{lem:three-region-transfer}]
Write $r=\|v\|$ and, for $r>0$, $e=v/r$; recall from~\prettyref{lem:pseudo-huber-derivatives} that $\nabla^2\psi_\mu(v)$ has tangential eigenvalue $1/s$ and radial eigenvalue $\mu^2/s^3$ with $s=\sqrt{r^2+\mu^2}$. We split into the three regions of the scales, and within the first two into a small-weight branch $\alpha<\alpha_\mu$ and a large-weight branch $\alpha\geq\alpha_\mu$. In each region $q^\sharp\in\alpha\partial\|v^\sharp\|$ by construction, so $(v^\sharp,q^\sharp)$ is a graph point; we verify the two bounds in~\prettyref{eq:tr-bounds}.

Case I: sticking core, $r\leq r_c$. Here $s^2=r^2+\mu^2\leq\mu^{8/5}+\mu^2\leq2\mu^{8/5}$ for small $\mu$, so $s\leq C\mu^{4/5}$ and the radial eigenvalue obeys $\mu^2/s^3\geq c\mu^2/\mu^{12/5}=c\mu^{-2/5}$.

Large weight $\alpha\geq\alpha_\mu$. Since the radial eigenvalue is the smallest eigenvalue of $\nabla^2\psi_\mu(v)$,~\prettyref{eq:tr-hyp} gives $C_E\geq\alpha\,z^T\nabla^2\psi_\mu(v)z\geq\alpha_\mu\cdot c\mu^{-2/5}\|z\|^2=c\mu^{-1/5}\|z\|^2$, whence $\|z\|\leq C\mu^{1/10}$. Take $v^\sharp=0$ and $q^\sharp=\alpha\nabla\psi_\mu(v)=\alpha v/s$; as $\|q^\sharp\|=\alpha r/s<\alpha$, this is a strict-stick graph point, so by~\prettyref{lem:weighted-norm-clarke}\ref{it:wnc-strict} the covector $n^\sharp=(\alpha\nabla^2\psi_\mu(v)z,\,0)$ is a Clarke normal. Its first block equals the target's, so the normal error is $\|(0,-z)\|=\|z\|\leq C\mu^{1/10}$, and the graph displacement is $\|v^\sharp-v\|+0=r\leq\mu^{4/5}$.

Small weight $\alpha<\alpha_\mu$. If $z\neq0$ set $e_z=z/\|z\|$, else let $e_z$ be any unit vector. Take $v^\sharp=0$, $q^\sharp=\alpha e_z$, and $n^\sharp=(\alpha\nabla^2\psi_\mu(v)z,-z)$. Since $-z\in\mathrm{span}\{e_z\}$, \prettyref{lem:weighted-norm-clarke}\ref{it:wnc-boundary} (for $\alpha>0$) or its $\alpha=0$ degenerate case makes $n^\sharp$ a Clarke normal, so the normal error is $0$. The graph displacement is $\|v^\sharp-v\|+\|q^\sharp-\alpha\nabla\psi_\mu(v)\|\leq r+\alpha(\|e_z\|+\|v/s\|)\leq\mu^{4/5}+2\alpha\leq\mu^{4/5}+2\mu^{1/5}\leq C\mu^{1/5}$.

Case II: transition, $r_c<r\leq r_s$. Here $r>\mu^{4/5}\geq\mu$, so $s\leq\sqrt{2}\,r_s\leq C\mu^{2/5}$ and $r\leq s$.

Large weight $\alpha\geq\alpha_\mu$. Take $v^\sharp=0$, $q^\sharp=\alpha e$; a stick-boundary graph point. The force displacement is $\|\alpha e-\alpha v/s\|=\alpha(1-r/s)=\alpha\frac{\mu^2}{s(s+r)}\leq\bar\Lambda\frac{\mu^2}{2r^2}\leq\bar\Lambda\frac{\mu^2}{2\mu^{8/5}}=C\mu^{2/5}$, and $\|v^\sharp-v\|=r\leq\mu^{2/5}$, so the graph displacement is $\leq C\mu^{2/5}\leq C\mu^{1/5}$. Split $z=z_\parallel+z_\perp$ with $z_\parallel=(e^Tz)e$ and $z_\perp=P_e z$; the tangential eigenvalue $1/s$ bounds $C_E\geq\alpha z^T\nabla^2\psi_\mu(v)z\geq\alpha\frac1s\|z_\perp\|^2\geq\alpha_\mu\frac{1}{C\mu^{2/5}}\|z_\perp\|^2=c\mu^{-1/5}\|z_\perp\|^2$, so $\|z_\perp\|\leq C\mu^{1/10}$. Take $n^\sharp=(\alpha\nabla^2\psi_\mu(v)z,-z_\parallel)$; since $-z_\parallel\in\mathrm{span}\{e\}$,~\prettyref{lem:weighted-norm-clarke}\ref{it:wnc-boundary} makes it a Clarke normal, with error from the target $\|(0,-z_\parallel+z)\|=\|z_\perp\|\leq C\mu^{1/10}$.

Small weight $\alpha<\alpha_\mu$. Use the same small-weight construction as Case~I: $v^\sharp=0$, $q^\sharp=\alpha e_z$, $n^\sharp=(\alpha\nabla^2\psi_\mu(v)z,-z)$. The normal error is $0$, the force displacement $\leq2\alpha\leq2\mu^{1/5}$, and $\|v^\sharp-v\|=r\leq\mu^{2/5}$, so the graph displacement is $\leq C\mu^{1/5}$.

Case III: sliding, $r>r_s$. Here $r>\mu^{2/5}\gg\mu$. Take $v^\sharp=v$, $q^\sharp=\alpha e$, and the sliding-branch normal $n^\sharp=(\tfrac{\alpha}{r}P_e z,-z)$, which by~\prettyref{lem:weighted-norm-clarke}\ref{it:wnc-slide} (or its $\alpha=0$ case) is a Clarke normal. The velocity displacement is $0$ and the force displacement is $\alpha(1-r/s)=\alpha\frac{\mu^2}{s(s+r)}\leq\bar\Lambda\frac{\mu^2}{2r^2}\leq\bar\Lambda\frac{\mu^2}{2\mu^{4/5}}=C\mu^{6/5}$. For the normal error, use the split $\nabla^2\psi_\mu(v)=\frac1s P_e+\frac{\mu^2}{s^3}ee^T$
\begin{ResizedEq}
\label{eq:tr-slide-normal}
&\alpha\nabla^2\psi_\mu(v)z-\frac{\alpha}{r}P_e z\\
=&\alpha\Big(\frac1s-\frac1r\Big)P_e z+\alpha\frac{\mu^2}{s^3}(e^Tz)e,
\end{ResizedEq}
whose norm is at most $\alpha\frac{\mu^2}{2r^3}\|z\|+\alpha\frac{\mu^2}{s^3}\|z\|\leq\frac32\bar\Lambda C_z\frac{\mu^2}{r^3}\leq\frac32\bar\Lambda C_z\frac{\mu^2}{\mu^{6/5}}=C\mu^{4/5}$ (using $|\tfrac1s-\tfrac1r|\leq\mu^2/(2r^3)$ and $s\geq r$); the second block matches the target exactly, so the normal error is $\leq C\mu^{4/5}$.

Taking the maximum over the three regions, the graph displacement is at most $C\mu^{1/5}=\omega_g(\mu)$ and the normal error at most $C\mu^{1/10}=\omega_n(\mu)$ (the exponents $1/5$ and $1/10$ being the smallest arising, hence dominant as $\mu\downarrow0$), uniformly over all data satisfying~\prettyref{eq:tr-hyp}.
\end{proof}

The fixed-weight normal $(a^\sharp,b^\sharp)$ produced by~\prettyref{lem:three-region-transfer} carries no weight-normal component; the pullback of~\prettyref{lem:normal-pullback} needs its extension to the master graph $\widehat{\mathcal{G}}$ of~\prettyref{lem:master-weighted-norm-clarke}, together with a bound on the mismatch from the smooth weight-normal $\tau_\mu\triangleq\nabla\psi_\mu(v)^Tz$.
\begin{lemma}[Weight-normal bookkeeping]
\label{lem:tau-bookkeeping}
Under the hypotheses and construction of~\prettyref{lem:three-region-transfer}, let $(v^\sharp,q^\sharp)$ and its fixed-weight Clarke normal $(a^\sharp,b^\sharp)\triangleq n^\sharp$ be the graph point and normal produced there, and set $\tau_\mu\triangleq\nabla\psi_\mu(v)^Tz$. Then there is a scalar $\tau^\sharp$ with
\begin{ResizedEq}
\label{eq:tau-membership}
(\tau^\sharp,a^\sharp,b^\sharp)\in N^C_{\widehat{\mathcal{G}}}(\alpha,v^\sharp,q^\sharp),
\end{ResizedEq}
the graph-force bound $\|q^\sharp-\alpha\nabla\psi_\mu(v)\|\leq2\alpha$, and, with $\alpha_\mu=\mu^{1/5}$ as in~\prettyref{lem:three-region-transfer},
\begin{enumerate}[label=(\alph*),ref=(\alph*),leftmargin=*]
\item\label{it:tau-large} if $\alpha\geq\alpha_\mu$, then $|\tau^\sharp-\tau_\mu|\leq C\mu^{1/10}$;
\item\label{it:tau-small} if $0\leq\alpha<\alpha_\mu$, then $|\tau^\sharp-\tau_\mu|\leq C$,
\end{enumerate}
with $C=C(C_z,C_E,\bar\Lambda)$ as in~\prettyref{lem:three-region-transfer}.
\end{lemma}
\begin{proof}[Proof of~\prettyref{lem:tau-bookkeeping}]
The graph-force bound is immediate: $q^\sharp\in\alpha\partial\|v^\sharp\|$ gives $\|q^\sharp\|\leq\alpha$, and $\|\alpha\nabla\psi_\mu(v)\|=\alpha\|v\|/s\leq\alpha$, so $\|q^\sharp-\alpha\nabla\psi_\mu(v)\|\leq2\alpha$. For the weight-normal we follow the regions and weight branches of~\prettyref{lem:three-region-transfer}, selecting $\tau^\sharp$ by~\prettyref{lem:master-weighted-norm-clarke}; throughout write $g_\mu\triangleq\nabla\psi_\mu(v)=v/s$ ($\|g_\mu\|\leq1$) and, for $r=\|v\|>0$, $e=v/r$, so $\tau_\mu=g_\mu^Tz=(r/s)\,e^Tz$.

Large weight in the sticking core (Case~I, $\alpha\geq\alpha_\mu$): $q^\sharp$ is a strict-stick point with $b^\sharp=0$, so~\prettyref{lem:master-weighted-norm-clarke}\ref{it:mwnc-strict} forces $\tau^\sharp=0$ and~\prettyref{eq:tau-membership} holds ($a^\sharp$ free). The core estimate of~\prettyref{lem:three-region-transfer} gives $\|z\|\leq C\mu^{1/10}$, hence $|\tau^\sharp-\tau_\mu|=|\tau_\mu|\leq\|z\|\leq C\mu^{1/10}$.

Large weight in the transition (Case~II, $\alpha\geq\alpha_\mu$): $q^\sharp=\alpha e$ is a stick-boundary point with $b^\sharp=-z_\parallel=-(e^Tz)e$, so~\prettyref{lem:master-weighted-norm-clarke}\ref{it:mwnc-boundary} gives~\prettyref{eq:tau-membership} with $\tau^\sharp=e^Tz$. Then
\begin{ResizedEq}
|\tau^\sharp-\tau_\mu|
=&|e^Tz|\Big(1-\frac{r}{s}\Big)
=|e^Tz|\frac{\mu^2}{s(s+r)}\\
\leq&C_z\frac{\mu^2}{2r^2}
\leq C_z\frac{\mu^2}{2\mu^{8/5}}
=C\mu^{2/5}
\leq C\mu^{1/10},
\end{ResizedEq}
using $r>r_c=\mu^{4/5}$ and $s\geq r$.

Sliding (Case~III, any $\alpha$): $(v^\sharp,q^\sharp)=(re,\alpha e)$ with $v^\sharp=re\neq0$. For $\alpha>0$ this is a sliding point with $b^\sharp=-z$, so~\prettyref{lem:master-weighted-norm-clarke}\ref{it:mwnc-slide} gives~\prettyref{eq:tau-membership} with $\tau^\sharp=e^Tz$ (its velocity block $\tfrac{\alpha}{r}P_ez$ equals $a^\sharp$). For $\alpha=0$ the point is $(0,re,0)$ with nonzero velocity, so~\prettyref{lem:master-weighted-norm-clarke}\ref{it:mwnc-zero-v} applies: the selection $\tau^\sharp=e^Tz$, $a^\sharp=0$, $b^\sharp=-z$ lies on the boundary of the corrected cone since $\tau^\sharp+(b^\sharp)^Te=e^Tz-z^Te=0\leq0$, giving~\prettyref{eq:tau-membership}. In either case, as above $|\tau^\sharp-\tau_\mu|=|e^Tz|\tfrac{\mu^2}{s(s+r)}\leq C_z\tfrac{\mu^2}{2r^2}\leq C_z\tfrac{\mu^2}{2\mu^{4/5}}=C\mu^{6/5}\leq C\mu^{1/10}$, using $r>r_s=\mu^{2/5}$; this discharges~\ref{it:tau-large} in the sliding region and, being $\leq C$, also~\ref{it:tau-small} there.

Small weight in the sticking core/transition (Cases~I-II, $0\leq\alpha<\alpha_\mu$): the construction takes $v^\sharp=0$, $q^\sharp=\alpha e_z$ with $e_z=z/\|z\|$ ($z\neq0$; any unit vector if $z=0$), and $n^\sharp=(\alpha\nabla^2\psi_\mu(v)z,-z)$. For $\alpha>0$ this is a stick-boundary point with $e=e_z$ and $b^\sharp=-z=-\|z\|e_z$, so~\prettyref{lem:master-weighted-norm-clarke}\ref{it:mwnc-boundary} gives~\prettyref{eq:tau-membership} with $\tau^\sharp=e_z^Tz=\|z\|$; for $\alpha=0$ the point is the origin $(0,0,0)$, where~\prettyref{lem:master-weighted-norm-clarke}\ref{it:mwnc-origin} gives the full space $N^C_{\widehat{\mathcal{G}}}(0,0,0)=\R\times\R^n\times\R^n$, so~\prettyref{eq:tau-membership} holds with $\tau^\sharp=0$. In either case $|\tau^\sharp|\leq\|z\|\leq C_z$ and $|\tau_\mu|\leq\|z\|\leq C_z$, so $|\tau^\sharp-\tau_\mu|\leq C$, proving~\ref{it:tau-small}.
\end{proof}
\begin{lemma}[Vertex velocity preimage and right inverse]
\label{lem:vertex-preimage}
For each vertex term $m\in\mathcal{T}$ and each fixed choice of its separating-plane data, the relative velocity of~\prettyref{eq:tangent-velocity} is affine in the damping-active velocity $\dot\vartheta_i^D$ at fixed $\theta_{\gamma(i)}$, $v_\parallel^m=C^m(\theta_{\gamma(i)})\,\dot\vartheta_i^D=A^m\dot{\theta}_i-B^m\varpi$ with $C^m\triangleq[\,A^m\ \ {-}B^m\,]$, and $\partial v_\parallel^m/\partial\dot{\varsigma}^{jh,j'h'}=-I$ for the plane-translation block of the pair carrying $m$. On the positive-weight branch $\alpha^m=\eta[\Lambda^i]^m>0$, hence $V_R^{jh,j'h'}(\theta_{\gamma(i)})>0$, the minimizing plane is unique and the coefficient maps $A^m,B^m,C^m$ are $C^1$ in $\theta_{\gamma(i)}$. Pointwise for any fixed coefficient representative, $C^m$ has full row rank, with the constant right inverse $\iota^m$ injecting $w\mapsto-w$ into that plane-translation block (all other velocity components zero),
\begin{ResizedEq}
&C^m\iota^m=I\\
&\|\iota^m\|_F=\sqrt{3}\\
&D_{\theta_{\gamma(i)}}\iota^m=0\\
&C^m(C^m)^T\succeq I.
\end{ResizedEq}
Consequently the damping-active velocity $\dot\vartheta^{D,m\sharp}\triangleq \dot\vartheta^D+\iota^m(v_\parallel^{m\sharp}-v_\parallel^m)$ realizes any target, $C^m \dot\vartheta^{D,m\sharp}=v_\parallel^{m\sharp}$, with
\begin{ResizedEq}
\label{eq:vertex-preimage-bound}
\|\dot\vartheta^{D,m\sharp}-\dot\vartheta^D\|\leq\|v_\parallel^{m\sharp}-v_\parallel^m\|,
\end{ResizedEq}
and the primitive force is recovered pointwise from any lifted force $q=(C^m)^T\zeta$ by the uniformly bounded inverse $\zeta=(C^m(C^m)^T)^{-1}C^m q$. On the positive-weight branch this recovery map is $C^1$ in $\theta_{\gamma(i)}$. Zero-weight incidences require neither chart inversion nor primitive-force recovery and are handled separately; no globally $C^1$ selection of inactive separating-plane coefficients is asserted. The construction may depend on $m$.
\end{lemma}
\begin{proof}[Proof of~\prettyref{lem:vertex-preimage}]
By~\prettyref{eq:tangent-velocity}, $v_\parallel^m=P^{jh,j'h'}\dot{X}^{m}-[\bar{n}^{jh,j'h'}\times X^{m}\dot{\omega}^{jh,j'h'}+\dot{\varsigma}^{jh,j'h'}]$, which at fixed previous configuration $\theta_{\gamma(i)}$ (hence fixed $P^{jh,j'h'},\bar{n}^{jh,j'h'},X^m$) is affine in $\dot\vartheta_i^D$, with Jacobians $A^m\triangleq\partial v_\parallel^m/\partial\dot{\theta}_i$ and $-B^m\triangleq\partial v_\parallel^m/\partial\varpi$, and in particular derivative $-I$ in the plane-translation velocity $\dot{\varsigma}^{jh,j'h'}$; thus $B^m=[\,I\ \ \bar{n}^{jh,j'h'}\times X^m\,]$ acts on the relevant pair block of $\varpi$. On the active side, the rigid-body vertex position and its velocity Jacobian are smooth in maximal coordinates, while the normalized plane normal $\bar n^{jh,j'h'}$ and tangent projector $P^{jh,j'h'}$ are $C^1$ by~\prettyref{lem:active-plane}. Hence the displayed formulas make $A^m$ and $B^m$ $C^1$ functions of $\theta_{\gamma(i)}$, and therefore $C^m=[\,A^m\ \ {-}B^m\,]$ is $C^1$ there.

The map $\iota^m$ injecting $w\mapsto-w$ into the three-dimensional $\dot{\varsigma}^{jh,j'h'}$ block satisfies $C^m\iota^m=-B^m(0,-w,0)=w$, i.e. $C^m\iota^m=I$. As a signed coordinate injection it has $\|\iota^m\|_F=\sqrt{3}$ and, separately, $\|\iota^m w\|_2=\|w\|_2$ for every $w\in\mathbb R^3$; being fixed, $D_{\theta_{\gamma(i)}}\iota^m=0$. Since $B^m(B^m)^T=I+(\bar{n}^{jh,j'h'}\times X^m)(\bar{n}^{jh,j'h'}\times X^m)^T\succeq I$, also $C^m(C^m)^T=A^m(A^m)^T+B^m(B^m)^T\succeq I$. Define the recovery map
\begin{ResizedEq}
R^m\triangleq(C^m(C^m)^T)^{-1}C^m.
\end{ResizedEq}
Then
\begin{ResizedEq}
R^m(R^m)^T=(C^m(C^m)^T)^{-1},
\end{ResizedEq}
and therefore
\begin{ResizedEq}
\|R^m\|_F^2=\operatorname{tr}\left((C^m(C^m)^T)^{-1}\right)\leq3,
\end{ResizedEq}
so $\|R^m\|_F\leq\sqrt{3}$. For a lifted force $q=(C^m)^T\zeta$, direct calculation gives $R^m q=\zeta$, while
\begin{ResizedEq}
\|q\|_2^2=\zeta^TC^m(C^m)^T\zeta\geq\|\zeta\|_2^2.
\end{ResizedEq}
Thus primitive-force recovery retains its constant-one vector estimate without using an induced matrix norm. Moreover, on the positive-weight branch $C^m$ is $C^1$, its Gram matrix remains positive definite, and matrix inversion is smooth on the positive-definite cone; consequently $R^m$ is $C^1$ there.

Finally $\dot\vartheta^{D,m\sharp}=\dot\vartheta^D+\iota^m(v_\parallel^{m\sharp}-v_\parallel^m)$ alters only the plane-translation block, by $-(v_\parallel^{m\sharp}-v_\parallel^m)$, so $C^m \dot\vartheta^{D,m\sharp}=v_\parallel^m+(v_\parallel^{m\sharp}-v_\parallel^m)=v_\parallel^{m\sharp}$. The coordinate-injection isometry gives $\|\dot\vartheta^{D,m\sharp}-\dot\vartheta^D\|=\|v_\parallel^{m\sharp}-v_\parallel^m\|$, proving~\prettyref{eq:vertex-preimage-bound}.
\end{proof}

\begin{lemma}[Chart and exact Clarke-normal pullback]
\label{lem:exact-pullback}
Fix a vertex term $m$, write $x=\theta_{\gamma(i)}$, $\Damp^{V,m}(x,\dot\vartheta^D)=\eta[\Lambda^i]^m(x)\,\|C^m(x)\dot\vartheta^D\|$; $\alpha\triangleq\eta[\Lambda^i]^m$, and split $\dot\vartheta^D=(\widehat{\dot\vartheta^D},u)$ with $u=\dot{\varsigma}^{jh,j'h'}$ the plane-translation block, so $C^m\dot\vartheta^D=\widehat A^m\widehat{\dot\vartheta^D}-u$ with a constant $-I$ block in $u$ (\prettyref{lem:vertex-preimage}). On the active side $\alpha(x)>0$, define the local force-range set
\begin{equation}
\mathcal M^m\triangleq\{(x,\dot\vartheta^D,q):q\in\operatorname{range}((C^m(x))^T)\}.
\end{equation}
It is locally a $C^1$ embedded submanifold of the ambient $(x,\dot\vartheta^D,q)$-space, and the map
\begin{ResizedEq}
&\Xi^m:(x,\widehat{\dot\vartheta^D},v,\zeta)\ \longmapsto\ \big(x,\ \dot\vartheta^D,\ (C^m(x))^T\zeta\big)\\
&\dot\vartheta^D=\TWOC{\widehat{\dot\vartheta^D}}{\widehat A^m\widehat{\dot\vartheta^D}-v},
\end{ResizedEq}
is locally a $C^1$ diffeomorphism from its source-coordinate neighborhood onto an open subset of $\mathcal M^m$. Its restriction carries $\mathcal{S}^m\triangleq\{(x,\widehat{\dot\vartheta^D},v,\zeta):(\alpha(x),v,\zeta)\in\widehat{\mathcal{G}}\}$ exactly onto $\gph\partial_{\dot\vartheta^D}\Damp^{V,m}$. For any master normal $(\tau,a,b)\in N^C_{\widehat{\mathcal{G}}}(\alpha(x),v,\zeta)$ and any $s$ with $C^m(x)s=b$,
\begin{ResizedEq}
\label{eq:np-pullback}
n_q=&s\\
n_{\dot\vartheta^D}=&(C^m)^Ta\\
n_x=&D_x\alpha^T\tau+D_x[C^m\dot\vartheta^D]^Ta-D_x[(C^m)^T\zeta]^Ts
\end{ResizedEq}
(with $\dot\vartheta^D,\zeta$ held fixed in the two $x$-derivatives) defines an ambient Clarke normal $(n_x,n_{\dot\vartheta^D},n_q)\in N^C_{\gph\partial_{\dot\vartheta^D}\Damp^{V,m}}(x,\dot\vartheta^D,(C^m)^T\zeta)$. The chart and pullback are constructed on the active side $\alpha(x)>0$; they remain valid at active points where $D_x\alpha(x)=0$, in which case the weight-pullback term $D_x\alpha^T\tau$ vanishes. The zero-weight case $\alpha(x)=0$, for which $D_x\alpha(x)=0$ and the weighted friction term is handled separately, requires no active-side coefficient chart.
\end{lemma}
\begin{proof}[Proof of~\prettyref{lem:exact-pullback}]
Work in a positive-weight neighborhood on which $C^m$ is $C^1$. By~\prettyref{lem:vertex-preimage}, $C^m$ has constant full row rank and
\begin{equation}
R^m(x)\triangleq(C^m(C^m)^T)^{-1}C^m
\end{equation}
is $C^1$. Hence
\begin{equation}
\Phi^m(x,\dot\vartheta^D,\zeta)\triangleq(x,\dot\vartheta^D,(C^m(x))^T\zeta)
\end{equation}
is injective and has the $C^1$ inverse $(x,\dot\vartheta^D,q)\mapsto(x,\dot\vartheta^D,R^m(x)q)$ on its image. Its differential is injective (the first two blocks determine their variations and $(C^m)^T\delta\zeta=0$ implies $\delta\zeta=0$), so after restricting the neighborhood, $\Phi^m$ is a $C^1$ embedding with image $\mathcal M^m$. This proves the asserted embedded-manifold property without selecting coefficients outside the active neighborhood.

The chain rule for the weighted norm composed with the surjective linear map $C^m$ gives $\partial_{\dot\vartheta^D}\Damp^{V,m}(x,\cdot)=(C^m)^T\big(\alpha(x)\partial\|\cdot\|\big)(C^m\dot\vartheta^D)$. Thus $(x,\dot\vartheta^D,q)$ belongs to the damping graph iff $q=(C^m)^T\zeta$ and $(\alpha(x),C^m\dot\vartheta^D,\zeta)\in\widehat{\mathcal{G}}$. The map $\Xi^m$ takes its full source-coordinate neighborhood into $\mathcal M^m$ and inverts there by reading off $x,\widehat{\dot\vartheta^D}$, setting $v=C^m\dot\vartheta^D$, and recovering $\zeta=R^m q$. Both maps are $C^1$, so $\Xi^m$ is the claimed relative diffeomorphism, and the preceding graph characterization proves $\Xi^m(\mathcal S^m)=\gph\partial_{\dot\vartheta^D}\Damp^{V,m}$ locally.

It remains to establish the source normal without invoking a general Clarke-tangent chain rule. We verify directly on the three positive-weight strata that
\begin{ResizedEq}
&(D_x\alpha^T\tau,0,a,b)\in N^C_{\mathcal S^m}\\
&\forall(\tau,a,b)\in N^C_{\widehat{\mathcal G}}(\alpha(x),v,\zeta).
\end{ResizedEq}
At a strict-stick point, $v=0$ and $\|\zeta\|<\alpha(x)$. Locally $\mathcal S^m=\{v=0\}$ with $x,\widehat{\dot\vartheta^D},\zeta$ free, while~\prettyref{lem:master-weighted-norm-clarke}\ref{it:mwnc-strict} gives $(\tau,a,b)=(0,a,0)$. Hence $(D_x\alpha^T\tau,0,a,b)=(0,0,a,0)$ is a Clarke normal to $\mathcal S^m$.

At a sliding point, write $v=re$, $r>0$, and $\zeta=\alpha(x)e$. Locally $\mathcal S^m$ is the smooth manifold
\begin{equation}
\mathfrak{F}(x,v,\zeta)\triangleq\zeta-\alpha(x)\frac{v}{\|v\|}=0,
\end{equation}
with $\widehat{\dot\vartheta^D}$ free; it is regular because $D_\zeta \mathfrak{F}=I$. By~\prettyref{lem:master-weighted-norm-clarke}\ref{it:mwnc-slide}, for some $z$ one has $(\tau,a,b)=(e^Tz,\frac{\alpha}{r}P_ez,-z)$. Direct differentiation gives
\begin{equation}
D \mathfrak{F}(x,v,\zeta)^T(-z)=\left(D_x\alpha^T\tau,\ a,\ b\right),
\end{equation}
so $(D_x\alpha^T\tau,0,a,b)$ is a regular, hence Clarke, normal to $\mathcal S^m$.

Finally, at a stick-boundary point $v=0$, $\zeta=\alpha(x)e$,~\prettyref{lem:master-weighted-norm-clarke}\ref{it:mwnc-boundary} gives $(\tau,a,b)=(t,a,-te)$. Strict-stick source points $(x,\widehat{\dot\vartheta^D},0,(\alpha(x)-\varepsilon_k)e)$ with $\varepsilon_k\downarrow0$ supply the limiting normals $(0,0,a,0)$. Sliding source points $(x,\widehat{\dot\vartheta^D},r_ke,\alpha(x)e)$ with $r_k\downarrow0$, using $z=te$ in the sliding formula and $P_e(te)=0$, supply the limiting normals $(D_x\alpha^Tt,0,0,-te)$. Since the Clarke normal cone is the closed convex cone generated by limiting normals, it contains their sum $(D_x\alpha^Tt,0,a,-te)$. Thus the required source-normal membership holds on every positive-weight stratum, including active points where $D_x\alpha=0$.

Applying tangent functoriality to the relative $C^1$ diffeomorphism $\Xi^m$ and to its relative inverse gives
\begin{equation}
D\Xi^m\,T^C_{\mathcal S^m}=T^C_{\gph\partial_{\dot\vartheta^D}\Damp^{V,m}},
\end{equation}
where both sides are ambient tangent vectors lying in $T\mathcal M^m$. For a source tangent variation, its image obeys $C^m\delta\dot\vartheta^D=\delta v-(D_xC^m)\dot\vartheta^D\,\delta x$ and $\delta q=(D_x(C^m)^T\zeta)\delta x+(C^m)^T\delta\zeta$. Pairing~\prettyref{eq:np-pullback} with the image variation,
\begin{ResizedEq}
\label{eq:np-pairing}
&\langle n_x,\delta x\rangle+\langle n_{\dot\vartheta^D},\delta\dot\vartheta^D\rangle+\langle n_q,\delta q\rangle\\
=&\langle n_x,\delta x\rangle+\langle a,C^m\delta\dot\vartheta^D\rangle+\langle s,(D_x(C^m)^T\zeta)\delta x\rangle+\langle C^m s,\delta\zeta\rangle\\
=&\langle D_x\alpha^T\tau,\delta x\rangle+\langle a,\delta v\rangle+\langle b,\delta\zeta\rangle=\langle(\tau,a,b),(D_x\alpha\,\delta x,\delta v,\delta\zeta)\rangle,
\end{ResizedEq}
using $C^m s=b$ and the definition of $n_x$. Thus the displayed pullback pairs nonpositively with every graph Clarke tangent whenever its force block $s$ solves $C^m s=b$ and the corresponding pairing is represented by a source tangent.

For completeness, we now verify ambient polarity and the arbitrary-$s$ assertion explicitly. Set
\begin{equation}
s_0\triangleq(C^m)^T(C^m(C^m)^T)^{-1}b,
\qquad C^m s_0=b.
\end{equation}
The pullback with $s=s_0$ pairs nonpositively with every graph Clarke tangent by the relative tangent equality and~\prettyref{eq:np-pairing}. The tangent space of $\mathcal M^m$, obtained by differentiating $\Phi^m$, consists of
\begin{equation}
(\delta x,\delta\dot\vartheta^D,(D_x(C^m)^T\zeta)\delta x+(C^m)^T\delta\zeta).
\end{equation}
For every $k\in\ker C^m$, the ambient covector
\begin{equation}
\kappa(k)\triangleq\big(-D_x[(C^m)^T\zeta]^Tk,\ 0,\ k\big)
\end{equation}
pairs to zero with each such tangent: the two moving-range derivative terms cancel and the remainder is $\langle C^m k,\delta\zeta\rangle=0$. Finally, any $s$ satisfying $C^m s=b$ decomposes as $s=s_0+k$ with $k\in\ker C^m$. Component by component, replacing $s_0$ by $s$ in~\prettyref{eq:np-pullback} adds exactly $\kappa(k)$: it adds $k$ to the force block, nothing to the velocity block, and $-D_x[(C^m)^T\zeta]^Tk$ to the $x$-block. Consequently the pullback for arbitrary $s$ has the same nonpositive pairing with every graph Clarke tangent as the canonical pullback. Direct Clarke polarity, rather than an ambient normal-transport or normal-cone sum rule, proves the claimed ambient normality. The argument is uniform across the positive-weight strata, and at a degenerate active point $D_x\alpha=0$ the weight term simply vanishes.
\end{proof}

Physically meaningful Clarke stationarity is per contact vertex and lives in the complete damping-active physical-plus-plane velocities $\dot\vartheta_i^D$, whereas~\prettyref{alg:CITO} solves the modified full-trajectory~\prettyref{pb:CITO-M} instantiated with the reduced CRDS~\prettyref{pb:PSR-CRDS}. We therefore work throughout in the damping-active velocity $\dot\vartheta^D$, keeping the plane-velocity damper $\EpsDamp\|\varpi\|^2$ exact and external: it is never smoothed away, never replaced by a graph normal, and never merged into the friction term set $\mathcal{T}$. We first record the full velocity Hessian and its common Schur response, then reconstruct the plane velocities and lift the reduced multiplier, identify the smooth full per-vertex normal and transfer it to a nearby full Clarke normal, and only afterwards normalize the complete full certificate.

We record the reduced-block setup used by the per-vertex normal transfer and the Schur lift. By the envelope identity of~\prettyref{lem:smooth-DRp} (Case~II), on the active side the reduced velocity gradient is $G(\theta_{\gamma(i)},\dot{\theta}_i)=\nabla_{\dot{\theta}_i}\Damp_{\EpsFri,\LHamilton}^V(\theta_{\gamma(i)},\dot{\theta}_i,\varpi^\star)$ with $\varpi^\star$ the unique plane-velocity minimizer, so, writing each vertex velocity as $v_\parallel^m=A^m\dot{\theta}_i-B^m\varpi$ (\prettyref{lem:vertex-preimage}) with $A^m,B^m$ the $C^1$ Jacobians of~\prettyref{eq:tangent-velocity},
\begin{ResizedEq}
\label{eq:abd-G}
G=\eta\sum_{m\in\mathcal{T}}[\Lambda^i]^m\,(A^m)^T\nabla\psi_\mu(v_\parallel^m)\Big|_{\varpi=\varpi^\star}.
\end{ResizedEq}
Assemble the friction velocity Hessian over $(\dot{\theta}_i,\varpi)$: with $J^m\triangleq[\,A^m\ \ {-}B^m\,]$,
\begin{ResizedEq}
\label{eq:abd-Psi}
\Psi_i\triangleq\eta\sum_{m}[\Lambda^i]^m (J^m)^T\nabla^2\psi_\mu(v_\parallel^m)J^m+\diag(0,2\EpsDamp I),
\end{ResizedEq}
with sub-blocks $\Psi_i^{\theta\theta},\Psi_i^{\theta\varpi},\Psi_i^{\varpi\varpi}$ and $\Psi_i^{\varpi\varpi}\succeq2\EpsDamp I$ (\prettyref{lem:boundary-flatness}). Differentiating the reduced value through the envelope (implicit-function theorem, licit since $\Psi_i^{\varpi\varpi}\succ0$) identifies the reduced blocks with Schur complements of $\Psi_i$ over the $\varpi$-block: with the $\varpi$-response $\zeta^\star(\xi)\triangleq-(\Psi_i^{\varpi\varpi})^{-1}\Psi_i^{\varpi\theta}\xi$,
\begin{ResizedEq}
\label{eq:abd-schur-blocks}
\delta_i K_i=&\Psi_i^{\theta\theta}-\Psi_i^{\theta\varpi}(\Psi_i^{\varpi\varpi})^{-1}\Psi_i^{\varpi\theta}\\
\Gamma_i=&\Psi_i^{\theta\gamma}-\Psi_i^{\theta\varpi}(\Psi_i^{\varpi\varpi})^{-1}\Psi_i^{\varpi\gamma},
\end{ResizedEq}
where the superscript $\gamma$ denotes the explicit $\theta_{\gamma(i)}$-partial (velocity and $\varpi$ held fixed); the reduced damping energy has the variational form
\begin{ResizedEq}
\label{eq:abd-schur-var}
&\delta_i\,\xi^T K_i\xi=\min_{\zeta}\Big\{\eta\sum_{m}[\Lambda^i]^m w_m(\zeta)^T\nabla^2\psi_\mu(v_\parallel^m)w_m(\zeta)+2\EpsDamp\|\zeta\|^2\Big\}\\
&w_m(\zeta)\triangleq A^m\xi-B^m\zeta,
\end{ResizedEq}
attained at $\zeta=\zeta^\star(\xi)$; write $w_m\triangleq w_m(\zeta^\star)$.

Two structural facts of~\prettyref{eq:abd-Psi} are used repeatedly, and neither yields a vertexwise decomposition of the reduced dissipation. First, the full velocity Hessian is an honest sum of positive-semidefinite pieces plus the external damper,
\begin{ResizedEq}
\label{eq:abd-Psi-split}
&\Psi_i=\sum_{m\in\mathcal{T}}\Psi_i^m+\diag(0,2\EpsDamp I),\\
&\Psi_i^m\triangleq\eta[\Lambda^i]^m\,(J^m)^T\nabla^2\psi_\mu(v_\parallel^m)J^m\succeq0,
\end{ResizedEq}
in which $\diag(0,2\EpsDamp I)$ is a separate external term, not one of the friction pieces $\Psi_i^m$. Second, a reduced direction $\xi$ lifts to the full velocity through the common aggregate Schur response
\begin{ResizedEq}
\label{eq:abd-common-lift}
\widehat\xi\triangleq(\xi,\zeta^\star(\xi)),\qquad \zeta^\star(\xi)=-(\Psi_i^{\varpi\varpi})^{-1}\Psi_i^{\varpi\theta}\xi,
\end{ResizedEq}
which depends on the entire $\Psi_i$---all friction pieces and the damper---not on any individual $\Psi_i^m$. Reading off the $\dot{\theta}_i$- and earlier-configuration components of the action $\Psi_i^m\widehat\xi$ of each piece on this common lift defines the per-term blocks used as the target of the vertexwise transfer,
\begin{ResizedEq}
\label{eq:abd-perterm-blocks}
\delta_i K_i^m\xi\triangleq&\big[\Psi_i^m\widehat\xi\big]_{\dot{\theta}_i}=(A^m)^T\eta[\Lambda^i]^m\nabla^2\psi_\mu(v_\parallel^m)\,w_m,\\
(\Gamma_i^m)^T\xi\triangleq&\Big[D_x\big((C^m)^T\eta[\Lambda^i]^m\nabla\psi_\mu(C^m\dot\vartheta^D)\big)^T\widehat\xi\Big]_{\theta_{\gamma(i)}},
\end{ResizedEq}
with $w_m=C^m\widehat\xi=A^m\xi-B^m\zeta^\star$ and $C^m=[\,A^m\ \ {-}B^m\,]$. We emphasize what $K_i^m,\Gamma_i^m$ are not: an individual $\Psi_i^m$ need not have an invertible plane block, so they are not individual Schur complements of $\Psi_i^m$; they are not Hessians of any reduced per-vertex scalar (the reduced dissipation does not decompose vertexwise,~\prettyref{eq:abd-schur-var}); and they need not be positive-semidefinite. The reduced energy does not decompose vertexwise: evaluating~\prettyref{eq:abd-schur-var} at its minimizer,
\begin{ResizedEq}
\label{eq:abd-energy-split}
&\delta_i\,\xi^T K_i\xi\\
=&\sum_{m\in\mathcal{T}}\eta[\Lambda^i]^m\,w_m^T\nabla^2\psi_\mu(v_\parallel^m)w_m+2\EpsDamp\|\zeta^\star(\xi)\|^2,
\end{ResizedEq}
carries the extra external-damper term $2\EpsDamp\|\zeta^\star(\xi)\|^2$ beyond the friction contributions. The $\dot{\theta}_i$-block action is nonetheless additive: since $\Psi_i\widehat\xi=(\delta_i K_i\xi,0)$ (\prettyref{lem:schur-lift} below) while $\diag(0,2\EpsDamp I)\widehat\xi=(0,2\EpsDamp\zeta^\star)$, subtracting gives $\sum_{m}[\Psi_i^m\widehat\xi]_{\dot{\theta}_i}=\delta_i K_i\xi$, i.e.\ $\sum_m\delta_i K_i^m\xi=\delta_i K_i\xi$, so $K_i=\sum_m K_i^m$ as operators---all that the exact adjoint lift of~\prettyref{lem:abd-kkt-lift} requires. What fails is the energy/positive-semidefinite decomposition: the genuine per-vertex energies $\widehat\xi^T\Psi_i^m\widehat\xi=\eta[\Lambda^i]^m\,w_m^T\nabla^2\psi_\mu(v_\parallel^m)w_m\ge0$ are not the quadratic forms $\xi^T\delta_i K_i^m\xi$ of the (non-symmetric, non-PSD) operator blocks, and they sum to $\delta_i\xi^T K_i\xi-2\EpsDamp\|\zeta^\star\|^2\le\delta_i\xi^T K_i\xi$ (strict when $\zeta^\star(\xi)\neq0$), so the reduced energy cannot be bounded by a Cauchy-Schwarz over positive-semidefinite summands $K_i^m$---the flaw of the reduced route that this specialization avoids. Each per-vertex energy is one nonnegative summand of~\prettyref{eq:abd-energy-split}, so
\begin{ResizedEq}
\label{eq:abd-perterm-energy}
&0\le\widehat\xi^T\Psi_i^m\widehat\xi\le\delta_i\,\xi^T K_i\xi\\
&2\EpsDamp\|\zeta^\star(\xi)\|^2\le\delta_i\,\xi^T K_i\xi,
\end{ResizedEq}
the two bounds the compatibility lemma consumes.

\begin{lemma}[Unique plane-velocity reconstruction]
\label{lem:varpi-reconstruction}
At every returned timestep $i$, the predecessor satisfies $V_R(\theta_{\gamma(i)},\EpsFJ)\leq\LHamilton$, so the globally extended reduced damping agrees with the reduced local pseudo-Huber core. Let $\widehat\Damp_{\mu_{\EpsFri,\LHamilton}}^V(\theta_{\gamma(i)},\dot\theta_i,\varpi)$ denote the corresponding unreduced friction-plus-plane-damper objective. Its minimization over the concatenated separating-plane velocities $\varpi=(\varpi^{jh,j'h'})$ has a globally unique minimizer $\varpi^\star$, characterized by
\begin{equation}
\nabla_\varpi\widehat\Damp_{\mu_{\EpsFri,\LHamilton}}^V(\theta_{\gamma(i)},\dot\theta_i,\varpi^\star)=0.
\end{equation}
The minimizer depends locally $C^1$ on $(\theta_{\gamma(i)},\dot\theta_i)$ on each smooth coefficient region, matches the zero inactive branch across activation boundaries, and satisfies the envelope identity
\begin{ResizedEq}
\nabla_{\dot\theta_i}\Damp_{R,\EpsFri,\LHamilton}^V
=\nabla_{\dot\theta_i}\widehat\Damp_{\mu_{\EpsFri,\LHamilton}}^V(\theta_{\gamma(i)},\dot\theta_i,\varpi^\star).
\end{ResizedEq}
The minimization and minimizer decouple over convex-hull pairs. Here $\varpi^\star$ reconstructs the auxiliary block of the damping-active velocity $(\dot\theta_i,\varpi^\star)$; the complete finite-difference velocity additionally contains the separate $\dot p_i$ block.
\end{lemma}
\begin{proof}[Proof of~\prettyref{lem:varpi-reconstruction}]
Abbreviate $\mu=\mu_{\EpsFri,\LHamilton}$ and
\begin{equation}
\Psi(x,v,\varpi)\triangleq\widehat\Damp_{\mu}^V(x,v,\varpi),
\qquad x=\theta_{\gamma(i)},\quad v=\dot\theta_i.
\end{equation}
On the certified sublevel $V_R(x,\EpsFJ)\leq\LHamilton$, the mollifier in~\prettyref{eq:global-reduced-damping} is identically one, so~\prettyref{eq:reduced-damping-cutoff} and the local construction~\prettyref{lem:smooth-DRp} give
\begin{equation}
\Damp_{R,\EpsFri,\LHamilton}^V(x,v)=\min_\varpi\Psi(x,v,\varpi).
\end{equation}
Every pseudo-Huber friction term is nonnegative, while the plane-velocity damper is retained exactly; hence
\begin{equation}
\Psi(x,v,\varpi)\geq\EpsDamp\|\varpi\|^2.
\end{equation}
Thus $\Psi(x,v,\cdot)$ is coercive and attains a finite minimum. Moreover, with $B^m$ the velocity Jacobian of the $m$th relative-velocity map with respect to $\varpi$,
\begin{ResizedEq}
\nabla_{\varpi\varpi}^2\Psi
=&\eta\sum_m[\Lambda^i]^m(B^m)^T\nabla^2\psi_\mu(v_\parallel^m)B^m+2\EpsDamp I\\
\succeq&2\EpsDamp I,
\end{ResizedEq}
because the pseudo-Huber Hessians and friction weights are positive semidefinite and nonnegative, respectively. Therefore $\Psi(x,v,\cdot)$ is globally $2\EpsDamp$-strongly convex, its minimizer
\begin{equation}
\varpi^\star(x,v)\triangleq\argmin{\varpi}\Psi(x,v,\varpi)
\end{equation}
is globally unique, and unconstrained Fermat stationarity gives $\nabla_\varpi\Psi(x,v,\varpi^\star)=0$. The same Hessian bound yields the uniform inverse estimate
\begin{equation}
\left\|\left[\nabla_{\varpi\varpi}^2\Psi(x,v,\varpi)\right]^{-1}\right\|
\leq\frac{1}{2\EpsDamp}.
\end{equation}

On an active smooth coefficient region, the stationarity map $\Xi(x,v,\varpi)\triangleq\nabla_\varpi\Psi(x,v,\varpi)$ is jointly $C^1$: the optimal separating-plane data and contact weights are $C^1$ by~\prettyref{lem:active-plane}, every relative velocity is affine in $(v,\varpi)$ with those coefficients, and $\psi_\mu$ is smooth. Since $D_\varpi\Xi=\nabla_{\varpi\varpi}^2\Psi$ is invertible, the $C^1$ implicit function theorem gives local $C^1$ dependence of $\varpi^\star$ on $(x,v)$. If all forces vanish, $\Psi=\EpsDamp\|\varpi\|^2$ and $\varpi^\star=0$. At an activation boundary, the derivative estimates in the proof of~\prettyref{lem:boundary-flatness} give $D_x\varpi^\star=\mathcal{O}(\rho)$ and $D_v\varpi^\star=\mathcal{O}(\rho^3)$; these tend to zero and match the inactive branch. Thus the required differentiability holds throughout the certified sublevel.

The attained minimum satisfies
\begin{equation}
\Damp_{R,\EpsFri,\LHamilton}^V(x,v)=\Psi(x,v,\varpi^\star(x,v)).
\end{equation}
Differentiating with respect to $v=\dot\theta_i$ is now justified and gives
\begin{ResizedEq}
\nabla_v\Damp_{R,\EpsFri,\LHamilton}^V
=&\nabla_v\Psi(x,v,\varpi^\star)+D_v\varpi^\star(x,v)^T\nabla_\varpi\Psi(x,v,\varpi^\star)\\
=&\nabla_v\Psi(x,v,\varpi^\star),
\end{ResizedEq}
where the chain-rule term vanishes by inner stationarity. Finally, $\Psi=\sum_{jh,j'h'}\widehat\Damp_{\mu}^{V,jh,j'h'}$ and distinct pairs share no plane-velocity block, so both the minimization and the unique minimizer decouple pairwise, $\varpi^\star=(\varpi^{\star,jh,j'h'})$.
\end{proof}

\begin{lemma}[Schur-complement multiplier lift]
\label{lem:schur-lift}
Let $\Psi_i=\begin{pmatrix}\Psi_i^{\theta\theta}&\Psi_i^{\theta\varpi}\\\Psi_i^{\varpi\theta}&\Psi_i^{\varpi\varpi}\end{pmatrix}$ of~\prettyref{eq:abd-Psi} be the full velocity Hessian at the reconstructed minimizer, with $\Psi_i^{\varpi\varpi}\succeq2\EpsDamp I$. For a reduced physics multiplier $\mu_{\theta_i}$ set $\nu_i\triangleq-(\Psi_i^{\varpi\varpi})^{-1}\Psi_i^{\varpi\theta}\mu_{\theta_i}$ and $\widehat\mu_i\triangleq(\mu_{\theta_i},\nu_i)$. Then, with $S_i\triangleq\Psi_i^{\theta\theta}-\Psi_i^{\theta\varpi}(\Psi_i^{\varpi\varpi})^{-1}\Psi_i^{\varpi\theta}=\delta_i K_i$ the Schur complement (the reduced diagonal damping block up to the common $\delta_i$ scaling),
\begin{ResizedEq}
\label{eq:schur-identities}
&\Psi_i\widehat\mu_i=\TWOC{S_i\mu_{\theta_i}}{0}\\
&\widehat\mu_i^T\Psi_i\widehat\mu_i=\mu_{\theta_i}^T S_i\mu_{\theta_i}\\
&2\EpsDamp\|\nu_i\|^2\leq\mu_{\theta_i}^T S_i\mu_{\theta_i}.
\end{ResizedEq}
\end{lemma}
\begin{proof}[Proof of~\prettyref{lem:schur-lift}]
The lower block of $\Psi_i\widehat\mu_i$ is $\Psi_i^{\varpi\theta}\mu_{\theta_i}+\Psi_i^{\varpi\varpi}\nu_i=\Psi_i^{\varpi\theta}\mu_{\theta_i}-\Psi_i^{\varpi\varpi}(\Psi_i^{\varpi\varpi})^{-1}\Psi_i^{\varpi\theta}\mu_{\theta_i}=0$, and the upper block is $\Psi_i^{\theta\theta}\mu_{\theta_i}+\Psi_i^{\theta\varpi}\nu_i=(\Psi_i^{\theta\theta}-\Psi_i^{\theta\varpi}(\Psi_i^{\varpi\varpi})^{-1}\Psi_i^{\varpi\theta})\mu_{\theta_i}=S_i\mu_{\theta_i}$, proving the first identity; pairing with $\widehat\mu_i$ gives $\widehat\mu_i^T\Psi_i\widehat\mu_i=\mu_{\theta_i}^T S_i\mu_{\theta_i}$. For the third, decompose $\Psi_i=G_i+\diag(0,2\EpsDamp I)$ with the friction Gram part $G_i\triangleq\eta\sum_m[\Lambda^i]^m(J^m)^T\nabla^2\psi_\mu J^m\succeq0$; then, using the second identity, $2\EpsDamp\|\nu_i\|^2=\widehat\mu_i^T\diag(0,2\EpsDamp I)\widehat\mu_i=\widehat\mu_i^T\Psi_i\widehat\mu_i-\widehat\mu_i^T G_i\widehat\mu_i=\mu_{\theta_i}^T S_i\mu_{\theta_i}-\widehat\mu_i^T G_i\widehat\mu_i\leq\mu_{\theta_i}^T S_i\mu_{\theta_i}$, the last step by $\widehat\mu_i^T G_i\widehat\mu_i\geq0$. The bound uses the block-diagonal placement of the damper $2\EpsDamp I$ on top of the positive-semidefinite $G_i$; the abstract inequality $\Psi_i^{\varpi\varpi}\succeq2\EpsDamp I$ alone would not suffice, and the friction Gram contribution to $\Psi_i^{\varpi\varpi}$---though it may be large as $\mu\to0$---is harmless because it enters only through the nonnegative $\widehat\mu_i^T G_i\widehat\mu_i$.
\end{proof}

\begin{lemma}[Identification of the smooth per-vertex normal]
\label{lem:smooth-normal-identification}
Fix a vertex term $m$ and a reduced direction $\xi$, lifted to the full velocity direction $\xi^{\mathrm f}\triangleq(\xi,\zeta^\star(\xi))$ by the plane-velocity Schur response of~\prettyref{eq:abd-schur-blocks}, so $z\triangleq C^m\xi^{\mathrm f}=A^m\xi-B^m\zeta^\star=w_m$. The smooth master normal $(g_\mu^Tz,\ \alpha\nabla^2\psi_\mu(v)z,\ -z)$ (with $\alpha=\eta[\Lambda^i]^m$, $v=v_\parallel^m$, $g_\mu=\nabla\psi_\mu(v)$), generated by direct differentiation of the smooth per-term graph $\widehat{\mathcal{G}}_\mu^m=\{(x,\dot\vartheta^D,q):q=(C^m(x))^T\alpha(x)\nabla\psi_\mu(C^m(x)\dot\vartheta^D)\}$ with force block $-\xi^{\mathrm f}$, is the smooth full-velocity normal
\begin{ResizedEq}
\label{eq:smooth-full-normal}
n_\mu=\Big(D_x\big[(C^m)^T\alpha\,\nabla\psi_\mu(C^m\dot\vartheta^D)\big]^T\xi^{\mathrm f},\ (C^m)^T\alpha\nabla^2\psi_\mu(v)C^m\xi^{\mathrm f},\ -\xi^{\mathrm f}\Big),
\end{ResizedEq}
which we denote $N_{i,m}^{\mathrm{sm}}(\xi^{\mathrm f})\triangleq n_\mu$, the smooth full per-vertex normal targeted by~\prettyref{def:S-D-compatible-ABD}; its evaluation on the common aggregate lift $\xi^{\mathrm f}=\widehat\xi$ followed by block projection onto the reduced $\dot{\theta}_i$-coordinate recovers the per-term blocks~\prettyref{eq:abd-perterm-blocks},
\begin{ResizedEq}
\label{eq:smooth-normal-id}
n_\mu\ \longmapsto\ \big((\Gamma_i^m)^T\xi,\ \delta_i K_i^m\xi,\ -\xi\big).
\end{ResizedEq}
\end{lemma}
\begin{proof}[Proof of~\prettyref{lem:smooth-normal-identification}]
On each positive-weight region the smooth per-term dissipation $\eta[\Lambda^i]^m\psi_\mu(C^m\dot\vartheta^D)$ has velocity gradient $(C^m)^T\alpha\,\nabla\psi_\mu(C^m\dot\vartheta^D)$ (the summand of~\prettyref{eq:abd-G} under the envelope identity of~\prettyref{lem:smooth-DRp}, Case~II), so its force graph is the smooth manifold $\widehat{\mathcal{G}}_\mu^m\triangleq\{(x,\dot\vartheta^D,q):q=Q(x,\dot\vartheta^D)\}$ with $Q(x,\dot\vartheta^D)\triangleq(C^m(x))^T\alpha(x)\,\nabla\psi_\mu(C^m(x)\dot\vartheta^D)$. Across a zero-weight boundary,~\prettyref{lem:boundary-flatness} and~\ref{lem:pervertex-activation} give a $C^1$ zero extension of this weighted force map and its first derivative; no globally $C^1$ selection of an inactive separating plane or of $C^m$ itself is asserted. We identify $n_\mu$ by differentiating this smooth graph directly, rather than by the exact-graph pullback technique: when $v\neq0$ the base point $(\alpha,v,\alpha\nabla\psi_\mu(v))$ of the smooth master normal has force norm $\alpha\|v\|/s<\alpha$, so it does not lie on the exact weighted-norm graph $\widehat{\mathcal{G}}=\{q\in\alpha\partial\|v\|\}$ and the exact-graph pullback technique does not apply to it verbatim. On the smooth manifold $\widehat{\mathcal{G}}_\mu^m$ the regular, limiting, and Clarke normal cones coincide; its tangent space at $(x,\dot\vartheta^D,Q(x,\dot\vartheta^D))$ is $\{(\delta x,\delta\dot\vartheta^D,D_xQ\,\delta x+D_{\dot\vartheta^D}Q\,\delta\dot\vartheta^D)\}$, so a covector $(n_x,n_{\dot\vartheta^D},n_q)$ is Clarke normal iff $n_x=-(D_xQ)^Tn_q$ and $n_{\dot\vartheta^D}=-(D_{\dot\vartheta^D}Q)^Tn_q$. Taking the force block $n_q=-\xi^{\mathrm f}$ (so that $z=C^m\xi^{\mathrm f}=w_m$) gives $n_x=(D_xQ)^T\xi^{\mathrm f}$ and $n_{\dot\vartheta^D}=(D_{\dot\vartheta^D}Q)^T\xi^{\mathrm f}$. Since $D_{\dot\vartheta^D}Q=(C^m)^T\alpha\nabla^2\psi_\mu(v)C^m$, the velocity block is $(C^m)^T\alpha\nabla^2\psi_\mu(v)C^m\xi^{\mathrm f}$; and the product rule $D_x[(C^m)^T\alpha g_\mu]=D_x(C^m)^T\alpha g_\mu+(C^m)^T D_x\alpha\,g_\mu+(C^m)^T\alpha\,\nabla^2\psi_\mu(v)(D_xC^m)\dot\vartheta^D$ (with $g_\mu=\nabla\psi_\mu(v)$; the middle term retains the $D_x\alpha$ contribution even where $\alpha=0$) gives the $x$-block $D_x[(C^m)^T\alpha\nabla\psi_\mu(C^m\dot\vartheta^D)]^T\xi^{\mathrm f}$. Together these are~\prettyref{eq:smooth-full-normal}.

As an energy-consistency check for~\prettyref{eq:smooth-normal-id}, contract the velocity block of~\prettyref{eq:smooth-full-normal} with $\xi^{\mathrm f}$: since $z=C^m\xi^{\mathrm f}=w_m$,
\begin{ResizedEq}
(\xi^{\mathrm f})^T(C^m)^T\alpha\nabla^2\psi_\mu(v)C^m\xi^{\mathrm f}=\alpha\,w_m^T\nabla^2\psi_\mu(v)w_m,
\end{ResizedEq}
which is the $m$-th friction summand of the reduced damping energy~\prettyref{eq:abd-schur-var} evaluated at its aggregate minimizer $\zeta^\star$. This scalar identity is only the energy-consistency check; the following block calculation proves the vector projection identity. The reduction~\prettyref{eq:smooth-normal-id} is not a Schur complementation of the individual per-term Hessian (which need not have an invertible plane block) and is not the second derivative of any reduced per-vertex scalar; it is the evaluation of the full per-term velocity Hessian $(C^m)^T\alpha\nabla^2\psi_\mu(v)C^m=(J^m)^T\alpha\nabla^2\psi_\mu(v)J^m=\Psi_i^m$ on the common aggregate lift $\xi^{\mathrm f}=\widehat\xi$ followed by block projection: by the definition~\prettyref{eq:abd-perterm-blocks}, the $\dot{\theta}_i$-component of the velocity block $\Psi_i^m\xi^{\mathrm f}$ is $\delta_i K_i^m\xi$ and the earlier-configuration component of the $x$-block is $(\Gamma_i^m)^T\xi$, while the force block $-\xi^{\mathrm f}$ projects to $-\xi$ in the reduced coordinate. This is~\prettyref{eq:smooth-normal-id}.
\end{proof}

\begin{lemma}[Per-vertex normal transfer to the dissipation graph]
\label{lem:normal-pullback}
Fix a vertex term $m\in\mathcal{T}$, with $\Damp^{V,m}(x,\dot\vartheta^D)=\eta[\Lambda^i]^m(x)\|C^m(x)\dot\vartheta^D\|$ smoothed as $\eta[\Lambda^i]^m(x)\psi_\mu(C^m(x)\dot\vartheta^D)$, $x=\theta_{\gamma(i)}$. Assume the primitive graph point $(v^\sharp,q^\sharp)$ and Clarke normal $n^\sharp=(a^\sharp,b^\sharp)$ of~\prettyref{lem:three-region-transfer}, its weight-normal $\tau^\sharp$ of~\prettyref{lem:tau-bookkeeping}, the preimage $\dot\vartheta^{D,m\sharp}$ and right inverse $\iota^m$ of~\prettyref{lem:vertex-preimage}, and the per-vertex estimates of~\prettyref{lem:pervertex-activation}. Assume $x=\theta_{\gamma(i)}\in\mathcal K_H$, let the full velocity $\dot\vartheta_i^D$ range over a fixed bounded set, and let the tested lifted direction $\xi^{\mathrm f}$ range over a fixed bounded set. The localized estimates of~\prettyref{lem:pervertex-activation} give uniform bounds for $C^m$, active-side estimates for $D_xC^m,D_x\alpha^m$, and uniform bounds for the weighted derivative products used below. Then there is a full Clarke normal $n_i^{m\sharp}\in N^C_{\gph\partial_{\dot\vartheta^D}\Damp^{V,m}}(x,\dot\vartheta_i^{D,m\dagger},q_i^{m\dagger})$ in the full velocity $\dot\vartheta_i^D$, with $\dot\vartheta_i^{D,m\dagger}=\dot\vartheta_i^D+\iota^m(v^\sharp-v)$, $q_i^{m\dagger}=(C^m)^Tq^\sharp$, and force block $\approx-\xi^{\mathrm f}$, satisfying
\begin{ResizedEq}
\label{eq:np-error}
\big\|n_i^{m\sharp}-N_{i,m}^{\mathrm{sm}}(\xi^{\mathrm f})\big\|\leq C\mu^{1/30},
\end{ResizedEq}
where $N_{i,m}^{\mathrm{sm}}(\xi^{\mathrm f})$ is the smooth full-velocity normal of~\prettyref{eq:smooth-full-normal} and $C$ depends only on the fixed bounds above and is independent of $\mu,\EpsFri,\EpsStop$, the timestep, and $m$. No Schur reduction to a reduced graph normal is performed here. At $\eta[\Lambda^i]^m(x)=0$ the transfer error is zero.
\end{lemma}
\begin{proof}[Proof of~\prettyref{lem:normal-pullback}]
Abbreviate $\alpha\triangleq\eta[\Lambda^i]^m$, $C\triangleq C^m$, $v\triangleq v_\parallel^m$, $z\triangleq C^m\xi^{\mathrm f}$ with $\xi^{\mathrm f}=(\xi,\zeta^\star)$ the full lift of~\prettyref{lem:smooth-normal-identification}, and $g_\mu=\nabla\psi_\mu(v)$, $H_\mu=\nabla^2\psi_\mu(v)$, so the smooth blocks are $a_\mu\triangleq\alpha H_\mu z$, $b_\mu\triangleq-z$, $\tau_\mu\triangleq g_\mu^Tz$, $\zeta_\mu\triangleq\alpha g_\mu$. Build the full-velocity data
\begin{ResizedEq}
\dot\vartheta^{D,m\sharp}=&\dot\vartheta^D+\iota^m(v^\sharp-v)\\
s^\sharp=&-\xi^{\mathrm f}+\iota^m(b^\sharp+z)\\
n_q^\sharp=&s^\sharp\\
n_{\dot\vartheta^D}^\sharp=&C^Ta^\sharp,
\end{ResizedEq}
and $n_x^\sharp=D_x\alpha^T\tau^\sharp+D_x[C \dot\vartheta^{D,m\sharp}]^Ta^\sharp-D_x[C^Tq^\sharp]^Ts^\sharp$, so $C s^\sharp=-z+(b^\sharp+z)=b^\sharp$ and $\|s^\sharp+\xi^{\mathrm f}\|=\|b^\sharp+z\|\leq\omega_n(\mu)$. For $\alpha>0$, by~\prettyref{lem:exact-pullback}, $n_i^{m\sharp}\triangleq n^{m\sharp}=(n_x^\sharp,n_{\dot\vartheta^D}^\sharp,n_q^\sharp)\in N^C_{\gph\partial_{\dot\vartheta^D}\Damp^{V,m}}$ is a full Clarke normal exactly (the degenerate $\alpha=0$ case is treated at the end), with force ($q$-)block $\approx-\xi^{\mathrm f}$; we keep it in the full velocity and perform no Schur reduction. Thus~\prettyref{eq:np-error} is exactly $\|n^{m\sharp}-n_\mu\|\leq C\mu^{1/30}$ with $n_\mu=N_{i,m}^{\mathrm{sm}}(\xi^{\mathrm f})$ of~\prettyref{eq:smooth-full-normal} (the full smooth normal generated by $\xi^{\mathrm f}$), which we now prove. The graph-point closeness $\|\dot\vartheta^{D,m\sharp}-\dot\vartheta^D\|\leq\omega_g(\mu)$ is~\prettyref{eq:vertex-preimage-bound}.

Because $\iota^m$ has a constant $-I$ block and $D_x\iota^m=0$, the preimage adjustment changes only the plane-translation block, whose coefficient in $C$ is the constant $-I$, so
\begin{ResizedEq}
\label{eq:np-cancel}
D_x[C \dot\vartheta^{D,m\sharp}]=D_x[C\dot\vartheta^D]=(D_xC)\dot\vartheta^D,
\end{ResizedEq}
and no possibly large $D_xC$ multiplies the displacement $\dot\vartheta^{D,m\sharp}-\dot\vartheta^D$. Subtracting $n_\mu$ (the pullback of the smooth master normal, with $s_\mu=-\xi^{\mathrm f}$) from $n^{m\sharp}$ and using~\prettyref{eq:np-cancel},
\begin{ResizedEq}
n_x^\sharp-n_{x,\mu}=&\underbrace{D_x\alpha^T(\tau^\sharp-\tau_\mu)}_{T_1}+\underbrace{D_x[C\dot\vartheta^D]^T(a^\sharp-a_\mu)}_{T_2}\\
&+\underbrace{D_x[C^T(\zeta^\sharp-\zeta_\mu)]^T\xi^{\mathrm f}}_{T_3}-\underbrace{D_x[C^Tq^\sharp]^T(s^\sharp+\xi^{\mathrm f})}_{T_4},
\end{ResizedEq}
with $\zeta^\sharp\triangleq q^\sharp$; the $\dot\vartheta^D$- and $q$-blocks obey $\|n_{\dot\vartheta^D}^\sharp-n_{\dot\vartheta^D,\mu}\|=\|C^T(a^\sharp-a_\mu)\|\leq\|C\|\,\omega_n(\mu)$ and $\|n_q^\sharp-n_{q,\mu}\|=\|s^\sharp+\xi^{\mathrm f}\|\leq\omega_n(\mu)$, both $\mathcal{O}(\mu^{1/10})$ since $\|C\|$ is bounded. We bound $T_1$-$T_4$, splitting at the weight scale $\alpha_\mu=\mu^{1/5}$ and using~\prettyref{lem:pervertex-activation} and the graph-force bound $\|\zeta^\sharp-\zeta_\mu\|\leq2\alpha$ of~\prettyref{lem:tau-bookkeeping}.

$T_4$ is weight-uniform: $\|q^\sharp\|\leq\alpha$ (subgradient of a norm) and $\|s^\sharp+\xi^{\mathrm f}\|\leq\omega_n(\mu)$, so $\|T_4\|\leq(\alpha\|D_xC\|)\,\omega_n(\mu)\leq C\alpha^{2/3}\mu^{1/10}\leq C\mu^{1/10}$, using $\alpha\|D_xC\|\leq C\alpha^{2/3}\leq C$.

Small weight, $0\leq\alpha<\alpha_\mu$. Here $\|D_x\alpha\|\leq C\alpha^{1/3}\leq C\mu^{1/15}$ and $|\tau^\sharp-\tau_\mu|\leq C$ (\prettyref{lem:tau-bookkeeping}\ref{it:tau-small}), so $\|T_1\|\leq C\mu^{1/15}$. For $T_3$, $\|\zeta^\sharp-\zeta_\mu\|\leq2\alpha$ gives $\|T_3\|\leq2\alpha\|D_xC\|\leq2C\alpha^{2/3}\leq C\mu^{2/15}$. For $T_2$: in the sticking/transition construction $a^\sharp=\alpha H_\mu z=a_\mu$ exactly (\prettyref{lem:three-region-transfer}, Cases~I-II small weight), so $T_2=0$; in the sliding construction $\|a^\sharp-a_\mu\|\leq C\alpha\,\mu^{2}/r^3\leq C\alpha\mu^{4/5}$ (\prettyref{eq:tr-slide-normal}, $r>r_s=\mu^{2/5}$), whence $\|T_2\|\leq\|\dot\vartheta^D\|\,\|D_xC\|\,C\alpha\mu^{4/5}=C\|\dot\vartheta^D\|(\alpha\|D_xC\|)\mu^{4/5}\leq C\mu^{4/5}$.

Large weight, $\alpha\geq\alpha_\mu$. Here $\|D_xC\|\leq C \alpha_\mu^{-1/3}=C\mu^{-1/15}$ (\prettyref{eq:pv-threshold}) and $\|D_x\alpha\|\leq C$. By~\prettyref{lem:three-region-transfer} and~\ref{lem:tau-bookkeeping}\ref{it:tau-large}, $|\tau^\sharp-\tau_\mu|\leq C\mu^{1/10}$, $\|a^\sharp-a_\mu\|\leq\omega_n(\mu)=C\mu^{1/10}$, and $\|\zeta^\sharp-\zeta_\mu\|\leq\omega_g(\mu)=C\mu^{1/5}$. Hence, with $\|\dot\vartheta^D\|$ bounded,
\begin{ResizedEq}
\|T_1\|\leq&C\mu^{1/10}\\
\|T_2\|\leq&C\mu^{-1/15}\mu^{1/10}=C\mu^{1/30}\\
\|T_3\|\leq&C\mu^{-1/15}\mu^{1/5}=C\mu^{2/15}.
\end{ResizedEq}
Combining branches, $\|n_x^\sharp-n_{x,\mu}\|\leq C\mu^{1/30}$ (the large-weight $T_2$ term is dominant), and with the $\mathcal{O}(\mu^{1/10})$ $\dot\vartheta^D$- and $q$-blocks, $\|n^{m\sharp}-n_\mu\|\leq C\mu^{1/30}$, proving~\prettyref{eq:np-error}.

If $\alpha(x)=0$, nonnegativity and $C^1$ regularity of $\alpha$ give $D_x\alpha(x)=0$ (\prettyref{lem:pervertex-activation}), hence $\alpha(x')=o(\|x'-x\|)$ as $x'\to x$. Every point of $\gph\partial_{\dot\vartheta^D}\Damp^{V,m}$ near $(x,\dot\vartheta^D,0)$ has force $q'=(C^m(x'))^T\zeta'$ with $\|\zeta'\|\le\alpha(x')$ and $\|C^m\|$ bounded, so $\|q'\|\le\|C^m\|\,\alpha(x')=o(\|x'-x\|)$; thus along any graph curve through $(x,\dot\vartheta^D,0)$ the force component is first-order stationary, and the Clarke tangent cone is contained in $\{(\delta x,\delta\dot\vartheta^D,0)\}$. Its polar, the Clarke normal cone $N^C_{\gph\partial_{\dot\vartheta^D}\Damp^{V,m}}(x,\dot\vartheta^D,0)$, therefore contains $\{0\}\times\{0\}\times\R^{\dim q}$, so $(0,0,-\xi^{\mathrm f})$ is an exact Clarke normal (here $q_i^{m\dagger}=0$). The smooth target $N_{i,m}^{\mathrm{sm}}(\xi^{\mathrm f})$ has $x$- and velocity-blocks carrying the factor $\alpha$ or $D_x\alpha$, both zero at $x$, and force block $-\xi^{\mathrm f}$; hence $n_i^{m\sharp}=(0,0,-\xi^{\mathrm f})=N_{i,m}^{\mathrm{sm}}(\xi^{\mathrm f})$ and the transfer error is zero.
\end{proof}

The pulled-back smooth full-velocity normal of~\prettyref{lem:smooth-normal-identification} is the per-term target the vertexwise transfer must reproduce; we now state the appendix-local, full physical-plus-plane compatibility notion that the pseudo-Huber smoothing will be shown to satisfy. It is stated on the damping-active velocity $\dot\vartheta_i^D$ and uses full graph normals of $\gph\partial_{\dot\vartheta^D}\Damp^{V,m}$; it does not presuppose that the reduced dissipation decomposes vertexwise, and it keeps the plane damper exact and external.
\begin{definition}[Schur-lifted per-vertex normal compatibility, full physical-plus-plane]
\label{def:S-D-compatible-ABD}
Work in the complete damping-active velocity $\dot\vartheta_i^D$ associated with~\prettyref{pb:CRDS-ABD-extended}, with the per-vertex full Hessian $\Psi_i^m$ of~\prettyref{eq:abd-Psi-split}, the common Schur lift $\widehat\xi=(\xi,\zeta^\star(\xi))$ of~\prettyref{eq:abd-common-lift}, and the smooth full per-vertex normal $N_{i,m}^{\mathrm{sm}}(\widehat\xi)$ of~\prettyref{eq:smooth-full-normal}; the plane-velocity damper $\EpsDamp\|\varpi\|^2$ is retained exactly and externally and carries no replaceable graph normal. Fix the $\EpsFri$-independent energy constant $\bar C_E\triangleq(|\mathcal{T}|+2\EpsDamp)/\delta_{\min}$. A~\prettyref{def:S-D} smoothing is Schur-lifted normal-compatible if there is a modulus $\omega_N(\EpsFri)\downarrow0$ as $\EpsFri\downarrow0$, uniform over $i\in\mathcal{I}$, over $m\in\mathcal{T}$, over every previous configuration $x=\theta_{\gamma(i)}\in\mathcal K_H$, and over damping-active velocities in a fixed bounded set, such that whenever the tested reduced direction obeys the normalized bounds $\|\xi\|\leq1$ and $\xi^T K_i\xi\leq\bar C_E$, each friction vertex $m\in\mathcal{T}$ admits a nearby damping-active velocity $\dot\vartheta_i^{D,m\dagger}$, an exact force $q_i^{m\dagger}\in\partial_{\dot\vartheta^D}\Damp^{V,m}(\theta_{\gamma(i)},\dot\vartheta_i^{D,m\dagger})$, and a full Clarke normal
\begin{ResizedEq}
&n_i^{m\sharp}\in N^C_{\gph\partial_{\dot\vartheta^D}\Damp^{V,m}}(\theta_{\gamma(i)},\dot\vartheta_i^{D,m\dagger},q_i^{m\dagger})\quad(\text{force block }\approx-\widehat\xi),\\
&\|\dot\vartheta_i^{D,m\dagger}-\dot\vartheta_i^D\|+\|q_i^{m\dagger}-q_i^{m,\EpsFri}\|\leq C\mu(\EpsFri)^{1/5},\\
&\|n_i^{m\sharp}-N_{i,m}^{\mathrm{sm}}(\widehat\xi)\|\leq\omega_N(\EpsFri),
\end{ResizedEq}
where $q_i^{m,\EpsFri}=\nabla_{\dot\vartheta^D}\Damp_{\EpsFri,\LHamilton}^{V,m}$ is the smooth damping-active friction force. No vertexwise decomposition of the reduced dissipation is assumed; only the full graph normals and the aggregate first-order lift are used.
\end{definition}

\begin{lemma}[Per-vertex normal compatibility of the ABD smoothing]
\label{lem:abd-normal-compatible}
Distinguish the aggregate weighted-load ceiling $\bar\Lambda_{\mathrm{agg}}$ of~\prettyref{appen:penalty-proofs}, now taken over the compact buffer $\widehat{\mathcal K}_H$ (with $\eta\sum_{\mathcal{A}}\lambda_r\leq\bar\Lambda_{\mathrm{agg}}$ there), from the individual friction-weight ceiling $\bar\alpha\triangleq\max\{\alpha_i^m(x):x\in\widehat{\mathcal K}_H,\ i\in\mathcal I,\ m\in\mathcal{T}\}$, $\alpha_i^m(x)\triangleq\eta[\Lambda^i]^m(x)$; since each $\alpha_i^m(x)$ is one nonnegative summand of $\eta\sum_{\mathcal{A}}\lambda_r$, one has $\bar\alpha\leq\bar\Lambda_{\mathrm{agg}}$, so the aggregate ceiling bounds the individual one. Choose the smoothing scale
\begin{ResizedEq}
\label{eq:abd-mu-schedule}
\mu(\EpsFri)\leq\min\Big\{\tfrac{\delta_{\min}\EpsFri^2}{\bar\Lambda_{\mathrm{agg}}},\ c_{\mathrm{nc}}\EpsFri^{30}\Big\}
\end{ResizedEq}
with $c_{\mathrm{nc}}>0$ fixed; the first factor is the aggregate reduced closeness choice of~\prettyref{lem:smooth-DRp-Regularity} and is dropped in the zero-load case $\bar\Lambda_{\mathrm{agg}}=0$ (treated separately below). Then the pseudo-Huber smoothing $\eta[\Lambda^i]^m\psi_\mu(v_\parallel^m)$ satisfies the full physical-plus-plane compatibility condition~\prettyref{def:S-D-compatible-ABD}, with a single modulus $\omega_N(\EpsFri)=C\mu(\EpsFri)^{1/30}+C\EpsFri\downarrow0$ uniform over $i,m$, $x\in\mathcal K_H$, and damping-active velocities in a fixed bounded set.
\end{lemma}
\begin{proof}[Proof of~\prettyref{lem:abd-normal-compatible}]
Zero-load case $\bar\Lambda_{\mathrm{agg}}=0$. Then every nonnegative friction weight $\alpha^m\leq\bar\Lambda_{\mathrm{agg}}=0$ vanishes, so each $\Damp^{V,m}\equiv0$; its graph is $\R^n\times\{0\}$ and the transfer is exact with zero error (\prettyref{lem:weighted-norm-clarke} at $\alpha=0$, and~\prettyref{lem:normal-pullback} at $\eta[\Lambda^i]^m=0$). The unreduced plane damper $\EpsDamp\|\varpi\|^2$ is unaffected---only its reduced minimizer, value, and force vanish---and remains exact and external. Only the $c_{\mathrm{nc}}\EpsFri^{30}$ branch of~\prettyref{eq:abd-mu-schedule} is used, and~\prettyref{def:S-D-compatible-ABD} holds with $\omega_N=0$. Assume henceforth $\bar\Lambda_{\mathrm{agg}}>0$.

Fix $i,m$ and a reduced direction with $\|\xi\|\leq1$ and $\xi^T K_i\xi\leq\bar C_E$; form the common lift $\xi^{\mathrm f}=\widehat\xi=(\xi,\zeta^\star)$ of~\prettyref{eq:abd-common-lift} and push $\xi$ to the primitive direction $z\triangleq C^m\widehat\xi=A^m\xi-B^m\zeta^\star=w_m$, with weight $\alpha\triangleq\alpha_i^m(x)=\eta[\Lambda^i]^m(x)\leq\bar\alpha\leq\bar\Lambda_{\mathrm{agg}}$ and primitive velocity $v\triangleq v_\parallel^m$.

(R1) Primitive bounds from~\prettyref{eq:abd-schur-var}. We do not use any operator estimate $\|\zeta^\star\|\leq C\|\xi\|$. By~\prettyref{eq:abd-perterm-energy} every summand of the variational identity~\prettyref{eq:abd-energy-split} is nonnegative, so
\begin{ResizedEq}
&2\EpsDamp\|\zeta^\star\|^2\leq\delta_i\,\xi^T K_i\xi\leq\delta_{\max}\bar C_E\\
&\|\zeta^\star\|\leq\sqrt{\tfrac{\delta_{\max}\bar C_E}{2\EpsDamp}},
\end{ResizedEq}
whence, with $\|A^m\|,\|B^m\|$ uniformly bounded and $\|\xi\|\leq1$, $\|z\|=\|A^m\xi-B^m\zeta^\star\|\leq C_z$; and the per-vertex energy is one nonnegative summand, $\alpha\,z^T\nabla^2\psi_\mu(v)z=\widehat\xi^T\Psi_i^m\widehat\xi\leq\delta_i\,\xi^T K_i\xi\leq\delta_{\max}\bar C_E=:C_E$. These are exactly the hypotheses~\prettyref{eq:tr-hyp} of~\prettyref{lem:three-region-transfer} with $\EpsFri$-independent constants $C_z,C_E$, stated only as boundedness on the normalized energy set.

(R2) Nearby point, exact force, and full-force closeness.~\prettyref{lem:three-region-transfer} supplies a primitive graph point $(v^\sharp,q^\sharp)\in\gph(\alpha\partial\|\cdot\|)$ (so $q^\sharp\in\alpha\partial\|v^\sharp\|$ by construction) and Clarke normal $n^\sharp$ with $\|v^\sharp-v\|+\|q^\sharp-\alpha\nabla\psi_\mu(v)\|\leq\omega_g(\mu)=C\mu^{1/5}$;~\prettyref{lem:tau-bookkeeping} supplies its weight-normal $\tau^\sharp$. If $\alpha>0$, then the corresponding pair is active and~\prettyref{lem:vertex-preimage} realizes $v^\sharp$ by the damping-active velocity $\dot\vartheta_i^{D,m\dagger}=\dot\vartheta_i^D+\iota^m(v^\sharp-v)$ with $\|\dot\vartheta_i^{D,m\dagger}-\dot\vartheta^D\|\leq\|v^\sharp-v\|\leq C\mu^{1/5}$, and recovers the exact damping-active force $q_i^{m\dagger}=(C^m)^Tq^\sharp\in\partial_{\dot\vartheta^D}\Damp^{V,m}(\theta_{\gamma(i)},\dot\vartheta_i^{D,m\dagger})$. If $\alpha=0$, take $\dot\vartheta_i^{D,m\dagger}=\dot\vartheta_i^D$ and $q_i^{m\dagger}=q_i^{m,\EpsFri}=0$; the zero-weight branch of~\prettyref{lem:normal-pullback} supplies the exact normal with zero error. Lifting the primitive force closeness through the bounded $(C^m)^T$ gives the damping-active physical-plus-plane friction-force closeness
\begin{ResizedEq}
&\|q_i^{m\dagger}-q_i^{m,\EpsFri}\|\\
=&\|(C^m)^T(q^\sharp-\alpha\nabla\psi_\mu(v))\|\\
\leq&\|C^m\|\,\omega_g(\mu)\leq C\mu^{1/5},
\end{ResizedEq}
with $q_i^{m,\EpsFri}=(C^m)^T\alpha\nabla\psi_\mu(v)=\nabla_{\dot\vartheta^D}\Damp_{\EpsFri,\LHamilton}^{V,m}$. Choose $c_{\mathrm{nc}}>0$ once, uniformly, so that $Cc_{\mathrm{nc}}^{1/5}\leq\min\{1,\delta_{\min}\}$. Under~\prettyref{eq:abd-mu-schedule}, $\mu\leq c_{\mathrm{nc}}\EpsFri^{30}$ then gives, for $0<\EpsFri\leq1$ and $\delta_i\geq\delta_{\min}$,
\begin{ResizedEq}
C\mu^{1/5}
\leq Cc_{\mathrm{nc}}^{1/5}\EpsFri^{6}
\leq\min\{\EpsFri,\delta_{\min}\EpsFri\}
\leq\min\{\EpsFri,\delta_i\EpsFri\}.
\end{ResizedEq}
Thus the~\prettyref{def:S-D-compatible-ABD} graph-point and force-closeness bounds hold. The aggregate reduced~\prettyref{def:S-D} smoothing separately satisfies its own $\theta_i$-gradient-scale closeness (an aggregate reduced nearby point and a $\delta_i\EpsFri$ gradient discrepancy, via $\nabla_{\theta_i}(\delta_i\Damp_R^V)=\nabla_{\dot{\theta}_i}\Damp_R^V$) by~\prettyref{lem:smooth-DRp-Regularity} under $\mu\leq\delta_{\min}\EpsFri^2/\bar\Lambda_{\mathrm{agg}}$; that aggregate bound is what the reduced algorithm consumes and is not conflated with the per-vertex full-space force closeness established above.

Normal transfer.~\prettyref{lem:normal-pullback}---using the per-vertex activation estimates of~\prettyref{lem:pervertex-activation} and the smooth-normal identification of~\prettyref{lem:smooth-normal-identification}---transfers $n^\sharp$ to a full Clarke normal $n_i^{m\sharp}\in N^C_{\gph\partial_{\dot\vartheta^D}\Damp^{V,m}}(\theta_{\gamma(i)},\dot\vartheta_i^{D,m\dagger},q_i^{m\dagger})$, with force block $\approx-\xi^{\mathrm f}$, matching the smooth full normal $N_{i,m}^{\mathrm{sm}}(\xi^{\mathrm f})$ up to $\|n_i^{m\sharp}-N_{i,m}^{\mathrm{sm}}(\xi^{\mathrm f})\|\leq C\mu^{1/30}$. Under $\mu\leq c_{\mathrm{nc}}\EpsFri^{30}$, $\mu^{1/30}\leq c_{\mathrm{nc}}^{1/30}\EpsFri$, so $\|n_i^{m\sharp}-N_{i,m}^{\mathrm{sm}}(\xi^{\mathrm f})\|\leq C\mu^{1/30}+C\EpsFri=:\omega_N(\EpsFri)\downarrow0$. Uniformity over the finite $\mathcal{T}$ and the timesteps holds because $|\mathcal{T}|,\delta_{\min},\EpsDamp$, the mesh, and the~\prettyref{lem:pervertex-activation} constants are $\EpsFri$-independent, so $C_z,C_E,\bar C_E$ enter only through $\EpsFri$-independent quantities.
\end{proof}
\begin{replemma}{lem:abd-compatible}
\label{lem:abd-compatible-formal}
The pseudo-Huber smoothing of~\prettyref{lem:smooth-DRp}, with $\mu(\EpsFri)$ chosen as in~\prettyref{eq:abd-mu-schedule}, replaces each Euclidean norm in the finite friction sum ($|\mathcal{T}|$ fixed, $\EpsFri,\EpsStop$-independent) of~\prettyref{eq:damping-ABD} by its pseudo-Huber surrogate; the resulting aggregate reduced dissipation is a~\prettyref{def:S-D} smoothing (\prettyref{lem:smooth-DRp-Regularity}), the plane-velocity damper $\EpsDamp\|\varpi\|^2$ is retained exactly and externally, and each full per-vertex graph satisfies the full physical-plus-plane Schur-lifted normal-compatibility condition~\prettyref{def:S-D-compatible-ABD} on the compact configuration/full-velocity region used by the returned-trajectory transfer (\prettyref{lem:abd-normal-compatible}). It does not claim that the reduced dissipation decomposes vertexwise or that the per-term blocks are positive-semidefinite Hessians of reduced per-vertex scalars: the $K_i^m,\Gamma_i^m$ are the common-lift components~\prettyref{eq:abd-perterm-blocks}, so the $\dot{\theta}_i$-block action is additive ($\sum_m K_i^m=K_i$ as operators), but the reduced energy carries the extra external-damper term~\prettyref{eq:abd-energy-split} and does not split vertexwise, and the $K_i^m$ need not be positive semidefinite. This is the sense in which~\prettyref{lem:abd-compatible} holds.
\end{replemma}
\begin{proof}[Proof of~\prettyref{lem:abd-compatible}]
The friction part of~\prettyref{eq:damping-ABD} is the finite sum $\sum_{m\in\mathcal{T}}\eta[\Lambda^i]^m\|v_\parallel^m\|$ over the fixed vertex-plane incidences, smoothed termwise by~\prettyref{lem:smooth-DRp}. The aggregate reduced dissipation (obtained after the common minimization over the plane velocities) is a~\prettyref{def:S-D} smoothing by~\prettyref{lem:smooth-DRp-Regularity}; this is the closeness the reduced algorithm consumes and is not a separate per-term~\prettyref{def:S-D} claim. The quadratic plane damper is already smooth and is kept exact and external. The full physical-plus-plane per-vertex normal-compatibility claim is~\prettyref{lem:abd-normal-compatible}.
\end{proof}

\begin{lemma}[Exact unnormalized lift to the reduced-configuration/full-damping-velocity subsystem]
\label{lem:abd-kkt-lift}
Suppose the reduced smooth modified optimization problem obtained by applying~\prettyref{pb:CITO-M} to the reduced CRDS~\prettyref{pb:PSR-CRDS} admits a finite $\EpsStop$-KKT certificate in the sense of~\prettyref{def:eps-KKT}, with reduced physics multipliers $\mu_{\theta_i}$, control multiplier $\mu_\StackCtrl$, and auxiliary controls $u_i'$. Reconstruct the unique minimizing plane velocities $\varpi^\star$ by~\prettyref{lem:varpi-reconstruction} and lift $\mu_{\theta_i}\mapsto\widehat\mu_i=(\mu_{\theta_i},\nu_i)$ by~\prettyref{lem:schur-lift}. Then, with no normalization applied, the reduced certificate lifts exactly to the reduced-configuration/full-damping-velocity subsystem of the smooth full first-order relations: the physical force-balance residual is unchanged, the plane-velocity force-balance residual is exactly $0$, the physical adjoint equals the reduced adjoint, the plane-velocity costate residual is exactly $0$, and the control-stationarity and two-point control lines of~\prettyref{eq:eps-kkt-system} are unchanged. No claim about reconstruction, feasibility, or stationarity of the auxiliary plane-configuration variables is made.
\end{lemma}
\begin{proof}[Proof of~\prettyref{lem:abd-kkt-lift}]
The reconstructed $\varpi^\star$ satisfies $\nabla_\varpi\Damp_{\EpsFri,\LHamilton}^V=0$ (\prettyref{lem:varpi-reconstruction}), so the plane-velocity force-balance residual vanishes identically, while the envelope identity makes the physical force-balance residual coincide with the reduced one.

Differentiating the inner first-order condition $\nabla_\varpi\Damp_{\EpsFri,\LHamilton}^V(\theta_{\gamma(i)},\dot\theta_i,\varpi^\star)=0$ by the implicit-function theorem identifies the derivative of the reduced velocity gradient with the Schur complement $S_i=\delta_iK_i$ of~\prettyref{eq:abd-schur-blocks}. Since
\begin{ResizedEq}
D_{\theta_i}\dot\theta_i=I/\delta_i,
\qquad
D_{\theta_{\gamma(i)}}\dot\theta_i=-I/\delta_i,
\end{ResizedEq}
the coordinate pullback gives the current contribution $S_i\mu_{\theta_i}/\delta_i=K_i\mu_{\theta_i}$ and the velocity-mediated predecessor contribution $-S_i\mu_{\theta_i}/\delta_i=-K_i\mu_{\theta_i}$.

For the full lift,~\prettyref{eq:schur-identities} of~\prettyref{lem:schur-lift} gives $\Psi_i\widehat\mu_i=(S_i\mu_{\theta_i},0)$. Hence its pullback through $D_{\theta_i}\dot\theta_i=I/\delta_i$ gives the same current physical block $K_i\mu_{\theta_i}$, while its plane-velocity costate block vanishes exactly. The explicit predecessor-configuration partial at fixed $(\dot\theta_i,\varpi_i)$ is
\begin{ResizedEq}
\Psi_i^{\gamma\theta}\mu_{\theta_i}+\Psi_i^{\gamma\varpi}\nu_i
=&\left(\Psi_i^{\gamma\theta}-\Psi_i^{\gamma\varpi}(\Psi_i^{\varpi\varpi})^{-1}\Psi_i^{\varpi\theta}\right)\mu_{\theta_i}\\
=&\Gamma_i^T\mu_{\theta_i},
\end{ResizedEq}
where the last equality is the transpose of the coupling identity in~\prettyref{eq:abd-schur-blocks}, using symmetry of $\Psi_i^{\varpi\varpi}$ and Schwarz symmetry of the mixed derivatives. Adding the velocity-mediated predecessor term gives
\begin{ResizedEq}
\Gamma_i^T\mu_{\theta_i}-K_i\mu_{\theta_i}=(\Gamma_i-K_i)^T\mu_{\theta_i}.
\end{ResizedEq}
Thus, in the adjoint block for a fixed $\theta_i$, residual row $i$ contributes $K_i\mu_{\theta_i}$ and every immediate successor $j$ with $\gamma(j)=i$ contributes $(\Gamma_j-K_j)^T\mu_{\theta_j}$, exactly as in the reduced Jacobian transpose; second-neighbor successors carry only kinetic terms. Hence the physical adjoint of the damping-active subsystem equals the reduced adjoint, and its plane-velocity costate block vanishes.

The objective and controls do not depend on the reconstructed plane velocities, so the control-stationarity and two-point control lines of~\prettyref{eq:eps-kkt-system} are unchanged. No approximation is introduced. The auxiliary plane-configuration variables and their feasibility and stationarity relations are not part of this Schur lift and must be supplied separately.
\end{proof}
To state the theorem we specialize~\prettyref{def:eps-stationary-NS} to the full physical-plus-plane formulation~\prettyref{pb:CRDS-ABD-extended}, keeping the plane damper exact and external. Lift each friction force to a free variable over the damping-active velocity $\dot\vartheta_i^D$: with $q_i^m$ a stand-in for $\nabla_{\dot\vartheta^D}\Damp^{V,m}$, the full residual and full graph map, for $i\in\mathcal I$, are
\begin{ResizedEq}
\label{eq:abd-full-residual}
&\widehat F_i^{\mathrm{full}}\triangleq\TWOC{M\ddot{\theta}_i-f-g-c+\sum_m[q_i^m]_{\dot{\theta}}-\FPP{h^e}{\theta_i}^T\lambda_i^e-\FPP{h^i}{\theta_i}(\vartheta_i^\star)^T\lambda_i^i}{\sum_m[q_i^m]_{\varpi}+2\EpsDamp\varpi_i}\\
&Q_i^{\mathrm{full}}(\StackConf,\varpi,q_i)\triangleq(\theta_{\gamma(i)},\tfrac{\theta_i-\theta_{\gamma(i)}}{\delta_i},\varpi,q_i),\\
&[\nabla Q_i^{\mathrm{full}}]^T(a,b,d,c)=(\tfrac{b}{\delta_i},\,a-\tfrac{b}{\delta_i},\,d,\,c),
\end{ResizedEq}
the plane block of $\widehat F_i^{\mathrm{full}}$ carrying the exact external damper $2\EpsDamp\varpi_i$, and $[\nabla Q_i^{\mathrm{full}}]^T$ passing the $\varpi$-covector $d$ through unchanged.
\begin{definition}[full physical-plus-plane C-$\TWO{\EpsStop}{\EpsFJAll}$-FJ-stationarity]
\label{def:eps-stationary-NS-ABD}
A returned physical trajectory $(\StackConf,\StackCtrl)$, after fixing the finite extended base state $\vartheta_0^\star$ and reconstructing in increasing order over $i\in\mathcal I$ the extended configurations $\vartheta_i^\star=\langle\theta_i,\{p_i^\star,\varsigma_i^\star,\omega_i^\star\}\rangle$ and true full finite-difference velocities $\dot\vartheta_i^\star=(\vartheta_i^\star-\vartheta_{\gamma(i)}^\star)/\delta_i$, whose damping-active velocities $\dot\vartheta_i^{D\star}\triangleq\TWO{\dot\theta_i}{\varpi_i^\star}$ satisfy $\varsigma_i^\star=\varsigma_{\gamma(i)}^\star+\delta_i\dot \varsigma_i^\star$ and $\omega_i^\star=\omega_{\gamma(i)}^\star+\delta_i\dot\omega_i^\star$, satisfies the full force-parameterized homogeneous Clarke $\TWO{\EpsStop}{\EpsFJAll}$-FJ-stationarity condition of~\prettyref{pb:CRDS-ABD-extended} if there exist an objective multiplier $\mu_0\geq0$, lifted physics/plane multipliers $\widehat\mu_i=(\mu_{\theta_i},\nu_i)$, control multipliers $\mu_\StackCtrl$ (sub-blocks $\mu_{u_i}$), free-sign/nonnegative band multipliers, auxiliary controls $u_i'$, and, for every $i\in\mathcal I$ and $m\in\mathcal{T}$, a nearby full velocity $\dot\vartheta_i^{D,m\sharp}$, a friction force $q_i^{m\sharp}\in\partial_{\dot\vartheta^D}\Damp^{V,m}(\theta_{\gamma(i)},\dot\vartheta_i^{D,m\sharp})$, and a full Clarke graph normal $n_i^{m\sharp}\in N^C_{\gph\partial_{\dot\vartheta^D}\Damp^{V,m}}(\theta_{\gamma(i)},\dot\vartheta_i^{D,m\sharp},q_i^{m\sharp})$, such that
\begin{ResizedEq}
&\|\widehat F^{\mathrm{full}}(\StackConf,\varpi^\star,\StackCtrl,(q^{m\sharp}))\|_1\leq\EpsStop,\quad\|\dot\vartheta_i^{D,m\sharp}-\dot\vartheta_i^{D\star}\|\leq\EpsStop\ \ \forall m\\
&\Big\|\nabla_{(\StackConf,\varpi,q)}\big[\mu_0 O+\textstyle\sum_i\langle\widehat\mu_i,\widehat F_i^{\mathrm{full}}\rangle+\text{(band terms)}\big]+\\
&\textstyle\sum_{i\in\mathcal I}\sum_{m\in\mathcal T}[\nabla Q_i^{\mathrm{full}}]^T n_i^{m\sharp}\Big\|\leq\EpsStop\\
&0=\FPP{h^i}{p_i}(\vartheta_i^\star)^T\lambda_i^i\quad\forall i\in\mathcal I\\
&\left\|\mu_0\FPP{O}{\StackCtrl}+[\StatJacU]^T\mu_\StackConf+\mu_\StackCtrl\right\|\leq\EpsStop\\
&u_i\in\mathcal{D},\ u_i'\in\mathcal{D},\ \|u_i-u_i'\|\leq\EpsStop,\ \mu_{u_i}\in N_{\mathcal{D}}(u_i')\\
&|h^{e,j}(\theta_i)|\leq\EpsFJAll\\
&\min\!\big(h^{i,j}(\vartheta_i^\star),\lambda_i^{i,j}(\vartheta_i^\star),\lambda_i^{i,j}(\vartheta_i^\star)(\EpsFJAll-h^{i,j}(\vartheta_i^\star))\big)\geq0\ \ \forall j,
\end{ResizedEq}
together with the full-tuple nontriviality normalization
\begin{ResizedEq}
\label{eq:abd-ns-normalization}
1=&\mu_0^2+\sum_{i\in\mathcal I}\|\widehat\mu_i\|^2+\|\mu_\StackCtrl\|^2+\\
&\sum_{i\in\mathcal I}\big(\|\mu_i^e\|^2+\|\mu_i^h\|^2+\|\mu_i^\lambda\|^2+\|\mu_i^c\|^2\big)+\\
&\sum_{i\in\mathcal I}\sum_{m\in\mathcal{T}}\|n_i^{m\sharp}\|^2.
\end{ResizedEq}
Here the ``(band terms)'' are the equality/inequality Lagrangian contributions of~\prettyref{def:eps-stationary-NS}, selected zero in the transfer. The plane balance is part of the primal residual~\prettyref{eq:abd-full-residual}. The external damper enters the adjoint solely through the plane-balance term $2\EpsDamp\varpi_i$ of $\widehat F_i^{\mathrm{full}}$: its $\varpi$-gradient in $\nabla_{(\StackConf,\varpi,q)}\langle\widehat\mu_i,\widehat F_i^{\mathrm{full}}\rangle$ equals $2\EpsDamp\nu_i$, so the damper contribution appears exactly once, is retained exactly, and is not replaced by any graph normal (only the friction forces $q_i^m$ and their graph normals are replaced).
\end{definition}

Start from the exact lifted reduced-configuration/full-damping-velocity subsystem supplied conditionally by~\prettyref{lem:abd-kkt-lift}; the auxiliary-configuration relations required for the full certificate are reconstructed separately in Step~2 of the proof of~\prettyref{thm:abd-C-stationarity}. We now normalize the complete full tuple---not the reduced one---because the Schur lift, though linear, is not an isometry, so the lift of a normalized reduced certificate is not a normalized full certificate. For the lifted physics/plane multipliers $\widehat\mu_i=(\mu_{\theta_i},\nu_i)$ (\prettyref{lem:schur-lift}), the control multiplier $\mu_\StackCtrl$, and the smooth full per-vertex normals $N_{i,m}^{\mathrm{sm}}(\widehat\mu_i)$ of~\prettyref{eq:smooth-full-normal}, set the full-space certificate scale
\begin{ResizedEq}
\label{eq:S-full-scale}
(S_{\EpsFri}^{\mathrm{full}})^2\triangleq1+\sum_{i\in\mathcal I}\|\widehat\mu_i\|^2+\|\mu_\StackCtrl\|^2+\sum_{i\in\mathcal{I}}\sum_{m\in\mathcal{T}}\big\|N_{i,m}^{\mathrm{sm}}(\widehat\mu_i)\big\|^2,
\end{ResizedEq}
where $\widehat\mu_i=(\mu_{\theta_i},\nu_i)$ is included directly, so every lifted multiplier participates even when $\mathcal T=\varnothing$; its occurrence in the force block $-\widehat\mu_i$ of each $N_{i,m}^{\mathrm{sm}}$ is additional rather than a substitute for normalization. Divide every component of the full certificate by $S_{\EpsFri}^{\mathrm{full}}$: $\bar\mu_0=1/S_{\EpsFri}^{\mathrm{full}}\in(0,1]$, $\bar\mu_{\theta_i}=\mu_{\theta_i}/S_{\EpsFri}^{\mathrm{full}}$, and, by linearity of the lift, $\widehat{\bar\mu}_i=\widehat\mu_i/S_{\EpsFri}^{\mathrm{full}}=(\bar\mu_{\theta_i},\bar\nu_i)$, $\bar\mu_\StackCtrl=\mu_\StackCtrl/S_{\EpsFri}^{\mathrm{full}}$, and $\bar N_{i,m}^{\mathrm{sm}}=N_{i,m}^{\mathrm{sm}}(\widehat\mu_i)/S_{\EpsFri}^{\mathrm{full}}=N_{i,m}^{\mathrm{sm}}(\widehat{\bar\mu}_i)$.
\begin{lemma}[Full-space normalized bounds]
\label{lem:abd-full-normalized}
The full normalized certificate lies on the unit sphere,
\begin{ResizedEq}
\label{eq:abd-full-unit}
\bar\mu_0^2+\sum_{i\in\mathcal I}\|\widehat{\bar\mu}_i\|^2+\|\bar\mu_\StackCtrl\|^2+\sum_{i\in\mathcal I}\sum_{m\in\mathcal T}\|\bar N_{i,m}^{\mathrm{sm}}\|^2=1,
\end{ResizedEq}
and, for every $i\in\mathcal{I}$ and $m\in\mathcal{T}$, with $\bar C_E=(|\mathcal{T}|+2\EpsDamp)/\delta_{\min}$,
\begin{ResizedEq}
\label{eq:abd-full-bounds}
&\|\widehat{\bar\mu}_i\|\leq1,\quad \|\Psi_i^m\widehat{\bar\mu}_i\|\leq1,\quad \widehat{\bar\mu}_i^T\Psi_i^m\widehat{\bar\mu}_i\leq1,\\
&\|\bar\mu_{\theta_i}\|\leq1,\quad \bar\mu_{\theta_i}^T K_i\bar\mu_{\theta_i}\leq\bar C_E.
\end{ResizedEq}
In addition, at every positive-weight incidence $\alpha^m(\theta_{\gamma(i)})>0$,
\begin{ResizedEq}
&\|C^m(\theta_{\gamma(i)})\widehat{\bar\mu}_i\|\leq C_z\\
&C_z\triangleq\max\left\{1,\sup_{i,m,\,x\in\widehat{\mathcal K}_H:\,\alpha^m(x)>0}\|C^m(x)\|\right\}<\infty\\
&\sup\varnothing\triangleq0.
\end{ResizedEq}
All these bounds are $\EpsFri,\EpsStop$-independent. None of them invokes the reduced-energy route of the general nonsmooth transfer; they are read off the full unit-sphere certificate directly.
\end{lemma}
\begin{proof}[Proof of~\prettyref{lem:abd-full-normalized}]
\prettyref{eq:abd-full-unit} is immediate from~\prettyref{eq:S-full-scale}; it gives $\|\widehat{\bar\mu}_i\|\leq1$ for every $i$ and $\|\bar\mu_\StackCtrl\|\leq1$. Suppose first that $\mathcal{T}=\varnothing$. The split~\prettyref{eq:abd-Psi-split} then gives $\Psi_i=\diag(0,2\EpsDamp I)$ and $\Psi_i^{\varpi\theta}=0$. Hence the Schur identities~\prettyref{eq:schur-identities} of~\prettyref{lem:schur-lift} give $\bar\nu_i=0$ and $K_i=0$. Consequently $\|\widehat{\bar\mu}_i\|=\|\bar\mu_{\theta_i}\|\leq1$ and $\bar\mu_{\theta_i}^T K_i\bar\mu_{\theta_i}=0\leq\bar C_E$; every assertion quantified over $m\in\mathcal{T}$ is vacuous.

Suppose now that $\mathcal{T}\neq\varnothing$. Each $\bar N_{i,m}^{\mathrm{sm}}=N_{i,m}^{\mathrm{sm}}(\widehat{\bar\mu}_i)$ has, by~\prettyref{lem:smooth-normal-identification}, specifically~\prettyref{eq:smooth-full-normal}, velocity block $\Psi_i^m\widehat{\bar\mu}_i$ and force block $-\widehat{\bar\mu}_i$; since $\|\bar N_{i,m}^{\mathrm{sm}}\|\leq1$ by~\prettyref{eq:abd-full-unit}, both $\|\Psi_i^m\widehat{\bar\mu}_i\|\leq1$ and $\|\widehat{\bar\mu}_i\|\leq1$, whence $\widehat{\bar\mu}_i^T\Psi_i^m\widehat{\bar\mu}_i\leq\|\widehat{\bar\mu}_i\|\,\|\Psi_i^m\widehat{\bar\mu}_i\|\leq1$. At an incidence satisfying $\alpha^m(\theta_{\gamma(i)})>0$, the definition of $C_z$ and $\|\widehat{\bar\mu}_i\|\leq1$ additionally give $\|C^m\widehat{\bar\mu}_i\|\leq\|C^m\|\,\|\widehat{\bar\mu}_i\|\leq C_z$; no such coefficient bound is asserted or needed at a zero-weight incidence. The block bound $\|\bar\mu_{\theta_i}\|\leq\|\widehat{\bar\mu}_i\|\leq1$ follows from~\prettyref{eq:abd-full-unit}. Finally, by~\prettyref{eq:schur-identities} and~\ref{eq:abd-Psi-split}, $\delta_i\bar\mu_{\theta_i}^T K_i\bar\mu_{\theta_i}=\widehat{\bar\mu}_i^T\Psi_i\widehat{\bar\mu}_i=\sum_m\widehat{\bar\mu}_i^T\Psi_i^m\widehat{\bar\mu}_i+2\EpsDamp\|\bar\nu_i\|^2\leq|\mathcal{T}|+2\EpsDamp$ (each per-vertex energy $\leq1$ and $\|\bar\nu_i\|\leq\|\widehat{\bar\mu}_i\|\leq1$), so $\bar\mu_{\theta_i}^T K_i\bar\mu_{\theta_i}\leq(|\mathcal{T}|+2\EpsDamp)/\delta_{\min}=\bar C_E$ by~\prettyref{eq:norm-mesh}. The finiteness of $C_z$ follows directly from the uniform active-incidence coefficient bound in~\prettyref{lem:pervertex-activation} and finite incidence. All constants inherit $\EpsFri,\EpsStop$-independence from the fixed $|\mathcal{T}|,\delta_{\min},\EpsDamp$ and compact-region data.
\end{proof}
\begin{lemma}[Full-space dual nontriviality]
\label{lem:abd-nonzero-full}
Let $\widehat{\bar\mu}_\StackConf\triangleq(\widehat{\bar\mu}_i)_{i\in\mathcal I}$ and $\bar Z_{\mathrm{full}}\triangleq(\bar\mu_0,\widehat{\bar\mu}_\StackConf,\bar\mu_\StackCtrl,(\bar N_{i,m}^{\mathrm{sm}}))$, so $\|\bar Z_{\mathrm{full}}\|=1$ by~\prettyref{eq:abd-full-unit}. Suppose that for each $i,m$ a full Clarke normal $n_i^{m\sharp}$ satisfies $\|n_i^{m\sharp}-\bar N_{i,m}^{\mathrm{sm}}\|\leq\omega_N(\EpsFri)$, and set $Z^\sharp_{\mathrm{full}}\triangleq(\bar\mu_0,\widehat{\bar\mu}_\StackConf,\bar\mu_\StackCtrl,(n_i^{m\sharp}))$, replacing only the friction normals (the external damper contributes no replaceable normal). Then $\|Z^\sharp_{\mathrm{full}}-\bar Z_{\mathrm{full}}\|\leq\sqrt{N^\star}\,\omega_N(\EpsFri)$, and if $\sqrt{N^\star}\,\omega_N(\EpsFri)\leq1/2$ then $\|Z^\sharp_{\mathrm{full}}\|\geq1/2$.
\end{lemma}
\begin{proof}[Proof of~\prettyref{lem:abd-nonzero-full}]
At most $N^\star=|\mathcal{T}|(|G_0|+N_{\mathrm{sub}})$ friction-normal blocks differ, each by at most $\omega_N(\EpsFri)$, so $\|Z^\sharp_{\mathrm{full}}-\bar Z_{\mathrm{full}}\|\leq(\sum_{i,m}\omega_N(\EpsFri)^2)^{1/2}\leq\sqrt{N^\star}\,\omega_N(\EpsFri)$; the reverse triangle inequality with $\|\bar Z_{\mathrm{full}}\|=1$ gives $\|Z^\sharp_{\mathrm{full}}\|\geq1-\sqrt{N^\star}\,\omega_N(\EpsFri)\geq1/2$.
\end{proof}
\begin{lemma}[Full-space conic renormalization]
\label{lem:abd-final-beta-full}
Let $Z^\sharp_{\mathrm{full}}$ be the replaced full-space tuple of~\prettyref{lem:abd-nonzero-full}. Under $\sqrt{N^\star}\,\omega_N(\EpsFri)\leq1/2$, set $\beta\triangleq\|Z^\sharp_{\mathrm{full}}\|^{-1}\leq2$ and multiply the objective coefficient $\bar\mu_0$, the lifted physics/plane multipliers $\widehat{\bar\mu}_i$, the control multiplier $\bar\mu_\StackCtrl$, the (zero) band multipliers, and all full friction Clarke normals $n_i^{m\sharp}$ by $\beta$. Then the full-space nontriviality normalization~\prettyref{eq:abd-ns-normalization} holds exactly, every cone and sign condition is preserved, and every dual (adjoint) residual is multiplied by at most two. The primal residual line, including the plane balance $\sum_m [q_i^{m\sharp}]_\varpi+2\EpsDamp\varpi_i^\star$, is a primal quantity and is unaffected; the external damper $\diag(0,2\EpsDamp I)$ is retained exactly in the smooth adjoint and contributes no replaceable normal.
\end{lemma}
\begin{proof}[Proof of~\prettyref{lem:abd-final-beta-full}]
Clarke normal cones and the convex control cone $N_{\mathcal{D}}(u_i')$ are cones, and all sign and complementary-slackness conditions are preserved by multiplication by $\beta>0$. By~\prettyref{lem:abd-nonzero-full}, $\beta=\|Z^\sharp_{\mathrm{full}}\|^{-1}\leq2$, giving exact unit norm and the factor-two bound on every adjoint residual, which is linear in the scaled dual tuple.
\end{proof}

\begin{proof}[Proof of~\prettyref{thm:abd-C-stationarity}]
We assemble the concrete homogeneous Clarke-stationarity certificate for the full formulation~\prettyref{pb:CRDS-ABD-extended} in five steps, running the reduced solve, lifting unnormalized to the full space, and only then normalizing and transferring on the full tuple.

Step~1: reduced solve. The hypotheses inherited from~\prettyref{thm:CITO-convergence} supply its objective, startup, clamp, initialization, and algorithm-input requirements. The standing~\prettyref{ass:abd-mass} supplies the positive-volume and positive-density hypotheses of~\prettyref{lem:mass}, which establishes the physical full-rank premise used by~\prettyref{cor:SF-CRDS}. By~\prettyref{cor:SF-CRDS} and~\prettyref{lem:smooth-DRp}, \prettyref{pb:PSR-CRDS} is a smooth full-rank CRDS with a valid~\prettyref{def:S-D} smoothing at the scale $\mu(\EpsFri)$ of~\prettyref{eq:abd-mu-schedule}, so~\prettyref{thm:CITO-convergence} applies:~\prettyref{alg:CITO} returns a point carrying a finite $\EpsStop$-KKT certificate for the reduced modified problem in the sense of~\prettyref{def:eps-KKT}, with reduced physics multipliers $\mu_{\theta_i}$, control multiplier $\mu_\StackCtrl$, and auxiliary controls $u_i'$, and satisfying the per-step Hamiltonian bound~\prettyref{eq:cito-ham-bound}. Only the aggregate reduced~\prettyref{def:S-D} closeness of~\prettyref{lem:smooth-DRp-Regularity} is used here.

Localization. The Hamiltonian bound from Step~1 and nonnegativity of the reduced potential give $V_R(\theta_i,\EpsFJ)\leq\LHamilton$. Moreover,~\prettyref{lem:trajectory-bound} gives $\|\theta_i\|\leq\LPosBound(H)$ and $\|\dot\theta_i\|\leq\LVelBound(H)$ at every returned node and predecessor, so each previous configuration $\theta_{\gamma(i)}$ lies in $\mathcal K_H$. Reconstruct $\varpi^\star$ as in~\prettyref{lem:varpi-reconstruction}. Comparing its smoothed full dissipation with the value at $\varpi=0$, using nonnegativity of friction, $\psi_\mu(v)\leq\|v\|$, the aggregate load ceiling $\bar\Lambda_{\mathrm{agg}}$ on $\widehat{\mathcal K}_H$, the localized coefficient bound of~\prettyref{lem:pervertex-activation}, and $\|\dot\theta_i\|\leq\LVelBound(H)$, gives
\begin{ResizedEq}
\EpsDamp\|\varpi^\star\|^2\leq\sum_m\alpha^m\psi_\mu(A^m\dot\theta_i)\leq\bar\Lambda_{\mathrm{agg}}\max_m\|A^m\|\LVelBound(H).
\end{ResizedEq}
Thus the reconstructed damping-active velocities $\dot\vartheta_i^{D\star}$ obey one fixed bound $\|\dot\vartheta_i^{D\star}\|\leq R_{\dot\vartheta^D}<\infty$, uniform in $\EpsFri,\EpsStop$. This proves the configuration and full-velocity localization required by~\prettyref{def:S-D-compatible-ABD} before compatibility is invoked.

Step~2: full-configuration reconstruction and exact unnormalized reduced-to-full-velocity lift. Choose an attained conservative base value $p_0^\star$ at the fixed configuration $\theta_0$, fix finite $\varsigma_0^\star,\omega_0^\star$, and set $\vartheta_0^\star=\langle\theta_0,\{p_0^\star,\varsigma_0^\star,\omega_0^\star\}\rangle$. Then, for every optimized node $i\in\mathcal I$ in increasing order, choose attained conservative minimizers $p_i^\star$ as in~\prettyref{it:active-plane-attain} and use the unique minimizing plane velocity $\varpi_i^\star$ of~\prettyref{lem:varpi-reconstruction} to set
\begin{equation}
\varsigma_i^\star=\varsigma_{\gamma(i)}^\star+\delta_i\dot \varsigma_i^\star,
\qquad
\omega_i^\star=\omega_{\gamma(i)}^\star+\delta_i\dot\omega_i^\star,
\end{equation}
and define $\vartheta_i^\star=\langle\theta_i,\{p_i^\star,\varsigma_i^\star,\omega_i^\star\}\rangle$. At $i=\min\mathcal I$ these recursions are anchored at $\vartheta_0^\star$; no reconstruction crosses the fixed history interval $(\gamma(0),0)$. Its true full velocity is $\dot\vartheta_i^\star=<\dot\theta_i,\dot p_i^\star,\varpi_i^\star>$ with $\dot p_i^\star=(p_i^\star-p_{\gamma(i)}^\star)/\delta_i$, while $\dot\vartheta_i^{D\star}\triangleq\TWO{\dot\theta_i}{\varpi_i^\star}$. This is the trajectory-level application of the one-step construction~\prettyref{eq:nonsmooth-varpi-min} and~\ref{eq:lifted-D-subgradient}: $\varpi_i^\star$ determines auxiliary coordinate increments and is not substituted for the full configuration $\vartheta_i^\star$. The conservative envelope identities~\prettyref{eq:VR-envelope-proof} give
\begin{equation}
\label{eq:abd-plane-config-stationarity}
0=\FPP{V(\vartheta_i^\star,\EpsFJ)}{p_i}
=\FPP{h^i}{p_i}(\vartheta_i^\star)^T\lambda_i^i,
\end{equation}
while the feasibility and complementarity argument in the proof of~\prettyref{lem:eps-SR-FJ} gives the componentwise bands for the single extended map $h^i(\vartheta_i^\star)$. Thus auxiliary-configuration stationarity is exact and contributes no new tolerance.

With these full configurations and the unique minimizing plane velocities reconstructed, lift $\mu_{\theta_i}\mapsto\widehat\mu_i=(\mu_{\theta_i},\nu_i)$ (\prettyref{lem:schur-lift}). Applying~\prettyref{lem:abd-kkt-lift} to that finite reduced certificate gives an exact unnormalized smooth certificate for the reduced-configuration/full-velocity equations: the physical force-balance residual is unchanged, the plane-velocity residual is exactly $0$ at $\varpi^\star$, the physical adjoint equals the reduced one, and the plane-velocity costate residual is exactly $0$. Together with the exact full-configuration reconstruction and plane equation above, these are the required first-order relations before normalization.

Step~3: full-space normalization. Form all smooth full per-vertex normals $N_{i,m}^{\mathrm{sm}}(\widehat\mu_i)$ of~\prettyref{eq:smooth-full-normal}, define the full-space scale $S_{\EpsFri}^{\mathrm{full}}$ of~\prettyref{eq:S-full-scale}, and divide the complete full tuple by it. By~\prettyref{lem:abd-full-normalized} the certificate lies on the unit sphere with $\bar\mu_0=1/S_{\EpsFri}^{\mathrm{full}}\in(0,1]$, and the full lift obeys the $\EpsFri,\EpsStop$-independent bounds~\prettyref{eq:abd-full-bounds}: in particular $\|\bar\mu_{\theta_i}\|\leq1$ and $\bar\mu_{\theta_i}^T K_i\bar\mu_{\theta_i}\leq\bar C_E$, the hypotheses that~\prettyref{def:S-D-compatible-ABD} consumes, obtained directly from the full unit-sphere certificate rather than from the reduced-energy route of the general transfer. Dividing the assembled smooth adjoint relations of Step~2 by $S_{\EpsFri}^{\mathrm{full}}$ bounds the normalized full smooth adjoint residual by $\EpsStop/S_{\EpsFri}^{\mathrm{full}}\leq\EpsStop$.

Step~4: vertexwise full Clarke transfer and renormalization. If $\mathcal{T}=\varnothing$, then $N^\star=0$ and there are no graph points, nearby velocities, friction normals, or normal replacements. The replacement perturbation is exactly zero, so the normalized smooth certificate from the exact lift in Steps~2-3 is already the required full certificate: its physical residual is at most $\EpsStop$, its plane residual is zero, and its full adjoint residual is at most $\EpsStop$. The full-space dual vector is unchanged and has unit norm by~\prettyref{eq:abd-full-unit}, so no further renormalization is needed.

Suppose henceforth that $\mathcal{T}\neq\varnothing$. The localization paragraph above gives $\theta_{\gamma(i)}\in\mathcal K_H$ and $\|\dot\vartheta_i^D\|\leq R_{\dot\vartheta^D}$ at every returned incidence. Apply~\prettyref{lem:abd-normal-compatible} (that is,~\prettyref{def:S-D-compatible-ABD}) with $\xi=\bar\mu_{\theta_i}$ per vertex $m\in\mathcal{T}$: its common lift is $\widehat\xi=\widehat{\bar\mu}_i$ from Step~3, and each normalized smooth full normal $\bar N_{i,m}^{\mathrm{sm}}=N_{i,m}^{\mathrm{sm}}(\widehat{\bar\mu}_i)$ is replaced by a full Clarke normal $n_i^{m\sharp}\in N^C_{\gph\partial_{\dot\vartheta^D}\Damp^{V,m}}$ at a nearby full velocity, with $\|n_i^{m\sharp}-\bar N_{i,m}^{\mathrm{sm}}\|\leq\omega_N(\EpsFri)$ and full-force closeness $\|q_i^{m\sharp}-q_i^{m,\EpsFri}\|\leq C\mu^{1/5}\leq\EpsFri$. Only friction forces and friction normals are replaced; the external damper (the primal term $2\EpsDamp\varpi_i$ in $\widehat F_i^{\mathrm{full}}$, with adjoint contribution $\diag(0,2\EpsDamp I)\widehat{\bar\mu}_i$) stays exact.

Dynamics (primal) line. At every zero-weight incidence, the zero-weight branch of~\prettyref{lem:normal-pullback} gives $q_i^{m\sharp}=q_i^{m,\EpsFri}=0$, so such incidences contribute no force-transfer error. The full residual~\prettyref{eq:abd-full-residual} at the exact forces therefore differs from the Step-2 residual (physical $\leq\EpsStop$, plane $=0$) only through the positive-weight incidences, by $\sum_{i,m:\,\alpha^m(\theta_{\gamma(i)})>0}\|q_i^{m\sharp}-q_i^{m,\EpsFri}\|$ across the physical and plane blocks, so
\begin{ResizedEq}
&\|\widehat F^{\mathrm{full}}(\StackConf,\varpi^\star,\StackCtrl,(q^{m\sharp}))\|_1\\
\leq&\EpsStop+C_F^{\mathrm{full}}N^\star\EpsFri,
\end{ResizedEq}
with $C_F^{\mathrm{full}}\triangleq\max\{1,\sqrt{\max_i\dim\dot\vartheta_i^D}\,\delta_{\max}\}\,C_z$, which is $\EpsFri,\EpsStop$-independent by~\prettyref{eq:norm-mesh}, finite incidence, and the active-incidence coefficient bound of~\prettyref{lem:pervertex-activation}, and which satisfies $C_F^{\mathrm{full}}\ge1$ (hence $C_F^{\mathrm{full}}N^\star\ge1$ as $N^\star\ge1$); the plane balance $\sum_m[q_i^{m\sharp}]_\varpi+2\EpsDamp\varpi_i^\star$ is controlled by the same full-force closeness, the external damper entering exactly. For positive-weight incidences, the nearby-velocity residual $\|\dot\vartheta_i^{D,m\sharp}-\dot\vartheta_i^{D\star}\|\leq C\mu^{1/5}\leq\EpsFri\leq C_F^{\mathrm{full}}N^\star\EpsFri\leq\EpsStopNS$ holds by~\prettyref{lem:vertex-preimage} and the schedule~\prettyref{eq:abd-mu-schedule}; at zero weight choose the nearby velocity equal to $\dot\vartheta_i^{D\star}$, giving zero displacement. These are primal, unaffected by the dual normalization.

Configuration/plane (adjoint) line. Replacing $\bar N_{i,m}^{\mathrm{sm}}$ by $n_i^{m\sharp}$ changes the damping-active adjoint by $\sum_i\sum_m[\nabla Q_i^{\mathrm{full}}]^T(n_i^{m\sharp}-\bar N_{i,m}^{\mathrm{sm}})$; by the banded, finite-incidence structure ($\delta_i\geq\delta_{\min}>0$) the assembled full pullback operator obeys $\|\sum_{i,m}[\nabla Q_i^{\mathrm{full}}]^Tz_i^m\|\leq C_N^{\mathrm{full}}(\sum_{i,m}\|z_i^m\|^2)^{1/2}$ uniformly in $\EpsFri,\EpsStop$, with $C_N^{\mathrm{full}}=C_N^{\mathrm{full}}(\delta_{\min},\text{incidence})$ finite by~\prettyref{eq:sqp-jacobian-structure} (the added $\varpi$-passthrough contributing a bounded block). With~\prettyref{eq:abd-full-bounds} and $\|n_i^{m\sharp}-\bar N_{i,m}^{\mathrm{sm}}\|\leq\omega_N(\EpsFri)$,
\begin{ResizedEq}
&\Big\|\nabla_{(\StackConf,\varpi,q)}\big[\bar\mu_0 O+\textstyle\sum_i\langle\widehat{\bar\mu}_i,\widehat F_i^{\mathrm{full}}\rangle\big]+\textstyle\sum_{i\in\mathcal I}\sum_{m\in\mathcal T}[\nabla Q_i^{\mathrm{full}}]^T n_i^{m\sharp}\Big\|\\
\leq&\EpsStop+C_N^{\mathrm{full}}\sqrt{N^\star}\,\omega_N(\EpsFri).
\end{ResizedEq}
By~\prettyref{lem:abd-nonzero-full} the replaced full dual vector has $\|Z^\sharp_{\mathrm{full}}\|\geq1/2$ (using $\sqrt{N^\star}\,\omega_N(\EpsFri)\leq1/2$); the full conic renormalization $\beta\leq2$ of~\prettyref{lem:abd-final-beta-full} restores~\prettyref{eq:abd-ns-normalization} exactly and multiplies the adjoint residual by at most two, so the configuration/plane line holds at $2\EpsStop+2C_N^{\mathrm{full}}\sqrt{N^\star}\,\omega_N(\EpsFri)$.

Step~5: full bands, auxiliary-configuration stationarity, tolerance, and fixed-data uniform boundedness. The control KKT line of~\prettyref{def:eps-KKT} is unchanged by the exact lift~\prettyref{lem:abd-kkt-lift}; division by $S_{\EpsFri}^{\mathrm{full}}\geq1$ and the final factor $\beta\leq2$ of~\prettyref{lem:abd-final-beta-full} therefore give the control-stationarity bound $\leq2\EpsStop$. The two-point condition $u_i,u_i'\in\mathcal D$, $\|u_i-u_i'\|\leq\EpsStop$ is inherited from the source certificate, and $\mu_{u_i}\in N_{\mathcal D}(u_i')$ is preserved under both positive scalings because the normal cone is a cone. Clamp inactivity supplies a finite penalty level. The lifted feasibility argument of~\prettyref{lem:eps-SR-FJ} gives the equality band, the inequality/complementarity bands for every component of $h^i(\vartheta_i^\star)$. \prettyref{eq:abd-plane-config-stationarity} gives exact auxiliary-configuration stationarity, while the zero band multipliers make complementary slackness immediate. Collecting the decoupled bands, every line of~\prettyref{def:eps-stationary-NS-ABD} holds at the single scalar $\EpsStopNS$ of~\prettyref{eq:abd-tildeeps} with $C_N\triangleq C_N^{\mathrm{full}}$ and $C_F\triangleq C_F^{\mathrm{full}}$ the full-space constants: the configuration/plane line at $2\EpsStop+2C_N^{\mathrm{full}}\sqrt{N^\star}\,\omega_N(\EpsFri)$ and the dynamics line at $\EpsStop+C_F^{\mathrm{full}}N^\star\EpsFri$. When $N^\star=0$, both transfer perturbation terms vanish and the same formula evaluates to $\EpsStopNS=2\EpsStop$. For fixed $\EpsFJ$ and $\EpsFri>0$, common problem and algorithm data and the endpoint clamp choice from Step~8 of~\prettyref{thm:nonsmooth-transfer} give trajectory, prescribed-force, and penalty bounds uniform over $0<\EpsStop\leq\EpsFJ$; no uniformity as $\EpsFJ\downarrow0$ or $\EpsFri\downarrow0$ is asserted. The full dual certificate is normalized to the unit sphere~\prettyref{eq:abd-ns-normalization} by construction, so no bound on the unnormalized multipliers is claimed, and at sticking the full Clarke normals' velocity blocks may be large without affecting any certificate quantity. This proves~\prettyref{thm:abd-C-stationarity}.
\end{proof}

\end{document}